\documentclass[english,11pt,a4paper,reqno,twoside,final]{smfbook}

\usepackage{math-book-dsimon}

 \usepackage{hyperref}

\theoremstyle{plain}
\newtheorem{theo}{Theorem}
\newtheorem{prop}{Proposition}[chapter]
\newtheorem{lemm}{Lemma}[chapter]
\newtheorem{coro}{Corollary}[chapter]
\newtheorem{rema}{Remark}[chapter]
\newtheorem{assumption}{Assumption}[chapter]

\theoremstyle{definition}
\newtheorem{defi}{Definition}[chapter]

\usetikzlibrary{patterns}

\usepackage[inline]{showlabels}

\showlabelsinline

\DeclareMathOperator{\identity}{id}
\DeclareMathOperator{\Guill}{Guill}
\DeclareMathOperator{\PatternShapes}{PS}
\DeclareMathOperator{\BoxShapes}{BS}
\DeclareMathOperator{\PatternTypes}{PT}

\DeclareMathOperator{\Shell}{Sh}

\DeclareMathOperator{\admisscut}{AdmCut}
\DeclareMathOperator{\Ass}{Ass}
\DeclareMathOperator{\Com}{Com}
\DeclareMathOperator{\Mat}{Mat}
\DeclareMathOperator{\Diag}{Diag}

\DeclareMathOperator{\Int}{Int}
\DeclareMathOperator{\vol}{vol}

\DeclareMathOperator{\Wall}{Walls}
\DeclareMathOperator{\Cell}{Cells}
\DeclareMathOperator{\boundaryweights}{bw}
\DeclareMathOperator{\Meas}{Meas}

\newcommand{\Edges}[1]{\mathrm{Edges}(#1)}
\newcommand{\Faces}{\mathrm{Faces}}

\newcommand{\MarkovWeight}[1]{\mathsf{#1}}

\newcommand{\asixv}{\mathsf{a}}
\newcommand{\bsixv}{\mathsf{b}}
\newcommand{\csixv}{\mathsf{c}}
\newcommand{\Lsixv}{\mathsf{L}}
\newcommand{\Msixv}{\mathsf{M}}
\newcommand{\Ssixv}{\mathsf{S}}

\newcommand{\Multind}[2]{\Pi_{#1}^{#2}}
\newcommand{\twodtens}[1]{\bigotimes \begin{bmatrix} #1 \end{bmatrix}}

\newcommand{\patterntype}[1]{\mathrm{\mathbf{#1}}}

\tikzstyle{guillpart} = [scale=0.5,baseline={(current bounding box.center)}]
\tikzstyle{guillfill}=[gray!70!white,opacity=0.35]
\tikzstyle{guillsep} = [ultra thick]
\tikzstyle{centerline} = [baseline={(current bounding box.center)}]

\newcommand{\bulkrect}[3]{ \fill[gray, opacity=0.2] #1 rectangle #2; \draw[guillsep] #1 rectangle #2; \path #1 -- #2 node[midway] {#3}; }

\usepackage{xcolor}
\newcommand{\removable}[1]{}

\title{Operads and the Markov Property on the square lattice.}
\author{Damien \textsc{Simon}}
\address{LPSM (UMR~8001), Sorbonne Université, 4, place Jussieu, F~75252 Paris Cedex 05}
\email{damien.simon@sorbonne-universite.fr}

\date{\today}

\begin{abstract}
	Markov processes on the lattices with arbitrary dimension are omnipresent in statistical mechanics; however their algebraic description is complete only in dimension 1, for which linear algebra provides many tools complementary to the probabilistic approach: as an example, invariant measures are eigenvectors of the generator. In larger dimension, such algebraic tools are absent due to the more involved structure of boundaries. The present work fills this gap by providing a new and complete algebraic description of these models without any other assumption than the Markov property. In order to handle higher dimensions and higher products, the language of operads is used in order to focus on associativities and the geometric interpretation of the algebraic products. This formalism leads to a new parametrization of boundary conditions of Markov processes, with concrete computations. This parametrization is inspired matrix product states in the physics literature. Among others, a probabilistic application elaborated in this paper is the construction of translation-invariant infinite-volume Gibbs measures on the whole lattice by using Kolmogorov's extension; this provides a new alternative tool to the traditional analytical approaches of large size limits. Various models are considered as illustrations.
\end{abstract}

\begin{document}

\frontmatter

\maketitle

\mainmatter

\chapter{Introduction}\label{sec:intro}

\section{Context and overview}

The present work investigates the intrinsic algebraic structure of two-dimensional Markov processes in statistical mechanics and probability theory, introducing the use of operads and matrix product states.

Two-dimensional Markov processes on the lattice encompass most lattice nearest-neighbour models in statistical mechanics, such as the Ising model, discrete Gaussian fields, the six-vertex model, and the Potts models, as well as other frustrated models and discrete versions of statistical field theories. These models are known to have non-trivial phase diagrams and exhibit phase transitions in dimension two, whereas these phenomena are absent in dimension one. Over the last hundred years, these models have been studied using various techniques, with numerous successes. It is, however, remarkable that throughout this vast literature—ranging from the Bethe Ansatz to FKG inequalities—the Markov property of these models is rarely exploited in isolation; it often remains just a basis for intuition. The present work fills this gap by defining an algebraic framework where the Markov property is the elementary building block of computation, irrespective of any other properties the models might possess.

By the Markov property, we mean the property that random variables in different domains separated by curves are conditionally independent given the random variables on the boundaries of the domains. The landscape of literature in dimension two is strikingly different from that in dimension one. In dimension one, Markov processes—and especially their time-oriented variant, Markov chains—have become a classical theory that can be taught at an early stage. What makes the theory so attractive in dimension one is its close relation to linear algebra: probability laws are parametrized by matrices and vectors, marginal laws are computed through matrix products, invariant laws are calculated as Perron-Frobenius eigenvectors, and so on, as explained in any textbook on the subject \cite{Norris}.

Until now, the only direct connection with linear algebra in dimension two has been slicing a domain into strips along a given dimension and manipulating these strips to revert to the one-dimensional case using the formalism of the transfer matrix \cite{baxterbook}. The present work extends this approach by using not just strips but a variety of additional shapes, and by defining a new "linear algebra" language for Markov processes on the two-dimensional lattice $\setZ^2$.

\bigskip

To build an intuitive understanding of the present construction, we invite the reader to reconsider some classical computations in statistical mechanics. The primary goal of statistical mechanics is to describe the large-size behaviour of systems governed by probabilistic local rules. Consider in dimension one (or dimension two) your favourite model on a segment of size $n$ (or a rectangle of size $(p,q)$), your favourite local observable $U$, and your favourite boundary condition. The goal is to extract the limit $l$ of a sequence of expectation values $u_n = \Espi{n}{U}$ (or $u_{p,q} = \Espi{p,q}{U}$) when $U$ is far from the boundaries (in the "bulk"), where the expectation can be written as a large weighted sum over all possible configurations.

We begin with a simple metaphor at the undergraduate level in dimension one. The study of the limit $l$ of such a sequence $(u_n)_{n\in\setN}$ is often performed along one of the following three lines: \begin{itemize} \item A rare case where there is an explicit formula $u_n = f(n)$ simple enough to perform exact asymptotic expansions. \item A frequent case where the sequence $u$ allows a \emph{guess} of the limit $l$ and an analytic control of the remainder $u_n - l$. \item Another frequent case where $(u_n)$ is defined by a recursion $u_{n+1} = f(u_n)$, and the possible limits $l$ must satisfy the fixed-point equation $l = f(l)$, which provides a way of \emph{determining} $l$ (in this case, convergence to a specific $l$ still needs to be proven). \end{itemize}

When we switch to dimension two and sequences of expectations $(u_{p,q})$ associated with large boxes, the correspondence is as follows. The first case corresponds to exactly solvable models, for which exact formulas are obtained for any finite size. The second case corresponds to the traditional analytic approach, in which one guesses suitable approximate boundary conditions, proposes a suitable limit candidate, and then controls the error in the large-size limit using probabilistic inequalities.

The third case is particularly interesting because it provides concrete fixed-point equations to identify the possible limits. In probability theory, the objects of study are much more complex than simple numerical sequences and correspond to joint laws of random variables. The recurrence property can be encoded in terms of consistent marginal laws on subdomains: this is the core idea behind Kolmogorov's extension theorem.

In dimension one, it is straightforward to make the third recursive case concrete for Markov processes. A partition function on a segment of length $n$, or a transition law after $n$ time steps, takes the form of a power $\MarkovWeight{A}^n$ of an elementary matrix describing the interaction between two neighbours. Expectation values of observables localized far from the boundaries involve limits of these powers $\MarkovWeight{A}^n$. Under well-known assumptions (Perron-Frobenius), we have $\MarkovWeight{A}^n \simeq \Lambda^n v_R v_L^*$, but \emph{in practice}, we are interested in the precise values of $\Lambda$, $v_R$, and $v_L$ to approximate the power—not the reverse. The power of linear algebra is that these three objects satisfy simple, easy-to-solve equations $\MarkovWeight{A} v_R = \Lambda v_R$ and $v_L \MarkovWeight{A} = \Lambda v_L$. Returning to probability theory and infinite-volume measures is then made easy by considering boundary equations given by $v_L$ and $v_R$.

This work extends precisely and rigorously this approach to two-dimensional models by introducing a higher-level linear algebra formalism based on operads.

\bigskip

The main challenge in dimension two is the spatially extended nature of boundaries. In dimension one, any segment of arbitrary size $p$ has exactly two endpoints. In dimension two, a rectangle of size $(p,q)$ has a boundary of length $2p + 2q$, which depends on the size. Thus, boundary conditions on different rectangles do not lie in the same space and cannot be directly compared, making it impossible to satisfy any fixed-point property directly. We show in this work that a deeper use of the Markov property and its algebraic counterpart introduces a sub-structure in boundary conditions with elementary blocks, on which the fixed-point approach outlined above can be successfully implemented.

\bigskip

The Markov property is a probabilistic notion. However, the geometric notion of gluing domains along boundaries forms the foundation of the present algebraic approach. In dimension one, gluing two consecutive segments produces a longer segment, but the boundary remains unchanged, consisting of two points. On the probabilistic side, Markov processes encode this geometry in the factorization of conditional expectations with respect to the boundaries. On the algebraic side, this is encoded by the product of matrices associated with the segments. In dimension two, segments are replaced by rectangles. Gluing rectangles to form larger rectangles is more complicated: first, the side sizes along the gluing must coincide, and second, gluing can occur along either of the two directions. On the probabilistic side, the Markov property is still formulated in terms of factorizations of conditional expectations with respect to the boundaries. However, on the algebraic side, the situation is not as clear, and we aim to fill this gap. We show that the framework of coloured operads is the correct framework that encompasses both the structural properties of Markov processes and the necessary computational tools to derive concrete equations.

An interesting observation that motivated us for this work is that, quite often, boundary conditions are not considered as significant in statistical mechanics\footnote{except in special cases such as Dobrushin boundary conditions for spin models or Arctic Circle phenomena for combinatorially constrained models.}, even though boundaries are central to the Markov property. For example, periodic boundary conditions are typically preferred as they avoid the need for a particular choice of values. This contrast likely explains why the Markov property is often not studied in detail, as it would raise the difficulty of describing boundaries.

Most approaches in statistical mechanics describe partition functions of two-dimensional Markov processes on finite domains as "tensors" whose number of indices is given by the length of the boundary. Gluing domains then corresponds to tensor contractions. In this process of growing larger and larger domains, the dimension of tensors raises exponentially and many attempts have been done to control this raise of dimension. For example, transfer matrices with periodic boundary conditions \cite{baxterbook} transforms these tensors into fixed size square matrices to be multiplied. The present work takes seriously this collection of tensors of various sizes on its own through a system of colours in an operad that we call "guillotine operad"; more importantly, we surround it by the more important notion of \emph{associativity} in the operadic framework, and the associated fruitful notions of actions and modules. The exponential growth of dimension of tensors, although still present, is controlled by an algebraic notion of equalities \emph{up to morphisms}. 

\bigskip

The connection between the previous elementary sub-structure of boundary conditions on the probabilistic side and the associativity-based algebraic construction around tensors is made concrete by the use of so-called "matrix product states" (MPS), frequently used in theoretical or numerical physics \cite{MPSreview}. One of the central results of the present work is the description of a deep relation between Markov processes and matrix product states through the formalism of operads. This idea is illustrated below in Section~\ref{sec:nutshell} just below. 

On the probabilistic side, our construction fully embodies the construction of infinite-volume Gibbs measures for Markov processes through Kolmogorov's extension theorem and the asymptotic computation of the thermodynamical free energy density. In \cite{GibbsVelenik}, section 6.2, a short discussion is made about the difficulty of using this extension theorem. These problems are solved in the present paper in the main theorem, Theorem~\ref{theo:eigenROPErep:invmeas}. The MPS sub-structure introduced for boundary weights transpose the consistency constraints of the extension theorem to eigenvalue-like equations on a finite set of objects that define the MPS. These MPS objects replace in dimension two the one-dimensional building blocks $\Lambda$, $v_R$ and $v_R$ required in the study the large size powers $\MarkovWeight{A}^n$ and we show that they exhibit similar properties.

On the algebraic side, the formalism of "guillotine" operads developed here is reminiscent of other operads, such the little-discs or little-cubes $E_n$-operad \cite{mayoperad} and their boundary versions such as the Swiss-cheese operad \cite{VoronovSwissCheese}, the commutative and associative operads as well as double groupoids (\cite{darrick1,darrick2} with similar topological gluings as us). This is not a surprise since the discrete models of statistical mechanics are related to quantum field theory (topological or not). Two major changes of perspective are introduced. First, the lattice structure and the geometrical shapes, which control the number of indices of the tensors, impose colours \cite{YauColoredOperads} in the guillotine operad. These colours forbid, among others, phenomena such as Eckman-Hilton argument or hidden commutativity, which cannot hold for Markov processes. The second major change is the development of boundary "module"-like structures for an operad. It has been a striking fact for us that such a "module" structure has not been elaborated before.
Indeed, one of the goal of operads is the generalization of commutative, associative and Lie structures to "higher" products; however, notions of (bi)modules, spectra and reduction of endomorphisms at the heart of most commutative and associative structure do not have a suitable counterpart at the operadic level. The counterpart for concrete computations is that standard definitions of eigen-elements are still valid in a generalized sense where equalities only hold "up to morphisms", hence meeting other theories in higher algebra.

The present work is the first part of a trilogy of works and contains all the theoretical structure. The second and third parts \cite{BodiotSimon} and \cite{SimonSixV} are applications of all the theorems presented here to two of the most paradigmatic models of statistical mechanics: Gaussian fields on the lattice and the six-vertex model. Previews are presented in Chapter~\ref{sec:applications}. 

For translation-invariant Gaussian fields, rigorous solutions already exist by Fourier transform, of course. The present construction provides an alternative solution, also rigorous, where all the new definitions presented here can be realized concretely and illustrated for beginners; it also highlights that the concrete computations introduced here provides simple recursions, that can be solved either numerically or exactly. 

The six-vertex model, described for example in \cite{baxterbook}, is an exactly solvable model in the sense that, for all finite size with periodic boundary conditions, exact and explicit formulae can be written for all observables, usually using Bethe Ansatz. The study of the thermodynamic limit however suffers from two main difficulties, \emph{even if the limit itself is often already known}: first, for large $N$, the explicit formulae can be (very) difficult to simplify even in the asymptotic regime and, second, a very technical analytical step \cite{DuminilBethe,DuminilCondensation} of "condensation of the Bethe roots" is required. We summarize in Section~\ref{sec:appli:sixvertex} how the present construction realized in details in the third work \cite{SimonSixV} fully circumvents these two difficulties.

	\section{The theory in a nutshell}\label{sec:nutshell}

\subsection{Markov property, tensors, shapes and gluings}
In dimension one, we consider random variables $X_v$ on a subset $D=\{K,K+1,\ldots,K+p\}$ of consecutive integers of $\setZ$ with values in a finite set $S$ that satisfy translational invariance in law and the Markov property. The joint law is then given by
\[
\probi{D}{ (X_v)_{v\in D} = (x_v)_{v\in D} } = \frac{1}{Z_D^{(1D)}} g_l(x_{K}) g_r(x_{L}) \prod_{k=K}^{K+p-1} \MarkovWeight{A}_{1D}( x_k,x_{k+1} )
\]
There is a direct correspondence between the coupling weights and matrices. The partition function and marginal laws are obtained as products of linear matrices and actions on vectors. In particular, one has 
\begin{equation}
	\label{eq:intro:1DZ}
Z_D^{(1D)}=\scal{g_l}{\MarkovWeight{A}_{1D}^{p} g_r}.
\end{equation}
As announced before, linear algebra provides an asymptotic behaviour for large domains with $p\to \infty$ in terms of Perron-Frobenius eigenvectors (as well as spectral gaps using other eigenvectors).

In dimension two, we consider random variables $X_e$ on the edges of $\setZ^2$ taking values in a finite set $S$ and satisfying both spatial Markov property and translational invariance in law. Their joint law on the rectangle $R$ is then given a product of face weights
\[
\probi{R}{ (X_e)_{e\in \Edges{R}}=(x_e)_{e\in \Edges{R}} } = \frac{1}{Z_R} g_{R}((x_e)_{e\in\partial R}) \prod_{f \in \Faces(R) } \MarkovWeight{W}((x_e)_{e\in \partial f}) 
\]
Each face $[k,k+1]\times[l,l+1]$ is an elementary square and has four edges on its boundary. Thus, the Markov weight $\MarkovWeight{W}$ is a function
\[
\MarkovWeight{W}: S^4 \to\setR_+,
\] 
which can be represented as a tensor with four indices, one for each edge. The boundary weight $g_R$ describes the boundary condition and is a function
\[
g_R : S^{2p+2q} \to \setR_+
\]
where $p$ and $q$ are the horizontal and vertical sizes of the rectangle $R$. Precise probabilistic constructions are made in Section~\ref{sec:proba}. The partition function, as well as expectation values, can be written as:
\begin{equation}
	\label{eq:intro:2DZ}
Z_R = \sum_{(x_e)_{e\in\partial R}} g_R((x_e)_{e\in\partial R}) \MarkovWeight{W}^{[p,q]} \left( (x_e)_{e\in\partial R} \right)
\end{equation}
where $\MarkovWeight{W}^{[p,q]}$ is obtained by summing the product of weights $\MarkovWeight{W}$ over all possible values on the internal edges.
This expression replaces in dimension two the previous 1D expression \eqref{eq:intro:1DZ}. There are striking similarities and differences:
\begin{enumerate}
	\item on the probabilistic side, there is no distinction between dimensions one and two: marginalization is obtained by summing over values on boundary points (1D) or edges (2D). On the algebraic side, a weak point of view is to interpret weights and summations as tensors and tensor contraction. 
	
	\item in dimension one, for any value of the size $p$, the number of indices of $\MarkovWeight{A}^p$ is always \emph{two} and the number of indices of each of the \emph{two} factors of the boundary condition is always \emph{one}. In dimension two, the term $\MarkovWeight{W}^{[p,q]}$ as well as the boundary weight $g_R$ both have $2p+2q$ indices and thus depend on the rectangular size $(p,q)$. Moreover, the term $g_R$ is associated to a 1D boundary curve and is not made of two disjoint points any more.
 	
 	\item in dimension one, gluing consecutive segments builds $\MarkovWeight{A}^{p}$ out of $\MarkovWeight{A}^{p'}$ and $\MarkovWeight{A}^{p-p'}$ by a matricial product. In dimension two, gluing compatible rectangles horizontally (resp. vertically) builds $\MarkovWeight{W}^{[p,q]}$ out of $\MarkovWeight{W}^{[p',q]}$ (resp. $\MarkovWeight{W}^{[p,q']}$) and $\MarkovWeight{W}^{[p-p',q]}$ (resp. $\MarkovWeight{W}^{[p,q-q']}$), see Figure~\ref{fig:markovonthelattice}.
 	
	\item in dimension one, there is a natural action of $\MarkovWeight{A}$ on $g_l$ and $g_r$ and a notion of eigenvectors.
\end{enumerate}
More precisely, in point 2, we have $\MarkovWeight{A}^{p} \in \End(\setR^{|S|})$ and $\MarkovWeight{W}^{[p,q]} \in \ca{T}_{p,q}=\End(\setR^{|S|})^{\otimes p}\otimes \End(\setR^{|S|})^{\otimes q}$. The products in point $3$ in dimension two is a combination of tensor products and matrix products. Describing the spaces $\ca{T}_{p,q}$ and their products is made in Section~\ref{sec:operad} with the introduction of the guillotine operad to describe products $\ca{T}_{p',q}\otimes \ca{T}_{p-p',q} \to \ca{T}_{p,q}$ and $\ca{T}_{p,q'}\otimes \ca{T}_{p,q-q'} \to \ca{T}_{p,q}$. Number of indices of tensors are dealt with through the shape labelling. Focus is made on associativities in order to move forward to action on boundaries.

This is the content of Theorems~\ref{theo:canonicalexampleGuill} and \ref{theo:partitionfuncguillotop}. Theorem~\ref{theo:1D:removingcoloursonboundaries} is a technical tool to identify the present construction on size-labelled spaces to the standard linear algebra formalism in dimension $1$.

\subsection{Decomposing and recomposing boundaries}

A way of formulating the final algebraic technical purpose of the paper is to rewrite \eqref{eq:intro:2DZ} in the form of \eqref{eq:intro:1DZ}, with actions and scalar products. There are two approaches in contrary motions to convince oneself that invariant boundary weights $g_R$ must naturally acquire an internal structure given by matrix product states. The present paper essentially establishes that all the following interpretations nicely coincide.

From a probabilistic perspective, one can consider the following configuration of rectangles in $\setZ^2$:
\[
\begin{tikzpicture}[guillpart,yscale=1,xscale=1]
	\draw[thick] (0,0) rectangle (4,3);
	\draw[thick,dotted] (-2,1.5) rectangle  (1,3.5);
	\draw[thick,dotted] (2.,-1) rectangle (3.5,1.25);
	\draw[ultra thick] (0,1.5)--(0,3)--(1,3);
	\draw[ultra thick] (2,0)--(3.5,0);
	\node at (4,2.5) [anchor = west] {$R$};
	\node at (-1.5,1.5) [anchor = north] {$R'$};
	\node at (3.5,-0.75) [anchor = west] {$R''$};
\end{tikzpicture}
\]
Using repeatedly marginalizations, the Markov Property on the various rectangles and their intersections and unions in various orders, one gets quickly convinced that the boundary weight $g_R$ on $R$ can be obtained in various ways and must be made by "gluing" four parts: the two parts inside $R$ and $R''$ and the two parts joining $R'$ and $R''$ through the bottom left vertex or the top right vertex. Repeating this procedure leads intuitively to a representation 
\begin{equation}\label{eq:boundaryweight:quick}
	\begin{split}
		g_R&\left( 
		\begin{tikzpicture}[yscale=0.3,xscale=0.5,baseline={(current bounding box.center)}]
			\draw (0,0) rectangle (4,3);
			\node at (0.5,0) [anchor = north] {{\footnotesize $x_1$}};
			\node at (1.5,0) [anchor = north]{{\footnotesize $x_2$}};
			\node at (2.5,0) [anchor = north]{{\footnotesize $\ldots$}};
			\node at (3.5,0) [anchor = north]{{\footnotesize $x_p$}};
			\node at (0.5,3) [anchor = south] {{\footnotesize $y_1$}};
			\node at (1.5,3) [anchor = south]{{\footnotesize $y_2$}};
			\node at (2.5,3) [above]{{\footnotesize $\ldots$}};
			\node at (3.5,3) [above]{{\footnotesize $y_p$}};
			\node at (0,0.5) [anchor = east]{{\footnotesize  $w_1$ }};
			\node at (0,1.5) [anchor = east]{{\footnotesize  $\vdots$ }};
			\node at (0,2.5) [anchor = east]{{\footnotesize  $w_q$ }};
			\node at (4,0.5) [anchor = west]{{\footnotesize  $z_1$ }};
			\node at (4,1.5) [anchor = west]{{\footnotesize  $\vdots$ }};
			\node at (4,2.5) [anchor = west]{{\footnotesize  $z_q$ }};
		\end{tikzpicture}
		\right) \\
		&= 
		"\Tr_{\ca{W}}"(U_{WS} B_S(x_1)\ldots B_S(x_p)U_{SE} B_E(z_1)\ldots B_E(z_q)U_{EN} B_N(y_p)\ldots B_W(w_1))
	\end{split}
\end{equation}
where the boundary configuration is read counter-clockwise and, at each value, a multiplication by an operator is performed. This "holonomy" is built out of $4+4|S|$ elements (4 change of directions with operator $U_{ab}$ and $|S|$ operators on each of the four sides). Such a representation is known in the literature as a "matrix product space" \cite{MPSreview}. These $4+4|S|$ sub-objects must still be related to the "tensor" weights $\MarkovWeight{W}$: guidance is provided by marginalization formulae inherited from probability theory.

\begin{figure}
	\begin{center}
		\begin{tabular}{p{6cm}|p{6cm}}
			\textbf{Dimension one}	& \textbf{Dimension two} \\
			\hline
			R.v. on vertices & R.v. on edges \\
			\hline
			Weights $\MarkovWeight{A}(x_L,x_R)$ on edges as matrices & Weights $\MarkovWeight{W}(x_S,x_N,x_W,x_E)$ on faces as four-index tensors \\
			\hline
			\begin{tikzpicture}[guillpart,yscale=2,xscale=2]
				\draw[guillsep] (0,0)--(1,0);
				\node at (0,0) [circle, fill, inner sep =2pt] {};
				\node at (1,0) [circle, fill, inner sep =2pt] {};
				\node at (0.5,0) {$\MarkovWeight{A}$};
			\end{tikzpicture}
			&
			\begin{tikzpicture}[guillpart,yscale=0.5,xscale=0.5]
				\fill[guillfill] (0,0) rectangle (2,2);
				\draw[guillsep] (0,0) rectangle (2,2);
				\node at (1,1) {$\MarkovWeight{W}$};
			\end{tikzpicture}
			\\
			\hline
			Unique 1D \emph{free} gluing &  Two 2D gluings of \emph{compatible} rectangles
			\\
			\begin{tikzpicture}[guillpart,yscale=2,xscale=2]
				\draw[guillsep] (0,0)--(2.5,0);
				\node at (0.5,0) [below] {$p'$};
				\node at (1.75,0) [below] {$p-p'$};
				\node at (0,0) [circle, fill, inner sep =2pt] {};
				\node at (1,0) [circle, fill, inner sep =2pt] {};
				\node at (2.5,0) [circle, fill, inner sep =2pt] {};
			\end{tikzpicture}
			&
			\begin{tikzpicture}[guillpart,yscale=1,xscale=1]
				\fill[guillfill] (0,0) rectangle (3,2);
				\draw[guillsep] (0,0) rectangle (3,2);
				\draw[guillsep] (2,0)--(2,2);
			\end{tikzpicture} and 	
			\begin{tikzpicture}[guillpart,yscale=1,xscale=1]
				\fill[guillfill] (0,0) rectangle (3,2);
				\draw[guillsep] (0,0) rectangle (3,2);
				\draw[guillsep] (0,1)--(3,1);
			\end{tikzpicture}
			\\
			&
			\\
			\hline
			One gluing associativity & $3=2\times 1+1$ gluing associativities 
			\\
			\begin{tikzpicture}[guillpart,yscale=2,xscale=2]
				\draw[guillsep] (0,0)--(4,0);
				\node at (0.5,0) [below] {$p_1$};
				\node at (1.75,0) [below] {$p_2$};
				\node at (3.25,0) [below] {$p_3$};
				\node at (0,0) [circle, fill, inner sep =2pt] {};
				\node at (1,0) [circle, fill, inner sep =2pt] {};
				\node at (2.5,0) [circle, fill, inner sep =2pt] {};
				\node at (4.,0) [circle, fill, inner sep =2pt] {};
			\end{tikzpicture}
			& 
			\begin{tikzpicture}[guillpart,yscale=1,xscale=1]
				\fill[guillfill] (0,0) rectangle (3,2);
				\draw[guillsep] (0,0) rectangle (3,2);
				\draw[guillsep] (0.9,0)--(0.9,2);
				\draw[guillsep] (2.2,0)--(2.2,2);
			\end{tikzpicture},
			\begin{tikzpicture}[guillpart,yscale=1,xscale=1]
				\fill[guillfill] (0,0) rectangle (3,2);
				\draw[guillsep] (0,0) rectangle (3,2);
				\draw[guillsep] (0,1.5)--(3,1.5);
				\draw[guillsep] (0,0.75)--(3,0.75);
			\end{tikzpicture}
			and 
			\begin{tikzpicture}[guillpart,yscale=1,xscale=1]
				\fill[guillfill] (0,0) rectangle (3,2);
				\draw[guillsep] (0,0) rectangle (3,2);
				\draw[guillsep] (1.8,0)--(1.8,2);
				\draw[guillsep] (0,1)--(3,1);
			\end{tikzpicture}
			\\
			&
			\\
			\hline
			Powers of matrices & "Surface" powers 
			\\
			$\MarkovWeight{A}^n = 
			\begin{tikzpicture}[guillpart,yscale=1.2,xscale=1.2]
				\draw[guillsep] (0,0)--(4,0);
				\node at (0.5,0) [] {$\MarkovWeight{A}$};
				\node at (1.5,0) [] {$\MarkovWeight{A}$};
				\node at (2.5,0) [] {$\ldots$};
				\node at (3.5,0) [] {$\MarkovWeight{A}$};
				\node at (0,0) [circle, fill, inner sep =2pt] {};
				\node at (1,0) [circle, fill, inner sep =2pt] {};
				\node at (2,0) [circle, fill, inner sep =2pt] {};
				\node at (3,0) [circle, fill, inner sep =2pt] {};
				\node at (4,0) [circle, fill, inner sep =2pt] {};
			\end{tikzpicture}
			$
			&
			$\MarkovWeight{W}^{[p,q]} =
			\begin{tikzpicture}[guillpart,yscale=1,xscale=1]
				\fill[guillfill] (0,0) rectangle (4,3);
				\draw[guillsep] (0,0)--(4,0) (0,1)--(4,1) (0,2)--(4,2) (0,3)--(4,3);
				\draw[guillsep] (0,0)--(0,3) (1,0)--(1,3) (2,0)--(2,3) (3,0)--(3,3) (4,0)--(4,3);
				\node at (0.5,0.5) {{\footnotesize$\MarkovWeight{W}$}};
				\node at (1.5,0.5) {{\footnotesize$\MarkovWeight{W}$}};
				\node at (2.5,0.5) {{\footnotesize$\ldots$}};
				\node at (3.5,0.5) {{\footnotesize$\MarkovWeight{W}$}};
				\node at (0.5,1.5) {{\footnotesize$\MarkovWeight{W}$}};
				\node at (1.5,1.5) {{\footnotesize$\vdots$}};
				\node at (2.5,1.5) {{\footnotesize$\ddots$}};
				\node at (3.5,1.5) {{\footnotesize$\MarkovWeight{\vdots}$}};
				\node at (0.5,2.5) {{\footnotesize$\MarkovWeight{W}$}};
				\node at (1.5,2.5) {{\footnotesize$\MarkovWeight{W}$}};
				\node at (2.5,2.5) {{\footnotesize$\ldots$}};
				\node at (3.5,2.5) {{\footnotesize$\MarkovWeight{W}$}};
			\end{tikzpicture}
			$
			\\
			\hline
		\end{tabular}
	\end{center}	
	\caption{\label{fig:markovonthelattice}Building blocks of Markov processes in dimensions one and two and their geometrical interpretations and representations.}
\end{figure}

From an algebraic perspective, one must find on which type of structure the four-leg tensors $\MarkovWeight{W}$ may act while preserving associativities. In dimension one,
the matrix $\MarkovWeight{A}_{1D}$ has two indices and each of them is used to act on a row or column vector $g_l$ or $g_r$ on the left or on the right. There is still one index left and the result $\MarkovWeight{A}_{1D}g_r$ is still a one-index tensor, i.e. a vector, as $g_r$ itself. In dimension two, contracting the four-index tensor $\MarkovWeight{W}$ with a vector leads to a three-index tensor and things get worse with tensors $\MarkovWeight{W}^{[p,q]}$. In fact, this simple computation does not respect the 2D lattice structure and would be true for any tensor. The 2D lattice structure can be introduced in a shape-preserving way so that the guillotine operad can still be used.

In dimension $1$, matrices are attached to segments. Left and right eigenvectors can be attached to half-lines and scalars to the full line. Indeed, gluing a segment on half-line produces again a half-line in the same way as a matrix acting on a vector produces a vectors. The equation $v(i)=\sum_{j} \MarkovWeight{A}(i,j)u(j)$ is geometrically encoded by:
\[
\begin{tikzpicture}[guillpart,yscale=1,xscale=1]
	\node (A)at (0,0) [circle,fill,inner sep=2pt] {};
	\node at (A) [below] {{\small $i$}};
	\draw (A) -- node [midway, above] {$v$} (3,0); 
\end{tikzpicture}
=\begin{tikzpicture}[guillpart,yscale=1,xscale=1]
	\node (A) at (0,0) [circle,fill,inner sep=2pt] {};
	\node at (A) [below] {{\small $i$}};
	\draw (A) -- node [midway, above] {$\MarkovWeight{A}$} (1.5,0); 
	\node at (1.5,0) [circle,fill, inner sep= 1pt] {};
	\node at (1.5,0) [below] {{\footnotesize $(j)$}};
	\draw (1.5,0)-- node [midway, above] {$u$} (4,0); 
\end{tikzpicture}
\]
Scalar product corresponds the gluing of opposite half-lines. This is made rigorous through $\Guill_1$-structure in the present paper and Theorem~\ref{theo:1D:removingcoloursonboundaries}. In dimension two, a similar geometrical reasoning leads to half-strips and an a priori enigmatic expression
\begin{equation}\label{eq:intro:enigmatichs}
\begin{tikzpicture}[guillpart,yscale=2,xscale=1]
	\fill[guillfill] (0,0) rectangle (3,1);
	\draw[guillsep] (3,0)--(0,0)--(0,1)--(3,1);
	\node at (1.5,0.5) {$B'$};
	\node at (0,0.5) [left] {{\small $i$}};
	\node at (1.5,0) {{\small $(k,?)$}}; \node at (1.5,1) {{\small $(l,?')$}};
\end{tikzpicture}
=
\begin{tikzpicture}[guillpart,yscale=2,xscale=1.5]
	\fill[guillfill] (0,0) rectangle (4,1);
	\draw[guillsep] (4,0)--(0,0)--(0,1)--(4,1);
	\draw[guillsep] (1.5,0)--(1.5,1); \node at (1.5,0.5) {{\footnotesize $(j)$}};
	\node at (2.5,0.5) {$B$};
	\node at (0.75,0.5) {$\MarkovWeight{W}$};
	\node at (0,0.5) [left] {{\small $i$}};
	\node at (2.5,0) {{\small $?$}}; \node at (2.5,1) {{\small $?'$}};
	\node at (0.75,0) {{\small $k$}}; \node at (0.5,1) {{\small $l$}};
\end{tikzpicture}
\end{equation}
encoding the action of one of the index of $\MarkovWeight{W}$ on an element $B$ on its right. What should replaced the question marks will be discussed later.

The half-strip vision above is a direct extension of the 1D picture to dimension two. However, dimension two is a bit more than twice the dimension one, as it can be seen for example in the square associativity in the guillotine operad. Geometric considerations easily convince oneself that corner elements depicted as:
\[
\begin{tikzpicture}[guillpart,yscale=0.8,xscale=0.8]
	\fill[guillfill] (0,0) rectangle (2,2);
	\draw[guillsep] (2,0)--(0,0)--(0,2);
	\node at (1,1) {$U$};
\end{tikzpicture}
\]
can be introduced and glued with the half-strip element $B$ as their top half-strip equivalents. One already remarks that, contrary to dimension one, the initial weight $\MarkovWeight{W}$ cannot be glued on a corner element...

The guillotine operad of Section~\ref{sec:operad} gives a precise meaning to such diagrams through Theorems~\ref{theo:extendedguilloperads}, \ref{theo:operadpointedversion} and \ref{theo:canonicalboundarystructure}. This leads to the important Theorem~\ref{theo:stability} that concentrates this discussion in the diagrammatic equation:
\begin{equation}\label{eq:intro:quickROPErep}
		g_R\left( 
		\begin{tikzpicture}[yscale=0.3,xscale=0.5,baseline={(current bounding box.center)}]
			\draw (0,0) rectangle (4,3);
			\node at (0.5,0) [anchor = north] {{\footnotesize $x_1$}};
			\node at (1.5,0) [anchor = north]{{\footnotesize $x_2$}};
			\node at (2.5,0) [anchor = north]{{\footnotesize $\ldots$}};
			\node at (3.5,0) [anchor = north]{{\footnotesize $x_p$}};
			\node at (0.5,3) [anchor = south] {{\footnotesize $y_1$}};
			\node at (1.5,3) [anchor = south]{{\footnotesize $y_2$}};
			\node at (2.5,3) [above]{{\footnotesize $\ldots$}};
			\node at (3.5,3) [above]{{\footnotesize $y_p$}};
			\node at (0,0.5) [anchor = east]{{\footnotesize  $w_1$ }};
			\node at (0,1.5) [anchor = east]{{\footnotesize  $\vdots$ }};
			\node at (0,2.5) [anchor = east]{{\footnotesize  $w_q$ }};
			\node at (4,0.5) [anchor = west]{{\footnotesize  $z_1$ }};
			\node at (4,1.5) [anchor = west]{{\footnotesize  $\vdots$ }};
			\node at (4,2.5) [anchor = west]{{\footnotesize  $z_q$ }};
		\end{tikzpicture}
		\right) 
		= \begin{tikzpicture}[guillpart,yscale=1.2,xscale=3.]
			\draw[guillsep, dotted] (0,0) rectangle (5,5);
			\draw[guillsep] (1,1) rectangle (4,4);
			\draw[guillsep] 	(1,0)--(1,1)
			(2,0)--(2,1)
			(3,0)--(3,1)
			(4,0)--(4,1);
			\draw[guillsep] 	(1,5)--(1,4)
			(2,5)--(2,4)
			(3,5)--(3,4)
			(4,5)--(4,4);
			\draw[guillsep] 	(0,1)--(1,1)
			(0,2)--(1,2)
			(0,3)--(1,3)
			(0,4)--(1,4);
			\draw[guillsep] 	(5,1)--(4,1)
			(5,2)--(4,2)
			(5,3)--(4,3)
			(5,4)--(4,4);
			\node at (0.5,0.5) { $U_{SW}$ };
			\node at (4.5,0.5) { $U_{SE}$ };
			\node at (0.5,4.5) { $U_{NW}$ };
			\node at (4.5,4.5) { $U_{NE}$ };
			\node at (1.5,0.5) { $B_S(x_1)$ };
			\node at (2.5,0.5) { $\ldots$ };
			\node at (3.5,0.5) { $B_S(x_P)$ };
			\node at (1.5,4.5) { $B_N(y_1)$ };
			\node at (2.5,4.5) { $\ldots$ };
			\node at (3.5,4.5) { $B_N(y_P)$ };
			\node at (0.5,1.5) { $B_W(w_1)$ };
			\node at (0.5,2.5) { $\vdots$ };
			\node at (0.5,3.5) { $B_W(w_Q)$ };
			\node at (4.5,1.5) { $B_E(z_1)$ };
			\node at (4.5,2.5) { $\vdots$ };
			\node at (4.5,3.5) { $B_E(z_Q)$ };
		\end{tikzpicture}
\end{equation}
which is the operadic formulation of \eqref{eq:boundaryweight:quick} and its properties. Such an operadic MPS representation (a ROPErep in our language) encodes simultaneously the boundary nature of $g_R$, its internal structure dictated the Markov property as well as a probabilistic stability property: if such a weight $g_R$ is chosen for a rectangle $R$ then all the marginals on small rectangles of the Markov process have boundary weights with a similar ROPErep structure.

\begin{figure}
\begin{center}
\begin{tabular}{p{6cm}|p{6cm}}	
\textbf{Dimension one} & \textbf{Dimension two}
\\
\hline
Boundary conditions $g(x,x')$ on the two end of a segment & Boundary conditions $g_R((x_e)_{e\in\partial R})$ on a whole closed curve.
\\
\hline
Left-right factorization $g(x,x')=g_l(x)g_r(x')$
&
ROPErep \eqref{eq:boundaryweight:quick} and \eqref{eq:intro:quickROPErep}
\\
\hline
Half-line interpretation & Half-strip interpretation 
\\
$\begin{tikzpicture}[guillpart,yscale=1,xscale=1.5]
	\draw[guillsep] (0,0)--(1,0);
	\draw[guillsep,dashed] (-0.5,0)--(0,0);
	\node at (1,0) [circle,fill,inner sep=2pt] {};
	\node at (0.5,0) [above] {$g_l$}; 
%	\node at (1,0) [anchor=west] {$x$};
\end{tikzpicture}$
and 
$\begin{tikzpicture}[guillpart,yscale=1,xscale=-1.5]
	\draw[guillsep] (0,0)--(1,0);
	\draw[guillsep,dashed] (-0.5,0)--(0,0);
	\node at (1,0) [circle,fill,inner sep=2pt] {};
	\node at (0.5,0) [above] {$g_r$}; 
%	\node at (1,0) [anchor=east] {$x'$};
\end{tikzpicture}$
& 
$\begin{tikzpicture}[guillpart,scale=1.5]
	\fill[guillfill] (0,0) rectangle (1,1);
	\draw[guillsep] (0,0)--(1,0)--(1,1)--(0,1);
	\node at (0.5,0.5) {$B_W$};
\end{tikzpicture}$
, 
$\begin{tikzpicture}[guillpart,scale=1.5,rotate=90]
	\fill[guillfill] (0,0) rectangle (1,1);
	\draw[guillsep] (0,0)--(1,0)--(1,1)--(0,1);
	\node at (0.5,0.5) {$B_S$};
\end{tikzpicture}$
,
$\begin{tikzpicture}[guillpart,scale=1.5,rotate=-90]
	\fill[guillfill] (0,0) rectangle (1,1);
	\draw[guillsep] (0,0)--(1,0)--(1,1)--(0,1);
	\node at (0.5,0.5) {$B_N$};
\end{tikzpicture}$
,
$\begin{tikzpicture}[guillpart,scale=1.5,rotate=-180]
	\fill[guillfill] (0,0) rectangle (1,1);
	\draw[guillsep] (0,0)--(1,0)--(1,1)--(0,1);
	\node at (0.5,0.5) {$B_E$};
\end{tikzpicture}$
\\
\hline
Left and right actions & Four actions
\\
$Ag_r = \begin{tikzpicture}[guillpart,yscale=1,xscale=-1.5]
	\draw[guillsep] (0,0)--(2,0);
	\draw[guillsep,dashed] (-0.5,0)--(0,0);
	\node at (1,0) [circle,fill,inner sep=1.5pt] {};
	\node at (2,0) [circle,fill,inner sep=2pt] {};
	\node at (0.5,0) [above] {$g_r$};
	\node at (1.5,0) [above] {$\MarkovWeight{A}$};
	%	\node at (1,0) [anchor=west] {$x$};
\end{tikzpicture}$, etc.
&
$\begin{tikzpicture}[guillpart,yscale=1.5,xscale=1.5]
	\fill[guillfill] (0,0) rectangle (2,1);
	\draw[guillsep] (2,0)--(0,0)--(0,1)--(2,1) (1,0)--(1,1);
	\node at (1.5,0.5) {$B_E$};
	\node at (0.5,0.5) {$\MarkovWeight{W}$};
\end{tikzpicture}$, etc.
\\
\hline
& Corners and two actions of half-strips for each one \\
{\footnotesize\emph{No equivalent}}
& 
$\begin{tikzpicture}[guillpart,yscale=1.5,xscale=2]
	\fill[guillfill] (0,0) rectangle (2,1);
	\draw[guillsep] (0,1)--(0,0)--(2,0) (1,0)--(1,1);
	\node at (0.5,0.5) {$B_N$};
	\node at (1.5,0.5) {$U_{NE}$};
\end{tikzpicture}$,
$\begin{tikzpicture}[guillpart,yscale=1.5,xscale=2]
	\fill[guillfill] (0,0) rectangle (2,1);
	\draw[guillsep] (0,0)--(2,0)--(2,1) (1,0)--(1,1);
	\node at (0.5,0.5) {$U_{NW}$};
	\node at (1.5,0.5) {$B_N$};
\end{tikzpicture}$, etc.
\\
\hline
One scalar product & Multiple pairings 
\\
$\begin{tikzpicture}[guillpart,yscale=1,xscale=1.5]
	\draw[guillsep] (0,0)--(2,0);
	\draw[guillsep,dashed] (-0.5,0)--(0,0) (2,0)--(2.5,0);
	\node at (1,0) [circle,fill,inner sep=2pt] {};
	\node at (0.5,0) [above] {$g_l$}; 
	\node at (1.5,0) [above] {$g_r$};
\end{tikzpicture}$
&
$\begin{tikzpicture}[guillpart,yscale=1.5,xscale=2]
	\fill[guillfill] (0,0) rectangle (2,1);
	\draw[guillsep] (0,1)--(2,1) (0,0)--(2,0) (1,0)--(1,1);
	\node at (0.5,0.5) {$B_W$};
	\node at (1.5,0.5) {$B_E$};
\end{tikzpicture}$, $\begin{tikzpicture}[guillpart,yscale=1.5,xscale=2]
\fill[guillfill] (0,0) rectangle (2,1);
\draw[guillsep]  (0,0)--(2,0) (1,0)--(1,1);
\node at (0.5,0.5) {$U_{NW}$};
\node at (1.5,0.5) {$U_{NE}$};
\end{tikzpicture}$, etc.
\\
\hline
\end{tabular}
\end{center}
	\caption{\label{fig:nutshell:boundarystructures} Correspondence of the basic building blocks  for boundary conditions in dimension one and two}
\end{figure}

\subsection{Towards diagonalization and large size limits}

In dimension one, the next step from \eqref{eq:intro:1DZ} towards large size limits consists in diagonalizing the matrix $\MarkovWeight{A}_{1D}$. In dimension two, working directly on boundary weights is hard and does not lead very far. Relating by marginalization two boundary weights on two rectangles that differ horizontally by a size $1$ on the right is given from probability theory by:
\begin{align*}
g_{p,q}\left( 
\begin{tikzpicture}[scale=0.3,baseline={(current bounding box.center)}]
	\draw (0,0) rectangle (4,3);
	\node at (0.5,0) [anchor = north] {{\footnotesize $x_1$}};
	\node at (1.5,0) [anchor = north]{{\footnotesize $x_2$}};
	\node at (2.5,0) [anchor = north]{{\footnotesize $\ldots$}};
	\node at (3.5,0) [anchor = north]{{\footnotesize $x_p$}};
	\node at (0.5,3) [anchor = south] {{\footnotesize $y_1$}};
	\node at (1.5,3) [anchor = south]{{\footnotesize $y_2$}};
	\node at (2.5,3) [above]{{\footnotesize $\ldots$}};
	\node at (3.5,3) [above]{{\footnotesize $y_p$}};
	\node at (0,0.5) [anchor = east]{{\footnotesize  $w_1$ }};
	\node at (0,1.5) [anchor = east]{{\footnotesize  $\vdots$ }};
	\node at (0,2.5) [anchor = east]{{\footnotesize  $w_q$ }};
	\node at (4,0.5) [anchor = west]{{\footnotesize  $z_1$ }};
	\node at (4,1.5) [anchor = west]{{\footnotesize  $\vdots$ }};
	\node at (4,2.5) [anchor = west]{{\footnotesize  $z_q$ }};
\end{tikzpicture}
\right) & = \sum_{(x,y,\gr{z})\in S^{2+q}}
g_{p+1,q}\left( 
\begin{tikzpicture}[scale=0.3,baseline={(current bounding box.center)}]
	\draw (0,0) rectangle (5,3);
	\draw[dotted] (4,0)--(4,3);
	\node at (0.5,0) [anchor = north] {{\footnotesize $x_1$}};
	\node at (1.5,0) [anchor = north]{{\footnotesize $x_2$}};
	\node at (2.5,0) [anchor = north]{{\footnotesize $\ldots$}};
	\node at (3.5,0) [anchor = north]{{\footnotesize $x_p$}};
	\node at (4.5,0) [anchor = north]{{\footnotesize $x$}};
	\node at (0.5,3) [anchor = south] {{\footnotesize $y_1$}};
	\node at (1.5,3) [anchor = south]{{\footnotesize $y_2$}};
	\node at (2.5,3) [above]{{\footnotesize $\ldots$}};
	\node at (3.5,3) [above]{{\footnotesize $y_p$}};
	\node at (4.5,3) [above]{{\footnotesize $y$}};
	\node at (0,0.5) [anchor = east]{{\footnotesize  $w_1$ }};
	\node at (0,1.5) [anchor = east]{{\footnotesize  $\vdots$ }};
	\node at (0,2.5) [anchor = east]{{\footnotesize  $w_q$ }};
	\node at (5,0.5) [anchor = west]{{\footnotesize  $z'_1$ }};
	\node at (5,1.5) [anchor = west]{{\footnotesize  $\vdots$ }};
	\node at (5,2.5) [anchor = west]{{\footnotesize  $z'_q$ }};
\end{tikzpicture}
\right)
T^v_q( \gr{z},\gr{z}' ; x,y )
\\
T^v_q( \gr{w}, \gr{z} | x,x' ) &= \sum_{ \substack{\gr{v}\in S^{q+1} \\ v_0=x, v_{P}=x'} } \prod_{l=1}^Q \MarkovWeight{W}(v_{l-1},v_l,w_{l},z_l)
\end{align*}
where $T_Q^v$ is a vertical transfer matrix, already hard to compute by itself. The difference of perimeters between the two rectangles make impossible any eigenvalue-like relation between $g_{p,q}$ and $g_{p+1,q}$.

On the contrary, \eqref{eq:intro:enigmatichs} naturally leads to an eigenvalue like equation 
\begin{equation}\label{eq:intro:eigenlike}
\Lambda \begin{tikzpicture}[guillpart,yscale=1,xscale=1]
	\fill[guillfill] (0,0) rectangle (2,1);
	\draw[guillsep] (2,0)--(0,0)--(0,1)--(2,1);
	\node at (1.,0.5) {$B$};
%	\node at (0,0.5) [left] {{\small $i$}};
%	\node at (1.5,0) {{\small $(k,?)$}}; \node at (1.5,1) {{\small $(l,?')$}};
\end{tikzpicture}
"\overset{\Phi}{=}"
\begin{tikzpicture}[guillpart,yscale=1,xscale=1]
	\fill[guillfill] (0,0) rectangle (4,1);
	\draw[guillsep] (4,0)--(0,0)--(0,1)--(4,1);
	\draw[guillsep] (1.5,0)--(1.5,1);
	% \node at (1.5,0.5) {{\footnotesize $(j)$}};
	\node at (2.5,0.5) {$B$};
	\node at (0.75,0.5) {$\MarkovWeight{W}$};
%	\node at (0,0.5) [left] {{\small $i$}};
%	\node at (2.5,0) {{\small $?$}}; \node at (2.5,1) {{\small $?'$}};
%	\node at (0.75,0) {{\small $k$}}; \node at (0.5,1) {{\small $l$}};
\end{tikzpicture}
\end{equation}
provided that the half-line structure represented here by question marks is elucidated and introduces suitable identifications between $(l,?')$ and $?'$. Again the structural insight provided by the guillotine operad and its associativities provides the solution to this identification problem on half-lines using "up-to-morphisms" definitions of eigenvectors. 

This technical part and its probabilistic applications corresponds to all the definitions of Chapter~\ref{sec:invariantboundaryelmts}. The three probabilistic consequences are Theorems~\ref{theo:reductionofmorphisms:canostruct} and \ref{theo:eigencorner:realizationMarkov}, which produces a short set of equations to be solved for the elements $B_a$ and $U_{ab}$ above, and the important Theorem~\ref{theo:eigenROPErep:invmeas}, which deduces the construction of infinite-volume translation-invariant Gibbs measures from the generalized eigen-elements.

In \eqref{eq:intro:eigenlike}, no transfer matrix is computed nor any power $\MarkovWeight{W}^{[p,q]}$. The result of such equations for half-strips and their equivalent on corners however provides the eigenvalue $\Lambda$ and the element $B_a$ and $U_{ab}$ which partly control the large size powers $\MarkovWeight{[p,q]}$ in the same way as $\Lambda$ and the eigenvectors $v_L$ and $v_R$ control $\MarkovWeight{A}^p$ in dimension one.

\section{A complete map of the paper}

The present paper is long essentially because many notions had to be defined from scratch. We organize it in chapters by first separating the material coming from different fields ---probability theory, algebra, matrix product states from physics--- and then assembling them step by step. The book is divided into six chapters after this introduction, which contains complements about the Markov property in the literature and alternative approaches, and a final conclusion.

Chapter~\ref{sec:proba} is a probabilistic introduction to Markov processes in dimension one and dimension two, whose main purposes are to fix notations, to describe the geometries we consider, to formulate the \textbf{Markov property} and its consequences and finally to present the non-trivial questions that are usually of interest in the theory of Markov processes and in statistical mechanics, together with their standard approaches. It does not contain any new material. Some computations are presented with a bit more detail than required to help the algebraist reader, who may not be at ease with probability theory.

Chapter~\ref{sec:operad} is algebraic and introduces briefly the concept of coloured operads in Section~\ref{sec:operadiclanguage} and then focuses in Section~\ref{sec:guillotinebasics} on the coloured operads of interest for Markov processes on the square lattice: the \textbf{guillotine operads} introduced in Definition~\ref{def:guillotineoperad:main}. In particular, the structure of the guillotine operad is shown to be generated by elementary so-called guillotine cuts and the three fundamental associativity rules~\eqref{eq:guill2:listassoc}, whose identification enlightens all the subsequent constructions. \textbf{Theorems~\ref{theo:canonicalexampleGuill} and \ref{theo:partitionfuncguillotop}} establish the connection between the partition functions of Markov processes and a suitable algebra over the guillotine operad. From a practical probabilist point of view, these two theorems are just interesting algebraic reformulations of classical computations but, on their own, do not contain any new direct computational material; they only open the way to the next sections.

Section~\ref{sec:boundaryguill2} is the algebraic core of the new constructions of the present work. It \emph{extends} the guillotine operad with \emph{boundary spaces} on the boundaries of the rectangles, in the same way as left and right modules extend associative algebras. From an algebraic perspective, this extension relies on the introduction of suitable colours and actions of the previous guillotine algebras on these additional spaces. The two cornerstones of this construction are to keep the associativity rules~\eqref{eq:guill2:listassoc} during this extension and to understand that boundaries may themselves carry a guillotine operad structure of lower dimension compatible with this extension. This section is quite long since many new objects have to be introduced, one for each shape of Figure~\ref{fig:admissiblepatterns} with their own structure and subtleties. We shortly revisit the dimension one case with Theorem~\ref{theo:1D:removingcoloursonboundaries}. The structures in dimension two case are summarized in \textbf{Theorems~\ref{theo:extendedguilloperads} and \ref{theo:operadpointedversion}}. Again, the reader will be invited to always think about the operations on these objects as gluing of geometric shapes which can be realized by matrix products and tensor products.

All along Section~\ref{sec:operad}, we develop \emph{graphical notations} for the elements of the guillotine operad acting on elements of the corresponding algebras. These notations are not just shorthand notations but deeply reflect the geometric nature of the objects and make computations easier; this should be put at the same level of rigour and simplicity as writing $ab$ or $abc$ the products of elements $m(a,b)$ or $m(a,m(b,c))=m(m(a,b),c)$ in an associative algebra (in dimension one).

The final section~\ref{sec:canonicalboundarystructure} of Section~\ref{sec:operad} presents through \textbf{Theorem~\ref{theo:canonicalboundarystructure}} a canonical boundary structure of extended guillotine algebra for two-dimensional Markov processes. This is a necessary tool for Chapter~\ref{sec:boundaryalgebra} and it also illustrates how basic factorized boundary conditions (which are common in the statistical mechanics literature) fit into this operadic framework and already make necessary a more general class of boundary conditions described in its full generality in the following Chapter~\ref{sec:boundaryalgebra}.

Chapter~\ref{sec:boundaryalgebra} introduces "matrix product states" to describe boundary weights of Markov processes, see \eqref{eq:boundaryweight:quick} above, and provides them with a guillotine operadic structure. We then describe the advantages of such a description: \begin{itemize} 
	\item a global boundary weight (which is associated with a dimension one curve) is built out of elementary objects associated with elementary segments and corners; \item there is an "associative" structure on such boundary elements, corresponding to a guillotine algebra in dimension one; 
	\item canonical actions of a two-dimensional guillotine algebra on such a one-dimensional-like guillotine algebra may be defined. 
\end{itemize} 
All the basic properties of such matrix product representations are first presented in Section~\ref{sec:introROPE} with the fundamental Definitions~\ref{def:ROPE}, \ref{def:ROPErep:FD} and \ref{def:ROPErep:SSS} of \textbf{ROPE and ROPErep}. A major result is \textbf{Theorem~\ref{theo:stability}} of Section~\ref{sec:ROPEstability} which states that ROPE representations of boundary weights as matrix products are \emph{stable} under the action of the face probability weights of the model seen as elements of a two-dimensional guillotine algebra acting on its boundary spaces. From a probabilistic point of view, this corresponds to the stability of ROPE representations of boundary weights when considering marginal laws of Markov processes on smaller rectangles. Up to our knowledge, no such result exists for matrix product states and we believe that such a structural stability may be fundamental and explain many of their observed wonderful computational properties. The proof of the theorem is in particular a constructive proof, so that the actions are explicitly written.

Chapter~\ref{sec:invariantboundaryelmts} contains the key definition of generalized eigen-elements on boundaries, targeting the construction of infinite-volume Gibbs measures on the probability side. It adds to the classical definition of eigenvectors morphisms of guillotine operads. These new eigen-elements then become the local building blocks of invariant boundary weights and lead to infinite-volume Gibbs measures. These definitions are, to the best of our knowledge, new and combine in a very intriguing way various modern algebraic notions: traditional eigenvalue and eigenvectors of linear algebra, coloured operadic products and actions and identification of algebras only up to morphisms. Several paragraphs are devoted to the exploration of various aspects of such a new structure, in particular it contains a Yang-Baxter-like equation which appears in Section~\ref{sec:commutuptomorph} as a higher-analogue of a classical exercise in linear algebra.

As a naive illustration, we present in Section~\ref{sec:trivialfactorizedcase} two simple models, that are trivial from the point of view of statistical mechanics but already exhibit the main non-trivial features of our definitions. For more general models, we present in Theorems~\ref{theo:reductionofmorphisms:canostruct} and \ref{theo:eigencorner:realizationMarkov} the new set of \textbf{concrete equations} ---in a generic language--- which have to be solved to compute the boundary eigen-elements introduced in the definitions. Section~\ref{sec:gibbsfromeigen} reaches the final probabilistic goal of the paper: \textbf{Theorem~\ref{theo:eigenROPErep:invmeas}} shows how \emph{infinite-volume Gibbs measures} arise from boundary eigen-elements up to morphisms.

Chapter~\ref{sec:applications} describes some applications and proofs of concept of the present approach. It provides concrete realizations of all the abstract definitions on various models, less trivial than the ones of Section~\ref{sec:trivialfactorizedcase} and thus shows how this new formalism allows to make concrete computations of statistical mechanics, as illustrated in \textbf{Theorem~\ref{theo:onedimmarginaloutofinvROPErep}}. It contains summaries of the two associated papers \cite{BodiotSimon} (through the summarizing Theorem~\ref{theofromBS:Gaussian}) and \cite{SimonSixV} about Gaussian models and the six-vertex model.

All the framework of the previous sections is presented in a purely algebraic form with finite sets, finite sums and finite-dimensional spaces and only the two-dimensional case in order to have clear and intuitive graphical notations to handle the new notions. This is sufficient for numerous interesting discrete models of statistical mechanics. However, it is interesting to have, on one side, more \textbf{general measurable state spaces} ---already for real-valued Gaussian processes for example--- and, on the other side, to \textbf{raise the spatial dimension of the lattice}. All these generalizations are presented in Chapter~\ref{sec:generalizations}. Theorems~\ref{theo:higherdim:canostruct} and \ref{theo:canonicalexampleGuill:continuous} present part of these generalizations.

Finally, Chapter~\ref{sec:openquestions} is a conclusion that contains various remarks and lists open questions and connections with other topics.

	\section{Remarks on the Markov property}
	
\subsection{About the Markov property in statistical mechanics.}

Most models of statistical mechanics on lattices are defined by local interactions, i.e. their probabilistic weights are made of factors that depend only on neighbouring (in the sense of the lattice) random variables; this corresponds to Gibbs random fields. This is also the case of quantum and statistical fields theories (in discrete or, more often, continuous space), for which most Lagrangians or Hamiltonians are local. From a probabilistic point of view, the consequence of this locality for these models is the Markov property, which says roughly that, conditionally on what happens on boundaries of domains with disjoint interiors, random variables on the interior of these domains are independent. 

In dimension 1, often interpreted as time instead of space when an orientation is chosen, laws of Markov chains are described by stochastic matrices and initial laws seen as left vectors, as already said.

For arbitrary graphs (not necessarily embedded in some $\setR^d$), Hammersley-Clifford's theorem provides a factorization of the law of a Markov process over complete sub-graphs, with an interaction which can be encoded in tensors. However, at this level of generality, one shall not expect for this combinatorial description to be related to nice geometrical and algebraic tools as in dimension 1 (no translation invariance, no notion of asymptotic free energy density, no products excepted tensor contractions etc.).

In dimension two, most models satisfying the Markov property are traditionally studied by using other methods that avoid the use of the Markov property in its \emph{full} generality. We mention, without any order, the following traditional tools relevant to our paper:
\begin{itemize}
\item the transfer matrix method of statistical mechanics \cite{baxterbook}, which splits the dimension $d$ into $(d-1)+1$ and maps the model to a one-dimensional Markov process with a (much) larger state space, either with periodic boundary conditions or arbitrary fixed ones.
\item infinite volume Gibbs measures \cite{GibbsGeorgii,GibbsVelenik}, which help to hide careful description of boundaries behind nice nestings of conditional measures (DLR approach) or through analytical approximation tools and Riesz representation theorem.
\item the theory of integrable systems \cite{baxterbook}, which uses additional features of some models to develop exact finite size formulae based on representation theory of quantum groups.
\end{itemize}

When teaching Markov chains to any non-specialized audience, the traditional way is to start with the discrete space $\setN$ and finite state spaces $S$, which can be done with a weak background both in probability theory and in linear algebra. Most geometric and algebraic features are already present at this level. Generalizing then to measurable state spaces and to the real line $\setR_+$ requires much more subtle measure theory and analytic tools on semi-groups but, excepted for particular properties such as explosion times for Feller diffusions, this generalization does not add new algebraic intuitions. This is why in the present paper, we focus only on \emph{finite state spaces and the square lattice $\setZ^d$}: our algebraic tools are sufficiently general to be then formulated for more general state spaces and continuous spaces provided that analytical tools may enrich the algebraic tools with suitable norms and continuities. This analytic remark may however not be an easy part for two main reasons: first, the theory of semi-groups in dimension 1 is already subtle and, second, one may expect that functions should be replaced by distributions in dimension larger than 2 since this is the standard behaviour of scaling limits of most models of statistical mechanics. We believe then that splitting the study of larger-dimensional Markov processes between an algebraic part in the discrete setting and analytic tools specific to the continuous setting may be a good method.

When considering the historical treatment of boundary conditions in statistical mechanics, the French word "désamour" sounds like a good summary. Hidden behind the credo of universality, a frequent reflex is to consider periodic boundary conditions or trivial ones, in order to reduce the number of parameters and other technicalities. This is nice for many computations but it is contrary to the principle of the Markov property, which makes boundaries appear everywhere by disconnecting domains and eliminating any periodicity. The present paper is essentially based on this fundamental remark and the will to put boundaries back in the light in a natural way. 

In dimension two and larger, boundaries of domains are still connected and made of segments and corners. Such boundaries are extended objects but may still be built as gluing of elementary local objects as in \eqref{eq:boundaryweight:quick}. The leitmotiv of the following sections is the following: all the global objects ---partition functions, free energy densities, invariant boundary conditions, infinite-volume Gibbs measure--- should be built out of local objects related to elementary geometric shapes ---rectangle of size $(1,1)$, segments of length $1$, corners--- whose gluings lead to operadic constructions on the local objects. Actually, such a philosophy pre-exists in statistical mechanics: this is the idea between matrix product states \cite{MPSreview} which has inspired part of the present work.

\subsection{Relations with other approaches}\label{sec:intro:otherapproaches}
	
The present paper considers the Markov property on its own and does not require any further assumptions on the lattice models. However, lattice models of statistical mechanics have a long story and exhibit many additional features. It is interesting to keep them in mind while reading the present work. In particular, as often as possible, we try in the following sections to point and discuss similarities and differences between our constructions and similar existing ones.

\subsubsection{Kolmogorov's extension theorem and its alternatives}

In dimension one (where no phase transition occurs), Kolmogorov's extension theorem is trivial to use and is thus used in practice since the required consistencies are directly related to the Perron-Frobenius eigenvectors on the two boundary vertices.

Most successes in larger dimensions so far are obtained through the analytical way: it consists in the choice of an \emph{arbitrary guessed}  boundary condition (for example, a fixed value all around the boundaries) ---which breaks consistency between increasing domains--- for large domains, followed by the proof of the convergence of local observables to some limit values when the domains grow to the whole lattice $\setZ^2$ and finally completed by a use of Riesz' theorem to obtain the infinite-volume Gibbs measures. After restriction from $\setZ^2$ to a finite domain, the restricted boundary conditions are only approximated by the ones chosen arbitrarily. We now discuss the advantages and disadvantages of this method:
\begin{itemize}
	\item guessing approximate boundary conditions is easy when the phase diagram is known and simple to understand (for example the Ising model) but it requires a physical intuition sometimes hard to acquire when subtle correlations enter the game or when the phase diagram is still unknown;
	
	\item convergence of local observables often requires a prior understanding of the physical content ---to separate relevant contributions from remainders going to zero--- as well as ad hoc inequalities (FKG ones for example); 
	
	\item the final measure following Riesz's existence theorem and its finite-domain restrictions may not have an easy-to-use structure and allows to produce precise computations of observables only through the approximation scheme above.
\end{itemize}
The present approach through Kolmogorov's extension theorem completely shifts the advantages and disadvantages above and thus provide an interesting different point of view:
\begin{itemize}
	\item there is no guessing in the finite domain boundary weights $g_R$: only a structure Ansatz \eqref{eq:boundaryweight:quick} and eigenvector-like equations to solve for the operators in \eqref{eq:boundaryweight:quick} (see Sections~\ref{sec:invariantboundaryelmts} and \ref{sec:applications}); 
	
	\item once these eigenvector-like equations are explicitly solved or existence of solutions is proved, finite domains laws are explicitly known in terms of the operators in \eqref{eq:boundaryweight:quick} and computations of correlation functions can be made without approximations, see Section~\ref{sec:appli:correlfn};
	
	\item all the difficulty is concentrated in the study of the eigenvector-like or fixed-point-like equation of Section~\ref{sec:invariantboundaryelmts} and there are mainly two ways of studies. The first one is, of course, to reuse all the guessings used in the analytical way when available and translate them into informations about the operators in \eqref{eq:boundaryweight:quick}. In absence of sufficient heuristics, the second way is to consider the eigenvector-like equations on their own and try to solve them, numerically or exactly, or at least extract interesting features.
\end{itemize}
Obtaining concrete equations to study on their own ---independent of any other heuristics or tools--- offers an interesting point of view.

\subsubsection{Matrix Product States and Integrability} Matrix product states (MPS) (see \cite{BlytheEvans,MPSreview} among others), similar to the ROPEreps introduced in Section~\ref{sec:boundaryalgebra}, have become a classical tool in quantum and statistical mechanics for various reasons. Heuristically, one may observe that whenever one manages to make a matrix product state technique work, from numerical density matrix renormalization group (DMRG) \cite{DMRGMPS} to exactly solvable models such as ASEP \cite{DEHP,ASEPBrownianExc,LargeDevDensityASEP}, it leads to significant and interesting results. Since the seminal work \cite{DEHP}, MPS have been used as an Ansatz (the "matrix Ansatz") in exact computations of invariant measures of stochastic processes on the line, such as the ASEP, for example. However, are MPS a consequence of integrability? This is unclear, and we only know that integrability can be a great help in determining or constructing Matrix Ansätze, as discussed in \cite{AlcarazLazo,GolinelliMallick,Tetrahedron}. The present paper may suggest that, in this case, MPS are more a consequence of the Markov property, which is common in integrable models, and integrability only provides a rich landscape to work in. As presented in Section~\ref{sec:invariantboundaryelmts}, computing boundary eigen-elements is similar to finding eigen-spaces of endomorphisms and is a difficult task in general: as with traditional transfer matrices, integrability should help to find these elements.

\subsubsection{Phase Transitions.} In one-dimensional statistical mechanics, where invariant boundary conditions are simply Perron-Frobenius eigenvectors, it has been known since Peierls that no phase transition can occur. However, in higher dimensions, phase transitions are observed in many models. This implies that the definitions of eigen-elements up to morphisms, given in Definitions~\ref{def:eigenalgebrauptomorphims} and \ref{def:cornereigensemigroups}, allow for a broader variety of behaviours. In particular, at the phase transition point of various two-dimensional Markovian models, conformal invariance is expected to appear; recent works by Smirnov and collaborators introduce discrete complex analysis on the lattice to describe this emergence of conformal invariance. Once again, one may expect a nice mathematical combination of these tools with our operadic approach to boundary weights.

\subsubsection{(Topological) (Quantum) Field Theory.} As briefly mentioned, lattice Markov processes with a discrete state space are toy models of Field Theory, featuring the same algebraic structure adapted to discrete lattices, without the analytical complications of distributional fields and associated renormalization procedures. Many recent works (see for example those by K.~Costello \cite{Costello,CostelloGwilliam}) provide solid algebraic foundations for Field Theory, whether topological or not, perturbative or not, based on various operadic structures and factorization algebras, which are very close to the structures present in this paper. However, the purposes are very different: while in Field Theory, the very first question is the one of rigorously \emph{defining} most models, in our case the definitions are immediate and our major focus is on the ability to \emph{compute} efficiently thermodynamic limits, correlation functions, and critical exponents. For the record, the starting point of the present work was recognizing the structure of MPS for boundary weights of two-dimensional Markov models on the lattice in figure 3 of \cite{GinotQFTfacto}, and the desire to formalize this correspondence in the simplest case of discrete space and finite state space.

\subsubsection{Integrable Toolbox and Exactly Solvable Models.} A striking feature of various constructions below is their resemblance to other constructions that appear in the study of exactly solvable models, though always with a different perspective and distinct mathematical constructions. However, we do not require any hint of integrability. The literature on exactly solvable systems already contains glimpses of some aspects of the construction described below, but without fully recognizing their absolute nature, establishing complete relations with each other, or extracting the fundamental role of generalized associativities. This last role of associativities is more easily seen here, free from the complexity of the Yang-Baxter equation \cite{JIMBO_1989}. Below, we list a variety of existing results directly related to the present paper: 
\begin{itemize} 
	\item Zamolodchikov-Faddeev algebras \cite{ZFFaddeev,ZFZamolodchikov} are already known to be related to matrix Ansätze \cite{vanicatZF} and correspond to algebras of operators on which Yang-Baxter face weights given by $R$-matrices act nicely (in the unitary case). We expect that there may be a relationship between ZF algebras and boundary algebras up to morphisms, as defined in Section~\ref{sec:invariantboundaryelmts}, but this may be subtle and is postponed to a second work in progress after the present one. 
	
	\item The corner elements in Chapters \ref{sec:operad}, \ref{sec:boundaryalgebra}, and \ref{sec:invariantboundaryelmts} are nearly identical to Baxter corner transfer matrices for integrable systems \cite{BaxterCTM_1981,BaxterCTM_2007,baxterbook}, and diagrams similar to our half-plane guillotine products are already present in Baxter's works and related studies. Some of the so-called "vertex operators" also share similarities. 
	
	\item Corner transfer matrices, as used in the corner transfer matrix renormalization group (CTMRG) \cite{Nishino_1996} (and more generally matrix product states in the density matrix renormalization group (DMRG) \cite{DMRGMPS}), show very good approximation properties under truncation (no rigorous proof yet, only strong heuristics), and there may be strong connections with the morphisms introduced in most of the definitions in Chapter~\ref{sec:invariantboundaryelmts}. 
	
	\item To the best of our knowledge, however, there is no systematic study of the relation between Zamolodchikov-Faddeev algebras and Baxter corner transfer matrices, though the present paper suggests that such a relationship may exist. 
	
	\item More generally, we do not know of any equivalent formulation of our definitions up to morphisms \ref{def:eigenalgebrauptomorphims}, \ref{def:cornereigensemigroups}, and \ref{def:commutuptomorph} in the literature on integrable systems. One possible reason is that most works around Baxter’s studies deal with strips and half-plane geometries (see Figure~\ref{fig:admissiblepatterns} below), or more generally, doubly infinite one-dimensional geometries. \item Drinfeld's construction of $R$-matrices relies on Hopf algebras, whose product and coproduct correspond, in the end, to each dimension of the two-dimensional geometry but breaks the symmetry between both. Our construction restores the symmetry and is amenable to generalization to larger dimensions, but since we do not focus on integrable systems, there is no such thing as the Yang-Baxter equation or quantum groups at this point. \end{itemize} We choose to separate our presentation as much as possible from the theory of integrable systems, despite various pointwise similarities. We rather think that various techniques have their origin in the Markov property, and the goal of this paper is to push this perspective as far as possible. The observation that some constructions may be nicely embodied in exactly solvable systems using the Yang-Baxter equation is postponed to a future paper in preparation.

	\section{Acknowledgements}
	
We first thank Gr\'egory Ginot for a gentle introduction to factorization algebras and $E_n$-structures at the origin of the present work. We also thank Najib Idrissi for all his answers to my numerous questions about operads and coloured operads as well as for pointing us the hidden commutativity phenomenon and references about the interchange relation. We also thank Dominique Manchon for encouragement and his insistence on Eckmann-Hilton argument. On the probabilistic side, we thank Loren Coquille for her questions about Gibbs measures and Vincent Beffara for a interesting discussion "on the road" about the computational complexity of finding ground states of Hamiltonians in dimension strictly larger than $1$. We thank Thierry L\'evy for a careful reading and interesting suggestions of presentation. We thank warmfully \'Emilien Bodiot for the weekly discussions about this work, his numerous readings and relevant comments at all stages of the writing.

	\section{Notations and theorems}
	
\paragraph*{Markov weights.}
During the paper, we will often switch between different dimensions. For example, weights for laws of Markov processes are matrices in dimension one and live on edges, whereas they are tensors and live on faces in dimension two (and other types of useful objects live on edges). Thus, we will write all the Markov weights, regardless of the dimension, as
\[
\MarkovWeight{A}, \MarkovWeight{B}, \MarkovWeight{W}, \MarkovWeight{\Omega}, \MarkovWeight{R}, \text{etc.}
\]
with capital "sans serif" letters.

\paragraph*{Rectangles.}
All the rectangles are written with a capital $R$ and have sizes written \[(P,Q), (P',Q'), (P_1,Q_1), \text{etc.} \] where $p$ is the horizontal length and $q$ the vertical one.

Variables associated to boundaries of rectangles are always called $x$ on the bottom ("South") of the rectangle, $y$ on the top ("North"), $w$ on the left ("West") and $z$ on the right ("East"), as illustrated in \eqref{eq:boundaryordering},  excepted inside summations on the boundary between two rectangles where these names cannot match.

\paragraph*{Random variables and state spaces.} 
Random variables are written with capital letters such as $X$. Discrete spaces are also written with a capital $S$ and $S_1$ (resp. $S_2$) is the state space associated to horizontal (resp. vertical) edges.

\paragraph*{Partition functions.} Regardless of the dimension, a partition function is always written\[
Z_{D}(\MarkovWeight{W}_\bullet ; x) \text{ or } Z_{D}^{\boundaryweights}(\MarkovWeight{W}_\bullet ; g)
\]
with a capital $Z$. The index $D$ indicates the domain. The first argument is the collection of weights associated to the faces of the domain $D$. In the first (resp. second) case, the second argument is the configuration of variables (resp. the boundary weight) on the boundary of the domain $D$. Boundary weights are always written with a letter $g$ or $G$.

\paragraph*{Algebras over operads.} All the algebras over (coloured) operads have names with capital curly letters $\ca{A}$, $\ca{B}$, $\ca{T}$, etc. Whenever this is clear from context, collections of spaces $(\ca{A}_i)_{i\in I}$ will be written $\ca{A}_I$ to simplify notations. The only other appearance of such a curly notations is for some $\sigma$-algebras in Section~\ref{sec:proba}.

\paragraph*{List of theorems.}
\begin{itemize}
	\item Theorem~\ref{theo:canonicalexampleGuill}: tensor algebra of matrices as guillotine algebras,
	\item Theorem~\ref{theo:partitionfuncguillotop}: partition functions as guillotine algebra,
	\item Theorem~\ref{theo:1D:removingcoloursonboundaries}: morphism from $\Guill_1$ to $\Ass$ through boundaries,
	
	\item Theorems~\ref{theo:extendedguilloperads} and \ref{theo:operadpointedversion}: guillotine operad with a corner and guillotine operad with pointed doubly-infinite shapes,
	
	\item Theorem~\ref{theo:canonicalboundarystructure}: canonical boundary guillotine algebra for Markov processes,
	
	\item Theorem~\ref{theo:stability}: structural stability of ROPEreps of boundary weights of Markov processes,
	
	\item Theorems~\ref{theo:reductionofmorphisms:canostruct} and \ref{theo:eigencorner:realizationMarkov}: reduction of morphisms on half-strips and corners for the canonical structure associated to a Markov process,
	
	\item Theorem~\ref{theo:eigenROPErep:invmeas}: invariant measures from boundary ROPE eigen-structure, from local to global, using Kolmogorov's extension theorem,
	
	\item Theorem~\ref{theo:onedimmarginaloutofinvROPErep}: one-dimensional marginals from invariant ROPEs,
	
	\item Theorem~\ref{theofromBS:Gaussian}: application to the Gaussian case summarizing \cite{BodiotSimon},
	
	\item Theorem~\ref{theo:higherdim:canostruct}: canonical linear algebra structure on hyper-cubic lattices,
	
	\item Theorem~\ref{theo:canonicalexampleGuill:continuous}: guillotine algebra for Markov processes in continuous space.
\end{itemize}

\chapter{Probabilistic formalism for Markov processes in dimensions one and two}\label{sec:proba}

		\section{Geometrical notations for the \texorpdfstring{$\setZ^2$}{Z2} and \texorpdfstring{$\setZ$}{Z} lattices}
		\subsection{The \texorpdfstring{$\setZ^2$}{Z2} lattice.}
We consider $\setZ^2$ as a planar graph with vertices $\mathcal{V}=\setZ^2$, oriented edges given by
\begin{align*}
\Edges{\setZ^2}_{\mathrm{or.}} &= E^{\uparrow}_{\mathrm{or.}} \cup {E}^{\downarrow}_{\mathrm{or.}} \cup {E}^{\leftarrow}_{\mathrm{or.}} \cup {E}^{\rightarrow}_{\mathrm{or.}}
\end{align*}
with the following types of oriented edges in the four possible directions
\begingroup
\allowdisplaybreaks
\begin{align*}
\Edges{\setZ^2}^{\uparrow}_{\mathrm{or.}} &= \{((k,l),(k,l+1)) ; (k,l)\in\setZ^2 \}
\\
\Edges{\setZ^2}^{\downarrow}_{\mathrm{or.}} &= \{((k,l),(k,l-1)) ; (k,l)\in\setZ^2 \}
\\
\Edges{\setZ^2}^{\leftarrow}_{\mathrm{or.}} &= \{((k,l),(k-1,l)) ; (k,l)\in\setZ^2 \}
\\
\Edges{\setZ^2}^{\rightarrow}_{\mathrm{or.}} &= \{((k,l),(k+1,l)) ; (k,l)\in\setZ^2 \}
\end{align*}
It is convenient to consider also unoriented edges. To this purpose, we introduce the projection $\pi : \setZ^2\times \setZ^2 \to P(\setZ^2)$, $(e_1,e_2)\mapsto \{e_1,e_2\}$, the horizontal and vertical unoriented edges
\begin{align*}
	\Edges{\setZ^2}_{h}&= \pi(\Edges{\setZ^2}^{\leftarrow}_{\mathrm{or.}}) = \pi( \Edges{\setZ^2}^{\rightarrow}_{\mathrm{or.}})
	\\
	\Edges{\setZ^2}_{v}&= \pi(\Edges{\setZ^2}^{\downarrow}_{\mathrm{or.}}) = \pi( \Edges{\setZ^2}^{\uparrow}_{\mathrm{or.}})
	\\
	\Edges{\setZ^2} &= \Edges{\setZ^2}_{h} \cup \Edges{\setZ^2}_{v}
\end{align*}
For oriented edges, we introduce the source and target maps $s,t: \mathcal{E}_{\mathrm{or.}} \to \setZ^2$ such that $s((a,b))=a$ and $t((a,b))=b$. We also introduce the set of faces
\[
\Faces(\setZ^2) = \{ [k,k+1]\times [l,l+1] ; (k,l)\in\setZ^2 \}
\]
Here the segments are either considered as subsets of $\setZ$ as $[a,b]=\{a,a+1,\ldots,b\}$ or as usual segments of $\setR$, without any further consequence for the results. For each $i\in \{N,W,S,E\}$ (to be interpreted as North, West, South and East), we introduce oriented boundaries $b_i : \Faces(\setZ^2) \to \Edges{\setZ^2}_{\mathrm{or.}}$ defined by:
\begin{align*}
b_N([k,k+1]\times [l,l+1]) &= ( (k+1,l+1),(k,l+1) ) \in \Edges{\setZ^2}_\mathrm{or.}^{\leftarrow}
\\
b_W([k,k+1]\times [l,l+1]) &= ( (k,l+1),(k,l) ) \in \Edges{\setZ^2}_\mathrm{or.}^{\downarrow}
\\
b_S([k,k+1]\times [l,l+1]) &= ( (k,l),(k+1,l) ) \in \Edges{\setZ^2}_\mathrm{or.}^{\rightarrow}
\\
b_E([k,k+1]\times [l,l+1]) &= ( (k+1,l),(k+1,l+1) ) \in \Edges{\setZ^2}_\mathrm{or.}^{\uparrow}
\end{align*}
and unoriented boundaries by $\pi\circ b_i$, that we will still write $b_i$ when there is no doubt. We also introduce the reversal map $r(((a,b),(a',b')))=((a',b'),(a,b))$ on $\Edges{\setZ^2}_\mathrm{or.}$.

\begin{defi}[domain of $\setZ^2$]
A \emph{domain} $D$ of $\setZ^2$ is a finite subset of $\Faces(\setZ^2)$. We introduce the set of oriented edges and unoriented edges:
\begin{align*}
\Edges{D}_\mathrm{or.} &= b_N(D) \cup b_W(D) \cup b_S(D) \cup b_E(D) \subset \Edges{\setZ^2}_{\mathrm{or.}}
\\
\Edges{D} &= \pi(\Edges{D}_\mathrm{or.}) \subset \Edges{\setZ^2}
\end{align*}
The \emph{oriented boundary} $\partial_\mathrm{or.} D$ of $D$ and the \emph{unoriented boundary} $\partial D$ are defined as:
\begin{align*}
\partial_\mathrm{or.} D &= \left\{
e\in \Edges{D}_\mathrm{or.} ; r(e) \notin \Edges{D}_\mathrm{or.} \right\}
\\
\partial D &= \pi(\partial_{\mathrm{or.}} D) \subset \Edges{\setZ^2}
\end{align*}
The set of \emph{vertices} of $D$ is the set $V(D)$ of points of $\setZ^2$ in the union of the faces $\cup_{f\in D} f$. The \emph{area} of $D$ is given by $\card(D)$ and its \emph{perimeter} is the cardinal of $\partial D$. 
\end{defi}
\endgroup

The internal unoriented edges of $D$ are elements of the set:
\[
\Int(D) = \Edges{D} \setminus \partial D
\]

The oriented boundary of a domain has an orientation, that is, for any $e\in \partial_\mathrm{or.} D$, there is a unique successor $\mathrm{succ}(e)\in \partial_\mathrm{or.} D$ such that $s(\mathrm{succ}(e))=t(e)$ and, if there are multiple solutions to this equation (and only in this case), $\mathrm{succ}(e)$ has an orientation turned counter-clockwise with respect to $e$ (i.e. if $e\in \Edges{\setZ^2}_\mathrm{or.}^\uparrow$ (resp. $\leftarrow$, $\downarrow$, $\rightarrow$), then $\mathrm{succ}(e)\in \Edges{\setZ^2}_\mathrm{or.}^{\leftarrow}$ (resp. $\downarrow$, $\rightarrow$, $\uparrow$)). A domain $D$ is connected if and only if $\mathrm{succ}$ has a unique cycle. Every cycle of $\mathrm{succ}$ is the oriented boundary of a connected component of $D$.

\begin{figure}
\begin{center}
\begin{tikzpicture}[scale=0.75]
% ENSEMBLE D
\begin{scope}
\node at (0.,3.) {$D$};
\foreach \x/\y in {0/0,1/0,2/0,3/0,0/1,1/1,2/1,3/1,1/2,2/2,3/2,2/-1} 
{
	\begin{scope}[xshift=\x cm,yshift=\y cm]
		\fill[guillfill] (0,0)--(1,0)--(1,1)--(0,1)--cycle ;
		\draw[guillsep] (0,0)--(1,0)--(1,1)--(0,1)--cycle ;
	\end{scope}
}
\end{scope}

\draw[->] (4.5,1) -- node [midway,above] {$b_N,b_W,b_S,b_E$} (7.5,1);

% ENSEMBLE E
\begin{scope}[xshift=8cm]
\node at (2.,3.5) {$\Edges{D}_{\mathrm{or.}}$};
\foreach \x/\y in {0/0,1/0,2/0,3/0,0/1,1/1,2/1,3/1,1/2,2/2,3/2,2/-1} 
{
	\begin{scope}[xshift=\x cm,yshift=\y cm]
		\draw[->,thick] (0.1,0.07) to (0.9,0.07);
		\draw[->,thick] (0.93,0.1) to (0.93,0.9);
		\draw[->,thick] (0.9,0.93) to (0.1,0.93);
		\draw[->,thick] (0.07,0.9) to (0.07,0.1);
	\end{scope}
}
\end{scope}

% ENSEMBLE Eunor
%\begin{scope}
%\foreach \x/\y in {0/0,1/0,2/0,3/0,0/1,1/1,2/1,3/1,1/2,2/2,3/2,2/-1} 
%{
%	\begin{scope}[xshift=\x cm,yshift=\y cm]
%		\draw[thick] (0.1,0.) to (0.9,0.0);
%		\draw[thick] (1,0.1) to (1,0.9);
%		\draw[thick] (0.9,1) to (0.1,1);
%		\draw[thick] (0.,0.9) to (0.,0.1);
%	\end{scope}
%}
%\end{scope}

\draw[->] (10,-1.5) -- node [midway,right] {$s,t$} (10,-2.5);
// ENSEMBLE V
\begin{scope}[xshift=8cm,yshift=-6cm]
\node at (0.,3.) {$V(D)$};
\foreach \x/\y in {0/0,1/0,2/0,3/0,0/1,1/1,2/1,3/1,1/2,2/2,3/2,2/-1} 
{
	\begin{scope}[xshift=\x cm,yshift=\y cm]
		\node at (0,0) [fill,circle,inner sep=0.5mm] {};
		\node at (1,0) [fill,circle,inner sep=0.5mm] {};
		\node at (1,1) [fill,circle,inner sep=0.5mm] {};
		\node at (0,1) [fill,circle,inner sep=0.5mm] {};
	\end{scope}
}
\end{scope}

\draw[->] (2,-1.5) -- node [midway,right] {$\partial$} (2,-2.5);

\begin{scope}[yshift=-6cm]
\node at (0.,3.) {$\partial D$};
\foreach \xs/\ys/\xe/\ye in {0/0/1/0, 1/0/2/0, 2/0/2/-1, 2/-1/3/-1, 3/-1/3/0, 3/0/4/0, 4/0/4/1, 4/1/4/2, 4/2/4/3, 4/3/3/3, 3/3/2/3, 2/3/1/3, 1/3/1/2, 1/2/0/2, 0/2/0/1, 0/1/0/0}
{
	\draw[thick] (\xs,\ys) to (\xe,\ye);
}
\end{scope}
\end{tikzpicture}
\end{center}
\caption{The structure of a domain $D$ of $\setZ^2$.}\label{fig:exampledomain}
\end{figure}

We now introduce the notions of partition and shell of a partition of a domain, which are fundamental for the definition of Markov processes.

\begin{defi}[partition of a domain, shell of a partition]
Let $D$ be a domain of $\setZ^2$. A \emph{partition} of $D$ with size $n$ is a $n$-uplet of non-empty domains $(D_i)_{1\leq i\leq n}$ such that $D_i\cap D_j=\emptyset$ for $i\neq j$ and $\cup_i D_i=D$. The \emph{shell} of a partition $(D_i)_{1\leq i\leq n}$ is the set $\Shell(D_\bullet)=\cup_{1\leq i\leq n} \partial D_i$.
\end{defi}

These two notions illustrated in Figure~\ref{fig:examplePartitionAndShell} are at the heart of the geometrical interpretation of the Markov property: conditioning on a shell introduces independence between the elements of the partition.

\begin{figure}
\begin{center}
\begin{tikzpicture}[scale=0.75]
% ENSEMBLE D
\begin{scope}
\node at (0.,3.) {$D_\bullet$};
\foreach \x/\y in {-0.1/0.1, 1.1/0, 2.1/0, 3.1/0, -0.1/1.1, 0.9/1.1, 2/1.2, 3.1/1, 1/2.2, 2/2.2, 3/2.2, 2.1/-1} 
{
	\begin{scope}[xshift=\x cm,yshift=\y cm]
		\fill[guillfill] (0,0)--(1,0)--(1,1)--(0,1)--cycle ;
		\draw[guillsep] (0,0)--(1,0)--(1,1)--(0,1)--cycle ;
	\end{scope}
}
\end{scope}

\begin{scope}[xshift=7cm]
\node at (0.,3.) {$\Shell(D_\bullet)$};
\foreach \xs/\ys/\xe/\ye in {0/0/1/0, 1/0/2/0, 2/0/2/-1, 2/-1/3/-1, 3/-1/3/0, 3/0/4/0, 4/0/4/1, 4/1/4/2, 4/2/4/3, 4/3/3/3, 3/3/2/3, 2/3/1/3, 1/3/1/2, 1/2/0/2, 0/2/0/1, 0/1/0/0, 1/0/1/1, 1/1/3/1, 3/1/3/2, 3/2/4/2, 2/2/1/2, 2/1/2/2}
{
	\draw[thick] (\xs,\ys) to (\xe,\ye);
}
\end{scope}
\end{tikzpicture}

\end{center}
\caption{Partition $D_\bullet$ of size $3$ of the domain $D$ of $\setZ^2$ described in figure~\eqref{fig:exampledomain} and its shell $\Shell(D_\bullet)$. We observe that $\Shell(D_\bullet)$ contains all the boundaries $\partial D_i$ and the boundary $\partial D$.}
\label{fig:examplePartitionAndShell}
\end{figure}

\subsection{Reminder about the one-dimensional lattice.} We make the choice of introducing the simpler case of the $\setZ$ lattice after the more complicated case of the $\setZ^2$ lattice, in order to connect the notions to the traditional ones in the theory of one-dimensional stochastic processes.

In dimension one, there are only vertices $\setZ$, oriented edges and unoriented ones $\{\{k,k+1\};k\in\setZ\}$. A domain $D$ is a finite collection of edges and $V(D)$ is its set of vertices (extremities of the edges). The boundary $\partial D$ of a domain $D$ is the set of vertices that belong to exactly one edge of the domain. A connected domain is a set of consecutive edges $\{(k,k+1); K_0\leq k <K_1\}$, with a two-point boundary $\{K_0,K_1\}$. The shell of a partition $(D_i)_{1\leq i\leq n}$ of a domain $D$ is the union of the boundaries $\partial D_i$ and it is a subset of $\setZ$. A domain $D$ is connected if and only if $\card(\partial D)=2$, or equivalently $D$ is of the form $\{ k\in\setZ; a\leq k\leq b\}$ for some (unique) integers $a,b\in\setZ$.

	\section{Structural property of the laws of Markov processes on \texorpdfstring{$\setZ^2$}{Z2}}
	
		\subsection{From the Hammersley-Clifford's theorem for graphs to regular lattices.}
	
In the discrete space setting, the Markov property may be defined on arbitrary unoriented graphs, for which edges provide the notion of neighbourhood of a vertex. The Markov property then states that the conditional law of a process on any given vertex $x$ w.r.t. to the other vertices depends only on the vertices in the neighbourhood of $x$. The Hammersley-Clifford's theorem then states that the law of the whole process can be factorized over the complete sub-graphs (or cliques) of the initial graph, where each factor depends only on the value of the process on the vertices in the corresponding complete sub-graph. 

This factorization structure appears everywhere in the formulae below but, at this level of generality, there may not be so much additional structure and no hope to formulate nice generalization of the one-dimensional setting as well as to recover various key notions of statistical mechanics in the Euclidean space. Moreover, it is not easy to translate it easily in the continuous space setting where the notion of neighbourhood is hidden in the infinitesimal structure of the model. 

Keeping in mind that our final purpose is to find an algebraic description of Markov processes common to discrete as well as to continuous space, we prefer here the point of view of a Markov property formulated in a \emph{global} perspective, i.e. by considering partitions of domains instead of neighbourhood dependencies. 

From a technical point of view, the system of colours introduced in Section~\ref{sec:operad}, which is a key feature of the present work, cannot be adapted in an as much powerful tool for the general case of graphs.

		\subsection{A reminder on one-dimensional time-oriented and spatial non-oriented Markov processes}

In order to fix some notations and introduce easy notions, we first start with a gentle reminder on one-dimensional Markov processes.	

In dimension one, the sets $\setZ$ or $\setN$ can either be interpreted as time or space depending on the context. This distinction leads to two different definitions of one-dimensional Markov processes: in the first case of time interpretation, the sets $\setZ$ and $\setN$ are oriented and there is the traditional definition of Markov chains through stochastic matrices and generators in the framework of filtrations (which definition involves the orientation). However, since the orientation is lost in most two-dimensional cases of statistical mechanics, we rather focus on a definition of one-dimensional spatial (non-oriented) Markov processes on domains of $\setZ$.

We will produce below two equivalent definitions: a purely probabilistic one in terms of conditional expectations related to the geometry of domains and a second one with specific probability laws involving algebraic objects (matrices here).

\subsubsection{Description by the law for a finite set of possible values}\label{para:lawfiniteset}

\begin{defi}[definition of a 1D Markov process with a finite set of values from its law]
\label{def:MarkovOneDimFirst}
Let $S$ be a finite set (endowed with its complete $\sigma$-algebra). Let $D$ be a domain of $\setZ$. An 
$S$-valued Markov process $X=(X_v)_{v\in V(D)}$ on $D$ is a collection 
of $S$-valued random variables $X_v$ (on a probability space $(\Omega,\ca{F},\Prob)$) 
indexed by the set of vertices of $D$ such that there exist a collection $(\MarkovWeight{A}_e)_{e\in D}$ of matrices indexed by $S\times S$ such that the conditional law given the boundary values is given by:
\begin{equation}\label{eq:MarkovLawDimOne}
\probc{ (X_v)_{v\in V(D)}=(x_v)_{v\in V(D)} }{ \ca{F}_{\partial D} } = \frac{1}{Z_D(\MarkovWeight{A}_\bullet;(x_v)_{v\in \partial D})} \prod_{v\in \partial D} \indic{X_v=x_b} \prod_{e\in D} \MarkovWeight{A}_e(x_{s(e)},x_{t(e)}) 
\end{equation}
where $\ca{F}_{\partial D}= \sigma(X_v; v \in \partial D)$ and $Z_D(\MarkovWeight{A}_\bullet;(x_v)_{v\in \partial D})$ is a positive real number, called the \emph{partition function} on the domain $D$.
\end{defi}

We first introduce a closed expression of the partition function, which plays a fundamental role in statistical mechanics and in the algebraic approach of Section~\ref{sec:operad}.
\begin{lemm}
With the same notations as in Definition~\ref{def:MarkovOneDimFirst} and Proposition~\ref{prop:MarkovOneDimSecond}, if the domain $D$ is the union of $C$ connected components $D_c=\{L_c,L_c+1,\ldots,R_c\}$ with $L_1<R_1<L_2<\ldots<L_C<R_C$,
it holds:
\begin{align*}
Z_D(\MarkovWeight{A}_\bullet; (x_v)_{v\in\partial D}) &= \prod_{c=1}^C 
Z_{D_c}(\MarkovWeight{A}_\bullet; x_{L_c},x_{R_c} )
\\
Z_{D_c}(\MarkovWeight{A}_\bullet; x_{L_c}, x_{R_c} ) &= \left(\MarkovWeight{A}_{(L_c,L_c+1)} \MarkovWeight{A}_{(L_c+1,L_c+2)}\ldots \MarkovWeight{A}_{(R_c-1,R_c)} \right)(x_{L_c},x_{R_c})
\end{align*}
where the product between matrices $\MarkovWeight{A}$ and $\MarkovWeight{B}$ is the standard matrix product
\begin{equation}\label{eq:matrixproduct:discrete}
	(\MarkovWeight{A}\MarkovWeight{B})(x,y)=\sum_{z\in S} \MarkovWeight{A}(x,z)\MarkovWeight{B}(z,y)
\end{equation}
\end{lemm}
\begin{proof}
This is a consequence of the simple fact that $\Espc{1}{\ca{F}_{\partial D}}=1$. Summing over the possible values of the r.v. $X_v$ produces, for each vertex $v$ not on the boundary, a matrix product. On the boundary, the conditioning selects the correct value $X_v$.
\end{proof}

\begin{lemm}[extraction of marginals]\label{lemma:marginallawDimOne}
Given a strictly increasing sequence $(K_i)_{0\leq i\leq n+1}$ of points in a connected domain $D=\{K,\ldots,L\}$ with $K_0=K$ and $K_{n+1}=L$, the marginal law of $(X_{K_i})_{1\leq i\leq n}$ is given by:
\begin{equation}\label{eq:law:marginalindimensionone}
\probc{ \prod_{i=1}^n \indic{X_{K_i}=x_i} }{\ca{F}_{\partial D} } = 
\frac{ Z_{D_1}(\MarkovWeight{A}_\bullet; X_K,x_{1}) Z_{D_2}(\MarkovWeight{A}_\bullet; x_1,x_2)\ldots  Z_{D_{n+1}}(\MarkovWeight{A}_\bullet; x_n,X_L) }{Z_{D}(\MarkovWeight{A}_\bullet;X_K,X_L)}
\end{equation}
where $D_i=\{K_{i-1},\ldots, K_i\}$. As a consequence, for functions $(h_i)_{1\leq i\leq n}$ from $S$ to $\setR$, one has:
\begin{equation}
	\Espc{ \prod_{i=1}^n f(X_{K_i}) }{\ca{F}_{\partial D}}
	= \frac{(Z_1 D_{h_1} Z_2 D_{h_2} \ldots Z_n D_{h_n} Z_{n+1})_{X_K,X_L}}{ Z_{X_K,X_L} }
\end{equation}
where $Z_i$ is the matrix with coefficients $(Z_i)_{x,y}=Z_{D_i}(A_\bullet; x,y)$ and $D_{h}$ is the diagonal matrix with coefficients $(D_h)_{x,x}=h(x)$.
\end{lemm}
The proof is left as a simple exercise in probability theory. This formula already illustrates the heart of Section~\ref{sec:operad} in dimension 2: observables are pointwise objects represented by a commutative algebra (of diagonal matrices), weights are objects related to edges and form an associative algebra related to the gluing of consecutive edges, boundary weights are left and right vectors to be acted on by the previous matrices. 

\begin{rema}[a fundamental structural observation] \label{rema:fundamentalremark} There is a deep structural property in formula~\eqref{eq:law:marginalindimensionone}, which will appear in the same way in larger dimensions and in the continuous case: the conditional law  \eqref{eq:law:marginalindimensionone} has the same structure as \eqref{eq:MarkovLawDimOne} excepted that the matrices $A_e(x_k,x_{k+1})$ attached to edges joining $k$ and $k+1$ are replaced by matrices $Z_{D_{k}}(\MarkovWeight{A}_\bullet; x_{k-1},x_k)$ attached to domains joining $X_{K_{k-1}}$ and $X_{K_k}$, which are matrices obtained by an algebraic combination of the matrices $A_e$ inside the domain $D_k$.
\end{rema}

\begin{rema}One also observes that the matrices $\MarkovWeight{A}_e$ are not uniquely defined. Let $(d_v)_{v\in V(D)}$ be diagonal matrices. Then, any global change $\MarkovWeight{A}_e \mapsto d_{s(e)} \MarkovWeight{A}_e d_{t(e)}^{-1}$ (with the corresponding change on the partition functions) lets the conditional probability \eqref{eq:MarkovLawDimOne} unchanged. Moreover, if, for some $v\in V(D)$, there is $x\in S$ such that $X_v\neq x$ a.s. then any coefficient $\MarkovWeight{A}_{(v-1,v)}(\cdot,x)$ or $\MarkovWeight{A}_{(v,v+1)}(x,\cdot)$ may be set equal to zero. 
\label{rema:gaugeinvariance1}
\end{rema}

An important property, trivially deduced from the previous lemmata, is the stability under restriction.
\begin{prop}\label{prop:onedim:markovrestrict}
Let $(X_v)_{v\in V(D)}$ be an $S$-valued 1D Markov process on a domain $D$. For any subdomain $D'\subset D$, the restriction $(X_v)_{v\in V(D')}$ is an $S$-valued 1D Markov process with the same weights $(\MarkovWeight{A}_e)_{e\in D'}$ inherited from $D$.
\end{prop}

Lemma~\ref{lemma:marginallawDimOne} and Proposition~\ref{prop:onedim:boundaryweight:restrict} together allow one to cut and join a Markov process defined on a domain $D$ to and from Markov processes on partitions $(D_i)_i$ of $D$. This geometric interpretation is formalized algebraically in Section~\ref{sec:operad}.

\subsubsection{Description from the conditional law for a measurable set of values}

All the previous examples can be generalized to the case where $(S,\ca{S})$ is a measurable space with small analytical modifications. We choose to present it here to reach a higher degree of generality and to show that the algebraic structure, which will be the main tool of the next sections, remains the same, in particular from a geometric point of view. However, the parallel with standard finite-dimensional linear algebra is somewhat more hidden and the formulation of the models requires additional analytical tools, e.g. reference measures on the vertices. Since the main purpose of the paper is to generalize the notion of eigenvectors to higher dimensions, the reader is encouraged to start by thinking with a finite state space in order to avoid technical analytical complications.

Switching from time-oriented Markov chains to spatial Markov chains requires a more subtle discussion on measures, transition measures, densities and transition densities. In the time-oriented case, there are one initial measure on $(S,\ca{S})$ at the initial time and transition laws $S\times\ca{S} \to S$ on the time increments (edges): increasing time by one combines an initial measure with a transition law to produce a new measure. If we are given a final time $N$, expectation values are obtained by combining the final measure with a measurable function (the observable) through an integral. In the spatial unoriented case, we now have two "initial" measures, one on the left and one on the right; an edge weight can be used to move the left initial measure to a new measure on the right neighbouring point or to move the right initial measure to a new measure on the left neighbouring point. When collapsing a segment to a single point from the left and from the right, we are left with two measures on the same point and there is a priori no reason to collapse it to a single number (measures may be mutually singular).

This problem can be overcome in the following way:
\begin{itemize}
	\item any point of $\setZ$ is endowed with a reference measure $\mu_x$ on $(S,\ca{S})$
	\item any edge $e$ of $\setZ$ is endowed with a measurable function $\MarkovWeight{A}_e : S\times S \to \setR_+$, which will serve as a transition density.
\end{itemize} 
Given two measurable functions $h_k,h_{k+1} : S \to \setR_+$ associated respectively to the points $k$ and $k+1$ and a measurable weight $\MarkovWeight{A}$ on $e=\{k,k+1\}$, we may define two new density functions associated respectively to the points $k+1$ and $k$
\begin{align*}
	(h_k \MarkovWeight{A})(y) &= \int_S h_k(x) A(x,y) d\mu_k(x)
	\\
	(\MarkovWeight{A} h_{k+1})(x) &= \int_S  A(x,y) h_{k+1}(y) d\mu_{k+1}(y)
\end{align*}
The difference with the time-oriented case are thus the following:
\begin{itemize}
	\item a transition law $\MarkovWeight{A}_{\{k,k+1\}} : S\times\ca{S} \to\setR_+ $ is assumed to have a density with the reference measure $\mu_{k}$ and the edge weight $\MarkovWeight{A}_{\{k,k+1\}}$ is now this density;
	\item an initial law on a point $k$ is assumed to have a density $h^L_k$ w.r.t. $\mu_k$;
	\item there is also a right initial law on the right which is also assumed to have a density $h^R_{k}$ w.r.t. $\mu_k$
	\item given left and right initial laws $h^L_k$ and $h^R_k$ on the point $k$, it defines a collapsed measure $h^L_k(x) h^R_k(x) d\mu(x)$ on the point $k$.
\end{itemize}

The density hypothesis w.r.t. to reference measures on points leads to an easy generalization of the previous section to Markov processes with values in $(S,\ca{S})$. Moreover, in practice, the reference measures $\mu_k$ are taken to be all equal in the interesting case of homogeneous processes.

\begin{prop}\label{prop:markov:1Dmeasurable}
	All the results of paragraph \ref{para:lawfiniteset} still hold provided the following replacements are made:
	\begin{itemize}
		\item $S$ is a measurable space $(S,\ca{S})$;
		\item matrices $\MarkovWeight{A}_e$ are replaced by measurable functions $S\times S\to \setR_+$;
		\item there is a family of $\sigma$-finite measures $(\mu_k)_{k\in\setZ}$, called reference measures;
		\item the product of matrices \eqref{eq:matrixproduct:discrete} between two weights $\MarkovWeight{A}_{\{k,k+1\}}$ and $\MarkovWeight{A}_{\{k+1,k+2\}}$ is replaced by
		\[
			(\MarkovWeight{A}_{\{k,k+1\}}\MarkovWeight{A}_{\{k+1,k+2\}})(x,y) = \int_S \MarkovWeight{A}_{\{k,k+1\}}(x,z)\MarkovWeight{A}_{\{k+1,k+2\}}(z,y) d\mu(z);
		\]
		\item the probability on the l.h.s. of \eqref{eq:MarkovLawDimOne} has a density w.r.t. $\otimes_{k\in\partial D} \mu_k$ given by the r.h.s. of \eqref{eq:MarkovLawDimOne};
		\item the density $(A_e)$ and the reference measures $(\mu_k)$ must satisfy the following finite partition function condition so that the r.h.s. exists: for any connected domain $D$ with extremities $k$ and $l$ and $\mu_k\otimes \mu_l$-almost any $(x,y)$ in $S\times S$, $Z_D(\MarkovWeight{A}_\bullet; x ,y) < +\infty$.
	\end{itemize}
\end{prop}

The previous case of finite $S$ can be recovered by considering $\ca{S}=\ca{P}(S)$ and by taking all the reference measures $\mu_k$ equal to the counting measure: all measurability assumptions become trivial as well as the density hypothesis and the finite partition function condition.

\subsubsection{Description by the Markov property}

A key property of the Markov processes, often used as a definition, is the Markov property, which has the advantage of avoiding the precise form of the law. Excepted for small sizes, it is equivalent to the previous definition.

\begin{prop}[characterization of a 1D Markov processes from conditional expectation]\label{prop:MarkovOneDimSecond}
With the same notations as Definition~\ref{def:MarkovOneDimFirst},
Let $S$ be a finite set (endowed with its complete $\sigma$-algebra) or a measurable set with a $\sigma$-algebra $\ca{S}$ and reference measures $(\mu_k)$. Let $D$ be a domain of $\setZ$ and 
 $X=(X_v)_{v\in V(D)}$ an $S$-valued Markov process on $D$. Then, for any partition $(D_i)_{1\leq i \leq n}$ of $D$ and any bounded real random variables $(U_i)_{1\leq i\leq n}$ such that, for each $1\leq i\leq n$, $U_i$ is $\sigma(X_v;v\in V(D_i))$-measurable, it holds:
\begin{equation}\label{eq:MarkovDimOne:CondExp}
\Espc{ \prod_{i=1}^n U_i }{\sigma(X_v;v\in\Shell(D_\bullet))}
= \prod_{i=1}^n \Espc{ U_i }{ \sigma(X_v; v\in \partial D_i)}
\end{equation}
almost everywhere. If the connected components of $D$ have a cardinal larger than $3$, then this property is equivalent to Definition~\ref{def:MarkovOneDimFirst}.
\end{prop}

This definition has the advantage of being immediately generalizable to the two-dimensional case \emph{mutatis mutandis} and hide the precise expression of the probability law. One also checks that it does not use any orientation of space in the conditional expectations.

\begin{proof}[Quick proof of property \ref{prop:MarkovOneDimSecond}]
Obtaining \eqref{eq:MarkovDimOne:CondExp} from the definition is a simple computation of conditional expectations using the previous lemmata.

We now assume \eqref{eq:MarkovDimOne:CondExp} and prove that the conditional law are given by \eqref{eq:MarkovLawDimOne}. The proof follows the same pattern as the proof of the Hammersley-Clifford theorem. We choose to present it in more details here since it gives a particular role to the associativity of the gluing of elements of partitions and nesting of partitions: the algebraic counterpart is precisely the content of Section~\ref{sec:operad}. 

We consider only the case of a connected domain $D=\{(k,k+1); K\leq k<L\}$ with $V(D)=\{K,K+1,\ldots,L\}$ and we establish the proof by recursion on the domain size. 

If $\card(D)=2$, then \eqref{eq:MarkovDimOne:CondExp} is a tautology and do not put any restriction on the joint law of the three variables $(X_K,X_{K+1},X_{K+2})$. If $\card(D)=3$, we consider the two partitions:
\begin{align*}
D^{(1)}_1 &=\{(K,K+1),(K+1,K+2)\}
& 
D^{(1)}_2 &=\{(K+2,K+3)\}
\\
D^{(2)}_1 &=\{(K,K+1)\}
& 
D^{(2)}_2 &=\{(K+1,K+2),(K+2,K+3)\}
\end{align*}
and we thus obtain:
\begin{align*}
\Espc{\prod_{j=1}^{2}\indic{X_{K+j}=x_{K+j}}}{\ca{F}_{\Shell(D^{(1)}_\bullet)}} 
&=  \Espc{\indic{X_{K+1}=x_{K+1}}}{X_K,X_{K+2}} \indic{X_{K+2}=x_{K+2}}
\\
&= a_1(x_{K+1}; X_K, x_{K+2}) \indic{X_{K+2}=x_{K+2}} 
\\
\Espc{\prod_{j=1}^{2}\indic{X_{K+j}=x_{K+j}}}{\ca{F}_{\partial D}} 
&=  a_1(x_{K+1}; X_K, x_{K+2})\Espc{\indic{X_{K+2}=x_{K+2}}}{\ca{F}_{\partial D}} 
\\
&= a_1(x_{K+1}; X_K, x_{K+2})b_1(x_{K+2}; X_K,X_{K+3}) 
\end{align*}
The same computation with the partition $D^{(2)}_\bullet$ provides the representation:
\[
\Espc{\prod_{j=1}^{2}\indic{X_{K+j}=x_{K+j}}}{\ca{F}_{\partial D}} 
= a_2(x_{K+2}; x_{K+1}, X_{K+3})b_2(x_{K+1}; X_K,X_{K+3}) 
\]
By restricting the functions to the only values $x\in S$ on each vertex with a non-zero probability, identifications of the previous two representations with $a_1$, $b_1$, $a_2$ and $b_2$ fixes the necessary factorization
\begin{align*}
u(x;y,z) = \alpha_u(x,y)\beta_u(y,z)\gamma_u(z,x)
\end{align*}
for each factor $u\in \{a_1,a_2,b_1,b_2\}$. There are additional relations between the factors of each functions but they are irrelevant here. We thus obtain the expected factorization:
\begin{align*}
\Espc{\prod_{j=1}^{2}\indic{X_{K+j}=x_{K+j}}}{\ca{F}_{\partial D}}  = \frac{\MarkovWeight{A}(X_K,x_{K+1}) \MarkovWeight{A}(X_{K+1},X_{K+2}) \MarkovWeight{A}(x_{K+2},X_{K+3})}{Z(X_K,X_{K+3})}
\end{align*}
from the previous factorizations.

We now assume that $D$ has $n> 3$ edges (and $n+1>4$ vertices) and that the equivalence between Definition~\ref{def:MarkovOneDimFirst} and Proposition~\ref{prop:MarkovOneDimSecond} holds for any process on a domain $D'$ with $3<m<n$ edges.

It is an easy to check that, if \eqref{eq:MarkovDimOne:CondExp} holds for a process $(X_v)_{v\in V(D)}$, then it also holds for any restriction $(X_v)_{v\in V(D')}$ with $D'\subset D$. 
We now consider the partition $D^{(1)}_1 = \{(K,K+1)\}$ and $D^{(1)}=\{(k,k+1); K+1\leq k<L\}$ and obtain the factorization:
\begin{align*}
\Espc{\prod_{j=K}^{L}\indic{X_{j}=x_{j}}}{\ca{F}_{\Shell(D^{(1)}_\bullet)}} 
&=  \Espc{\prod_{j=K+1}^{L}\indic{X_{j}=x_{j}}}{ \ca{F}_{\partial D^{(1)}_2} } \indic{X_{K}=x_{K}}
\\
&= \frac{  \indic{X_K=x_K} \indic{X_{K+1}=x_{K+1}} \indic{X_L=x_L} }{Z_{D^{(1)}_2}(\MarkovWeight{A}^{(1)}_\bullet; X_{K+1}, X_L)}
\prod_{k=K+1}^{L-1} \MarkovWeight{A}^{(1)}_{(k,k+1)}(x_k,x_{k+1})
\end{align*}
from the recursion hypothesis, which provides the existence of matrices $(\MarkovWeight{A}^{(1)}_{(k,k+1)})_{K+1\leq k< L}$ such that the conditional law in the r.h.s. is given by the \eqref{eq:MarkovLawDimOne}. We thus obtain finally
\begin{align*}
\Espc{\prod_{j=K}^{L}\indic{X_{j}=x_{j}}}{\ca{F}_{\partial D}} 
&= a_1(x_{K+1}; x_K,x_L)  \frac{  \indic{X_K=x_K}  \indic{X_L=x_L} }{Z_{D^{(1)}_2}(\MarkovWeight{A}^{(1)}_\bullet; x_{K+1}, x_L)}
\prod_{k=K+1}^{L-1} \MarkovWeight{A}^{(1)}_{(k,k+1)}(x_k,x_{k+1})
\\
a_1(x_{K+1}; x_K,x_L) &= \Espc{\indic{X_{K+1}=x_{K+1}}}{X_K=x_k,X_L=x_L}
\end{align*}
We still have to factorize the expression of $a_1$ as in the case of size $3$. This can be done by considering the partition $D^{(2)}_1=\{(k,k+1); K\leq k<L-1 \}$, $D^{(2)}_2=\{(L-1,L)\}$. The same computations produce the representation:
\begin{align*}
\Espc{\prod_{j=K}^{L}\indic{X_{j}=x_{j}}}{\ca{F}_{\partial D}} 
&= a_2(x_{L-1}; x_K,x_L)  \frac{  \indic{X_K=x_K}  \indic{X_L=x_L} }{Z_{D^{(2)}_1}(\MarkovWeight{A}^{(2)}_\bullet; x_{K}, x_{L-1})}
\prod_{k=K}^{L-2} \MarkovWeight{A}^{(2)}_{(k,k+1)}(x_k,x_{k+1})
\\
a_2(x_{L-1}; x_K,x_L) &= \Espc{\indic{X_{L-1}=x_{L-1}}}{X_K=x_k,X_L=x_L}
\end{align*}
Considering the equality between the two representations restricted to values $x_k$ with non-zero probability implies factorizations $a_i(y;x_K,x_L)=\alpha_i(x_K,y)\beta_i(y,X_L)\gamma_i(X_K,X_L)$ and thus the conditional law~\eqref{eq:MarkovLawDimOne}.

The generalization to disconnected domains $D$ partitioned into disconnected domain is made easy from the following fact: if a domain $D$ is the union of $C$ connected components $D'_c$, $1\leq c\leq C$, then any r.v. $U$ measurable w.r.t. $\ca{F}_D$ can be written as a finite sum $U=\sum_{j} \prod_{c=1}^{C} U_{c}^{(j)}$ where $U^{(j)}_c$ is $\ca{F}_{D'_c}$-measurable, due to the finite cardinal of $S$. By splitting a domain $D$ and the elements of any partition $D_\bullet$ of $D$ into finite numbers of connected components, linearity of conditional expectations reduces all the computations to the previous case.
\end{proof}

It is interesting to remark that \eqref{eq:MarkovDimOne:CondExp} does not give any restriction for a domain of size $2$, whereas any restriction of a Markov process on a large domain to a subdomain of size $2$ induces a factorization of the probability law as in \eqref{eq:MarkovLawDimOne}.		

\subsubsection{The role of boundary weights.} In Definition~\ref{def:MarkovOneDimFirst} as well as in Proposition~\ref{prop:MarkovOneDimSecond}, it is interesting to notice that the full law of the process is not written, although it is required in practice for any expectation values and correlation functions. Going from $\probc{\bullet}{\ca{F}_D}$ to $\prob{\bullet}$ on a connected domain $D=\{L,\ldots,R\}$ requires averaging over the two boundary variables $X_L$ and $X_R$. Excepted for particular periodic boundary conditions, one expects the two variables to exhibit correlations only through the variables in the strict interior of $D$ and we may consider a factorized law
\[
\prob{ X_L=x_L, X_R=x_r } = \frac{1}{Z'} u_{D,1}(x_L) Z_D(\MarkovWeight{A}_\bullet; x_L,x_R) u_{D,2}(x_R)
\]
with two boundary vectors $u_{D,1}$ and $u_{D,2}$ and a normalization $Z'=\scal{u_{D,1}}{Z_D(\MarkovWeight{A}_\bullet;\cdot,\cdot)u_{D,2}}$. This corresponds to a complete probability law:
\begin{equation}\label{eq:onedim:lawwithboundary}
\prob{ (X_v)_{v\in V(D)}=(x_v)_{v\in V(D)} }
=
\frac{1}{Z'} u_{D,1}(x_L) \left(
 \prod_{e\in D} \MarkovWeight{A}_e(x_{s(e)},x_{t(e)}) 
 \right) u_{D,2}(x_R)
\end{equation}
in which each variable appears exactly twice. In the case of a measurable state space with reference measures, the previous equation should be interpreted as a density with respect to the product of the reference measures on each point of the domain.

The following proposition illustrates how linear algebra governs the evolution of the boundary weights under the restrictions described in Proposition~\ref{prop:onedim:markovrestrict}.

\begin{prop}\label{prop:onedim:boundaryweight:restrict}
Following the notations of Proposition~\ref{prop:onedim:markovrestrict}, if the process $(X_v)_{v\in V(D)}$ on a connected domain $D=\{L,L+1,\ldots,R\}$ have a boundary weight given by \eqref{eq:onedim:lawwithboundary}, then the restriction to a connected subdomain $D'=\{L',L'+1,\ldots,R'\}$ have the same edge weights $\MarkovWeight{A}_\bullet$ and boundary weights given:
\begin{align*}
	u_{D',1}(y) = \sum_{x\in S}u_{D,1}(x) Z_{\{L,\ldots,L'\}}(\MarkovWeight{A}_\bullet; x,y)
	\\
	u_{D',2}(x) = \sum_{y\in S} Z_{\{R',\ldots,R\}}(\MarkovWeight{A}_\bullet; x,y)u_{D,2}(y)
\end{align*}
when $S$ is finite. In the case of a measurable state space with reference measures, the previous equations hold for densities and the summation are integrals over the reference measures at the given points.
\end{prop}
This illustrates the well-known fact that, in dimension one, the edge weights are elementary matrices whose products are the partition functions, which themselves act as matrices on the boundary weights seen as vectors. The next Section~\ref{sec:prob:twodim} presents the two-dimensional case, for which the description of the boundary weights is more intricate and is the key point of the present paper.

Boundary weights $u_{D,i}$ play an important role in the study of thermodynamic limits where $D$ is a large connected domain, whenever the edge weights $\MarkovWeight{A}_e$ are all the same.

\begin{prop}\label{prop:onedim:kolmoextension}
If all the edge weights $\MarkovWeight{A}_e$ are constant, equal to some $\MarkovWeight{A}$, and for any domain $D=\{K,\ldots,L\}$, the boundary weights $u_{D,1}$ and $u_{D,2}$ are all equal to the Perron-Frobenius left and right eigenvectors $v_{1}$ and $v_2$, with eigenvalue $\Lambda$, then one can define a process $(X_v)_{v\in\setZ}$, invariant under translation, whose marginal laws on $D$ are given by~\eqref{eq:onedim:lawwithboundary}. 
\end{prop}
\begin{proof}
It is an elementary computation to check that $Z'=c\Lambda^{K-L}$ for a domain $D=\{K,\ldots,L\}$ where $c$ is a normalization constant that depends on the two vectors $v_1$ and $v_2$. The existence of the process  the whole $\setZ$ is then a direct consequence of Kolmogorov's extension theorem. 
\end{proof}

From a computational perspective, boundary weights, and more generally, eigenvectors of the edge weight $\MarkovWeight{A}$, describe the correlation functions and are thus very important despite their absence in the Markov property. Phrased differently, we may say that the Markov property is a first step to reduce most computations of expectations to quantities that then depend only on boundary weights.

To illustrate this fact, one has, in the context of proposition of Proposition~\ref{prop:onedim:kolmoextension}, the following expectations:
\begin{align*}
	\Esp{\indic{X_v=x}} &= \frac{1}{c}v_1(x)v_2(x) 
	\\
	\Var\left(\indic{X_v=x} \indic{X_{v'}=x'}\right) &\sim c' \left(\frac{\Lambda'}{\Lambda}\right)^{d(v,v')} 
\end{align*}
as $d(v,v')\to\infty$ if $\Lambda'$ is the second (non-degenerate eigenvalue) of $\MarkovWeight{A}$.

	\subsection{Two-dimensional spatial processes}\label{sec:prob:twodim}

We now present the two-dimensional case and try to be as close as possible in notations to the one-dimensional case in order to illustrate the new difficulties. In order to gain some place, we formulate it in the general case of measurable state spaces $(S_1,\ca{S}_1)$ and $(S_2,\ca{S}_2)$ endowed with reference measures on the edges of $\setZ^2$.

\subsubsection{Definition by the law}

The switch from dimension one to dimension two is made easy from proposition~\eqref{prop:MarkovOneDimSecond} since the notion of shell of a partition of a domain also exists in dimension two. The main difference is that the process $X$ must be defined on the edges of the domain (and in dimension $d$, the process should be defined on the $(d-1)$-dimensional objects). In order to keep the structure of the presentation simple and avoid the discussion of domains with two faces, we first start with conditional laws and proceed then with factorization of the conditional expectations.

\begin{defi}[definition of a 2D Markov process from its law]
\label{def:MarkovTwoDimFirst}
Let $(S_1,\ca{S}_1)$ and $(S_2,\ca{S}_2)$ be two measurable sets. Let $D$ be a domain of $\setZ^2$. Let $(\mu_e)_{e\in\Edges{D}}$ be $\sigma$-finite measures on $S_1$ for horizontal edges and on $S_2$ for vertical edges.

An $S$-valued Markov process $X=(X_e)_{e\in \Edges{D}}$ on $D$ is a collection of random variables (on a probability space $(\Omega,\ca{F},\Prob)$) indexed by the set of unoriented edges of $D$ such that:
\begin{enumerate}[(i)]
\item for any horizontal edge $e\in \Edges{D}_h$ (resp. $\Edges{D}_v$), $X_e$ takes values in $S_1$ (resp. $S_2$),
\item there exists a collection $(\MarkovWeight{W}_f)_{f\in D}$ (written $\MarkovWeight{W}_\bullet$ whenever $D$ is clear from context) of measurable functions $\MarkovWeight{W}_f:S_1\times S_1\times S_2\times S_2 \to \setR_+$ such that the conditional law of the process $(X_e)_{e\in\Int(D)}$ given the boundary values $(X_e)_{e\in\partial D}$ admits, almost surely, a conditional density $p$ w.r.t.~the product measure $\otimes_{e\in \Int(D)} \mu_e$ given, for any sequence $(x_e)_{e\in \Edges{D}}$, by
\begin{equation}\label{eq:MarkovLawDimTwo}
p( (x_e)_{e\in\Int(D)} \vert (x_e)_{e\in\partial D} ) = \frac{1}{Z_D(\MarkovWeight{W}_\bullet;x_{\partial D})} \prod_{f\in D} \MarkovWeight{W}_f( x_{b_S(f)},x_{b_N(f)},x_{b_W(f)},x_{b_E(f)}  ) 
\end{equation}
where $Z_D(\MarkovWeight{W}_\bullet;x_{\partial D})$ is a \emph{finite} positive real number, called the \emph{partition function} on the domain $D$.
\end{enumerate}
\end{defi}

It is an easy exercise to check from the normalization of the probability measure that the partition function is given by:
\begin{equation}\label{eq:2Dpartitionfunc:sum}
Z_D(\MarkovWeight{W}_\bullet; x_{\partial D}) = 
\int    \prod_{f\in D} \MarkovWeight{W}_f( x_{b_S(f)},x_{b_N(f)},x_{b_W(f)},x_{b_E(f)}  ) 
\otimes_{e\in \Int(D)} \mu_e(dx_e)
\end{equation}
where each $x_e$ with $e\in\Int(D)$ appears exactly twice in the product of weights $\MarkovWeight{W}_f$, either as one in the north edge of a face and one in the south edge of another face, or one as a west edge and one as a East edge. In dimension one, this summation was interpreted as matrix products. In dimension two, it is the purpose of Section~\ref{sec:operad} to propose a similar powerful algebraic structure.

As in dimension one, boundary weights are absent from the Markov property but are present in the complete law of the process. 

\begin{defi}
A boundary \emph{measure} on a domain $D$ of a Markov process $(X_e)\in\Edges{D}$ is a measure $\mu_{\partial D}$ on $S_1^{\Edges{D}_h } \times S_2^{\Edges{D}_v}$ such that the marginal law of the process $(X_e)\in\Edges{D}$ is given by
\[
\frac{Z_D(\MarkovWeight{W}_\bullet; (x_e)_{e\in\partial D})}{ Z_D^{\boundaryweights}(\MarkovWeight{W}_\bullet; \mu_{\partial D}) }d\mu_{\partial D}( (x_e)_{e\in\partial D} )
\]
with
\begin{equation}\label{eq:proba:ZboundaryfromZdet:meas}
	Z_D^{\boundaryweights}(\MarkovWeight{W}_\bullet;\mu_{\partial D}) = \int   Z_D(\MarkovWeight{W}_\bullet;(x_e)_{e\in\partial D})
	d\mu_{\partial D}( (x_e)_{e\in\partial D} )
\end{equation}
A boundary \emph{weight} on a domain $D$ of a Markov process $(X_e)\in\Edges{D}$ is a measurable function
\[
g: S_1^{\Edges{D}_h } \times S_2^{\Edges{D}_v} \to \setR_+ 
\]
such that the measure  $g \otimes_{e\in\partial D} \mu_e$ is a boundary measure of the process $(X_e)\in\Edges{D}$.
\end{defi}
For a boundary weight $g_D$, we introduce the partition function
\begin{equation}\label{eq:proba:ZboundaryfromZdet}
	Z_D^{\boundaryweights}(\MarkovWeight{W}_\bullet;g_D) = \int  g_D((x_e)_{e\in\partial D}) Z_D(\MarkovWeight{W}_\bullet;(x_e)_{e\in\partial D})
	\otimes_{e\in\partial D} \mu_e(dx_e)
\end{equation}

We choose to introduce boundary measures in addition to boundary weights only for the convenience of being able to describe deterministic boundary conditions, which are useful in practice: indeed, it corresponds to Dirac measures on the boundary values and it may not admit a density w.r.t. the reference measures $\mu_e$, which may be diffuse. Besides this particular case, it has only a limited interest: for any rectangle $R=[a,b]\times [c,d]$, the restriction of the Markov process (see Lemma~\ref{lemma:proba:Markovrestriction} below) to any sub-rectangles admits a boundary weights w.r.t. the reference measures $(\mu_e)$.

As in dimension one, each $(x_e)$ appears exactly twice: once in the partition function and once in the boundary weights $g_D$. In dimension one, we have seen in~\eqref{eq:onedim:lawwithboundary} that the boundary weight generically have a nice factorized structure in terms of two boundary left and right vectors, but the one-dimensional structure of the boundary of a connected domain in dimension two prevents such a simple situation: proposing an nice algebraic description of boundary weights $g_{D}$ in dimension two is precisely the subject of Section~\ref{sec:boundaryalgebra} using the so-called matrix product states or matrix Ans\"atze that have appeared in the physics literature.

\begin{rema}
\label{rema:gaugeinvariance2}
Given a two-dimensional Markov processes, the face weights are uniquely defined only up to multiplication by inverse diagonal matrices on each edge, i.e. for two faces $f_1$ and $f_2$ such that $b_N(f_1)=b_S(f_2)$, for any function $c: S_1\to \setR_+^*$ (resp. $S_2\to\setR_+^*$), the face weights
\begin{align*}
\MarkovWeight{W}'_{f_1}(x,y,w,z) &= \MarkovWeight{W}_{f_1}(x,y,w,z) c(y)		& &\text{(resp. $\MarkovWeight{W}_{f_1}(x,y,w,z) c(z)$)}
\\
\MarkovWeight{W}'_{f_2}(x,y,w,t) &= \MarkovWeight{W}_{f_2}(x,y,w,z) c(x)^{-1} & &\text{(resp. $\MarkovWeight{W}_{f_2}(x,y,w,z) c(w)^{-1}$)}
\end{align*}
define the same Markov process conditioned on its boundaries. Additional factors $c$ must be added to the boundary weights $g_D$ to have the same full law. This observation will be related to the commutative edge structure introduced in Section~\ref{sec:eckmanhilton}.
\end{rema}

\begin{rema}
In dimension one, edge weights are functions of two discrete variables and hence can be identified to matrices (with two indices). In dimension two on the square lattice, face weights are functions of four discrete variables and hence can be identified to tensors with four indices. Matrix products then become tensor contractions as in \eqref{eq:2Dpartitionfunc:sum}. However, the core of Section~\ref{sec:operad} is to embed these operations in an operadic framework: the focus will then be shifted from a computational definition of tensor contractions to a deeper understanding of the composition rules with their associativities and actions on auxiliary spaces.
\end{rema}

We follow closely our previous presentation of the one-dimensional results by introducing some technical lemmata.

\begin{lemm}[stability under restriction]\label{lemma:proba:Markovrestriction}
Let $(X_e)_{e\in\Edges{D}}$ be a two-dimensional Markov process on a domain $D$ of $\setZ^2$ with boundary weight $g_D$. For any subdomain $D'\subset D$, the process $(X_e)_{e\in\Edges{D'}}$ is again a Markov process, whose weight faces are the same (up to restrictions) as for $(X_e)$ and whose boundary weight $g_{D'}$ is given:
\begin{equation}\label{eq:twodim:weightforarestriction}
g_{D'}(x_{\partial D'}) = \int  Z_{D\setminus D'}(\MarkovWeight{W}_\bullet;x_{\partial (D\setminus D')}) g_D(x_{\partial D})
\otimes_{e \in  \partial D\setminus \partial D'} \mu_e(dx_e)
\end{equation}
\end{lemm}
The proof is a simple computation from the law of the process. This is the higher-dimensional equivalent of Lemma~\ref{lemma:proba:Markovrestriction} and Proposition~\ref{prop:onedim:boundaryweight:restrict} with the major difference that the lack of immediate structure of the boundary weights prevents an easy interpretation of \eqref{eq:twodim:weightforarestriction} as an action of matrices on vectors. 

\removable{
\begin{lemm}[factorization over connected components]
Let $D$ be a domain made of $C$ connected components $(D_c)_{1\leq c\leq C}$. Then, for any boundary condition $(x_e)_{e\in \partial D}$, the following factorization holds:
\[
Z_D(\MarkovWeight{W}_\bullet;(x_e)_{e\in\partial D}) = \prod_{c=1}^C
Z_{D_c}(\MarkovWeight{W}_\bullet;(x_e)_{e\in \partial D_c})
\]
\end{lemm}
The purpose of this last lemma is to reduce most computations to the case of connected domains, whose boundary is a unique one-dimensional broken line.
}

\subsubsection{The Markov property.}

We now introduce the key Markov property, which has the advantage, as in dimension one, of characterizing Markov processes without referring to the details of their laws.

\begin{prop}[characterization of a 2D Markov processes from conditional expectation]\label{prop:MarkovTwoDimSecond}
With the same notations as Definition~\ref{def:MarkovTwoDimFirst}, we introduce, for any subset $E$ of $\Edges{D}$, the $\sigma$-algebra $\ca{F}_E$ generated by the family of r.v. $(X_e)_{e\in E}$. Then, for any partition $\mathbf{D}=(D_i)_{1\leq i \leq n}$ of $D$ and any bounded real random variables $(U_i)_{1\leq i\leq n}$ such that, for each $1\leq i\leq n$, $U_i$ is $\ca{F}_{\Edges{D_i}}$-measurable, it holds:
\begin{equation}\label{eq:MarkovDimTwo:CondExp}
\Espc{ \prod_{i=1}^n U_i }{ \ca{F}_{\Shell(\mathbf{D})}}
= \prod_{i=1}^n \Espc{ U_i }{ \ca{F}_{\partial D_i} }
\end{equation}
Moreover, if all the connected components of $D$ have a cardinal larger than $3$, then this property is equivalent to Definition~\ref{def:MarkovOneDimFirst}.
\end{prop}
\begin{proof}
The proof is exactly the same as in the one-dimensional case and is close to the one of the Hammersley-Clifford's theorem. We first assume that $D$ and all the elements $D_i$ of the partition are connected: the non-connected case is obtained by decomposing all the r.v. $U_i$ as finite linear combinations of products of r.v. localized on the connected components. 

Assuming Definition~\ref{def:MarkovTwoDimFirst}, the proposition is a simple consequence of the previous lemmata and summation over values of r.v. on edges to compute all the required conditional expectations. 

When the size of the connected domain $D$ is smaller than $2$, \eqref{eq:MarkovDimTwo:CondExp} is a tautology and does not provide any factorization of the law. When the size of the connected domain is larger than $3$, the same recursion over the domain sizes as in the one-dimensional case holds and we only briefly overview the computations. 

For any connected domain $D$ of size larger than $3$, one can always find a partition $(\mathbf{D}^{(3)}_i)_{i\in \{1,2,3\}})$ such that:
\begin{enumerate}[(i)]
\item  $\mathbf{D}^{(3)}_1$ and  $\mathbf{D}^{(3)}_2$ are connected,
\item  $\mathbf{D}^{(3)}_3=\{f_0\}$ contains a single face,
\item the two distinct partitions $(\mathbf{D}^{(a)}_i)_{i\in \{1,2\}})$ with $a\in\{1,2\}$ defined by
\begin{align*}
\mathbf{D}^{(1)}_1 &= \mathbf{D}^{(3)}_1 \cup \mathbf{D}^{(3)}_3
&
\mathbf{D}^{(1)}_2 &= \mathbf{D}^{(3)}_2 
\\
\mathbf{D}^{(2)}_1 &= \mathbf{D}^{(3)}_1 
&
\mathbf{D}^{(2)}_2 &= \mathbf{D}^{(3)}_2  \cup \mathbf{D}^{(3)}_3
\end{align*}
as illustrated in Figure~\ref{fig:twodimcondexptofactorization} are connected.
\end{enumerate}

\begin{figure}
\begin{center}
\begin{tikzpicture}
\begin{scope}
\node at (1.,1.) {$\mathbf{D}^{(3)}_1$};
\node at (3.,1.) {$\mathbf{D}^{(3)}_2$};
\node at (2.5,2.5) {$f_0$};
\foreach \xs/\ys/\xe/\ye in {0/0/0/2, 0/2/1/2, 1/2/1/3, 1/3/4/3, 4/3/4/0, 4/0/3/0, 3/0/3/-1, 3/-1/2/-1/, 2/-1/2/0, 2/0/0/0, 2/0/2/3, 2/2/3/2, 3/2/3/3}
{
	\draw[thick] (\xs,\ys) to (\xe,\ye);
}
\end{scope}
\begin{scope}[xshift=7cm,yshift=2cm,scale=0.75]
\node at (1.,1.) {$\mathbf{D}^{(1)}_1$};
\node at (3.,1.) {$\mathbf{D}^{(1)}_2$};
\foreach \xs/\ys/\xe/\ye in {0/0/0/2, 0/2/1/2, 1/2/1/3, 1/3/4/3, 4/3/4/0, 4/0/3/0, 3/0/3/-1, 3/-1/2/-1/, 2/-1/2/0, 2/0/0/0, 2/0/2/2, 2/2/3/2, 3/2/3/3}
{
	\draw[thick] (\xs,\ys) to (\xe,\ye);
}
\end{scope}
\begin{scope}[xshift=7cm,yshift=-2cm,scale=0.75]
\node at (1.,1.) {$\mathbf{D}^{(2)}_1$};
\node at (3.,1.) {$\mathbf{D}^{(2)}_2$};
\foreach \xs/\ys/\xe/\ye in {0/0/0/2, 0/2/1/2, 1/2/1/3, 1/3/4/3, 4/3/4/0, 4/0/3/0, 3/0/3/-1, 3/-1/2/-1/, 2/-1/2/0, 2/0/0/0, 2/0/2/3,}
{
	\draw[thick] (\xs,\ys) to (\xe,\ye);
}
\end{scope}
\draw[->] (4.5, 1.5) -- (6.5,3);
\draw[->] (4.5, 1.) -- (6.5,-1);
\end{tikzpicture}
\end{center}
\caption{Three partitions of a domain $D$ as considered in the proof of Proposition~\ref{prop:MarkovTwoDimSecond}: the finest on the left with a single face element, the two coarsest on the right obtained by adding the single face to one of the other elements.}
\label{fig:twodimcondexptofactorization}
\end{figure}

Let $(x_e)_{e\in\Edges{D}}$ be a fixed sequence of values in $S_1$ and $S_2$. The conditional expectation~\eqref{eq:MarkovDimTwo:CondExp} gives for any $a\in\{1,2\}$
\[
\Espc{\prod_{e\in\Edges{D}} \indic{X_e=x_e} }{\ca{F}_{\Shell(\mathbf{D}^{(a)})}} = \left(\prod_{e \in \partial D} \indic{X_e=x_e} \right)V^{(a)}
 \prod_{i\in\{1,2\}} U_i^{(a)}
 \]
with the following definitions of $V^{(a)}$ and $U_i^{(a)}$: 
\begin{align*}
V^{(a)}&=\prod_{e\in\Shell(\mathbf{D}^{(a)})\setminus \partial D}\indic{X_e=x_e}
&
U^{(a)}_i &=\Espc{\prod_{e\in\Edges{\mathbf{D}^{(a)}_i}} \indic{X_e=x_e} }{\ca{F}_{\partial \mathbf{D}^{(a)}_i}}  
\end{align*}
The r.v. $U_i^{(a)}$ can be written as measurable functions $\Upsilon^{(a)}_i$ of $(X_e)_{e\in\partial D \cap \partial \mathbf{D}^{(a)}_i}$ and $(X_e)_{e\in \partial \mathbf{D}^{(a)}_i \setminus \partial D}$. The r.v $\Espc{V^{(a)}}{\ca{F}_{\partial D}}$ can also be written as a measurable function $\Psi^{(a)}$ of $(X_e)_{e\in \partial D}$. So we obtain:
\begin{equation}
\Espc{\prod_{e\in\Edges{D}} \indic{X_e=x_e} }{\ca{F}_{\partial D}} = \left(\prod_{e \in \partial D} \indic{X_e=x_e} \right)
\Psi^{(a)}(x_{\partial D})
 \prod_{i\in\{1,2\}} \Upsilon^{(a)}_i\left( x_{\partial D \cap \partial \mathbf{D}^{(a)}_i}, x_{\partial \mathbf{D}^{(a)}_i  \setminus \partial D}\right)
\end{equation}
The r.h.s. does not depend on the label $a$ and the same argument as in the one-dimensional case (with a careful look at the relative positions of the face $\{f_0\}$ w.r.t. the boundary) provides the existence of a factorization of $\Psi^{(a)}$ over the fine partitions $\mathbf{D}^{(3)}$ with a weight $W_{f_0}$. For size $3$, all the other terms also factorize. For sizes larger than $3$, a quick recursion provides a product structure for the other $\Upsilon^{(a)}_i$.\end{proof}

The factorization in the Markov property based on shells of partitions and boundary conditions is closely related to the operads introduced in Section~\ref{sec:operad}, which will be restricted to rectangles and partitions of rectangles into rectangles through guillotine cuts, in particular it exhibits the same "associativity" in the gluing of elements of partitions, which is, finally, the key feature of the algebraic description of Markov processes as seen below.

\subsection{Finite state space, linear algebra, and some examples.}

Introducing one-dimensional discrete-time and discrete-state-space Markov chains is a classical academic presentation with a low level of mathematics essentially based on linear algebra of matrices and expurgated from any analytical measure theory. However, all the algebraic ingredients are already there. Since the purpose present paper is to describe the algebraic structures behind two-dimensional Markov processes, it is interesting to first formulate all the constructions with discrete state spaces.

When $S_1$ and $S_2$ are finite sets with cardinals $s_1$ and $s_2$, they are canonically endowed with their complete $\sigma$-algebra $\ca{S}_i=\ca{P}(S_i)$ and the counting reference measure $\mu(A)=\card(A)$. In this case, any measure on $S_1^k\times S_2^l$ has a density w.r.t. $\mu_1^{\otimes k}\otimes \mu_{2}^{\otimes l}$ and all the integration steps in the previous formulae are replaced by sums over the sets $S_1$ or $S_2$.

Given a fixed enumeration of each $S_i$ that we will never write, there is an identification between functions $f:S_i\to \setR$ (resp. $\setC$ and the vector space $V(S_i)=\setR^{s_i}$ (resp. $\setC^{s_i})$ through the canonical basis. The notation $V(S_i)$ is introduced to remove the appearance of the base field and to remind the initial space. A face weight $\MarkovWeight{W} : S_1\times S_1\times S_2\times S_2$ is identified to an element $\ha{\MarkovWeight{W}}$ of $\End(V(S_1)) \otimes \End(V(S_2))$ (or $V(S_1)\otimes V(S_1)\otimes V(S_2)\otimes V(S_2)$) through the canonical bases of each space. 

For a rectangular domain $R$ with non-zero horizontal and vertical lengths $P$ and $Q$, a partition function $Z_R(\MarkovWeight{W}; \cdot)$ is a function $S_1^{P}\times S_1^{P}\times S_2^Q\times S_2^Q \to \setR_+$ where the first (resp. second, third, fourth) space corresponds to the $P$ (resp. $P$, $Q$, $Q$) variables on the South (resp. North, West, East) boundary side oriented from left to right (resp. left to right, bottom to top, bottom to top). Hence, $Z_R(\MarkovWeight{W},\cdot)$ is canonically identified with an element of \begin{equation}\label{eq:matrixspaceforZ}\End(V(S_1))^{\otimes P}\otimes \End(V(S_2))^{\otimes Q}\end{equation} where the $i$-th space of the first (resp. second) tensor product corresponds to the $i$-th vertical (resp. horizontal) row, which joins the $i$-th South (resp. West) variable with the $i$-th North variable (resp. East). These identifications are chosen so that gluing horizontally two rectangles with sizes $(P_1,Q)$ and $(P_2,Q)$ along the East side of the first one and the West side of the second one produces a global partition function obtained from the partition functions on the two rectangles by multiplying on $\End(V(S_2))^{\otimes Q}$ and tensorizing $\End(V(S_1))^{\otimes P_1}\otimes \End(V(S_1))^{\otimes P_2}$: this is a trivial observation from the interpretation of the summations in \eqref{eq:2Dpartitionfunc:sum} as matrix products. At first sight, this can be interpreted as contracting suitable tensors to produce new tensors but such an interpretation completely misses the underlying geometric structure of gluings as it will be seen in Section~\ref{sec:operad}. 

A boundary weight $g$ on such a rectangle $R$ with sizes $P$ and $Q$ is, in the same way as $Z_R(\MarkovWeight{W};\cdot)$, a function $S_1^p\times S_1^p\times S_2^q\times S_2^q\to\setR_+$ but may be rather seen as an element of:
\begin{equation}\label{eq:matrixspaceforg}
	(V(S_1)^*\otimes V(S_1))^{\otimes P} \otimes 	(V(S_2)^*\otimes V(S_2))^{\otimes Q}
\end{equation}
with the same correspondence with the rows and the columns of the rectangles as for $Z_R(\MarkovWeight{W};\cdot)$. Then, \eqref{eq:proba:ZboundaryfromZdet}, whose l.h.s. is just a number, can be interpreted as the contraction of each factor of $g$ in \eqref{eq:matrixspaceforg} with the corresponding factor of $Z_R(\MarkovWeight{W};\cdot)$ in \eqref{eq:matrixspaceforZ} through the linear map
\begin{align*}
V(S_i)^*\otimes V(S_i) \otimes \End(V(S_i)) & \to \setR  \\
	a\otimes b\otimes M & \mapsto \scal{a}{Mb}
\end{align*}
Proper operadic definitions of such spaces, multiplications and pairings will be described in Section~\ref{sec:canonicalguill2:markov}.

This framework already encompasses most models of statistical mechanics such as the Ising model and the six-vertex model (two states) and more generally all spin models and the Potts model.

\section{A quick overview of traditional studies of two-dimensional Markov processes}
	\subsection{Large domains in dimension two: asymptotic theorems, Gibbs measures and boundary conditions}\label{sec:gibbs}
	
In most cases, from large-time asymptotics of Markov chains to statistical mechanics in any dimension, one is interested in the study of Markov process on large domains (large time for chains), the so-called thermodynamic limit in statistical mechanics. For  simplicity, we restrict here our discussion to the homogeneous (i.e. all the weights (edge in 1D, face in 2D) are equal) and irreducible (all the configurations have a positive probability) case. For (time-oriented) Markov chains, this corresponds to convergence to the invariant law, ergodic theorems (and their fluctuations), scaling limits (such as Donsker theorem  for random walks). In statistical mechanics, one is interested in the large distance asymptotics of multi-points correlation functions or in the construction of the models on infinite domains through Gibbs measures, which may be a hard task in dimension larger than 1 due to the potential presence of phase transitions.

In dimension one, the Perron-Frobenius eigenvectors, which have positive coordinates, control the thermodynamic limit through the easy Proposition~\ref{prop:onedim:kolmoextension} as well as expectation of one-point observables. The second fact is that such eigenvectors are easy to compute since they are given by an eigenvalue/eigenvector system of equations, for which many algebraic tools are available. It is thus remarkable that large size asymptotic computations are reduced to a purely local eigenvector computation $Au_2=\Lambda u_2$ and $u_1A=\Lambda u_1$. For larger dimensions, to the best of our knowledge, there is no clear, general and simple description of invariant boundary conditions for boundary weights $g_D$. The main obstacle is that, whereas in dimension one boundaries are local, boundaries in larger dimensions are extended global objects with various shapes.

As nicely described in \cite{GibbsVelenik} (chapter 6, section 6.2), defining Gibbs measure with Kolmogorov's extension theorem requires to have a set of \emph{compatible} boundary weights $(g_{D})_D$ for any finite domains, in order to build \emph{compatible} probability measures. In dimension larger than 2, this apparent difficulty replaced by the analytic formalism of states or the DLR formalism. 

These are precisely these preliminary remarks that have motivated us for the constructions of Sections~\ref{sec:operad} and \ref{sec:boundaryalgebra}, in which we propose, for the simpler geometry of rectangles, an \emph{algebraic} construction of \emph{compatible boundary weights} $(g_D)_D$ made of \emph{local} objects determined by \emph{local} equations. We do not know yet whether our construction may provide all the reasonable Gibbs states for a given Markov model but, every time a boundary structure is built, defining an associated Gibbs state is immediate through Kolmogorov's extension theorem. Moreover, as in the one-dimensional case, the boundary structures we provide bring a good help in the computation of multi-point correlation functions.

	\subsection{Some classical computations: other steps towards algebra}\label{sec:classicalcomputations}
	
Before starting the description of the algebraic framework, we present some classical computations in dimension two, which help to solve the previous questions without our formalism.

\subsubsection{The generic behaviour of the partition function.} For an homogeneous model (i.e. all the weights $\MarkovWeight{W}_f$ in dimension two (resp. $\MarkovWeight{A}_e$ in dimension one) are equal to a fixed weight $\MarkovWeight{W}$ (resp. $\MarkovWeight{A}$), the partition function $Z^{2D}_{M,N}$ (resp. $Z^{1D}_{N}$) on a rectangle of sizes $(M,N)$ (resp. on a segment of length $N$) is generically expected to scale as:
\begin{align*}
\frac{1}{MN} \log Z^{2D}_{M,N}  &= f^{2D} + o(1) 
&
\frac{1}{N} \log Z^{1D}_{N}  &= f^{1D} + o(1) 
\end{align*}
when $M$ and $N$ are large. The number $f^{nD}$ is called the \emph{free energy density} and is one of the first quantities of interest.

In dimension one with simple irreducibility hypotheses, $f^{1D}=\log \Lambda$ is the Perron-Frobenius eigenvalue of the transition matrix $\MarkovWeight{A}$ and The associated eigenvectors determine the averages $\Esp{\indic{X_v=x}}$ when $v$ is far from the boundaries. We will see how this is generalized to larger dimensions in Section~\ref{sec:invariantboundaryelmts}.s

\subsubsection{Transfer matrix in dimension two.} A traditional approach to two-dimensional models of statistical mechanics is the so-called transfer matrix formalism (see \cite{baxterbook} for an overview of methods and exact computations). For rectangular domains of size $(M,N)$, the idea is to divide a domain into $N$ (resp. $M$) horizontal (resp. vertical) strips of length $M$ (resp. $N$) and to consider the gluing of successive strips as a one-dimensional model with weight $T_M$ (resp. $T'_N$). In the transverse direction, one must choose either periodic boundary conditions or independent values on the boundary edges: in both cases, it does not correspond to the "natural" boundary conditions coming from infinity but one expects that their effect vanishes in the large size limit.

 The drawback is that the one-dimensional strip weight $T_M$ is a matrix of large size $|S_1|^M\times |S_1|^M$ which may be hard to diagonalize (even for so-called exactly solvable models...). In this case, one has $f^{2D} = \lim_{M\to\infty} (\log \Lambda_M)/M$ where $\Lambda_M$ is the Perron-Frobenius eigenvalue of $T_M$, or similarly $f^{2D} =\lim_{N\to\infty} (\log \Lambda'_N)/N$ for $T'_N$.

Moreover the symmetry breaking between the two directions during the computation may also be a problem hard to encompass (but absent from our construction below): it is not clear why the previous two expressions of $\Lambda_M$ and $\Lambda'_N$ should lead to the same value. 

Finally, the question of the boundary conditions that shall be put on the transverse boundaries, is not trivial: it is usually avoided through periodic boundary conditions and it is not always satisfactory, since large size limits have to be considered first.

The classical book of Baxter \cite{baxterbook} is one of the most common introductions to these approaches.

\subsubsection{Exactly solvable "integrable" systems.} Diagonalization of transfer matrices $T_N$ is a hard task and computations are made partially possible for so-called "integrable systems" (see again \cite{baxterbook} for an overview), which are based on algebraic constructions such as R-matrices that help to build complete families of operators that commute with the transfer matrices. However, most computations use very elaborated algebra to obtain very partial results on very specific models, which appear to be like miracles.

\subsubsection{Matrix Ans\"atze and matrix product states.}

In exactly solvable models, an alternative to Bethe Ans\"atze methods for space-time models such as the asymmetric exclusion process is the so-called Matrix Ansatz, as introduced in \cite{DEHP}. Despite its generalizations (see for example \cite{BlytheEvans}), it powerful computational implications (see for example \cite{ASEPBrownianExc,LargeDevDensityASEP,CrampeRagoucySimon},) and its connections to other domains such as combinatorics (see \cite{CorteelWilliams}), the Matrix Ansatz method remains a mysterious object, hard to construct and to apply widely. In particular, it was not clear to us whether it is directly related to integrable systems or not. 

Besides these theoretical perspectives, matrix product states have become over the past few fifteen years a very powerful approximation and numerical tool in physics, as nicely explained in \cite{MPSreview} and \cite{DMRGMPS}, without attracting a large interest from mathematicians. We deeply think that there must be nice mathematical grounds for these successes. 

Our constructions of section \ref{sec:operad} and \ref{sec:boundaryalgebra} are essentially inspired from Matrix Ans\"atze and matrix product states and lift them to a \emph{general structural property} of two-dimensional Markov processes. In this spirit, integrable systems are just a helping tool to build them. we also hope that the tools proposed in Section~\ref{sec:boundaryalgebra} may help to understand the nice approximation property of matrix product states and provide later new algorithmic tools.

\chapter[Coloured operads for Markov processes]{The algebraic framework of coloured operads for Markov processes}\label{sec:operad}

	\section{The language of operads}\label{sec:operadiclanguage}
		
\subsection{Definitions}

\subsubsection{Without colours.}
We give here a brief definition of operads, mostly for the probabilist reader: all the definitions are classical algebraic definitions. We emphasize on the less used notion of coloured operad, since the colours will play a fundamental role in the next constructions and are not just a required technical specification. The interested reader may refer to \cite{mayoperad} for a historical starting point and to \cite{stasheffoperad} for a global overview.

\begin{defi}[operad]
A (symmetric) operad is a sequence of sets $(P(n))_{n\in\setN}$ (an element of $P(n)$ is called an $n$-ary operation) such that:
\begin{enumerate}[(i)]
\item for any $n\in\setN$ and any finite sequence $(k_i)_{1\leq i\leq n}$ in $\setN$, there is a \emph{composition function}
\[
\begin{split}
\circ : P(n)\times P(k_1)\timesdots P(k_n) & \to P(k_1+\ldots+ k_n) \\
(\alpha,\alpha_1,\ldots,\alpha_n) &\mapsto \alpha\circ(\alpha_1,\ldots,\alpha_n)
\end{split}
\]
\item there is an element $\identity \in P(1)$ (called the identity) such that, for any $\alpha\in P(n)$,
\begin{align*}
\alpha\circ(\identity,\ldots,\identity) &= \alpha \\
\identity \circ\,\alpha &= \alpha
\end{align*}
\item for any $n\in\setN$, any $(k_i)_{1\leq i\leq n}$ and any $((l_{i,j})_{1\leq j\leq k_i})_{1\leq i\leq n}$, the composition functions satisfy the following associativity rules for any operations $\alpha\in P(n)$, $\alpha_{i}\in P(k_i)$ and $\alpha_{i,j}\in P(l_{i,j})$:
\begin{equation}\label{eq:operad:compoasso}
\begin{split}
(\alpha\circ &(\alpha_1,\ldots,\alpha_n))\circ (\alpha_{1,1},\ldots,\alpha_{1,l_{1,k_1}},\alpha_{2,1},\ldots,\alpha_{n,l_{n,k_n}}) \\
&= \alpha \circ (
\alpha_1\circ(\alpha_{1,1},\ldots,\alpha_{1,k_1}),
\ldots,
\alpha_n\circ(\alpha_{n,1},\ldots,\alpha_{n,k_n}))
\end{split}
\end{equation}
\item there is an action of the symmetric group $\Sym{n}$ on $P(n)$ such that \begin{equation}\label{eq:operad:symmaction}\alpha_\sigma\circ(\alpha_1,\ldots,\alpha_n)=\alpha\circ(\alpha_{\sigma^{-1}(1)},\ldots,\alpha_{\sigma^{-1}(n)})
\end{equation}
\end{enumerate}
\end{defi}

\begin{rema}[Cartesian products versus tensor products]
	In the case where the spaces $(P(n))$ are vector spaces, we may require the compositions and the actions of the symmetric group to be multilinear and thus lift all the Cartesian products to tensor products. If it is not the case, we may still consider the vector spaces generated by the elements of $P(n)$. 
\end{rema}

\begin{defi}[algebra over an operad]\label{def:operad:algebra}
Let $P=(P(n))_{n\in\setN}$ be an operad. An algebra over $P$ is a set $\ca{A}$ such that:
\begin{enumerate}[(i)]
\item for any $n\in\setN$, there are maps $P(n)\times \ca{A}^n \to \ca{A}$, $(\alpha,u_1,\ldots,u_n)\mapsto m_\alpha(u_1,\ldots,u_n)$,
\item for any $u\in \ca{A}$, $m_{\id}(u)=u$
\item for any $(\alpha,\alpha_1,\ldots,\alpha_n)\in P(n)\times P(k_1)\timesdots P(k_n)$ and any sequence $(u_i)_{1\leq i\leq k_1+\ldots+k_n}$, it holds:
\begin{equation}
m_{\alpha\circ(\alpha_1,\ldots,\alpha_n)}(u_1,\ldots,u_{k_1+\ldots+k_n}) =
m_\alpha(v_1,\ldots,v_n)
\end{equation}
with elements $v_j\in\ca{A}$ given by the following partial maps
\begin{equation}
	v_j=m_{\alpha_j}(u_{k_1+\ldots+k_{j-1}+1},\ldots,u_{k_1+\ldots+k_j})
\end{equation}
\item for any $\alpha\in P(n)$, any $\sigma\in\Sym{n}$ and any $(u_1,\ldots,u_n)\in \ca{A}$, \begin{equation}m_{\alpha^\sigma}(u_1,\ldots,u_n)=m_\alpha(u_{\sigma^{-1}(1)},\ldots,u_{\sigma^{-1}(n)})
\end{equation}
\end{enumerate} 
\end{defi}

Operads can be viewed as the common structure behind any sets with internal laws sharing the same axioms, such as commutativity, associativity, etc. A map $m_\alpha :\ca{A}^n\to\ca{A}$ is called a \emph{product} We provide some examples below.

\subsubsection{With colours.}

In many cases, products, and hence compositions, may not be defined for all the objects in an algebra $\ca{A}$ over an operad $(\ca{P}_n)$. We thus have to split further the space $\ca{A}$ according to some criteria which lead to the notion of \emph{coloured} operad, in which the spaces are indexed by so-called \emph{colours}. The interested reader may refer to \cite{YauColoredOperads} for a general overview; we recall here only the basic definitions required by our construction.

\begin{defi}[coloured operad]
Let $C$ be a set, called the set of \emph{colours}. A coloured operad over $C$ is a sequence of sets $(P(c;c_1,\ldots,c_n)_{(c,c_1,\ldots,c_n)\in C^{n+1}})_{n\in\setN}$ such that:
\begin{enumerate}[(i)]
\item for any $n\in\setN$, any colours $(c,c_1,\ldots,c_n)\in C^{n+1}$ and any finite sequence of colours $((c_i,c_1^{(i)},\ldots,c_{k_i}^{(i)}))_{1\leq i\leq n}$, there is a \emph{composition function}:
\begin{equation}\label{eq:defcompositionmaps}
\begin{split}
\circ : P(c;c_1,\ldots,c_n) \times \left( 
\prod_{i=1}^n P(c_i;c^{(i)}_1,\ldots,c^{(i)}_{k_i})
\right)
& \to P(c;c^{(1)}_1,\ldots,c^{(1)}_{k_1},c^{(2)}_1,\ldots,c^{(n)}_{k_n}) \\
(\alpha,\alpha_1,\ldots,\alpha_n) & \mapsto \alpha\circ (\alpha_1,\ldots,\alpha_n)
\end{split}
\end{equation}
\item for any $c\in C$, there is an element $\id_c\in P(c;c)$ such that, for $\alpha\in P(c;c_1,\ldots,c_n)$, $\alpha\circ(\id_{c_1},\ldots,\id_{c_n})=\alpha$ and $\id_c\circ \alpha=\alpha$;
\item \eqref{eq:operad:compoasso} holds for any maps $\alpha$, $\alpha_i$ and $\alpha_{i,j}$ with suitable colours for the maps and the composition to be defined;
\item for any $P(c;c_1,\ldots,c_n)$, there is an action of the subgroup of of permutations $\sigma\in\Sym{n}$ such that, for any $\alpha\in P(c;c_1,\ldots,c_n)$,  $\alpha^\sigma\in P(c;c_{\sigma(1)},\ldots,c_{\sigma(n)})$ and such that \eqref{eq:operad:symmaction} holds.
\end{enumerate}
\end{defi}

\begin{defi}[algebra over a coloured operad]\label{def:colouroperad:algebra}
An \emph{algebra over a coloured operad} $(P(c;c_1,\ldots,c_n))$ over a set of colours $C$ is a collection of sets $(\ca{A}(c))_{c\in C}$ such that, for any $n\in\setN$, there are maps 
\begin{equation}\label{eq:defproductmaps}
\begin{split}
P(c;c_1,\ldots,c_n)\times \ca{A}(c_1) \timesdots \ca{A}(c_n) & \to \ca{A}(c) \\
(\alpha,u_1,\ldots,u_n) &\mapsto m_\alpha(u_1,\ldots,u_n)
\end{split}
\end{equation}
satisfying the same conditions as in Definition~\ref{def:operad:algebra}, up to colour compatibility conditions.
\end{defi}

In practice, the colours encode what type of operations can be made between objects of different natures (i.e. elements of various spaces $\ca{A}_c$) and, if a space $P(c;c_1,\ldots,c_n)$ is empty (this is the case below), this means that the corresponding objects cannot be combined with some operations or products. A typical example corresponds to rectangular matrices whose colours are their shape $c=(n,m)$, for which multiplications exists only if the number of columns on a matrix matches the number of lines of the next one.

\subsubsection{Notation.} We defined partial compositions $\circ_i : P(n)\times P(m) \to P(n+m-1)$ by $\alpha\circ_i \beta= \alpha\circ(\id,\ldots,\id,\beta,\id,\ldots,\id)$ with $\beta$ in $i$-th position. Specifying axioms on partial compositions or the complete composition functions $\circ$ is fully equivalent. The same notations may be adopted for coloured operads whenever the colours of the various elements are compatible with composition.

\subsubsection{Linear structure and tensor products.}

In all the previous definitions, the sets $P(c;c_1,\ldots,c_n)$ and $\ca{A}_c$ are generic sets without extra structure. In most cases below (but not all), they are vector spaces. In this case, we require the composition maps $\circ$ in \eqref{eq:defcompositionmaps} as well as the product maps $m_\alpha$ in \eqref{eq:defproductmaps} to be \emph{linear}: they are thus canonically extended to the tensor products that replaces the Cartesian products in \eqref{eq:defcompositionmaps} and \eqref{eq:defproductmaps}. Excepted if explicitly mentioned all the algebras over an operad that appear below have a linear structure and only tensor products will be used.

		\subsection{Related examples in one- and two-dimensional geometries}
		
The coloured operad of the next section relevant for two-dimensional Markov processes has a geometric structure, which is the case of other well-known operads, and mixes several features of these operads. In order to inspire intuition to the reader, we present here briefly these basic cases.

\subsubsection{The operad \texorpdfstring{$\mathrm{Com}$}{Com} and \texorpdfstring{$\Ass$}{Ass}: zero- and one-dimensional view} The associative operad $\Ass$ (for "associative") is the operad corresponding to monoids (set case) and algebras (vector space case). It consists of objects $\Ass(n) = \{m_n\}\times \Sym{n}$ (for $n\geq 1$) such that $(m_n,\tau) = (m_n,\id)^\tau$ (action of the symmetric group). We set by convention $m_n:= (m_n,\id)$ and define the composition functions by $m_n(m_{k_1},\ldots,m_{k_n})= m_{k_1+\ldots+k_n}$. One can check that it is indeed generated by $\id\in\Ass(1)$ and $m_2\in\Ass(2)$ with the associativity condition 
\begin{equation}\label{eq:associativity1D}
m_2(m_2,\id)=m_2(\id,m_2)
\end{equation} An algebra $\ca{A}$ over $\Ass$ is just a monoid or an algebra where $m_2$ corresponds to the usual product $m_2(a,b)=ab$. 

The one-dimensional view emerges from the fact that, up to the action of symmetric group to relabel the arguments, $m_n$ is an operation that maps a one-dimensional oriented sequence $(a_1,\ldots,a_n)\to m_n(a_1,\ldots,a_n)$ and $m_{n+1}$ can be obtained from $m_n$ by either left or right multiplication. If one sees the elements $a_i$ as attached to some segments, the multiplication $m_n$ corresponds to the successive gluing of the segments in the correct order.

The operad $\mathrm{Com}$ (for "commutative") is even simpler and corresponds to $\mathrm{Com}_n=\{m_n\}$ and a trivial action of $\Sym{n}$. This corresponds to a zero-dimensional view since the multiplication $(a_1,\ldots,a_n)\to m_n(a_1,\ldots,a_n)$ does not take into account the ordering of the elements $a_i$, in the same manner as gluing $n$ points on a single point erases any order.

These dimensional interpretations will become much clearer for larger dimensions in the next sections.

\subsubsection{The little \texorpdfstring{$d$}{d}-cubes operad \texorpdfstring{$E_d$}{Ed}.} For any $d\geq 1$, the little $d$-cubes operad $(E_d(n))_{n\geq 1}$ is made of elements $\alpha \in E_d(n)$ such that $\alpha=(R_1,\ldots,R_n)$ where the $(R_i)$ are rectangles $\times_{1\leq i\leq d}[u_i,v_i]\subset [0,1]^d$ with disjoint interiors. The action of the symmetric group is the permutation of the rectangles. The composition laws are given by:
\[(R_1,\ldots,R_n)\circ_i (R'_1,\ldots,R'_m) = 
(R_1,\ldots,R_{i-1},s_{R_n}(R'_1),\ldots,s_{R_n}(R'_m),R_{i+1},\ldots,R_n)\]
where $s_R(R')$ is the image of $R'$ by the unique affine transformation $(x_1,\ldots,x_d)\mapsto (\alpha_1 x_1+\beta_1,\ldots,\alpha_d x_d+\beta_d)$ sending $[0,1]^d$ to $R$.

Such operads are well-studied in algebraic topology from the point of view of homotopy. Moreover, up to homotopies, rectangles can be replaced by any convex shapes, without changing many algebraic topological properties. The topological nature of the operad, which is a problem for us (see below), is related to the presence of the maps $s_R$ that discards the size of the rectangles during compositions: the next section presents a closely related operad where the introduction of colours related to sizes replace the scaling transformation $s_R$. 

It is interesting to notice that the operad $E_1$ is related to $\Ass$ through the mapping $(R_1,\ldots,R_n)\mapsto (m_n,\sigma)$ where $\sigma$ is the permutation that orders the rectangles $R_1$,\ldots,$R_n$ from left to right.

	\section{The guillotine partition operad for Markov processes}	\label{sec:guillotinebasics}

We choose to start with the two-dimensional case since it is more graphical and there are less possible confusions with other well-known operads and the geometrical nature appears in a clearer way.

		\subsection{Definition and elementary properties in the two-dimensional case}

The present paper deals with discrete space Markov processes on $\setZ^2$ as well as some continuous extensions to processes on $\setR^2$ with additional analytical considerations as presented in Section~\ref{sec:higherdim} below. From an operadic point of view, all the content of the next definitions (excepted the discrete examples on $\setZ^2$ obtained as tensor products of matrices) may be formulated in $\setZ^2$ or in $\setR^2$. We thus provide directly general definitions within two abstract spaces $\setP$ and $\setL$ ($P$ stands for points and $L$ for lengths) which can be chosen as
\begin{itemize}
	\item $\setP=\setZ$, $\setL=\setN_0$ and $\setL^* =\setN_1$ if the space on which the Markov processes are defined is $\setZ^2$
	\item $\setP=\setR$, $\setL=\setR_+$ and $\setL^* =\setR_+^*$ if the space on which the Markov processes are defined is $\setR^2$
\end{itemize}
In both cases, we have the following properties:
\begin{enumerate}[(i)]
	\item there is a commutative addition law $\setL \times \setL \mapsto \setL$ with neutral element $0$;
	\item there is a total order $\leq$ on $\setP$;
	\item there is a total order $\leq$ on $\setL$ such that, for all $x,y\in\setL$, one has $x\leq x+y$ with an equality if and only if $y=0$;
	\item for any $(a,b)\in \setP$, such that $a<b$, there exists $u\in \setL$ such that $b=a+u$ (by abuse of notation, we write it $u=b-a$).
\end{enumerate}
and they are enough for all the definitions related to rectangles and guillotine cuts that we now introduce.

\subsubsection{Guillotine cuts and partitions.}
A segment in $\setP$ is a set $[u,v] = \{w\in\setP ; u \leq w \leq v\}$.

A rectangle in $\setP^2$ is a set $[u_1,v_1]\times [u_2,v_2]$ with $u_1\leq u_2$ and $v_1\leq v_2$. A degenerate rectangle corresponds to an equality case $u_1=v_1$ (vertical segment) or $u_2=v_2$ (horizontal segment). A doubly degenerate rectangle corresponds to $u_1=v_1=u$ and $u_2=v_2=v$, i.e. this is just the singleton $\{(u,v)\}$.

In relation with the geometric considerations of the previous section, we have that a rectangle $R=[u_1,v_1]\times [u_2,v_2]$ where $u_1$, $u_2$, $v_1$, $v_2$ are integers in $\setZ$ with $u_i < v_i$ corresponds to the domain (set of faces) $D_R=\{ [k_1,k_1+1]\times [k_2,k_2+1]; u_i\leq k_i< v_i\}$. If, for some $u_1=v_1$ or $u_2=v_2$, then the domain does not contain any face and only edges and vertices.

\begin{defi}[guillotine cut]
Let $R=[u_1,v_1]\times [u_2,v_2]$ be a rectangle in $\setP^2$. A horizontal (resp. vertical) guillotine cut on $R$ is a pair $(1,w)$ (resp. $(2,w)$) with $w\in\setP$ such that $u_1\leq w\leq v_1$ (resp. $u_2\leq w\leq v_2$). The cut is strict if the inequalities are strict.
\end{defi}
Prescribing a guillotine cut $(i,w)$ divides $R$ into two rectangles 
$R_1=[u_1,v_1]\times [u_2,w]$ and $R_2=[u_1,v_1]\times [w,v_2]$ for an horizontal cut and $R_1=[u_1,w]\times [u_2,v_2]$ and $R_2=[w,v_1]\times [u_2,v_2]$ for a vertical one, such that the union of the two rectangles is $R$ and they are not overlapping, i.e. their intersection does not contain any face (only edges and vertices) and has zero area (we prefer not to talk about interiors in order to encompass the case of degenerate rectangles). 

The \emph{shape} of a rectangle (degenerate or not) $[u_1,v_1]\times [u_2,v_2]$ is the element $(v_1-u_1,v_2-u_2)\in\setL^2$.

\begin{defi}[guillotine partition]\label{def:guillotinepartition2D}
Let $R=[u_1,v_1]\times [u_2,v_2]$ be a rectangle in $\setP^2$. A guillotine partition of size $n$ of $R$ is a finite $n$-uplet of rectangles $(R_1,\ldots,R_n)$ such that:\begin{enumerate}[(i)]
\item \label{item:hyp:disjointint} for all $i\neq j$, the rectangles $R_i$ and $R_j$ do not overlap (i.e. their intersection has zero area),
\item \label{item:hyp:union} $R=\cup_{i=1}^n R_i$,
\item \label{item:hyp:existsguillcut}if $n\geq 2$, there exists a permutation $\sigma\in \Sym{n}$, a guillotine cut $(i,w)$ of $R$ that divides $R$ into two rectangles $S$ and $T$ and an index $1\leq k\leq n-1$ such that $(R_{\sigma(1)},\ldots,R_{\sigma(k)})$ and $(R_{\sigma(k+1)},\ldots,R_{\sigma(n)})$ are respective guillotine partitions of $S$ and $T$.
\end{enumerate}
\end{defi}

This definition is recursive in the size of the partition in the third item. The triplet $(\sigma,(i,w),k)$ may also be not unique: one can imagine a square cut into four squares with two cuts, one horizontal and one vertical, in the middle of the edges.

If the rectangles are not degenerate, the last item defines uniquely $k$; on the contrary, if a rectangle $R_i$ is degenerate and coincide with the guillotine cut, then it can be added either in the $S$ part or in the $T$ part.

Due tu the first point, only degenerate rectangles can appear several times in a guillotine partition.

\begin{figure}
\begin{center}
\begin{tikzpicture}[scale=0.6]
\draw (0,0) rectangle (4,5);
\draw (0,0) rectangle (2,2.5) node [pos=0.5] {$1$};
\draw (0,2.5) rectangle (2,5) node [pos=0.5] {$3$};
\draw (2,0) rectangle (4,2) node [pos=0.5] {$5$};
\draw (2,2) rectangle (3,5) node [pos=0.5] {$2$};
\draw (3,2) rectangle (4,4) node [pos=0.5] {$4$};
\draw (3,4) rectangle (4,5) node [pos=0.5] {$6$};
\draw[ultra thick] (2,0) -- (2,5);
\end{tikzpicture}
\hskip 1cm
\begin{tikzpicture}[scale=0.6]
\draw (0,0) rectangle (4,5);
\draw (0,0) rectangle (1,3) node [pos=0.5] {$1$};
\draw (1,0) rectangle (4,1) node [pos=0.5] {$2$};
\draw (0,3) rectangle (2,5) node [pos=0.5] {$3$};
\draw (2,1) rectangle (4,5) node [pos=0.5] {$4$};
\draw (1,1) rectangle (2,3) node [pos=0.5] {$5$};
\end{tikzpicture}
\end{center}
\caption{Examples of a guillotine partition with non-degenerate rectangles on the left and a non-guillotine partition on the right. The first guillotine cut (third item of Definition~\ref{def:guillotinepartition2D}) on the left is vertical and displayed in thick.\label{fig:guillotineexamples}}
\end{figure}

\begin{defi}[equivalence under translations]\label{def:guillpart:equivtranslat}
Two guillotine partitions $(R_1,\ldots,R_n)$ and $(S_1,\ldots,S_m)$ of two rectangles $R$ and $S$ are equivalent if $n=m$ and there exists a translation $\theta$ on $\setP^2$ such that $\theta(R_i)=S_i$ for all $1\leq i\leq n$ (and thus $\theta(R)=S$).

A \emph{guillotine partition class} with shape $((p,q),(p_1,q_1),\ldots,(p_n,q_n))$ is an equivalence class $[\rho]$ under this translation equivalence of guillotine partitions $\rho=(R_1,\ldots,R_n)$ of a rectangle $R$ such that $R$ has shape $(p,q)$ and each $R_i$ has shape $(p_i,q_i)$.
\end{defi}

\begin{lemm}[composition of guillotine partitions by nesting]\label{lemma:guillotinecompo}
If $\rho=(R_1,\ldots,R_n)$ is a guillotine partition of $R$ and, for any $1\leq i\leq n$, $\rho_i=(R^{(i)}_1,\ldots,R^{(i)}_{k_i})_{1\leq i\leq n}$ is a guillotine partition of $R_i$, then \[
(R^{(1)}_{1},\ldots,R^{(1)}_{k_1},\ldots,R^{(n)}_{1},\ldots,R^{(n)}_{k_n})\] is again a guillotine partition of $R$, called the composition $\rho\circ(\rho_1,\ldots,\rho_n)$ of the partitions. This composition (or nesting) extends to guillotine partitions \emph{classes} by considering composable representatives of the equivalence classes.
\end{lemm}
\begin{proof}
This can be proved quickly by a recursion from the definitions. A detailed proof will be given in a more general setting for Lemma~\ref{lemma:equivandcompo:extendedguill:notranslation}.
\end{proof}
Such an example of nesting is displayed in Figure~\ref{fig:compoguillotinepartitions}.

\subsubsection{The guillotine operad: definition and notations.}
We now introduce a coloured operad whose colours are precisely the set of rectangle shapes (hence the acronym $\patterntype{r}$):
\begin{equation}
\PatternShapes(\patterntype{r}) = \setL^2
\end{equation}
i.e. sizes (horizontal and vertical) of rectangles of $\setP^2$. The acronym $\PatternShapes$ stands for "pattern shapes" and other geometries will be considered below. It includes degenerate rectangles: a particular discussion of degenerate rectangles is presented in Section~\ref{sec:eckmanhilton}.

\begin{defi}[guillotine partition operad]\label{def:guillotineoperad:main}
The guillotine partition operad $\Guill_2$ on $\setP^2$ is the coloured operad defined as follows:
\begin{enumerate}[(i)]
	\item the set of colours is $\PatternShapes(\patterntype{r})$;
	\item the set $\Guill_2((p,q);(p_1,q_1),\ldots,(p_n,q_n))$ of $n$-ary operations is the set of guillotine partition equivalence classes with shape $((p,q);(p_1,q_1),\ldots,(p_n,q_n))$;
	\item the composition functions are given by the nestings defined in Lemma~\ref{lemma:guillotinecompo};
	\item the symmetric group acts by permuting the rectangles of the partitions.
\end{enumerate} 
\end{defi}
Checking the axioms of an operad is  almost trivial from the previous Lemma~\ref{lemma:guillotinecompo}. One also remarks that the same definition may be extended to the continuous setting directly.

Depending on the context, we may write the action of a guillotine partition $\rho=(R_1,\ldots,R_n)$ on elements $(a_1,\ldots,a_n)$ of a $\Guill_2$-algebra $(\ca{A}_{p})_{p\in\PatternShapes(\patterntype{r})}$ 
\begin{itemize}
\item either by $m_{\rho}(a_1,\ldots,a_n)$ to recall the similarity with products and to stick the previous notations,
\item or simply by the drawing of $\rho$ applied on $(a_1,\ldots,a_n)$ as in the example:
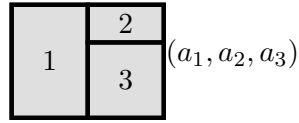
\[
\begin{tikzpicture}[guillpart,xscale=2]
\fill[guillfill] (0,0) rectangle (2,3);
\draw[guillsep] (0,0) rectangle (2,3);
\draw[guillsep] (0,0) rectangle (1,3) node [pos=0.5] {$1$};
\draw[guillsep] (1,0) rectangle (2,2) node [pos=0.5] {$3$};
\draw[guillsep] (1,2) rectangle (2,3) node [pos=0.5] {$2$};
\end{tikzpicture}(a_1,a_2,a_3)
\]
to illustrate compositions as in Figure~\ref{fig:compoguillotinepartitions}
\item or simply by the drawing of $\rho$ with the $a_i$ inside the corresponding rectangles, as in the example:
\[
\begin{tikzpicture}[guillpart,xscale=2]
\fill[guillfill] (0,0) rectangle (2,3);
\draw[guillsep] (0,0) rectangle (2,3);
\draw[guillsep] (0,0) rectangle (1,3) node [pos=0.5] {$a_1$};
\draw[guillsep] (1,0) rectangle (2,2) node [pos=0.5] {$a_3$};
\draw[guillsep] (1,2) rectangle (2,3) node [pos=0.5] {$a_2$};
\end{tikzpicture}
\]
\end{itemize}

\begin{figure}
\begin{align*}
\begin{tikzpicture}[guillpart]
\fill[guillfill] (0,0) rectangle (4,5);
\draw[guillsep] (0,0) rectangle (4,5);
\draw[guillsep] (0,0) rectangle (2,2.5) node [pos=0.5] {$1$};
\draw[guillsep] (0,2.5) rectangle (2,5) node [pos=0.5] {$3$};
\draw[guillsep] (2,0) rectangle (4,2) node [pos=0.5] {$5$};
\draw[guillsep] (2,2) rectangle (3,5) node [pos=0.5] {$2$};
\draw[guillsep] (3,2) rectangle (4,4) node [pos=0.5] {$4$};
\draw[guillsep] (3,4) rectangle (4,5) node [pos=0.5] {$6$};
\end{tikzpicture}
&=  \sigma\left(\;
\begin{tikzpicture}[guillpart]
\fill[guillfill] (0,0) rectangle (4,5);
\draw[guillsep] (0,0) rectangle (4,5);
\draw[guillsep] (0,0) rectangle (2,2.5) node [pos=0.5] {$1$};
\draw[guillsep] (0,2.5) rectangle (2,5) node [pos=0.5] {$2$};
\draw[guillsep] (2,0) rectangle (4,2) node [pos=0.5] {$3$};
\draw[guillsep] (2,2) rectangle (3,5) node [pos=0.5] {$4$};
\draw[guillsep] (3,2) rectangle (4,4) node [pos=0.5] {$5$};
\draw[guillsep] (3,4) rectangle (4,5) node [pos=0.5] {$6$};
\end{tikzpicture}\;\right)
, \qquad \sigma = \begin{pmatrix} 1 & 2 & 3 & 4 & 5 & 6 \\ 1 & 4 & 2 & 5 & 3 & 6\end{pmatrix}\in\Sym{6}
\\
\begin{tikzpicture}[guillpart]
\fill[guillfill] (0,0) rectangle (4,5);
\draw[guillsep] (0,0) rectangle (4,5);
\draw[guillsep] (0,0) rectangle (2,2.5) node [pos=0.5] {$1$};
\draw[guillsep] (0,2.5) rectangle (2,5) node [pos=0.5] {$2$};
\draw[guillsep] (2,0) rectangle (4,2) node [pos=0.5] {$3$};
\draw[guillsep] (2,2) rectangle (3,5) node [pos=0.5] {$4$};
\draw[guillsep] (3,2) rectangle (4,4) node [pos=0.5] {$5$};
\draw[guillsep] (3,4) rectangle (4,5) node [pos=0.5] {$6$};
\end{tikzpicture}
&=  
\begin{tikzpicture}[guillpart]
\fill[guillfill] (0,0) rectangle (4,5);
\draw[guillsep] (0,0) rectangle (4,5);
\draw[guillsep] (0,0) rectangle (2,5) node [pos=0.5] {$1$};
\draw[guillsep] (2,0) rectangle (4,2) node [pos=0.5] {$2$};
\draw[guillsep] (2,2) rectangle (4,5) node [pos=0.5] {$3$};
\end{tikzpicture}
\circ \left(\;
\begin{tikzpicture}[guillpart]
\fill[guillfill] (0,0) rectangle (2,5);
\draw[guillsep] (0,0) rectangle (2,5);
\draw[guillsep] (0,0) rectangle (2,2.5) node [pos=0.5] {$1$};
\draw[guillsep] (0,2.5) rectangle (2,5) node [pos=0.5] {$2$};
\end{tikzpicture}
\;, \id,\;
\begin{tikzpicture}[guillpart]
\fill[guillfill] (0,0) rectangle (2,3);
\draw[guillsep] (0,0) rectangle (2,3);
\draw[guillsep] (0,0) rectangle (1,3) node [pos=0.5] {$1$};
\draw[guillsep] (1,0) rectangle (2,2) node [pos=0.5] {$2$};
\draw[guillsep] (1,2) rectangle (2,3) node [pos=0.5] {$3$};
\end{tikzpicture}
\;\right)
\end{align*}
\caption{Composition rules of guillotine partitions: the guillotine partition on the left of Figure~\ref{fig:guillotineexamples} can be seen as the composition of two more elementary guillotine partitions.}\label{fig:compoguillotinepartitions}
\end{figure}

\begin{prop}[orbit of a partition under permutations] Given a guillotine partition $\rho=(R_1,\ldots,R_n)$ and a permutation $\sigma$, the permuted partition $(R_{\sigma^{-1}(i)})_{1\leq i\leq n}$ is equal to $\rho$ if and only if, for every non-degenerate cycle $c:i_1\mapsto i_2\mapsto \ldots \mapsto i_k \mapsto i_1$ of $\sigma$ with $k\geq 2$, all the rectangles $R_{i_l}$, $1\leq l\leq k$, are equal to a same \emph{degenerate} rectangle $R_c$.
\end{prop}

One must be careful that the relative positions of the rectangles matter and not only their size. In particular, we have that the partition of $[0,2p]\times[0,q]$ into $R_1=[0,p]\times[0,q]$ and $R_2=[p,2p]\times[0,q]$ is different from the permuted permutation $(R_2,R_1)$ since $R_2$ and $R_1$ are different although having the same size.

\begin{prop}
	We see that that any non-degenerate rectangle can appear only once due to item (i) of Definition~\ref{def:guillotinepartition2D} and thus $\sigma$ can only permute identical degenerate rectangles.
\end{prop}
This proposition has strong consequences in terms of commutativity as it will be seen below.

\subsubsection{Fundamental scalar example.}\label{par:trivialoperad} The simplest example of an algebra over the operad $\Guill_2$ is given by the algebras $\ca{A}(p,q)=\setK$ where $\setK$ is any commutative field with composition rules  \[
m_{(R_1,\ldots,R_n)}(u_1,\ldots,u_n)=u_1\ldots u_n \]
for any guillotine partition $(R_1,\ldots,R_n)$ and any scalars $u_i$. In this case, neither the relative positions nor the size of the rectangles play a role and it illustrates only the $n$-arity of such a partition. This canonical trivial $\Guill_2$-operads will be simply written $(\setK_{p,q})_{(p,q)\in\PatternShapes(\patterntype{r})}$ or $\setK_{\PatternShapes(\patterntype{r})}$ to distinguish it from the field $\setK$ itself.  This example will reappear in the Section~\ref{sec:eckmanhilton} about the Eckmann-Hilton argument and in the proof of Theorem~\ref{theo:stability} and is closely related to the the notion of eigenvalue of Section~\ref{sec:boundaryalgebra}.

\subsubsection{Generators with arity two.}

\begin{prop}\label{prop:guill2generators}
The operad $\Guill_2$ is generated by the action of the symmetric group, the identities $\id \in \Guill_2((p,q);(p,q))$ and the elementary guillotine partitions $([0,r]\times[0,q],[p-r,p]\times[0,q])$ and $([0,p]\times[0,s],[0,p]\times [s,q])$ for any $(p,q)\in\PatternShapes(\patterntype{r})$, $0\leq r \leq p$ and $0\leq s \leq q$.
\end{prop}
Since the colours are always clear from context, one chooses to introduce the simplified notations:
\begin{subequations}
\label{eq:defgenerators}
\begin{align}
m_{WE}^{r,p-r | q}&=m_{([0,r]\times[0,q],[p-r,p]\times[0,q])} \in \Guill_2((p,q);(r,q),(p-r,q))
\\
m_{SN}^{p | s,q-s}&=m_{([0,p]\times[0,s],[0,p]\times [s,q])}\in \Guill_2((p,q);(p,s),(p,q-s))
\end{align}
\end{subequations}
for the products in an algebra over the operad $\Guill_2$. The indices $WE$ (resp. $SN$) stands for West-East (resp. South-North) in order to illustrate the first rectangle is on the left (resp. at the bottom) and the second rectangle on the right (resp. at the top). The exponents reflect the sizes of the rectangles. Whenever the sizes are clear from context (or arbitrary), we may choose to drop the exponents as in the following proposition.

\begin{proof}
The proof is obtained simply by induction on the arity of the operations, using the recursivity (third item) of Definition~\ref{def:guillotinepartition2D} to reduce by a guillotine cut a $n$-arity operation to a $2$-arity operation applied to the guillotine sub-partitions, until reaching identities.
\end{proof}

\begin{prop}[generalized associativity of the generators]\label{prop:guill2elemassoc}
All the relations satisfied by the elementary guillotine partitions $m_{WE}$ and $m_{SN}$ are generated by the following relations
\begin{subequations}
\label{eq:guill2:listassoc}
\begin{align}
m_{WE} \circ (\id, m_{WE} ) &= m_{WE}\circ (m_{WE},\id)
\label{eq:guill2:horizassoc}
\\
m_{SN} \circ (\id, m_{SN} ) &= m_{SN}\circ (m_{SN},\id)
\label{eq:guill2:vertassoc}
\\
m_{WE} \circ (m_{SN}, m_{SN} ) &= m_{SN}\circ (m_{WE},m_{WE})
\label{eq:guill2:interchangeassoc}
\end{align}
\end{subequations}
whenever the colours are compatible with the composition rules.
\end{prop}
\begin{proof}
The existence of such relations is essentially related to the recursive nature definition of guillotine partitions~\ref{def:guillotinepartition2D} and the non-uniqueness of the first guillotine cut in point \eqref{item:hyp:existsguillcut}. Any guillotine partition can be encoded with a binary tree with guillotine cuts on the vertices and rectangles on the edges: this representation is not unique and we must quotient these tree representations by an equivalence relation defined by fundamental "moves".

For a guillotine partition $\rho=(R_1,\ldots,R_n)$ of $R$, we introduce $\admisscut(\rho)$ as the set of guillotine cuts such that \eqref{item:hyp:existsguillcut} of \ref{def:guillotinepartition2D} is valid. For every cut $C\in\admisscut(\rho)$, we write $(\rho_1(\rho,C),\rho_2(\rho,C))$ the guillotine partitions of $R$ obtained from the cut (the two rectangles may be ordered canonically from left to right for a vertical cut and from bottom to top for a horizontal one). 

Any composition of a guillotine partition from two smaller guillotine partitions has to use a first cut $C\in\admisscut$. Any cut $D\in\admisscut(\rho))$ different from $C$ thus defines:
\begin{itemize}
\item either one cut $D \in \admisscut(\rho_1(\rho,C))$  if $D$ and $C$ have the same direction and $D$ is at the left or bottom of $C$,
\item either one cut $D \in \admisscut(\rho_2(\rho,C))$  if $D$ and $C$ have the same direction and $D$ is at the right or top of $C$,
\item or two cuts $D$, one in $\admisscut(\rho_1(\rho,C))$ and a second one in $\admisscut(\rho_2(\rho,C))$ if $C$ and $D$ have different directions.
\end{itemize}
We first treat the first case. For any $\rho$, for any pair $(C,D)$ of distinct elements of $\admisscut(\rho)$, we must have equivalence of the following representation:
\begin{align*}
\rho&= \begin{tikzpicture}[guillpart]
\fill[guillfill] (0,0) rectangle (2,1);
\draw[guillsep] (0,0) rectangle (2,1);
\draw[guillsep] (1,0)--(1,1);
\node at (0.5,0.5) {$1$};
\node at (1.5,0.5) {$2$};
%\bulkrect{(0,0)}{(1,1)}{$1$};
%\bulkrect{(1,0)}{(2,1)}{$2$};
\node at (1,0) [below] {$C$};
\end{tikzpicture} 
\circ \left(
	\begin{tikzpicture}[guillpart]
		\fill[guillfill] (0,0) rectangle (2,1);
		\draw[guillsep] (0,0) rectangle (2,1);
		\draw[guillsep] (1,0)--(1,1);
		\node at (0.5,0.5) {$1$};
		\node at (1.5,0.5) {$2$};
	%\bulkrect{(0,0)}{(1,1)}{$1$};
	%\bulkrect{(1,0)}{(2,1)}{$2$};
	\node at (1,0) [below] {$D$};
	\end{tikzpicture}
	\circ\Big( \rho_1(\rho_1(\rho,C),D), \rho_2(\rho_1(\rho,C),D)\Big)
	,\rho_2(C)
\right) 
\\
&= \begin{tikzpicture}[guillpart]
\fill[guillfill] (0,0) rectangle (2,1);
\draw[guillsep] (0,0) rectangle (2,1);
\draw[guillsep] (1,0)--(1,1);
\node at (0.5,0.5) {$1$};
\node at (1.5,0.5) {$2$};
%\bulkrect{(0,0)}{(1,1)}{$1$};
%\bulkrect{(1,0)}{(2,1)}{$2$};
\node at (1,0) [below] {$D$};
\end{tikzpicture} 
\circ \left(
\rho_1(D),
	\begin{tikzpicture}[guillpart]
	\fill[guillfill] (0,0) rectangle (2,1);
	\draw[guillsep] (0,0) rectangle (2,1);
	\draw[guillsep] (1,0)--(1,1);
	\node at (0.5,0.5) {$1$};
	\node at (1.5,0.5) {$2$};	
%	\bulkrect{(0,0)}{(1,1)}{$1$};
%	\bulkrect{(1,0)}{(2,1)}{$2$};
	\node at (1,0) [below] {$C$};
	\end{tikzpicture}
	\circ \Big( \rho_1(\rho_2(\rho,D),C), \rho_2(\rho_2(\rho,D),C)\Big)
\right) 
\end{align*}
It is then easy to check geometrically that \begin{align*}
\rho_1(\rho_1(\rho,C),D) &=\rho_1(D) 
\\
\rho_2(\rho_1(\rho,C),D) &=  \rho_1(\rho_2(\rho,D),C)
\\
\rho_2(C) &= \rho_2(\rho_2(\rho,D),C)
\end{align*}
Since $\rho$ is arbitrary, for any guillotine partitions $a,b,c$ with compatible sizes, it holds
\[
\begin{tikzpicture}[guillpart]
	\fill[guillfill] (0,0) rectangle (2,1);
	\draw[guillsep] (0,0) rectangle (2,1);
	\draw[guillsep] (1,0)--(1,1);
	\node at (0.5,0.5) {$1$};
	\node at (1.5,0.5) {$2$};
%\bulkrect{(0,0)}{(1,1)}{$1$};
%\bulkrect{(1,0)}{(2,1)}{$2$};
\node at (1,0) [below] {$C$};
\end{tikzpicture} 
\circ \left(
	\begin{tikzpicture}[guillpart]
		\fill[guillfill] (0,0) rectangle (2,1);
		\draw[guillsep] (0,0) rectangle (2,1);
		\draw[guillsep] (1,0)--(1,1);
		\node at (0.5,0.5) {$1$};
		\node at (1.5,0.5) {$2$};
%	\bulkrect{(0,0)}{(1,1)}{$1$};
%	\bulkrect{(1,0)}{(2,1)}{$2$};
	\node at (1,0) [below] {$D$};
	\end{tikzpicture}
	\circ\Big( a, b\Big)
	,c
\right) 
=
\begin{tikzpicture}[guillpart]
	\fill[guillfill] (0,0) rectangle (2,1);
	\draw[guillsep] (0,0) rectangle (2,1);
	\draw[guillsep] (1,0)--(1,1);
	\node at (0.5,0.5) {$1$};
	\node at (1.5,0.5) {$2$};
%\bulkrect{(0,0)}{(1,1)}{$1$};
%\bulkrect{(1,0)}{(2,1)}{$2$};
\node at (1,0) [below] {$D$};
\end{tikzpicture} 
\circ \left(
a,
	\begin{tikzpicture}[guillpart]
		\fill[guillfill] (0,0) rectangle (2,1);
		\draw[guillsep] (0,0) rectangle (2,1);
		\draw[guillsep] (1,0)--(1,1);
		\node at (0.5,0.5) {$1$};
		\node at (1.5,0.5) {$2$};
%	\bulkrect{(0,0)}{(1,1)}{$1$};
%	\bulkrect{(1,0)}{(2,1)}{$2$};
	\node at (1,0) [below] {$C$};
	\end{tikzpicture}
	\circ \Big(b,c\Big)
\right) 
\]
which is exactly equation~\eqref{eq:guill2:horizassoc}. The same computations hold for the other cases.

One also checks that there are no other type of redundancies in the binary tree representations of a guillotine partitions and that such "moves" on all vertices of a tree can always lead a canonical binary tree with vertical cuts ordered from left to right first followed by horizontal cuts ordered from bottom to top.
\end{proof}

\begin{rema}[double-monoid structure with interchange] Equations \eqref{eq:guill2:listassoc} already appear in the literature \cite{Kock,BagherzadehBremner} in the \emph{absence} of colours on the operad and they lead to the interesting phenomenon of hidden commutativity. However, hidden commutativity relations do \emph{not} appear in our context: this is trivial to see in the relevant construction of the next Section~\ref{sec:canonicalguill2:markov}). Indeed, the proofs of hidden commutativity in the previous reference requires associative moves that are incompatible with the colour compatibility requirements of the operad $\Guill_2$.
\end{rema}

\begin{rema}[similarities with the $E_2$ operad] The guillotine partition operad has a geometrical origin similar to the little $2$-cubes operad corresponding to the extraction of sub-rectangles of a larger rectangles. However, the scaling $s_R$ that appears in the little $2$-cubes operad is replaced here by a colour composability criterium in $\PatternShapes(\patterntype{r})$ (only a part of the translations remains through the equivalence classes), which however changes dramatically what can be done with this operad but is however unavoidable in order to be able to formulate the notion of eigenvalue in Section~\ref{sec:boundaryalgebra}. Everything is formulated here in a discrete setting but one can see that the palette of colours is related to the Riemannian structure that is needed to speak about lengths, surfaces and volumes.
\end{rema}

\subsection{Back to dimension 1: some definitions and remarks}
\subsubsection{Definition}
It is also interesting to consider the one-dimensional equivalent $\Guill_1$ of $\Guill_2$ and reinterpret already known results. This operad $\Guill_1$ may look trivial at first sight or close to existing ones but reveals to be interesting from the points of view of eigenvalues (see Section~\ref{sec:boundaryalgebra}) and commutativity. It also appears many times as relevant sub-operads of $\Guill_2$.

The rectangles are replaced by segments $[u,v]$ with $(u,v)\in\setP$ and $u\leq v$, rectangle shapes by length $p\in\setL$.

\begin{defi}[one-dimensional "guillotine" partitions]
\label{def:oneguillotpartition}
A 1D guillotine partitions of a segment $[u,v]$ with $(u,v)\in\setP^2$  with $u\leq v$ is a sequence of segments $\ca{S}=([u_i,v_i])_{1\leq i\leq n}$ of $[0,p]$ such that the interiors $(u_i,v_i)$ are disjoint and $\bigcup_{1\leq i\leq n} [u_i,v_i]= [u,v]$.
\end{defi}
Compared to the two-dimensional case, the shapes are two simple to require the notion of guillotine cut, which would correspond here to the choice of a single point $w\in [u,v]$. In particular, if all the degenerate segments (i.e. $u_i<v_i$) are distinct, there exists a unique ordering permutation $\sigma\in\Sym{n}$ such that $u_{\sigma(1)}=u$, $v_{\sigma(i)}=u_{\sigma(i+1)}$ for $1\leq i\leq n-1$ and $v_\sigma(n)=v$. 

One may also define composition rules in the same way as in the two-dimensional case, by first considering equivalence classes under translations and then nestings. The operad we obtain corresponds to the $\mathrm{Ass}$ operad with additional colours related to the size of the segments.

\begin{defi}[$\Guill_1$ operad]
\label{def:Guill1operad}
The $\Guill_1$ operad is the operad with colours in $\setL$ whose $n$-ary operations of $\Guill_1(p;p_1,\ldots,p_n)$ are the equivalence classes under translations of partitions $([u_i,v_i])_{1\leq i\leq n}$ with $v_i-u_i=p_i$ of the segment $[0,p]$. The action of the symmetric group is the permutation of the segments $[u_i,v_i]$.
\end{defi}

In particular, $\Guill_1$ is generated by the permutations, the identities and $2$-ary partitions, which satisfy an associative property properly decorated by the length of the segments.

\subsubsection{Relation with \texorpdfstring{$\Ass$}{Ass}}
An algebra over this operad $\Guill_1$ is given by a collection of spaces $(\ca{A}_p)_{p\in\setL^*}$ and products $m_\rho:\ca{A}_{p_1}\otimes \ldots \ca{A}_{p_n} \to \ca{A}_p$ for any partition $\rho\in \Guill_1(p;p_1,\ldots,p_n)$ with $p=\sum_{i=1}^n p_i$. One checks easily that this operad is just a length colourization of $\Ass$ with the following construction. Let $\ha{\ca{A}}=\bigoplus_{p\in\setL^*} \ca{A}_p$. One may define products $\ha{m}_n:\ca{A}^{\otimes n}\to \ca{A}$ such that, for any $(a_1,\ldots,a_n)\in \ca{A}_{p_1}\otimesdots \ca{A}_{p_n}$, it holds:
\begin{equation}
\ha{m}_n(a_1,\ldots,a_n)=m_{([0,p_1],[p_1,p_1+p_2],\ldots,[p-p_n,p])}(a_1,\ldots,a_n) \in \ca{A}_p
\end{equation}
where $p=p_1+\ldots+p_n$. One checks then easily that $\ha{\ca{A}}$ is an algebra over the operad $\mathrm{Ass}$. 

Similarly, given an algebra $\ca{A}$ over the operad $\mathrm{Ass}$, we may always define $\ca{A}_p=\ca{A}$ and obtain trivially an algebra over $\Guill_1$. Such tricks are however impossible in larger dimensions.

We will see in Section~\ref{sec:onedim:extendedlengthass} that a second deeper link between $\Guill_1$-algebras and $\Ass$-algebra exists as soon as a generalization of left or right modules is introduced.

\subsubsection{Null length, commutativity and unit elements}

A particular sub-operad of $\Guill_1$ is given by the sets $\Guill_1(0;0,\ldots,0)$ ($n$ zeroes) with $n\in\setN$. From the geometric point of view, for any integer $n$, this set is a singleton and it corresponds to the guillotine partition of the degenerate segment $\{u\}$ into $n$ degenerate segments $\{u\}$. This partition is invariant under the action of $\Sym{n}$ and thus it is isomorphic to the commutative operad $\Com$. We see that the presence of colours and their geometric interpretation provides a nice framework to illustrate the relation between $\Ass$ and $\Com$.

Moreover, one may be interested in unit elements in a $\Guill_1$-algebra $\ca{A}$. Given a $2$-ary product $m\in \Guill_2(p_1+p_2;p_1,p_2)$, a unit element $e$ should at least satisfy $m(e,a)=a$ and thus $p_2$ and $p_1+p_2$ must be equal; hence the first size $p_1$ has to be zero. Unit elements necessarily belongs to $\ca{A}_0$. We have seen above that $\ca{A}_0$ is necessarily commutative and thus unit elements necessarily commute. A deeper discussion about unit elements is presented in Section~\ref{sec:eckmanhilton}.

\subsection{Relevant sub-operads of \texorpdfstring{$\Guill_2$}{Guill2}.} 
\paragraph*{Restriction on sizes}

Given an element $(\kappa_1,\kappa_2)\in \setL^2$, we may consider only the guillotine partitions with rectangles of sizes $(p,q)$ such that $p\geq \kappa_1$ and $q\geq \kappa_2)$ and one observes that it defines a sub-operad $\Guill_2^{(\leq \kappa_1,\leq \kappa_2)}$ of $\Guill_2$. The large inequalities on sizes can as well be replaced by strict inequalities. In particular, degenerate rectangles can be removed from the definition of guillotine partitions by considering $\Guill_2^{(> 0,>0)}$.

\paragraph*{Sub-\texorpdfstring{$\Guill_1$}{Guill1}-operads}

One checks that, for any fixed $p\in\setN_1$, the sets \begin{equation}\Guill_1^{(p),\mathrm{vert}}(q;q_1,\ldots,q_n) =  \Guill_2((p,q);(p,q_1),\ldots,(p,q_n)),
\label{eq:suboperadGuill1Guill2}
\end{equation}
which contain only partitions produced by horizontal guillotine cuts,
form an operad isomorphic to $\Guill_1$. This corresponds to the fact that the elementary products $m_{SN}$ satisfies the coloured associativity~\eqref{eq:guill2:vertassoc}, which is the building block of the operad $\Guill_1$. The same holds for $\Guill_1^{(q),\mathrm{horiz}}$ obtained by vertical guillotine cuts, with the elementary products $m_{WE}$. 

Thus, the operad $\Guill_2$ can be seen as collections of vertical and horizontal $\Guill_1$ operads with additional interchange relation~\eqref{eq:guill2:interchangeassoc}, in the spirit of \cite{Kock,BagherzadehBremner}. This remark will become fundamental in Section~\ref{sec:boundaryguill2} for the generalization of left and right modules to dimension two.

\paragraph*{Commutative sub-operads}

From the discussion above about the relation between $\Com$ and the null size partitions in $\Guill_1$, we may extract a $\Com$ sub-operad of $\Guill_2$ for every segment. 

\begin{prop}
	Let $S=\{u_1\}\times [u_2,v_2]$ (resp. $[u_1,v_1]\times \{u_2\}$) be a vertical (resp. horizontal) segment of $\setP^2$ with size $(0,p_2)$ (resp. $(p_1,0)$). The sub-operad $\Guill_2( (0,p_2) ; (0,p_2),\ldots,(0,p_2) )$ ($n$ times), $n\in\setN$, which inherits the horizontal product,
    is isomorphic to $\Com$ and is noted $\Com_2(S)$, as well as the sub-operad $\Guill_2( (p_1,0) ; (p_1,0),\ldots,(p_1,0) )$ ($n$ times) with the vertical product.
\end{prop}

Gluing two consecutive segments $S_1$ and $S_2$ in a given direction introduces the use of the second transverse product and provides morphisms from two algebras over $\Com_2(S_1)$ and $\Com_2(S_2)$ respectively to an algebra over $\Com_2(S_1\cup S_2)$. Such an example is omnipresent in statistical mechanics and corresponds to the commutative algebra of observables over the edge variables as it will be seen below.

		\subsection{An algebra over guillotine partitions for two-dimensional Markov processes.}\label{sec:canonicalguill2:markov}
	
We now introduce an algebra over $\Guill_2$ that is fundamental for the probabilistic formalism of two-dimensional Markov processes since it encodes algebraically the gluing property of domains at the heart of the Markov property.

\subsubsection{A canonical finite-dimensional \texorpdfstring{$\Guill_2$}{Guill2}-algebra for discrete state spaces and discrete space \texorpdfstring{$\setP=\setZ$}{P=Z}.}

\begin{theo}[tensor algebra of matrices as $\Guill_2$-algebra] 
\label{theo:canonicalexampleGuill}
Let $V_1$ and $V_2$ be two finite-dimensional vector spaces with canonical bases $(e^{(1)}_k)_{k}$ and $(e^{(2)}_k)$. For any $(p,q)\in \setN_1\times\setN_1$, let $T_{p,q}(V_1,V_2)$ be the vector space defined by 
\begin{equation}
\label{eq:defTpqV1V2}
T_{p,q}(V_1,V_2) = 
\begin{cases}
\End(V_1)^{\otimes p}\otimes \End(V_2)^{\otimes q} & \text{for $p>0$ and $q>0$} \\
\Diag_{e^{(1)}}(V_1)^{\otimes p} & \text{for $p>0$ and $q=0$} \\
\Diag_{e^{(2)}}(V_2)^{\otimes q}  & \text{for $p=0$ and $q>0$} \\
 \setK & \text{for $p=0$ and $q=0$} \\
\end{cases}
\end{equation}
where $\Diag_{e^{(i)}}(V_i)$ is the set of operators on $V_i$ that are diagonal in the basis $(e^{(i)}_k)$. Let the elementary products $m_{WE}$ and $m_{SN}$ be defined by:
\begin{align*}
m_{WE}^{p_1,p_2|q}: T_{p_1,q}(V_1,V_2)\otimes T_{p_2,q}(V_1,V_2) 
& \to T_{p_1+p_2,q}(V_1,V_2)
\\
\left(\bigotimes_{i=1}^{p_1} A_i  \otimes \bigotimes_{k=1}^q B_k\right)
\otimes 
\left(\bigotimes_{i=1}^{p_2} A'_i  \otimes \bigotimes_{k=1}^q B'_k\right)
&\mapsto 
\left(\bigotimes_{i=1}^{p_1} A_i\otimes \bigotimes_{i=1}^{p_2} A'_i\right)\otimes \bigotimes_{k=1}^q (B_k B'_k)
\\
\\
m_{SN}^{p|q_1,q_2}: T_{p,q_1}(V_1,V_2)\otimes T_{p,q_2}(V_1,V_2) 
& \to T_{p,q_1+q_2}(V_1,V_2)
\\
\left(\bigotimes_{i=1}^{p} A_i  \otimes \bigotimes_{k=1}^{q_1} B_k\right)
\otimes 
\left(\bigotimes_{i=1}^{p} A'_i  \otimes \bigotimes_{k=1}^{q_2} B'_k\right)
&\mapsto 
\bigotimes_{i=1}^{p} (A_i A'_i)\otimes \left(\bigotimes_{k=1}^{q_1} B_k\otimes \bigotimes_{k=1}^{q_2} B'_k\right)
\end{align*} 
where, by convention the tensor products are just absent when a size $p_i$ or $q_i$ is equal to $0$. Then, the sequence $(T_{p,q}(V_1,V_2))_{(p,q)\in \PatternShapes(\patterntype{r})}$ is an algebra over the operad $\Guill_2$ generated by the elementary products $m_{WE}^{p_1,p_2|q}$ and $m_{SN}^{p|q_1,q_2}$.
\end{theo}

\begin{proof}
We first check that the invariance of the products under permutations of $\Sigma{2}$ when degenerate rectangles appear with multiplicity in ensured by the commutative structure of the sets of diagonal operators $\Diag_{e^{(i)}}(V_i)$ and of the base field $\setK$.	
	
From propositions~\ref{prop:guill2generators} and  \ref{prop:guill2elemassoc}, it is sufficient to show that the elementary products satisfy the associativity relations~\eqref{eq:guill2:listassoc}: this is an elementary exercise from the definitions above. 

However, we propose a more algebraic perspective by observing that, for any vector space $V$, $\End(V)^{\otimes p}$ can be endowed with two associative products. The first one is a coloured identity and the other is not coloured and relies on the usual product on $\End(V)$:
\begin{subequations}
\label{eq:productandconcat:assoc}
\begin{align}
\bullet:\End(V)^{\otimes p_1} \otimes \End(V)^{\otimes p_2} 
&\to \End(V)^{\otimes p_1+ p_2}
\\
\left(\bigotimes_{i=1}^{p_1} A_i\right) \otimes \left(\bigotimes_{i=1}^{p_2} A'_i \right) & \mapsto \bigotimes_{i=1}^{p_1} A_i \otimes \bigotimes_{i=1}^{p_2} A'_i
\\
m: \End(V)^{\otimes p} \otimes \End(V)^{\otimes p} 
&\to \End(V)^{\otimes p}
\\
\left(\bigotimes_{i=1}^{p} A_i\right)\otimes \left(\bigotimes_{i=1}^{p} A'_i\right) &\mapsto \bigotimes_{i=1}^{p} (A_i A'_i)
\end{align} 
Both products are associative, the first in $\Guill_1$ and the second in $\Ass$. These two products moreover satisfy an interchange relation
\begin{equation}\label{eq:interchangeforcanostruct}
m(A\bullet A',B\bullet B') = m(AB)\bullet m(A'B')
\end{equation}
\end{subequations}
as soon as the colours are compatible for the two products. The product  $m_{WE}$ (resp. $m_{SN}$) is thus inherited from the the product $\bullet$ (resp. $m$) on $(\End(V_1)^{\otimes p})_{p}$ and the product $m$ (resp. $\bullet$) on $(\End(V_2)^{\otimes p})_p$. This is then an elementary exercise that the combination of associativities of $\bullet$ and $m$ together with \eqref{eq:interchangeforcanostruct} provides a $\Guill_2$ structure and the associativities \eqref{eq:guill2:listassoc}.
\end{proof}

\subsubsection{Operadic view on partition functions of Markov processes.}

A partition function as introduced in the probability law \eqref{eq:MarkovLawDimTwo} on a rectangular domain can be embedded in a suitable space $T_{p,q}(V_1,V_2)$ via a matrix representation. This description can indeed be lifted to the operadic level as illustrated by the following structure theorem.

In the transcription from probabilistic notations to algebraic ones, we introduce the following notations:
\begin{subequations}
\label{eq:boundary:matrixcoeff:corresp}
\begin{itemize}
	\item for any non-degenerate rectangle of size $(p,q)$, we identify bijectively any function $F : S_1^p\times S_1^p \times S_2^q\times S_2^q \to \setK$ on the boundary edge spaces to an element $\ha{F} \in T_{p,q}(V(S_1),V(S_2))$ through
	\begin{equation}
	\ha{F} \coloneqq \sum_{(x,y,w,z)} F(x,y,w,z) E^{(p)}_{x,y}\otimes E^{(q)}_{w,z}
	\end{equation}
	with $E^{(p)}_{x,y} = E_{x_1,y_1}\otimesdots E_{x_p,y_p}$ and $E_{i,j}$ is the elementary matrix with $1$ at position $(i,j)$ and $0$ elsewhere.
	
	\item for a degenerate rectangle of size $(p,0)$ (resp. $(0,q)$), we identify bijectively any function $F:S_1^p \to \setK$ (resp. $S_2^q\to\setK$) to the diagonal matrix $\ha{F}\in T_{p,0}(V(S_1),V(S_2))$ through
	\begin{equation}
	\ha{F} \coloneqq \sum_{x\in S_1^p} F(x) E^{(q)}_{x,x}	\qquad \left(\text{resp.\,} \sum_{x\in S_2^q} F(x) E^{q}_{x,x}  \right)
	\end{equation}
\end{itemize}
\end{subequations}

\begin{theo}[partition functions as  $\Guill_2$-algebra]
	\label{theo:partitionfuncguillotop}
Let $S_1$ and $S_2$ be two finite sets and let $V(S_1)=\setR^{S_1}$ and $V(S_2)=\setR^{S_2}$ (with their canonical bases) be their associated finite-dimensional vector spaces. Let $R=[P_1,P_2]\times [Q_1,Q_2]$ be a rectangle in $\setZ^2$. Let $(X_e)_{e\in \Edges{R}}$ be a $(S_1,S_2)$-valued Markov process on $R$ with face weights $(\MarkovWeight{W}_f)_{f\in R}$.

For any sub-rectangle $R'\subset R$ (with shape $(p,q)$) and any guillotine partition $(R'_1,\ldots,R'_n)$ of $R'$, it holds (under the previous identification):
\begin{equation}\label{eq:compopartitionfunctions}
\ha{Z}_{R'}(\MarkovWeight{W}_\bullet) = m_{R'_1,\ldots,R'_n}(\ha{Z}_{R'_1}(\MarkovWeight{W}_\bullet),\ldots, \ha{Z}_{R'_n}(\MarkovWeight{W}) )
\end{equation}
in the canonical $\Guill_2$-algebra $(T_{p,q}(V(S_1),V(S_2)))_{(p,q)\in \PatternShapes(\patterntype{r})}$, with $\ha{Z}_{f}(\MarkovWeight{W}_\bullet) = \ha{\MarkovWeight{W}}_f$ for $R=f$ an elementary face. 

Moreover, for any collection of functions $(h_e)_{e\in \Edges{R}}$ from $S_1$ or $S_2$ to $\setR$, it holds:
	\begin{equation}\label{eq:observablesexpection}
		\Espc{\prod_{e\in E}h_e(X_e)}{ (X_e)_{e\in\partial R} }
		= \frac{1}{
			Z_R(\MarkovWeight{W}_\bullet ; (X_e)_{e\in\partial R} )
		}
		\begin{tikzpicture}[guillpart]
			\begin{scope}[yscale=2.3,xscale=2.6]
				\fill[guillfill] (0,0) rectangle (4,3);
				\draw[guillsep] (0,0) rectangle (4,3);
				\draw[guillsep] (0,0) rectangle node {$\ha{\MarkovWeight{W}}_{\diamond}$} (1,1);
				\draw[guillsep] (1,0) rectangle node {$\ha{\MarkovWeight{W}}_{\diamond}$} (2,1);
				\draw[guillsep] (2,0) rectangle node {$\ldots$} (3,1);
				\draw[guillsep] (3,0) rectangle node {$\ha{\MarkovWeight{W}}_{\diamond}$} (4,1);
				\draw[guillsep] (0,1) rectangle node {$\vdots$} (1,2);
				\draw[guillsep] (1,1) rectangle node {$\vdots$} (2,2);
				\draw[guillsep] (2,1) rectangle node {$\vdots$} (3,2);
				\draw[guillsep] (3,1) rectangle node {$\vdots$} (4,2);
				\draw[guillsep] (0,2) rectangle node {$\ha{\MarkovWeight{W}}_{\diamond}$} (1,3);
				\draw[guillsep] (1,2) rectangle node {$\ha{\MarkovWeight{W}}_{\diamond}$} (2,3);
				\draw[guillsep] (2,2) rectangle node {$\ldots$} (3,3);
				\draw[guillsep] (3,2) rectangle node {$\ha{\MarkovWeight{W}}_{\diamond}$} (4,3);
				\foreach \x in {0.5,1.5,...,3.5}
				{
					\foreach \y in {0,1,...,3}
					\node at (\x,\y) {$\ha{h}_{\diamond}$};
				}
				\foreach \x in {0,1,...,4} 
				{
					\foreach \y in {0.5,1.5,...,2.5}
					\node at (\x,\y) {$\ha{h}_{\diamond}$};
				}
			\end{scope}
		\end{tikzpicture}_{(X_e)_{e\in\partial R}}
	\end{equation}
	where each symbol $\diamond$ has to be replaced by the corresponding face or edge $D$ on an edge $e$. The last subscript in the r.h.s. corresponds to the extraction of the matrix element associated to the boundary condition in the correspondence~\eqref{eq:boundary:matrixcoeff:corresp}.
\end{theo}

This theorem reflects the well-known gluing properties of partitions functions in statistical mechanics using the Markov property. The proof below reinterprets these properties in an operadic language. For example, the associativity corresponds to Fubini's theorem. The theorem may look useless at first sight in practice: this is however a first step towards the much more interesting results of Section~\ref{sec:boundaryalgebra}, where specific algebra-motivated boundary conditions will enter the framework. It provides an algebraic fully-two-dimensional counterpart to the Markov property, in the same way as linear algebra is related to one-dimensional Markov processes.

\begin{proof}
The proof relies on a recursion on the number of elements of the guillotine partitions and on suitable identifications of partial sums. For $n=1$, the identification between $Z_{R'}(\MarkovWeight{W}_\bullet;\cdot)$ and $\ha{Z}_{R'}(\MarkovWeight{W}_\bullet)$ is a definition and there are no computations involved.

For $n=2$ and a vertical guillotine cut dividing a rectangle $R'$ of size $(p,q)$ into two rectangles $R'_1$ of size $(p_1,q)$ and $R'_2$ of size $(p_2,q)$ with $p_1+p_2=p$, we have from the definition of partition functions: 
\begin{align*}
&\ha{Z}_{R'}(\MarkovWeight{W}_\bullet)_{x_1, y_1; \ldots; x_p, y_p; w_1, z_1; \ldots; w_q, z_q}=
\sum_{(x^{int}_e)_{e\in\Int(R')}} \prod_{f\in R'} \MarkovWeight{W}_f(x_{\partial f}) 
\\
&= \sum_{(x^{int}_e)_{e\in \Int(R')}} \left(\prod_{f\in R'_1} \MarkovWeight{W}_f(x_{\partial f})  \right) \left(\prod_{f\in R'_2} \MarkovWeight{W}_f(x_{\partial f} )  \right)
\end{align*}
and we use here and below the following abuses of notations: $x_E$ is the sequence $(x_e)_{e\in E}$, $x_{\partial f}$ is the $4$-uplet of variables on the four sides of $f$ in the order South, North, West and East; each $x_e$ is either one of the $x^{int}_e$ if it lies strictly inside the rectangle or one of the $x_i$, $y_i$, $w_i$ or $z_i$ if it lies on the boundary.

Now, the set of edges $\Edges{R'}\setminus\partial R'$ is the disjoint union of the three sets  $\Edges{R'_1}\setminus\partial R'_1$, $\Edges{R'_2}\setminus\partial R'_2$ and the edges $(x^{cut}_k)_{1\leq k \leq q}$ on the cut between $R'_1$ and $R'_2$ (enumerated from bottom to top). We thus have, from the definition of the partition functions:
\begin{align*}
&\ha{Z}_{R'}(\MarkovWeight{W}_\bullet)_{x_1, y_1; \ldots; x_p, y_p; w_1, z_1; \ldots; w_q, z_q} 
\\
&= \sum_{(x^{cut}_k)} \left( \sum_{x^{int}_{\Edges{R'_1}\setminus \partial R'_1}} \prod_{f\in R'_1} \MarkovWeight{W}_f(x_{\partial f})\right)\left( \sum_{x^{int}_{ \Edges{R'_2}\setminus \partial R'_2}} \prod_{f\in R'_2} \MarkovWeight{W}_f(x_{\partial f})\right)
\\
&= \sum_{(x^{cut}_k)}
\ha{Z}_{R'_1}(\MarkovWeight{W}_\bullet)_{x_1, y_1; \ldots; x_{p_1}, y_{p_1}; w_1, x^{cut}_1; \ldots; w_q, x^{cut}_q} 
\ha{Z}_{R'_2}(\MarkovWeight{W}_\bullet)_{x_{p_1+1}, y_{p_1+1}; \ldots; x_p, y_p; x^{cut}_1, z_1; \ldots; x^{cut}_q, z_q} 
\end{align*}
and we recognize the definition of the composition law $m_{WE}$ on $(T_{p,q})$ in its matrix version. The same computation holds for a vertical cut.

We now assume that, for a fixed integer $n\geq 3$, \eqref{eq:compopartitionfunctions} holds for any guillotine partitions of size $k\leq n$. For any guillotine partition $\rho=(R_1,\ldots,R_{n+1})$ of size $n+1$ of a rectangle $R'$, we first consider a cut $C \in \admisscut{\rho}$, which divides $R'$ into $R'_1$ and $R'_2$ with inherited guillotine partitions $\rho_1=(R_i)_{i\in I_1}$ and $\rho_2=(R_{i})_{i\in I_2}$ where $(I_1,I_2)$ is a non-trivial partition of $\{1,\ldots,n+1\}$. The case $n=2$ treated before leads to 
\[
	\ha{Z}_{R'}(\MarkovWeight{W}_\bullet) = m_{(R'_1,R'_2)}( \ha{Z}_{R'_1}(\MarkovWeight{W}_\bullet), \ha{Z}_{R'_2}(\MarkovWeight{W}_\bullet) )
\]
The strong recursion hypothesis provides now
\[
\ha{Z}_{R'_j}(\MarkovWeight{W}_\bullet) = m_{\rho_j}\left( (\ha{Z}_{R_i}(\MarkovWeight{W}_\bullet))_{i\in I_j} \right)
\]
and thus, from the operadic composition of guillotine partitions, we obtain
\begin{align*}
\ha{Z}_{R'}(\MarkovWeight{W}_\bullet) &= 
m_{(R'_1,R'_2)}\left(m_{\rho_1}\left( (\ha{Z}_{R_i}(\MarkovWeight{W}_\bullet))_{i\in I_1} \right),
m_{\rho_2}\left( (\ha{Z}_{R_i}(\MarkovWeight{W}_\bullet))_{i\in I_2} \right) \right)
\\
&= 
m_\rho( \ha{Z}_{R_1}(\MarkovWeight{W}_\bullet),\ldots \ha{Z}_{R_{n+1}}(\MarkovWeight{W}_\bullet) )
\end{align*}
and the proof is complete.
\end{proof}

\subsubsection{A similar view on the one-dimensional case.} We rephrase here the well-known one-dimensional situation in the framework of the operad $\Guill_1$. Let $S$ be a finite state space and $V(S)=\setR^S$ the associated vector space. In this case, the relevant algebra over $\Guill_1$ is simply $\ca{A}_p = \End(V(S))$ for $p>0$ and $\ca{A}_0=\Diag(V(S))$ (this is frequent case, see Theorem~\ref{theo:1D:removingcoloursonboundaries} and Section~\ref{sec:eagerandlazyROPEs}). Contrary to the two-dimensional case, the spaces $\ca{A}_p$ do not depend on the colour and thus the correspondence with $\mathrm{Ass}$ is complete. The products $m_n$ are simply the usual products of $n$ matrices from left to right. 

Let $(X_k)_{0\leq k \leq p}$ be a $S$-valued Markov process on $\{0,1,\ldots,p\}$ with weights $(\MarkovWeight{A}_k)_{0\leq k<p}$. Each weight $\MarkovWeight{A}_k$ defines a matrix $\MarkovWeight{A}_k\in\ca{A}_1$ defined entry by entry by $\MarkovWeight{A}_k(i,j)$. The partition function on $\{k_1,k_1+1,\ldots,k_2\}$ with boundary conditions $(i,j)$ is then given by:
\begin{equation}
Z_{k_1,k_2}(\MarkovWeight{A}_\bullet;i,j) = (\MarkovWeight{A}_{k_1}\ldots \MarkovWeight{A}_{k_2-1})(i,j)
\end{equation}
It is then a simple exercise to check that for any (discrete guillotine) partition $([u_i,v_i])_{1\leq i\leq n}$ of $[0,p]$ with ordering permutation $\sigma$ the following identity holds:
\begin{equation}
Z_{0,p}(\MarkovWeight{A}_\bullet;i,j) =\left( Z_{u_{\sigma(1)},v_{\sigma(1)}}(\MarkovWeight{A}_\bullet;\cdot,\cdot)\ldots Z_{u_{\sigma(1)},v_{\sigma(1)}}(\MarkovWeight{A}_\bullet,\cdot,\cdot)\right)(i,j)
\end{equation}
where the products hold in $\End(V(S))$. This is the one-dimensional counterpart of Theorem~\ref{theo:partitionfuncguillotop} with traditional linear algebra instead of operads.

\begin{table}
	\begin{center}
		\begin{tabular}{|p{4.5cm}|p{4.5cm}|p{4.5cm}|}
			\hline
			\textbf{Probability} & \textbf{Discrete} & \textbf{Continuous} (Section~\ref{sec:generalizations}) \\
			\hline
			State space for a segment of length $p$ & 
			Finite sets $S^p$, $p\in\setN_1$ & 
			Segment state space $(S_p,\ca{S}_p,\nu_p)_{p\in\setL^*}$ \\
			\hline
			Boundary weight on state space $S$ &
			Vector in $V(S)=\setR^{|S|}$ & 
			(linear space) Measurable function in $\Meas(S)$ \\
			\hline
			Observable on state space $S$ &
			Diagonal matrix on $V(S)$ & 
			(algebra) Measurable function in $\Meas(S)$ \\
			\hline
			Joint functions on a product $S\times T$ &
			Tensor product $V(S)\otimes V(T)$ & 
			Measurable function in $\Meas(S\times T)$ \\
			\hline
			Expectation $\sum_{x\in S} f(x)g(x) $ &
			Canonical scalar product on $V(S)$ & 
			Integration \textbf{(I)} $\int_S f(x)g(x)d\nu(x)$\\
			\hline
			Weight $w: S\times S \to \setR $ & 
			$\End( V(S) )$ & 
			Weight in $\Meas(S\times S)$ \\
			\hline
			Marginalization $\sum_{y\in S} f(x,y)g(y,z)$ & 
			Composition in $\End(V(S))$ & 
			Composition \textbf{(I)} $\int_{S} f(x,y) g(y,z) d\nu(y)$ \\
			\hline
			Reorganizing sums &
			Associativity of $\circ$& 
			Fubini's theorem \\
			\hline 
			Boundary marginalization $\sum_{y\in S} f(x,y)g(y)$ & 
			Action of $\End(V(S))$ on $V(S)$ & 
			Composition \textbf{(I)} $\int_{S} f(x,y) g(y) d\nu(y)$ \\
			\hline
			Reorganizing sums &
			Axiom of an action & 
			Fubini's theorem \\
			\hline
		\end{tabular}
	\end{center}
	\caption{Correspondence between the situation of finite state space and discrete space (linear algebra with tensor products) and general state space with discrete state space or continuous state space (measure theory with Fubini's theorem, as in Section~\ref{sec:generalizations}). The symbol \textbf{(I)} requires additional analytical restrictions for the integrals to be defined: $\sigma$-finite measures, positive $\ov{\setR}_+$-valued functions or essentially bounded functions or $L^2$ functions depending on the situation.}\label{tab:equiv:vectorspace:segmentstatespace}
\end{table}

\subsection{Surface powers in the discrete setting and regular partitions.}
\label{par:surfacepowers}  

The addition of the colours in $\PatternShapes(\patterntype{r})$ for the generators $m_{WE}^{p_1,p_2|q}$ and $m_{SN}^{p|q_1,q_2}$ gives access to the notion of \emph{surface power} of an element. 

In the commutative case, there exist unique maps $\setK \to \setK_{p,q}$, $u\mapsto u^{[p,q]}$ such that $\setK\to \ca{A}_{1,1}$ is the identity and \begin{equation}\label{eq:surfacepowers}
m_{(R_1,\ldots,R_n)}(u^{[p_1,q_1]},\ldots,u^{[p_n,q_n]}) = u^{[p,q]}
\end{equation}
where $(R_1,\ldots,R_n)$ is a guillotine partition of $[0,p]\times[0,q]$ and the $(p_i,q_i)$s are the shapes of the $R_i$. These maps are precisely given by $u^{[p,q]}=u^{pq}$. The quantity $pq$ is the \emph{area} of the rectangle and justifies the name of \emph{surface power}. This is the first situation where the metric notion of area of a rectangle appears and it will reappear later for free energy densities in statistical mechanics and generalization of eigenvalues in Section~\ref{sec:invariantboundaryelmts}.

For a more general $\Guill_2$-operad $(\ca{A}_{p,q})_{(p,q)\in\setN^2}$, we may define similar \emph{surface powers} by assigning to any element $u\in\ca{A}_{p,q}$, a collection of elements $u^{[n,m]}\in\ca{A}_{np,mq}$ such that $u^{[1,1]}=u$ and such that it satisfies the same equation~\eqref{eq:surfacepowers}.

\begin{prop}[partition function of homogeneous process]\label{prop:partitionfunc:homogeneous}
	Under the same notations as for Theorem~\ref{theo:partitionfuncguillotop}, if all the face weights $\MarkovWeight{W}_{f}$  are all equal to the same face weight $\MarkovWeight{W}$, identified with $\ha{\MarkovWeight{W}}\in T_{1,1}(V(S_1),V(S_2))$, then the partition function $\ha{Z}_{R}(\ha{W})$ on a rectangle $R$ with shape $(p,q)$ is equal to the following surface power
	\[
	\ha{Z}_{R}(\ha{\MarkovWeight{W}}) = \ha{\MarkovWeight{W}}^{[p,q]}
	\]
\end{prop}
\begin{proof}
	It follows from a double recursion on the rectangle sizes $p$ and $q$: the recursion step is obtained using~\eqref{eq:compopartitionfunctions} and the definition of surface powers.
\end{proof}

\subsection{Dihedral group and generalization of opposite algebras.} In dimension one, it is a well-known fact that any algebra structure on a space also defines an opposite algebra structure obtained by defining the opposite product $m^{\mathrm{op}}(a,b)=m(b,a)$, which corresponds to the second element in the cardinal two set $\Ass_2$. This construction can be directly extended to the coloured version $\Guill_1$: an algebra $(\ca{A}_p)_{p\in\setL^*}$ over this operad admits an opposite product with the same colour compatibility conditions as the initial one. It corresponds geometrically to the isometry that flips the real line $x\mapsto -x$.

In dimension two, the situation is a bit more involved and we may consider transformations under which equations~\eqref{eq:guill2:listassoc} behave well. The first one is horizontal reversal corresponding to considering $m_{WE}^{\mathrm{op}}$ instead of $m_{WE}$ while keeping the same vertical $m_{SN}$ products: in this case, the two sides of \eqref{eq:guill2:horizassoc} are permuted, \eqref{eq:guill2:vertassoc} and \eqref{eq:guill2:vertassoc} are unchanged. The second one is vertical reversal, under which $m_{SN}^{\mathrm{op}}$ is changed to $m_{SN}^{\mathrm{op}}$ while $m_{WE}$ is unchanged. For these two operations, the colour compatibility conditions remain unchanged.

There is a new symmetry that can be introduced corresponding to diagonal reversal, under which a rectangle of shape $(p,q)$ is mapped to the rectangle of of shape $(q,p)$. Under this transformation, $m_{WE}$ becomes $m_{SN}$ and $m_{SN}$ becomes $m_{WE}$: it switches equations~\eqref{eq:guill2:horizassoc} and \eqref{eq:guill2:vertassoc} and switches the two sides of~\eqref{eq:guill2:interchangeassoc}. However, the colour compatibility conditions change and, under this diagonal reversal, an algebra $(\ca{A}_{p,q})$ over this coloured operad is an algebra $(\ca{A}'_{p,q})$ over the new coloured operad up to the identification $\ca{A}'_{p,q}=\ca{A}_{q,p}$.

These new symmetries induce an action of the full \emph{dihedral group} on an algebra over the coloured operad $\Guill_2$. In the same way as opposite algebras in dimension one are related to $*$-involution and transpositions (hence opening the route to $C^*$-algebra, $U(n)$ groups, etc), it would be very interesting to investigate mathematical consequences of this action of the dihedral group in dimension two.

\subsection{A discussion about units}
	\label{sec:eckmanhilton}
	
This section is devoted to a discussion that lies further from our main purposes of Section~\ref{sec:boundaryalgebra}. However, we think that it sheds some light on the relations between the $\Guill_2$ coloured operad and its topological counterparts (double semi-groups, $E_2$-operad, etc).

\subsection{A review on units for \texorpdfstring{$\Ass$}{Ass} and guillotine operads.} In dimension one, as often as possible, algebras with a unit are considered. From the operadic point of view, it corresponds to the specification of a new operation $e$ of arity $0$ in $\Ass_0$ such that the composition $m\circ(\id\otimes e)=m\circ(e\otimes \id)=\id$ maps the $2$-arity product $m$ to the unique $1$-arity element, which is the identity $\id$. At the level of algebras over $\Ass$, the $0$-arity corresponds, for an associative algebra $\ca{A}$, to a (commutative) morphism $e : \setK \to \ca{A}$, traditionally written $e(\lambda)=\lambda 1_{\ca{A}}$.

When introducing colours associated to segment lengths in the operad $\Guill_1$, the neutrality condition requires a unit $e(\lambda)$ to belong to $\Guill_1(0;\emptyset)$, i.e. degenerate segments: the absence of colours, denoted by $\emptyset$ after the semicolon ";", corresponds to null arity. This interpretation of identities as objects associated to points, i.e. $0$-dimensional objects has in fact a deeper geometric meaning: it is related to the fact that $\Guill_0$ could be defined as the operad $\Com$ as seen previously. 

In the canonical case $\ca{A}_0=\Diag( V(S)$ and $\ca{A}_p = \End(V(S)) \simeq \Mat_{|S|}(\setR)$, for $p\geq 1$, for a discrete-space Markov process, the product $m^{(0)}$ on $\ca{A}_0=\Diag( V(S)$  is commutative and this space contains the identity matrix $I$ and hence the unit $e(\lambda)=\lambda I$. Moreover it is related to the gauge conjugation of Markov transition matrices by diagonal positive matrices, which may be seen as the arbitrary gluing or not of infinitely thin points between segments.

\subsection{Dimension two: Eckmann-Hilton argument vs. coloured \texorpdfstring{$\Guill_2$}{Guill2} operad.}

The question of unit in $E_2$-operad-like structures with associativity relations similar to \eqref{eq:guill2:listassoc} is well-known under the name of Eckmann-Hilton argument \cite{eckmannhilton}: \emph{without colours}, if there exists a common neutral element $e$ for two products $m_{WE}$ and $m_{SN}$ that satisfy \eqref{eq:guill2:listassoc}, then both products are equal and commutative! Indeed, we have
\begin{align*}
m_{WE}(a,b) &= m_{WE}(m_{SN}(a,e),m_{SN}(e,b)) = m_{SN}(m_{WE}(a,e),m_{WE}(e,b)) = m_{SN}(a,b),
\\
m_{WE}(a,b) &= m_{WE}(m_{SN}(e,a),m_{SN}(b,e)) = m_{SN}(m_{WE}(e,b),m_{WE}(a,e)) = m_{SN}(b,a).
\end{align*}
The mechanism is similar to the one behind hidden commutativity already mentioned \cite{BagherzadehBremner,Kock} and is broken by the colour compatibility constraints in the $\Guill_2$ operad, which makes the previous equalities impossible to write. One also observes trivially from the canonical example $(\ca{T}_{p,q})$ of Theorem~\ref{theo:canonicalexampleGuill} that such a mechanism cannot hold.

However, in the presence of the colours $\PatternShapes(\patterntype{r})$, there is still a generalized notion of units obtained by mimicking the one-dimensional case and it leads to the interesting construction of a whole family of generalized units. It is inspired by the particular case of $(\ca{T}_{p,q})$ in which identities in the endomorphisms are natural candidates to build units.

In the present case, unit elements are necessarily are necessarily elements of spaces $\ca{A}_{p,q}$ associated to \emph{degenerate} rectangles, i.e. $p=0$ and $q=0$ and, in this case, they act as unit elements only for the products in the degenerate direction of the rectangles; in particular, the units may have a non-trivial product structure in the non-degenerate direction.  Unit elements necessarily belongs to the spaces associated to \emph{boundary spaces}, i.e. spaces $\ca{A}_{p,q}$ with $p=0$ and $q=0$. We however do not know how far this observation has consequences.

\begin{figure}
\begin{center}
\begin{tabular}{|l|c|}
\hline
 \text{Vertical units:} &
 $\begin{tikzpicture}[guillpart,scale=1.75]
 	\bulkrect{(0,0)}{(1,1)}{$a_{p,q}$};
 \end{tikzpicture}
 \otimes 
  \begin{tikzpicture}[guillpart,scale=1.75]
 	\bulkrect{(0,0)}{(0,1)}{};
 	\node at (0,0.5) [right] {$e_{0,q}$};
 \end{tikzpicture}
  \mapsto m_{WE}^{p,0 | q}(a_{p,q},e_{0,q}) = 
  \begin{tikzpicture}[guillpart,scale=1.75]
 	\bulkrect{(0,0)}{(1,1)}{$a_{p,q}$};
 \end{tikzpicture}$
 \\
 \hline
 \text{Horizontal units:} &
$\begin{matrix} \begin{tikzpicture}[guillpart,scale=1.75]
 	\bulkrect{(0,0)}{(1,1)}{$a_{p,q}$};
 \end{tikzpicture}
 \\
 \otimes 
 \\
  \begin{tikzpicture}[guillpart,scale=1.75]
 	\bulkrect{(0,0)}{(1,0)}{};
 	\node at (0.5,0) [below] {$e_{p,0}$};
 \end{tikzpicture}
 \end{matrix}
  \mapsto m_{SE}^{p | 0,q}(a_{p,q},e_{0,q}) = 
  \begin{tikzpicture}[guillpart,scale=1.75]
 	\bulkrect{(0,0)}{(1,1)}{$a_{p,q}$};
 \end{tikzpicture}$
 \\
 \hline
  \text{Concatenation of units:} &
$\begin{tikzpicture}[guillpart,scale=1.75]
 	\bulkrect{(0,0)}{(1.33,0)}{};
 	\node at (0.66,0) [below] {$e_{p_1,0}$};
 \end{tikzpicture}
 \otimes \begin{tikzpicture}[guillpart,scale=1.75]
 	\bulkrect{(0,0)}{(0.66,0)}{};
 	\node at (0.33,0) [below] {$e_{p_2,0}$};
 \end{tikzpicture}
  \mapsto m_{WE}^{p_1,p_2|0}(e_{p_1,0},e_{p_2,0}) = 
 \begin{tikzpicture}[guillpart,scale=1.75]
 	\bulkrect{(0,0)}{(2,0)}{};
 	\node at (1,0) [below] {$e_{p_1+p_2,0}$};
 \end{tikzpicture}$
 \\
 \hline
 \text{Superposition of units:} &
$\begin{matrix}\begin{tikzpicture}[guillpart,scale=1.75]
 	\bulkrect{(0,0)}{(1,0)}{};
 	\node at (0.5,0) [below] {$e_{p,0}$};
 \end{tikzpicture}
\\ \otimes \\
\begin{tikzpicture}[guillpart,scale=1.75]
 	\bulkrect{(0,0)}{(1,0)}{};
 	\node at (0.5,0) [below] {$e_{p,0}$};
 \end{tikzpicture}
 \end{matrix}
 \mapsto m_{SN}^{p|0,0}(e_{p_1,0},e_{p_2,0}) = 
 \begin{tikzpicture}[guillpart,scale=1.75]
 	\bulkrect{(0,0)}{(1,0)}{};
 	\node at (0.5,0) [below] {$e_{p,0}$};
 \end{tikzpicture}$
\\
\hline
 \text{Horiz. gluing with a point}
 & 
 $\begin{tikzpicture}[guillpart,scale=1.75]
 	\bulkrect{(0,0)}{(1.33,0)}{};
 	\node at (0.66,0) [below] {$e_{p,0}$};
 \end{tikzpicture}
 \otimes \begin{tikzpicture}[guillpart,scale=1.75]
 	\node (P) at (0.33,0) [inner sep=0.5mm,circle,fill] {};
 	\node at (P) [below] {$e_{0,0}$};
 \end{tikzpicture}
  \mapsto m_{WE}^{p,0|0}(e_{p,0},e_{0,0}) = 
 \begin{tikzpicture}[guillpart,scale=1.75]
 	\bulkrect{(0,0)}{(1.33,0)}{};
 	\node at (0.66,0) [below] {$e_{p,0}$};
 \end{tikzpicture}
 $
\\
\hline
 \end{tabular}
 \end{center}
\caption{Graphical description of the operations on coloured units for degenerate rectangles. Indices $(p,q)$ on the elements only indicate the space $\ca{A}_{p,q}$ to which they belong. We illustrate only and randomly some of the directions but the list is much longer.}
\label{fig:graphicalunits}
\end{figure}

	\section{More colours on boundaries: higher modules for \texorpdfstring{$\Guill_2$}{Guill2}-algebras}\label{sec:boundaryguill2}
	
The algebraic heart of the paper is the present section, from which section \ref{sec:invariantboundaryelmts} will follow. The previous sections may be seen as a reorganization of classical spaces, the $T_{p,q}(V(S_1),V(S_2))$, in terms of guillotine partitions of rectangles where products correspond to partial sums on boundary variables, as in the one-dimensional case of matrices. The major novelty in dimension two that does not exist in dimension one is the non-trivial geometric structure of the boundary of a rectangle, which is a closed curve with a length that depends on the rectangle. Moving from two-point boundaries in dimension 1 to one-dimensional arbitrary large boundaries in dimension two require a deep generalization of the notion of left and right modules and the step is performed below.

		\subsection{The one-dimensional case: a reminder and  a geometric interpretation}\label{sec:onedim:extendedlengthass}
		
\subsubsection{A first observation.} This section starts with a striking remark: whereas the theory of left and right modules and bi-modules is omnipresent in the literature about associative algebras and linear algebra, we did not find any appealing operadic descriptions for it. There exists a notion of module over an operad \cite{fressemoduleoveroperad} but it is different from the notion of module over an algebra on an operad. There also exists the particular case of Swiss-Cheese operad \cite{VoronovSwissCheese} related to $E_2$ operad, which is closer to the construction presented below but it does not contain all the ingredients ---such as the colours, edges and corners--- and is of topological nature.
		
\subsubsection{Left modules over an algebra.} The definition of a left module $\ca{M}$ over an algebra $\ca{A}$ with product $m:\ca{A}\otimes\ca{A}\to\ca{A}$, corresponds to the existence of a morphism $m':\ca{A} \otimes \ca{M} \to\ca{M}$ such that
\[
m' \circ (\id_{\ca{A}} \otimes m') = m'\circ(m\otimes \id_\ca{M})
\]
on $\ca{A}\otimes\ca{A}\otimes\ca{M}$ and follows directly from the associativity condition~\eqref{eq:associativity1D} on $m$ by replacing in a systematic way (in all compositions of arbitrary arity of $m$ and $m'$) the space $\ca{A}$ on the right by the set $\ca{M}$. From an operadic point of view, this could be seen by considering two colours, one for the algebra and one for the module, with suitable compatibility conditions; however, without further justification, this ad-hoc construction may be hard to generalize. The module colour is absorbing in the sense that, every time an operation has a module argument, then the result is again in the module.

The fact that there are exactly two notions (left and right) is indeed deeply related to the one-dimensional nature of $\mathrm{Ass}$.

\subsubsection{Switching to the \texorpdfstring{$\Guill_1$}{Guill1} operad.} Incorporating left and right modules can be easily done by enriching the notion of one-dimensional guillotine partitions and segment lengths with \emph{infinite shapes}. We introduce now definitions generalizing definitions~\ref{def:oneguillotpartition} and \ref{def:Guill1operad}.

\begin{defi}[left and right extended 1D guillotine partitions]
A right-extended (resp. left-extended) 1D guillotine partition is a guillotine partition of a discrete half-line $[u,+\infty)$ (resp. $(-\infty,u]$) with $u\in\setP$, i.e. a finite sequence of segments or half-lines $\ca{S}=([u_i,v_i])_{1\leq i\leq n}$ such that \begin{enumerate}[(i)]
	\item for all $1\leq i\leq n$, $u_i\leq v_i$ and $u_i\in\setP$ and $v_i\in \setP\sqcup\{+\infty\}$ (resp. $v_i\in\setP\sqcup\{-\infty\}$),
	\item for any distinct $i$ and $j$, the intersection of $[u_i,v_j]$ and $[u_j,v_j]$ has Lebesgue-negligible,
	\item $\cup_{1\leq i\leq n}[u_i,v_i]=[u,+\infty)$ (resp. $(-\infty,u]$).
\end{enumerate}
\end{defi}

Most properties of 1D guillotine partitions remain the same --- such as the existence of an ordering permutation $\sigma$ (with uniqueness excepted in case of degenerate segments with multiplicity), and the notion of compositions --- and one can check that there is exactly one $i\in\{1,\ldots,n\}$ such that $v_i=+\infty$. Proceeding to the definition of the extended operad $\Guill_{1,R}$ (resp. $\Guill_{1,L}$) requires then:\begin{itemize}
\item  the introduction of a new colour $\infty_R$ (resp. $\infty_L$), which corresponds to half-lines $[u,+\infty)$ (resp. $(-\infty,u]$), to the set of colours $\setL^*$ with the addition rule $\infty_R+ p =\infty_R$ for any $p\in\setL^*$,
\item  the introduction of right-extended 1D guillotine partitions to the set of guillotine partitions.
 \end{itemize}

One may remark that $\infty_R+\infty_R$ is not required to be defined since it makes no sense to glue together two half-lines $[u,+\infty)$ and $[u',+\infty)$: on the algebraic side, it corresponds to the fact that is no product nor action inside a module. The value $\infty_R$ is absorbent w.r.t the finite values in $\setL$, which corresponds algebraically to the fact that acting by elements of an algebra on an element of a module produces again an element of the module.

An algebra over the operad $\Guill_{1,R}$ now consists in spaces $(\ca{A}_c)_{c\in\setL^* \cup \{\infty_R\}}$ such that $ (\ca{A}_c)_{c\in\setL^*} $ provides an algebra over $\Guill_1$ and a new space $\ca{M}=\ca{A}_{\infty_R}$ (for the new colour) that plays precisely the role of left modules (the correspondence is clear in the constant case $\ca{A}_c=\ca{A}$ for $c\in\setL^*$ where $\ca{A}$ is an algebra): the action $\ca{A}_c\otimes \ca{M} \to \ca{M}$ corresponds to the partition $([0,c],[c,+\infty))$ of $[0,+\infty)$.

\subsubsection{From \texorpdfstring{$\Guill_1$}{Guill_1}- to \texorpdfstring{$\Ass$}{Ass}-algebras using boundary spaces}

We have already seen that any $\Ass$-algebra can be lifted to a $\Guill_1$-algebra and we now consider the converse operation: we state an important representation theorem in practice which establishes a relation between $\Guill_1$-algebras (i.e. with length colours) and standard $\Ass$-algebra using representation theory. Its consequences will reappear in Section~\ref{sec:eagerandlazyROPEs} in the context of fully-extended $\Guill_2$-algebras and again in Section~\ref{sec:invariantboundaryelmts} in the context of practical computations of generalized eigen-elements.

\begin{theo}[morphism from $\Guill_1$ to $\Ass$ through boundaries]\label{theo:1D:removingcoloursonboundaries}
Let $(\ca{A}_c)_{c\in\setL^*\cup\{\infty_R\}}$ be a $\Guill_{1,R}$-algebra. The set $\Hom(\ca{A}_{\infty_R},\ca{A}_{\infty_R})$ of linear maps on $\ca{A}_{\infty_R}$ is an $\Ass$-algebra. We define maps $\pi_p : \ca{A}_p \to \Hom(\ca{A}_{\infty_R},\ca{A}_{\infty_R})$ by 
\[
\pi_p(a)u = m_{[0,p],[p,+\infty)}(a,u) 
\]
for all $a\in\ca{A}_{p}$ and $u\in\ca{A}_{\infty_R}$. Then for all $n\in\setN^*$, for all $p_1,\ldots,p_n\in\setL^*$ with $P=\sum_{i=1}^n p_i$, for all sequence $(a_i)_{1\leq i\leq n}$ such that $a_i\in\ca{A}_{p_i}$ for all $1\leq i\leq n$,
\[
\pi_{P}\left(
m_{[0,p_1],\ldots, [P-p_n,P]}(a_1,\ldots,a_n)
\right)
=\pi_{p_1}(a_1)\circ \pi_{p_2}(a_2) \circ \ldots \circ \pi_{p_n}(a_n)
\]
where $\circ$ is the usual associative composition rule on $\Hom(\ca{A}_{\infty_R},\ca{A}_{\infty_R})$.
\end{theo}
\begin{proof}
This a simple exercise with the definitions of the two operads $\Guill_1$ and $\Ass$. For all $u\in\ca{A}_{\infty_R}$, we write the same quantity in two ways. We first have, for the guillotine partition $\rho=([0,p_1], \ldots, [P-p_n,P], [P,+\infty))$:
\begin{align*}
m_{\rho}(a_1,\ldots,a_n,u)&=m_{[0,p_1],[p_1,+\infty)}\left(
a_1, m_{[p_1,p_2],\ldots, [P-p_n,P],[P,+\infty)}(a_2,\ldots,a_n,u)
\right)
\\
&=\pi_{p_1}(a_1) m_{[p_1,p_2],\ldots, [P-p_n,P],[P,+\infty)}(a_2,\ldots,a_n,u)
\end{align*}
On the other hand, using the $\Guill_1$-associativity properties, we have
\begin{align*}
m_{[0,p_1],\ldots, [P-p_n,P],[P,+\infty)}(a_1,\ldots,a_n,u)
&= m_{[0,P],[P,+\infty)}\left( m_{[0,p_1],\ldots, [P-p_n,P]}(a_1,\ldots,a_n), u\right)
\\
&= \pi_P\left(m_{[0,p_1],\ldots, [P-p_n,P]}(a_1,\ldots,a_n)\right) u
\end{align*}
By a recursion on $p$, we obtain the expected result by combining the previous equation
\end{proof}

The palette of colours may then appear useless at first sight in dimension one, since, given a left or right module, one may switch from $\Guill_1$ to $\Ass$-algebras. However, this redundancy is the key point for the higher dimensional generalizations, as seen below, both to formulate the Markov property (going from $\Ass$ to $\Guill_1$) in Section~\ref{sec:canonicalguill2:markov} and compute boundary eigen-structures (going back from $\Guill_1$ to $\Ass$) in Section~\ref{sec:invariantboundaryelmts}.

\subsubsection{Fully extended 1D guillotine partitions.} We overview quickly line-extended 1D guillotine partitions defined as guillotine partition of the full line $\setR$. Besides the colour $\infty_R$ and $\infty_L$, there is an additional colour $\infty_{LR}$, defined to be $\infty_{LR} = \infty_R+\infty_L$ and an operad $\Guill_{1,LR}$ made of guillotine partitions of segments, of the two types of half-lines and of the full line.

An algebra over this operad consists in spaces $(\ca{A}_c)_{c\in\setL^*}$, a left-module $\ca{A}_{\infty_R}$, with action $m_R: \ca{A}_p\otimes\ca{A}_{\infty_R} \to \ca{A}_{\infty_R}$, a right-module $\ca{A}_{\infty_L}$, with action $m_R: \ca{A}_p\otimes\ca{A}_{\infty_R} \to \ca{A}_{\infty_R}$, and a vector space $\ca{A}_{\infty_{LR}}$ with a linear map $b: \ca{A}_{\infty_L} \otimes  \ca{A}_{\infty_R} \to  \ca{A}_{\infty_{LR}}$ such that 
\begin{equation}
\label{eq:pairing1D}
b(m_L(u,a),v)=b(u,m_R(a,v))
\end{equation}
where $a\in\ca{A}_p$ for any $p\in\setL^*$. This last relation is the counterpart of the associativity~\eqref{eq:associativity1D} where left and right spaces are replaced by the modules. The additional linear map $b$ corresponds to the equivalence class under translation made of the partitions \[
B=( (-\infty,u], [u,+\infty) )
\]
for any $u\in\setP$ of the infinite line
and the previous associativity condition~\eqref{eq:pairing1D} corresponds to the two decompositions of the following partition:
\begin{align*}
( (-\infty,0],[0,c], [c,+\infty) ) &= B\circ \left(\id_{\infty_L}, ([0,c], [c,+\infty)) \right) \\
 &= B\circ \left( ((-\infty,0], [0,c]),\id_{\infty_R} \right)
\end{align*}

\subsubsection{The canonical construction for one-dimensional discrete-space Markov processes.}
We have seen that, for a $S$-valued Markov process in dimension one, the $\Guill_1$-algebra is given by $\ca{A}_p=\End(V(S))\simeq \Mat_{|S|}(\setR)$. We may now consider the extension
\begin{align*}
\ca{A}_{\infty_L}&=V(S)^*
&
\ca{A}_{\infty_R} &= V(S)
& \ca{A}_{\infty_{LR}} &= \setR
\end{align*}
with the actions $m_R$ and $m_L$ defined as the natural actions on $\End(V(S))$ on $V(S)$ and its dual. The pairing $b$ is just the action of the dual space $V(S)^*$ on the initial space $V(S)$ with $b(u,v)=u(v)$ (or the canonical scalar product $\scal{u}{v}$ if $V(S)^*$ is identified to $V(S)$).

It makes now possible to consider boundary conditions $b_L \in \ca{A}_{\infty_L}$ and $b_R\in \ca{A}_{\infty_R}$ and the partition function
\[
Z^{\text{boundary}}_{\{1,\ldots,n\}}(\MarkovWeight{A}_\bullet; b_L, b_R) = \sum_{(x_1,x_n)\in S^2} Z_{\{1,\ldots,n\}}(\MarkovWeight{A}_\bullet; x_1,x_n) b_L(x_1) b_R(x_n) \in\ca{A}_{\infty_{LR}}
\]
which is now a real number, and not a function of the boundary variables due to the summations that can be interpreted as the action of the matrix $Z_{D}(\MarkovWeight{A}_\bullet; x_1,x_n) \in\ca{A}_n $.

There is an interesting fact about the "infinite" colours $\infty_L$ and $\infty_R$ in terms of probabilistic models. Boundary conditions $b_L$ and $b_R$ are often considered as "coming from infinity" in Gibbs measures on the whole set $\setZ$, or at least as weights that summarize the contribution from outer parts of the process on the segment.

\subsection{Two-dimensional guillotine partitions of infinite domains: basic notations.}		

The previous construction of left and right modules for algebras as 1D guillotine partitions extended to infinite shapes can be generalized in a straightforward way to the two-dimensional case by considering half-infinite or infinite rectangles. This is the purpose of the present section. Once the geometry is well understood, there is no additional technical tool to develop; however, there is a wide variety of (half)-infinite extensions to introduce and each of them has its own algebraic interpretation: the following definitions may thus look like a long list but cannot be shortened.

The two-dimensional patterns are obtained as Cartesian products of the previous one-dimensional patterns, i.e. segments, half-lines and the full line. The shapes will thus be elements of $(\setL^*\cup\{\infty_L,\infty_R,\infty_{LR}\})^2$.

\begin{defi}[two-dimensional admissible patterns]\label{def:admissiblepatterns}
	Two-dimensional admissible patterns are the elementary geometric shapes obtained from successive guillotine cuts in the plane $\setP^2$. They are given the following names and descriptions:
	\begin{itemize}
		\item 
as already introduced, a \emph{rectangle} in $\setP^2$ is a set $[u_1,v_1]\times [u_2,v_2]$. Its shape is the couple $(v_1-u_1, v_2-u_2)\in \setL^*\times \setL^*$.

\item A North \emph{half-strip} is a subset of $\setP^2$ defined by $[u_1,v_1]\times [u_2,+\infty)$. Its shape is the couple $(v_1-u_1,\infty_R)$. A South (resp. West, East) half-strip is defined in the same way where the value $-\infty$ (resp. $-\infty$, $+\infty$ replaces $u_2$ (resp. $u_1$, $v_1$) in the definition of a rectangle and has a shape $(v_1-u_1, \infty_L)$ (resp. $(\infty_L, v_2-u_2)$, $(\infty_R, v_2-u_2)$).

\item A South-North \emph{strip} (resp. West-East strip) is a subset defined by $[u_1,v_1]\times \setP$ (resp. $\setP\times [u_2,v_2]$). Its shape is the couple $(v_1-u_1,\infty_{LR})$ (resp. $(\infty_{LR},v_2-u_2)$).

\item A North-West \emph{corner} (resp. West-South, South-East, East-North corner) is a subset defined by $(-\infty,v_1]\times [u_2,+\infty)$ (resp. $(-\infty,v_1]\times (-\infty,v_2]$, $[u_1,+\infty)\times(-\infty,v_2]$, $[u_1,+\infty)\times [u_2,+\infty)$). Its shape is the couple $(\infty_L,\infty_R)$ (resp. $(\infty_L,\infty_L)$, $(\infty_R,\infty_L)$, $(\infty_R,\infty_R)$).

\item A North \emph{half-plane} (resp. West, South and East half-plane) is a subset defined by $\setP\times [u_2,+\infty)$ (resp. $(-\infty,v_1]\times\setP$, $\setP\times (-\infty,v_2]$, $[u_1,+\infty)\times\setP$. Its shape is $(\infty_{LR},\infty_R)$ (resp. $(\infty_L,\infty_{LR})$, $(\infty_{LR},\infty_L)$ and $(\infty_{R},\infty_{LR})$).

\item The full plane is the full plane $\setP^2$. Its shape is $(\infty_{LR},\infty_{LR})$.
	\end{itemize}
\end{defi}

Figure~\ref{fig:admissiblepatterns} provides a visual catalogue of such patterns, grouped by types. We also define the following set $\PatternTypes$ of pattern types:
\begin{align*}
\PatternTypes = \{ & \patterntype{{r}}, \patterntype{{hs}}_N, \patterntype{{hs}}_W, \patterntype{{hs}}_S, \patterntype{{hs}}_E, 
 \patterntype{{s}}_{SN}, \patterntype{{s}}_{WE}, \\ &
 \patterntype{{c}}_{NW}, \patterntype{{c}}_{SW}, \patterntype{{c}}_{SE}, \patterntype{{c}}_{NE},\patterntype{{hp}}_{N}, \patterntype{{hp}}_{S}, \patterntype{{hp}}_{W},\patterntype{{hp}}_{E}, \patterntype{{fp}} \}
\end{align*}
where the names are made of the initials of the pattern and its directions.

\begin{figure}
\begin{center}
\begin{tikzpicture}
%Plane
\begin{scope}[yshift=9cm]
	\begin{scope}
	\fill[guillfill] (0,0) rectangle (1,1);
	\node at (0.5,0.5) {$\setP^2$};	
	\node at (0,0.5) [anchor=east] {Full plane};
		\node at (0.5,0.) [anchor = north] {$(\infty_{WE},\infty_{SN})$};
	\end{scope}
\end{scope}

%Half-planes
\begin{scope}[yshift=7cm]
	\begin{scope}[xshift= 4.cm]
	\fill[guillfill] (0,0) rectangle (1,1);
	\draw[guillsep] (0,1) -- (1,1);
	%\node at (0.5,0) [] {$-\infty$};
	%\node at (0,0.5) [] {$-\infty$};
	%\node at (1,0.5) [] {$+\infty$};
		\node at (0.5,0.) [anchor = north] {$(\infty_{WE},\infty_S)$};
	\end{scope}
	
	\begin{scope}[xshift= 1.5cm]
	\fill[guillfill] (0,0) rectangle (1,1);
	\draw[guillsep] (0,0) -- (1,0);
	%\node at (0,0.5) [] {$-\infty$};
	%\node at (0.5,1) [] {$+\infty$};
	%\node at (1,0.5) [] {$+\infty$};
		\node at (0.5,0.) [anchor = north] {$(\infty_{WE},\infty_N)$};
	\end{scope}
	
	\begin{scope}[xshift= -1.5cm]
	\fill[guillfill] (0,0) rectangle (1,1);
	\draw[guillsep] (1,0) -- (1,1);
	%\node at (0.5,0) [] {$-\infty$};
	%\node at (0.5,1) [] {$+\infty$};
	%\node at (0,0.5) [] {$-\infty$};
		\node at (0.5,0.) [anchor = north] {$(\infty_W,\infty_{SN})$};
	\end{scope}
	
	\begin{scope}[xshift= -4cm]
	\fill[guillfill] (0,0) rectangle (1,1);
	\draw[guillsep] (0,0) -- (0,1);
	%\node at (0.5,0) [] {$-\infty$};
	%\node at (0.5,1) [] {$+\infty$};
	%\node at (1,0.5) [] {$+\infty$};
	\node at (0,0.5) [anchor = east] {Half-planes};
		\node at (0.5,0.) [anchor = north] {$(\infty_E,\infty_{SN})$};
	\end{scope}
\end{scope}

%Strips and corners
\begin{scope}[yshift=4.5cm]
	
	\node at (-2.5,1.5) {Corners};
	\draw (-5.5,1.1)-- (-5.5,1.3) -- (1.5,1.3) -- (1.5,1.1);

	\node at (5.,1.5) {Strips};
	\draw (3.5,1.1)-- (3.5,1.3) -- (6.5,1.3) -- (6.5,1.1);

	\begin{scope}[xshift= 5.5cm]
	\fill[guillfill] (0,0) rectangle (1,1);
	\draw[guillsep] (0,0) -- (0,1);
	\draw[guillsep] (1,1) -- (1,0);
	%\node at (0.5,0) [] {$-\infty$};
	%\node at (0.5,1) [] {$+\infty$};
		\node at (0.5,0.) [anchor = north] {$(p,\infty_{SN})$};
	\end{scope}

	\begin{scope}[xshift= 3.5cm]
	\fill[guillfill] (0,0) rectangle (1,1);
	\draw[guillsep] (0,0) -- (1,0); 
	\draw[guillsep] (1,1) -- (0,1);
	%\node at (0,0.5) [] {$-\infty$};
	%\node at (1,0.5) [] {$+\infty$};
		\node at (0.5,0.) [anchor = north] {$(\infty_{WE},q)$};
	\end{scope}
	
	\begin{scope}[xshift= 0.5cm]
	\fill[guillfill] (0,0) rectangle (1,1);
	\draw[guillsep]  (0,1) -- (1,1) -- (1,0);
	%\node at (0,0.5) [] {$-\infty$};
	%\node at (0.5,0) [] {$-\infty$};
		\node at (0.5,0.) [anchor = north] {$(\infty_W,\infty_S)$};
	\end{scope}

	\begin{scope}[xshift= -1.5cm]
	\fill[guillfill] (0,0) rectangle (1,1);
	\draw[guillsep] (1,1) -- (1,0) -- (0,0);
	%\node at (0,0.5) [] {$-\infty$};
	%\node at (0.5,1) [] {$+\infty$};
		\node at (0.5,0.) [anchor = north] {$(\infty_W,\infty_N)$};
	\end{scope}

	\begin{scope}[xshift= -3.5cm]
	\fill[guillfill] (0,0) rectangle (1,1);
	\draw[guillsep] (1,0) -- (0,0) -- (0,1);
	%\node at (0.5,1) [] {$+\infty$};
	%\node at (1,0.5) [] {$+\infty$};
		\node at (0.5,0.) [anchor = north] {$(\infty_E,\infty_N)$};
	\end{scope}
	
	\begin{scope}[xshift= -5.5cm]
	\fill[guillfill] (0,0) rectangle (1,1);
	\draw[guillsep] (0,0) -- (0,1) -- (1,1);
	%\node at (1,0.5) [] {$+\infty$};
	%\node at (0.5,0) [] {$-\infty$};
		\node at (0.5,0.) [anchor = north] {$(\infty_E,\infty_S)$};
	\end{scope}
\end{scope}

%half strips
\begin{scope}[yshift=2cm]
	\begin{scope}[xshift= 3cm]
	\fill[guillfill] (0,0) rectangle (1,1);
	\draw[guillsep] (0,0) -- (0,1) -- (1,1) -- (1,0);
	%\node at (0.5,0) [] {$-\infty$};
		\node at (0.5,0.) [anchor = north] {$(p,\infty_S)$};
	\end{scope}
	
	\begin{scope}[xshift= -3cm]
	\fill[guillfill] (0,0) rectangle (1,1);
	\draw[guillsep]  (0,1) -- (1,1) -- (1,0)-- (0,0);
	%\node at (0,0.5) [] {$-\infty$};
	\node at (-0.5,0.5) [anchor = east] {Half-strips};
		\node at (0.5,0.) [anchor = north] {$(\infty_W,q)$};
	\end{scope}

	\begin{scope}[xshift= 1cm]
	\fill[guillfill] (0,0) rectangle (1,1);
	\draw[guillsep] (1,1) -- (1,0) -- (0,0) -- (0,1);
	%\node at (0.5,1) [] {$+\infty$};
	\node at (0.5,0.) [anchor = north] {$(p,\infty_N)$};
	\end{scope}

	\begin{scope}[xshift= -1cm]
	\fill[guillfill] (0,0) rectangle (1,1);
	\draw[guillsep] (1,0) -- (0,0) -- (0,1)--(1,1);
	%\node at (1,0.5) [] {$+\infty$};
	\node at (0.5,0.) [anchor = north] {$(\infty_E,q)$};
	\end{scope}
\end{scope}

\begin{scope}
	\fill[guillfill] (0,0) rectangle (1,1);
	\draw[guillsep] (0,0) rectangle (1,1);
	\node at (0,0.5) [anchor = east] {Rectangles};
	\node at (0.5,0.) [anchor = north] {$(p,q)$};
\end{scope}
\end{tikzpicture}
\end{center}
\caption{\label{fig:admissiblepatterns}Five types of admissible patterns, with their pattern shape below. Any pattern in a line may appear in the guillotine partitions of some of the patterns above this line. The set $\setP$ stands for $\setZ^2$ in the discrete cases and $\setR^2$ in the continuous case. All finite sizes $p$ and $q$ are elements of $\setL^*$ ($\setN_1$ in the discrete case and $\setR_+^*$ in the continuous setting).}
\end{figure}

\begin{defi}[guillotine partition of an admissible pattern]
A guillotine partition of size $n$ of an admissible pattern $[u_1,v_1]\times[u_2,v_2]$ (here, any of the two boundary $u_1$ and $u_2$ may be replaced by $-\infty$ as well as $v_1$ and $v_2$ by $+\infty$) is a finite $n$-uplet $(R_1,\ldots,R_n)$ of admissible patterns satisfying the same three conditions as in Definition~\ref{def:guillotinepartition2D}.
\end{defi}

One remarks that there is a hierarchy of pattern types: rectangles can only contain rectangles, half-strips may contain half-strips and rectangles but not corners, corners can contain corners, half-strips and rectangles, strips can contain half-strips in the same directions but not in the other one, etc. We define a partial order $\preceq$ on the set of $\PatternTypes$ by $u\preceq v$ if and only if there exists a pattern of type $u$ that can be included in a pattern of type $v$. For every type $u\in\PatternTypes$, we introduce $D(u) = \{ v\in \PatternTypes; v \preceq u \}$. For example, one has $D(\patterntype{r})=\{\patterntype{r}\}$, $D(\patterntype{s}_{WE})=\{\patterntype{r},\patterntype{hs}_W,\patterntype{hs}_E,\patterntype{s}_{WE} \}$ and $D(\patterntype{fp})=\PatternTypes$.

As pointed in Definition~\ref{def:admissiblepatterns}, an important change is now the addition of more "infinite" values to the rectangle shapes. Indeed, rectangle shapes are obtained as quotient classes of rectangles under translation equivalence; quotient classes of patterns under translation equivalence require the introduction of the following palette of \emph{pattern shapes}:
\begin{equation}
\PatternShapes(\patterntype{fp}) = \big(\setL^* \sqcup \{\infty_W,\infty_E,\infty_{WE}\}\big) \times 
\big(\setL^* \sqcup \{\infty_S,\infty_N,\infty_{SN}\}\big)
\end{equation} 
where the additional values are dictated by the locations of $\pm \infty$ on the boundaries in the intuitive way. It is just a rewriting of the shapes of Definition~\ref{def:admissiblepatterns} with the substitution $L \text{ (left)}\leftrightarrow S$ (South) or $W$ (West) and $R \text{ (right)}\leftrightarrow N$ (North) or $E$ (East) depending on the coordinate in order to make notations easier to interpret.

For any pattern type $u\in \PatternTypes$, we thus introduce the corresponding subset of colours:
\[
\PatternShapes(u) = \{ (p,q)\in \PatternShapes  ; \exists x\text{ of type }v\in D(u) \text{ and with shape } (p,q) \} 
\]
and observe that $u\preceq v$ implies $\PatternShapes(u) \subset \PatternShapes(v)$.

\subsection{The extended guillotine operad, part one: boundary algebras and corner double modules.}\label{sec:extendedguillpartI}	

We may again define equivalence of guillotine partitions under translations exactly as in Definition~\ref{def:guillpart:equivtranslat}. 

\begin{lemm}[composition of extended guillotine partitions]
\label{lemma:equivandcompo:extendedguill:notranslation}
For any pattern type \[u\in\{\patterntype{r},\patterntype{hs}_N,\patterntype{hs}_S,\patterntype{hs}_W,\patterntype{hs}_E,\patterntype{c}_{NW},\patterntype{c}_{WS},\patterntype{c}_{SE},\patterntype{c}_{EN} \},
\] i.e. any pattern type with a pattern size that does not contain the value $\infty_{SN}$ nor $\infty_{WE}$, any equivalence class $[\rho]$ of extended guillotine partitions of type $u$ with representative $\rho=(R_1,\ldots,R_n)$ of size $n$ and shapes $((p,q);(p_1,q_1),\ldots,(p_n,q_n))$ all in $\PatternShapes(u)$ and any sequence $([\rho_i])_{1\leq i\leq n}$ of equivalence classes of extended guillotine partitions with representant $\rho_i=(R_1^{(i)},\ldots,R^{(i)}_{k_i})$ of type $u_i\in D(u)$ of size $k_i$ and shapes $((p_i,q_i);(p'_1,q'_1),\ldots,(p'_{k_i},q'_{k_i}))$ all in $\PatternShapes(v)$, the composition $[\rho]\circ([\rho_1],\ldots,[\rho_n])$ is well-defined as the equivalence class of the guillotine partition
\begin{equation}\label{eq:compoextendedguillpart}
\rho \circ (\rho_1,\ldots,\rho_n) = \left(\theta_1(R_1^{(1)}),\ldots,\theta_1(R_{k_1}^{(1)}),\ldots,\theta_n(R_1^{(n)}),\ldots,\theta_1(R_{k_n}^{(n)}) \right)
\end{equation}
where $\theta_i$ is the unique translation that maps $\cup_{j=1}^{k_i} R_j^{(i)}$ to $R_i$.

(A pointed versions of this lemma for the other pattern types is discussed for Theorem~\ref{theo:operadpointedversion} below.)
\end{lemm}

One remarks immediately that the particular case $u=\patterntype{r}$ is the previous Lemma~\ref{lemma:guillotinecompo}.

\begin{proof}
We first prove that, given two patterns $P$ and $P'$ with the same pattern type $u\in\{\patterntype{r},\patterntype{hs}_N,\patterntype{hs}_S,\patterntype{hs}_W,\patterntype{hs}_E,\patterntype{c}_{NW},\patterntype{c}_{WS},\patterntype{c}_{SE},\patterntype{c}_{EN} \}$ and the same pattern shape $(p,q)\in\PatternShapes(u)$, there exists a unique translation $\theta$ of $\setR^2$ such that $\theta(P)=P'$. The common feature of all these pattern types is that, for both axis, at least one extremity of the pattern is finite and, among all the patterns with same type, the finite extremity can be chosen to be in the same directions. For example, all the north-west corners can be written $(-\infty,w_1]\times[w_2,\infty)$ with both numbers $w_i$ finite. Using these finite extremities of $P$ and $P'$, we may define a translation $\theta(x_1,x_2)=(x_1+w'_1-w_1,x_2+w'_2-w_2)$.  Both patterns have the same shape if and only if $\theta(P)=P'$ (which can be checked using the opposite extremities). Moreover, this translation is unique.

Given extended partitions $\rho$ and $\rho_i$ with compatible shapes as stated in the lemma, then there is a unique way of shifting each partition $\rho_i$ of the pattern $\cup_{j=1}^{k_i} R_j^{(i)}$ to $R_i$ with a translation $\theta_i$. The composition $\rho\circ(\rho_1,\ldots,\rho_n)$ is well-defined and produces a new guillotine partition of the initial pattern $\cup_{j=1}^n R_j$ through \eqref{eq:compoextendedguillpart}.

Given equivalence classes $[\rho]$ and $[\rho_i]$ with compatible shapes, it is straightforward to check that any choice of representatives $\rho$ and $\rho_i$ of the equivalence classes leads to a composition $\rho \circ (\rho_1,\ldots,\rho_n)$ in the same equivalence class.
\end{proof}
		
Following the previous notations, we may now defined a whole hierarchy of coloured operads extending $\Guill_2$.
\begin{theo}[Guillotine operad with a corner]\label{theo:extendedguilloperads}
For any pattern type \[u\in
\{\patterntype{r},\patterntype{hs}_N,\patterntype{hs}_S,\patterntype{hs}_W,\patterntype{hs}_E,\patterntype{c}_{NW},\patterntype{c}_{SW},\patterntype{c}_{SE},\patterntype{c}_{NE} \},
\]
any $n\in\setN^*$ and any $((p,q),(p_1,q_1),\ldots,(p_n,q_n))\in\PatternShapes(u)^{n+1}$, let \[\Guill_2^{(u)}((p,q);(p_1,q_1),\ldots,(p_n,q_n))\] be the set of equivalence classes of guillotine partitions $(R_1,\ldots,R_n)$ of a pattern with shape $(p,q)$ into $n$ patterns $R_i$ with shapes $(p_i,q_i)$. Then, the collection of sets \[
\left(\left( \Guill_2^{(u)}(c;c_1,\ldots,c_n) \right)_{(c,c_1,\ldots,c_n)\in \PatternShapes(u)^{n+1}} \right)_{n\in\setN^*}
\]
with the compositions of equivalence classes of guillotine partitions defines a coloured operad $\Guill_2^{(u)}$, which admits, for any $w\in D(u)$, $\Guill_2^{(w)}$ as a sub-operad. 

(Pointed versions of this theorem for the other pattern types are discussed in Theorem~\ref{theo:operadpointedversion}).
\end{theo}

\begin{proof}
This is a direct consequence of Lemma~\ref{lemma:equivandcompo:extendedguill:notranslation} and the proof is a simple exercise of verification of the axioms of an operad using the definition \eqref{eq:compoextendedguillpart} of the compositions.

The sub-operad structure is a simple consequence of the partial order $\preceq$ defined previously on the pattern shapes.
\end{proof}

The previous definitions and properties are just verification of axioms: it becomes now interesting to see how all this structure fits with the classical structures of algebra and how this leads to interesting structural consequences for two-dimensional Markov processes. The next paragraphs exhibit other sub-operad of the various $\Guill_2^{(u)}$, which are the key ingredients for our next sections.

The extension of the $\Guill_2^{(\patterntype{r})}$-operad to half-infinite shapes does not change various other properties, in particular the generators. The following property holds, although we do not provide a proof which would be a copy \emph{mutatis mutandis} of the proofs of properties~\ref{prop:guill2generators} and \ref{prop:guill2elemassoc}.

\begin{prop}[generators of extended $\Guill_2$-algebras]
\label{prop:extendedgeneratorsI}
For any pattern type $u\in
\{\patterntype{r},\patterntype{hs}_N,\patterntype{hs}_S,\patterntype{hs}_W,\patterntype{hs}_E,\patterntype{c}_{NW},\patterntype{c}_{SW},\patterntype{c}_{SE},\patterntype{c}_{NE} \}$,
the $\Guill_2^{(u)}$ is generated by the same elementary partitions as in~\eqref{eq:defgenerators} where the values $r,s$ or $p-r,q-s$ can now take infinite values in $\{\infty_L,\infty_R\}$. The same fundamental associativity relations~\eqref{eq:guill2:listassoc} hold.
\end{prop}

It is standard in linear algebra to call $m_{WE}^{\infty_L,r' | q}$ and $m_{SN}^{p|\infty_L,s'}$ right actions and $m_{WE}^{r,\infty_R | q}$ and $m_{SN}^{p|s,\infty_R}$ left actions instead of products, as it will be seen below with much more details.

\subsubsection{A closer look to "boundary algebras".} Before switching back to the more general picture, we choose to discuss a little bit more the particular case $\Guill_2^{(u)}$ where $u\in \patterntype{hs}_E$, i.e. East half-strips are considered. We have $D(\patterntype{hs_E})=\{\patterntype{r},\patterntype{hs}_E\}$ and $\PatternShapes(\patterntype{hs}_E) = \PatternShapes(\patterntype{r}) \sqcup \left( \{\infty_E\}\times \setL^* \right)$.

Moreover, there is a trivial correspondence between $\{\infty_E\}\times \setL^*$ and the set of colours $\setL^*$ of $\Guill_1$. This correspondence is in fact deeper and holds between half-strips $[u_1,+\infty)\times [u_2,v_2]$ and segments $[u_2,v_2]$ of $\setP$. Horizontal guillotine cuts in East half-strips produces only East half-strips. Purely horizontal East half-strips partitions into half-strips are thus identified with one-dimensional guillotine partitions through the correspondence.

The following proposition deals with the four cases, North, West, South and East by symmetry.

\begin{prop}\label{prop:algebrasonsides}
An algebra over $\Guill_2^{(\patterntype{hs}_E)}$ (resp. $\Guill_2^{(\patterntype{hs}_N)}$, $\Guill_2^{(\patterntype{hs}_W)}$ and $\Guill_2^{(\patterntype{hs}_S)}$) is given by a collection of spaces $(\ca{A}_c)_{c\in\PatternShapes(\patterntype{hs}_E)}$ (resp. $c\in\PatternShapes(\patterntype{hs}_N)$, $c\in\PatternShapes(\patterntype{hs}_W)$ and $c\in\PatternShapes(\patterntype{hs}_S)$), which we choose to see as two collections, one indexed by $\PatternShapes(\patterntype{r})$ and one, written $(\ca{B}_n)$, indexed by $\setL^*$ with the correspondence $\ca{B}_r = \ca{A}_{\infty_E,r}$ (resp. $\ca{A}_{r,\infty_N}$, $\ca{A}_{\infty_W,r}$ and $\ca{A}_{r,\infty_S}$) with $r\in\setL^*$. Then, 
\begin{enumerate}[(i)]
\item $(\ca{A}_c)_{c\in\PatternShapes(\patterntype{r})}$ is an algebra over $\Guill_2^{(\patterntype{r})}$;
\item $(\ca{B}_r)_{r\in\setL^*}$ is an algebra over $\Guill_1$, that we call a \emph{boundary algebra};
\item for every $r\in\setL^*$, $(\ca{A}_{p,r})_{p\in\setL^*}$ is an algebra over $\Guill_1$ (resp. $(\ca{A}_{r,q})_{q\in\setL^*}$, $(\ca{A}_{p,r})_{p\in\setL^*}$ and $(\ca{A}_{r,q})_{q\in\setL^*}$) and $\ca{B}_r$ is a left (resp. left, right, right) module for this algebra in the sense of $\Guill_{1,R}$ (resp. $\Guill_{1,R}$, $\Guill_{1,L}$ and $\Guill_{1,L}$).
\end{enumerate}
\end{prop}
\begin{proof}
The first point is just Theorem~\ref{theo:extendedguilloperads} with $\patterntype{r} \preceq \patterntype{hs}_E$. The second point is a direct consequence of the identification of half-strips with one-dimensional segments. The first part of the third point has already been seen in \eqref{eq:suboperadGuill1Guill2} and the second part follows in the same way by a direct check of axioms of $\Guill_{1,R}$. For example, the composition structure of guillotine partitions implies (for the East case) the commutativity of the following diagrams, for any finite numbers $(r_1,r_2,r_3)\in(\setL^*)^3$:
\newsavebox{\hse}
\savebox{\hse}{
	\begin{tikzpicture}[guillpart]
	\fill [guillfill] (0,0)--(2,0)--(2,2)--(0,2)--(0,0);
	\draw [guillsep] (2,0)--(0,0)--(0,2)--(2,2);
	\draw [guillsep] (0,1) -- (2,1);
	\node at (1,0.5) {$1$};
	\node at (1,1.5) {$2$};
	\end{tikzpicture}    
}
\newsavebox{\hsebis}
\savebox{\hsebis}{
	\begin{tikzpicture}[guillpart]
	\fill [guillfill] (0,0)--(2,0)--(2,2)--(0,2)--(0,0);
	\draw [guillsep] (2,0)--(0,0)--(0,2)--(2,2);
	\draw [guillsep] (0,1) -- (2,1);
	\node at (1,0.5) {$2$};
	\node at (1,1.5) {$3$};
	\end{tikzpicture}    
}
\newsavebox{\hseter}
\savebox{\hseter}{
	\begin{tikzpicture}[guillpart]
	\fill [guillfill] (0,0)--(2,0)--(2,3)--(0,3)--(0,0);
	\draw [guillsep] (2,0)--(0,0)--(0,3)--(2,3);
	\draw [guillsep] (0,1) -- (2,1);
	\draw [guillsep] (0,2) -- (2,2);
	\node at (1,0.5) {$1$};
	\node at (1,1.5) {$2$};
	\node at (1,2.5) {$3$};
	\end{tikzpicture}    
}
\begin{equation}
\begin{tikzpicture}[baseline={(current bounding box.center)}]
\matrix (m) [matrix of math nodes,row sep=3em,column sep=4em,minimum width=2em]
  {
	\ca{B}_{r_1}\otimes \ca{B}_{r_2}\otimes \ca{B}_{r_3}
     &  &
	\ca{B}_{r_1+r_2}\otimes \ca{B}_{r_3}
     \\
	\ca{B}_{r_1}\otimes \ca{B}_{r_2+r_3}     
     & &
     \ca{B}_{r_1+r_2+r_3} \\};
\path[-stealth]
    (m-1-1) edge node [left] {$\id_1\otimes$\usebox{\hsebis}} (m-2-1)
    (m-1-1) edge node [above] {\usebox{\hse}$\otimes\id_3$} (m-1-3)
    (m-2-1) edge node [below] {\usebox{\hse}} (m-2-3)
    (m-1-3) edge node [right] {\usebox{\hse}} (m-2-3)
    (m-1-1) edge node [above right] {$\Phi$} (m-2-3);
\end{tikzpicture}
\qquad\text{where }
\Phi= \usebox{\hseter}
\label{eq:interchangedecomp:ex1}
\end{equation}
\end{proof}

For the East case again, rephrased in terms of horizontal and vertical products $m_{WE}$ and $m_{SN}$ related through the interchange relation \eqref{eq:guill2:interchangeassoc}, we thus see that $(\ca{B}_r)_{r\in\setL^*}$ carries a double structure of an algebra (with products similar to $m_{SN}$) and of a module with an action of $(\ca{A}_c)_{c\in\PatternShapes(\patterntype{r})}$ compatible with the horizontal products $m_{WE}$. The interchange relation \eqref{eq:guill2:interchangeassoc} thus states that this action is in fact a morphism of algebra. We use the same words of algebras, modules and actions for algebras over the \emph{coloured} $\Guill_1$ and $\Guill_{1,R}$ as for algebras over the operad $\Ass$ since both structures are identical excepted the colours and since, in practice, this will work in the same way in the computations of the next sections.

\subsubsection{A closer look to corner (double) modules.}

The $\Guill_2$ boundary zoology also contains corners. We consider here the case of South-East corners. In this case, we have $D(\patterntype{c}_{SE})=\{\patterntype{r},\patterntype{hs}_E,\patterntype{hs}_{S},\patterntype{c}_{SE}\}$ for shapes generated by guillotine partitions of a corner as illustrated by the following example of guillotine partition:
\[
\begin{tikzpicture}[guillpart,scale=1.5]
\draw[guillsep] (0,0) -- (4,0);
\draw[guillsep] (0,0) -- (0,-3);
\fill[guillfill] (0,0) -- (4,0) -- (4,-3) -- (0,-3);
\draw[guillsep] (0,-1) -- (4,-1);
\draw[guillsep] (1,-1) -- (1,-3);
\draw[guillsep] (2,0) -- (2,-3);
\node at (1,-0.5) {$\patterntype{r}$};
\node at (3,-0.5) {$\patterntype{hs}_E$};
\node at (0.5,-2) {$\patterntype{hs}_S$};
\node at (1.5,-2) {$\patterntype{hs}_S$};
\node at (3,-2) {$\patterntype{c}_{SE}$};
\end{tikzpicture}
\]
and the pattern shapes are given by $\PatternShapes(\patterntype{c}_{SE})= (\setL^*\sqcup \{\infty_E\})\times (\setL^*\sqcup \{\infty_S\})$.

An algebra $(\ca{A}_c)_{c\in\PatternShapes(\patterntype{c}_{SE})}$  over the operad $\Guill_2^{(\patterntype{c}_{SE})}$ consists now in four objects:
\begin{itemize}
\item an algebra $(\ca{A}_c)_{c\in\PatternShapes(\patterntype{r})}$ over $\Guill_2^{(\patterntype{r})}$,
\item an East boundary algebra $(\ca{B}^E_r)_{r\in\setL^*}$ over $\Guill_1$ using  $\ca{B}_r^E = \ca{A}_{\infty_E,r}$,
\item a South boundary algebra $(\ca{B}^S_r)_{r\in\setL^*}$ over $\Guill_1$ using  $\ca{B}_r^S = \ca{A}_{r,\infty_S}$,
\item a new space $\ca{M}^{SE}= \ca{A}_{\infty_E,\infty_S}$, called a \emph{corner double module}.
\end{itemize}

The boundary algebras are the ones described in the previous paragraph and correspond to the sub-operads $\Guill_2^{(\patterntype{hs}_E)}$ and $\Guill_2^{(\patterntype{hs}_S)}$. It is now interesting to explore the algebraic structure of $\ca{M}^{SE}$.

Considering only vertical guillotine cuts identifies $(\ca{B}^S_r)_{r\in\setL^*}$ and $\ca{M}^{SE} =: \ca{B}^S_{\infty_R}$ to an algebra over the extended guillotine operad $\Guill_{1,R}$ and thus $\ca{M}^{SE}$ appears to be a left module over $(\ca{B}^S_r)_{r\in\setL^*}$. The action corresponds to the guillotine partitions \[
a \triangleright m =
\begin{tikzpicture}[scale=0.5,baseline={(current bounding box.center)}]
\draw[guillsep] (0,0) -- (3,0);
\draw[guillsep] (0,0) -- (0,-2);
\fill[guillfill] (0,0) -- (3,0) -- (3,-2) -- (0,-2);
\draw[guillsep] (1,0) -- (1,-2);
\node at (0.5,-1) {$a$};
\node at (2,-1) {$m$};
\end{tikzpicture}
\]

Considering now only horizontal guillotine cuts of corners identifies $(\ca{B}^E_r)_{r\in\setL^*}$ and $\ca{M}^{SE} =: \ca{B}^E_{\infty_L}$ to an algebra over the extended guillotine operad $\Guill_{1,L}$ and thus $\ca{M}^{SE}$ appears to be a right module over $(\ca{B}^E_r)_{r\in\setL^*}$. The action corresponds to guillotine partition \[
m \triangleleft b =
\begin{tikzpicture}[scale=0.5,baseline={(current bounding box.center)}]
\draw[guillsep] (0,0) -- (2,0);
\draw[guillsep] (0,0) -- (0,-3);
\fill[guillfill] (0,0) -- (2,0) -- (2,-3) -- (0,-3);
\draw[guillsep] (0,-1) -- (2,-1);
\node at (1,-0.5) {$b$};
\node at (1,-2) {$m$};
\end{tikzpicture}
\]

Is thus a corner double module $\ca{M}^{SE}$ a bimodule? In whole generality, no. But it looks like. The interchange relation~\eqref{eq:guill2:interchangeassoc} corresponds here to the commutativity of the diagram:
\newsavebox{\squarediag}
\savebox{\squarediag}{
	\begin{tikzpicture}[scale=0.55,baseline={(current bounding box.center)}]
	\fill [guillfill] (0,0)--(2,0)--(2,2)--(0,2)--(0,0);
	\draw [guillsep] (2,2)--(0,2)--(0,0);
	\draw [guillsep] (0,1) -- (2,1);
	\draw [guillsep] (1,0) -- (1,2);
	\node at (0.5,0.5) {$2$};
	\node at (0.5,1.5) {$1$};
	\node at (1.5,0.5) {$3$};
	\node at (1.5,1.5) {$4$};
	\end{tikzpicture}    
}

\newsavebox{\vertsouthbox}
\savebox{\vertsouthbox}{
	\begin{tikzpicture}[scale=0.55,baseline={(current bounding box.center)}]
	\fill [guillfill] (0,0)--(1,0)--(1,2)--(0,2)--(0,0);
	\draw [guillsep] (0,0)--(0,2)--(1,2)--(1,0);
	\draw [guillsep] (0,1) -- (1,1);
	\node at (0.5,0.5) {$2$};
	\node at (0.5,1.5) {$1$};
	\end{tikzpicture}    
}

\newsavebox{\vertsoutheastbox}
\savebox{\vertsoutheastbox}{
	\begin{tikzpicture}[scale=0.55,baseline={(current bounding box.center)}]
	\fill [guillfill] (0,0)--(1,0)--(1,2)--(0,2)--(0,0);
	\draw [guillsep] (0,0)--(0,2)--(1,2);
	\draw [guillsep] (0,1) -- (1,1);
	\node at (0.5,0.5) {$3$};
	\node at (0.5,1.5) {$4$};
	\end{tikzpicture}    
}

\newsavebox{\verteastsouthbox}
\savebox{\verteastsouthbox}{
	\begin{tikzpicture}[scale=0.55,baseline={(current bounding box.center)}]
	\fill [guillfill] (0,0)--(1,0)--(1,2)--(0,2)--(0,0);
	\draw [guillsep] (0,0)--(0,2)--(1,2);
	\draw [guillsep] (0,1) -- (1,1);
	\node at (0.5,0.5) {$1$};
	\node at (0.5,1.5) {$2$};
	\end{tikzpicture}    
}

\newsavebox{\horizeastbox}
\savebox{\horizeastbox}{
	\begin{tikzpicture}[scale=0.55,baseline={(current bounding box.center)}]
	\fill [guillfill] (0,0)--(0,1)--(2,1)--(2,0)--(0,0);
	\draw [guillsep] (0,0)--(0,1)--(2,1);
	\draw [guillsep] (1,0) -- (1,1);
	\node at (0.5,0.5) {$1$};
	\node at (1.5,0.5) {$2$};
	\end{tikzpicture}    
}

\newsavebox{\horizonefourbox}
\savebox{\horizonefourbox}{
	\begin{tikzpicture}[scale=0.55,baseline={(current bounding box.center)}]
	\fill [guillfill] (0,0)--(0,1)--(2,1)--(2,0)--(0,0);
	\draw [guillsep] (2,0)--(0,0)--(0,1)--(2,1);
	\draw [guillsep] (1,0) -- (1,1);
	\node at (0.5,0.5) {$1$};
	\node at (1.5,0.5) {$4$};
	\end{tikzpicture}    
}
\newsavebox{\horiztwothreebox}
\savebox{\horiztwothreebox}{
	\begin{tikzpicture}[scale=0.55,baseline={(current bounding box.center)}]
	\fill [guillfill] (0,0)--(0,1)--(2,1)--(2,0)--(0,0);
	\draw [guillsep] (0,0)--(0,1)--(2,1);
	\draw [guillsep] (1,0) -- (1,1);
	\node at (0.5,0.5) {$2$};
	\node at (1.5,0.5) {$3$};
	\end{tikzpicture}    
}

\begin{equation}
\begin{tikzpicture}[baseline={(current bounding box.center)}]
\matrix (m) [matrix of math nodes,row sep=3em,column sep=4em,minimum width=2em]
  {
	\ca{A}_{p,q}\otimes \ca{B}^S_{p}\otimes \ca{M}^{SE} \otimes \ca{B}^E_{q} 
     &  &
	\ca{M}^{SE}\otimes \ca{B}^E_{q}
     \\
	\ca{B}^S_{p}\otimes \ca{M}^{SE}     
     & &
     \ca{M}^{SE}      \\};
\path[-stealth]
    (m-1-1) edge node [left] {\usebox{\vertsouthbox}$\otimes$\usebox{\vertsoutheastbox}} (m-2-1)
    (m-1-1) edge node [above] {\usebox{\horizonefourbox}$\otimes$\usebox{\horiztwothreebox}} (m-1-3)
    (m-2-1) edge node [below] {\usebox{\horizeastbox}} (m-2-3)
    (m-1-3) edge node [right] {\usebox{\verteastsouthbox}} (m-2-3)
    (m-1-1) edge node [above right] {$\Psi$} (m-2-3);
\end{tikzpicture}
\qquad\text{where }
\Psi= \usebox{\squarediag}
\label{eq:interchangedecomp:ex2}
\end{equation}
In the case of an uncoloured version of $\Guill_2$ similar to $E_2$-like operad with an Eckmann-Hilton unit $e$, the previous diagram specified with $e\in\ca{A}$ in first position would precisely endow $\ca{M}^{SE}$ with a $\ca{B}^S$-$\ca{B}^E$-bimodule structure. However, we have seen here that it is not the case. This question of bimodule property will reappear later in the probabilistic context however.

With respect to the palette of colours $\PatternShapes(\patterntype{hs}_{S})$ and $\PatternShapes( \patterntype{hs}_{E})$ on the boundary South and East $\Guill_1$-algebras, the corners play a special "erasing" role: indeed, on one hand, colours are \emph{required} in order to have well-defined South and East actions of the $\Guill_2^{(\patterntype{r})}$-algebras and, on the other hand, Theorem~\ref{theo:1D:removingcoloursonboundaries} can be applied since on each South and East side for which the corner introduces a $\Guill_{1,R}$- and $\Guill_{1,R}$-structure so that both side $\Guill_1$-algebras $\ca{B}^S_\bullet$ and $\ca{B}^E_\bullet$ can be represented as sub-$\Ass$-algebras of $\Hom(\ca{M}_{SE},\ca{M}_{SE})$. This heuristics on colours reappear at various stages in the next sections with always the same general approach: colours on a boundary shape allows for actions of higher dimensional objects whereas "corners" (or more generally a shape with one more $\infty_{L}$ or $\infty_{R}$ value, i.e. a lower-dimensional shape) allows for mappings to simpler operads with less colours.

\begin{rema}
We have dealt here with the case of South-East corner but the same structure holds for the three other corners \emph{mutatis mutandis} by rotating the tensor products in the diagrams.
\end{rema}

\subsubsection{A useful notation for algebras.}

For any pattern type $\patterntype{x}$ listed in Proposition~\ref{prop:extendedgeneratorsI}, an algebra is a collection of spaces $(\ca{A}_{p,q})_{(p,q)\in \PatternShapes(\patterntype{x})}$: in order to lighten indices, we choose to write such a collection of spaces as:
\begin{equation}\label{eq:shortnotation:guillalgebra}
\ca{A}_{\PatternShapes(\patterntype{x})}
\end{equation}
in order to recall the set to which indices belong, since considering extensions and restrictions of such collections will become frequent in the next sections.

\subsection{The extended operads, part two: the pointed versions of strips, half-planes and full plane.}	\label{sec:extendedguillpartII}

In the previous cases, we use a equivalence relation on guillotine partitions based on translational invariance. The idea is that, given a specified corner of an admissible pattern and a shape $(p,q)$, one can reconstruct the admissible pattern. The translation of an admissible pattern to another one with the same shape is then easily defined from the reference corners in a unique way. In absence of corners, the situation is more subtle and requires to break the translation invariance of the corner-less shapes.

We illustrate this translational ambiguity in the case of \emph{strips}, through the example of West-East strips, in order to see the need for pointed versions of extended partitions. We consider the following two guillotine partitions:
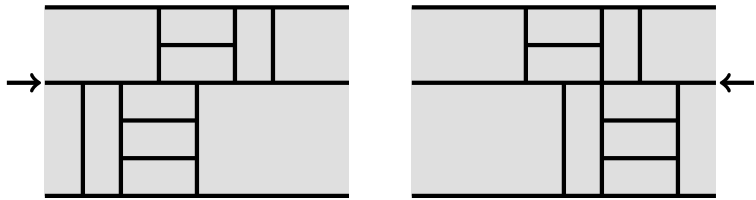
\[
\begin{tikzpicture}[guillpart]
\fill[guillfill] (0,0) rectangle (8,5);
\draw[guillsep] (0,0) -- (8,0);
\draw[guillsep] (0,3) -- (8,3);
\draw[guillsep] (0,5) -- (8,5);

\draw[guillsep] (1,0)--(1,3);
\draw[guillsep] (2,0)--(2,3);
\draw[guillsep] (4,0)--(4,3);
\draw[guillsep] (2,1)--(4,1);
\draw[guillsep] (2,2)--(4,2);

\draw[guillsep] (3,3)--(3,5);
\draw[guillsep] (5,3)--(5,5);
\draw[guillsep] (6,3)--(6,5);
\draw[guillsep] (3,4)--(5,4);
\draw[->,ultra thick] (-1,3) -- (-0.1,3);
\end{tikzpicture}
\qquad
\begin{tikzpicture}[guillpart]
\fill[guillfill] (0,0) rectangle (8,5);
\draw[guillsep] (0,0) -- (8,0);
\draw[guillsep] (0,3) -- (8,3);
\draw[guillsep] (0,5) -- (8,5);
\begin{scope}[xshift=3cm];
\draw[guillsep] (1,0)--(1,3);
\draw[guillsep] (2,0)--(2,3);
\draw[guillsep] (4,0)--(4,3);
\draw[guillsep] (2,1)--(4,1);
\draw[guillsep] (2,2)--(4,2);
\end{scope}
\draw[guillsep] (3,3)--(3,5);
\draw[guillsep] (5,3)--(5,5);
\draw[guillsep] (6,3)--(6,5);
\draw[guillsep] (3,4)--(5,4);

\draw[->,ultra thick] (9,3) -- (8.1,3);
\end{tikzpicture}
\]
Considering the horizontal cut pointed by the arrows, both of them would be obtained by the same composition of the two same guillotine partitions of a strip and there is no reason for the two guillotine partitions to produce similar operadic operations. 

However, there is an easy trick to fix this composition problem through the introduction of \emph{pointed} guillotine partitions of these patterns, with an adequate notion of change of base points coherent with the previous constructions. Moreover, these pointings of guillotine partitions have a deeper computational and physical interpretations that we choose to delay to other works in preparation since, in the next sections, this problem will not appear directly: thus we only give for completeness a brief overview of how it works, with only the main definitions and theorems, without spending lines of text about the possible interpretations.

\begin{defi}\label{def:pointedguillotine}
A \emph{pointed} guillotine partition of a strip $S=\setP\times [u_2,v_2]$ (resp. a strip $S=[u_1,v_1]\times \setP$, a half-plane $S$ of any type) is a pair $(z,\rho)$ where $z\in \setP$ and $\rho$ is an extended guillotine partition of $S$. A \emph{pointed} guillotine partition of the full plane is a pair $((z_1,z_2),\rho)$ where $(z_1,z_2)\in \setP^2$ and $\rho$ is a guillotine partition of $\setP^2$.

By convention, for a given pointed partition $(z,(R_1,\ldots,R_n))$, any sub-pattern $R_i$ with a type $u \in \{ 
 \patterntype{{s}}_{SN}, \patterntype{{s}}_{WE}, 
 \patterntype{{hp}}_{N}, \patterntype{{hp}}_{S}, \patterntype{{hp}}_{W},\patterntype{{hp}}_{E}, \patterntype{{fp}} \}$ may be canonically considered as a pointed pattern $(z,R_i)$.

The \emph{pattern type} of pointed horizontal strip (resp. horizontal strips, north, west, south, east half planes and the full plane) is noted $\patterntype{s}_{WE}^*$ (resp. $\patterntype{s}_{SN}^*$,
$\patterntype{hp}_N^*$,$\patterntype{hp}_W^*$,$\patterntype{hp}_S^*$,$\patterntype{hp}_E^*$ and $\patterntype{fp}^*$).
\end{defi}
For horizontally (resp. vertically) infinite strips and half-planes, $z$ has to be interpreted as a $x$-coordinate (resp. a $y$-coordinate). We choose to represent it graphically by a dotted line between opposite boundaries and thick points at the intersection between dotted lines or between a dotted line and a finite boundary (West-East strip on the left, East half-plane in the middle and full plane on the right):
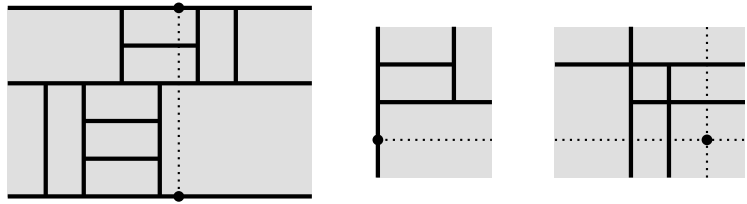
\begin{equation*}
\begin{tikzpicture}[guillpart]
\fill[guillfill] (0,0) rectangle (8,5);
\draw[guillsep] (0,0) -- (8,0);
\draw[guillsep] (0,3) -- (8,3);
\draw[guillsep] (0,5) -- (8,5);

\draw[guillsep] (1,0)--(1,3);
\draw[guillsep] (2,0)--(2,3);
\draw[guillsep] (4,0)--(4,3);
\draw[guillsep] (2,1)--(4,1);
\draw[guillsep] (2,2)--(4,2);

\draw[guillsep] (3,3)--(3,5);
\draw[guillsep] (5,3)--(5,5);
\draw[guillsep] (6,3)--(6,5);
\draw[guillsep] (3,4)--(5,4);

\node at (4.5,0) [circle, fill, inner sep=0.5mm] {};
\node at (4.5,5) [circle, fill, inner sep=0.5mm] {};
\draw[thick, dotted] (4.5,0)--(4.5,5);
\end{tikzpicture}
\qquad
\begin{tikzpicture}[guillpart]
\fill[guillfill] (0,0) rectangle (3,4);
\draw[guillsep] (0,2) -- (3,2);
\draw[guillsep] (0,3) -- (2,3);
\draw[guillsep] (2,2) -- (2,4);
\draw[guillsep] (0,0)--(0,4);

\draw[thick,dotted] (0,1)--(3,1);
\node at (0,1) [circle,fill,inner sep=0.5mm] {};
\end{tikzpicture}
\qquad
\begin{tikzpicture}[guillpart]
\fill[guillfill] (0,0) rectangle (5,4);
\draw[guillsep] (0,3) -- (5,3);
\draw[guillsep] (2,0) -- (2,4);
\draw[guillsep] (3,0) -- (3,3);
\draw[guillsep] (2,2)--(5,2);
\draw[thick,dotted] (4,0)--(4,4);
\draw[thick,dotted] (0,1)--(5,1);
\node at (4,1) [circle,fill,inner sep=0.5mm] {};
\end{tikzpicture}
\end{equation*}

The pattern shapes do not vary and we set $\PatternShapes(\patterntype{x}^*)=\PatternShapes(\patterntype{x})$ for any pattern type $\patterntype{x}\in\PatternTypes$ with at least one doubly-infinite dimension. The set of descendants $D(\patterntype{x}^*)$ of such a pointed pattern contains the same patterns as $D(\patterntype{x})$, excepted that patterns with at least one doubly-infinite dimensions are replaced by their pointed versions. One has for example:
\[
D(\patterntype{hp}_N^*)= \{ \patterntype{hp}_N^*, \patterntype{s}_{WE}^*, \patterntype{c}_{NW},\patterntype{c}_{NE}, \patterntype{hs}_{N}, \patterntype{hs}_{W}, \patterntype{hs}_{E}, \patterntype{r} \}
\]

\begin{lemm}\label{lemm:translationforpointedpatterns}
Given two pointed patterns $(z,M)$ and $(z',M')$ of this type $\patterntype{x}^* \in \{  
\patterntype{s}_{WE}^*, \patterntype{hp}_N^*,\patterntype{hp}_S^*
\}$ (resp. $\in \{ \patterntype{s}_{SN}^* ,\patterntype{hp}_W^*,\patterntype{hp}_E^*\}$ and $\in \{\patterntype{fp}^*\}$), the following statements are equivalent:
\begin{enumerate}[(i)]
\item \label{item:sameshape} $M$ and $M'$ have the same shape
\item \label{item:existstransl}there exists a translation $\theta: (x_1,x_2)\mapsto (x_1+a_1,x_2+a_2)$ of $\setP^2$ such that $M'=\theta(M)$ and $z'= z+a_1$ (resp. $z'=z+a_2$ and $z'=\theta(z)$) 
\end{enumerate}
If both statements are true, then the translation is uniquely defined.
\end{lemm}
\begin{proof}
From \eqref{item:existstransl} to \eqref{item:sameshape}, it is trivial since translations are isometries: the finite lengths in the shapes are unchanged under translations and one also checks that it is also the case for the infinite values.

We now assume \eqref{item:sameshape}. Along each axis, either there is a finite boundary at coordinates $w$ and $w'$ in the same direction for both patterns $M$ and $M'$ (left or right, top or bottom, depending on the pattern type), or the $z$ and $z'$ correspond a coordinate along this axis. We then define for each axis $a_i = w'-w$ or $a_i = z'-z$. This defines the shifts of the translation $\theta$ and, in case of ambiguity in the choice of the finite boundary, \eqref{item:sameshape} assures that the $a_i$ do not depend on this choice. It is then easy to check that $\theta(M)=M'$ by enumerating cases. Moreover, the previous choice of the $a_i$ is unique. 
\end{proof}

We may now strictly follow the precedent construction of Section~\ref{sec:extendedguillpartI} and define equivalence classes of pointed patterns under translations and compositions of such equivalence classes: the previous Lemma~\ref{lemm:translationforpointedpatterns} provides the unique translations required at the beginning of the proof of Lemma~\ref{lemma:equivandcompo:extendedguill:notranslation}. We do not provide the detailed proof since it is just another case study with the same arguments.

As an example, the previous pointed guillotine partition of a strip can be written unambiguously as the following composition:
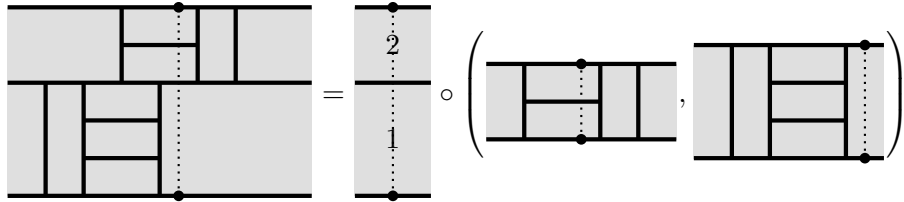
\begin{equation*}
\begin{tikzpicture}[guillpart]
\fill[guillfill] (0,0) rectangle (8,5);
\draw[guillsep] (0,0) -- (8,0);
\draw[guillsep] (0,3) -- (8,3);
\draw[guillsep] (0,5) -- (8,5);

\draw[guillsep] (1,0)--(1,3);
\draw[guillsep] (2,0)--(2,3);
\draw[guillsep] (4,0)--(4,3);
\draw[guillsep] (2,1)--(4,1);
\draw[guillsep] (2,2)--(4,2);

\draw[guillsep] (3,3)--(3,5);
\draw[guillsep] (5,3)--(5,5);
\draw[guillsep] (6,3)--(6,5);
\draw[guillsep] (3,4)--(5,4);

\node at (4.5,0) [circle, fill, inner sep=0.5mm] {};
\node at (4.5,5) [circle, fill, inner sep=0.5mm] {};
\draw[thick, dotted] (4.5,0)--(4.5,5);
\end{tikzpicture}
=
\begin{tikzpicture}[guillpart]
\fill[guillfill] (0,0) rectangle (2,5);
\draw[guillsep] (0,0)--(2,0);
\draw[guillsep] (0,3)--(2,3);
\draw[guillsep] (0,5)--(2,5);
\node at (1,1.5) {$1$};
\node at (1,4) {$2$};
\node at (1,0) [circle, fill, inner sep=0.5mm] {};
\node at (1,5) [circle, fill, inner sep=0.5mm] {};
\draw[thick, dotted] (1,0)--(1,5);
\end{tikzpicture}
\circ
\left(
\begin{tikzpicture}[guillpart]
\fill[guillfill] (0,0) rectangle (5,2);
\draw[guillsep] (0,0)--(5,0);
\draw[guillsep] (0,2)--(5,2);
\draw[guillsep] (1,0) -- (1,2);
\draw[guillsep] (3,0) -- (3,2);
\draw[guillsep] (4,0) -- (4,2);
\draw[guillsep] (1,1) -- (3,1);
\node at (2.5,0) [circle, fill, inner sep=0.5mm] {};
\node at (2.5,2) [circle, fill, inner sep=0.5mm] {};
\draw[thick, dotted] (2.5,0)--(2.5,2);
\end{tikzpicture}
,
\begin{tikzpicture}[guillpart]
\fill[guillfill] (0,0) rectangle (5,3);
\draw[guillsep] (0,0)--(5,0);
\draw[guillsep] (0,3)--(5,3);
\draw[guillsep] (2,1)--(4,1);
\draw[guillsep] (2,2)--(4,2);
\draw[guillsep] (1,0)--(1,3);
\draw[guillsep] (2,0)--(2,3);
\draw[guillsep] (4,0)--(4,3);
\node at (4.5,0) [circle, fill, inner sep=0.5mm] {};
\node at (4.5,3) [circle, fill, inner sep=0.5mm] {};
\draw[thick, dotted] (4.5,0)--(4.5,3);
\end{tikzpicture}
\right)
\end{equation*}
where we have canonically pointed the two sub-patterns with the initial point $z$. The first sub-pattern can now be written unambiguously:
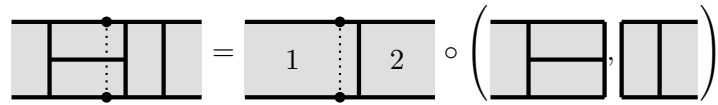
\begin{equation*}
\begin{tikzpicture}[guillpart]
\fill[guillfill] (0,0) rectangle (5,2);
\draw[guillsep] (0,0)--(5,0);
\draw[guillsep] (0,2)--(5,2);
\draw[guillsep] (1,0) -- (1,2);
\draw[guillsep] (3,0) -- (3,2);
\draw[guillsep] (4,0) -- (4,2);
\draw[guillsep] (1,1) -- (3,1);
\node at (2.5,0) [circle, fill, inner sep=0.5mm] {};
\node at (2.5,2) [circle, fill, inner sep=0.5mm] {};
\draw[thick, dotted] (2.5,0)--(2.5,2);
\end{tikzpicture}
=
\begin{tikzpicture}[guillpart]
\fill[guillfill] (0,0) rectangle (5,2);
\draw[guillsep] (0,0)--(5,0);
\draw[guillsep] (0,2)--(5,2);
\draw[guillsep] (3,0) -- (3,2);
\node at (2.5,0) [circle, fill, inner sep=0.5mm] {};
\node at (2.5,2) [circle, fill, inner sep=0.5mm] {};
\node at (1.25,1) {$1$};
\node at (4,1) {$2$};
\draw[thick, dotted] (2.5,0)--(2.5,2);
\end{tikzpicture}
\circ \left(
\begin{tikzpicture}[guillpart]
\fill[guillfill] (0,0) rectangle (3,2);
\draw[guillsep] (0,0)--(3,0);
\draw[guillsep] (0,2)--(3,2);
\draw[guillsep] (1,0) -- (1,2);
\draw[guillsep] (3,0) -- (3,2);
\draw[guillsep] (1,1) -- (3,1);
\end{tikzpicture}
,
\begin{tikzpicture}[guillpart]
\fill[guillfill] (3,0) rectangle (5,2);
\draw[guillsep] (3,0)--(5,0);
\draw[guillsep] (3,2)--(5,2);
\draw[guillsep] (3,0) -- (3,2);
\draw[guillsep] (4,0) -- (4,2);
\end{tikzpicture}
\right)
\end{equation*}
where the two sub-patterns are not doubly-infinite and do not require any pointing, since the translations may be defined relatively to the guillotine cut.

The previous Theorem~\ref{theo:extendedguilloperads} can then be generalized to the pointed patterns directly.

\begin{theo}[Guillotine operad with pointed doubly-infinite shapes]\label{theo:operadpointedversion}
For any pattern type \[u^*\in
\{\patterntype{s}_{WE}^*,\patterntype{s}_{SN}^*,\patterntype{hp}_N^*,\patterntype{hp}_W^*,\patterntype{hp}_S^*,\patterntype{hp}_{E}^*,\patterntype{fp}^* \},
\] the collection of equivalence classes of
\emph{pointed} guillotine partitions of patterns $v\in D(u^*)$, coloured with their shapes in $\PatternShapes(u^*)$, defines a coloured operad $\Guill_2^{(u^*)}$ which admits, for any $w\in D(u^*)$, $\Guill_2^{(w)}$ as a sub-operad. 
\end{theo}
Again, the proof is just an exercise on operads using the previous definitions and mimics exactly the unpointed case.

We discuss a bit more the consequences of the base points on the structure of these operads. The base points do not appear in the palette $\PatternShapes(u^*)=\PatternShapes(u)$ and thus an algebra over $\Guill_2^{(u^*)}$ is a collection of sets $(\ca{A}_c)_{c \in \PatternShapes(u)}$. The base points appear only in the products between these sets and, in particular, it changes slightly the generators as explained in the following proposition. 

\begin{prop}\label{prop:extendedgeneratorsII}
The $\Guill_2^{(\patterntype{fp}^*)}$ is generated by the following generators:
\begin{itemize}
\item the generators of the sub-operads $\Guill_2^{(u)}$ for the unpointed pattern types \[u\in\{\patterntype{r},\patterntype{hs}_N,\patterntype{hs}_W,\patterntype{hs}_S,\patterntype{hs}_E,\patterntype{c}_{SW},\patterntype{c}_{SE},\patterntype{c}_{NW},\patterntype{c}_{NE}\},\]
\item the set of pointed partitions of size $2$ indexed by the scalar $x$ and $y$ for strips
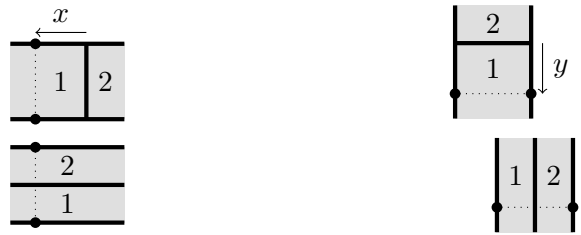
\begin{align*}
\begin{tikzpicture}[guillpart]
	\fill[guillfill] (3,0)--(0,0)--(0,2)--(3,2)--(3,0);
	\draw[guillsep] (3,0)--(0,0);
	\draw[guillsep] (0,2)--(3,2);
	\draw[guillsep] (2,0)--(2,2);
	\node at (1.33,1.) {$1$};
	\node at (2.5,1.) {$2$};	
	\node at (0.66,0) [circle, fill, inner sep=0.5mm] {};
	\node at (0.66,2) [circle, fill, inner sep=0.5mm] {};
	\draw [dotted] (0.66,0) -- (0.66,2);
		\draw [->] (2,2.3) -- node [midway, above] {$x$} (0.66,2.3);
\end{tikzpicture}
&
&
\begin{tikzpicture}[guillpart,rotate=90]
	\fill[guillfill] (3,0)--(0,0)--(0,2)--(3,2)--(3,0);
	\draw[guillsep] (3,0)--(0,0);
	\draw[guillsep] (0,2)--(3,2);
	\draw[guillsep] (2,0)--(2,2);
	\node at (1.33,1.) {$1$};
	\node at (2.5,1.) {$2$};	
	\node at (0.66,0) [circle, fill, inner sep=0.5mm] {};
	\node at (0.66,2) [circle, fill, inner sep=0.5mm] {};
	\draw [dotted] (0.66,0) -- (0.66,2);
		\draw [->] (2,-0.3) -- node [midway, right] {$y$} (0.66,-0.3);
\end{tikzpicture}
\\
\begin{tikzpicture}[guillpart]
	\fill[guillfill] (3,0)--(0,0)--(0,2)--(3,2)--(3,0);
	\draw[guillsep] (3,0)--(0,0);
	\draw[guillsep] (0,2)--(3,2);
	\draw[guillsep] (0,1)--(3,1);
	\node at (1.5,0.5) {$1$};
	\node at (1.5,1.5) {$2$};	
	\node at (0.66,0) [circle, fill, inner sep=0.5mm] {};
	\node at (0.66,2) [circle, fill, inner sep=0.5mm] {};
	\draw [dotted] (0.66,0) -- (0.66,2);
\end{tikzpicture}
&
&
\begin{tikzpicture}[guillpart,rotate=90]
	\fill[guillfill] (2.5,0)--(0,0)--(0,2)--(2.5,2)--(2.5,0);
	\draw[guillsep] (2.5,0)--(0,0);
	\draw[guillsep] (0,2)--(2.5,2);
	\draw[guillsep] (0,1)--(2.5,1);
	\node at (1.5,0.5) {$2$};
	\node at (1.5,1.5) {$1$};	
	\node at (0.66,0) [circle, fill, inner sep=0.5mm] {};
	\node at (0.66,2) [circle, fill, inner sep=0.5mm] {};
	\draw [dotted] (0.66,0) -- (0.66,2);
\end{tikzpicture}
\end{align*}
\item the corner pairings and strip action for the South half-planes
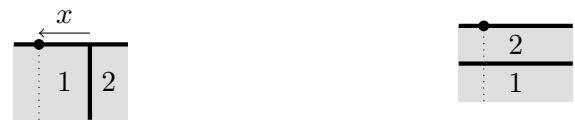
\begin{align*}
\begin{tikzpicture}[guillpart]
	\fill[guillfill] (3,0)--(0,0)--(0,2)--(3,2)--(3,0);
	\draw[guillsep] (0,2)--(3,2);
	\draw[guillsep] (2,0)--(2,2);
	\node at (1.33,1.) {$1$};
	\node at (2.5,1.) {$2$};	
	\node at (0.66,2) [circle, fill, inner sep=0.5mm] {};
	\draw [dotted] (0.66,0) -- (0.66,2);
		\draw [->] (2,2.3) -- node [midway, above] {$x$} (0.66,2.3);
\end{tikzpicture}
&
&
\begin{tikzpicture}[guillpart]
	\fill[guillfill] (3,0)--(0,0)--(0,2)--(3,2)--(3,0);
	\draw[guillsep] (0,2)--(3,2);
	\draw[guillsep] (0,1)--(3,1);
	\node at (1.5,0.5) {$1$};
	\node at (1.5,1.5) {$2$};	
	\node at (0.66,2) [circle, fill, inner sep=0.5mm] {};
	\draw [dotted] (0.66,0) -- (0.66,2);
\end{tikzpicture}
\end{align*}
indexed by the distance between the cut and the base point $x$,
and the equivalent ones for the three other half-planes.
\item the half-planes pairings for the full plane, indexed by the distance between the guillotine cut and the base point:
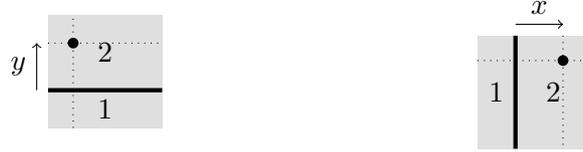
\begin{align*}
\begin{tikzpicture}[guillpart]
	\fill[guillfill] (0,0) rectangle (3,3);
	\draw[guillsep] (0,1)--(3,1);
	\node at (1.5,0.5) {$1$};
	\node at (1.5,2) {$2$};
	\draw [dotted] (0.66,0) -- (0.66,3);
	\draw [dotted] (0,2.25) -- (3,2.25);
	\draw [->] (-0.3,1) -- node [midway, left] {$y$} (-0.3,2.25);
	\node at (0.66,2.25) [circle, fill, inner sep=0.5mm] {};
\end{tikzpicture}
&
&
\begin{tikzpicture}[guillpart,rotate=-90]
	\fill[guillfill] (0,0) rectangle (3,3);
	\draw[guillsep] (0,1)--(3,1);
	\node at (1.5,0.5) {$1$};
	\node at (1.5,2) {$2$};
	\draw [dotted] (0.66,0) -- (0.66,3);
	\draw [dotted] (0,2.25) -- (3,2.25);
	\draw [->] (-0.3,1) -- node [midway, above] {$x$} (-0.3,2.25);
	\node at (0.66,2.25) [circle, fill, inner sep=0.5mm] {};
\end{tikzpicture}
\end{align*}
\end{itemize}
These generators are submitted to the same associativity relations~\eqref{eq:guill2:listassoc} where finite boundaries should be replaced by infinite boundaries when needed and where base points shall be added.
\end{prop}

We first comment the notation with a scalar $x$ over an arrow present in the previous drawings such as:
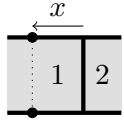
\[\begin{tikzpicture}[guillpart]
	\fill[guillfill] (3,0)--(0,0)--(0,2)--(3,2)--(3,0);
	\draw[guillsep] (3,0)--(0,0);
	\draw[guillsep] (0,2)--(3,2);
	\draw[guillsep] (2,0)--(2,2);
	\node at (1.33,1.) {$1$};
	\node at (2.5,1.) {$2$};	
	\node at (0.66,0) [circle, fill, inner sep=0.5mm] {};
	\node at (0.66,2) [circle, fill, inner sep=0.5mm] {};
	\draw [dotted] (0.66,0) -- (0.66,2);
		\draw [->] (2,2.3) -- node [midway, above] {$x$} (0.66,2.3);
\end{tikzpicture}\]
Such a drawing corresponds to the equivalence class of the pointed guillotine partition of the strip $\setP\times [v,v+q]$ with pointing $z\in\setP$ into two half-strips $(-\infty,u]\times [v,v+q]$ and $[u,+\infty)\times [v,v+q]$. Vertical translations change the value of $v$ only and thus the equivalence class is characterized only by the vertical coordinate $q$ (the shape); horizontal translations change simultaneously $z$ and $u$ and the equivalence class is thus characterized by the difference $x=z-u$. The previous drawing illustrates this fact: the sizes $\infty_{WE}$ and $q$ should be read on the pattern shape and the additional parameter $x$ is displayed by an arrow.

No additional coordinate $x$ appear in the drawing
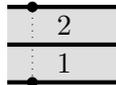
\[\begin{tikzpicture}[guillpart]
	\fill[guillfill] (3,0)--(0,0)--(0,2)--(3,2)--(3,0);
	\draw[guillsep] (3,0)--(0,0);
	\draw[guillsep] (0,2)--(3,2);
	\draw[guillsep] (0,1)--(3,1);
	\node at (1.5,0.5) {$1$};
	\node at (1.5,1.5) {$2$};	
	\node at (0.66,0) [circle, fill, inner sep=0.5mm] {};
	\node at (0.66,2) [circle, fill, inner sep=0.5mm] {};
	\draw [dotted] (0.66,0) -- (0.66,2);
\end{tikzpicture}
\]
because it corresponds to the pointed guillotine partition of $\setR\times [v,v+q]$ with base point $z$ into $\setP\times [v,v+q_1]$ and $\setP\times [v+q_1,v+q]$ with $0<q_1<q$. Vertical translations change $v$ but keep $q_1$ and $q$, which are the vertical parameters of the pattern shapes. Horizontal translations shift $x$ only and thus $x$ disappears in the equivalence classes.

\begin{proof}
The proof is the exact copy of the proof of Proposition~\ref{prop:guill2elemassoc} since the definition of guillotine partitions from guillotine cuts remains the same when the domain are infinite. The introduction of base points to break some of the translational invariance only adds new generators indexed by the distance between the base points and the cuts of partitions. The case study is long but easy.
\end{proof}

In practical computations, the pointings may play a subtle role which will briefly appear in Section~\ref{sec:canonicalboundarystructure} with shifts in the graduation of the spaces; it may however deserve a wider study which goes too far beyond the scope of the present paper and we present only the basic properties required to formulate our theorems.

\subsubsection{An example: west-east strips.}
As an example, we discuss in detail the case of West-East strips. We have in this case:
\begin{align*}
D(\patterntype{s}_{WE}^*)&=\{\patterntype{r},\patterntype{hs}_W,\patterntype{hs}_E,\patterntype{s}_{WE}^*\}
\\
\PatternShapes(\patterntype{s}_{WE}^*)&=(\setL^*\sqcup\{\infty_W,\infty_E,\infty_{WE}\})\times\setN_1
\end{align*}
An algebra $(\ca{A}_c)_{c\in\PatternShapes(\patterntype{s}_{WE})}$ over $\Guill_2^{(\patterntype{s}_{WE}^*)}$ contains an algebra $(\ca{A}_c)_{c\in\PatternShapes(\patterntype{r})}$ over $\Guill_2^{(\patterntype{r})}$ and two opposite West and East algebras $\ca{B}^W_r=\ca{A}_{\infty_W,r}$ and $\ca{B}^E_r=\ca{A}_{\infty_E,r}$ (with $r\in\setL^*$) , as well as a new class of objects $\ca{S}^{WE}_r=\ca{A}_{\infty_{WE},r}$, which we call a \emph{strip algebra}, due to the following property.

\begin{prop}\label{prop:stripalgebra}
$\ca{S}^{WE}_\bullet$ is an algebra over $\Guill_1$ and there exist morphisms of $\Guill_1$-algebras $(h_{r,x})_{r\in\setL^*,x\in\setP}$ such that $h_{r,x}:\ca{B}^W_r\otimes \ca{B}^E_r \to \ca{S}^{WE}_r$ (given explicitly in the proof below), which we call \emph{pointed pairings}, such that the following diagrams are commutative for any finite $p,q\in\setL^*$:
\newsavebox{\boxUp}
\savebox{\boxUp}{
	\begin{tikzpicture}[scale=0.5,baseline={(current bounding box.center)}]
	\fill[guillfill] (3,0)--(0,0)--(0,2)--(3,2)--(3,0);
	\draw[guillsep] (3,0)--(0,0)--(0,2)--(3,2);
	\draw[guillsep] (1,0)--(1,2);
	\node at (0.5,1) {$2$};
	\node at (2,1) {$3$};
	\end{tikzpicture}
}
\newsavebox{\boxRight}
\savebox{\boxRight}{
	\begin{tikzpicture}[scale=0.5,baseline={(current bounding box.center)}]
	\fill[guillfill] (3,0)--(0,0)--(0,2)--(3,2)--(3,0);
	\draw[guillsep] (3,0)--(0,0);
	\draw[guillsep] (0,2)--(3,2);
	\draw[guillsep] (2,0)--(2,2);
	\node at (1.33,1.) {$1$};
	\node at (2.5,1.) {$2$};
	
	\node at (0.66,0) [circle, fill, inner sep=0.5mm] {};
	\node at (0.66,2) [circle, fill, inner sep=0.5mm] {};
	\draw [dotted] (0.66,0) -- (0.66,2);
		\draw [->] (2,2.3) -- node [midway, above] {$x$} (0.66,2.3);
	\end{tikzpicture}
}
\newsavebox{\boxDown}
\savebox{\boxDown}{
	\begin{tikzpicture}[scale=0.5,baseline={(current bounding box.center)}]
	\fill[guillfill] (4,0)--(0,0)--(0,2)--(4,2)--(4,0);
	\draw[guillsep] (4,0)--(0,0);
	\draw[guillsep] (0,2)--(4,2);
	\draw[guillsep] (3,0)--(3,2);
	\node at (2.33,1) {$1$};
	\node at (3.5,1) {$2$};
	\node at (0.66,0) [circle, fill, inner sep=0.5mm] {};
	\node at (0.66,2) [circle, fill, inner sep=0.5mm] {};
	\draw [dotted] (0.66,0) -- (0.66,2);
	\draw [->] (3,2.3) -- node [midway, above] {$x+p$} (0.66,2.3);	
	\end{tikzpicture}
}
\newsavebox{\boxLeft}
\savebox{\boxLeft}{
	\begin{tikzpicture}[scale=0.5,baseline={(current bounding box.center)}]
	\fill[guillfill] (0,0)--(2,0)--(2,2)--(0,2)--(0,0);
	\draw[guillsep] (0,0)--(2,0)--(2,2)--(0,2);
	\draw[guillsep] (1,0)--(1,2);
	\node at (0.5,1) {$1$};
	\node at (1.5,1) {$2$};
	\end{tikzpicture}
}
\newsavebox{\boxDiag}
\savebox{\boxDiag}{
	\begin{tikzpicture}[scale=0.5,baseline={(current bounding box.center)}]
	\fill[guillfill] (0,0)--(5,0)--(5,2)--(0,2)--(0,0);
	\draw[guillsep] (0,0)--(5,0);
	\draw[guillsep] (0,2)--(5,2);
	\draw[guillsep] (2,0)--(2,2);
	\draw[guillsep] (3,0)--(3,2);
	\node at (1.33,1) {$1$};
	\node at (2.5,1) {$2$};
	\node at (4,1) {$3$};
	\node at (0.66,0) [circle, fill, inner sep=0.5mm] {};
	\node at (0.66,2) [circle, fill, inner sep=0.5mm] {};
	\draw [dotted] (0.66,0) -- (0.66,2);
	\draw [->] (2,2.3) -- node [midway, above] {$x$} (0.66,2.3);
	\end{tikzpicture}
}
\begin{equation}
\begin{tikzpicture}[baseline={(current bounding box.center)}]
\matrix (m) [matrix of math nodes,row sep=3em,column sep=4em,minimum width=2em]
  {
	\ca{B}^W_{q}\otimes \ca{A}_{p,q} \otimes \ca{B}^E_{q} 
     &  &
	 \ca{B}^W_{q} \otimes \ca{B}^E_{q}
     \\
    &  &  
     \\
	\ca{B}^W_{p} \otimes \ca{B}^E_{q}
     & &
     \ca{S}^{WE}      \\};
\path[-stealth]
    (m-1-1) edge node [left] {\usebox{\boxLeft}$\otimes\id_3$} (m-3-1)
    (m-1-1) edge node [above] {$\id_1\otimes$\usebox{\boxUp}} (m-1-3)
    (m-3-1) edge node [below] {\usebox{\boxDown}$=h_{q,x+p}$} (m-3-3)
    (m-1-3) edge node [right] {\usebox{\boxRight}$=h_{q,x}$} (m-3-3)
    (m-1-1) edge node [above right] {$\Xi$} (m-3-3);
\end{tikzpicture}
\label{eq:commutdiag:pairing}
\end{equation}
where $\Xi$ is the map associated to the pointed guillotine partition:
\begin{equation*}
\Xi = \usebox{\boxDiag}
\end{equation*}
\end{prop}
\begin{proof}
The $\Guill_1$-algebra is obtained by considering only horizontal guillotine cuts of pointed strips into pointed strips. Since a purely horizontal guillotine partition $(R_i)_{1\leq i\leq n}$ is given by $R_i=\setP\times [u_2,v_2]$, it is associated to the one-dimensional partition $([u_i,v_i])_{1\leq i\leq n}$. $\ca{S}^{WE}_\bullet=(\ca{S}^{WE}_r)_{r\in\setL^*}$ is thus an algebra over $\Guill_1$ with associative products given by, for any integer $m\geq 1$ and $r_i\in\setL^*$, by
\[
\begin{tikzpicture}[scale=0.5,baseline={(current bounding box.center)}]
\fill[guillfill] (3,0)--(0,0)--(0,5) -- (3,5) -- (3,0);
\draw[guillsep] (0,5)--(3,5);
\draw[guillsep](0,0)--(3,0);
\draw[guillsep] (0,1) -- (3,1);
\draw[guillsep] (0,1) -- (3,1);
\draw[guillsep] (0,2) -- (3,2);
\draw[guillsep] (0,4) -- (3,4);
\node at (2,0.5) {$1$};
\node at (2,1.5) {$2$};
\node at (2,3) {$\vdots$};
\node at (2,4.5) {$m$};
\node at (1,0) [circle, fill, inner sep=0.5mm] {};
\node at (1,5) [circle, fill, inner sep=0.5mm] {};
\draw [dotted] (1,0) -- (1,5);
\end{tikzpicture}
: \ca{S}^{WE}_{r_1}\otimesdots \ca{S}^{WE}_{r_m} \to \ca{S}^{WE}_{r_1+\ldots+r_m}
\]
We remark that, because of the quotient by translation equivalence, the position of the base line plays no role except the alignment of successive strips.

The relation with the two opposite algebras $\ca{B}^W_\bullet$ and $\ca{B}^E_\bullet$ is the existence of morphisms of algebras:
\savebox{\boxRight}{
	\begin{tikzpicture}[scale=0.5,baseline={(current bounding box.center)}]
	\fill[guillfill] (3,0)--(0,0)--(0,2)--(3,2)--(3,0);
	\draw[guillsep] (3,0)--(0,0);
	\draw[guillsep] (0,2)--(3,2);
	\draw[guillsep] (2,0)--(2,2);
	\node at (1.33,1.) {$1$};
	\node at (2.5,1.) {$2$};
	
	\node at (0.66,0) [circle, fill, inner sep=0.5mm] {};
	\node at (0.66,2) [circle, fill, inner sep=0.5mm] {};
	\draw [dotted] (0.66,0) -- (0.66,2);
		\draw [->] (2,2.3) -- node [midway, above] {$x$} (0.66,2.3);
	\end{tikzpicture}
}
\[
h_{q,x}=\usebox{\boxRight}: \ca{B}^{W}_{q}\otimes \ca{B}^{E}_{q} \to \ca{S}^{WE}_{q}
\]
The morphism property is obtained from the two possible decomposition (interchange relation again !) of the partition:
\[
\begin{tikzpicture}[guillpart]
\fill[guillfill] (0,0)--(5,0)--(5,2)--(0,2)--(0,0);
\draw[guillsep] (0,0)--(5,0);
\draw[guillsep](0,1)--(5,1);
\draw[guillsep](0,2)--(5,2);
\draw[guillsep] (3,0) -- (3,2);
\node at (2.,0.5) {$1$};
\node at (4,0.5) {$2$};
\node at (2.,1.5) {$3$};
\node at (4,1.5) {$4$};
\node at (1,0) [circle, fill, inner sep=0.5mm] {};
\node at (1,2) [circle, fill, inner sep=0.5mm] {};
\draw [dotted] (1,0) -- (1,2);
		\draw [->] (3,2.3) -- node [midway, above] {$x$} (1,2.3);
\end{tikzpicture}
: (\ca{B}^{W}_{r_1}\otimes \ca{B}^{E}_{r_1})\otimes (\ca{B}^{W}_{r_2}\otimes \ca{B}^{E}_{r_2}) \to \ca{S}^{WE}_{r_1+r_2}
\]
similarly to equations~\eqref{eq:interchangedecomp:ex1} and \eqref{eq:interchangedecomp:ex2}.
\end{proof}

We will not discuss in detail the case of half-planes and of the full plane since it works in the same way but we recollect in the following property all the important facts.
\begin{prop}\label{prop:halfplanemoduleactions}
An algebra $(\ca{A}_c)_{c\in\PatternShapes(\patterntype{hp}_N^*)}$ (resp. $\patterntype{hp}_S^*$, $\patterntype{hp}_W^*$ and $\patterntype{hp}_E^*$) over $\Guill_2^{(\patterntype{hp}_N^*)}$ (resp. $\patterntype{hp}_S^*$, $\patterntype{hp}_W^*$ and $\patterntype{hp}_E^*$) is specified by
\begin{itemize}
\item a $\Guill_2^{(\patterntype{r})}$-algebra $(\ca{A}_c)_{c\in\PatternShapes(\patterntype{r})}$,
\item three boundary algebras $(\ca{B}^{i}_r)_{r\in\setL^*}$ with $i\in\{W,N,E\}$ (resp. $\{W,S,E\}$, $\{S,W,N\}$ and $\{S,E,N\}$),
\item a strip algebra $(\ca{B}^{WE}_q)_{q\in\setL^*}$ with a pairing from $\ca{B}^W_\bullet$ and $\ca{B}^E_\bullet$ to it,
\item two corner double modules $\ca{M}^{j}$ with $j\in\{NW,NE\}$ (resp. $\{SW,SE\}$, $\{SW,NW\}$ and $\{SE,NE\}$) with suitable actions of the boundary algebras,
\item an additional space $\ca{H}^{N}$ (resp. $\ca{H}^S$, $\ca{H}^W$ and $\ca{H}^E$), called a \emph{half-plane module}, such that:
\begin{enumerate}[(i)]
\item there exist pointed pairings
\[
\begin{tikzpicture}[guillpart]
\fill[guillfill] (0,0)--(3,0)--(3,1)--(0,1)--(0,0);
\draw[guillsep] (0,0)--(3,0);
\draw[guillsep] (2,0) -- (2,1);
\node at (0.66,0.5) {$1$};
\node at (2.5,0.5) {$2$};
	\node (A1) at (1.33,0) [circle, fill, inner sep=0.5mm] {};
	\node (A2) at (1.33,1) [circle, fill, inner sep=0.5mm] {};
	\draw [dotted] (A1) -- (A2);
		\draw [->] (2,1.3) -- node [midway, above] {$x$} (1.33,1.3);
\end{tikzpicture}
: \ca{M}^{NW} \otimes \ca{M}^{NE} \to \ca{H}^{N}
\]
satisfying a commutative diagram similar to~\eqref{eq:commutdiag:pairing} (resp. similar pairings).
\item $\ca{S}^{WE}_\bullet$ and $\ca{H}^N$ form a algebra over $\Guill_{1,R}$ (i.e. $\ca{H}^N$ is a left module over $\ca{S}^{WE}_\bullet$) such that the following diagram 
is commutative:
\newsavebox{\boxUpA}
\savebox{\boxUpA}{
	\begin{tikzpicture}[guillpart]
	\fill [guillfill] (0,0)--(2,0)--(2,1)--(0,1)--(0,0);
	\draw [guillsep] (0,0)--(2,0);
	\draw [guillsep] (0,1)--(2,1);
	\draw [guillsep] (1.33,0) -- (1.33,1);
	\node at (0.75,0.5) {$1$};
	\node at (1.66,0.5) {$2$};
	\node (A1) at (0.33,0) [circle, fill, inner sep=0.5mm] {};
	\node (A2) at (0.33,1) [circle, fill, inner sep=0.5mm] {};
	\draw [dotted] (A1) -- (A2);
	\end{tikzpicture}    
}
\newsavebox{\boxUpB}
\savebox{\boxUpB}{
	\begin{tikzpicture}[guillpart]
	\fill [guillfill] (0,0)--(2,0)--(2,1)--(0,1)--(0,0);
	\draw [guillsep] (0,0)--(2,0);
	\draw [guillsep] (1.33,0) -- (1.33,1);
	\node at (0.75,0.5) {$3$};
	\node at (1.66,0.5) {$4$};
	\node (A1) at (0.33,0) [circle, fill, inner sep=0.5mm] {};
	\coordinate (A2) at (0.33,1) ;
	\draw [dotted] (A1) -- (A2);
	\end{tikzpicture}    
}
\newsavebox{\boxLeftA}
\savebox{\boxLeftA}{
	\begin{tikzpicture}[scale=0.5,baseline={(current bounding box.center)}]
	\fill [guillfill] (0,0)--(1,0)--(1,2)--(0,2)--(0,0);
	\draw [guillsep] (1,0)--(1,2);
	\draw [guillsep] (0,1) -- (1,1);
	\draw [guillsep] (0,0) -- (1,0);
	\node at (0.5,0.5) {$1$};
	\node at (0.5,1.5) {$3$};
	\end{tikzpicture}    
}
\newsavebox{\boxLeftB}
\savebox{\boxLeftB}{
	\begin{tikzpicture}[scale=0.55,baseline={(current bounding box.center)}]
		\fill [guillfill] (0,0)--(1,0)--(1,2)--(0,2)--(0,0);
	\draw [guillsep] (0,0)--(0,2);
	\draw [guillsep] (0,1) -- (1,1);
	\draw [guillsep] (0,0) -- (1,0);
	\node at (0.5,0.5) {$2$};
	\node at (0.5,1.5) {$4$};
	\end{tikzpicture}    
}
%\newsavebox{\boxRight}
\savebox{\boxRight}{
	\begin{tikzpicture}[scale=0.5,baseline={(current bounding box.center)}]
	\fill [guillfill] (0,0)--(1.33,0)--(1.33,2)--(0,2)--(0,0);
	\draw [guillsep] (0,1) -- (1.33,1);
	\draw [guillsep] (0,0) -- (1.33,0);
	\node at (0.75,0.5) {$1$};
	\node at (0.75,1.5) {$3$};
	\node (A1) at (0.33,0) [circle, fill, inner sep=0.5mm] {};
	\coordinate (A2) at (0.33,2) ;
	\draw [dotted] (A1) -- (A2);

	\end{tikzpicture}    
}
%\newsavebox{\boxDown}
\savebox{\boxDown}{
	\begin{tikzpicture}[scale=0.55,baseline={(current bounding box.center)}]
		\fill [guillfill] (0,0)--(2,0)--(2,1)--(0,1)--(0,0);
	\draw [guillsep] (0,0)--(2,0);
	\draw [guillsep] (1.33,0) -- (1.33,1);
	\node at (0.75,0.5) {$1$};
	\node at (1.66,0.5) {$2$};
	\node (A1) at (0.33,0) [circle, fill, inner sep=0.5mm] {};
	\coordinate (A2) at (0.33,1) ;
	\draw [dotted] (A1) -- (A2);
	\end{tikzpicture}    
}
%\newsavebox{\boxDiag}
\savebox{\boxDiag}{
	\begin{tikzpicture}[scale=0.5,baseline={(current bounding box.center)}]
	\fill[guillfill] (0,0)--(2,0)--(2,2)--(0,2)--(0,0);
	\draw[guillsep] (0,0)--(2,0);
	\draw[guillsep] (0,1)--(2,1);
	\draw[guillsep] (1.33,0)--(1.33,2);
	\node at (0.75,0.5) {$1$};
	\node at (1.75,0.5) {$2$};
	\node at (0.75,1.5) {$3$};
	\node at (1.75,1.5) {$4$};
	\node (A1) at (0.33,0) [circle, fill, inner sep=0.5mm] {};
	\coordinate (A2) at (0.33,2) ;
	\draw [dotted] (A1) -- (A2);
	\end{tikzpicture}    
}
\begin{equation}
\begin{tikzpicture}[baseline={(current bounding box.center)}]
\matrix (m) [matrix of math nodes,row sep=3em,column sep=4em,minimum width=2em]
  {
	\ca{B}^W_{q}\otimes \ca{B}^E_{q}\otimes \ca{M}^{NW} \otimes \ca{M}^{NE}
     &  &
	\ca{S}^{WE}\otimes \ca{H}^N
     \\
	\ca{M}^{NW}\otimes \ca{M}^{NE}     
     & &
     \ca{H}^{N}      \\};
\path[-stealth]
    (m-1-1) edge node [left] {\usebox{\boxLeftA}$\otimes$\usebox{\boxLeftB}} (m-2-1)
    (m-1-1) edge node [above] {\usebox{\boxUpA}$\otimes$\usebox{\boxUpB}} (m-1-3)
    (m-2-1) edge node [below] {\usebox{\boxDown}} (m-2-3)
    (m-1-3) edge node [right] {\usebox{\boxRight}} (m-2-3)
    (m-1-1) edge node [above right] {$\Upsilon$} (m-2-3);
\end{tikzpicture}
\end{equation}
where the diagonal map is the pointed guillotine partition map
\begin{equation*}
\Upsilon= \usebox{\boxDiag}
\end{equation*}
(resp. similar $\Guill_{1,R}$- or $\Guill_{1,L}$-algebra structure with similar diagrams with rotated guillotine partitions).
\end{enumerate}
\end{itemize}
\end{prop}
The full proof is quite lengthy and uses the same ingredients as before, i.e. purely vertical or horizontal guillotine cuts for associativity and actions diagrams and square partitions similar to $\Phi$, $\Psi$ and $\Upsilon$ for commutative diagrams.

From now on, we also use for the pointed pattern type $\patterntype{x}^*$ of the previous proposition the same short notation $\ca{A}_{\PatternShapes(\patterntype{x}^*)}$ for the collection $(\ca{A}_{p,q})_{(p,q)\in\PatternShapes(\patterntype{x}^*)}$ as already introduced in~\eqref{eq:shortnotation:guillalgebra}.

Finally, we arrive at least at the full plane case which contains, without any surprise, objects associated to all the objects of figure~\ref{fig:admissiblepatterns}: a $\Guill_2^{(\patterntype{r})}$-algebra for rectangles, four boundary algebras, four corner double modules, two strip algebras, four half-plane modules and a last space associated to the pointed full plane obtained from the two pointed pairings from two opposite half-plane modules with a suitable compatibility condition similar to~\eqref{eq:commutdiag:pairing}.

\subsubsection{Change of base points.} In many situations, in particular the one of the next Section~\ref{sec:canonicalboundarystructure}, the pairings $h_{q,x}$ of eq.~\eqref{eq:commutdiag:pairing} and their analogues for all the doubly-infinite shapes which require base points can be described as a composition $\theta_{q,x}\circ h_{q,0}$ where:
\begin{itemize}
\item  the objects $h_{q,0}$ correspond to single-cut generators $m_{WE}^{r,r'|q}$ and $m_{SN}^{p|s,s'}$ with a base point placed exactly on the cut.
\item the translator operators $\theta_{q,x}$ act as a representation of the translation group on the spaces of the algebra with $\theta_{q,x} = \theta_{q,1}^{\circ x}$ (resp. $\theta_{q,-1}^{\circ -x}$) for positive (resp. negative) $x$. 
\end{itemize}
In the simple purely rectangular geometries considered below, these distinctions play a minor role but we expect them to be more important for generalized geometries.

\subsubsection{Summary.} Spaces of an algebra over a fully extended operad $\Guill_2^{(u)}$ are associated to any rectangle with the shape
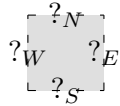
\begin{equation*}
\begin{tikzpicture}[scale=0.5,baseline={(current bounding box.center)}]
\fill[guillfill] (0,0) rectangle (2,2);
\draw[loosely dashed] (0,0)--(2,0);
\draw[loosely dashed] (0,0)--(0,2);
\draw[loosely dashed] (2,0)--(2,2);
\draw[loosely dashed] (0,2)--(2,2);
\node at (1,0) {$?_S$};
\node at (1,2) {$?_N$};
\node at (0,1) {$?_W$};
\node at (2,1) {$?_E$};
\end{tikzpicture}
\end{equation*}
where each of the four boundaries are either solid (finite bound) or empty (infinite bound). Products or actions are associated to simple guillotine partitions (where, for each $a\in\{N,S,W,E\}$, every $?_a$ is replaced by the same type of boundary):
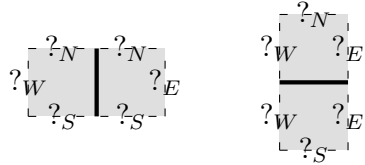
\begin{equation*}
\begin{tikzpicture}[scale=0.45,baseline={(current bounding box.center)}]
\fill[guillfill] (0,0) rectangle (4,2);
\draw[loosely dashed] (0,0)--(4,0);
\draw[loosely dashed] (0,0)--(0,2);
\draw[loosely dashed] (4,0)--(4,2);
\draw[loosely dashed] (0,2)--(4,2);
\node at (1,0) {$?_S$};
\node at (3,0) {$?_S$};
\node at (1,2) {$?_N$};
\node at (3,2) {$?_N$};
\node at (0,1) {$?_W$};
\node at (4,1) {$?_E$};
\draw[guillsep] (2,0)--(2,2);
\end{tikzpicture}
\qquad 
\begin{tikzpicture}[scale=0.45,baseline={(current bounding box.center)}]
\fill[guillfill] (0,0) rectangle (2,4);
\draw[loosely dashed] (0,0)--(2,0);
\draw[loosely dashed] (0,0)--(0,4);
\draw[loosely dashed] (2,0)--(2,4);
\draw[loosely dashed] (0,4)--(2,4);
\node at (1,0) {$?_S$};
\node at (1,4) {$?_N$};
\node at (0,1) {$?_W$};
\node at (0,3) {$?_W$};
\node at (2,1) {$?_E$};
\node at (2,3) {$?_E$};
\draw[guillsep] (0,2)--(2,2);
\end{tikzpicture}
\end{equation*}
The associativity of each product, the compatibility between a product and an action and the compatibility between actions and a pairing correspond all to the commutative diagrams obtained from the two possible decompositions of the simple horizontal and vertical guillotine partitions:
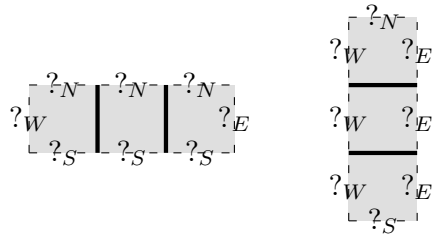
\begin{equation*}
\begin{tikzpicture}[scale=0.45,baseline={(current bounding box.center)}]
\fill[guillfill] (0,0) rectangle (6,2);
\draw[loosely dashed] (0,0)--(6,0);
\draw[loosely dashed] (0,0)--(0,2);
\draw[loosely dashed] (6,0)--(6,2);
\draw[loosely dashed] (0,2)--(6,2);
\draw[guillsep] (2,0)--(2,2);
\draw[guillsep] (4,0)--(4,2);
\node at (1,0) {$?_S$};
\node at (3,0) {$?_S$};
\node at (5,0) {$?_S$};
\node at (1,2) {$?_N$};
\node at (3,2) {$?_N$};
\node at (5,2) {$?_N$};
\node at (0,1) {$?_W$};
\node at (6,1) {$?_E$};
\end{tikzpicture}
\qquad 
\begin{tikzpicture}[scale=0.45,baseline={(current bounding box.center)}]
\fill[guillfill] (0,0) rectangle (2,6);
\draw[loosely dashed] (0,0)--(2,0);
\draw[loosely dashed] (0,0)--(0,6);
\draw[loosely dashed] (2,0)--(2,6);
\draw[loosely dashed] (0,6)--(2,6);
\draw[guillsep] (0,2)--(2,2);
\draw[guillsep] (0,4)--(2,4);
\node at (1,0) {$?_S$};
\node at (1,6) {$?_N$};
\node at (0,1) {$?_W$};
\node at (0,3) {$?_W$};
\node at (0,5) {$?_W$};
\node at (2,1) {$?_E$};
\node at (2,3) {$?_E$};
\node at (2,5) {$?_E$};
\end{tikzpicture}
\end{equation*}
Compatibilities between transverse products, actions and pairings all correspond to the two possible decompositions of the guillotine partitions:
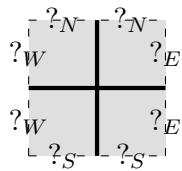
\begin{equation*}
\begin{tikzpicture}[scale=0.45,baseline={(current bounding box.center)}]
\fill[guillfill] (0,0) rectangle (4,4);
\draw[loosely dashed] (0,0)--(4,0);
\draw[loosely dashed] (0,0)--(0,4);
\draw[loosely dashed] (4,0)--(4,4);
\draw[loosely dashed] (0,4)--(4,4);
\draw[guillsep] (2,0)--(2,4);
\draw[guillsep] (0,2)--(4,2);
\node at (1,0) {$?_S$};
\node at (3,0) {$?_S$};
\node at (1,4) {$?_N$};
\node at (3,4) {$?_N$};
\node at (0,1) {$?_W$};
\node at (0,3) {$?_W$};
\node at (4,1) {$?_E$};
\node at (4,3) {$?_E$};
\end{tikzpicture}
\end{equation*}

In the case of doubly infinite domains, base points must be added to these diagrams to complete the set of generators and associativity conditions.

		\section[Canonical boundary structure for discrete-space Markov processes]{Canonical boundary structure for two-dimensional discrete-space Markov processes}\label{sec:canonicalboundarystructure}

Theorems~\ref{theo:canonicalexampleGuill} and \ref{theo:partitionfuncguillotop} (and their counterpart~\ref{theo:canonicalexampleGuill:continuous} below for continuous spaces) introduce a $\Guill_2^{(\patterntype{r})}$-algebra related to two-dimensional Markov processes. As seen in Section~\ref{sec:proba}, boundary conditions play an essential role in probabilistic applications to two-dimensional Markov chain. The present section is devoted to the introduction of canonical boundary algebras, corner double modules, strip algebras, half-plane modules and full plane scalar field, thus extending Theorem~\ref{theo:canonicalexampleGuill} to the extended operad $\Guill_2^{(\patterntype{fp}^*)}$.

To this purpose, all the constructions will rely on the construction of direct sums/disjoint unions of elementary spaces related to the spaces $T_{p,q}(V(S_1),V(S_2))$ where some of spaces of endomorphisms are replaced by vector spaces and their duals so that they may admit actions of some of the tensorial factors of $T_{p,q}(V(S_1),V(S_2))$.

\begin{rema}
	We warn the reader that this canonical boundary structure is \emph{not} suitable to describe \emph{any} type of boundary conditions (not even some relevant ones) but only provides a nice framework for proofs such as the the proof of Theorem~\ref{theo:stability} and the description of linear maps in Theorem~\ref{theo:reductionofmorphisms:canostruct}. Suitable spaces for interesting non-trivial boundary conditions are provided in sections~\ref{sec:boundaryalgebra} and \ref{sec:applications}.
\end{rema}

\subsection{Disjoint union (sum) of spaces}
\subsubsection{Basic properties}
We will consider the following sets associated to the various geometric shapes:
\begin{itemize}
	\item $(\setL^*,\leq)$ for half-strips,
	\item $(\setL^*\times \setL^*,\preccurlyeq)$ for corners, strips and half-planes where the partial order is given by $(a,b)\preccurlyeq(a',b')$ if and only if $a\leq a'$ and $b\leq b'$.
\end{itemize}
On the corners, each factor $\setL^*$ corresponds to one of the directions of the two sides of the corners. On half-strips and half-planes, each factor $\setL^*$ corresponds to one of the two directions on the axis where the shape admits full translation invariance.

Let $I$ be a set and let $(E_i)_{i\in I}$ a collection of sets indexed by $I$. The disjoint union 
\[
\bigsqcup_{i\in I} E_i
\]
is defined as the universal set $X$ (up to isomorphism) endowed with canonical injection $\iota_i : E_i \to X$, such that, for any space $Y$ and any collection of morphisms $(\phi_i)_{i\in I}$ with $\phi_i : E_i\to Y$, there exists a morphism $\Phi: X\to Y$ such that $\phi_i = \Phi\circ \iota_i$. It may be concretely realized as the set $\cup_{i\in I} \{i\}\times E_i$ with $\iota_i(x) = (i,x)$.

First, the disjoint union of vector spaces or algebras $(E_i)_{i\in I}$ are again vector spaces and algebras. In this case, it is traditionally written
\[
\bigoplus_{i\in I} E_i
\]
and an element of it is a \emph{finite} linear combination of elements of the spaces $E_i$. We will write it $a=(a_i)_{i\in I}$ with $a_i=0_{E_i}$ except for a finite number of them. In case of algebras, we have a product given by the component-wise products $a\cdot b= (a_i b_i)_{i\in I}$. In case of modules $(F_i)_{i\in I}$ over algebras $(E_i)_{i\in I}$, we have in the same way an action of $\oplus_{i} E_i$ on $\oplus_i F_i$ defined component-wise.

A categorical point of view is to consider the previous construction as a colimit with a functor $I\to C$ where $C$ is a category (sets, vector spaces or algebras for us) and $I$ is endowed with the pre-order given by equality (the only arrows are the identity arrows). This remark may look abstract nonsense at this point but it induces an interesting structure when combined with Section~\ref{sec:invariantboundaryelmts}. Indeed, in all the cases below, the index sets will be $\setL^*$ and $\setL^*\times \setL^*$, which may equipped with the total order $\leq$ (inherited $\setN$ or $\setR_+$) for the first and with the order $\preccurlyeq$ for the second, with $(p,q)\preccurlyeq (p',q')$ if and only if $p\leq p'$ and $q\leq q'$. The disjoint unions below are considered as limits or colimits with respect to $\setL^*$ or $\setL^*\otimes\setL^*$ \emph{without} any consideration of these partial orders. Taking them into account may produce instead projective or direct limits: this will be the case in particular in Section~\ref{sec:invariantboundaryelmts} through the maps $\psi_p^{\MarkovWeight{F},S,q,r}$ (see \eqref{eq:def:actionmorphism:psi}) and $\Phi_p^{S,q,r}$ (see Definition~\ref{def:eigenalgebrauptomorphims}) and it leads to further interesting algebraic and probabilistic interpretations discussed below.

\subsubsection{A useful notation.}

In most cases below, we will consider the case $I=\setN$ or $\setR_+$ (i.e. $\setL$ in our notation) and spaces $E_i$ with a bilinear map $E_i\times E_j\to E_{i+j}$. For finite state space and discrete space, this will correspond to the canonical map
$(V\otimes W^{\otimes i})\otimes (V'\otimes W^{\otimes j}) \to V\otimes V'\otimes W^{\otimes i+j}$
where $W$ is a vector space or an algebra of endomorphisms of a vector space. In this case, these additional maps define two actions on the direct sums over $\setL$, one on the left and on the right. In all the spaces below, all the direct sums will be equipped with only one of them and we decide to incorporate this information on the notations for easier reading by using a $(L)$ or $(R)$ as exponent on the direct sum symbol.

We have canonical maps, for each $r\in\setL$,
\begin{align*}
	\left( {\bigoplus_{p\in\setL}}^{(L)}  V\otimes W^{\otimes p} \right) \otimes (V'\otimes W^{\otimes r}) & \mapsto {\bigoplus_{p\in\setL}}^{(L)} V\otimes V'\otimes W^{\otimes p}
	\\
	(v\otimes \iota_p(w) \otimes (v'\otimes \iota_r(w') &\mapsto v\otimes v' \otimes \iota_{p+q}(w\otimes w')
	\\
	(V'\otimes W^{\otimes r})\otimes \left( {\bigoplus_{p\in\setL}}^{(L)}  V\otimes W^{\otimes p} \right) & \mapsto {\bigoplus_{p\in\setL}}^{(R)} V'\otimes V\otimes W^{\otimes p}
	\\
	(v'\otimes \iota_r(w') \otimes (v\otimes \iota_p(w) &\mapsto v'\otimes v \otimes \iota_{p+q}(w'^\otimes w)
\end{align*}
where $w \in W^{\otimes p}$, $w'\in W^{\otimes r'}$ and thus $w\otimes w' \in W^{p+r}$.

\subsection{The new coloured spaces on boundaries, part one (without base points).}

Theorem~\ref{theo:canonicalexampleGuill} provides a $\Guill_2^{(\patterntype{r})}$-structure for partition functions of 2D Markov processes. As seen in Section~\ref{sec:proba}, the complete laws of the process also involve boundary weights, which become a key object for the description of Gibbs measures on the whole space. The aim of the next sections is to fit these boundary weights into the framework of extended guillotine operads by suitable local-to-global constructions. To this purpose, we first present a canonical construction of a $\Guill_2^{(\patterntype{fp}^*)}$-algebra, which extends the $\Guill_2^{(\patterntype{r})}$-algebra of Theorem~\ref{theo:canonicalexampleGuill} and provides the most elementary structure to formulate Theorem~\ref{theo:stability}. This toolbox allows one to implement actions on boundaries on a computer through traditional linear algebra and illustrate some interesting phenomena typical of the dimension two.

We provide here only a description for the discrete space $(\setP,\setL)=(\setZ,\setN)$, which highlights best the algebraic nature of the various objects and compare them to other standard constructions. The generalization to continuous space is rejected to Section~\ref{sec:generalizations} since it requires additional measure-theoretical tools.

\subsubsection{Spaces for half-strips and corners from direct limits.}

All through the constructions below, we systematically assume that $E^{\otimes 0}=\setK$ for any vector space or algebra $E$ on a field $\setK$. For a discrete state space and discrete space, we introduce the following spaces for any $p,q\geq 1$:
\begingroup
\allowdisplaybreaks
\begin{subequations}
\label{eq:canonicalboundaryalg:FD}
\begin{align}
T_{p,\infty_S}(V(S_1),V(S_2)) &= (V(S_1)^*)^{\otimes p} \otimes  {\bigoplus_{q\in\setN}}^{(L)} \End(V(S_2))^{\otimes q}
\\
T_{p,\infty_N}(V(S_1),V(S_2)) &= V(S_1)^{\otimes p} \otimes  {\bigoplus_{q\in\setN}}^{(R)} \End(V(S_2))^{\otimes q} 
\\
T_{\infty_W,q}(V(S_1),V(S_2)) &= {\bigoplus_{p\in\setN}}^{(L)} \End(V(S_1))^{\otimes p} \otimes  (V(S_2)^*)^{\otimes q}
\\
T_{\infty_E,q}(V(S_1),V(S_2)) &=  {\bigoplus_{p\in\setN}}^{(R)} \End(V(S_1))^{\otimes p} \otimes V(S_2)^{\otimes q} 
\end{align}
as extended versions of the previous spaces $T_{p,q}(V(S_1),V(S_2))$ on the half-strips. For degenerate half-strips, we use instead
\begin{align}
T_{0,\infty_S}(V(S_1),V(S_2)) &= \setK \otimes  {\bigoplus_{q\in\setN}}^{(L)} \Diag(V(S_2))^{\otimes q}
\end{align}
and equivalent formulae in the three other directions. For the corners, we introduce the spaces:
\begin{align} 
T_{\infty_W,\infty_S}(V(S_1),V(S_2)) &= 
\left( {\bigoplus_{p\in\setN}}^{(L)} (V(S_1)^*)^{\otimes p} \right) \otimes
\left( {\bigoplus_{q\in\setN}}^{(L)} (V(S_2)^*)^{\otimes q}\right)
\\
T_{\infty_W,\infty_N}(V(S_1),V(S_2)) &= 
\left( {\bigoplus_{p\in\setN}}^{(L)} V(S_1)^{\otimes p} \right) \otimes 
\left( {\bigoplus_{q\in\setN}}^{(R)} (V(S_2)^*)^{\otimes q}\right)
\\
T_{\infty_E,\infty_S}(V(S_1),V(S_2)) &=
\left( {\bigoplus_{p\in\setN}}^{(R)} (V(S_1)^*)^{\otimes p} \right) \otimes 
\left( {\bigoplus_{q\in\setN}}^{(L)} V(S_2)^{\otimes q}\right)
\\
T_{\infty_E,\infty_N}(V(S_1),V(S_2)) &= 
\left( {\bigoplus_{p\in\setN}}^{(R)} V(S_1)^{\otimes p} \right) \otimes 
\left( {\bigoplus_{q\in\setN}}^{(R)} V(S_2)^{\otimes q}\right)  
\end{align}
\end{subequations}
For corners, there are two direct sums over $\setL^*$: depending on the situation, it may be advantageous to describe them as a single direct sums over $\setL^*\times\setL^*$ following a trivial identification.

The left and right direct sums as well as the choice of the spaces among $V$, $V^*$ and $\End(V)$ can be read directly from the pattern shapes as follows. We first canonically identify $\End(V)$ with $V^*\times V$ and interpret $V^*$ as associated to a bottom or a left side of a rectangle. Thus, for finite $p$ and $q$, $T_{p,q}(V(S_1),V(S_2))$ is identified to
\[
(V(S_1)^*)^{\otimes p}\otimes V(S_1)^{\otimes p} \otimes (V(S_2)^*)^{\otimes q} \otimes V(S_2)^{\otimes q}
\]
where the factors correspond to South, North, West and East sides respectively.

Sending a size of the pattern shape to infinity (left or right) in one of the four directions has the following consequences on the shapes and the definition of the spaces:
\begin{itemize}
	\item an edge, perpendicular to the infinite direction, is moved towards infinity and disappears: the corresponding factor $V^{\otimes p}$ in the previous equation is removed.
	\item edge(s) parallel to the infinite direction, grow(s) to infinity: the tensor product(s) w.r.t. this size is replaced by a direct sums indexed by all the sizes, i.e. $\setL$.
\end{itemize}
For segment state spaces, the space $S_r^{(i)}$ that corresponds to the disappearing edge is removed and the the factor(s) that correspond(s) to the infinite length edges becomes part of a direct sum over all the possible finite size.

\subsubsection{Products between half-strips and rectangles from morphisms between direct structures.}

All the products between half-strips and corners spaces as well as actions of rectangles involve products defined on direct sums. 

\begin{lemm}[products for half-strips]
Let $V$ and $W$ be two vectors spaces. For any $p_1,p_2,q\in \setL^*$, let $M_{p_1,p_2,q}$ be the product defined by:
\begin{align*}
M_{p_1,p_2,q}:	\left(V^{\otimes p_1}\otimes \End(W)^{\otimes q}\right)\otimes \left(V^{\otimes p_2}\otimes \End(W)^{\otimes q}\right)
& \to V^{\otimes p_1+p_2}\otimes \End(W)^{\otimes q} 
\\
\left( u \otimes A\right)\otimes \left(v \otimes B\right) & \mapsto (u\otimes v) \otimes AB
\end{align*}
This collection of products defines canonically a product over the direct sums
\[
\begin{split}
M_{p_1,p_2}: &	\left(V^{\otimes p_1}\otimes {\bigoplus_{q\in\setL}}^{(L|R)} \End(W)^{\otimes q}\right)\otimes \left(V^{\otimes p_2}\otimes {\bigoplus_{q\in\setL}}^{(L|R)} \End(W)^{\otimes q}\right) 
\\
& \to V^{\otimes p_1+p_2}\otimes {\bigoplus_{q\in\setL}}^{(L|R)} \End(W)^{\otimes q}
\end{split} 
\]
\end{lemm}
By considering suitable spaces $V$ and $W$ among $V(S_1)$, $V(S_2)$ and their duals in each case, we use this lemma to define on the four sides the following elementary 
guillotine products;
\begingroup
\allowdisplaybreaks
\begin{subequations}
\label{eq:productboundaryalgebra}
\begin{align*}
%		m_{WE}^{p_1,p_2|\infty_S}=
	\begin{tikzpicture}[guillpart]
			\fill[guillfill] (0,0) rectangle (2,1);
			\draw[guillsep] (0,0)--(0,1)--(2,1)--(2,0)  (1,0)--(1,1)  ;
			\node at (0.5,0.5) {$1$};
			\node at (1.5,0.5) {$2$};
	\end{tikzpicture}
		: T_{p_1,\infty_S}(V(S_1),V(S_2)) \otimes T_{p_2,\infty_S}(V(S_1),V(S_2))  &\to T_{p_1+p_2,\infty_S}(V(S_1),V(S_2)) 
		\\
	\begin{tikzpicture}[guillpart]
		\fill[guillfill] (0,0) rectangle (2,1);
			\draw[guillsep] (0,1)--(0,0)--(2,0) --(2,1)  (1,0)--(1,1)  ;
			\node at (0.5,0.5) {$1$};
			\node at (1.5,0.5) {$2$};
	\end{tikzpicture}
		: T_{p_1,\infty_N}(V(S_1),V(S_2)) \otimes T_{p_2,\infty_N}(V(S_1),V(S_2))  &\to T_{p_1+p_2,\infty_N}(V(S_1),V(S_2)) 
		\\
	\begin{tikzpicture}[guillpart,rotate=-90]
			\fill[guillfill] (0,0) rectangle (2,1);
			\draw[guillsep] (0,0)--(0,1)  (1,0)--(1,1)  (2,0)--(2,1)  (0,1)--(2,1);
			\node at (0.5,0.5) {$2$};
			\node at (1.5,0.5) {$1$};
	\end{tikzpicture}
		: T_{\infty_W,q_1}(V(S_1),V(S_2)) \otimes T_{\infty_W,q_2}(V(S_1),V(S_2))  &\to T_{\infty_W,q_1+q_2}(V(S_1),V(S_2)) 
		\\
		\begin{tikzpicture}[guillpart,rotate=90]
			\fill[guillfill] (0,0) rectangle (2,1);
			\draw[guillsep] (0,0)--(0,1)  (1,0)--(1,1)  (2,0)--(2,1)  (0,1)--(2,1);
			\node at (0.5,0.5) {$1$};
			\node at (1.5,0.5) {$2$};
		\end{tikzpicture}
	: T_{\infty_E,q_1}(V(S_1),V(S_2)) \otimes T_{\infty_E,q_2}(V(S_1),V(S_2))  &\to T_{\infty_E,q_1+q_2}(V(S_1),V(S_2)) 
\end{align*}
\end{subequations}
\endgroup
It is then a simple exercise to check that each of these products satisfies the $\Guill_1$-associativity inherited from the associativity of the tensor product and of the composition of endomorphisms.

A similar construction can be performed for the gluing of a rectangle on a half-strip, excepted that in this case the tensor product is between endomorphisms and the composition of endomorphisms is replaced by the action of endomorphisms on a vector space.

\begin{lemm}[actions on rectangles on half-strips]
	Let $V$ and $W$ be two vectors spaces. Let $M_{p,q,r}^R$ and $M_{p,q,r}^L$ be the products defined by:
	\begin{align*}
		M^R_{p,q,r}:	\left(\End(V)^{\otimes p}\otimes \End(W)^{\otimes q}\right)\otimes \left(V^{\otimes p}\otimes \End(W)^{\otimes r}\right)
		& \to V^{\otimes p}\otimes \End(W)^{\otimes q+r} 
		\\
		\left( A \otimes B\right)\otimes \left(v \otimes C\right) & \mapsto (Av) \otimes (B\otimes C)
		\\
		M^L_{p,q,r}:	\left(\End(V)^{\otimes p}\otimes \End(W)^{\otimes q}\right)\otimes \left({V^*}^{\otimes p}\otimes \End(W)^{\otimes r}\right)
		& \to {V^*}^{\otimes p}\otimes \End(W)^{\otimes q+r} 
		\\
		\left( A \otimes B\right)\otimes \left(v \otimes C\right) & \mapsto (vA) \otimes (C\otimes B)
	\end{align*}
	This collection of products define the following two products over the left- and right-oriented direct sums 
	\begin{align*}
	M_{p,q}^R:	\left(\End(V)^{\otimes p}\otimes \End(W)^{\otimes q}\right)\otimes \left(V^{\otimes p}\otimes {\bigoplus_{r\in\setL}}^{(R)} \End(W)^{r}\right)
	& \to V^{\otimes p}\otimes {\bigoplus_{r\in\setL}}^{(R)} \End(W)^{\otimes r} 
		\\
	M_{p,q}^L:	\left(\End(V)^{\otimes p}\otimes \End(W)^{\otimes q}\right)\otimes \left({V^*}^{\otimes p}\otimes {\bigoplus_{r\in\setL}}^{(L)} \End(W)^{\otimes r}\right)
	& \to {V^*}^{p}\otimes {\bigoplus_{r\in\setL}}^{(L)} \End(W)^{\otimes r}
	\end{align*}
\end{lemm}
From this lemma applied to suitable spaces $V$ and $W$ among $V(S_1)$, $V(S_2)$ and their duals, we directly obtain the following four new actions:
\begingroup
\allowdisplaybreaks
\begin{subequations}\label{eq:actionboundaryalgebra}
	\begin{align}
		\begin{tikzpicture}[guillpart]
			\fill[guillfill] (0,0) rectangle (1,2);
			\draw[guillsep] (0,0)--(0,2)--(1,2)--(1,0)  (0,1)--(1,1)  ;
			\node at (0.5,0.5) {$2$};
			\node at (0.5,1.5) {$1$};
		\end{tikzpicture}
		: T_{p,q}(V(S_1),V(S_2)) \otimes T_{p,\infty_S}(V(S_1),V(S_2))  &\to T_{p,\infty_S}(V(S_1),V(S_2)) 
		\\
		\begin{tikzpicture}[guillpart]
			\fill[guillfill] (0,0) rectangle (1,2);
			\draw[guillsep] (0,2)--(0,0)--(1,0)--(1,2)  (0,1)--(1,1)  ;
			\node at (0.5,0.5) {$1$};
			\node at (0.5,1.5) {$2$};
		\end{tikzpicture}
		: T_{p,q}(V(S_1),V(S_2)) \otimes T_{p,\infty_N}(V(S_1),V(S_2))  &\to T_{p,\infty_N}(V(S_1),V(S_2)) 
		\\
		\begin{tikzpicture}[guillpart]
			\fill[guillfill] (0,0) rectangle (2,1);
			\draw[guillsep]  (0,0)--(2,0)--(2,1)--(0,1) (1,0)--(1,1);
			\node at (0.5,0.5) {$2$};
			\node at (1.5,0.5) {$1$};
		\end{tikzpicture}
		: T_{p,q}(V(S_1),V(S_2)) \otimes T_{\infty_W,q}(V(S_1),V(S_2))  &\to T_{\infty_W,q}(V(S_1),V(S_2)) 
		\\
		\begin{tikzpicture}[guillpart]
			\fill[guillfill] (0,0) rectangle (2,1);
			\draw[guillsep] (2,0)--(0,0)--(0,1)--(2,1) (1,0)--(1,1);
			\node at (0.5,0.5) {$1$};
			\node at (1.5,0.5) {$2$};
		\end{tikzpicture}
		: T_{p,q}(V(S_1),V(S_2)) \otimes T_{\infty_E,q}(V(S_1),V(S_2))  &\to T_{\infty_E,q}(V(S_1),V(S_2)) 
	\end{align}
\end{subequations}
\endgroup
Checking the axiom of an $\Guill_1$ action as well as the square associativity is then again an easy exercise and the proof follows from the one made for the spaces $(T_{p,q})$ with finite $p$ and $q$.

\subsubsection{Actions on corners from morphisms between direct sums.}

The same constructions of consistent products on direct sums are generalized to corner spaces: the only difference is the treatment of direct sums over the product set $\setL^*\times\setL^*$.

\begin{lemm}
	Let $V$ and $W$ be two vector spaces. Let $A_{r,p,q}$ be the map
	\begin{align*}
		A_{r,p,q} : \left( V^{\otimes r} \otimes \End(W)^{\otimes q} \right) \otimes \left( V^{\otimes p} \otimes W^{\otimes q} \right) 
		& \to \left( V^{\otimes p+r} \otimes W^{\otimes q} \right) 
		\\
		(u \otimes A) \otimes (v\otimes w) 
		& \mapsto (u\otimes v)\otimes  Aw
	\end{align*} 
	This collection of products define the following products:
	\begin{align*}
		A_{r} : \left( V^{\otimes r} \otimes {\bigoplus_{q\in\setL}}^{(R)} \End(W)^{\otimes q} \right)
		\otimes \left( {\bigoplus_{(p,q)\in\setL^2}}^{(R,R)} V^{\otimes p}\otimes W^{\otimes q}\right) \to {\bigoplus_{(p,q)\in\setL^2}}^{(R,R)} V^{\otimes p} \otimes W^{\otimes q} 
	\end{align*}
Switching the order of the tensor products and the action of $\End(W)$ on $W^*$ instead of $W$ provides similar extensions for direct sums over $\setL^2$ with orientations $(L,R)$, $(R,L)$ and $(L,L)$.
\end{lemm}

By taking $V$ and $W$ equal to the spaces $V(S_1)$, $V(S_2)$ and their duals, we then obtain the following eight actions (two for each corner since there may be two adjacent half-strips):
\begingroup
\allowdisplaybreaks
\begin{subequations}\label{eq:canonical:actionsoncorners}
	\begin{align*}
		\begin{tikzpicture}[guillpart]
			\fill[guillfill] (0,0) rectangle (2,1);
			\draw[guillsep] (0,0)--(0,1)--(2,1)   (1,0)--(1,1)  ;
			\node at (1.5,0.5) {$2$};
			\node at (0.5,0.5) {$1$};
		\end{tikzpicture}
		: T_{p,\infty_S}(V(S_1),V(S_2)) \otimes T_{\infty_E,\infty_S}(V(S_1),V(S_2))  &\to T_{\infty_E,\infty_S}(V(S_1),V(S_2)) 
	\\
	\begin{tikzpicture}[guillpart]
		\fill[guillfill] (0,0) rectangle (1,2);
		\draw[guillsep] (0,0)--(0,2)--(1,2)   (0,1)--(1,1)  ;
		\node at (0.5,0.5) {$2$};
		\node at (0.5,1.5) {$1$};
	\end{tikzpicture}
	: T_{\infty_E,q}(V(S_1),V(S_2)) \otimes T_{\infty_E,\infty_S}(V(S_1),V(S_2))  &\to T_{\infty_E,\infty_S}(V(S_1),V(S_2)) 
	\\
	\begin{tikzpicture}[guillpart]
		\fill[guillfill] (0,0) rectangle (2,1);
		\draw[guillsep] (0,1)--(2,1)--(2,0)   (1,0)--(1,1)  ;
		\node at (0.5,0.5) {$2$};
		\node at (1.5,0.5) {$1$};
	\end{tikzpicture}
	: T_{p,\infty_S}(V(S_1),V(S_2)) \otimes T_{\infty_W,\infty_S}(V(S_1),V(S_2))  &\to T_{\infty_W,\infty_S}(V(S_1),V(S_2)) 
	\\
	\begin{tikzpicture}[guillpart]
		\fill[guillfill] (0,0) rectangle (1,2);
		\draw[guillsep] (0,2)--(1,2)--(1,0)   (0,1)--(1,1)  ;
		\node at (0.5,0.5) {$2$};
		\node at (0.5,1.5) {$1$};
	\end{tikzpicture}
	: T_{\infty_W,q}(V(S_1),V(S_2)) \otimes T_{\infty_W,\infty_S}(V(S_1),V(S_2))  &\to T_{\infty_W,\infty_S}(V(S_1),V(S_2)) 
	\\
		\begin{tikzpicture}[guillpart]
		\fill[guillfill] (0,0) rectangle (2,1);
		\draw[guillsep] (0,1)--(0,0)--(2,0)   (1,0)--(1,1)  ;
		\node at (1.5,0.5) {$2$};
		\node at (0.5,0.5) {$1$};
	\end{tikzpicture}
	: T_{p,\infty_N}(V(S_1),V(S_2)) \otimes T_{\infty_E,\infty_N}(V(S_1),V(S_2))  &\to T_{\infty_E,\infty_N}(V(S_1),V(S_2)) 
	\\
	\begin{tikzpicture}[guillpart]
		\fill[guillfill] (0,0) rectangle (1,2);
		\draw[guillsep] (1,0)--(0,0)--(0,2)   (0,1)--(1,1)  ;
		\node at (0.5,0.5) {$2$};
		\node at (0.5,1.5) {$1$};
	\end{tikzpicture}
	: T_{\infty_E,N}(V(S_1),V(S_2)) \otimes T_{\infty_E,\infty_N}(V(S_1),V(S_2))  &\to T_{\infty_E,\infty_N}(V(S_1),V(S_2)) 
	\\
	\begin{tikzpicture}[guillpart]
		\fill[guillfill] (0,0) rectangle (2,1);
		\draw[guillsep] (0,0)--(2,0)--(2,1)   (1,0)--(1,1)  ;
		\node at (0.5,0.5) {$2$};
		\node at (1.5,0.5) {$1$};
	\end{tikzpicture}
	: T_{p,\infty_N}(V(S_1),V(S_2)) \otimes T_{\infty_W,\infty_N}(V(S_1),V(S_2))  &\to T_{\infty_W,\infty_N}(V(S_1),V(S_2)) 
	\\
	\begin{tikzpicture}[guillpart]
		\fill[guillfill] (0,0) rectangle (1,2);
		\draw[guillsep] (0,0)--(1,0)--(1,2)   (0,1)--(1,1)  ;
		\node at (0.5,0.5) {$1$};
		\node at (0.5,1.5) {$2$};
	\end{tikzpicture}
	: T_{\infty_W,q}(V(S_1),V(S_2)) \otimes T_{\infty_W,\infty_N}(V(S_1),V(S_2))  &\to T_{\infty_W,\infty_N}(V(S_1),V(S_2)) 
	\end{align*}
\end{subequations}
\endgroup
Checking again the various associativity from the fundamental associativities of tensor products and actions of endomorphisms on vectors is again left as an exercise.

\subsection{The new coloured spaces on boundaries, part two (with a base point).}

We now define the canonical spaces for Markov processes when the pattern shapes have at least a doubly-infinite directions, i.e. for a pattern in $\{\patterntype{s}_{WE},\patterntype{s}_{SN},\patterntype{hp}_{S}^*, \patterntype{hp}_{N}^*, \patterntype{hp}_{W}^*, \patterntype{hp}_{E}^*, \patterntype{fp}^* \}$. The difference with the previous spaces is that the morphisms now require a further information, the base point. In the present case, this will take the form of an additional action of the translation group on the various spaces. 

Previously, for each non-doubly-infinite shape $\patterntype{u}$, the corresponding space gave rise to a direct sum over $\setL$ of tensor products with a left or right action. This may interpreted geometrically as fixing a reference side (the one not sent to infinity) and the index $p\in \setL$ plays the role of a distance to this point. In the present case, there is no canonical reference point and one must choose one arbitrarily (our so-called base point). This means that the tensor product $V^{\otimes p}$ may be divided arbitrarily into two parts $V^{\otimes p_1}\otimes V^{\otimes p_2}$ (if the base point lies in-between) and then direct sums may be considered over both $p_1$ and $p_2$. This last case does not however contain all the cases and there is no action of the translation group: a space like
\[
{\bigoplus_{p_1\in\setL}}^{(L)}
{\bigoplus_{p_2\in\setL}}^{(R)}
V^{p_1}\otimes V^{p_2}
\]
may be enough to contain gluing of two corners or half-strips with a base point placed on the guillotine cut but there is no action of the translation group on it.

Geometrically, the missing case corresponds to the case of a base point taken outside the segment $[0,p]$. In this case, distances to the base point shall be taken algebraically with signs, while keeping the information that the left extremities lies on the left of the right extremities. We must therefore use a direct sum over 
\begin{equation}\label{eq:indexingsetfortranslationinvar}
\setP^2_{\leq} = \left\{ (x_L,x_R)\in\setP^2 ; x_L\leq x_R \right\}
\end{equation}
of tensor products $V^{\otimes (x_R-x_L)}$. In particular, there is a canonical embedding map:
$\setL\times \setL \to \setP^2_\leq$, $(p,q) \mapsto (-p,q)$. There is then a new informal rule to obtain the new spaces from the canonical spaces $T_{p,q}(V(S_1),V(S_2))$ on rectangles: any doubly-infinite line (with size $\infty_{SN}$ or $\infty_{WE}$) on these pattern type gives rise to a direct sum over $(x_L,x_R)\setP^2_\leq$ of the tensor products $V^{\otimes p}$ with $p=x_R-x_L$.

We define the following spaces, both for finite state space in discrete space with $p,q\geq 1$:
\begingroup
\allowdisplaybreaks
\begin{subequations}
\label{eq:defcanonicalspaces:doublyinfinitepatterns:FD}
\begin{align}
	T_{\infty_{WE},0}( V(S_1), V(S_2) ) &=
	{\bigoplus_{(x_L,x_R)\in\setP^2_\leq}}
	\Diag(V(S_1))^{\otimes x_R-x_L}
	\\
	T_{0,\infty_{SN}}( V(S_1), V(S_2) ) &=
	{\bigoplus_{(y_L,y_R)\in\setP^2_\leq}}
	\Diag(V(S_2))^{\otimes y_R-y_L}
	\\
	T_{\infty_{WE},q}( V(S_1), V(S_2) ) &=
	{\bigoplus_{(x_L,x_R)\in\setP^2_\leq}}
	\End(V(S_1))^{\otimes x_R-x_L}
	\\
	T_{p,\infty_{SN}}( V(S_1), V(S_2) ) &=
	{\bigoplus_{(y_L,y_R)\in\setP^2_\leq}}
	\End(V(S_2))^{\otimes y_R-y_L}
	\\
	T_{\infty_{WE},\infty_S}( V(S_1), V(S_2) ) &=
	{\bigoplus_{(x_L,x_R)\in\setP^2_\leq}}
	(V(S_1)^*)^{\otimes x_R-x_L}
	\\
	T_{\infty_{WE},\infty_N}( V(S_1), V(S_2) ) &=
	{\bigoplus_{(x_L,x_R)\in\setP^2_\leq}}
	V(S_1)^{\otimes x_R-x_L}
	\\
	T_{\infty_{W},\infty_{SN}}( V(S_1), V(S_2) ) &=
	{\bigoplus_{(y_L,y_R)\in\setP^2_\leq}}
	(V(S_2)^*)^{\otimes y_R-y_L}
	\\
	T_{\infty_{E},\infty_{SN}}( V(S_1), V(S_2) ) &=
	{\bigoplus_{(y_L,y_R)\in\setP^2_\leq}}
	V(S_2)^{y_R-y_L}
	\\
	T_{\infty_{WE},\infty_{SN}}( V(S_1), V(S_2) ) &= \setK
\end{align}
\end{subequations}
\endgroup
All these spaces have the same structure of a direct sum 
over $\setP^2_\leq$ where two spaces indexed by $(u_L,u_R)$ and $(v_L,v_R)$ are the same as soon as $v_R-v_L=u_R-u_L$. We first define an action of translations on  such directs sums and prove directly the following lemma.

\begin{lemm}\label{lemm:canonical:translationinvariance}
	Let $(E_{u,v})_{(u,v)\in \setP^2_\leq}$ be a collection of sets such that $v'-u'=v-u$ implies $E_{u,v}=E_{u',v'}$. The group of translations of $\setP$ acts on $(E_{u,v})$ in the following way: the translation $\theta_k : x\in \setP \mapsto x+k$ acts on $(E_{u,v})$ through
	\[
		\theta_{k}( (a_{u,v})_{(u,v)\in\setP^2_\leq} ) = (a_{k+u,k+v})_{(u,v)\in\setP^2_\leq}
	\]
\end{lemm}

We may now define pairings from the previous shapes, half-strips and corners, which produce strips and half-plane.

\begin{lemm}[pairings from half-strips to strips and from corners to half-spaces]\label{eq:canonical:halfstriptostrip}
Let $V$, $V'$ and $W$ be vector spaces and $\scal{\cdot}{\cdot}: V'\otimes V\to\setK$ a bilinear form. Let $(H_{q_1,q_2})$ be the collection of maps defined by
\begin{align*}
	H_{q_1,q_2}: \left(V' \otimes W^{\otimes q_1}\right) \otimes \left(V \otimes W^{\otimes q_2}\right) &\to W^{\otimes q_1+ q_2} \\
	 (a\otimes w_1)\otimes (b\otimes w_2) & \mapsto \scal{a}{b} (w_1\otimes w_2)
\end{align*}
These maps define component-wise maps for each $x\in\setP$:
\begin{align*}
	H^{(x)}: \left(V' \otimes {\bigoplus_{q_1\in\setL}}^{(L)} W^{\otimes q_1}\right) \otimes \left(V \otimes {\bigoplus_{q_2\in\setL}}^{(R)} W^{\otimes q_2}\right) 
	& \to 
	{\bigoplus_{(u_L,u_R)\in\setP^2_\leq}}
	W^{\otimes u_R-u_L}
	\\
	(a\otimes w_{q_1})\otimes (b\otimes w_{q_2}) \mapsto 
	\theta_x\left([H_{q_1,q_2}(a\otimes w_{q_1})\otimes (b\otimes w_{q_2})]_{\eta(q_1,q_2} \right)
\end{align*}
where $[x]_{u,v}$ indicates that the element $x\in W^{v-u}$ is placed in the $(u,v)$-component of the direct sum (and 0 are put in the other spaces).
\end{lemm}
Applying this lemma to the case $V=V(S_i)^{\otimes r}$, $V'=(V(S_i)^*)^{\otimes r}$ and $W=\End(V(S_j))$, we obtain the products with a base point placed at the location of the guillotine cut and then, by left composition by $\theta_x$, we obtain all the following $\Guill_2$-products:
\begin{subequations}
\label{eq:canonical:pairingofstrips}
\begin{align}
	\begin{tikzpicture}[guillpart,yscale=0.65,xscale=0.75,baseline={(current bounding box.center)}]
		\fill[guillfill] (3,0)--(0,0)--(0,2)--(3,2)--(3,0);
		\draw[guillsep] (3,0)--(0,0);
		\draw[guillsep] (0,2)--(3,2);
		\draw[guillsep] (2,0)--(2,2);
		\node at (1.33,1.) {$1$};
		\node at (2.5,1.) {$2$};	
		\node at (0.66,0) [circle, fill, inner sep=0.5mm] {};
		\node at (0.66,2) [circle, fill, inner sep=0.5mm] {};
		\draw [dotted] (0.66,0) -- (0.66,2);
		\draw [->] (2,2.3) -- node [midway, above] {$x$} (0.66,2.3);
	\end{tikzpicture}
	: T_{\infty_{W},q}(V(S_1),V(S_2))\otimes T_{\infty_E,q}(V(S_1),V(S_2))  &\to T_{\infty_{WE},q}(V(S_1),V(S_2))
	\\
	\begin{tikzpicture}[guillpart,baseline={(current bounding box.center)},rotate=90,scale=0.75]
		\fill[guillfill] (3,0)--(0,0)--(0,2)--(3,2)--(3,0);
		\draw[guillsep] (3,0)--(0,0);
		\draw[guillsep] (0,2)--(3,2);
		\draw[guillsep] (2,0)--(2,2);
		\node at (1.33,1.) {$1$};
		\node at (2.5,1.) {$2$};	
		\node at (0.66,0) [circle, fill, inner sep=0.5mm] {};
		\node at (0.66,2) [circle, fill, inner sep=0.5mm] {};
		\draw [dotted] (0.66,0) -- (0.66,2);
		\draw [->] (2,2.3) -- node [midway, left] {$x$} (0.66,2.3);
	\end{tikzpicture}
	: T_{p,\infty_{S}}(V(S_1),V(S_2))\otimes T_{p,\infty_N}(V(S_1),V(S_2))  &\to T_{p,\infty_{SN}}(V(S_1),V(S_2))
\end{align}
\end{subequations}
Applying this lemma to the case $V={\bigoplus_{p\in\setL}}^{L|R} V(S_i)^{\otimes p}$ and $V'={\bigoplus_{p\in\setL}}^{L|R} (V(S_i)^*)^{\otimes p}$ with a bilinear form obtained by 
\[
B( (v'_p)_p , (v_p)_p ) = \sum_{p\in\setL} \scal{v'_p}{v_p}
\]
(which is well defined since the direct sums contain only a finite number of non-zero terms), we obtained the following elementary products associated to the gluing between consecutive corners that produces a half-plane, after use of a translation $\theta_x$ to shift the base point off the guillotine cut:
\begingroup
\allowdisplaybreaks
\begin{subequations}
	\label{eq:canonical:pairingofcorners}
	\begin{align}
		\begin{tikzpicture}[guillpart,yscale=0.75,baseline={(current bounding box.center)}]
			\fill[guillfill] (0,0)--(0,2)--(3,2)--(3,0);
			%\draw[guillsep] (3,0)--(0,0);
			\draw[guillsep] (0,2)--(3,2);
			\draw[guillsep] (2,0)--(2,2);
			\node at (1.33,1.) {$1$};
			\node at (2.5,1.) {$2$};	
			\node at (0.66,2) [circle, fill, inner sep=0.5mm] {};
			\draw [dotted] (0.66,0) -- (0.66,2);
			\draw [->] (2,2.3) -- node [midway, above] {$x$} (0.66,2.3);
		\end{tikzpicture}
		: T_{\infty_{W},\infty_S}(V(S_1),V(S_2))\otimes T_{\infty_E,\infty_S}(V(S_1),V(S_2))  &\to T_{\infty_{WE},\infty_S}(V(S_1),V(S_2))
		\\
		\begin{tikzpicture}[guillpart,baseline={(current bounding box.center)},rotate=90,scale=0.75]
			\fill[guillfill] (3,0)--(0,0)--(0,2)--(3,2)--(3,0);
			%\draw[guillsep] (3,0)--(0,0);
			\draw[guillsep] (0,2)--(3,2);
			\draw[guillsep] (2,0)--(2,2);
			\node at (1.33,1.) {$1$};
			\node at (2.5,1.) {$2$};	
			%\node at (0.66,0) [circle, fill, inner sep=0.5mm] {};
			\node at (0.66,2) [circle, fill, inner sep=0.5mm] {};
			\draw [dotted] (0.66,0) -- (0.66,2);
			\draw [->] (2,2.3) -- node [midway, left] {$x$} (0.66,2.3);
		\end{tikzpicture}
		: T_{\infty_E,\infty_{S}}(V(S_1),V(S_2))\otimes T_{\infty_E,\infty_N}(V(S_1),V(S_2))  &\to T_{\infty_E,\infty_{SN}}(V(S_1),V(S_2))
	\\
		\begin{tikzpicture}[guillpart,yscale=0.75,baseline={(current bounding box.center)}]
		\fill[guillfill] (3,0)--(0,0)--(0,2)--(3,2)--(3,0);
		\draw[guillsep] (3,0)--(0,0);
		%\draw[guillsep] (0,2)--(3,2);
		\draw[guillsep] (2,0)--(2,2);
		\node at (1.33,1.) {$1$};
		\node at (2.5,1.) {$2$};	
		\node at (0.66,0) [circle, fill, inner sep=0.5mm] {};
		%\node at (0.66,2) [circle, fill, inner sep=0.5mm] {};
		\draw [dotted] (0.66,0) -- (0.66,2);
		\draw [->] (2,-0.3) -- node [midway, below] {$x$} (0.66,-0.3);
	\end{tikzpicture}
	: T_{\infty_{W},\infty_N}(V(S_1),V(S_2))\otimes T_{\infty_E,\infty_N}(V(S_1),V(S_2))  &\to T_{\infty_{WE},\infty_N}(V(S_1),V(S_2))
	\\
	\begin{tikzpicture}[guillpart,baseline={(current bounding box.center)},rotate=90,scale=0.75]
		\fill[guillfill] (3,0)--(0,0)--(0,2)--(3,2)--(3,0);
		\draw[guillsep] (3,0)--(0,0);
		%\draw[guillsep] (0,2)--(3,2);
		\draw[guillsep] (2,0)--(2,2);
		\node at (1.33,1.) {$1$};
		\node at (2.5,1.) {$2$};	
		\node at (0.66,0) [circle, fill, inner sep=0.5mm] {};
		%\node at (0.66,2) [circle, fill, inner sep=0.5mm] {};
		\draw [dotted] (0.66,0) -- (0.66,2);
		\draw [->] (2,-0.3) -- node [midway, right] {$x$} (0.66,-0.3);
	\end{tikzpicture}
	: T_{\infty_W,\infty_{S}}(V(S_1),V(S_2))\otimes T_{\infty_W,\infty_N}(V(S_1),V(S_2))  &\to T_{\infty_W,\infty_{SN}}(V(S_1),V(S_2))
	\end{align}
\end{subequations}
\endgroup

We may now combine doubly-infinite shapes together in a trivial way.
\begin{lemm}
	The spaces $(T_{\infty_{WE},q}(V(S_1),V(S_2)))$ (resp. $(T_{p,\infty_{SN}}(V(S_1),V(S_2)))$) are $\Guill_1$-algebras as direct sums of $\Guill_1$-algebras $\End(V)^{\otimes p}$.
\end{lemm}
This provides the definition of the following two products associated to the gluing of strips:
\begin{subequations}
	\label{eq:canonical:gluingstrips}
\begin{align}
	\begin{tikzpicture}[guillpart,baseline={(current bounding box.center)},xscale=0.8]
		\fill[guillfill] (2.5,0)--(0,0)--(0,2)--(2.5,2)--(2.5,0);
		%\draw[guillsep] (2.5,0)--(0,0);
		\draw[guillsep] (0,2)--(2.5,2);
		\draw[guillsep] (0,1)--(2.5,1);
		\draw[guillsep] (0,0)--(2.5,0);
		\node at (1.5,0.5) {$1$};
		\node at (1.5,1.5) {$2$};	
		\node at (0.66,0) [circle, fill, inner sep=0.5mm] {};
		\node at (0.66,2) [circle, fill, inner sep=0.5mm] {};
		\draw [dotted] (0.66,0) -- (0.66,2);
	\end{tikzpicture}	
	: T_{\infty_{WE},q_1}(V(S_1),V(S_2)) \otimes T_{\infty_{WE},q_2}(V(S_1),V(S_2)) &\to 
	T_{\infty_{WE},q_1+q_2}(V(S_1),V(S_2))
	\\
	\begin{tikzpicture}[guillpart,baseline={(current bounding box.center)},rotate=-90,xscale=0.8]
		\fill[guillfill] (2.5,0)--(0,0)--(0,2)--(2.5,2)--(2.5,0);
		\draw[guillsep] (2.5,0)--(0,0);
		\draw[guillsep] (0,2)--(2.5,2);
		\draw[guillsep] (0,1)--(2.5,1);
		\node at (1.5,0.5) {$1$};
		\node at (1.5,1.5) {$2$};	
		\node at (0.66,0) [circle, fill, inner sep=0.5mm] {};
		\node at (0.66,2) [circle, fill, inner sep=0.5mm] {};
		\draw [dotted] (0.66,0) -- (0.66,2);
	\end{tikzpicture}	
	: T_{p_1,\infty_{SN}}(V(S_1),V(S_2)) \otimes T_{p_2,\infty_{SN}}(V(S_1),V(S_2)) &\to 
	T_{p_1+p_2,\infty_{SN}}(V(S_1),V(S_2))
\end{align}
\end{subequations}

\begin{lemm}
	The spaces $T_{\infty_{WE},q}(V(S_1),V(S_2))$ (resp. $T_{p,\infty_{SN}}(V(S_1),V(S_2))$) act trivially on $T_{\infty_{WE},\infty_S}(V(S_1),V(S_2))$ and $T_{\infty_{WE},\infty_N}(V(S_1),V(S_2))$ (resp. $T_{\infty_{W},\infty_{SN}}$ and $T_{\infty_{E},\infty_{SN}}$) through the component-wise action of $\End(V(S_1))^{\otimes x_R-x_L}$ on $(V(S_1)^*)^{\otimes x_R-x_L}$ and $V(S_1)^{\otimes x_R-x_L}$ (resp. $\End(V(S_2))^{\otimes y_R-y_L}$ on $(V(S_2)^*)^{\otimes y_R-y_L}$ and $V(S_2)^{\otimes y_R-y_L}$.
\end{lemm}
This lemma provides the definition of the following products associated to the gluing of a strip on a half-plane:
\begingroup
\allowdisplaybreaks
\begin{subequations}
	\label{eq:canonical:gluingstripsonhalfplanes}
\begin{align}
\begin{tikzpicture}[guillpart,baseline={(current bounding box.center)},xscale=0.8]
	\fill[guillfill] (2.5,0)--(0,0)--(0,2)--(2.5,2)--(2.5,0);
	%\draw[guillsep] (2.5,0)--(0,0);
	\draw[guillsep] (0,2)--(2.5,2);
	\draw[guillsep] (0,1)--(2.5,1);
	\node at (1.5,0.5) {$1$};
	\node at (1.5,1.5) {$2$};	
	%\node at (0.66,0) [circle, fill, inner sep=0.5mm] {};
	\node at (0.66,2) [circle, fill, inner sep=0.5mm] {};
	\draw [dotted] (0.66,0) -- (0.66,2);
\end{tikzpicture}	
: T_{\infty_{WE},\infty_S}(V(S_1),V(S_2)) \otimes T_{\infty_{WE},q}(V(S_1),V(S_2)) &\to 
T_{\infty_{WE},\infty_S}(V(S_1),V(S_2))
\\
\begin{tikzpicture}[guillpart,baseline={(current bounding box.center)},xscale=0.8]
	\fill[guillfill] (2.5,0)--(0,0)--(0,2)--(2.5,2)--(2.5,0);
	\draw[guillsep] (2.5,0)--(0,0);
	%\draw[guillsep] (0,2)--(2.5,2);
	\draw[guillsep] (0,1)--(2.5,1);
	\node at (1.5,0.5) {$1$};
	\node at (1.5,1.5) {$2$};	
	\node at (0.66,0) [circle, fill, inner sep=0.5mm] {};
	%\node at (0.66,2) [circle, fill, inner sep=0.5mm] {};
	\draw [dotted] (0.66,0) -- (0.66,2);
\end{tikzpicture}	
: T_{\infty_{WE},q}(V(S_1),V(S_2)) \otimes T_{\infty_{WE},\infty_N}(V(S_1),V(S_2)) &\to 
T_{\infty_{WE},\infty_N}(V(S_1),V(S_2))
\\
\begin{tikzpicture}[guillpart,baseline={(current bounding box.center)},rotate=-90]
	\fill[guillfill] (2.5,0)--(0,0)--(0,2)--(2.5,2)--(2.5,0);
	\draw[guillsep] (2.5,0)--(0,0);
	%\draw[guillsep] (0,2)--(2.5,2);
	\draw[guillsep] (0,1)--(2.5,1);
	\node at (1.5,0.5) {$1$};
	\node at (1.5,1.5) {$2$};	
	\node at (0.66,0) [circle, fill, inner sep=0.5mm] {};
	%\node at (0.66,2) [circle, fill, inner sep=0.5mm] {};
	\draw [dotted] (0.66,0) -- (0.66,2);
\end{tikzpicture}	
: T_{p,\infty_{SN}}(V(S_1),V(S_2)) \otimes T_{\infty_{E},\infty_{SN}}(V(S_1),V(S_2)) &\to 
T_{\infty_{E},\infty_{SN}}(V(S_1),V(S_2))
\\
\begin{tikzpicture}[guillpart,baseline={(current bounding box.center)},rotate=-90]
	\fill[guillfill] (2.5,0)--(0,0)--(0,2)--(2.5,2)--(2.5,0);
	%\draw[guillsep] (2.5,0)--(0,0);
	\draw[guillsep] (0,2)--(2.5,2);
	\draw[guillsep] (0,1)--(2.5,1);
	\node at (1.5,0.5) {$1$};
	\node at (1.5,1.5) {$2$};	
	%\node at (0.66,0) [circle, fill, inner sep=0.5mm] {};
	\node at (0.66,2) [circle, fill, inner sep=0.5mm] {};
	\draw [dotted] (0.66,0) -- (0.66,2);
\end{tikzpicture}	
: T_{\infty_{W},\infty_{SN}}(V(S_1),V(S_2)) \otimes T_{p,\infty_{SN}}(V(S_1),V(S_2)) &\to 
T_{\infty_{W},\infty_{SN}}(V(S_1),V(S_2))
\end{align}
\end{subequations}
\endgroup

We finally provide the final pairing of half-plane spaces to produce scalar numbers associated to the whole plane.
\begin{lemm}
	The canonical pairing of spaces $(V(S_i)^*)^{\otimes p}\otimes V(S_i)^{\otimes p} \to \setK$ on each component of the direct sums in half-plane spaces provides, after finite summation over the index $p$, the two products
	\begin{subequations}
		\label{eq:defcanonical:gluingdoublyinfinite}
		\begin{align}
			\begin{tikzpicture}[guillpart,baseline={(current bounding box.center)},xscale=0.8]
				\fill[guillfill] (2.5,0)--(0,0)--(0,2)--(2.5,2)--(2.5,0);
				%\draw[guillsep] (2.5,0)--(0,0);
				%\draw[guillsep] (0,2)--(2.5,2);
				\draw[guillsep] (0,1)--(2.5,1);
				\node at (1.5,0.5) {$1$};
				\node at (1.5,1.5) {$2$};	
				%\node at (0.66,0) [circle, fill, inner sep=0.5mm] {};
				\node at (0.66,1) [circle, fill, inner sep=0.5mm] {};
				\draw [dotted] (0.66,0) -- (0.66,2);
			\end{tikzpicture}	
			: T_{\infty_{WE},\infty_S}(V(S_1),V(S_2)) \otimes T_{\infty_{WE},\infty_N}(V(S_1),V(S_2)) &\to 
			T_{\infty_{WE},\infty_{SN}}(V(S_1),V(S_2))
			\\
			\begin{tikzpicture}[guillpart,baseline={(current bounding box.center)},rotate=-90,xscale=0.8]
				\fill[guillfill] (2.5,0)--(0,0)--(0,2)--(2.5,2)--(2.5,0);
				%\draw[guillsep] (2.5,0)--(0,0);
				%\draw[guillsep] (0,2)--(2.5,2);
				\draw[guillsep] (0,1)--(2.5,1);
				\node at (1.5,0.5) {$1$};
				\node at (1.5,1.5) {$2$};	
				%\node at (0.66,0) [circle, fill, inner sep=0.5mm] {};
				\node at (0.66,1) [circle, fill, inner sep=0.5mm] {};
				\draw [dotted] (0.66,0) -- (0.66,2);
			\end{tikzpicture}	
			: T_{\infty_{W},\infty_{SN}}(V(S_1),V(S_2)) \otimes T_{\infty_{E},\infty_{SN}}(V(S_1),V(S_2)) &\to 
			T_{\infty_{WE},\infty_{SN}}(V(S_1),V(S_2))
			\\
			(a'_{u,v})_{(u,v)\in\setP^2_\leq} \otimes (a_{u,v})_{(u,v)\in\setP^2_\leq}
			&\mapsto
			\sum_{(u,v)\in\setP^2_\leq} \scal{a'_{u,v}}{a_{u,v}}
		\end{align}
	\end{subequations}
\end{lemm}

\subsection{A global structure theorem for canonical boundary spaces..}

We have defined in this section a whole set of spaces $(T_{p,q}(V(S_1),V(S_2)))_{(p,q)\in\PatternShapes(\patterntype{fp}^*)}$  with suitable products, actions and pairings. The following theorem summarizes all these constructions in an operadic way.

\begin{theo}[canonical boundary guillotine algebra for Markov processes]\label{theo:canonicalboundarystructure}
The collection of spaces 
$(T_{p,q}(V(S_1),V(S_2)))_{(p,q)\in \PatternShapes(\patterntype{fp}^*)}$  introduced 
in the previous sections with their elementary guillotine products define a $\Guill_2^{(\patterntype{fp}^*)}$-algebra.
\end{theo}
\begin{proof}
The extended guillotine operad $\Guill_2^{(\patterntype{fp}^*)}$ is generated by the fundamental guillotine partitions, pointed or not depending on whether there is an doubly-infinite dimension, given by a single guillotine cut, which satisfy the associativities \eqref{eq:guill2:horizassoc}, \eqref{eq:guill2:vertassoc} and \eqref{eq:guill2:interchangeassoc} and their pointed versions.   This corresponds to properties \ref{prop:guill2generators}, \ref{prop:guill2elemassoc}, the first extended version  \ref{prop:extendedgeneratorsI} and the second extended version~\ref{prop:extendedgeneratorsII}.

Thus, it is enough to verify that the products, actions and pairings defined previously on the spaces $T_{p,q}(V_1,V_2)$ satisfy these associativity relations. This is then just an exercise on tensor products of vectors spaces and endomorphisms once we recognize that the usual and concatenation products used in equations
\eqref{eq:productandconcat:assoc} are extended in a straightforward way to these spaces.
\end{proof}

Whenever no doubt on the spaces $V(S_1)$ and $V(S_2)$ nor on the segment state space $S^{(1)}_\bullet$ and $S^{(2)}_\bullet$ arises, e.g. if one considers a fixed given Markov process with its face weights, this $\Guill_2^{(\patterntype{fp}^*)}$-algebra of spaces will be denoted as follows for simplicity, using again notation~\eqref{eq:shortnotation:guillalgebra}:
	\begin{equation}\label{notation:shortTpq}
		\ca{T}_{\PatternShapes(\patterntype{fp}^*)} =
		 (T_{p,q}(V(S_1),V(S_2)))_{(p,q)\in\PatternShapes(\patterntype{fp}^*)} 
	\end{equation}

\subsection{A first simple probabilistic consequence: 2D Markov processes with factorized boundary weights.}\label{sec:factorizedweights:partone}

\subsubsection{From algebra to probability}
Such a collection of canonical boundary spaces $(\ca{T}_{p,q})_{(p,q)\in\PatternShapes(\patterntype{fp}^*)\setminus\PatternShapes(\patterntype{r})}$ as defined previously can not describe generic boundary conditions of arbitrary Markov processes. The reason is that all a generic element of the direct sums is a finite linear combination of finite dimensional spaces: it can not describe the whole complexity of a two-dimensional model defined as a Gibbs measure on the whole plane and, moreover, it would imply identities on correlation functions that would probably have been already identified in the literature of statistical mechanics.

The restriction to finite linear combinations however describes nicely marginals of Markov processes on a \emph{finite} domain with boundary weights that can be factorized over the edges of the boundary, as it will be seen in the following proposition. 

This is a first step towards Theorem~\ref{theo:stability} and already illustrate the key features of the next sections~\ref{sec:boundaryalgebra} and \ref{sec:invariantboundaryelmts}, which will contain additional algebraic features for more general boundary weights.

\subsubsection{Partition functions with boundary weights factorized over the edges.}

\begin{prop}\label{prop:partitionfuncfactorweights} 
Let $S_1$ and $S_2$ be two finite sets and let $V(S_1)=\setR^{S_1}$ and $V(S_2)=\setR^{S_2}$ (with their canonical bases) be their associated finite-dimensional vector spaces. Let $R=[P_1,P_2]\times [Q_1,Q_2]$ be a rectangle in $\setZ^2$ with shape $(P,Q)=(P_2-P_1,Q_2-Q_1)$. Let $(X_e)_{e\in \Edges{R}}$ be a $(S_1,S_2)$-valued Markov process on $R$ with face weights $(\MarkovWeight{W}_{x,y})_{P_1\leq x < Q_1,P_2\leq y< Q_2}$ and a boundary weight $G_R : (x_e)_{e\in\Edges{\partial R}} \mapsto \setR_+$ given by a factorized expression
\begin{equation}
G_R( (x_e)_{e\in\Edges{\partial R}} ) =\prod_{e\in\Edges{\partial R}} \gamma_e(x_e)
\end{equation}
where $\gamma_e$ is a function from $S_1$ (resp. $S_2$) to $\setR_+$ if $e$ is an horizontal (resp. vertical) edge.

Then, the partition function $Z_{R}^{\boundaryweights}(\MarkovWeight{W}_\bullet;\Phi)$ is a real number given by:
\begin{equation}\label{eq:partitionfuncwithfactorbw}
Z_R^{\boundaryweights}(\MarkovWeight{W}_\bullet;G_R) =
\begin{tikzpicture}[guillpart]
\begin{scope}[scale=1.65 ];
\fill[guillfill] (0,0) rectangle (6,5);
\draw[guillsep] (0,1) -- (6,1)
			 (0,4) -- (6,4)
			 (1,0) -- (1,5)
			 (5,0) -- (5,5);
\draw[guillsep] (2,0)--(2,1)	(3,0)--(3,1) (4,0)--(4,1);
\draw[guillsep] (2,4)--(2,5)	(3,4)--(3,5) (4,4)--(4,5);
\draw[guillsep] (0,2)--(1,2) (0,3)--(1,3) (0,4)--(1,4);
\draw[guillsep] (5,2)--(6,2) (5,3)--(6,3) (5,4)--(6,4);
\node at (0.5,0.5) {$u_{SW}$}; 
\node at (5.5,0.5) {$u_{SE}$};
\node at (0.5,4.5) {$u_{NW}$};
\node at (5.5,4.5) {$u_{NE}$};

\node at (1.5,0.5) {$\ha{\gamma}_{e_1^S}$};
\node at (2.5,0.5) {$\ha{\gamma}_{e_2^S}$};
\node at (3.5,0.5) {$\dots$};
\node at (4.5,0.5) {$\ha{\gamma}_{e_P^S}$};

\node at (1.5,4.5) {$\ha{\gamma}_{e_1^N}$};
\node at (2.5,4.5) {$\ha{\gamma}_{e_2^N}$};
\node at (3.5,4.5) {$\dots$};
\node at (4.5,4.5) {$\ha{\gamma}_{e_P^N}$};

\node at (0.5,1.5) {$\ha{\gamma}_{e_1^W}$};
\node at (0.5,2.5) {$\vdots$};
\node at (0.5,3.5) {$\ha{\gamma}_{e_Q^W}$};

\node at (5.5,1.5) {$\ha{\gamma}_{e_1^E}$};
\node at (5.5,2.5) {$\vdots$};
\node at (5.5,3.5) {$\ha{\gamma}_{e_Q^E}$};

\node at (3,2.5) {$Z_R(\MarkovWeight{W}_\bullet;\cdot)$};

\node at (4,2) [circle,fill,inner sep=0.5mm] {};
\draw[dotted] (4,0)--(4,5) (0,2)--(6,2);
\end{scope}
\end{tikzpicture}
\end{equation}
in the canonical $\Guill_2^{(\patterntype{fp}^*)}$-algebra $T(V(S_1),V(S_2))_{\PatternShapes(\patterntype{fp})}$
where:
\begin{itemize}
\item for each $a\in\{N,W,S,E\}$, $(e^{a}_i)_{i}$ is a left to right (or bottom to top) enumeration of the edges of the side $a$,
\item the base point is arbitrary and plays no role,
\item the partition function $Z_R(\MarkovWeight{W}_\bullet,\cdot)$ is the element of $T_{P,Q}(V(S_1),V(S_2))$ given by Theorem~\ref{theo:partitionfuncguillotop},
\item the corner elements are given by $u_{ab}=1$ for $(a,b)\in\{S,N,\}\times\{W,E\}$, where the $1$ are placed in the $(0,0)$-component, equal to $\setK$, of the direct sums over $(p,q)\in\setN\times\setN$ of \eqref{eq:canonicalboundaryalg:FD}
and the boundary elements are given by:
\begin{align*}
\ha{\gamma}_{e^S_i} &= \gamma_{e^S_i} \otimes 1 
&
\ha{\gamma}_{e^N_i} &= \gamma_{e^N_i} \otimes 1
\\
\ha{\gamma}_{e^W_i} &=  1 \otimes\gamma_{e^W_i} 
&
\ha{\gamma}_{e^E_i} &=  1\otimes\gamma_{e^E_i} 
\end{align*}
where the functions $\gamma_{\bullet}$ are identified with the corresponding vector in $V(S_i)$ and where the $1$ is placed in the $0$-component, equal to $\setK$ of the direct sums over $p\in\setN$ or $q\in\setN$ of \eqref{eq:canonicalboundaryalg:FD}.
\end{itemize}
\end{prop}
\begin{proof}
The proof is a simple exercise in evaluation of tensor products in the previous spaces. We first consider the marginal law of the process on its boundaries:
We already have, from Section~\ref{sec:prob:twodim}, by definition:
\begin{equation}\label{eq:relationbetweenpartitionfunc}
Z_R^{\boundaryweights}(\MarkovWeight{W}_\bullet;G_R) = \sum_{(x_e)_{e\in\Edges{R}}} Z_R(\MarkovWeight{W}_\bullet;(x_e)_{e\in\Edges{R}} ) \prod_{e\in\Edges{R}} \gamma_e(x_e)
\end{equation}
which we will reinterpret as a representation of the extended guillotine operad $\Guill_2^{(\patterntype{fp}^*)}$.

We now evaluate the r.h.s. of~\eqref{eq:partitionfuncwithfactorbw} using Theorem~\ref{theo:canonicalboundarystructure} and the generators of the products in the canonical structure. We consider first the boundary algebra terms and we have:
\begin{equation*}
\begin{tikzpicture}[guillpart,yscale=0.7]
\begin{scope}[scale=1.8]
\fill[guillfill] (0,0) rectangle (4,1);
\draw[guillsep] (0,1) -- (4,1)
		(0,0) -- (0,1)	(1,0)--(1,1)	  (2,0)--(2,1)  (3,0)--(3,1)  (4,0)--(4,1);
\node at (0.5,0.5) {$\ha{\gamma}_{e_1^S}$};
\node at (1.5,0.5) {$\ha{\gamma}_{e_2^S}$};
\node at (2.5,0.5) {$\dots$};
\node at (3.5,0.5) {$\ha{\gamma}_{e_P^S}$};
\end{scope}
\end{tikzpicture}
= \left( \gamma_{e_1^S} \otimesdots \gamma_{e_P^S} \right)
\otimes 1
\end{equation*}
from the products~\eqref{eq:productboundaryalgebra} and similar formulae hold on the three other sides. We then obtain
\begin{equation*}
\begin{tikzpicture}[guillpart,yscale=0.7]
\begin{scope}[scale=1.8]
\fill[guillfill] (-1,0) rectangle (5,1);
\draw[guillsep] (-1,1) -- (5,1)
		(0,0) -- (0,1)	(1,0)--(1,1)	  (2,0)--(2,1)  (3,0)--(3,1)  (4,0)--(4,1);
\node at (0.5,0.5) {$\ha{\gamma}_{e_1^S}$};
\node at (1.5,0.5) {$\ha{\gamma}_{e_2^S}$};
\node at (2.5,0.5) {$\dots$};
\node at (3.5,0.5) {$\ha{\gamma}_{e_P^S}$};
\node at (-0.5,0.5) {$u_{SW}$};
\node at (4.5,0.5) {$u_{SE}$};
\node at (0.05,1) [circle, fill, inner sep=0.5mm] {};
\draw[dotted] (0.05,0)--(0.05,1);
\end{scope}
\end{tikzpicture}
= \gamma_{e_1^S} \otimesdots \gamma_{e_P^S} \in [T_{\infty_{WE},\infty_S}(V(S_1),V(S_2))]_{0,P}
\end{equation*}
where $[x]_{0,P}$ means that the element is placed in the $(0,P)$-component of the direct sum over $\setP^2_\leq$ of \eqref{eq:defcanonicalspaces:doublyinfinitepatterns:FD} since the base point is placed on the left vertical guillotine cut. Another choice of the base point would shift this index.
Similar results hold on the North, West and East side.

Contracting the West and East side with the partition function $Z_R(\MarkovWeight{W}; \cdot)$ seen as an element of $\ca{T}_{P,Q})$ as in Theorem~\ref{theo:partitionfuncguillotop} now provide an element of $\End(V(S_1))^{\otimes P}$ given by
\begin{equation}
\begin{split}
\begin{tikzpicture}[guillpart]
\begin{scope}[xscale=1.6,yscale=1.6];
\fill[guillfill] (0,1) rectangle (6,4);
\draw[guillsep] (0,1) -- (6,1)
			 (0,4) -- (6,4)
			 (1,1) -- (1,4)
			 (5,1) -- (5,4);
\draw[guillsep] (0,2)--(1,2) (0,3)--(1,3) (0,4)--(1,4);
\draw[guillsep] (5,2)--(6,2) (5,3)--(6,3) (5,4)--(6,4);
\node at (0.5,1.5) {$\ha{\gamma}_{e_1^W}$};
\node at (0.5,2.5) {$\vdots$};
\node at (0.5,3.5) {$\ha{\gamma}_{e_Q^W}$};
\node at (5.5,1.5) {$\ha{\gamma}_{e_1^E}$};
\node at (5.5,2.5) {$\vdots$};
\node at (5.5,3.5) {$\ha{\gamma}_{e_Q^E}$};
\node at (3,2.5) {$Z_R(\MarkovWeight{W}_\bullet;\cdot)$};
\node at (1.05,1.) [circle,fill,inner sep=0.5mm] {};
\node at (1.05,4.) [circle,fill,inner sep=0.5mm] {};
\draw[dotted] (1.05,1)--(1.05,4);
\node at (1.5,1) [anchor=north] {$x_{e_1^S}$};
\node at (2.5,1) [anchor=north] {$x_{e_2^S}$};
\node at (3.5,1) [anchor=north] {$\dots$};
\node at (4.5,1) [anchor=north] {$x_{e_P^S}$};
\node at (1.5,4) [anchor=south] {$x_{e_1^N}$};
\node at (2.5,4) [anchor=south] {$x_{e_2^N}$};
\node at (3.5,4) [anchor=south] {$\dots$};
\node at (4.5,4) [anchor=south] {$x_{e_P^N}$};
\end{scope}
\end{tikzpicture}
= 
\sum_{\substack{(x_{e^W_i})_{1\leq i\leq Q}\\ (x_{e^E_i})_{1\leq i\leq Q}} }
Z_R(\MarkovWeight{W}_\bullet;(x_e)_{e\in \partial R}) \prod_{i=1}^Q \gamma_{e^W_i}(x_{e_i^W}) \prod_{j=1}^Q \gamma_{e^W_j}(x_{e_j^W})
\end{split}
\label{eq:strippartitionfuncwithfactor}
\end{equation}
where the indices $x_{e^{a}_i}$ are just here to write the element $Z_R(\MarkovWeight{W}_\bullet;\cdot)\in T_{P,Q}(V_1,V_2)$ as a matrix and not an endomorphism, so that we get an nice identification with the probabilistic formula~\eqref{eq:relationbetweenpartitionfunc}.

Placing the base point on the left vertical guillotine cut places this element of $\End(V(S_1))^{\otimes P}$ in the component $(0,P)\in\setP^2_\leq$ of the direct sum of $\ca{T}_{\infty_{WE},Q}$.

Now, gluing together the South side, the previous equation and the North side with the partition 
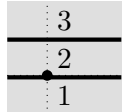
\[
\begin{tikzpicture}[guillpart]
\fill[guillfill] (0,0) rectangle (3,3);
\draw[guillsep] (0,1) -- (3,1)  (0,2)--(3,2);
\node at (1.05,1.05) [circle, fill, inner sep=0.5mm] {};
\draw[dotted] (1.05,0)--(1.05,3) (0,1.05)--(3,1.05);
\node at (1.5,0.5) {$1$};
\node at (1.5,1.5) {$2$};
\node at (1.5,2.5) {$3$};
\end{tikzpicture}
\]
we obtain the expression \eqref{eq:relationbetweenpartitionfunc}. The reader may easily check that any other choice of the base point would just produce a shift in the tensor products irrelevant at the the end.

We finally remark that we made a particular choice of decomposition of the guillotine partition~\eqref{eq:partitionfuncwithfactorbw} but, using generalized associativities, any other decomposition gives the same result.
\end{proof}

\subsubsection{Boundary weights for marginal processes.}

The same mechanism combined with Proposition~\ref{lemma:proba:Markovrestriction} provides the following proposition for restrictions to sub-rectangles: we do not provide the proof since Theorem~\ref{theo:stability} will generalize it much further. We however encourage the reader to try to write this proof as an exercise on the extended $\Guill_2^{(\patterntype{fp}^*)}$-operad.

\begin{prop}\label{prop:factorizedweights:marginal}
Under the same notations as in Proposition~\ref{prop:partitionfuncfactorweights}, let $R'\subset R$ be a sub-rectangle $[P'_1,P'_2]\times [Q'_1,Q'_2]$ with shape $(P',Q')$. The process $(X_e)_{e\in \Edges{R'}}$ is a Markov process with a boundary weight $g_{R'}$ such that, for any $Y\in\ca{T}_{P',Q'}(V(S_1),V(S_2))$, 
\begin{equation}\label{eq:MPSrestrictfromindep}
\sum_{(x_e)_{e\in\Edges{\partial R'}}} g_{R'}((x_e)) Y((x_e)) 
=
\begin{tikzpicture}[guillpart]
\begin{scope}[scale=2.];
\fill[guillfill] (0,0) rectangle (6,5);
\draw[guillsep] (0,1) -- (6,1)
			 (0,4) -- (6,4)
			 (1,0) -- (1,5)
			 (5,0) -- (5,5);
\draw[guillsep] (2,0)--(2,1)	(3,0)--(3,1) (4,0)--(4,1);
\draw[guillsep] (2,4)--(2,5)	(3,4)--(3,5) (4,4)--(4,5);
\draw[guillsep] (0,2)--(1,2) (0,3)--(1,3) (0,4)--(1,4);
\draw[guillsep] (5,2)--(6,2) (5,3)--(6,3) (5,4)--(6,4);
\node at (0.5,0.5) {$u'_{SW}$}; 
\node at (5.5,0.5) {$u'_{SE}$};
\node at (0.5,4.5) {$u'_{NW}$};
\node at (5.5,4.5) {$u'_{NE}$};

\node at (1.5,0.5) {$\ha{\gamma}'_{e_{P'_1+1}^S}$};
\node at (2.5,0.5) {$\ha{\gamma}'_{e_{P'_1+2}^S}$};
\node at (3.5,0.5) {$\dots$};
\node at (4.5,0.5) {$\ha{\gamma}'_{e_{P'_2}^S}$};

\node at (1.5,4.5) {$\ha{\gamma}'_{e_{P'_1+1}^N}$};
\node at (2.5,4.5) {$\ha{\gamma}'_{e_{P'_1+2}^N}$};
\node at (3.5,4.5) {$\dots$};
\node at (4.5,4.5) {$\ha{\gamma}'_{e_{P'_2}^N}$};

\node at (0.5,1.5) {$\ha{\gamma}'_{e_{Q'_1+1}^W}$};
\node at (0.5,2.5) {$\vdots$};
\node at (0.5,3.5) {$\ha{\gamma}'_{e_{Q'_2}^W}$};

\node at (5.5,1.5) {$\ha{\gamma}'_{e_{Q'_1+1}^E}$};
\node at (5.5,2.5) {$\vdots$};
\node at (5.5,3.5) {$\ha{\gamma}'_{e_{Q'_2}^E}$};

\node at (3,2.5) {$Y$};

\node at (1,1) [circle,fill,inner sep=0.5mm] {};
\draw[dotted] (1,0)--(1,5) (0,1)--(6,1);
\end{scope}
\end{tikzpicture}
\end{equation}
with the following boundary elements:
\begin{align}\label{eq:MPSeltsfromindep}
u'_{SW} &= 
\begin{tikzpicture}[guillpart,xscale=2,yscale=1.6]
	\fill[guillfill] (0,0) rectangle (4,4);
	\draw[guillsep] (0,1) -- (4,1) (1,0)--(1,4)
		(2,0) -- (2,1) (3,0)--(3,1) (4,0)--(4,4)
		(0,2) -- (1,2) (0,3)--(1,3) (0,4)--(4,4);
	\node at (2.5,2.5) { $Z_{[P_1,P'_1]\times [Q_1,Q'_1]}$ };
	\node at (0.5,0.5) { $u_{SW}$ };
	\node at (1.5,0.5) { $\ha{\gamma}_{e_1^S}$ };
	\node at (2.5,0.5) { \dots };
	\node at (3.5,0.5) { $\ha{\gamma}_{e_{P_1}^S}$ };
	\node at (0.5,1.5) { $\ha{\gamma}_{e_1^E}$ };
	\node at (0.5,2.5) { \dots };
	\node at (0.5,3.5) { $\ha{\gamma}_{e_{Q_1}^E}$ };
\end{tikzpicture}
&
\ha{\gamma}'_{e_i^S} &=
\begin{tikzpicture}[guillpart,yscale=1.6,xscale=5]
	\fill[guillfill] (0,0) rectangle (1,4);
	\draw[guillsep] (0,0)--(0,4) (1,0)--(1,4) (0,1)--(1,1) (0,4)--(1,4);
	\node at (0.5,2.5) { $Z_{[i-1,i]\times [Q_1,Q'_1]}$ };
	\node at (0.5,0.5) { $\ha{\gamma}_{e_i^S}$ };
\end{tikzpicture}
\end{align}
(where the partition functions have simplified notations to indicate only their domain) and similar rotated guillotine partitions for the other corners and sides.
\end{prop}

\begin{rema}[matricial description of the boundary weights] Writing down carefully the new boundary elements \eqref{eq:MPSeltsfromindep} that appear in the expression of the boundary weight $g_{R'}$ in terms of elements of tensor products of endomorphisms provides the following representation
	\begin{equation}
\begin{split}
g_{R'}((x_e)) = 
\Tr\Big(&
A^{S}_1(x^S_1)\ldots A^{S}_{P'}(x^S_{P'}) U_{SE} 
 A^{E}_1(x^E_1)\ldots A^{S}_{Q'}(x^E_{Q'}) U_{NE} 
\\
& A^{N}_{P'}(x^N_{P'})\ldots A^{S}_1(x^N_1) U_{NW} 
A^{W}_{Q'}(x^W_{Q'})\ldots A^{S}_1(x^W_1) U_{SW} 
\Big)
\end{split}
\end{equation}
where the $A^{a}_i(x)$ and $U_{ab}$ are respectively square and rectangular matrices. Such a representation is called a \emph{matrix product state} and is the object of study of the next Section~\ref{sec:boundaryalgebra}.
\end{rema}

\subsection{Some comments about relations with transfer matrices, cylinder geometry and Hopf algebras}	\label{sec:guill2:cylindershape}		

Most traditional computations in two-dimensional statistical mechanics rely on the transfer matrix, which consists in dividing a rectangular into horizontal (or vertical strips) in order to reduce the study of a 2D Markov process to a 1D Markov process with a larger state space, which corresponds to the sequence of values along a line. This corresponds in our formalism to the sub-operad of guillotine partitions with only horizontal guillotine cuts, i.e. to the $\Guill_1$-operad generated by the $m_{SN}^{P|q_1,q_2}$ for some fixed $P\geq 1$. Within this perspective, it is interesting to remark that it may not be natural to consider boundary conditions described by spaces $(\ca{T}_{p,\infty_S})_{p\in\setN_1}$, since the products on these spaces are described by extensions of the products $m_{WE}$, which have been forgotten in the construction of the transfer matrix.

Traditionally, two types of boundary conditions on the vertical sides are considered to build transfer matrices: factorized or periodic.

\subsubsection{Factorized boundary weights on the vertical side.} The factorized case corresponds exactly to the situation of propositions~\ref{prop:partitionfuncfactorweights} and \ref{prop:factorizedweights:marginal}. Factorized boundary weights may not be the most interesting boundary weights in statistical mechanics models. However, there are frequent in many exact results from statistical mechanics since they induce no external correlations between the successive horizontal lines. 

Indeed, the partition function with factorized weights on the West and East side obtained in~\eqref{eq:strippartitionfuncwithfactor} in the proof of Proposition~\ref{prop:halfplanemoduleactions} combined with Theorem~\ref{theo:partitionfuncguillotop} gives the factorization:
\[
\begin{tikzpicture}[guillpart]
\begin{scope}[xscale=2,yscale=1.5];
\fill[guillfill] (0,1) rectangle (5,4);
\draw[guillsep] (0,1) -- (5,1)
			 (0,4) -- (5,4)
			 (1,1) -- (1,4)
			 (4,1) -- (4,4);
\draw[guillsep] (0,2)--(1,2) (0,3)--(1,3) (0,4)--(1,4);
\draw[guillsep] (4,2)--(5,2) (4,3)--(5,3) (4,4)--(5,4);
\node at (0.5,1.5) {$\ha{\gamma}_{e_1^W}$};
\node at (0.5,2.5) {$\vdots$};
\node at (0.5,3.5) {$\ha{\gamma}_{e_Q^W}$};
\node at (4.5,1.5) {$\ha{\gamma}_{e_1^E}$};
\node at (4.5,2.5) {$\vdots$};
\node at (4.5,3.5) {$\ha{\gamma}_{e_Q^E}$};
\node at (2.5,2.5) {$Z_R(\MarkovWeight{W}_\bullet;\cdot)$};
\node at (1.05,1.) [circle,fill,inner sep=0.5mm] {};
\node at (1.05,4.) [circle,fill,inner sep=0.5mm] {};
\draw[dotted] (1.05,1)--(1.05,4) ;
\end{scope}
\end{tikzpicture}
= \begin{tikzpicture}[guillpart]
\begin{scope}[xscale=1.5,yscale=1.5]
\fill[guillfill] (0,1) rectangle (5,4);
\draw[guillsep] (0,1) -- (5,1)
			(0,2) -- (5,2)
			(0,3) -- (5,3)
			 (0,4) -- (5,4);
 \node at (1.,1.) [circle,fill,inner sep=0.5mm] {};
  \node at (1.,4.) [circle,fill,inner sep=0.5mm] {};
\draw[dotted] (1.,1)--(1.,4);
\node at (2.5,1.5) {$T_1$};
\node at (2.5,2.5) {$\vdots$};
\node at (2.5,3.5) {$T_Q$};
\end{scope}
\end{tikzpicture}
\]
where the r.h.s.~is just the $m_{SN}$-products of $Q$ elements of the strip algebra $\ca{S}_{WE}$ and $T_i$ is the usual transfer matrix of the $i$-th line (with some "diagonal" boundary conditions from the vertical weights $\phi_e$) given by:
\[
T_i = 
\begin{tikzpicture}[guillpart]
\begin{scope}[scale=1.7]
\fill[guillfill] (0,0) rectangle (5,1);
\draw[guillsep] (0,0)--(5,0) (0,1)--(5,1);
\node at (1.05,1.) [circle,fill,inner sep=0.5mm] {};
\node at (1.05,0) [circle,fill,inner sep=0.5mm] {};
\draw[dotted] (1.05,0)--(1.05,1);
\draw[guillsep] (1,0)--(1,1) (4,0)--(4,1) (2,0)--(2,1) (3,0)--(3,1);
\node at (0.5,0.5) {$\ha{\phi}_{e_i^W}$};
\node at (4.5,0.5) {$\ha{\phi}_{e_i^E}$};
\node at (1.5,0.5) {$w_{1,i}$};
\node at (2.5,0.5) {$\dots$};
\node at (3.5,0.5) {$w_{P,i}$};
\end{scope}
\end{tikzpicture}
\in \ca{T}_{\infty_{WE},1}(V(S_1),V(S_2))
\]
and the transfer matrix $T_i$ is trivially identified to an element $\End(V(S_1))^{\otimes P}$.

\subsubsection{The periodic case on vertical sides: the cylinder.}

Periodic boundary conditions correspond to the identification of the opposite vertical sides of rectangles to build a cylinder. There are major technical advantages of this geometry:
\begin{itemize}
\item most often, the (double) thermodynamic limit of a large radius and large length cylinder coincide with the one obtained by a Gibbs measure on the whole $\setZ^2$;
\item it "simplifies" the computations through the study of a one-dimensional system;
\item in some cases, one may use Fourier transform in the circular dimension;
\item it removes the need to specify often-irrelevant boundary conditions on the vertical edges.
\end{itemize}
However, it hides the fundamental ingredients of a two-dimensional geometry, namely the guillotine partitions and the presence of two coherent products which are fundamental for boundary conditions given by matrix product states as described in Section~\ref{sec:boundaryalgebra}. From a more physical perspective on Markov processes, the cylinder geometry is only a computational tool to extract bulk properties in the thermodynamic limit, which should corresponds to invariant boundary weights of finite rectangles (embedded in either large cylinders of the full plane with adequate Gibbs measure), as studied in the next sections.

It is interesting to see however how periodic boundary conditions may fit in our framework. Once again, the colour palette of the $\Guill_2^{(\patterntype{r})}$-operad may be completed to describe periodicity. To this purpose, we may introduce new pattern shapes $\patterntype{cyl}_{NS}^*$, $\patterntype{cyl}_{WE}^*$ and $\patterntype{torus}^*$), associated colours $(p,q,c)\in \setN\times \setN \times \{h,v,hv\}$ where the last element indicates which direction is periodic (horizontal, vertical or both) and the first two elements indicate the sizes of the cylinder. Guillotine partitions of cylinders cut cylinders into either cylinders or rectangles.  Translational invariance along the periodic dimensions requires a pointing in these directions, as already seen for the doubly-infinite planar case of Section~\ref{sec:extendedguillpartII}. 

The canonical spaces for these shapes are trivially given by
\begin{align*}
\ca{T}_{(p,q,h)} &= \End( V(S_1))^{\otimes p}
&
\ca{T}_{(p,q,v)} &= \End( V(S_2))^{\otimes q}
\end{align*} and
$\ca{T}_{(p,q,hv)}= \setR$.
For example we consider a cylinder $[0,p] \times [0,q]$ with shape $(p,q,h)$ (i.e. points $(0,y)$ and $(p,y)$ are identified) cut into two rectangles $[0,r] \times [0,q]$ and $[r,p]\times [0,q]$ with shapes $(r,q)$ and $(p-r,q)$ and with base point $x\in [0,p]$, the associated product is defined by:
\begin{align*}
\ca{T}_{r,q} \otimes \ca{T}_{p-r,q} & \to \ca{T}_{(p,q,h)}
\\
\left[ 
\left(\bigotimes_{i=1}^r A_i\right) \otimes
\left(\bigotimes_{j=1}^q B_j\right)
\right]
\otimes 
\left[ 
\left(\bigotimes_{i=1}^{p-r} A'_i\right) \otimes
\left(\bigotimes_{j=1}^q B'_j\right)
\right]
&
\mapsto
\prod_{j=1}^q \Tr_{V(S_2)}(B_jB'_j) \theta_x\left(\bigotimes_{i=1}^p A''_i\right) 
\end{align*}
where the element $A''_i$ is defined as $A_i$ if $1\leq i\leq r-x$ and $A''_{i-r}$ if $r+1\leq i\leq p$. The translation $\theta_x$ corresponds to a circular shift of the indices in the sequence $(A''_i)$. This type of formulae with traces guarantees the generalized associativity relations~\eqref{eq:guill2:listassoc} adapted to nested guillotine partitions of cylinders and tori.

One checks that, for fixed $p\in\setN_1$, the spaces $(\ca{T}_{p,q,h})_{q\in\setN_1}$ have a $\Guill_1$-structure related  to the gluing of cylinders with the same radius $p$.

\chapter[Matrix product states for boundary weights]{Matrix product states for boundary weights of two-dimensional Markov processes}\label{sec:boundaryalgebra}

As illustrated in propositions~\ref{prop:partitionfuncfactorweights} and \ref{prop:factorizedweights:marginal} and subsequent remarks, the fact that the $\Guill_2^{(\patterntype{r})}$-algebra structure of partition functions can be nicely extended to $\Guill_2^{(\patterntype{fp}^*)}$-algebras in which sides and corners are equipped with algebra and modules structure is very reminiscent to the so-called \emph{matrix product states}, omnipresent for several years in quantum and statistical mechanics, both in theory and in numerics. The present section studies the interplay between such matrix product states, 2D Markov processes and extended $\Guill_2$-operads.

The map of this section is the following: we first introduce matrix product states on their own and progressively reinterpret them in the language of $\Guill_2$-operads and arrive at the final Definition~\ref{def:ROPErep:FD}. We then formulate and prove the main stability Theorem~\ref{theo:stability} for Markov processes with matrix product states as boundary weights. We then explore the concrete computational consequences of this theorem in order to put in evidence the lack of a suitable notion of eigen-objects (values, vectors or matrices) in the framework of operads. We postpone such generalized algebraic definitions of eigen-objects, which are the key results of the paper, to Section~\ref{sec:invariantboundaryelmts}.
 
	\section{Rectangular generalization of matrix product weights: an overview}\label{sec:introROPE}
		\subsection{A quick overview of matrix product representations of functions}

Let $S$ be a finite set. The idea of Matrix Ans\"atze or matrix product states is to represent a function $f: S^n \to \setC$ by a set of matrices $(A(x))_{x\in S}$ acting on some vector space $\ca{H}$ and two vectors $v\in\ca{H}$ and $u\in\ca{H}^*$ such that, for all $(x_1,\ldots,x_n)\in S^n$,
\begin{equation}\label{eq:segmentmatrixansatz}
f(x_1,\ldots,x_n) = \scal{u}{A(x_1)\ldots A(x_n) v}
\end{equation}
Formulated as such, this may not look so interesting since such a representation always exists; this can be seen by taking sufficiently large matrices to encode all the words $(x_1,\ldots,x_n)$. 

This becomes much more interesting when one considers families $(f_N)_{N\in\setN}$ of functions with $f_N: S^N \to \setC$ represented by a \emph{common} and \emph{homogeneous} family of matrices $(A(x))_{x\in S}$. In particular, this allows to build "global" functions on a arbitrary large set $S^N$ out of a finite set of "local" building blocks: the matrices and the two vectors. This problematic of comparing functions on different sets made of the same building blocks is precisely at the heart of the description of boundary weights for Markov process and also at the heart of other various topics in theoretical physics.

A well-known example is in particular the asymmetric simple exclusion process (ASEP) with reservoirs, whose invariant law as as well as other eigenvectors are made of such Matrix Ans\"atze: this is well illustrated in \cite{DEHP,ASEPBrownianExc,LargeDevDensityASEP,BlytheEvans,CrampeRagoucySimon,vanicatZF}. This is also the case in quantum mechanics and solid state physics as illustrated in \cite{MPSreview} and \cite{DMRGMPS} (for numerical approaches). Moreover, in many cases where the models under study exhibit nice algebraic properties, these matrices and vectors can be related to well-known algebraic structures, such as quantum groups or tridiagonal matrices: in particular, this may be of great help for computations of asymptotics of correlation functions.

More generally, one also encounters representations $f(x_1,\ldots,x_n)=\omega(A(x_1)\ldots A(x_n))$ where $\omega$ is a state (such as the trace for example) and the precise complete topological framework is often not specified in practical computations. The main idea, which is useful in practice, is to map a one-dimensional sequence of values $(x_1,\ldots,x_n)$ to an operator $A_1(x_1)\ldots A_n(x_n)$ such that the operator associated to concatenated sequences is the product of the operators associated to each sequence. The extension to cyclic sequences $(x_k)_{k\in\setZ/N\setZ}$ does not change this main idea but uses only tracial states $\omega$. 

In the same way, extensions to rectangular geometries may be made by considering four Matrix Ans\"atze on the four sides of a rectangle with a state made of a trace (since the curve is closed) and elements of bimodules on the corners. The next section is devoted to rigorous operadic definitions to describe the following situation. For any sequences $x=(x_i)$ and $y=(y_j)$ of $p$ elements in a set $S_1$ and any sequences $z=(z_k)$ and $w=(w_l)$ of $q$ elements in a set $S_2$, we wish to think about the four sequences as placed around a rectangle (following the axis orientation) and we wish to represent a function $f$ --~such as a boundary weight $g$ of a two-dimensional Markov process~-- on decorated rectangles as:
\begin{subequations}
\label{eq:firstrectangularMPS}
\begin{align}
\label{eq:boundaryordering}
f(x,y,w,z) &=
f\left( 
\begin{tikzpicture}[scale=0.5,baseline={(current bounding box.center)}]
\draw (0,0) rectangle (4,3);
\node at (0.5,0) [anchor = north] {$x_1$};
\node at (1.5,0) [anchor = north]{$x_2$};
\node at (2.5,0) [anchor = north]{$\ldots$};
\node at (3.5,0) [anchor = north]{$x_p$};
\node at (0.5,3) [anchor = south] {$y_1$};
\node at (1.5,3) [anchor = south]{$y_2$};
\node at (2.5,3) [above]{$\ldots$};
\node at (3.5,3) [above]{$y_p$};
\node at (0,0.5) [anchor = east]{ $w_1$ };
\node at (0,1.5) [anchor = east]{ $\vdots$ };
\node at (0,2.5) [anchor = east]{ $w_q$ };
\node at (4,0.5) [anchor = west]{ $z_1$ };
\node at (4,1.5) [anchor = west]{ $\vdots$ };
\node at (4,2.5) [anchor = west]{ $z_q$ };
\end{tikzpicture}
\right) 
\\
&= 
\Tr_{\ca{H}_W}\left( U_{WS} \mathbf{A}_S(x)
U_{SE} \mathbf{A}_E(z)U_{EN} \mathbf{A}_N(y) U_{NW} \mathbf{A}_W(w) \right)
\end{align}
with the following definition of the operators on each side:
\begin{align}
\mathbf{A}_S(x) &= A_S(x_1)\ldots A_S(x_p)  
&
\mathbf{A}_N(z) &= A_N(y_p)\ldots A_N(y_1)  
\\
\mathbf{A}_E(y) &= A_E(z_1)\ldots A_E(z_q)  
&
\mathbf{A}_W(y) &= A_W(w_q)\ldots A_W(w_1)  
\end{align}
\end{subequations}
where $A_a(x)\in \Hom(\ca{H}_a,\ca{H}_a)$ are linear maps on some vector spaces $\ca{H}_a$ with $a\in\{N,W,S,E\}$ and $U_{ab}\in \Hom(\ca{H}_a,\ca{H}_b)$ are linear maps between these spaces. One remarks that the product above is evaluated counter-clockwise (this will not be necessarily the case in the next sections). 

		\subsection{Algebraic and probabilistic remarks on boundary weights on the discrete space \texorpdfstring{$\setZ^2$}{Z2}.}
		
\subsubsection{Some notations.}
A boundary weight on a non-degenerate rectangle $R$ of size $(p,q)\in\setN^*\times\setN^*$, as  in \eqref{eq:proba:ZboundaryfromZdet}, is a function $g_R$ mapping a sequence $(x_e)_{e\in\partial R}$ with values in $S_1$ and $S_2$ (horizontal and vertical) to $\setR_+$. Any such sequence is associated to an element of $S_1^p\times S_1^p\times S_2^q\times S_2^q$ where the first (resp. second, third, fourth) set corresponds to the elements of the sequence on the South (resp. North, West, East) read from left to right (resp. left to right, bottom to top, bottom to top). In order to keep notations as clear as possible, we will systematically use the graphical notation
\begin{equation}\label{eq:notation:seqonrectboundary}
g\left( 
\begin{tikzpicture}[scale=0.5,baseline={(current bounding box.center)}]
\draw (0,0) rectangle (4,3);
\node at (0.5,0) [anchor = north] {$x_1^S$};
\node at (1.5,0) [anchor = north]{$x_2^S$};
\node at (2.5,0) [anchor = north]{$\ldots$};
\node at (3.5,0) [anchor = north]{$x_p^S$};
\node at (0.5,3) [anchor = south] {$x_1^N$};
\node at (1.5,3) [anchor = south]{$x_2^N$};
\node at (2.5,3) [above]{$\ldots$};
\node at (3.5,3) [above]{$x_p^N$};
\node at (0,0.5) [anchor = east]{ $x_1^W$ };
\node at (0,1.5) [anchor = east]{ $\vdots$ };
\node at (0,2.5) [anchor = east]{ $x_q^W$ };
\node at (4,0.5) [anchor = west]{ $x_1^E$ };
\node at (4,1.5) [anchor = west]{ $\vdots$ };
\node at (4,2.5) [anchor = west]{ $x_q^E$ };
\end{tikzpicture}
\right),
\end{equation}
every time such a geometric interpretation of a sequence $(x_e)_{e\in\partial R}$, seen as an element of $S_1^p\times S_1^p \times S_2^q\times S_2^q$ is used.
This notation was already used in \eqref{eq:boundaryordering} with a different counter-clockwise orientation of edges: the present left-to-right and bottom-to-top ordering will be more suited to guillotine cuts.

In order to describe marginal one-dimensional probability laws, we may also consider the case of segments seen as degenerate rectangles with sizes $(p,0)$ or $(0,q)$. In this case, a boundary weight $g$ is simply a function $S_1^p\to\setR_+$ (left-to-right orientation) or $S_2^q\to\setR_+$ (bottom-to-top orientation) of the configurations along the segments.

\subsubsection{Boundary weights as elements of a dual space for non-degenerate rectangles.}\label{par:dualviewonweights} A first good algebraic interpretation is to see such a boundary weight $g$ on a \emph{non-degenerate} rectangle $R$ of size $(p,q)$ as an element of the dual $\ca{T}_{p,q}(V(S_1),V(S_2))^*$ (that we will write again $\ca{T}_{p,q}^*$ for a short notation). Indeed, \eqref{eq:proba:ZboundaryfromZdet} maps linearly a function $Z_R(\MarkovWeight{W}_\bullet;\cdot)$ of the boundary values to a scalar number $Z_R^{\boundaryweights}(\MarkovWeight{W}_\bullet)$ and the function $Z_R(\MarkovWeight{W}_\bullet;\cdot)$ can be described (using Theorem~\ref{theo:partitionfuncguillotop}) as an element of $\ca{T}_{p,q}$. 

In the present case with finite spaces $S_1$ and $S_2$, the spaces $\ca{T}_{p,q}$ are finite-dimensional and have a canonical basis labelled by the elements of the sets $S_i$ and thus could be arbitrarily identified with their duals $\ca{T}_{p,q}^*$. However, this is not a good idea, neither from a probabilistic point of view, nor from an algebraic one. On the probabilistic side, boundary weights most often are used to describe marginal laws of a system inside a larger one by resuming contributions outside the small system. On the algebraic side, already in dimension one for associative algebras, states are often obtained from representations (see for example the GNS representation of $C^*$-algebras). Following this idea, boundary weights should be associated to full plane algebras $\ca{T}_{\patterntype{fp}^*}$ and built out of the boundary spaces $\ca{T}_{\patterntype{fp}^*\setminus\patterntype{r}}$: such an idea induces naturally "matrix product" representations of the weights, as seen below in a more formal way.

Section~\ref{sec:operad} already provides a description of interesting elements of the dual spaces $\ca{A}_{p,q}^*$ of a $\Guill_2(\patterntype{r})$-algebra, that we will exploit. Assuming that we are given a full plane extension $(\ca{A}_{p,q})_{(p,q)\in\PatternShapes(\patterntype{fp}^*)}$ of the algebra $\ca{A}_{\PatternShapes(\patterntype{r})}$ with the requirement 
\begin{equation}\label{eq:requirement:ROPErepfordual}
\ca{A}_{\infty_{WE},\infty_{SN}}= \setK
\end{equation} there is then a natural inclusion:
\[
\iota : \ca{A}_{p,\infty_S}\otimes \ca{A}_{p,\infty_N} \otimes \ca{A}_{\infty_W,q} \otimes \ca{A}_{\infty_E,q} \otimes
\left(
\bigotimes_{(a,b)\in\{S,N\}\times \{W,E\} }
\ca{A}_{\infty_b,\infty_a } \right) \to \ca{A}_{p,q}^*
\]
given by:
\begin{equation}\label{eq:dualityfromboundarystruct}
\iota\left(A_S^{(p)}\otimes A_N^{(p)}\otimes A_W^{(q)}\otimes A_E^{(q)}\otimes \bigotimes_{a,b} U_{ab}\right)(w) = 
\begin{tikzpicture}[guillpart,scale=2]
	\fill[guillfill] (0,0) rectangle (3,3) ;
	\draw[guillsep] (1,0)--(1,3) (2,0)--(2,3) (0,1)--(3,1) (0,2)--(3,2);
	\node at (1.1,1.15) [circle,fill,inner sep=0.5mm] {};
	\draw[dashed] (0,1.15) -- (3,1.15)
		(1.1,0) -- (1.1,3);
	\node at (0.5,0.5) { $U_{SW}$ };
	\node at (1.5,0.5) { $A^{(p)}_{S}$ };
	\node at (2.5,0.5) { $U_{SE}$ };
	\node at (0.5,1.5) { $A^{(p)}_{W}$ };
	\node at (1.5,1.5) { $w$ };
	\node at (2.5,1.5) { $A^{(q)}_{E}$ };
	\node at (0.5,2.5) { $U_{NW}$ };
	\node at (1.5,2.5) { $A^{(p)}_{N}$ };
	\node at (2.5,2.5) { $U_{NE}$ };
\end{tikzpicture}\in\setK
\end{equation}
In particular, this remark invites one to consider boundary weights $g \in \ca{T}_{p,q}^*$ as given by inclusion of elements of a suitable boundary representation $\ca{A}_{\PatternShapes(\patterntype{fp}^*)\setminus \PatternShapes(\patterntype{r})}$ of $\ca{T}_{\PatternShapes(\patterntype{r})}$ into the dual space of $\ca{T}_{p,q}^*$, especially if the elements $A^{(p)}_S$, $A^{(p)}_N$, $A^{(q)}_W$ and $A^{(q)}_E$ above are given by products of $p$ or $q$ elements $A^{(1)}_a \in \ca{A}_{1,\infty_a}$ (or $\ca{A}_{\infty_a,1}$ on the vertical sides), related by $\Guill_1$-product structure to the sequence of boundary elements $(x_e)_{e\in\partial R}$.

\subsubsection{Boundary weights as elements of a dual space for segments (degenerate rectangles).} 

The configuration on a horizontal segment of length $p$ is simply a sequence of $p$ elements of $S_1$, i.e. an element of $S_1^p$. A boundary weight on a segment is a function $g: S_1^p\to \setR$, which does \emph{not} correspond to a boundary weight on a non-degenerate rectangle $(p,q)$ with $q=0$: this latter case would correspond to a function $g:S_1^p\times S_1^p \to \setR$ which would rather corresponds to a \emph{slit} in the plane, each copy of $S_1$ corresponding to a side of the slit.

Describing a segment boundary weight $g:S_1^p \to \setR$ as a dual element remains valid provided one considers the commutative spaces $\ca{T}_{p,q}$ when $p$ or $q$ vanishes. Indeed, $\ca{T}_{p,0}= \Diag_{|S_1|}(\setR)^{\otimes p}$ is canonically identified with the tensor $p$-power of $\setR^{|S_1|}$, whose dual is identified to itself and then identified to the set of functions $S_1^p\to\setR$ through the canonical basis. 

There is thus a deep connection between the transverse $0$-dimensionality of segments, the commutativity induced on the gluing of such segments along this direction, the restriction to diagonal matrices instead of full matrices, the duality with weights $g : S_i^k \to \setR$ (and not $S_i^k\times S_i^k\to\setR$) due to the null width.

From the point of view of extended $\Guill_2^{(\patterntype{fp}^*)}$-algebra, \eqref{eq:dualityfromboundarystruct} should be replaced in the case of horizontal segments by:
\[
	\iota\left(A_S\otimes A_N\otimes \bigotimes_{a,b} U_{ab}\right)(D_h^{(p)}) = 
	\begin{tikzpicture}[guillpart,scale=2]
		\fill[guillfill] (0,0) rectangle (3,2) ;
		\draw[guillsep] (1,0)--(1,2) (2,0)--(2,2) (0,1)--(3,1) ;
		\node at (1.1,1.15) [circle,fill,inner sep=0.5mm] {};
		\draw[dashed] (1.1,0)--(1.1,2) (0,1.15)--(3,1.15);
		\node at (0.5,0.5) { $U_{SW}$ };
		\node at (1.5,0.5) { $A^{(p)}_{S}$ };
		\node at (2.5,0.5) { $U_{SE}$ };
		\node at (1.5,1.) { $D_h^{(p)}$ };
		\node at (0.5,1.5) { $U_{NW}$ };
		\node at (1.5,1.5) { $A^{(p)}_{N}$ };
		\node at (2.5,1.5) { $U_{NE}$ };
	\end{tikzpicture}\in\setK
\]

		\section{Guillotine framework for rectangular boundaries}
	
\subsection{Rectangular Operator Product Environments (ROPE)}	

\subsubsection{Definition of ROPE as particular \texorpdfstring{$\Guill_2$}{Guill2}-algebras}
A more direct operadic approach based on the previous Section~\ref{sec:operad} and more adapted to the operad $\Guill_2$ is chosen here. It consists into guillotine cuts of the boundary of the rectangles, which cut opposite sides at the same coordinate.

Instead of going through all the same steps as in Section~\ref{sec:operad},
a clever approach, which will acquire a deeper meaning in the next Section~\ref{sec:invariantboundaryelmts} relies on the introduction of the trivial algebra $(\setK_{p,q})_{(p,q)\in\PatternShapes(\patterntype{r})}$ over the guillotine operad of Section~\ref{par:trivialoperad}, for which all the spaces are equal to a field $\setK$ and all the products are equal to the commutative product of the field $\setK$.

\begin{defi}[rectangular operator product environment (ROPE)]
\label{def:ROPE}
We consider here \emph{both} cases $(\setP,\setL)=(\setZ,\setN)$ (discrete space) and $(\setR,\setR_+)$ (continuous space).
A rectangular operator product environment (ROPE) over a field $\setK$ is an algebra $\ca{B}_{\PatternShapes(\patterntype{fp}^*)}=(\ca{B}_{p,q})_{(p,q)\in\PatternShapes(\patterntype{fp}^*)}$  over the extended operad $\Guill_2^{(\patterntype{fp}^*)}$ such that:
\begin{enumerate}[(i)]
\item the sub-algebra $\ca{B}_{\PatternShapes(\patterntype{r})}$ is the trivial $\Guill_2^{(\patterntype{r})}$-algebra with $\ca{B}_{p,q}=\setK$ for all $(p,q)\in\PatternShapes(\patterntype{r})$;
\item  the action of $\setK_{\PatternShapes(\patterntype{r})}$ on the boundary algebras $\ca{B}_{p,\infty_S}$, $\ca{B}_{p,\infty_N}$, $\ca{B}_{\infty_W,q}$ and $\ca{B}_{\infty_E,q}$ with $p,q\in\setL^*$
with $(p,q)\notin \PatternShapes(\patterntype{r})$ 
is given by scalar multiplication;
\item $\ca{B}_{\infty_{SN},\infty_{WE}}=\setK$.
\end{enumerate} 
\end{defi}

\subsubsection{Remark on the name.} We may have called it \emph{rectangular matrix product structure} to stick to the idea of \emph{matrix product states} but, in practice, the matrices that will appear are most often infinite dimensional and we prefer the word \emph{operator}. In the case of segment state spaces and continuous space where products are given by integration of kernels, this is even more justified. Moreover, the switch from "structure" to "environment" is made to make the acronym "ROPE", which sounds well adapted to a closed one-dimensional object such as the boundary of a rectangle or a two-dimensional domain.

\subsubsection{Additional graphical notations.} For all segments on the sides, \begin{tikzpicture}[scale=0.3]
\draw (1,0)--(0,0)--(0,1);
\end{tikzpicture}-shapes (and rotated ones), \begin{tikzpicture}[scale=0.3]
\draw (0,0.8)--(0,0)--(1,0)--(1,0.8);
\end{tikzpicture}-shapes (and rotated ones) and the full rectangle, we introduce the following shorthand notations (we display only examples only for one side and one or two corners):
\begin{subequations}
\label{eq:ROPEcorresp}
\begin{equation}
\label{eq:ROPEcorresp1}
\begin{tikzpicture}[guillpart,xscale=2,yscale=1.5]
\draw[guillsep] (0,1)--(4,1);
\draw[guillsep,dotted] (0,0)--(4,0);
\draw[guillsep] (1,0)--(1,1)  (2,0)--(2,1) (3,0)--(3,1) (0,0)--(0,1) (4,0) -- (4,1);
\node at (0.5,0.5) { $A_S^{(1)}$ };
\node at (1.5,0.5) { $A_S^{(2)}$ };
\node at (2.5,0.5) { $\ldots$ };
\node at (3.5,0.5) { $A_S^{(p)}$ };
\end{tikzpicture}
:=
\begin{tikzpicture}[guillpart,xscale=2,yscale=1.5]
\fill[guillfill] (0,0) rectangle (4,1);
\draw[guillsep] (0,1)--(4,1);
\draw[guillsep] (1,0)--(1,1)  (2,0)--(2,1) (3,0)--(3,1) (4,0)--(4,1) (0,0)--(0,1);
\node at (0.5,0.5) { $A_S^{(1)}$ };
\node at (1.5,0.5) { $A_S^{(2)}$ };
\node at (2.5,0.5) { $\ldots$ };
\node at (3.5,0.5) { $A_S^{(p)}$ };
\end{tikzpicture}
\end{equation}
\begin{equation}
\label{eq:ROPEcorresp2}
\begin{tikzpicture}[guillpart,xscale=2,yscale=1.5]
\draw[guillsep] (0,1)--(5,1);
\draw[guillsep,dotted] (0,1)--(0,0)--(5,0)--(5,1);
\draw[guillsep] (1,0)--(1,1)  (2,0)--(2,1) (3,0)--(3,1) (4,0)--(4,1);
\node at (0.5,0.5) { $U_{SW}$ };
\node at (1.5,0.5) { $A_S^{(1)}$ };
\node at (2.5,0.5) { $\ldots$ };
\node at (3.5,0.5) { $A_S^{(p)}$ };
\node at (4.5,0.5) { $U_{SE}$ };
\end{tikzpicture}
:=
\begin{tikzpicture}[guillpart,xscale=2,yscale=1.5]
\fill[guillfill] (0,0) rectangle (5,1);
\draw[guillsep] (0,1)--(5,1);
\draw[guillsep] (1,0)--(1,1)  (2,0)--(2,1) (3,0)--(3,1) (4,0)--(4,1);
\node at (0.5,0.5) { $U_{SW}$ };
\node at (1.5,0.5) { $A_S^{(1)}$ };
\node at (2.5,0.5) { $\ldots$ };
\node at (3.5,0.5) { $A_S^{(p)}$ };
\node at (4.5,0.5) { $U_{SE}$ };
  \node at (1.05,1.) [circle,fill,inner sep=0.5mm] {};
 \draw[dotted] (1.05,0)--(1.05,1);
\end{tikzpicture}
\end{equation}
\begin{equation}
\label{eq:ROPEcorresp3}
\begin{tikzpicture}[guillpart,xscale=2,yscale=1.5]
\draw[guillsep, dotted] (0,4)--(0,0)--(4,0);
\draw[guillsep] (1,4)--(1,1) -- (4,1);
\draw[guillsep] 	(1,0)--(1,1)
				(2,0)--(2,1)
				(3,0)--(3,1)
				(4,0)--(4,1);
\draw[guillsep] 	(0,1)--(1,1)
				(0,2)--(1,2)
				(0,3)--(1,3)
				(0,4)--(1,4);
\node at (0.5,0.5) { $U_{SW}$ };
\node at (1.5,0.5) { $A_S^{(1)}$ };
\node at (2.5,0.5) { $\ldots$ };
\node at (3.5,0.5) { $A_S^{(p)}$ };
\node at (0.5,1.5) { $A_W^{(1)}$ };
\node at (0.5,2.5) { $\vdots$ };
\node at (0.5,3.5) { $A_W^{(q)}$ };
\end{tikzpicture}
:=
\begin{tikzpicture}[guillpart,xscale=2,yscale=1.5]
\fill[guillfill] (0,0) rectangle (4,4);
\draw[guillsep] (1,1) rectangle (4,4);
\draw[guillsep] 	(1,0)--(1,1)
				(2,0)--(2,1)
				(3,0)--(3,1)
				(4,0)--(4,1);
\draw[guillsep] 	(0,1)--(1,1)
				(0,2)--(1,2)
				(0,3)--(1,3)
				(0,4)--(1,4);
\node at (0.5,0.5) { $U_{SW}$ };
\node at (1.5,0.5) { $A_S^{(1)}$ };
\node at (2.5,0.5) { $\ldots$ };
\node at (3.5,0.5) { $A_S^{(p)}$ };
\node at (0.5,1.5) { $A_W^{(1)}$ };
\node at (0.5,2.5) { $\vdots$ };
\node at (0.5,3.5) { $A_W^{(q)}$ };
\node at (2.5,2.5) {$1_\setK^{(p,q)}$};
\end{tikzpicture}
\end{equation}
\begin{equation}
\label{eq:ROPEcorresp4}
\begin{tikzpicture}[guillpart,xscale=2,yscale=1.5]
\draw[guillsep, dotted] (0,4)--(0,0)--(5,0)-- (5,4);
\draw[guillsep] (1,4)--(1,1) -- (4,1)-- (4,4);
\draw[guillsep] 	(1,0)--(1,1)
				(2,0)--(2,1)
				(3,0)--(3,1)
				(4,0)--(4,1);
\draw[guillsep] 	(0,1)--(1,1)
				(0,2)--(1,2)
				(0,3)--(1,3)
				(0,4)--(1,4);
\draw[guillsep] 	(5,1)--(4,1)
				(5,2)--(4,2)
				(5,3)--(4,3)
				(5,4)--(4,4);
\node at (0.5,0.5) { $U_{SW}$ };
\node at (4.5,0.5) { $U_{SE}$ };
\node at (1.5,0.5) { $A_S^{(1)}$ };
\node at (2.5,0.5) { $\ldots$ };
\node at (3.5,0.5) { $A_S^{(p)}$ };
\node at (0.5,1.5) { $A_W^{(1)}$ };
\node at (0.5,2.5) { $\vdots$ };
\node at (0.5,3.5) { $A_W^{(q)}$ };
\node at (4.5,1.5) { $A_E^{(1)}$ };
\node at (4.5,2.5) { $\vdots$ };
\node at (4.5,3.5) { $A_E^{(q)}$ };
\end{tikzpicture}
:=
\begin{tikzpicture}[guillpart,xscale=2,yscale=1.5]
\fill[guillfill] (0,0) rectangle (5,4);
\draw[guillsep] (1,1) rectangle (4,4);
\draw[guillsep] 	(1,0)--(1,1)
				(2,0)--(2,1)
				(3,0)--(3,1)
				(4,0)--(4,1);
\draw[guillsep] 	(0,1)--(1,1)
				(0,2)--(1,2)
				(0,3)--(1,3)
				(0,4)--(1,4);
\draw[guillsep] 	(4,1)--(5,1)
				(4,2)--(5,2)
				(4,3)--(5,3)
				(4,4)--(5,4);
\node at (0.5,0.5) { $U_{SW}$ };
\node at (4.5,0.5) { $U_{SE}$ };
\node at (1.5,0.5) { $A_S^{(1)}$ };
\node at (2.5,0.5) { $\ldots$ };
\node at (3.5,0.5) { $A_S^{(p)}$ };
\node at (0.5,1.5) { $A_W^{(1)}$ };
\node at (0.5,2.5) { $\vdots$ };
\node at (0.5,3.5) { $A_W^{(q)}$ };
\node at (4.5,1.5) { $A_E^{(1)}$ };
\node at (4.5,2.5) { $\vdots$ };
\node at (4.5,3.5) { $A_E^{(q)}$ };
\node at (2.5,2.5) {$1_\setK^{(p,q)}$};
 \node at (1.05,4.) [circle,fill,inner sep=0.5mm] {};
 \draw[dotted] (1.05,0)--(1.05,4);
\end{tikzpicture}
\end{equation}
\begin{equation}
\label{eq:ROPEcorresp5}
\begin{tikzpicture}[guillpart,xscale=2,yscale=1.5]
\draw[guillsep, dotted] (0,0) rectangle (6,5);
\draw[guillsep] (1,1) rectangle (5,4);
\draw[guillsep] 	(1,0)--(1,1)
				(2,0)--(2,1)
				(3,0)--(3,1)
				(4,0)--(4,1)
				(5,0)--(5,1);
\draw[guillsep] 	(1,5)--(1,4)
				(2,5)--(2,4)
				(3,5)--(3,4)
				(4,5)--(4,4)
				(5,5)--(5,4);
\draw[guillsep] 	(0,1)--(1,1)
				(0,2)--(1,2)
				(0,3)--(1,3)
				(0,4)--(1,4);
\draw[guillsep] 	(5,1)--(6,1)
				(5,2)--(6,2)
				(5,3)--(6,3)
				(5,4)--(6,4);
\node at (0.5,0.5) { $U_{SW}$ };
\node at (5.5,0.5) { $U_{SE}$ };
\node at (0.5,4.5) { $U_{NW}$ };
\node at (5.5,4.5) { $U_{NE}$ };
\node at (1.5,0.5) { $A_S^{(1)}$ };
\node at (2.5,0.5) { $A_S^{(2)}$ };
\node at (3.5,0.5) { $\ldots$ };
\node at (4.5,0.5) { $A_S^{(p)}$ };
\node at (1.5,4.5) { $A_N^{(1)}$ };
\node at (2.5,4.5) { $A_N^{(2)}$ };
\node at (3.5,4.5) { $\ldots$ };
\node at (4.5,4.5) { $A_N^{(p)}$ };
\node at (0.5,1.5) { $A_W^{(1)}$ };
\node at (0.5,2.5) { $\vdots$ };
\node at (0.5,3.5) { $A_W^{(q)}$ };
\node at (5.5,1.5) { $A_E^{(1)}$ };
\node at (5.5,2.5) { $\vdots$ };
\node at (5.5,3.5) { $A_E^{(q)}$ };
\end{tikzpicture}
:=
\begin{tikzpicture}[guillpart,xscale=2,yscale=1.5]
\fill[guillfill] (0,0) rectangle (6,5);
\draw[guillsep] (1,1) rectangle (5,4);
\draw[guillsep] 	(1,0)--(1,1)
				(2,0)--(2,1)
				(3,0)--(3,1)
				(4,0)--(4,1)
				(5,0)--(5,1);
\draw[guillsep] 	(1,5)--(1,4)
				(2,5)--(2,4)
				(3,5)--(3,4)
				(4,5)--(4,4)
				(5,5)--(5,4);
\draw[guillsep] 	(0,1)--(1,1)
				(0,2)--(1,2)
				(0,3)--(1,3)
				(0,4)--(1,4);
\draw[guillsep] 	(5,1)--(6,1)
				(5,2)--(6,2)
				(5,3)--(6,3)
				(5,4)--(6,4);
\node at (0.5,0.5) { $U_{WS}$ };
\node at (5.5,0.5) { $U_{SE}$ };
\node at (0.5,4.5) { $U_{NW}$ };
\node at (5.5,4.5) { $U_{EN}$ };
\node at (1.5,0.5) { $A_S^{(1)}$ };
\node at (2.5,0.5) { $A_S^{(2)}$ };
\node at (3.5,0.5) { $\ldots$ };
\node at (4.5,0.5) { $A_S^{(p)}$ };
\node at (1.5,4.5) { $A_N^{(1)}$ };
\node at (2.5,4.5) { $A_N^{(2)}$ };
\node at (3.5,4.5) { $\ldots$ };
\node at (4.5,4.5) { $A_N^{(p)}$ };
\node at (0.5,1.5) { $A_W^{(1)}$ };
\node at (0.5,2.5) { $\vdots$ };
\node at (0.5,3.5) { $A_W^{(q)}$ };
\node at (5.5,1.5) { $A_E^{(1)}$ };
\node at (5.5,2.5) { $\vdots$ };
\node at (5.5,3.5) { $A_E^{(q)}$ };
\node at (3,2.5) {$1_\setK^{(p,q)}$};
 \node at (1.05,3.95) [circle,fill,inner sep=0.5mm] {};
 \draw[dotted] (1.05,0)--(1.05,5);
 \draw[dotted] (0,3.95)--(6,3.95);
\end{tikzpicture}
\end{equation}
\end{subequations}
where the notation $1^{(p,q)}_\setK$ stands for $1_\setK$ seen as an element of $\ca{B}_{p,q}=\setK$ for a correct colour compatibility.

\subsubsection{Relation with matrix product states.}
Before exploring further the properties of rectangular operator product environment, we formulate, as an example, the relations of these objects with traditional matrix product states or matrix Ans\"atze.

\begin{prop}[relation with matrix product states]\label{prop:linkwithMPS}
For any $a\in \{S,N,W,E\}$, let $\ca{H}_a$ be a real Hilbert space (and we identify it with its dual). We define an associated ROPE through
\begingroup
\allowdisplaybreaks
\begin{align*}
\ca{B}_{p,\infty_S} &= \End(\ca{H}_S) &
\ca{B}_{p,\infty_N} &= \End(\ca{H}_N) 
\\
\ca{B}_{\infty_W,q} &= \End(\ca{H}_W) &
\ca{B}_{\infty_E,q} &= \End(\ca{H}_E) 
\\
\ca{B}_{\infty_W,\infty_S} &= \ca{H}_W\otimes \ca{H}_{S} \simeq \Hom(\ca{H}_W,\ca{H}_S)
&
\ca{B}_{\infty_E,\infty_S} &= \ca{H}_E\otimes \ca{H}_{S} \simeq \Hom(\ca{H}_E,\ca{H}_S)
\\
\ca{B}_{\infty_W,\infty_S} &= \ca{H}_W\otimes \ca{H}_{N} \simeq \Hom(\ca{H}_W,\ca{H}_N)
&
\ca{B}_{\infty_E,\infty_S} &= \ca{H}_E\otimes \ca{H}_{N} \simeq \Hom(\ca{H}_E,\ca{H}_N)
\\
\ca{B}_{p,\infty_{SN}} &= \End({\ca{H}_S})\otimes \End(\ca{H}_N)
&
\ca{B}_{\infty_{WE},q} &= \End({\ca{H}_W})\otimes \End(\ca{H}_E)
\\
\ca{B}_{\infty_b,\infty_{SN}} &= \ca{H}_S \otimes \ca{H}_N \simeq \Hom(\ca{H}_N,\ca{H}_S)
&
\ca{B}_{\infty_{WE},\infty_a} &= \ca{H}_W \otimes \ca{H}_E \simeq \Hom(\ca{H}_E,\ca{H}_W)
\\
\ca{B}_{\infty_b,\infty_a} &= \ca{H}_b\otimes \ca{H}_a \simeq \Hom(\ca{H}_b,\ca{H}_a)
&
\ca{B}_{\infty_{WE},\infty_{SN}} &= \setR
\end{align*}
with the usual product of endomorphism, left actions on East and North and right actions on West and South. The tensor products of Hilbert spaces may completed to Hilbert spaces if needed. The following correspondence holds between graphical representations~\eqref{eq:ROPEcorresp} and matrix product states:
\begin{align*}
\eqref{eq:ROPEcorresp1} &= A_S^{(1)} A_S^{(2)} \ldots A_S^{(p)}
\\
\eqref{eq:ROPEcorresp2} &= U_{SW}^r A_S^{(1)} A_S^{(2)} \ldots A_S^{(p)} U_{SE}
\\
\eqref{eq:ROPEcorresp3} &= (A_W^{(1)} \ldots A_W^{(q)} )^t U_{SW}^r A_S^{(1)} \ldots A_S^{(p)}
\\
\eqref{eq:ROPEcorresp4} &= (A_W^{(1)} \ldots A_W^{(q)} )^t U_{SW}^r A_S^{(1)} \ldots A_S^{(p)} U_{SE} A_E^{(1)}\ldots A_E^{(q)}
\\
\eqref{eq:ROPEcorresp5} &= \Tr_{\ca{H}_W}\left( (A_W^{(1)} \ldots A_W^{(q)} )^t U_{SW}^r A_S^{(1)} \ldots A_S^{(p)} U_{SE} A_E^{(1)}\ldots A_E^{(q)} U_{NE}^r (A_N^{(1)}\ldots A_N^{(q)})^t U_{NW} \right)
\end{align*}
\endgroup
where transpositions appear on the West and North sides (from axis orientation to counter-clockwise one) and scalar products between identification of $\ca{H}_a$ with its dual and the products between $U_{ab}$ and $U_{bc}$ (the exponent $\bullet^r$ reflects the passage from left to right vectors).
\end{prop}
The proof is left to the reader as an exercise with the definition of a ROPE.

\subsubsection{Bimodule properties on the corners.} Definition~\ref{def:ROPE} \emph{already} includes some of the associativities mentioned previously : $\Guill_1$-algebras on sides, "modules" on corners, etc., seen as sub-operads of the $\Guill_2^{(\patterntype{fp}^*)}$ operad. In particular, the bimodule property on a corner is a consequence of the Eckmann-Hilton property already discussed in Section~\ref{sec:eckmanhilton}, by using the following interchange relation~\eqref{eq:guill2:interchangeassoc} on the following object
\begin{align*}
\begin{tikzpicture}[guillpart,xscale=1.8,yscale=1.3]
%\fill[guillfill] (0,0)--(2,0)--(2,1)--(1,1)--(1,2)--(0,2)--(0,0);
\draw[guillsep] (0,2)--(1,2)--(1,0);
\draw[guillsep] (2,0)--(2,1)--(0,1);
\draw[guillsep,dotted] (2,0)--(0,0)--(0,2);
\node at (0.5,1.5) { $A_W$ };
\node at (1.5,.5) { $A_S$ };
\node at (.5,.5) {$U_{SW}$ };
\end{tikzpicture}
 := 
 \begin{tikzpicture}[guillpart,xscale=1.8,yscale=1.3]
\fill[guillfill] (0,0)--(2,0)--(2,2)--(0,2)--(0,0);
\draw[guillsep] (0,2)--(2,2)--(2,0);
\draw[guillsep] (2,1)--(0,1);
\draw[guillsep] (1,2)--(1,0);
\node at (0.5,1.5) { $A_W$ };
\node at (1.5,.5) { $A_S$ };
\node at (.5,.5) {$U_{SW}$ };
\node at (1.5,1.5) {$1$ };
\end{tikzpicture}
\end{align*}

\subsubsection{A remark on pointed products and double-infinite shapes.}

Definition~\ref{def:ROPE} of a ROPE already implied the bimodule property of the corner element. One checks that it has a second consequence: the use of a scalar multiplication for the action of $\setK_{p,q}$ on a boundary element makes irrelevant the use of pointed versions of guillotine partition as in Definition~\ref{def:pointedguillotine} as stated in the following proposition.

\begin{prop}\label{prop:ROPE:irrelevantpointing}
Following the notations of Definition~\ref{def:ROPE}, for any $x\in \setZ$, $q\in\setN^*$, any $u\in \ca{B}_{\infty_W,q}$, any $v\in\ca{B}_{\infty_E,q}$ and any $p\in\setN$,
\begin{align*}
	\begin{tikzpicture}[scale=0.5,baseline={(current bounding box.center)}]
	\fill[guillfill] (0,0)--(5,0)--(5,2)--(0,2)--(0,0);
	\draw[guillsep] (0,0)--(5,0);
	\draw[guillsep] (0,2)--(5,2);
	\draw[guillsep] (2,0)--(2,2);
	\draw[guillsep] (3,0)--(3,2);
	\node at (1.33,1) {$u$};
	\node at (2.5,1) {$1_p$};
	\node at (4,1) {$v$};
	\node at (0.66,0) [circle, fill, inner sep=0.5mm] {};
	\node at (0.66,2) [circle, fill, inner sep=0.5mm] {};
	\draw [dotted] (0.66,0) -- (0.66,2);
	\draw [->] (2,2.3) -- node [midway, above] {$x$} (0.66,2.3);
	\end{tikzpicture}
	&= 
		\begin{tikzpicture}[scale=0.5,baseline={(current bounding box.center)}]
	\fill[guillfill] (0,0)--(5,0)--(5,2)--(0,2)--(0,0);
	\draw[guillsep] (0,0)--(5,0);
	\draw[guillsep] (0,2)--(5,2);
	\draw[guillsep] (3,0)--(3,2);
	\node at (1.33,1) {$u$};
	\node at (4,1) {$v$};
	\node at (0.66,0) [circle, fill, inner sep=0.5mm] {};
	\node at (0.66,2) [circle, fill, inner sep=0.5mm] {};
	\draw [dotted] (0.66,0) -- (0.66,2);
	\draw [->] (3,2.3) -- node [midway, above] {$x-p$} (0.66,2.3);
	\end{tikzpicture}
	\\
&=		\begin{tikzpicture}[scale=0.5,baseline={(current bounding box.center)}]
	\fill[guillfill] (0,0)--(5,0)--(5,2)--(0,2)--(0,0);
	\draw[guillsep] (0,0)--(5,0);
	\draw[guillsep] (0,2)--(5,2);
	\draw[guillsep] (2,0)--(2,2);
	\node at (1.33,1) {$u$};
	\node at (4,1) {$v$};
	\node at (0.66,0) [circle, fill, inner sep=0.5mm] {};
	\node at (0.66,2) [circle, fill, inner sep=0.5mm] {};
	\draw [dotted] (0.66,0) -- (0.66,2);
	\draw [->] (2,2.3) -- node [midway, above] {$x$} (0.66,2.3);
	\end{tikzpicture}
\end{align*}
where $1_p$ stands for $1\in\setK$ seen as the element of $\setK_{p,q}$. Similar equations hold for vertical strips, half-planes and the full plane.
\end{prop}
\begin{proof}
The two equalities are obtained by using the associativity of pointed horizontal guillotine partitions to obtain the two decompositions of the left hand-side. We obtain for the first line:
\begin{align*}
	\begin{tikzpicture}[scale=0.5,baseline={(current bounding box.center)}]
	\fill[guillfill] (0,0)--(5,0)--(5,2)--(0,2)--(0,0);
	\draw[guillsep] (0,0)--(5,0);
	\draw[guillsep] (0,2)--(5,2);
	\draw[guillsep] (2,0)--(2,2);
	\draw[guillsep] (3,0)--(3,2);
	\node at (1.33,1) {$u$};
	\node at (2.5,1) {$1_p$};
	\node at (4,1) {$v$};
	\node at (0.66,0) [circle, fill, inner sep=0.5mm] {};
	\node at (0.66,2) [circle, fill, inner sep=0.5mm] {};
	\draw [dotted] (0.66,0) -- (0.66,2);
	\draw [->] (2,2.3) -- node [midway, above] {$x$} (0.66,2.3);
	\end{tikzpicture}
	=
	\begin{tikzpicture}[scale=0.5,baseline={(current bounding box.center)}]
	\fill[guillfill] (0,0)--(5,0)--(5,2)--(0,2)--(0,0);
	\draw[guillsep] (0,0)--(5,0);
	\draw[guillsep] (0,2)--(5,2);
	\draw[guillsep] (3,0)--(3,2);
	\node at (1.33,1) {$1$};
	\node at (4,1) {$2$};
	\node at (0.66,0) [circle, fill, inner sep=0.5mm] {};
	\node at (0.66,2) [circle, fill, inner sep=0.5mm] {};
	\draw [dotted] (0.66,0) -- (0.66,2);
	\draw [->] (3,2.3) -- node [midway, above] {$x-p$} (0.66,2.3);
	\end{tikzpicture}\left(	   
	   \begin{tikzpicture}[scale=0.5,baseline={(current bounding box.center)}]
	\fill[guillfill] (0,0) rectangle (3,2);
	\draw[guillsep] (0,0)--(3,0);
	\draw[guillsep] (0,2)--(3,2);
	\draw[guillsep] (3,0)--(3,2);
	\draw[guillsep] (2,0)--(2,2);
	\node at (1.,1) {$u$};
	\node at (2.5,1) {$1_p$};
	\end{tikzpicture}
	   ,
	   v
	\right)
	=\begin{tikzpicture}[scale=0.5,baseline={(current bounding box.center)}]
	\fill[guillfill] (0,0)--(5,0)--(5,2)--(0,2)--(0,0);
	\draw[guillsep] (0,0)--(5,0);
	\draw[guillsep] (0,2)--(5,2);
	\draw[guillsep] (3,0)--(3,2);
	\node at (1.33,1) {$1$};
	\node at (4,1) {$2$};
	\node at (0.66,0) [circle, fill, inner sep=0.5mm] {};
	\node at (0.66,2) [circle, fill, inner sep=0.5mm] {};
	\draw [dotted] (0.66,0) -- (0.66,2);
	\draw [->] (3,2.3) -- node [midway, above] {$x-p$} (0.66,2.3);
	\end{tikzpicture}\left(	   
		u,v
	\right)
\end{align*}
by using the fact that the action of $1$ on $u$ leaves $u$ invariant. A similar computation for the second line.
\end{proof}

Following this result, pointings will often be dropped in the next equations whenever guillotine partitions are applied to ROPE elements, in order to make notations a bit lighter. 

One also checks that the irrelevance of pointing product already appear in the construction of Proposition~\ref{prop:linkwithMPS} through the use of the same scalar product on Hilbert spaces regardless of the position of the base point.

\subsubsection{Constant and universal ROPEs}\label{sec:eagerandlazyROPEs}

The presence of a colour palette for the guillotine operads introduce on the boundaries as many spaces as possible lengths. Proposition~\ref{prop:linkwithMPS} is a particular case for which the spaces $\ca{B}_{p,\infty_S}$ do not depend on the colours and are all equal to a common space $\ca{B}_S$ (and the same is true on all sides). Each gluing $\ca{B}_{p,\infty_S}\otimes \ca{B}_{p',\infty_S} \to \ca{B}_{p+p',\infty_S}$ corresponds then a bilinear product on $\ca{B}_S\otimes \ca{B}_S$ to itself. When evaluating an element such as \eqref{eq:ROPEcorresp2}, all the vertical guillotine cuts correspond to the concrete evaluation of a product on $\ca{B}_S$. 

On the other hand, with the \emph{same} objects, we may build the following alternative ROPE: we define $\ca{B}'_{p,\infty_S} = \ca{B}_S^{\otimes p}$ and the product $\ca{B}'_{p,\infty_S}\otimes \ca{B}'_{p',\infty_S}\to \ca{B}'_{p+p',\infty_S}$ is now simply the identity on tensor products. The product on $\ca{B}_{S}$ may then be performed only globally during the action on corners through the map:
\begin{align*}
\ca{B}_S^{\otimes p} \otimes \ca{B}_{\infty_E,\infty_S} &\to \ca{B}_{\infty_E,\infty_S}
\\
\left(\bigotimes_{i=1}^p a_i \right)\otimes b &\mapsto 
a_1\ldots a_p b
\end{align*}

From a categorical point of view, the second case corresponds to a universal construction above the operadic product structure. From a computational point of view, the first case corresponds to \emph{eager} evaluation (products are performed as soon as two elements are glued together) whereas the second case to \emph{lazy} evaluation (the whole products are evaluated only when actions on boundary spaces are required). We will stick to this computational vocabulary since one may imagine various intermediate cases. The discussion is important since it will control the number of equations and unknowns in Section~\ref{sec:invariantboundaryelmts} through definitions~\ref{def:eigenalgebrauptomorphims} and the following ones.

\begin{defi}\label{def:eagerandlazy}
	We consider only the discrete case $(\setP,\setL)=(\setZ,\setN)$. A ROPE $\ca{B}_{\PatternShapes(\patterntype{fp}^*)}$ is \emph{eager} if there exists four spaces $\ti{\ca{B}}_{a}$, $a\in\{S,N,W,E\}$ with associative products $m_a : \ti{\ca{B}}_{a}^{\otimes 2}\to\ti{\ca{B}}_{a} $ such that for all $r\in\setL^*$
	\begin{align*}
		\ca{B}_{r,\infty_S} &= \ti{\ca{B}}_{S}
		&
		\ca{B}_{r,\infty_N} &= \ti{\ca{B}}_{N}
		\\
		\ca{B}_{\infty_W,r} &= \ti{\ca{B}}_{W}
		&
		\ca{B}_{\infty_E,r} &= \ti{\ca{B}}_{E}
	\end{align*}
	and all the $\Guill_1$-products $m_{SN}$ (resp. $m_{WE}$) are given by the products $m_S$ and $m_N$ ($m_W$ and $m_E$).
	
	A ROPE $\ca{B}_{\PatternShapes(\patterntype{fp}^*)}$ is \emph{lazy} if there exists four spaces $\ti{\ca{B}}_{a}$, $a\in\{S,N,W,E\}$ such that for all $r\in\setL^*$
	\begin{align*}
		\ca{B}_{r,\infty_S} &= \ti{\ca{B}}_{S}^{\otimes r}
		&
		\ca{B}_{r,\infty_N} &= \ti{\ca{B}}_{N}^{\otimes r}
		\\
		\ca{B}_{\infty_W,r} &= \ti{\ca{B}}_{W}^{\otimes r}
		&
		\ca{B}_{\infty_E,r} &= \ti{\ca{B}}_{E}^{\otimes r}
	\end{align*}
	and all the $\Guill_1$-products $m_{SN}$ (resp. $m_{WE}$) are given by the tensor products on these spaces.
\end{defi}

\begin{prop}[universality of lazy ROPES]
	For any eager ROPE on spaces $\ti{\ca{B}}_{a}$ with associative products $m_a$, there exists a corresponding lazy ROPE with a morphism of $\Guill_2$-algebra from the lazy ROPE to the eager ROPE. 
\end{prop}

There is no fundamental difference in usage between both since, for ROPEs, the full plane space $\ca{B}_{\infty_{WE},\infty_{SN}}$ is the scalar field $\setK$: this implies necessarily that at some stage of gluings of shapes the products and pairings on spaces are necessarily performed. 

However, for computational perspective, lazy ROPES are bad in practice for the same reasons as Theorem~\ref{theo:stability} are bad: as the colour increases, the "size" of the objects increase dramatically and cannot be stored without simplification. The presence of a product allows for a reduction of the "size" such that objects belonging to different spaces $\ca{B}_{p,\infty_S}$ can be compared in some sense. A concrete realization of these ideas is explained in Section~\ref{sec:computational-considerations}. 

Eager ROPEs may look quite restrictive at first sight since they erase any dependence on the colours, which are however fundamental in the guillotine framework. However, there are fundamental for computations (see again Section~\ref{sec:computational-considerations}) and have a natural interpretation in terms of representation theory. A boundary side alone without corners---let's say South for this discussion--- is made of a $\Guill_1$-algebra and colours are necessary to separate how the $\Guill_2^{(\patterntype{r})}$-algebra  $\ca{A}_{\PatternShapes(\patterntype{r})}$ acts on it. As soon as a corner is introduced, Theorem~\ref{theo:1D:removingcoloursonboundaries} shows that, from the point of view of the corner alone, length colours can be erased from the computation. For ROPEs, since the rectangular part $\ca{B}_{\PatternShapes(\patterntype{r})}$ is trivial and does not depend on the colours, it is natural to use Theorem~\ref{theo:1D:removingcoloursonboundaries} to \emph{represent} all the spaces $\ca{B}_{p,\infty_S}$ as elements of $\Hom(\ca{B}_{\infty_E,\infty_S}, \ca{B}_{\infty_E,\infty_S})$, which then defines an eager ROPE, which is a \emph{final} object in its categorical interpretation: in any practical computation of a boundary weight value, all the elements of the boundary spaces are necessarily represented at some point on some vector spaces in order to obtain a scalar number at the end.

The concept of lazy ROPEs is however interesting since it shows the importance of corners and pairings. Internal products may (eager) or may (lazy) not be performed in the side $\Guill_1$-algebras since morphism of ROPEs can be used; but these products have be computed as soon as the $\Guill_1$-algebras are extended to the left and to the right and this is essentially due to the fact that $\infty_L$, $\infty_R$ and $\infty_{LR}$ are absorbing colours with "smaller" structures in order to achieve total reduction to the scalar field for $\ca{B}_{\infty_{WE},\infty_{SN}}=\setK$. 

The important conclusion behind this discussion is that constraints on the objects emerge only when products,  actions or pairings are computationally performed, i.e. every time a map $\ca{A}\otimes \ca{B} \to \ca{C}$ leads to a space $\ca{C}$ intuitively "smaller" than $\ca{A}\otimes \ca{B}$. For finite matrices, "smaller" corresponds to a smaller dimension; for infinite objects, the situation is more subtle as it may be seen in concrete models such as \cite{BodiotSimon}.

\subsection{ROPE representations of functions.}
\subsubsection{Definition of ROPEreps.}
We may now proceed and use ROPEs to obtain interesting parametrizations of functions over boundaries of rectangles decorated by sequences in two sets $S_1$ and $S_2$, such as boundary weights of two-dimensional Markov processes. The following definition rewrites \emph{matrix product state} or \emph{matrix Ans\"atze} in the language of $\Guill_2$-algebras and ROPEs.

\begin{defi}[rectangular operator product environment representation (ROPErep) for finite state space and discrete space]\label{def:ROPErep:FD}
Let $S_1$ and $S_2$ be two finite sets. Let $I$ be a set and $(p_i,q_i)_{i\in I}$ a sequence of integers in $\setN\times\setN\setminus\{(0,0)\}$. Let $(f_i)_{i\in I}$ be a collection of functions such that $f_i: S_1^{p_i}\times S_1^{p_i}\times S_2^{q_i}\times S_2^{q_i} \to \setK$ if $p_i\neq 0$ and $q_i\neq 0$ (boundaries of non-degenerate rectangles) and $f_i : S_1^{p_i} \to \setK$ if $q_i=0$ and $f_i:S_2^{q_i}\to\setK$ if $p_i=0$ (rectangles degenerating on segments). A \emph{homogeneous rectangular operator product environment representation} (ROPErep) of $(f_i)_{i\in I}$ over a ROPE $\ca{B}_{\PatternShapes(\patterntype{fp}^*)}$ consists in:
\begin{enumerate}[(i)]
\item for every $a\in \{N,S\}$, a function $S_1 \to \ca{B}_{1,\infty_a}$, $x\mapsto A_{a,1}(x)$.
\item for every $b\in \{W,E\}$, a function $ S_2 \to \ca{B}_{\infty_b,1}$, $x\mapsto A_{b,1}(x)$.
\item for every $(a,b)\in \{N,S\}\times\{W,E\}$, elements $U_{a,b}\in \ca{B}_{\infty_a,\infty_b}$,
\end{enumerate}
such that:
\begin{itemize}
\item for all $i\in I$ with $p_i>0$ and $q_i>0$, for all $(x,y,w,z)\in S_1^{p_i}\times S_1^{p_i}\times S_2^{q_i}\times S_2^{q_i}$, it holds:
\begin{subequations}
\label{eq:completeRMPR}
\begin{equation}
f_i\left(
\begin{tikzpicture}[scale=0.45,baseline={(current bounding box.center)}]
\draw (0,0) rectangle (4,3);
\node at (0.5,0) [anchor = north] {$x_1$};
\node at (1.5,0) [anchor = north]{$x_2$};
\node at (2.5,0) [anchor = north]{$\ldots$};
\node at (3.5,0) [anchor = north]{$x_{p_i}$};
\node at (0.5,3) [anchor = south] {$y_1$};
\node at (1.5,3) [anchor = south]{$y_2$};
\node at (2.5,3) [above]{$\ldots$};
\node at (3.5,3) [above]{$y_{p_i}$};
\node at (0,0.5) [anchor = east]{ $w_1$ };
\node at (0,1.5) [anchor = east]{ $\vdots$ };
\node at (0,2.5) [anchor = east]{ $w_{q_i}$ };
\node at (4,0.5) [anchor = west]{ $z_1$ };
\node at (4,1.5) [anchor = west]{ $\vdots$ };
\node at (4,2.5) [anchor = west]{ $z_{q_i}$ };
\end{tikzpicture}
\right) = 
\begin{tikzpicture}[guillpart,yscale=1.35,xscale=3.1]
\draw[guillsep, dotted] (0,0) rectangle (6,5);
\draw[guillsep] (1,1) rectangle (5,4);
\draw[guillsep] 	(1,0)--(1,1)
				(2,0)--(2,1)
				(3,0)--(3,1)
				(4,0)--(4,1)
				(5,0)--(5,1);
\draw[guillsep] 	(1,5)--(1,4)
				(2,5)--(2,4)
				(3,5)--(3,4)
				(4,5)--(4,4)
				(5,5)--(5,4);
\draw[guillsep] 	(0,1)--(1,1)
				(0,2)--(1,2)
				(0,3)--(1,3)
				(0,4)--(1,4);
\draw[guillsep] 	(5,1)--(6,1)
				(5,2)--(6,2)
				(5,3)--(6,3)
				(5,4)--(6,4);
\node at (0.5,0.5) { $U_{SW}$ };
\node at (5.5,0.5) { $U_{SE}$ };
\node at (0.5,4.5) { $U_{NW}$ };
\node at (5.5,4.5) { $U_{NE}$ };
\node at (1.5,0.5) { $A_{S,1}(x_1)$ };
\node at (2.5,0.5) { $A_{S,1}(x_2)$ };
\node at (3.5,0.5) { $\ldots$ };
\node at (4.5,0.5) { $A_{S,1}(x_{p_i})$ };
\node at (1.5,4.5) { $A_{N,1}(y_1)$ };
\node at (2.5,4.5) { $A_{N,1}(y_2)$ };
\node at (3.5,4.5) { $\ldots$ };
\node at (4.5,4.5) { $A_{N,1}(y_{p_i})$ };
\node at (0.5,1.5) { $A_{W,1}(w_1)$ };
\node at (0.5,2.5) { $\vdots$ };
\node at (0.5,3.5) { $A_{W,1}(w_{q_i})$ };
\node at (5.5,1.5) { $A_{E,1}(z_1)$ };
\node at (5.5,2.5) { $\vdots$ };
\node at (5.5,3.5) { $A_{E,1}(z_{q_i})$ };
\end{tikzpicture}
\end{equation}
\item for all $i\in I$ with $q_i=0$, for all $x\in S_1^{p_i}$,
\begin{equation}
	f_i(x) = \begin{tikzpicture}[guillpart,yscale=1.35,xscale=3.1]
		\draw[guillsep, dotted] (0,0) rectangle (6,2);
		%\draw[guillsep] (1,1) rectangle (5,4);
		\draw[guillsep] 	(1,0)--(1,2)
		(2,0)--(2,2)
		(3,0)--(3,2)
		(4,0)--(4,2)
		(5,0)--(5,2);
		\draw[guillsep] (0,1)--(6,1);
		\node at (0.5,0.5) { $U_{SW}$ };
		\node at (5.5,0.5) { $U_{SE}$ };
		\node at (0.5,1.5) { $U_{NW}$ };
		\node at (5.5,1.5) { $U_{NE}$ };
		\node at (1.5,0.5) { $A_{S,1}(x_1)$ };
		\node at (2.5,0.5) { $A_{S,1}(x_2)$ };
		\node at (3.5,0.5) { $\ldots$ };
		\node at (4.5,0.5) { $A_{S,1}(x_{p_i})$ };
		\node at (1.5,1.5) { $A_{N,1}(x_1)$ };
		\node at (2.5,1.5) { $A_{N,1}(x_2)$ };
		\node at (3.5,1.5) { $\ldots$ };
		\node at (4.5,1.5) { $A_{N,1}(x_{p_i})$ };
	\end{tikzpicture}
\end{equation}
\item for all $i\in I$ with $p_i=0$, for all $w\in S_2^{q_i}$,
\begin{equation}
	f_i(w) = \begin{tikzpicture}[guillpart,yscale=1.35,xscale=3.3]
		\draw[guillsep, dotted] (0,0) rectangle (2,5);
		\draw[guillsep] (1,0)--(1,5);
		\draw[guillsep] 	(0,1)--(2,1)
		(0,2)--(2,2)
		(0,3)--(2,3)
		(0,4)--(2,4);
		\node at (0.5,0.5) { $U_{SW}$ };
		\node at (1.5,0.5) { $U_{SE}$ };
		\node at (0.5,4.5) { $U_{NW}$ };
		\node at (1.5,4.5) { $U_{NE}$ };
		\node at (0.5,1.5) { $A_{W,1}(w_1)$ };
		\node at (0.5,2.5) { $\vdots$ };
		\node at (0.5,3.5) { $A_{W,1}(w_{q_i})$ };
		\node at (1.5,1.5) { $A_{E,1}(w_1)$ };
		\node at (1.5,2.5) { $\vdots$ };
		\node at (1.5,3.5) { $A_{E,1}(w_{q_i})$ };
	\end{tikzpicture}
\end{equation}
\end{subequations}
\end{itemize}
\end{defi}

As stated in Proposition~\ref{prop:linkwithMPS}, whenever the rectangular matrix product environment is the one of linear algebra, i.e. $\ca{B}_{p,\infty_a}= \End(V_a)$, $\ca{B}_{\infty_b,q}=\End(V_b)$ with suitable intertwiners $U_{ab}$ and suitable transpositions (on the North and on the West) to obtain a counter-clockwise product, the previous representation \eqref{eq:firstrectangularMPS} corresponds to usual matrix product states and matrix Ans\"atze \eqref{eq:firstrectangularMPS}.
Section \ref{sec:basicRMPR} below provides various illustrations of such ROPEreps of boundary weights of some probabilistic models.

\subsubsection{Local-to-global principle and back}

As already announced, homogeneous ROPEreps allows one to build functions on arbitrary large power sets with cardinal  $|S_1|^{2p_i} |S_2|^{2q_i}$ out of a \emph{finite} number of objects the ROPE spaces, in a \emph{local-to-global} construction: gluing of neighbouring segments labelled with values $x_i$ corresponds to a product in the ROPE. The values of the functions $f_i$ are obtained after gluing all pieces of the boundary together. 

One may wonder how to build a ROPErep for a given set of functions $(f_i)$. This would correspond to a global-to-local construction and this is highly non-trivial in generic cases to the best of our knowledge.

Given a single function $f: S_1^p\times S_1^p\times S_2^q \times S_2^q \to \setR$, it is easy to see that a finite-dimensional eager ROPE $\ca{B}_{p,\infty_{a}}=\ca{B}_{\infty_b,q}=\End(W)$ and a tracial state can always be found by considering a space $W$ sufficiently large to register all the possible arguments of the function. Moreover, a lot of ROPEreps can be found for this function $f$. But this is not the spirit of ROPEreps: the idea is to find a \emph{common} underlying structure for an arbitrary large or infinite set of functions. To fulfil this purpose, the general idea is rather to find common properties to all the functions in the collection and encode them in the ROPErep.

The idea can be illustrated in the following very basic example. Let $A$ be a fixed subset of $S_1^{r}$ for some integer $r\geq 1$. If all the functions $f_i$ satisfy for all $y\in S_1^{p_i}$, $w,z\in S_2^{q_i}$, for all $0\leq k\leq p_i-r$, for all $x\in S_1^k$ and $x'\in S_1^{p_i-r}$,
\[
\sum_{a \in A} \alpha_a f_i( xax',y,w,z)=0
\]
for some scalar coefficients $\alpha : A\to\setK$, where $x a x'$ is the concatenation of the three sequences $x$, $a$ and $x'$, one may wish to encode it directly at the level of the ROPErep by requiring the following equality 
\[
\sum_{a \in A} \alpha_a A_{S,1}(a_1)A_{S,1}(a_2)\ldots A_{S,1}(a_r) = 0
\]
in the boundary space $\ca{B}_{r,\infty_S}$.

On the probabilistic side, the collection of functions $f_i$ may be provided by the boundary weights of Markov processes obtained from a translation-invariant Gibbs measure on $\setZ^2$. The structure of a ROPErep should then reflect  all the remarkable identities on boundary weights deduced from the Gibbs measure (identities on correlation functions, correlation lengths, etc.). We will see below in Section~\ref{sec:invariantboundaryelmts} that it will be preferable to first compute explicitly the ROPEreps by solving suitable equations and then deduce the Gibbs measure and its properties from the ROPErep following the local-to-global principle.

\subsubsection{Non-unicity of ROPEreps and morphisms of guillotine algebras.}

A very important point is that the ROPErep of a collection of functions may \emph{not} be unique and this non-uniqueness will allow us for crucial reductions in Section~\ref{sec:invariantboundaryelmts}. We now provide an easy way of building various ROPEreps from a given one.

To this purpose, we recall the notion of morphisms of algebras over a coloured operad in the special case of $\Guill_2^{(\patterntype{fp}^*)}$-algebras. This definition is classical in the theory of operads.

\begin{defi}[morphisms of $\Guill_2$-algebras]\label{def:guill2morphisms}
Let $\ca{A}_{\PatternShapes(\patterntype{fp}^*)}$ and $\ca{A}'_{\PatternShapes(\patterntype{fp}^*)}$ be two $\Guill_2^{(\patterntype{fp}^*)}$-algebras. A collection of maps $\Phi_{\PatternShapes(\patterntype{fp}^*)}=(\Phi_{p,q})_{(p,q)\in \PatternShapes(\patterntype{fp}^*)}$ is a morphism of $\Guill_2^{(\patterntype{fp}^*)}$-algebras if and only if:
\begin{enumerate}[(i)]
\item for all $(p,q)\in\PatternShapes(\patterntype{fp}^*)$, $\Phi_{p,q}$ is a linear map from $\ca{A}_{p,q}$ to $\ca{A}'_{p,q}$
\item for all $(p,q)\in\PatternShapes(\patterntype{fp}^*)$ and for any decompositions $p=p_1+p_2$ and $q=q_1+q_2$, for all $(a,b)\in \ca{A}_{p,q_1}\otimes \ca{A}_{p,q_2}$ and for all $(c,d)\in\ca{A}_{p_1,q}\otimes\ca{A}_{p_2,q}$,
\begin{align*}
\Phi_{p,q}\left(
 \begin{tikzpicture}[guillpart,yscale=1.6,xscale=2.]
 	\fill[guillfill] (0,0) rectangle (1,2);
 	\draw[guillsep] (0,1)--(1,1);
 	\draw[dashed] (0,0)--(1,0)--(1,2)--(0,2)--(0,0);
 	\node at (0.5,0.5) { $a$ };
 	\node at (0.5,1.5) { $b$ };
 \end{tikzpicture}_{\ca{A}}
\right)
&=
 \begin{tikzpicture}[guillpart,yscale=1.6,xscale=2.7]
 	\fill[guillfill] (0,0) rectangle (1,2);
 	\draw[guillsep] (0,1)--(1,1);
 	\draw[dashed] (0,0)--(1,0)--(1,2)--(0,2)--(0,0);
 	\node at (0.5,0.5) { $\Phi_{p,q_1}(a)$ };
 	\node at (0.5,1.5) { $\Phi_{p,q_2}(b)$ };
 \end{tikzpicture}_{\ca{A}'}
&
\Phi_{p,q}\left(
 \begin{tikzpicture}[guillpart,yscale=1.6,xscale=2,rotate=-90]
 	\fill[guillfill] (0,0) rectangle (1,2);
 	\draw[guillsep] (0,1)--(1,1);
 	\draw[dashed] (0,0)--(1,0)--(1,2)--(0,2)--(0,0);
 	\node at (0.5,0.5) { $c$ };
 	\node at (0.5,1.5) { $d$ };
 \end{tikzpicture}_{\ca{A}}
\right)
&=
 \begin{tikzpicture}[guillpart,yscale=1.6,xscale=3.,rotate=-90]
 	\fill[guillfill] (0,0) rectangle (1,2);
 	\draw[guillsep] (0,1)--(1,1);
 	\draw[dashed] (0,0)--(1,0)--(1,2)--(0,2)--(0,0);
 	\node at (0.5,0.5) { $\Phi_{p_1,q}(c)$ };
 	\node at (0.5,1.5) { $\Phi_{p_2,q}(d)$ };
 \end{tikzpicture}_{\ca{A}'}
\end{align*}
\end{enumerate}
with the following arithmetic on the colours: $p+q$ is the usual addition for finite $p$ and $q$, $r+\infty_a=\infty_a$ for finite $r$ and $a\in\{S,N,W,E\}$, $\infty_S+\infty_N =\infty_{SN}$ and $\infty_{W}+\infty_E =\infty_{WE}$. The dashed lines correspond to either a full line or no line depending on whether $p$ or $q$ is finite. The index $\ca{A}$ or $\ca{A}'$ indicates in which algebra the products are evaluated.
\end{defi}

Such morphisms of $\Guill_2^{(\patterntype{fp}^*)}$-algebra with some additional restrictions help to build alternative ROPEreps of a function from an initial ROPErep.

\begin{prop}[equivalent ROPE representation of a function]\label{prop:ROPErepuptomorphism}
Let $\ca{B}_{\PatternShapes(\patterntype{fp}^*)}$ and $\ca{B}'_{\PatternShapes(\patterntype{fp}^*)}$ be two ROPEs. Let $(f_i)_{i\in I}$ be a collection of functions $f_i:S_1^{p_i}\times S_1^{p_i}\times S_2^{q_i}\times S_2^{q_i} \to \setK$ (resp. $S^{(1)}_{p_i}\times S^{(1)}_{p_i}\times S^{(2)}_{q_i}\times S^{(2)}_{q_i}$) with a ROPErep \eqref{eq:completeRMPR} over $\ca{B}$. If $\Phi_{\PatternShapes(\patterntype{fp}^*)}=(\Phi_{p,q})_{(p,q)\in \PatternShapes(\patterntype{fp}^*)}$ is a morphism of $\Guill_2^{(\patterntype{fp}^*)}$-algebras from $\ca{B}$ to $\ca{B}'$ such that, for any $\lambda\in\setK$,
\begin{enumerate}[(i)]
\item for any $(p,q)\in\PatternShapes(\patterntype{r})$, $\Phi_{p,q}(\lambda)=\lambda$
\item $\Phi_{\infty_{WE},\infty_{SN}}(\lambda)=\lambda$, 
\end{enumerate}
then the collection of functions $(f_i)_{i\in I}$ admits a ROPErep over $\ca{B}'$ given by
\begin{equation}
\label{eq:morphismROPErep}
f\left(x,y,w,z\right)
 = 
\begin{tikzpicture}[guillpart,yscale=1.6,xscale=3.2]
\draw[guillsep, dotted] (0,0) rectangle (3,3);
\draw[guillsep] (1,1) rectangle (2,2);
\draw[guillsep] 	(1,0)--(1,1)
				(2,0)--(2,1);
\draw[guillsep] 	(1,3)--(1,2)
				(2,3)--(2,2);
\draw[guillsep] 	(0,1)--(1,1)
				(0,2)--(1,2);
\draw[guillsep] 	(2,1)--(3,1)
				(2,2)--(3,2);
\node at (0.5,0.5) { $\ti{U}_{WS}$ };
\node at (2.5,0.5) { $\ti{U}_{SE}$ };
\node at (0.5,2.5) { $\ti{U}_{NW}$ };
\node at (2.5,2.5) { $\ti{U}_{EN}$ };
\node at (1.5,0.5) { $\ti{A}_{S,p_i}(x)$ };
\node at (1.5,2.5) { $\ti{A}_{N,p_i}(y)$ };
\node at (0.5,1.5) { $\ti{A}_{W,q_i}(w_1)$ };
\node at (2.5,1.5) { $\ti{A}_{E,q_i}(z_1)$ };
\end{tikzpicture}
\end{equation}
with the new elements, for all $(a,b)\in\{S,N\}\times \{W,E\}$, for all $p,q\in\setL^*$, for all $x\in S_1^p$ (resp. $S^{(1)}_p$) and all $w\in S_2^q$ (resp. $S^{(2)}_q$), given by
\begin{align*}
\ti{U}_{ab} &= \Phi_{\infty_b,\infty_a}(U_{ab}) \\
\ti{A}_{a,p}(x) &= \Phi_{p,\infty_a}( A_{a,p}(x) )
&
\ti{A}_{b,q}(w) &= \Phi_{\infty_b,q}( A_{b,q}(w) )
\end{align*}
\end{prop}
\begin{proof}
As an element of $\setK=\ca{B}'_{\infty_{WE},\infty_{SN}} = \ca{B}_{\infty_{WE},\infty_{SN}}$, we simply write $f(x,y,w,z)=\Phi_{\infty_{WE},\infty_{SN}}(f(x,y,w,z))$. We now introduce \eqref{eq:ROPEcorresp5} in the argument of $\Phi_{\infty_{WE}}$ and 
use the morphism property of $(\Phi_{p,q})$ to split it into the expected form (reduction of the central $1_\setK^{(p,q)}$ uses the identity property of $\Phi_{p,q}$ for finite $p$ and $q$).
\end{proof}

This non-unique choice of ROPEreps of a given function is widely used in the next sections in order to build a notion of \emph{invariant boundary weights}. 

\subsubsection{Fundamental examples of reduction of dimension of a ROPErep} \label{sec:ROPErep:fundamentalexample} We first illustrate this property on a simple case, typical of linear algebra. We assume $\ca{B}_{p,q} = \Mat_{n_a,n_a}(\setR)$ when exactly one element among $p$ or $q$ is finite and $a\in\{S,N,W,E\}$ is the associated side. We also set $\ca{B}_{\infty_b,\infty_a}=\setR^{n_a}\otimes \setR^{n_b}$ on the corners with the traditional left and right actions of $\Mat_{n,n}(\setR)$ on $\setR^n$ and the euclidean scalar product. We now assume that, on each side, matrices $A_{a,p}(x)$ is similar an upper diagonal block structure 
\begin{equation}
A_{a,p}(x) = K\begin{pmatrix}
A_{a,p}^{1}(x)  & R_{a,p}^{1,2}(x) & \cdots & R_{a,p}^{1,r_a}(x) \\ 
0 & A_{a,p}^{2 }(x) & \cdots & R_{a,p}^{2,r_a}(x) \\
0 & 0 &  \ddots & \vdots\\
0 & 0 & 0 & A_{a,p}^{r_a}(x)
\end{pmatrix}K^{-1}
\end{equation}
where the $A_{a,p}^{k}(x)$ (resp. $R_{a,p}^{k,l}(x)$) have dimension $d_{a,k}\times d_{a,k}$ (resp. , $d_{a,k}\times d_{a,l}$) with dimensions $d_{a,k}$ that satisfy $\sum_{k=1}^{r_a}=n_a$. 

For each side $a\in\{S,N,W,E\}$, we now choose a specific sub-block $k_a$. We may define morphisms 
\begin{equation}
\Phi_a(A_{a,p}(x)) =  P_a K^{-1} A_{a,p}(x) K Q_a
\end{equation}
with matrices $P_a$ and $Q_a$ that have the block structure
\begin{align*}
P_a  &=  \begin{pmatrix} 0 &\cdots & 0 & I_{d_{a,k}} & P_a^{k_a+1} & \cdots & P_a^{r_a} \end{pmatrix} \in \mathrm{Mat}_{d_{k_a},n_a}(\setR)
\\
Q_a  &=\begin{pmatrix} Q_a^{1} \\ \vdots \\ Q_a^{k_a-1} \\ I_{d_{a,k}} \\ 0 \\ \vdots \\ 0 \end{pmatrix}\in \mathrm{Mat}_{n_a,d_{k_a}}(\setR)
\end{align*}
We have indeed $P_a Q_a = I_{d_a}$ and the morphism property \begin{equation}
\Phi_a  ( A_{a,p}(x_1,\ldots,x_p)A_{a,p'}(x_{p+1},x_{p+p'}))
= A_{a,p+p'}(x_1,\ldots,x_{p+p'})
\end{equation}
and its equivalent on segment state spaces with morphisms instead of concatenations. Using the morphisms $\Phi_a$ on each side provides the first part of a morphism of ROPEreps from a $n_a$-dimensional one to a $d_{a,k_a}$-dimensional one. We also observe that there is still a freedom in the choice of the blocks $(Q_a^{l})_{l<k_a}$ and $(P_a^{l})_{l>k}$.

In order to have a complete morphism of ROPErep, we still need to specify \emph{both} a particular structure for the corner elements and particular morphisms on corners in a coherent way. We focus here on the South-East element as an example and leave the reader write down the exact relations on the three other corners. Whenever the corner element $U_{SE}$ is of the form $(K_SQ_S\otimes {K_E^{-1}}^tP_E^t ) V_{SE} \equiv  K_SQ_SV_{SE} P_E K_E^{-1}$, then the morphism $\Phi_{\infty_E,\infty_S}$ defined by:
\begin{equation}
\Phi_{\infty_E,\infty_S}(u) = (P_S K_S^{-1}\otimes Q_E^t K_E^t)  u \equiv P_S K_S^{-1} u K_E Q_E
\end{equation}
with the identification $\setR^{n_S}\otimes \setR^{n_E} \equiv \Mat_{n_S,n_E}(\setR)$.
satisfy the morphism property:
\begin{equation}
\Phi_{\infty_E,\infty_S}\left( 
	A_{S,p}(x) U_{SE} A_{E,q}(z)
	\right) 
	=
	A_{S,p}^{k_S}(x)V_{SE} A_{E,q}^{k_E}(z)
\end{equation}

Defining consequently half-plane morphisms $\Phi_{\infty_{WE},\infty}(u) = Q_W u Q_E$ and strip morphisms in a similar way and choosing the corner elements similarly then provide a complete morphism of ROPErep such that matrices $A_a^{(i)}(x)$ are replaced by sub-blocks $A_a^{(i)}(x)$ and the corner elements $U_{ba}$ are replaced by the lower-dimensional $V_{ba}$. We also observe that the extraction of sub-blocks of the side matrices can be performed with various operators $Q_a$ and $P_a$, whose sub-blocks $(Q_a^{l})_{l<k_a}$ and $(P_a^{l})_{l>k}$ can be left free and are related to the choice of the corner elements $U_{ba}$ built out of the $V_{ba}$. 

This example of morphism of ROPEreps is the key picture to have in mind while reading definitions~\ref{def:eigenalgebrauptomorphims} and \ref{def:cornereigensemigroups} in their abstract form. In practice however, most ROPEreps are infinite dimensional: the previous finite-dimensional structure may be lifted to von Neumann algebras with suitable projectors and identification of factors and characters. Such a case appears naturally for Gaussian Markov models, see \cite{BodiotSimon}.

		\section{Relevant boundary conditions described by matrix product weights}\label{sec:basicRMPR}

This section is devoted to the description of some basic ROPEreps that correspond to classical boundary weights appearing in statistical mechanics as well as interesting generalizations, which appear naturally in the ROPEreps framework.

\subsection{Factorized boundary weights and constant boundary conditions.}\label{par:factorizedbw}

We consider first the case of factorized functions:
\[
g\left(
\begin{tikzpicture}[scale=0.5,baseline={(current bounding box.center)}]
\draw (0,0) rectangle (4,3);
\node at (0.5,0) [anchor = north] {$x_1$};
\node at (1.5,0) [anchor = north]{$x_2$};
\node at (2.5,0) [anchor = north]{$\ldots$};
\node at (3.5,0) [anchor = north]{$x_p$};
\node at (0.5,3) [anchor = south] {$y_1$};
\node at (1.5,3) [anchor = south]{$y_2$};
\node at (2.5,3) [above]{$\ldots$};
\node at (3.5,3) [above]{$y_p$};
\node at (0,0.5) [anchor = east]{ $w_1$ };
\node at (0,1.5) [anchor = east]{ $\vdots$ };
\node at (0,2.5) [anchor = east]{ $w_q$ };
\node at (4,0.5) [anchor = west]{ $z_1$ };
\node at (4,1.5) [anchor = west]{ $\vdots$ };
\node at (4,2.5) [anchor = west]{ $z_q$ };
\end{tikzpicture}
\right)= \prod_{i=1}^p u_S(x_i) \prod_{j=1}^p u_N(y_j) \prod_{k=1}^q u_W(w_k) \prod_{l=1}^q u_E(z_l)
\]
where $u_S, u_N : S_1 \to \setK$ and $u_W, u_E: S_2 \to \setK$ as already seen in Section~\ref{sec:factorizedweights:partone}. In this case, it is easy to see that the simplest choice $\ca{B}_{p,q} = \setK$ for any $(p,q)\in\PatternShapes(\patterntype{fp}^*)$ of rectangular matrix product structure (with all products equal to the product in $\setK$) provides a (inhomogeneous) ROPErep with one-dimensional matrices $A_{a,1}(x)= u_a(x) \in \setK$ for any $a\in\{S,N,W,E\}$ and $U_{ab}=1$ on each of the four corners. In Section~\ref{sec:factorizedweights:partone}, the same description was done by embedding the functions $u_{a}$ in the canonical structure $\ca{T}_{\PatternShapes(\patterntype{fp}^*)}$ and one may check indeed that there exists a trivial morphism of ROPEs between both representation.

In particular, this encompasses the case of constant functions $u_a = \indic{x_a}$ (fixed value on each side). Breaking homogeneity of the ROPEreps of such factorized functions would allow for fixed sequences of values on each side instead of a single value.

\subsection{Boundary conditions parametrized by a dynamical system or a hidden Markov chain.} 

The structure of matrix product states \eqref{eq:firstrectangularMPS} is very reminiscent of hidden Markov chains if the matrices $A_a(x)$ and the corner matrices have positive coefficients in some basis: the indices of the matrices are then interpreted as an internal state space. 

\begin{defi}\label{def:hiddenmarkovrect}
Let $S_{\mathrm{hidden}}$ be a countable space. A hidden Markov process on $\partial R$ is a couple of processes $((S_v)_{v\in V(\partial R)}, (X_e)_{e\in\partial R})$ such that:
\begin{enumerate}[(i)]
\item  $(S_v)_v$ is a process with values in $S_{\mathrm{hidden}}$, $(X_e)_{e}$ has values in $S_1$ (resp. $S_2$) for horizontal (resp. vertical) edges.
\item the law of $(S_v)_v$ is given, for any sequence $(s_v)_v$ by:
\[
\prob{ (S_v) = (s_v) } = \frac{1}{Z_R^{\mathrm{hidden}}} 
\prod_{ (v_1,v_2)\in \partial_{\mathrm{or.}}R }\MarkovWeight{T}_{\mathrm{side}(v_1,v_2)}( s_{v_1}, s_{v_2} )
\]
where $\MarkovWeight{T}_{a}: S_{\mathrm{hidden}}^2 \to \setR_+$, $a\in \{S,N,W,E\}$ are edge weights associated to the four sides and the function partition $Z_R^\mathrm{hidden}$ is the trace of the counter-clockwise product of the weights $\MarkovWeight{\Omega}_a$ and is assumed to be finite.

\item conditionally on the process $(S_v)$, the law of the process $(X_e)$ is given by \[
\probc{ (X_e)_e= (x_e)_e }{ (S_v)_v } = \prod_{e\in \partial_{\mathrm{or.}} R } \nu_{\mathrm{side}(e)}( x_e | S_{s(e)}, S_{t(e)} )
\]
where the $\nu_a( \cdot | s_1,s_2)$ are probability laws on $S_1$ for $a\in\{S,N\}$ and $S_2$ for $a\in \{W,E\}$.
\end{enumerate}
\end{defi}
This means that, conditionally on the process $(S_v)$, the r.v. $(X_e)$ are independent and the law of each $X_e$ depends only on the two values of $(S_v)$ on the extremities of the edge $e$. Because of the circular geometry of $\partial R$, the process $(S_v)$ should rather be called a Markov bridge.

\begin{prop}\label{prop:ROPErep:hiddenmarkov}
Under the same notations as the previous Definition~\ref{def:hiddenmarkovrect}, the ROPE defined according to Proposition~\ref{prop:linkwithMPS} with the spaces
\[
\ca{H}_W = \ca{H}_E = \ca{H}_S = \ca{H}_N = \setR^{|S_\mathrm{hidden}|},
\]
the ROPE elements defined by the following matrices (indices $\alpha$ and $\beta$ are elements of $S_\mathrm{hidden}$)
\begin{subequations}\label{eq:hiddenmarkovrates}
\begin{align}
A_S(x)_{\alpha,\beta} &= \nu_S(x | \alpha,\beta ) \MarkovWeight{T}_S(\alpha,\beta) 
& 
A_N(x)_{\alpha,\beta} &= \nu_N(x | \beta,\alpha ) \MarkovWeight{T}_N(\beta,\alpha) 
\\
A_W(x)_{\alpha,\beta} &= \nu_W(x | \beta,\alpha ) \MarkovWeight{T}_W(\beta,\alpha) 
& 
A_E(x)_{\alpha,\beta} &= \nu_E(x | \alpha,\beta ) \MarkovWeight{T}_E(\alpha,\beta),
\end{align}
\end{subequations}
and the corner bimodule elements $U_{ab} = \id_{\setR^{|S_\mathrm{hidden}|}}$
determine a complete homogeneous ROPE representation of the weights $(h_R)$ defined by the \emph{unnormalized} laws of the process $(X_e)$:
\[
h_R( (x_e) ) = Z_R^\mathrm{hidden} \prob{(X_e)=(x_e)}
\] 
\end{prop}
\begin{proof}
The marginal law of $(X_e)$ is obtained by summing the product of \eqref{eq:hiddenmarkovrates} over all the possible values of the hidden process $(S_e)$. The previous proposition is then a straightforward reinterpretation of the sums as product of matrices correctly indexed (and correctly transpose to recast counter-clockwise products to West-East and South-North products).
\end{proof}

This encompasses naturally the previous case of factorized weights by taking $S_\mathrm{hidden}$ with cardinal one and a trivial dynamics. This is also the case for the traditional (segment and not rectangle) Matrix Ansatz solution of the exclusion process \cite{DEHP}. It is an open question for us of whether Perron-Frobenius-like fundamental eigenvectors as defined in the next section can be rewritten systematically in the form of hidden Markov chains with suitable hidden state space $S_\mathrm{hidden}$.

	\section{Structural stability for boundary weights with a matrix product structure}\label{sec:ROPEstability}
		\subsection{The theorem and its proof}

We now have all the definitions to formulate and prove our main stability theorem below. 

\begin{theo}[structural stability of ROPEreps of boundary weights of Markov processes]\label{theo:stability}
Let $S_1$ and $S_2$ be two finite sets. Let $R=[X_1,X_2]\times [Y_1,Y_2]$ be a non-degenerate rectangle in $\setZ^2$ with shape $(P,Q)=(X_2-X_1,Y_2-Y_1)\in \setL^*\times\setL^*$. Let $(X_e)_{e\in \Edges{R}}$ be a $(S_1,S_2)$-valued Markov process on $R$ with face weights $(\MarkovWeight{W}_{x,y})_{X_1\leq x < X_2,Y_1\leq y< Y_2}$ and boundary weight $g_{R}: S_1^{P}\times S_1^{P} \times S_2^{Q}\times S_2^{Q} \to \setR_+$.

If the boundary weight $g_R$ admits a (not necessarily) rectangular operator product representation over a ROPE $\ca{B}_{\PatternShapes(\patterntype{fp}^*)}$ given by:
\begin{equation}\label{eq:ROPErepforboundaryweight}
g_R(x,y,w,z) = \begin{tikzpicture}[guillpart,yscale=1.6,xscale=3.4]
\draw[guillsep, dotted] (0,0) rectangle (5,5);
\draw[guillsep] (1,1) rectangle (4,4);
\draw[guillsep] 	(1,0)--(1,1)
				(2,0)--(2,1)
				(3,0)--(3,1)
				(4,0)--(4,1);
\draw[guillsep] 	(1,5)--(1,4)
				(2,5)--(2,4)
				(3,5)--(3,4)
				(4,5)--(4,4);
\draw[guillsep] 	(0,1)--(1,1)
				(0,2)--(1,2)
				(0,3)--(1,3)
				(0,4)--(1,4);
\draw[guillsep] 	(5,1)--(4,1)
				(5,2)--(4,2)
				(5,3)--(4,3)
				(5,4)--(4,4);
\node at (0.5,0.5) { $U_{SW}$ };
\node at (4.5,0.5) { $U_{SE}$ };
\node at (0.5,4.5) { $U_{NW}$ };
\node at (4.5,4.5) { $U_{NE}$ };
\node at (1.5,0.5) { $A_S^{(1)}(x_1)$ };
\node at (2.5,0.5) { $\ldots$ };
\node at (3.5,0.5) { $A_S^{(P)}(x_P)$ };
\node at (1.5,4.5) { $A_N^{(1)}(y_1)$ };
\node at (2.5,4.5) { $\ldots$ };
\node at (3.5,4.5) { $A_N^{(P)}(y_P)$ };
\node at (0.5,1.5) { $A_W^{(1)}(w_1)$ };
\node at (0.5,2.5) { $\vdots$ };
\node at (0.5,3.5) { $A_W^{(Q)}(w_Q)$ };
\node at (4.5,1.5) { $A_E^{(1)}(z_1)$ };
\node at (4.5,2.5) { $\vdots$ };
\node at (4.5,3.5) { $A_E^{(Q)}(z_Q)$ };
\end{tikzpicture}
\end{equation}
then there exists a larger ROPE $\ca{B}'_{\PatternShapes(\patterntype{fp}^*)}$ such that, for any non-degenerate sub-rectangle $R'=[X'_1,X'_2]\times[Y'_1,Y'_2]\subset R$ with shape $(P',Q')=(X'_2-X'_1,Y'_2-Y'_1)$, the restriction $(X_e)_{e\in\Edges{R'}}$ of $(X_e)$ to $R'$ is a Markov process with the same restricted face weights $(\MarkovWeight{W}_{x,y})_{(x,y)\in R'}$ and a boundary weight $g_{R'}$ that admits a ROPE representation over $\ca{B}'_{\PatternShapes(\patterntype{fp}^*)}$ given by
\begin{equation}\label{eq:stabtheo:newROPErep}
g_{R'}(x,y,w,z) = \begin{tikzpicture}[guillpart,yscale=1.65,xscale=3.4]
\draw[guillsep, dotted] (0,0) rectangle (5,5);
\draw[guillsep] (1,1) rectangle (4,4);
\draw[guillsep] 	(1,0)--(1,1)
				(2,0)--(2,1)
				(3,0)--(3,1)
				(4,0)--(4,1);
\draw[guillsep] 	(1,5)--(1,4)
				(2,5)--(2,4)
				(3,5)--(3,4)
				(4,5)--(4,4);
\draw[guillsep] 	(0,1)--(1,1)
				(0,2)--(1,2)
				(0,3)--(1,3)
				(0,4)--(1,4);
\draw[guillsep] 	(5,1)--(4,1)
				(5,2)--(4,2)
				(5,3)--(4,3)
				(5,4)--(4,4);
\node at (0.5,0.5) { ${V}_{SW}$ };
\node at (4.5,0.5) { ${V}_{SE}$ };
\node at (0.5,4.5) { ${V}_{NW}$ };
\node at (4.5,4.5) { ${V}_{NE}$ };
\node at (1.5,0.5) { $B_S^{(1)}(x_1)$ };
\node at (2.5,0.5) { $\ldots$ };
\node at (3.5,0.5) { $B_S^{(P')}(x_{P'})$ };
\node at (1.5,4.5) { $B_N^{(1)}(y_1)$ };
\node at (2.5,4.5) { $\ldots$ };
\node at (3.5,4.5) { $B_N^{(P')}(y_{P'})$ };
\node at (0.5,1.5) { $B_W^{(1)}(w_1)$ };
\node at (0.5,2.5) { $\vdots$ };
\node at (0.5,3.5) { $B_W^{(Q')}(w_{Q'})$ };
\node at (4.5,1.5) { $B_E^{(1)}(z_1)$ };
\node at (4.5,2.5) { $\vdots$ };
\node at (4.5,3.5) { $B_E^{(Q')}(z_{Q'})$ };
\end{tikzpicture}
\end{equation}
obtained from the ROPE representation of $g$ through the set of equations \eqref{eq:def:boundaryobjectsrestrictions} and \eqref{eq:decompboundaryeltoROPErep} below.

If, moreover, the face weights are constant, all equal to a face weight $\MarkovWeight{W}$ and the ROPErep of $g_R$ is homogeneous, i.e. all the operators $A_a^{(i)}(u_i)$ are equal to some $A_a(u_i)$, then the previous ROPErep of $g_{R'}$ is also homogeneous, i.e. the operators $B_a^{(j)}(u'_j)$ are equal to some operators $B_a(u'_j)$.
\end{theo}

\begin{rema}[change of perspective]
 	The perspective behind this theorem is that the boundary weights, which are global objects on the boundaries related by global operations (summations on the edges between $R'$ and $R$) and thus difficult to study, are built out (through operadic structures) of local objects related by local operations easier to manipulate. More generally, one is invited, on one side, to read the properties of global objects, such as (long-range) correlation functions, (large size) partition functions for example, in the properties of local objects, and, on the other side, to manipulate as long as possible local objects.
\end{rema}

\begin{rema}
	Equations \eqref{eq:def:boundaryobjectsrestrictions} and \eqref{eq:decompboundaryeltoROPErep} below are explicit and express the $B_a^{(j)}(u'_j)$ and $V_{ab}$ from the $A_a^{(i)}(u_i)$ and $U_{ab}$ using tensor products and linear algebra. This is only a first tool with some drawbacks listed below. Section~\ref{sec:invariantboundaryelmts} provides further tools to get around these drawbacks.
\end{rema}

\begin{figure}
\begin{center}
\begin{tikzpicture}[scale=0.8]
\coordinate (a30) at (0,0);
\coordinate (a01) at (10,0);
\coordinate (a12) at (10,8);
\coordinate (a23) at (0,8);
\coordinate (b30) at (2,2);
\coordinate (b01) at (7,2);
\coordinate (b12) at (7,7);
\coordinate (b23) at (2,7);

\fill[color=gray!40!white] (0,0) -- (10,0) -- (10,2) -- (0,2) -- cycle;
\fill[color=gray!40!white] (0,8) -- (10,8) -- (10,7) -- (0,7) -- cycle;
\fill[color=gray!40!white] (0,2) -- (2,2) -- (2,7) -- (0,7) -- cycle;
\fill[color=gray!40!white] (7,2) -- (10,2) -- (10,7) -- (7,7) -- cycle;

\draw[ultra thick,red] (a30) -- (a01) -- (a12) -- (a23) -- (a30);
\draw[thick] (b30) -- (b01) -- (b12) -- (b23) -- (b30);

\draw[dashed] (2,0) -- (2,2) -- (0,2);
\draw[dashed] (2.5,0) -- (2.5,2) ;
\draw[dashed] (3,0) -- (3,2) ;
\draw[dashed] (2,2.5) -- (0,2.5);
\draw[dashed] (2,3) -- (0,3);

\draw[dashed] (7,0) -- (7,2) -- (10,2);
\draw[dashed] (6,0) -- (6,2) ;
\draw[dashed] (6.5,0) -- (6.5,2);
\draw[dashed] (7,2.5) -- (10,2.5);
\draw[dashed] (7,3) -- (10,3);

\draw[dashed] (7,8) -- (7,7) -- (10,7);
\draw[dashed] (6.5,8) -- (6.5,7) ;
\draw[dashed] (6,8) -- (6,7) ;
\draw[dashed] (7,6.5) -- (10,6.5);
\draw[dashed] (7,6) -- (10,6);

\draw[dashed] (2,8) -- (2,7) -- (0,7);
\draw[dashed] (2.5,8) -- (2.5,7);
\draw[dashed] (3,8) -- (3,7);
\draw[dashed] (2,6.5) -- (0,6.5);
\draw[dashed] (2,6) -- (0,6);

\node at (a30) [anchor=north east] {$(X_1,Y_1)$};
\node at (a30) [inner sep=0.7mm,red,circle,fill] {};
\node at (a01) [anchor=north west] {$(X_2,Y_1)$};
\node at (a01) [inner sep=0.7mm,red,circle,fill] {};
\node at (a12) [anchor=south west] {$(X_2,Y_2)$};
\node at (a12) [inner sep=0.7mm,red,circle,fill] {};
\node at (a23) [anchor=south east] {$(X_1,Y_2)$};
\node at (a23) [inner sep=0.7mm,red,circle,fill] {};
\node at (b30) [anchor=south west] {$(X'_1,Y'_1)$};
\node at (b01) [anchor=south east] {$(X'_2,Y'_1)$};
\node at (b12) [anchor=north east] {$(X'_2,Y'_2)$};
\node at (b23) [anchor=north west] {$(X'_1,Y'_2)$};
\draw[<->] (0,-0.5) -- node [midway,below] {$N_1$} (2,-0.5); 
\draw[<->] (2,-0.5) -- node [midway,below] {$P'$} (7,-0.5); 
\draw[<->] (7,-0.5) -- node [midway,below] {$N_2$} (10,-0.5); 
\draw[<->] (-0.5,0) -- node [midway,left] {$M_1$} (-0.5,2); 
\draw[<->] (-0.5,2) -- node [midway,left] {$Q'$} (-0.5,7); 
\draw[<->] (-0.5,7) -- node [midway,left] {$M_2$} (-0.5,8); 
% nD-ALGEBRAS
\node at (4.5,1) [blue] {$\ldots B_S^{(i)}\ldots $};
\node at (1,1) [blue] {$V^{(SW)}$};
\node at (8.5,1) [blue] {$V^{(SE)}$};

\node at (1,4.5) [blue] {$\begin{matrix}\vdots \\ B_W^{(j)} \\\vdots\end{matrix}$};
\node at (8.5,4.5) [blue] {$\begin{matrix}\vdots \\ B_E^{(j)}\\\vdots\end{matrix}$};

\node at (4.5,7.5) [blue] {$\ldots B^{(i)}_N\ldots$};
\node at (1,7.5) [blue] {$V^{(NW)}$};
\node at (8.5,7.5) [blue] {$V^{(NE)}$};

\draw[<-,red] (10,4) -- (10.5,4) node [red,anchor=west] {$A^{(j)}_E$};
\draw[<-,red] (5,8) -- (5,8.5) node [red,anchor=south] {$A^{(i)}_N$};
\draw[<-,red] (0,4) -- (-1,4) node [red,anchor=east] {$A^{(j)}_W$};
\draw[<-,red] (5,0) -- (5,-1) node [red,anchor=north] {$A^{(i)}_S$};

\draw[<-,red] (10,0) -- (10.5,0.5) node [red, anchor = west] {$U^{(SE)}$};

\draw[<-,red] (10,8) -- (10.5,7.5) node [red, anchor = west] {$U^{(NE)}$};
\draw[<-,red] (0,8) -- (0.5,8.5) node [red, anchor = west] {$U^{(NW)}$};
\draw[<-,red] (0,0) -- (-1,0.5) node [red, anchor = east] {$U^{(SW)}$};

\end{tikzpicture}
\end{center}
\caption{\label{fig:nestedrectanglesparam}Parametrizations of the rectangles for the proof of Theorem \ref{theo:stability} (black: coordinates and lengths, red: elements of the ROPE $\ca{B}_{\PatternShapes(\patterntype{fp}^*)}$, blue: elements of the ROPE $\ca{B}'_{\PatternShapes(\patterntype{fp}^*)}$).}
\end{figure}

\begin{proof}
The proof relies on the definition of a larger $\Guill_2^{(\patterntype{fp}^*)}$-algebra $(\ca{E}_{p,q})$ which contains both the canonical structure $(\ca{T}_{p,q})$ of Section~\ref{sec:canonicalboundarystructure} and the ROPE $(\ca{B}_{p,q})$ and in which the statement of the theorem becomes trivial.

\begin{lemm}\label{lemm:fromROPErelelementtotensorprod}
	Let $(\ca{R}_{p,q})$ be a ROPE and $S$ a finite set. A function $\phi: S \to \ca{R}_{1,\infty_{S}}$ (resp. $S\to\ca{R}_{1,\infty_N}$, $S\to\ca{R}_{\infty_W,1}$ and $S\to\ca{R}_{\infty_E,1}$) defines injectively a canonical element of $V(S)^{\otimes r}\otimes \ca{R}_{r,\infty_{S}}$ (resp. ${V(S)^*}^{\otimes r}\otimes \ca{R}_{r,\infty_N}$, ${V(S)^*}^{\otimes r}\otimes \ca{R}_{\infty_W,r}$ and $V(S)^{\otimes r}\otimes \ca{R}_{\infty_E,r}$) for any $r\geq 1$ through
	\begin{align*}
	\sum_{(u_1,\ldots,u_r)\in S^r} e_{u_r}\otimesdots e_{u_r}\otimes 
	\begin{tikzpicture}[guillpart,xscale=2,yscale=1.5]
		\fill[guillfill] (0,0) rectangle (3,1);
		\draw[guillsep] (0,0)--(0,1)--(3,1)--(3,0) (1,0)--(1,1) (2,0)--(2,1);
		\node at (0.5,0.5) {$\phi(u_1)$};
		\node at (1.5,0.5) {$\dots$};
		\node at (2.5,0.5) {$\phi(u_r)$};
	\end{tikzpicture}
	\end{align*}
\end{lemm}

\begin{lemm}[canonical extension of a ROPErep]\label{lemm:ROPEcanonicalextension}
	Let $S_1$ and $S_2$ be two finite sets and $\ca{B}_{\PatternShapes(\patterntype{fp}^*)}$ be a ROPE. Let $\ca{B}'_{\PatternShapes(\patterntype{fp}^*)}$ be the collection of spaces defined by
	\begingroup
	\allowdisplaybreaks
	\begin{align*}
		\ca{B}'_{p,q} &= \setK 
		\\
		\ca{B}'_{p,\infty_S} &= \left(\oplus_{q\in\setN} \End(V(S_2))^{\otimes q} \right) \otimes \ca{B}_{p,\infty_S}
		&
		\ca{B}'_{p,\infty_N} &= \left(\oplus_{q\in\setN} \End(V(S_2))^{\otimes q} \right) \otimes \ca{B}_{p,\infty_N}
		\\
		\ca{B}'_{\infty_W,q} &= \left(\oplus_{p\in\setN} \End(V(S_1))^{\otimes p} \right) \otimes \ca{B}_{\infty_W,q}
		&
		\ca{B}'_{\infty_E,q} &= \left(\oplus_{p\in\setN} \End(V(S_1))^{\otimes p} \right) \otimes \ca{B}_{\infty_E,q}
		\\
		\ca{B}'_{\infty_W,\infty_S} &= \left(\oplus_{(p,q)\in\setN^2} {V(S_1)^*}^{\otimes p}\otimes {V(S_2)^*}^{\otimes q} \right) \otimes \ca{B}_{\infty_W,\infty_S}
		\\
		\ca{B}'_{\infty_E,\infty_S} &= \left(\oplus_{(p,q)\in\setN^2} {V(S_1)^*}^{\otimes p}\otimes {V(S_2)}^{\otimes q} \right) \otimes \ca{B}_{\infty_E,\infty_S}
		\\
		\ca{B}'_{\infty_W,\infty_N} &= \left(\oplus_{(p,q)\in\setN^2} {V(S_1)}^{\otimes p}\otimes {V(S_2)^*}^{\otimes q} \right) \otimes \ca{B}_{\infty_W,\infty_N}
		\\
		\ca{B}'_{\infty_E,\infty_N} &= \left(\oplus_{(p,q)\in\setN^2} {V(S_1)}^{\otimes p}\otimes {V(S_2)}^{\otimes q} \right) \otimes \ca{B}_{\infty_E,\infty_N}
	\end{align*}
	for finite $p$ and $q$ for bounded and one-sided infinite shapes and defined by
	\begin{align*}
		\ca{B}'_{\infty_{WE},q} &= \left( \oplus_{(p_1,p_2)\in\setP^2_\leq} \End(V(S_1))^{\otimes p_1}\otimes \End(V(S_1))^{\otimes p_2} \right)\otimes \ca{B}_{\infty_{WE},q}
		\\
		\ca{B}'_{\infty_{WE},\infty_S} &= \left( \oplus_{(p_1,p_2)\in\setP^2_\leq}  {V(S_1)^*}^{\otimes p_1}\otimes {V(S_1)^*}^{\otimes p_2} \right)\otimes \ca{B}_{\infty_{WE},\infty_S}
		\\
		\ca{B}'_{\infty_{WE},\infty_N} &= \left( \oplus_{(p_1,p_2)\in\setP^2_\leq} V(S_1)^{\otimes p_1}\otimes V(S_1)^{\otimes p_2} \right)\otimes \ca{B}_{\infty_{WE},\infty_N}
		\\
		\ca{B}'_{p,\infty_{SN}} &= \left( \oplus_{(q_1,q_2)\in\setP^2_\leq} \End(V(S_2))^{\otimes q_1}\otimes \End(V(S_2))^{\otimes q_2} \right)\otimes \ca{B}_{q,\infty_{SN}}
		\\
		\ca{B}'_{\infty_W,\infty_{SN}} &= \left( \oplus_{(q_1,q_2)\in\setP^2_\leq} {V(S_2)^*}^{\otimes q_1}\otimes {V(S_2)^*}^{\otimes q_2} \right)\otimes \ca{B}_{\infty_W,\infty_{SN}}
		\\
		\ca{B}'_{\infty_E,\infty_{SN}} &= \left( \oplus_{(q_1,q_2)\in\setP^2_\leq} V(S_2)^{\otimes q_1}\otimes V(S_2)^{\otimes q_2} \right)\otimes \ca{B}_{\infty_E,\infty_{SN}}
		\\
		\ca{B}'_{\infty_{WE},\infty_{SN}} &= \setK \otimes \ca{B}_{\infty_{WE},\infty_{SN}}\simeq \setK
	\end{align*}
	\endgroup
	for finite $p$ and $q$ and two-sided infinite shapes. The spaces are endowed with similar products as the canonical structures $\ca{T}_{\PatternShapes(\patterntype{fp}^*)}$ (excepted that there are no concatenation products on the missing spaces). The spaces $\ca{B}'_{\PatternShapes(\patterntype{fp}^*)}$ define a ROPE with a canonical morphism of ROPE $\Phi_{\PatternShapes(\patterntype{fp}^*)}: \ca{B}_{\PatternShapes(\patterntype{fp}^*)} \to \ca{B}'_{\PatternShapes(\patterntype{fp}^*)}$ defined for any element $A \in \ca{B}_{p,q}$ by $\Phi_{p,q}(A)= (1,0,0,\ldots)\otimes A$ where $1\in\setK$ is always placed in the $0$-components (resp. $(0,0)$-component) of the direct sums and $0$ in the other components for rectangles and half-strips (resp. for the other shapes).
\end{lemm}
\begin{proof}
	The proof follows exactly the same steps as for the canonical structure $\ca{T}_{p,q}$ with fewer terms (no concatenation of tensor products) or the same steps as for Proposition~\ref{prop:linkwithMPS} and is not repeated.
\end{proof}
\begin{lemm}\label{lemm:globalGuillalgebraEpq}
	Let $S_1$ and $S_2$ be two finite sets and $\ca{B}_{\PatternShapes(\patterntype{fp}^*)}$ be a ROPE. The spaces $(\ca{E}_{p,q})$ defined by 
	\[
	\ca{E}_{p,q} = \ca{T}_{p,q} \otimes \ca{B}_{p,q}
	\]
	form a $\Guill_2^{(\patterntype{fp}^*)}$-algebra as a tensor product of two $\Guill_2^{(\patterntype{fp}^*)}$-algebras and one has:
	\begingroup
	\allowdisplaybreaks
	\begin{align*}
		\ca{E}_{p,q} &= \ca{T}_{p,q} \otimes \ca{B}'_{p,q} \simeq \ca{T}_{p,q}\\
		\ca{E}_{p,\infty_S} &= {V(S_1)^*}^{\otimes p}\otimes \ca{B}'_{p,\infty_S} &
		\ca{E}_{p,\infty_N} &= V(S_1)^{\otimes p}\otimes \ca{B}'_{p,\infty_N} \\
		\ca{E}_{\infty_W,q} &= {V(S_2)^*}^{\otimes q}\otimes \ca{B}'_{\infty_W,q} &
		\ca{E}_{\infty_E,q} &= V(S_2)^{\otimes q}\otimes \ca{B}'_{\infty_E,q} 
	\end{align*}
	\endgroup
	for all finite integers $p$ and $q$ and $\ca{E}_{s} = \ca{B}'_{s}$
	for all other shapes.
\end{lemm}
\begin{proof}
	As a tensor product of two $\Guill_2^{(\patterntype{fp}^*)}$-algebras, $(\ca{E}_{p,q})$ is also a $\Guill_2^{(\patterntype{fp}^*)}$-algebra with products defined, for any (pointed or not) guillotine partition $\rho$ of size $n$, by:
	\begin{equation}
		m_{\rho}\left( \bigotimes_{k=1}^n t_{p_k,q_k}\otimes b_{p_k,q_k} \right) 
		= m_{\rho}\left( \bigotimes_{k=1}^n t_{p_k,q_k}\right) \otimes m_{\rho}\left(\bigotimes_{k=1}^n b_{p_k,q_k} \right) 
	\end{equation}
	for any sequences of elements $(t_{p_k,q_k})_k$ and $(b_{p_k,q_k})_k$ such that $t_{p_k,q_k} \in \ca{T}_{p_k,q_k}$ and $b_{p_k,q_k}\in \ca{B}_{p_k,q_k}$. 
	
	For rectangular shapes $(p,q)\in\PatternShapes(\patterntype{r})$, we have the following identification
	\[
	\ca{E}_{p,q} = \ca{T}_{p,q}\otimes \setR \simeq \ca{T}_{p,q}
	\]
	which is thus the same $\Guill_2^{(\patterntype{r})}$-algebra as the one associated to the Markov process.
\end{proof}

\textbf{Embeddings in $\ca{E}_{\PatternShapes(\patterntype{fp}^*)}$.}
We can now proceed to the proof of the theorem. First, we canonically embed any weight face $\MarkovWeight{w}\in\ca{T}_{1,1}$ in $\ca{E}_{\PatternShapes(\patterntype{fp}^*)}$ through
\[
\ha{\MarkovWeight{w}} = \MarkovWeight{w} \otimes 1 \in \ca{E}_{1,1}
\]

Using lemma \ref{lemm:fromROPErelelementtotensorprod} with $\ca{R}_{p,q}=\ca{B}'_{p,q}$ defined in lemma \ref{lemm:ROPEcanonicalextension} with the identification of Lemma~\ref{lemm:globalGuillalgebraEpq}, we canonically embed the ROPErep elements $A_a^{(i)}(u)$ and $U_{ab}$ in the larger $\Guill_2^{(\patterntype{fp}^*)}$-algebra   $\ca{E}_{\PatternShapes(\patterntype{fp}^*)}$. This corresponds to the following expressions
\begingroup
\allowdisplaybreaks
\begin{subequations}
\label{eq:fromROPEreptoEpq}
\begin{align}
\ha{A}_S^{(i)}  &= \sum_{x\in S_1} e_x^*\otimes (1,0,\ldots)  \otimes A_S^{(i)}(x) \in \ca{E}_{1,\infty_S} \label{eq:fromROPEreptoEpq:South}
\\
\ha{A}_N^{(i)}  &= \sum_{y\in S_1}e_y\otimes (1,0,\ldots) \otimes A_N^{(i)}(y) \in \ca{E}_{1,\infty_N}
\\
\ha{A}_W^{(i)}  &= \sum_{w\in S_2}e^*_w\otimes(1,0,\ldots)\otimes A_W^{(i)}(w) \in \ca{E}_{\infty_W,1}
\\
\ha{A}_E^{(i)}  &= \sum_{z\in S_2} e_z \otimes (1,0,\ldots) \otimes A_E^{(i)}(z) \in \ca{E}_{\infty_E,1}
\\
\ha{U}_{SW} &=  (1,0,\ldots)\otimes(1,0,\ldots) \otimes U_{SW} \in \ca{E}_{\infty_W,\infty_S}
\\
\ha{U}_{SE} &= (1,0,\ldots)\otimes(1,0,\ldots) \otimes U_{SE} \in \ca{E}_{\infty_E,\infty_S}
\\
\ha{U}_{NE} &= (1,0,\ldots)\otimes(1,0,\ldots)  \otimes U_{NE} \in \ca{E}_{\infty_E,\infty_N}
\\
\ha{U}_{NW} &= (1,0,\ldots)\otimes(1,0,\ldots) \otimes U_{NW} \in \ca{E}_{\infty_W,\infty_N}
\end{align}
\end{subequations}
\endgroup

\textbf{Boundary weights: from $g_R$ to $g_{R'}$.}
Following Lemma~\ref{lemma:proba:Markovrestriction}, the boundary weight $g_{R'}$ is obtained from $g_R$ by summation over the ranges of the r.v. on edges that do not appear in $R'$. However, it is easier (see paragraph~\ref{par:dualviewonweights}) to characterize $g_{R'}$ as the almost unique element of the dual $\ca{T}_{P',Q'}^*$ such that, for all $Y \in \ca{T}_{P',Q'}$,
\begin{align*}
\langle g_{R'}, Y \rangle_{P',Q'}
&= \langle g_{R}, \mathbf{Z}(Y) \rangle_{P,Q}
\\
\text{with } \mathbf{Z}(Y)
 &=
\begin{tikzpicture}[guillpart,xscale=6.,yscale=1.5]
\fill[guillfill] (0,0) rectangle (3,3);
\foreach \x in {0,1,...,3} {
	\draw[guillsep] (\x,0) -- (\x,3);
	\draw[guillsep] (0,\x) -- (3,\x);
}
\node at (0.5,0.5) { $Z_{[X_1,X'_1]\times [Y_1,Y'_1]}$ };
\node at (1.5,0.5) { $Z_{[X'_1,X'_2]\times [Y_1,Y'_1]}$ };
\node at (2.5,0.5) { $Z_{[X'_2,X_2]\times [Y_1,Y'_1]}$ };
\node at (0.5,1.5) { $Z_{[X_1,X'_1]\times [Y'_1,Y'_2]}$ };
\node at (1.5,1.5) { $Y$ };
\node at (2.5,1.5) { $Z_{[X'_2,X_2]\times [Y'_1,Y'_2]}$ };
\node at (0.5,2.5) { $Z_{[X_1,X'_1]\times [Y'_2,Y_2]}$ };
\node at (1.5,2.5) { $Z_{[X'_1,X'_2]\times [Y'_2,Y_2]}$ };
\node at (2.5,2.5) { $Z_{[X'_2,X_2]\times [Y'_2,Y_2]}$ };
\end{tikzpicture}\in \ca{T}_{P,Q}
\end{align*}
where $Z_{[a_1,a_2]\times [b_1,b_2]}$ is a shortcut notation for the partition functions $Z_{[a_1,a_2]\times [b_1,b_2]}(\MarkovWeight{W}_\bullet;\cdot) \in \ca{T}_{a_2-a_1,b_2-b_1}$. This correspondence is established by first considering $\Esp{\prod_{e\in \Edges{R'}}\indic{X_e=x_e}}$ for any sequence $(x_e)$ (hence the arbitrary choice of $\Omega$) and then computing this expectation value using both Proposition~\ref{prop:MarkovTwoDimSecond} and equation~\eqref{eq:proba:ZboundaryfromZdet}. As usual in probability theory, the "almost" unicity means here that $g_{R'}$ may take arbitrary values on boundary configurations $x_{\partial R'}$ that have a null probability. The previous equation is the two-dimensional equivalent of the one-dimensional situation where the boundary weights on a subsegment are obtained as action of the product of transition matrices on boundary vectors, excepted that the products are now replaced by products in the operad $\Guill_2$ using \eqref{eq:compopartitionfunctions}.

The evaluation $\langle \cdot, \cdot \rangle_{P,Q}$ of a dual element on an element of $\ca{T}_{P,Q}$ is just a summation over the values of boundary variables. Considering now the ROPErep \eqref{eq:ROPErepforboundaryweight} of $g_R$, it is easy to move from $\ca{T}_{\PatternShapes(\patterntype{fp}^*)}$ to $\ca{E}_{\PatternShapes(\patterntype{fp}^*)}$ using lemmata~\ref{lemm:ROPEcanonicalextension} and \ref{lemm:fromROPErelelementtotensorprod} and get:
\begin{align}\label{eq:weightinEpqbeforetensor}
\langle g_R, \mathbf{Z}(Y)\rangle_{P,Q} = & \sum_{x_{\partial R}} 
g_R(x_{\partial R}) \mathbf{Z}(Y)(x_{\partial R}) 
\\
= & \sum_{\substack{x,y\in S_1^P \\ w,z\in S_2^Q }}
\begin{tikzpicture}[guillpart,yscale=1.6,xscale=3.4]
\draw[guillsep, dotted] (0,0) rectangle (5,5);
\draw[guillsep] (1,1) rectangle (4,4);
\draw[guillsep] 	(1,0)--(1,1)
				(2,0)--(2,1)
				(3,0)--(3,1)
				(4,0)--(4,1);
\draw[guillsep] 	(1,5)--(1,4)
				(2,5)--(2,4)
				(3,5)--(3,4)
				(4,5)--(4,4);
\draw[guillsep] 	(0,1)--(1,1)
				(0,2)--(1,2)
				(0,3)--(1,3)
				(0,4)--(1,4);
\draw[guillsep] 	(5,1)--(4,1)
				(5,2)--(4,2)
				(5,3)--(4,3)
				(5,4)--(4,4);
\node at (0.5,0.5) { $U_{WS}$ };
\node at (4.5,0.5) { $U_{SE}$ };
\node at (0.5,4.5) { $U_{NW}$ };
\node at (4.5,4.5) { $U_{EN}$ };
\node at (1.5,0.5) { $A_S^{(1)}(x_1)$ };
\node at (2.5,0.5) { $\ldots$ };
\node at (3.5,0.5) { $A_S^{(P)}(x_P)$ };
\node at (1.5,4.5) { $A_N^{(1)}(y_1)$ };
\node at (2.5,4.5) { $\ldots$ };
\node at (3.5,4.5) { $A_N^{(P)}(y_P)$ };
\node at (0.5,1.5) { $A_W^{(1)}(w_1)$ };
\node at (0.5,2.5) { $\vdots$ };
\node at (0.5,3.5) { $A_W^{(Q)}(w_Q)$ };
\node at (4.5,1.5) { $A_E^{(1)}(z_1)$ };
\node at (4.5,2.5) { $\vdots$ };
\node at (4.5,3.5) { $A_E^{(Q)}(z_Q)$ };
\end{tikzpicture}\\
& \times 
 \begin{tikzpicture}[guillpart,yscale=1.6,xscale=2.8]
\fill[guillfill] (0,0) rectangle (5,5);
\draw[guillsep, dotted] (0,0) rectangle (5,5);
\draw[guillsep] (1,1) rectangle (4,4);
\draw[guillsep] 	(1,0)--(1,1)
				(2,0)--(2,1)
				(3,0)--(3,1)
				(4,0)--(4,1);
\draw[guillsep] 	(1,5)--(1,4)
				(2,5)--(2,4)
				(3,5)--(3,4)
				(4,5)--(4,4);
\draw[guillsep] 	(0,1)--(1,1)
				(0,2)--(1,2)
				(0,3)--(1,3)
				(0,4)--(1,4);
\draw[guillsep] 	(5,1)--(4,1)
				(5,2)--(4,2)
				(5,3)--(4,3)
				(5,4)--(4,4);
\node at (0.5,0.5) { $\tilde{u}_{WS}$ };
\node at (4.5,0.5) { $\tilde{u}_{SE}$ };
\node at (0.5,4.5) { $\tilde{u}_{NW}$ };
\node at (4.5,4.5) { $\tilde{u}_{EN}$ };
\node at (1.5,0.5) { $\tilde{e}^S_{x_1}$ };
\node at (2.5,0.5) { $\ldots$ };
\node at (3.5,0.5) { $\tilde{e}^S_{x_P}$ };
\node at (1.5,4.5) { $\tilde{e}^N_{y_1}$ };
\node at (2.5,4.5) { $\ldots$ };
\node at (3.5,4.5) { $\tilde{e}^N_{y_P}$ };
\node at (0.5,1.5) { $\tilde{e}^W_{w_1}$ };
\node at (0.5,2.5) { $\vdots$ };
\node at (0.5,3.5) { $\tilde{e}^W_{w_Q}$ };
\node at (4.5,1.5) { $\tilde{e}^E_{z_1}$ };
\node at (4.5,2.5) { $\vdots$ };
\node at (4.5,3.5) { $\tilde{e}^E_{z_Q}$ };
\node at (2.5,2.5) { $\mathbf{Z}(Y)$ };
\end{tikzpicture}
\end{align}
where we chose once and for all arbitrarily the base point $(X_1,Y_1)$ and skip it from the graphical representations.

Then, from the definition of a ROPE as an extension of the $\Guill_2^{(\patterntype{r})}$-algebra $\setR_{p,q}$ and from the definition of the embeddings~\eqref{eq:fromROPEreptoEpq}, the previous equation can be recast in the larger $\Guill_2^{(\patterntype{fp}^*)}$-algebra $(\ca{E}_{p,q})$ through
\begin{equation}\label{eq:weightinEpq}
\langle g_{R}, \mathbf{Z}(Y) \rangle_{P,Q}
=
\begin{tikzpicture}[guillpart,xscale=2,yscale=1.7]
\fill[guillfill] (0,0) rectangle (5,5);
\draw[guillsep] (0,1) -- (5,1) (0,4)--(5,4) (1,0)--(1,5) (4,0)--(4,5);
\draw[guillsep] (2,0)--(2,1) (3,0)--(3,1)
	(2,4)--(2,5) (3,4)--(3,5)
	(0,2)--(1,2) (0,3)--(1,3)
	(4,2)--(5,2) (4,3)--(5,3);
\node at (2.5,2.5) { $\mathbf{Z}(Y)\otimes 1$ };
\node at (0.5,0.5) { $\ha{U}_{SW}$ };
\node at (1.5,0.5) { $\ha{A}_S^{(1)}$ };
\node at (2.5,0.5) { $\ldots$ };
\node at (3.5,0.5) { $\ha{A}_S^{(P)}$ };
\node at (4.5,0.5) { $\ha{U}_{SE}$ };
\node at (0.5,4.5) { $\ha{U}_{NW}$ };
\node at (1.5,4.5) { $\ha{A}_N^{(1)}$ };
\node at (2.5,4.5) { $\ldots$ };
\node at (3.5,4.5) { $\ha{A}_N^{(P)}$ };
\node at (4.5,4.5) { $\ha{U}_{NE}$ };
\node at (0.5,1.5) { $\ha{A}_W^{(1)}$ };
\node at (0.5,2.5) { $\ldots$ };
\node at (0.5,3.5) { $\ha{A}_W^{(Q)}$ };
\node at (4.5,1.5) { $\ha{A}_E^{(1)}$ };
\node at (4.5,2.5) { $\ldots$ };
\node at (4.5,3.5) { $\ha{A}_E^{(Q)}$ };
\end{tikzpicture}
\end{equation}
Using the associativity relations~\eqref{eq:guill2:listassoc} in the larger structure $\ca{E}_{\PatternShapes(\patterntype{fp}^*}$, we may now compose the guillotine partitions in another way. We introduce the following guillotine partition $\rho$ of size $5+2P+2Q$ of the full plane defined by:
\begingroup
\allowdisplaybreaks
\begin{align*}
\rho_c &= [X'_1,X'_2]\times [Y'_1,Y'_2] 
\\
\rho_{SW} &= (-\infty,X'_1]\times (-\infty,Y'_1]
&
\rho_{SE} &= [X'_2,+\infty)\times (-\infty,Y'_1]
\\
\rho_{NW} &= (-\infty,X'_1]\times [Y'_2,+\infty)
&
\rho_{NE} &= [X'_2,+\infty)\times [Y'_2,+\infty)
\\
\rho_{S,i} &= [X'_1+i-1,X'_1+i] \times (-\infty,Y'_1]
&
\rho_{N,i} &= [X'_1+i-1,X'_1+i] \times [Y'_2,+\infty)
\\
\rho_{W,j} &= (-\infty,X'_1] \times [Y'_1+j-1,Y'_j]
&
\rho_{E,j} &= [X'_2,+\infty) \times [Y'_1+j-1,Y'_j]
\end{align*}
\endgroup
for $1\leq i\leq P'$ and $1\leq j\leq Q'$. As an example, $\rho_{S,i}$ contains both the half-strip $[X'_1+i-1,X'_1+i] \times (-\infty,Y_1]$ already present in~\eqref{eq:weightinEpq} in which it is associated to $\ha{A}_S^{(i+N_1)}$ and the squares that lie on top of this half-strip and below $[X'_1,X'_2]\times[Y'_1,Y'_2]$. More generally, every element of $\rho$ excepted $\rho_c$ contains elements of the partition in~\eqref{eq:weightinEpq} excepted the central rectangle for $\mathbf{Z}(\Omega)$ with additional squares cut in the central rectangle. It is an easy exercise to check that this is indeed a guillotine partition by identifying the correct sequence of cuts.

Switching to this new partition $\rho$ and using the operadic structure of $\ca{E}_{\PatternShapes(\patterntype{fp}^*)}$ we obtain the following representation of $g_{R'}$:
\begin{align}\label{eq:bwrestrictedaction}
\langle g_{R'}, Y \rangle_{P',Q'} =
\langle g_R, \mathbf{Z}(Y)\rangle_{P,Q}= 
\begin{tikzpicture}[guillpart,xscale=2,yscale=1.8]
\fill[guillfill] (0,0) rectangle (5,5);
\draw[guillsep] 
	(0,1)--(5,1)
	(0,4)--(5,4)
	(1,0)--(1,5)
	(4,0)--(4,5);
\draw[guillsep]
	(2,0)--(2,1)
	(3,0)--(3,1)
	(2,4)--(2,5)
	(3,4)--(3,5)
	(0,2)--(1,2)
	(0,3)--(1,3)
	(4,2)--(5,2)
	(4,3)--(5,3)
	(2,0)--(2,1)
	(3,0)--(3,1);
\node at (2.5,2.5) {$Y\otimes 1$};
\node at (0.5,0.5) { $\ha{V}_{SW}$ };
\node at (1.5,0.5) { $\ha{B}_S^{(1)}$ };
\node at (2.5,0.5) { $\ldots$ };
\node at (3.5,0.5) { $\ha{B}_S^{(P')}$ };
\node at (4.5,0.5) { $\ha{V}_{SE}$ };
\node at (0.5,4.5) { $\ha{V}_{NW}$ };
\node at (1.5,4.5) { $\ha{B}_N^{(1)}$ };
\node at (2.5,4.5) { $\ldots$ };
\node at (3.5,4.5) { $\ha{B}_N^{(P')}$ };
\node at (4.5,4.5) { $\ha{V}_{NE}$ };
\node at (0.5,1.5) { $\ha{B}_W^{(1)}$ };
\node at (0.5,2.5) { $\ldots$ };
\node at (0.5,3.5) { $\ha{B}_W^{(Q')}$ };
\node at (4.5,1.5) { $\ha{B}_E^{(1)}$ };
\node at (4.5,2.5) { $\ldots$ };
\node at (4.5,3.5) { $\ha{B}_E^{(Q')}$ };
\end{tikzpicture}
\end{align}
with the following definitions of boundary side elements for $1\leq i\leq P'$ and $1\leq j\leq Q'$ attached to the shapes in the partition $\rho$:
\begin{subequations}
\label{eq:def:boundaryobjectsrestrictions}
\begingroup
\allowdisplaybreaks
\begin{align}
\label{eq:def:boundarySouthobjectrestriction}
\ha{B}^{(i)}_S  &=
\begin{tikzpicture}[guillpart,xscale=9.,yscale=1.6]
\fill[guillfill] (0,0) rectangle (1,2);
\draw[guillsep] (0,0)--(0,2)--(1,2)--(1,0)
	(0,1)--(1,1)	;
\node at (0.5,0.5) { $\ha{A}_S^{(N_1+i)} $ };
\node at (0.5,1.5) { $Z_{[X'_1+i-1,X'_1+i]\times [Y_1,Y'_1]}\otimes 1$ };
\end{tikzpicture}
\\
\ha{B}^{(i)}_N  &=
\begin{tikzpicture}[guillpart,xscale=9,yscale=1.6]
\fill[guillfill] (0,0) rectangle (1,2);
\draw[guillsep] (0,2)--(0,0)--(1,0)--(1,2)
	(0,1)--(1,1)	;
\node at (0.5,1.5) { $\ha{A}_N^{(N_1+i)} $ };
\node at (0.5,0.5) { $Z_{[X'_1+i-1,X'_1+i]\times [Y'_2,Y_2]}\otimes 1$ };
\end{tikzpicture}
\\
\ha{B}^{(j)}_W &=
\begin{tikzpicture}[guillpart,xscale=2.7,yscale=1.6]
\fill[guillfill] (0,0) rectangle (6,1);
\draw[guillsep] (0,0)--(6,0)--(6,1)--(0,1)
	(1,0)--(1,1)	;
\node at (0.5,0.5) { $\ha{A}_W^{(M_1+j)} $ };
\node at (3.5,0.5) { $Z_{[X_1,X'_1]\times [Y'_1+i-1,Y'_1+i]}\otimes 1$ };
\end{tikzpicture}
\\
\ha{B}^{(j)}_E &=
\begin{tikzpicture}[guillpart,xscale=2.7,yscale=1.6]
\fill[guillfill] (0,0) rectangle (6,1);
\draw[guillsep] (6,0)--(0,0)--(0,1)--(6,1)
	(5,0)--(5,1)	;
\node at (5.5,0.5) { $\ha{A}_E^{(M_1+j)} $ };
\node at (2.5,0.5) { $Z_{[X'_2,X_2]\times [Y'_1+i-1,Y'_1+i]}\otimes 1$ };
\end{tikzpicture}
\end{align}
\endgroup
and the following definitions of boundary corner elements:
\begin{align}
\ha{V}_{SW} &=
\begin{tikzpicture}[guillpart,yscale=1.6,xscale=3.5]
\fill[guillfill] (0,0) rectangle (4,4);
\draw[guillsep]
	(1,0)--(1,4)		(0,1)--(4,1)
	(4,0)--(4,4)		(0,4)--(4,4);
\draw[guillsep]
	(2,0)--(2,1)		(3,0)--(3,1)
	(0,2)--(1,2)		(0,3)--(1,3);
\node at (0.5,0.5) { $\ha{U}_{SW}$ };
\node at (1.5,0.5) { $\ha{A}_S^{(1)}$ };
\node at (2.5,0.5) { $\ldots$ };
\node at (3.5,0.5) { $\ha{A}_S^{(N_1)}$ };
\node at (0.5,1.5) { $\ha{A}_W^{(1)}$ };
\node at (0.5,2.5) { $\ldots$ };
\node at (0.5,3.5) { $\ha{A}_W^{(M_1)}$ };
\node at (2.5,2.5) { $Z_{[X_1,X'_1]\times[Y_1,Y'_1]}\otimes 1$ };
\end{tikzpicture}
\end{align}
\begin{align}
\ha{V}_{SE} &=
\begin{tikzpicture}[guillpart,yscale=1.6,xscale=3.5,rotate=90]
\fill[guillfill] (0,0) rectangle (4,4);
\draw[guillsep]
	(1,0)--(1,4)		(0,1)--(4,1)
	(4,0)--(4,4)		(0,4)--(4,4);
\draw[guillsep]
	(2,0)--(2,1)		(3,0)--(3,1)
	(0,2)--(1,2)		(0,3)--(1,3);
\node at (0.5,0.5) { $\ha{U}_{SE}$ };
\node at (1.5,0.5) { $\ha{A}_E^{(1)}$ };
\node at (2.5,0.5) { $\ldots$ };
\node at (3.5,0.5) { $\ha{A}_E^{(M_1)}$ };
\node at (0.5,1.5) { $\ha{A}_S^{(P)}$ };
\node at (0.5,2.5) { $\ldots$ };
\node at (0.5,3.5) { $\ha{A}_S^{(P-N_2+1)}$ };
\node at (2.5,2.5) { $Z_{[X'_2,X_2]\times[Y_1,Y'_1]}\otimes 1$ };
\end{tikzpicture}
\end{align}
\begin{align}
\ha{V}_{NE} &=
\begin{tikzpicture}[guillpart,yscale=1.6,xscale=3.5,rotate=180]
\fill[guillfill] (0,0) rectangle (4,4);
\draw[guillsep]
	(1,0)--(1,4)		(0,1)--(4,1)
	(4,0)--(4,4)		(0,4)--(4,4);
\draw[guillsep]
	(2,0)--(2,1)		(3,0)--(3,1)
	(0,2)--(1,2)		(0,3)--(1,3);
\node at (0.5,0.5) { $\ha{U}_{NE}$ };
\node at (1.5,0.5) { $\ha{A}_N^{(P)}$ };
\node at (2.5,0.5) { $\ldots$ };
\node at (3.5,0.5) { $\ha{A}_N^{(P-N_2+1)}$ };
\node at (0.5,1.5) { $\ha{A}_E^{(Q)}$ };
\node at (0.5,2.5) { $\ldots$ };
\node at (0.5,3.5) { $\ha{A}_E^{(Q-M_2+1)}$ };
\node at (2.5,2.5) { $Z_{[X'_2,X_2]\times[Y'_2,Y_2]}\otimes 1$ };
\end{tikzpicture}
\end{align}
\begin{align}
\ha{V}_{NW} &=
\begin{tikzpicture}[guillpart,yscale=1.6,xscale=3.5,rotate=270]
\fill[guillfill] (0,0) rectangle (4,4);
\draw[guillsep]
	(1,0)--(1,4)		(0,1)--(4,1)
	(4,0)--(4,4)		(0,4)--(4,4);
\draw[guillsep]
	(2,0)--(2,1)		(3,0)--(3,1)
	(0,2)--(1,2)		(0,3)--(1,3);
\node at (0.5,0.5) { $\ha{U}_{NW}$ };
\node at (1.5,0.5) { $\ha{A}_W^{(Q)}$ };
\node at (2.5,0.5) { $\ldots$ };
\node at (3.5,0.5) { $\ha{A}_W^{(Q-M_2+1)}$ };
\node at (0.5,1.5) { $\ha{A}_N^{(1)}$ };
\node at (0.5,2.5) { $\ldots$ };
\node at (0.5,3.5) { $\ha{A}_N^{(N_1)}$ };
\node at (2.5,2.5) { $Z_{[X_1,X'_1]\times[Y'_2,Y_2]}\otimes 1$ };
\end{tikzpicture}
\end{align}
\end{subequations}

Equation~\eqref{eq:bwrestrictedaction} is of the same type as \eqref{eq:weightinEpq} where the boundary elements $\ha{A}_a^{(i)}$ and $\ha{U}_{ab}$ were built out of the canonical boundary structure and a ROPE following~\eqref{eq:fromROPEreptoEpq}. We may now use again Lemma~\ref{lemm:fromROPErelelementtotensorprod} to extract a ROPErep of $g_{R'}$ over the ROPE $\ca{B}'_{\PatternShapes(\patterntype{fp}^*)}$. In particular, seen as a function $S_1^{P'}\times S_1^{P'} \times S_2^{Q'}\times S_2^{Q'}\to \setR$ for which we find a new ROPErep, the value $g_{R'}(x,y,w,z)$ is obtained from~\eqref{eq:bwrestrictedaction} with an element $Y \in \ca{T}_{P',Q'}$ given by:
\begin{equation}\label{eq:Yforselectionofvalue}
	Y = \bigotimes_{k=1}^{P'} (e_{x_k}\otimes e^*_{y_k}) \otimes \bigotimes_{k=1}^{Q'} (e_{w_k}\otimes e^*_{z_k})
\end{equation}
Each vector in this tensor product will be contracted with the corresponding dual vector component in the definitions of $\ca{E}_{p,q}$ with $(p,q)$ among $(1,\infty_{a})$ and $(\infty_b,1)$ and thus the operator part of $\ca{T}_{p,q}\otimes \ca{B}_{p,q}$ will be extracted and will participate to the new ROPE $\ca{B}'_{p,q}$, as stated in Lemma~\ref{lemm:globalGuillalgebraEpq}. We may thus write any element $\ha{B}_a^{(i)}$
as
\begin{subequations}
\label{eq:decompboundaryeltoROPErep}
\begin{align}
	\ha{B}_S^{(i)} &= \sum_{x\in S_1} e_x^* \otimes B_S^{(i)}(x)	&
	\ha{B}_N^{(i)} &= \sum_{y\in S_1} e_y \otimes B_N^{(i)}(x)	\\
	\ha{B}_W^{(i)} &= \sum_{w\in S_2} e_w^* \otimes B_W^{(i)}(x) &
	\ha{B}_E^{(i)} &= \sum_{z\in S_2} e_z \otimes B_E^{(i)}(x)
\end{align}
\end{subequations}
where the $B_a^{(i)}(x)$ belong to the new ROPE $\ca{B}'_{\PatternShapes(\patterntype{fp}^*)}$. Corner elements $\ha{V}_{ab}\in\ca{E}_{\infty_a,\infty_b}$ are then directly identified to the corner element $V_{ab}\in\ca{B}'_{\infty_a,\infty_b}$.

This proves the theorem of structural stability of ROPEreps of boundary weights with explicit formulae for the ROPEreps of the boundary weights of smaller rectangles.

If the ROPErep of $g_R$ is homogeneous and the face weights are constant, i.e. the elements $A_a^{(i)}(u)$ do not depend on the position index $i$, then the expressions of $\ha{B}_a^{(i)}(x)$ also do not depend on the position index $i$ and the ROPErep of $g_{R'}$ is also homogeneous.
\end{proof}

The ROPE $\ca{B}'_{\PatternShapes(\patterntype{fp}^*)}$ is defined in such a way that \emph{all} the weights on sub-rectangles can be embedded in it; however, for any given sub-rectangle $R'$, one can find another ROPE $\ca{B}^{(R')}_{\PatternShapes(\patterntype{fp}^*)}$ that provides a simpler ROPErep of $g_{R'}$ ---simpler in the sense that it corresponds to subspaces of the spaces $\ca{B}'_{p,q}$.

\begin{coro}[null block elimination]\label{coro:stab:smallerROPE}
Using the same notation as in Theorem~\ref{theo:stability} and the parametrization of figure~\ref{fig:nestedrectanglesparam}, the boundary weight $g_{R'}$ admits a ROPE over over the ROPE $\ca{B}^{(R')}_{\PatternShapes(\patterntype{fp}^*)}$ defined by
\begin{align*}
	\ca{B}^{(R')}_{p,q} &= \setK \\
	\ca{B}^{(R')}_{p,\infty_S} &= \End(V(S_2))^{\otimes M_1} \otimes \ca{B}_{p,\infty_S} &
	\ca{B}^{(R')}_{p,\infty_N} &= \End(V(S_2))^{\otimes M_2} \otimes \ca{B}_{p,\infty_N} \\
	\ca{B}^{(R')}_{\infty_W,q} &= \End(V(S_1))^{\otimes N_1} \otimes \ca{B}_{\infty_W,q} &
	\ca{B}^{(R')}_{\infty_E,q} &= \End(V(S_1))^{\otimes N_2} \otimes \ca{B}_{\infty_E,q} \\
	\ca{B}^{(R')}_{\infty_W,\infty_S} &= {V(S_1)^*}^{N_1}\otimes {V(S_2)^*}^{\otimes M_1} \otimes \ca{B}_{\infty_W,\infty_S} &
	\ca{B}^{(R')}_{\infty_W,\infty_N} &= {V(S_1)}^{N_1}\otimes {V(S_2)^*}^{\otimes M_2} \otimes \ca{B}_{\infty_W,\infty_N} \\
	\ca{B}^{(R')}_{\infty_E,\infty_S} &= {V(S_1)^*}^{N_2}\otimes {V(S_2)}^{\otimes M_1} \otimes \ca{B}_{\infty_E,\infty_S} &	\ca{B}^{(R')}_{\infty_E,\infty_N} &= {V(S_1)}^{N_2}\otimes {V(S_2)}^{\otimes M_2} \otimes \ca{B}_{\infty_E,\infty_N}
\end{align*}
for the non-doubly-infinite shapes and similar formulae for doubly-infinite shapes.  
\end{coro}
\begin{proof}
	We consider in a closer way formulae \eqref{eq:def:boundarySouthobjectrestriction} and observe, from the structure of the canonical spaces $\ca{T}_{p,q}$ that their tensor component in the direct sum $\oplus_{r\in\setN} \End(V(S_a))^{\otimes r}$ are zero in all the spaces $\End(V(S_a))^{\otimes r}$ except for one of them with an index $r_0$ given the suitable $N_i$ or $M_i$ in the parametrization of figure~\ref{fig:nestedrectanglesparam}. One considers then the projection morphisms $\oplus_{r\in\setN} \End(V(S_a))^{\otimes r} \to \End(V(S_a))^{\otimes r_0}$ and checks that they define morphisms of ROPE from $\ca{B}'_\bullet\to\ca{B}^{(R')}_\bullet$ when restricted to the elements $\ha{B}_a^{(i)}(u)$ and extended in a similar way to all the other pattern shapes.
\end{proof}
No computation other than removal of null components is involved in this simplification: this has a practical consequence of encoding all the $\ha{B}_a^{(i)}(u)$ in more compact way, in particular for numerical computations. It also prepares the definitions of Section~\ref{sec:invariantboundaryelmts} by extracting the relevant components of the various ROPEreps of boundary weights on sub-rectangles.

		\subsection{Interest and drawbacks}
		
Theorem~\ref{theo:stability} is both a cornerstone of the operadic approach and only a mere starting point for further practical computations. We explain now why.

The main interest of Theorem~\ref{theo:stability} is to emphasize the fact that matrix product states and Matrix Ans\"atze, which are particular but typical cases of ROPEreps, fit perfectly with two-dimensional Markov processes: by "perfectly" we mean that, as soon as one believes that the operadic approach, such as the $\Guill_2^{(\patterntype{r})}$-operad, is a good approach then the algebraic framework suggests naturally that boundary conditions should be described by objects related extended operads such as the $\Guill_2^{(\patterntype{fp}^*)}$-operad. The stability of the ROPErep structure w.r.t. to restrictions as illustrated in Theorem~\ref{theo:stability} is just the translation of this fact.

A critical mind could argue that expressions~\eqref{eq:def:boundaryobjectsrestrictions} are just an abstract rewriting of partial summations over averaged r.v. in statistical mechanics in order to make appear the strip partition functions such as $Z_{[a,a+1]\times [b,c]}$ and rearrange the products and tensor products: as such, it may not bring anything new. This is perfectly true. However, we argue that this is already the case in dimension one: one may just sum with recognizing matrix products and linear algebra but the algebraic perspective brings a lot of other important tools (such as diagonalization, eigenvectors, etc). The next section \ref{sec:invariantboundaryelmts} illustrates this point by introducing well-adapted algebraic tools to study "invariant" boundary conditions, as discussed in Section~\ref{sec:gibbs}, from a new operadic point of view.

A pragmatic critical mind could also argue that \eqref{eq:decompboundaryeltoROPErep} only hides the complexity of practical computations behind useful abstract representations. The practical expressions \eqref{eq:def:boundaryobjectsrestrictions} of the $\ha{B}_a^{(i)}(u)$ describe these operators in tensor products such as $\End(V(S_a))^{\otimes R} \otimes \ca{B}_{1,\infty_a}$ where the initial operators $A_a^{(i)}(u)$ belong to the factor $\ca{B}_{1,\infty_a}$: as the relative distances $N_i$ and $M_i$ grow, the dimensions of the tensor products grow exponentially and, as such, are not suitable for numerical computations.We invite the reader to proceed immediately to Section~\ref{sec:invariantboundaryelmts} to see how we introduce natural additional algebraic tools to transform these formulae into relevant computational and practical tools. In particular, the examples presented in Section~\ref{sec:trivialfactorizedcase} illustrate in a brief and convincing way that this proof of the theorem of structural stability is only a small part of the present approach.

\section{Combining ROPEreps to build new ROPEreps.}\label{sec:combinationofROPEreps}

ROPEreps of functions and their morphisms present similarities with representation theory ---with boundary operators acting on corner spaces and morphisms between representations--- but, as it will be seen below, the probabilistic interpretation of functions on rectangles as boundary weights of Markov processes on arbitrarily large rectangles will be provided by Gibbs measures on infinite domains, which are known to have a simplicial structure. The purpose of this section is to present elementary operations on ROPEreps that fill the gaps between the algebraic and probabilistic structures.

\subsection{Tensorizing ROPEreps and products of functions}
\begin{prop}
	Let $\ca{B}_{\PatternShapes(\patterntype{fp}^*)}$ and $\ca{B}'_{\PatternShapes(\patterntype{fp}^*)}$ be two ROPEs. Let $(f_i)_{i\in I}$ and $(g_i)_{i\in I}$ be two sequences of functions such that:
	\begin{itemize}
		\item for all $i\in I$, both $f_i$ and $g_i$ are maps $S_1^{p_i}\times S_1^{p_i}\times S_2^{q_i}\times S_2^{q_i} \to \setK$;
		\item the collection $(f_i)_{i\in I}$ admits a ROPErep over $\ca{B}_{\PatternShapes(\patterntype{fp}^*)}$ with operators $A_{a,1}$ for $a\in\{S,N,W,E\}$ and corner elements $U_{a,b}$ with $a\in\{S,N\}$ and $b\in\{W,E\}$;
		\item the collection $(g_i)_{i\in I}$ admits a ROPErep over $\ca{B}'_{\PatternShapes(\patterntype{fp}^*)}$ with operators $A'_{a,1}$ for $a\in\{S,N,W,E\}$ and corner elements $U'_{a,b}$ with $a\in\{S,N\}$ and $b\in\{W,E\}$.
	\end{itemize} 
	The spaces \[
	\ca{B}''_{p,q}=\ca{B}_{p,q}\otimes \ca{B}'_{p,q}\]
	for $(p,q)\in \PatternShapes(\patterntype{fp}^*)$ admits a canonical ROPE structure and the collection of product of functions $(f_i g_i)_{i\in I}$ admits a ROPErep over $\ca{B}''_{\PatternShapes(\patterntype{fp}^*)}$ with boundary operators and corner elements given by
	\begin{align*}
		A''_{a,1}(x) &= A_{a,1}(x)\otimes A'_{a,1}(x)
		\\
		U''_{a,b} &= U_{a,b}\otimes U'_{a,b} 
	\end{align*}
\end{prop}
\begin{proof}
	The ROPE structure of tensor products of ROPE is directly inherited from the $\Guill_2$-structure of tensor products of $\Guill_2$-algebras and the canonical usual identification $\setK\otimes\setK\simeq \setK$. One then writes $f_i(u)g_i(u)$ as a product of two values in $\setK$ identified either as $\setK\otimes \setK$ where each copy of $\setK$ is identified to $\ca{B}_{\infty_{WE},\infty_{SN}}$ and $\ca{B}'_{\infty_{WE},\infty_{SN}}$ or as an element of $\ca{B}''_{\infty_{WE},\infty_{SN}}$. The tensor product of $\Guill_2$-algebras then allows for making the tensor products enter into each boundary components.
\end{proof}

From a purely algebraic perspective, this property is similar to the tensorization of representations of algebras and the product property of characters. From a probabilistic property, multiplying a boundary weight $f_i$ by a function $g_i$ corresponds to a new measure for the Markov process with a density $g_i$ (localized on the boundary) w.r.t to the measure with a boundary weight $f_i$.

\subsection{Linear combination of ROPEreps and functions}
As for representation theory, one may be interested in the decomposition of a tensor product of ROPEreps into a (direct) sum of more elementary "irreducible" ROPEreps with multiplicity. On the probabilistic side, the simplicial structure of Gibbs measures allows for a decomposition of a Gibbs measure as a positive barycentre of \emph{extremal} Gibbs measures. The ROPEreps of collection of functions encompasses simultaneously these situations.

\begin{prop}
	Let $K$ be a finite set. Let $(\ca{B}^{(\alpha)}_{\PatternShapes(\patterntype{fp}^*)})_{\alpha\in K}$ be a collection of ROPEs. Let $(f^{(\alpha)}_i)_{i\in I,\alpha\in K}$ be a collection of functions such that:
	\begin{itemize}
		\item for all $i\in I$,  there exists integers $(p_i,q_i)$ such that all the functions $f^{(\alpha)}_i$ with $\alpha\in K$ are maps $S_1^{p_i}\times S_1^{p_i}\times S_2^{q_i}\times S_2^{q_i} \to \setK$;
		\item for all $\alpha\in K$, the collection $(f^{(\alpha)}_i)_{i\in I}$ admits a ROPErep over $\ca{B}^{(\alpha)}_{\PatternShapes(\patterntype{fp}^*)}$ with operators $A^{(\alpha)}_{a,1}$ for $a\in\{S,N,W,E\}$ and corner elements $U^{(\alpha)}_{a,b}$ with $a\in\{S,N\}$ and $b\in\{W,E\}$;
	\end{itemize} 
	The spaces $\ha{\ca{B}}_{\PatternShapes(\patterntype{fp}^*)}$ defined, for any $(p,q)\in\PatternShapes(\patterntype{fp}^*)$, by
	\[
	\ha{\ca{B}}_{p,q} = \bigoplus_{\alpha \in K} \ca{B}^{(\alpha)}_{p,q}
	\]
	excepted for the full plane and rectangular space for which we still have $\ha{\ca{B}}_{p,q} = \setK$
	admits a ROPE structure where the products are inherited from the component-wise ROPE products of each $\ca{B}^{(\alpha)}$ with two exceptions:
	\begin{itemize}
		\item an element $\lambda\in\ha{\ca{B}}_{p,q}$ for $(p,q)\in\PatternShapes(\patterntype{r})$ is identified with the element $\bigoplus_{\alpha\in K} \lambda$ 
		\item products with a range in $\bigoplus_{\alpha\in K} \ha{\ca{B}}_{\infty_{WE},\infty_{SN}}$ are further composed with the map
		\begin{align*}
		\bigoplus_{\alpha\in K} \ha{\ca{B}}_{\infty_{WE},\infty_{SN}} \mapsto \setK \\
		(\lambda_\alpha)_{\alpha\in K} \mapsto \sum_{\alpha\in K} \lambda_\alpha
		\end{align*}
	\end{itemize}
	The collection of functions $(\sum_{\alpha \in K} f^{(\alpha)}_i)_{i\in K}$ admits a ROPErep over $\ha{\ca{B}}_{\PatternShapes(\patterntype{fp}^*)}$ given by elementary operators and corner elements:
	\begin{align*}
		\ha{A}_{a,1}(x) &= \bigoplus_{\alpha\in K} A^{(\alpha)}_{a,1}(x)
		&
		\ha{U}_{a,b} &= \bigoplus_{\alpha\in K} U^{(\alpha)}_{a,b}
	\end{align*}
\end{prop}
\begin{proof}
A convenient way of obtaining the result is the introduction of the intermediate $\Guill_2$-algebra
\[
\ca{B}^{\oplus}_{p,q} = \bigoplus_{\alpha\in K} \ca{B}^{(\alpha)}_{p,q}
\]
for any $(p,q)\in\PatternShapes(\patterntype{fp}^*)$. This is trivially a $\Guill_2$-algebra: it is sufficient to verify all the axioms component-wise. This is not a ROPE since the spaces for rectangles and the full plane are not $\setK$ but $\oplus_{\alpha\in K} \setK$. 

Given the collection of functions $(f_i)$, we may introduce the functions $F_i : S_1^{p_i}\times S_1^{p_i} \times S_2^{q_i}\times S_2^{q_i} \to \oplus_{\alpha \in K} \setK$ defined by:
\[
F_i(u) = (f_i^{(\alpha)}(u))_{\alpha\in\setK}
\]
We have component-wise a generalized (in the sense that $\setK$ is replaced by $\setK^{|K|}$) ROPErep of the collection of function $(F_i)$. The function $\sum_{i} f_i$ is then obtained from $F_i$ by summing the coordinates. 

We then build a map $\ca{B}^{\oplus}_{p,q}\to \ha{\ca{B}}_{p,q}$ on the boundaries in the following way:
\begin{itemize}
	\item all the products with target space $\ca{B}^{\oplus}_{\infty_{WE},\infty_{SN}}$ are left-composed with the summation of coordinates
	\item all the other spaces excepted the ones with pattern type $\patterntype{r}$ are trivially identified between the two $\Guill_2$-spaces.
\end{itemize}
One checks that all the associativity conditions required by $\Guill_2$-structures are not modified under this change. Incorporating rectangles can be done by a simple coproduct structure $\ha{\ca{B}}_{p,q}=\setK$ is sent to $\ca{B}^{\oplus}$ by a canonical coproduct inclusion that maps $\lambda$ to the constant sequence $(\lambda)_{\alpha\in K}$. Again, one checks easily that all the required axioms of a ROPE are satisfied.

Again, it is interesting to notice that such a coproduct construction already exists in representation theory \emph{mutatis mutandis}.
\end{proof}

\chapter[Operadic eigen-elements and invariant boundary conditions]{Operadic generalized eigen-elements and invariant boundary conditions}\label{sec:invariantboundaryelmts}

		\section{Preliminary remarks}
		
\subsubsection{General purpose of the construction}

The previous sections have introduced a boundary structure based on matrix product states for boundary weights of Markov processes on the plane, generalizing the classical matricial framework in dimension $1$. We show in the present chapter how the classical 1D  linear algebra properties of Markov processes in the thermodynamic limit can be adapted in dimension 2, provided that the boundary weights are replaced by the elements of ROPE representations presented in the previous chapter. The major novelty of the present chapter is that the operadic structure allows for new definitions "up to morphisms" of eigenvectors and it culminates with Theorem~\ref{theo:eigenROPErep:invmeas}, which builds infinite-volume Gibbs measures from a purely algebraic perspective.

\subsection{Homogeneity assumptions and two-dimensional semi-groups}

We now make further assumptions for the whole sections \ref{sec:invariantboundaryelmts}. The first assumption is the \emph{bulk homogeneity} in the face weights.

\begin{defi}[2D-semi-group]\label{def:2Dsemigroup}
	Let $\ca{A}_{\PatternShapes(\patterntype{r})}$ be a $\Guill_2^{(\patterntype{r})}$-algebra. A 2D-semigroup is a sequence of elements $(\MarkovWeight{F}_{p,q})_{(p,q)\in \PatternShapes(\patterntype{r})}$ such that,
	\begin{enumerate}[(i)]
		\item for all $(p,q)\in\PatternShapes(\patterntype{r})$, $\MarkovWeight{F}_{p,q}\in\ca{A}_{p,q}$
		\item for all $p_1,p_2,q\in \setL^*$,
		\begin{equation}\label{eq:horizsemigroup}
			\begin{tikzpicture}[guillpart,yscale=1.5,xscale=2]
				\fill[guillfill] (0,0) rectangle (2,1);
				\draw[guillsep] (0,0) rectangle (2,1);
				\draw[guillsep] (1,0)--(1,1);
				\node at (0.5,0.5) {$\MarkovWeight{F}_{p_1,q}$};
				\node at (1.5,0.5) {$\MarkovWeight{F}_{p_2,q}$};
			\end{tikzpicture}
			=\MarkovWeight{F}_{p_1+p_2,q}
		\end{equation}
		\item for all $p,q_1,q_2\in\setL^*$,
		\begin{equation}\label{eq:vertsemigroup}
			\begin{tikzpicture}[guillpart,yscale=1.5,xscale=2]
				\fill[guillfill] (0,0) rectangle (1,2);
				\draw[guillsep] (0,0) rectangle (1,2);
				\draw[guillsep] (0,1)--(1,1);
				\node at (0.5,0.5) {$\MarkovWeight{F}_{p,q_1}$};
				\node at (0.5,1.5) {$\MarkovWeight{F}_{p,q_2}$};
			\end{tikzpicture}
			=\MarkovWeight{F}_{p,q_1+q_2}
		\end{equation}
	\end{enumerate}
\end{defi}
In the case of discrete space, the collection of weights is completely characterized by $\MarkovWeight{F}_{1,1}$ from which all the elements are given as surface powers $\MarkovWeight{F}_{p,q} = \MarkovWeight{F}_{1,1}^{[p,q]}$. One also checks that a 2D-semi-group generates a canonical $\Guill_2^{(\patterntype{r})}$-sub-algebra of $\ca{A}_{\PatternShapes(\patterntype{r})}$ with one-dimensional spaces $\setK \MarkovWeight{F}_{p,q}$ for any $(p,q)\in \PatternShapes(\patterntype{r})$.

\begin{assumption}\label{assumption:bulkhomogeneity}
	All through the present Section~\ref{sec:invariantboundaryelmts}, we will fix once and for all a given 2D-semi-group $(\MarkovWeight{F}_{p,q})_{(p,q)\in\PatternShapes(\patterntype{r})}$, written shortly $\MarkovWeight{F}_\bullet$.
\end{assumption}

The situations covered by Assumption~\ref{assumption:bulkhomogeneity} is the case of all the models in statistical mechanics for which the Hamiltonian is local and translation invariant and all the discrete field theories with local and translation-invariant Hamiltonians. By definition, it does not encompass the case of disordered models.

\begin{rema}
	We will also assume that all the ROPEreps are homogeneous, following definitions~\ref{def:ROPErep:FD}. The description of non-translational-invariant Gibbs measures is not covered in this section; however, the previous sections are still valid in this case.
\end{rema} 

\subsection{Considerations about the dimension one and the passage to dimension two}

\subsubsection{Some remarks about the well-known one-dimensional case.}
Homogeneous Markov processes in dimension one correspond to Markov weights $(\MarkovWeight{A}_e)$ equal to some given weight $\MarkovWeight{A}$.
 
On any connected domain $D$, the boundary weight $g_D(x,y)$ is often factorized in a left and a right component through
\begin{equation}\label{eq:1D:localtoglobal}
g_D(x,y)	= u_{D,1}(x) u_{D,2}(y)
\end{equation}
This factorization of a global object (in the sense that $g_D$ depends on the whole domain) into two local objects (in the sense that each $u_{D,i}$ is associated to a single point) is the easiest example of a local-to-global approach but is so simple that the local-to-global aspect is useless. A case in which the boundary weights are not factorized corresponds to a non-trivial barycentre of extremal Gibbs measures.

The function partition on a domain of length $L$ (global object) is given by $Z_L=\scal{u_{D,1}}{\MarkovWeight{A}^Lu_{D,2}}$ (gluing of local objects) and more generally the one-point correlation functions are given by $\scal{u_{D,1}}{\MarkovWeight{A}^{L_1}\MarkovWeight{B}\MarkovWeight{A}^{L_2}u_{D,2}}/Z_L$ where $\MarkovWeight{B}$ is another matrix related to some observable. The asymptotic free energy density then corresponds to the largest (in modulus) eigenvalue and the associated Perron-Frobenius eigenvectors can be used to define the Gibbs measure on the whole space $\setZ$: this is how the linear algebra structure enters the game and participates to the study of the large size properties of the models through local equations. The eigenvalue equation 
\begin{equation}\label{eq:1D:eigenval}
	\MarkovWeight{A}u=\Lambda u
\end{equation} is an example of a local equation (an edge and a point), whereas the asymptotic $\scal{u_{D,1}}{\MarkovWeight{A}^Lu_{D,2}} \simeq \Lambda^L c$ is an example of local-to-global principle.

Using the Perron-Frobenius left and right eigenvectors, the partition functions is exactly equal to $Z_L=c \Lambda^L$ and the one-point correlation functions is simplified to the local equation $\scal{u^L_\Lambda}{\MarkovWeight{B}u^R_{\Lambda}}/c$ and no asymptotics is required.

More generally, we can make the following trivial remark, from a strictly computational perspective, the use of eigenvectors of $\MarkovWeight{A}$ provides the following well-known simplification:
\begin{subequations}
	\label{eq:1Deigenval:replacements}
\begingroup
\allowdisplaybreaks
\begin{align}
	u^L_\lambda \MarkovWeight{A}^k \MarkovWeight{B}  &= \lambda^k u^L_\lambda \MarkovWeight{B}
	\\
	\MarkovWeight{B} \MarkovWeight{A}^l u^R_{\mu}  &= \mu^l  \MarkovWeight{B} u^R_\mu
	\\
	\scal{u^L_\lambda}{\MarkovWeight{A}^k\MarkovWeight{B}\MarkovWeight{A}^l u^R_{\mu}} &= \lambda^{k}\mu^{l}\scal{u^L_\lambda}{\MarkovWeight{B}u^R_{\mu}}
\end{align}
\endgroup
\end{subequations}
which are almost always used to replace the l.h.s. by the r.h.s (but rarely the converse).

\subsubsection{Equivalent computations in dimension two and some wishes.}\label{par:requiredsimplificationrules}

Previous Section~\ref{sec:boundaryalgebra} generalizes a local-to-global factorization~\eqref{eq:1D:localtoglobal} through the notion of ROPE representations of rectangular boundary weights $g_R$, in such a way that face weights $\MarkovWeight{W}_{p,q}$ act in a natural way on these local-to-global factorizations through the theorem of structural stability~\ref{theo:stability}. At this point, the missing part is a generalization of the eigenvalue equation \eqref{eq:1D:eigenval} with a face weight  $\MarkovWeight{F}\in\ca{T}_{1,1}$ (or its associated 2D-semi-group) acting on boundary elements $A_{a}(x)$ and corner elements $U_{ab}$ in a ROPE. We first list some requirements about such a wished definition of 2D-eigen-elements.

A typical computation in probability and statistical mechanics consists is the one of correlation functions. This can be performed by introducing suitable elements along the edges as seen in Section~\ref{sec:eckmanhilton} or, if one wishes to avoid these commutative spaces associated to segments, by replacing some face weights $\MarkovWeight{F}$ by face weights $\MarkovWeight{G}_k$ containing both the weight $\MarkovWeight{F}$ and the local observables of interest. This corresponds to the classical trick of statistical mechanics of introduction of impurities to study the pure system. This is also the idea behind $E_2$-algebras which corresponds to the insertion of non-trivial elements associated to subdomain of a domain. Correlation functions require the evaluation of quantities such as
\[
\begin{tikzpicture}[guillpart,yscale=1.,xscale=1.6]
	\fill[guillfill] (0,0) rectangle (12,10);
	\foreach \x in {1,2,...,11} 
	{
		\draw[guillsep] (\x,0) -- (\x,10);
	}
	\foreach \y in {1,2,...,9} 
	{
		\draw[guillsep] (0,\y) -- (12,\y);
	}
	\node at (4.5,3.5) {$\MarkovWeight{G}_2$};
	\node at (5.5,2.5) {$\MarkovWeight{G}_1$};
	\node at (7.5,4.5) {$\MarkovWeight{G}_4$};
	\node at (5.5,6.5) {$\MarkovWeight{G}_3$};
	\foreach \x in {1.5,2.5,...,10.5}
	{
		\node at (\x,0.5) {$A_S$};
		\node at (\x,9.5) {$A_N$};
	}
	\foreach \y in {1.5,2.5,...,8.5} 
	{
		\node at (0.5,\y) {$A_W$};
		\node at (11.5,\y) {$A_E$};
	}
	\node at (0.5,0.5) {$U_{SW}$};
	\node at (0.5,9.5) {$U_{NW}$};
	\node at (11.5,0.5) {$U_{SE}$};
	\node at (11.5,9.5) {$U_{NE}$};
\end{tikzpicture}
\]
where all the internal elements except those written are equal to $\MarkovWeight{F}$ and the boundary weights are made of the homogeneous ROPErep obtained from an up-to-now-hypothetical translation-invariant Gibbs measure and up-to-now-hypothetical eigen-elements in a suitable sense. It is then natural to expect the previous quantity to be equal to
\[
\Lambda^{60} s_{W}^3 s_{E}^3 s_{S}^6 s_{N}^6 \begin{tikzpicture}[guillpart,xscale=1.6]
	\fill[guillfill] (3,1) rectangle (9,8);
	\foreach \x in {4,5,...,8} 
	{
		\draw[guillsep] (\x,1) -- (\x,8);
	}
	\foreach \y in {2,3,...,7} 
	{
		\draw[guillsep] (3,\y) -- (9,\y);
	}
	\node at (4.5,3.5) {$\MarkovWeight{G}_2$};
	\node at (5.5,2.5) {$\MarkovWeight{G}_1$};
	\node at (7.5,4.5) {$\MarkovWeight{G}_4$};
	\node at (5.5,6.5) {$\MarkovWeight{G}_3$};
	\foreach \x in {4.5,5.5,...,7.5}
	{
		\node at (\x,1.5) {$A_S$};
		\node at (\x,7.5) {$A_N$};
	}
	\foreach \y in {2.5,3.5,...,6.5} 
	{
		\node at (3.5,\y) {$A_W$};
		\node at (8.5,\y) {$A_E$};
	}
	\node at (3.5,1.5) {$U_{SW}$};
	\node at (3.5,7.5) {$U_{NW}$};
	\node at (8.5,1.5) {$U_{SE}$};
	\node at (8.5,7.5) {$U_{NE}$};
\end{tikzpicture}
\]
where the exponent $60$ corresponds to the number of removed faces and each $s_a$ corresponds to the number of removed elements $A_a$. Moreover, if one of the $\MarkovWeight{G}_k$ is equal to $\MarkovWeight{F}$, one expects that this simplification should be pursued further. In the same way, we also may have some rows full of $\MarkovWeight{F}$ but not all of them. Considering all this choices leads to the following wished local "replacement" rules (not complete, see full definitions in the following subsections):
\begin{itemize}
	\item any sequence of half-strips on any side such as
	\begin{subequations}
		\label{eq:2D:wishedreplacements}
	\begin{equation}\label{eq:2D:wishedreplacements:first}
		\begin{tikzpicture}[guillpart,xscale=1.5]
			\fill[guillfill] (0,0) rectangle (5,1);
			\draw[guillsep] (0,0)--(5,0)--(5,1)--(0,1) (1,0)--(1,1) (2,0)--(2,1) (3,0)--(3,1) (4,0)--(4,1);
			\node at (0.5,0.5) {$A_W$};
			\node at (1.5,0.5) {$\MarkovWeight{F}$};
			\node at (2.5,0.5) {$\MarkovWeight{F}$};
			\node at (3.5,0.5) {$\dots$};
			\node at (4.5,0.5) {$\MarkovWeight{F}$};
		\end{tikzpicture}
	\end{equation}
	or 
	\begin{equation}\label{eq:2D:wishedreplacements:second}
		\begin{tikzpicture}[guillpart,xscale=1.5]
			\fill[guillfill] (0,0) rectangle (6,1);
			\draw[guillsep] (0,0)--(6,0)--(6,1)--(0,1) (1,0)--(1,1) (2,0)--(2,1) (3,0)--(3,1) (4,0)--(4,1) (5,0)--(5,1);
			\node at (0.5,0.5) {$A_W$};
			\node at (1.5,0.5) {$\MarkovWeight{F}$};
			\node at (2.5,0.5) {$\MarkovWeight{F}$};
			\node at (3.5,0.5) {$\dots$};
			\node at (4.5,0.5) {$\MarkovWeight{F}$};
			\node at (5.5,0.5) {$\MarkovWeight{G}$};
		\end{tikzpicture}
	\end{equation}
	on the West with $k$ elements $\MarkovWeight{F}$ may be "replaced" by a sequence of reduced half-strips with a number $l\leq k$ of elements $\MarkovWeight{F}$, with an overall factor $\Lambda^{k-l}$, when surrounded by similar half-strip elements and suitable corner elements;
	
	\item any element on any corner such as
	\begin{equation}\label{eq:2D:wishedreplacements:corner}
		\begin{tikzpicture}[guillpart,yscale=1.25,xscale=1.6]
			\fill[guillfill] (0,0) rectangle (4,4);
			\draw[guillsep] (0,0)--(4,0)--(4,4) (1,0)--(1,4) (2,0)--(2,4) (3,0)--(3,4) (0,1)--(4,1) (0,2)--(4,2) (0,3)--(4,3);
			\node at (0.5,3.5) {$U_{NW}$};
			\node at (0.5,2.5) {$A_{W}$};
			\node at (0.5,1.5) {$\vdots$};
			\node at (0.5,0.5) {$A_{W}$};
			\node at (1.5,3.5) {$A_{N}$};
			\node at (2.5,3.5) {$\dots$};
			\node at (3.5,3.5) {$A_{N}$};
			\node at (1.5,0.5) {$\MarkovWeight{F}$};
			\node at (3.5,0.5) {$\MarkovWeight{F}$};
			\node at (1.5,2.5) {$\MarkovWeight{F}$};
			\node at (3.5,2.5) {$\MarkovWeight{F}$};
			\node at (2.5,0.5) {$\dots$};
			\node at (2.5,2.5) {$\dots$};
			\node at (2.5,1.5) {$\ddots$};
			\node at (1.5,1.5) {$\vdots$};
			\node at (3.5,1.5) {$\vdots$};
		\end{tikzpicture}
	\end{equation}
	on the North-West corner 
	with $k_W$ elements $A_W$ and $k_N$ elements $A_N$ may be "replaced" by a reduced corner element with $l_W<k_W$ and $l_N<k_N$ elements on each boundary, with an overall factor $\Lambda^{(k_W-l_W)(k_N-l_N)}s_W^{k_W-l_W}s_N^{k_N-l_N}$, when surrounded by suitable homogeneous half-strip elements;
	
	\item any row or column such as
	\begin{equation}\label{eq:2D:wishedreplacements:strip}
		\begin{tikzpicture}[guillpart,xscale=1.3]
			\fill[guillfill] (0,0) rectangle (10,1);
			\draw[guillsep] (0,0)--(10,0) (0,1)--(10,1) ;
			\foreach \x in {1,2,3,4,6,7,8,9} {
				\draw[guillsep] (\x,0)--(\x,1);
			}
			\node at (0.5,0.5) {$A_W$};
			\node at (9.5,0.5) {$A_E$};
			\foreach \x in {1.5,3.5,6.5,8.5} {
				\node at (\x,0.5) {$\MarkovWeight{F}$};	
			}
			\node at (2.5,0.5) {$\dots$};
			\node at (7.5,0.5) {$\dots$};
			\node at (5,0.5) {$\MarkovWeight{G}$};
		\end{tikzpicture}
	\end{equation}
	\end{subequations}		
	with $k$ elements $\MarkovWeight{F}$ on the left and $l$ on the right and \emph{any} element $\MarkovWeight{G}\in \ca{A}_{r,1}$ with arbitrary $r\in\setL$ may be "replaced" by a similar row or column with fewer elements $\MarkovWeight{F}$ on the left and on the right and an unchanged element $\MarkovWeight{G}$.
\end{itemize}

These replacements are similar to the ones in \eqref{eq:1Deigenval:replacements} but the main novelty here is that these "replacement" cannot proceed from simple equalities as in dimension one. This can be seen for Markov processes directly from the constructions in the proof of the theorem of structural stability~\ref{theo:stability}, in which the elements such as the ones in~\eqref{eq:2D:wishedreplacements} do not belong to the same spaces for different lengths $k$, $l$, etc. 

Instead of equalities, we shall rather introduced suitable "replacement" maps $\Phi_{\bullet}$ on the boundary spaces (half-strips, corners, strips, half-plane and plane) which maps the elements listed above to their simplified value with the correct multiplication factors, with the additional constraints that they must preserve the guillotine products and their associativities.

Before providing the general definitions, we start with a review of two trivial Markov processes where this new mechanism of morphisms of boundary eigen-structures instead of equalities~\eqref{eq:1Deigenval:replacements} is explicit and obvious. We also describe some generic features of morphisms related to the increase of dimension from $1$ to $2$ (and which can be generalized also to larger dimension trivially).

\section{A warm-up with two trivial Markov processes in dimension two}\label{sec:trivialfactorizedcase}

\subsection{The first example: independent horizontal and vertical lines.}

\subsubsection{Definitions}
We consider finite sets $S_1$ and $S_2$ and two square matrices $\MarkovWeight{A}$ and $\MarkovWeight{B}$, with respective sizes $|S_1|\times |S_1|$ and $|S_2|\times |S_2|$. We index them respectively by elements of $S_1$ and $S_2$. We moreover assume that $\MarkovWeight{A}$ and $\MarkovWeight{B}$ contain only positive entries and are irreducible, i.e. for any $(x,y)\in S_1\times S_1$ (resp. $S_2\times S_2)$, there exists $n\in\setN$ such that $(\MarkovWeight{A}^n)(x,y)>0$ (resp. $(\MarkovWeight{B}^n)(x,y)>0$).

We now introduce the face weight $\MarkovWeight{F}$ defined, for any $(x,y,w,z)\in S_1^2\times S_2^2$, by
\begin{equation}\label{eq:trivialmodelHV}
\MarkovWeight{F}(x,y,w,z) = \MarkovWeight{A}_{x,y} \MarkovWeight{B}_{w,z}
\end{equation}
and consider, on any rectangle $R=[x_1,x_2]\times [y_1,y_2]$, the two-dimensional Markov process $(X_e)_{e\in\Edges{R}}$ with all face weights equal to $\MarkovWeight{F}$. 

From a probabilistic point of view, this process can be decomposed into one-dimensional processes. We introduce, for any $x_1\leq x<x_2$, the process $(U^{(x)}_l)_{y_1\leq l \leq y_2}$ and, for any $y_1\leq y <y_2$, the processes $(V^{(y)}_k)_{x_1\leq k \leq x_2}$ defined by
\begin{align*}
U^{(x)}_l &= X_{ \{(x,l),(x+1,l)\} }
&
V^{(y)}_k &= X_{ \{(k,y),(k,y+1)\} }
\end{align*}
All the processes $U^{(x)}$ and $V^{(y)}$ are independent and, moreover, the processes $U^{(x)}$ (resp. $V^{(y)}$) are one-dimensional $S_1$-valued (resp. $S_2$-valued) Markov processes with edge weights equal to $\MarkovWeight{A}$ (resp. $\MarkovWeight{B}$). 

\subsubsection{Partition functions and boundary weights.}
The partition function on a rectangle $R=[x_1,x_2]\times [y_1,y_2]$ with shape $(P,Q)=(x_2-x_1,y_2-y_1)$ is given, for any $(x,y,w,z)\in S_1^{P}\times S_1^{P}\times S_2^{Q} \times S_2^{Q}$ (with the same boundary ordering as in \eqref{eq:boundaryordering}) by the product of the one-dimensional partitions functions:
\begin{equation}
Z_{R}(\MarkovWeight{F};x,y,w,z) = \prod_{i=1}^{P} (\MarkovWeight{A}^{L_2})(x_i,y_i) \prod_{j=1}^{Q} (\MarkovWeight{B}^{L_1})(w_j,z_j)
\end{equation}
Using the positivity and irreducibility hypotheses, Perron-Frobenius's theorem implies the existence of a real largest eigenvalue $\alpha$ (resp. $\beta$) of $\MarkovWeight{A}$ (resp. $\MarkovWeight{B}$) with left and right eigenvectors $a_l$ and $a_r$ (resp. $b_l$ and $b_r$) with strictly positive entries and the normalizations $\scal{a_l}{a_r}=\scal{b_l}{b_r}=1$. The double limit $P,Q\to\infty$ for rectangles with increasing sizes gives a free energy density 
\begin{equation}
\frac{1}{PQ}\log Z_{R} \to f= \log(\alpha\beta)
\label{eq:trivialmodel:freeenergy}
\end{equation}
for any "reasonable" sequence of boundary conditions (for example, constant ones).

Using the eigenvectors $a_l$, $a_r$, $b_l$ and $b_r$, it is easy to build a family of factorized boundary weights, for any rectangle $R$ with shape $(P,Q)$:
\begin{equation}\label{eq:horizvertmodel:boundaryweight}
G_{R}(x,y,w,z) = \prod_{i=1}^{P} a_l(x_i) a_r(y_i) \prod_{j=1}^{Q} b_l(w_j) b_r(z_j)
\end{equation}
such that:
\begin{enumerate}[(i)]
\item the partition function $Z^{\boundaryweights}_R(\MarkovWeight{F}; G_{R})$ is \emph{exactly} equal to $(\alpha\beta)^{L_1 L_2}$,
\item if $(X_e)_{e\in\Edges{R}}$ is a Markov process with boundary weight $G_{R}$, then, for \emph{any} sub-rectangle $R'$ of $R$, $(X_e)_{e\in\Edges{R'}}$ is a Markov process with boundary weight $G_{R'}$,
\item there exists, by Kolmogorov's extension theorem, a Markov process $(X_e)_{e\in \Edges{\setZ^2}}$ on the whole plane with marginal laws on any rectangle $R$ described by the boundary weight $G_R$, with free energy density equal to $f=\log(\alpha\beta)$.
\end{enumerate}

\subsubsection{ROPE representations.}
We now describe the algebraic counterpart in terms of ROPE representation. The boundary weights $G_R$ are given in a factorized form and thus the easiest homogeneous ROPE representation of the collection of weights $G_R$ is the trivial one introduce in section \ref{par:factorizedbw}:
\begin{equation}
\ca{B}_{p,q} = \setR
\end{equation}
for rectangles, half-strips and corners. The local operators in this ROPE are given by $A_S(x)= a_l(x)$, $A_N(y)=a_r(x)$, $A_W(w)= b_l(w)$ and $A_E(z)= b_r(z)$ and $U_{ab}=1$.

Following Theorem~\ref{theo:stability} and Corollary~\ref{coro:stab:smallerROPE}, the restriction from a rectangle $R=[x_1,x_2]\times[y_1,y_2]$ to a sub-rectangle $R'=[x'_1,x'_2]\times[y'_1,y'_2]$ produces boundary objects $B_S(x)$  in a ROPE $(\ca{B}^{(R')}_{p,q})$ (which depends on the relative positions $N_1$, $N_2$, $M_1$ and $M_2$ of the two rectangles in figure~\ref{fig:nestedrectanglesparam} but we skip them to have lighter notations) for the ROPErep of the marginal boundary weight $G_{R\to R'}$. Here, the ROPE $(\ca{B}^{(R')}_{p,q})$ is simply given on the boundary sides by:
\begingroup
\allowdisplaybreaks
\begin{align*}
\ca{B}^{(R')}_{p,\infty_S} &= \End(V(S_2))^{\otimes M_1}
&
\ca{B}^{(R')}_{p,\infty_N} &= \End(V(S_2))^{\otimes M_2}
\\
\ca{B}^{(R')}_{\infty_W,q} &= \End(V(S_1))^{\otimes N_1}
&
\ca{B}^{(R')}_{\infty_E,q} &= \End(V(S_1))^{\otimes N_2}
\\
\ca{B}^{(R')}_{\infty_W,\infty_S} &= {V(S_1)^*}^{\otimes N_1}\otimes {V(S_2)^*}^{\otimes M_1}
&
\ca{B}^{(R')}_{\infty_W,\infty_N} &= {V(S_1)}^{\otimes N_1}\otimes {V(S_2)^*}^{\otimes M_2}
\\
\ca{B}^{(R')}_{\infty_E,\infty_S} &= {V(S_1)^*}^{\otimes N_2}\otimes {V(S_2)}^{\otimes M_1}
&
\ca{B}^{(R')}_{\infty_E,\infty_N} &= {V(S_1)}^{\otimes N_2}\otimes {V(S_2)}^{\otimes M_2}
\end{align*}
and the operators $B_a(\bullet)$ and the corner elements can be evaluated directly as:
\begin{subequations}
\label{eq:trivialcase:boundaryaction}
\begin{align}
 B_S(x) &= \alpha^{M_1} a_l(x) \MarkovWeight{B}^{\otimes M_1} 
 &
 B_N(x) &= \alpha^{M_2} a_r(x) \MarkovWeight{B}^{\otimes M_2}
 \\
 B_W(x) &= \beta^{N_1} b_l(x) \MarkovWeight{A}^{\otimes N_1}
 &
 B_E(x) &= \beta^{N_2} b_r(x) \MarkovWeight{A}^{\otimes N_2}
 \\
 V_{SW} &= (\alpha\beta)^{N_1 M_1} a_l^{\otimes N_1}\otimes b_l^{\otimes M_1} 
 & 
 V_{NW} &= (\alpha\beta)^{N_1 M_2} a_r^{\otimes N_1}\otimes b_l^{\otimes M_2}
 \\
 V_{SE} &= (\alpha\beta)^{N_2 M_1} a_l^{\otimes N_2}\otimes b_r^{\otimes M_1}
 &
 V_{NE} &= (\alpha\beta)^{N_2 M_2} a_r^{\otimes N_2}\otimes b_r^{\otimes M_2}
\end{align}
\end{subequations}
\endgroup
The new objects $B_a(u)$ and $V_{ab}$ are much larger (in dimension) than the initial ones $A_a(u)$ and $U_{ab}$. However, they provide the \emph{same} boundary weight $G=_{R\to R'}\propto G_{R'}$ as~\eqref{eq:horizvertmodel:boundaryweight} up to a global constant $\Lambda^{PQ-P'Q'}$ where $(P,Q)$ and $(P',Q')$ are the shapes of $R$ and $R'$. This indicates that $B_a(u)$ and $V_{ab}$ are uselessly complicated for the computations of the boundary weight $G_{R'}$.  

Indeed, on each side, the operators $(B_a(x))_{x\in S_i}$ generate a commutative subalgebra of the corresponding $\End(V(S_i))^{\otimes K}$ and the $V(S_i)^{\otimes K}$-factor of $V_{ab}$ is a (Perron-Frobenius) eigenvector of the $(B_a(x))_{x\in S_i}$. All the products of elements of the ROPErep of $G_{R'}$ with elements $B_a(u)$ and $V_{ab}$ use indeed only a one-dimensional subspace of the whole space. Mapping this one-dimensional subspaces to the initial one-dimensional ROPE $\ca{B}_{p,q}$ then provides ROPE morphisms between the ROPEreps on $\ca{B}$ and $\ca{B}^{(R')}$.

In practice, we introduce now the rank one projector ${P}_\alpha$ (resp. $P_\beta$) onto the eigenspace of $\MarkovWeight{A}$ (resp. $\MarkovWeight{B}$) with eigenvalue $\alpha$ (resp. $\beta$) identified with $\setK$ itself. We then define the following linear maps on the tensor products on the half-strip elements:
\begin{subequations}
	\label{eq:horizvertmodel:morphisms}	
\begingroup
\allowdisplaybreaks
\begin{align}
\Phi^{S,N_2,0}_{p} : \ca{B}'_{p,\infty_S}  & \to \ca{B}_{p,\infty_S} \simeq \setK
&
\otimes_{i=1}^{N_2} C_i &\mapsto \prod_{i=1}^{N_2} ({P}_\beta C_i {P}_\beta)
\\
\Phi^{N,M_2,0}_{p} : \ca{B}'_{p,\infty_N}  & \to \ca{B}_{p,\infty_N} \simeq \setK
&
\otimes_{i=1}^{M_2} C_i &\mapsto \prod_{i=1}^{M_2} ({P}_\beta C_i {P}_\beta)
\\
\Phi^{W,N_1,0}_{q} : \ca{B}'_{\infty_W,q}  & \to \ca{B}_{\infty_W,q} \simeq \setK 
&
\otimes_{i=1}^{N_1} C_i &\mapsto \prod_{i=1}^{N_1} ({P}_\alpha C_i {P}_\alpha)
\\
\Phi^{E,M_1,0}_{q} : \ca{B}'_{\infty_E,q}  & \to \ca{B}_{\infty_E,q} \simeq \setK
&
\otimes_{i=1}^{M_1} C_i &\mapsto \prod_{i=1}^{M_1} ({P}_\alpha C_i {P}_\alpha)
\end{align}
where the notations correspond to Definition~\ref{def:eigenalgebrauptomorphims} below. We also introduce the following linear maps on the corner spaces:
\begin{align}
	\Phi_{\infty_W,\infty_S}^{SW,(N_1,M_1),(0,0)} : \ca{B}^{(R')}_{\infty_W,\infty_S}
	&\to
	\ca{B}^{\infty_W,\infty_S} \simeq \setK
	&
	u\otimes v &\mapsto uP_\alpha^{\otimes N_1}\otimes vP_\beta^{\otimes M_1}
\\
	\Phi_{\infty_W,\infty_N}^{NW,(N_1,M_2),(0,0)} : \ca{B}^{(R')}_{\infty_W,\infty_N}
	&\to
	\ca{B}^{\infty_W,\infty_N} \simeq \setK
	&
	u\otimes v &\mapsto P_\alpha^{\otimes N_1}u\otimes vP_\beta^{\otimes M_2}
\\
	\Phi_{\infty_E,\infty_S}^{ES,(N_2,M_1),(0,0)} : \ca{B}^{(R')}_{\infty_E,\infty_S}
	&\to
	\ca{B}^{\infty_E,\infty_S} \simeq \setK
	&
	u\otimes v &\mapsto uP_\alpha^{\otimes N_2}\otimes P_\beta^{\otimes M_1} v
\\
	\Phi_{\infty_E,\infty_N}^{EN,(N_2,M_2),(0,0)} : \ca{B}^{(R')}_{\infty_E,\infty_N}
	&\to
	\ca{B}^{\infty_E,\infty_N} \simeq \setK
	&
	u\otimes v &\mapsto P_\alpha^{\otimes N_2}u\otimes P_\beta^{\otimes M_2}v
\end{align}
\endgroup
\end{subequations}
where the notations are chosen in order to match with the ones of Proposition~\ref{prop:eigencorner:commutmorphisms} below.

We now have the following results about the relations between the two ROPEreps of the boundary weight $g_R$ and $g_{R'}$.
\begin{prop}
	For any relative positions $N_1,N_2,M_1,M_2$ of the rectangle $R'$ inside $R$, the restriction of the \emph{linear maps} $\Phi_\bullet^\bullet$ defined in equations~\eqref{eq:horizvertmodel:morphisms} to the subspaces of $\ca{B}^{(R')}_{\PatternShapes(\patterntype{fp}^*)}$ generated by the elements $B_a(u)$ and $V_{ab}$ lifts to \emph{morphism of ROPEs} such that:
	\begin{itemize}
		\item on the boundary side elements, one has
		\begin{subequations}
			\label{eq:horizvertlines:lambdamorph}
			\begin{align}
				\Phi^{M_1,0}_{1,\infty_S}( B_S(x) ) &= (\alpha\beta)^{M_1} A_S(x) 
				& 
				\Phi^{M_2,0}_{1,\infty_N}( B_N(x) ) &= (\alpha\beta)^{M_2} A_N(x) 
				\\
				\Phi^{N_1,0}_{\infty_W,1}( B_W(x) ) &= (\alpha\beta)^{N_1} A_W(x) 
				&
				\Phi^{N_2,0}_{\infty_E,1}( B_E(x) ) &= (\alpha\beta)^{N_2} A_E(x) 
			\end{align}

		\item on the corner elements, one has
			\begin{align}
				\Phi^{(N_1,M_1),(0,0)}_{\infty_W,\infty_S}(V_{SW}) &= \Lambda^{N_1M_1} U_{SW}
				&
				\Phi^{(N_1,M_2),(0,0)}_{\infty_W,\infty_N}(V_{WN}) &= \Lambda^{N_1M_2} U_{NW}
				\\
				\Phi^{(N_2,M_1),(0,0)}_{\infty_E,\infty_S}(V_{SE}) &= \Lambda^{N_2 M_1} U_{SE}
				&
				\Phi^{(N_2,M_2),(0,0)}_{\infty_E,\infty_N}(V_{NE}) &= \Lambda^{N_2M_2} U_{NE}
			\end{align}
		\end{subequations}
	
		\item the boundary weight $G_{R\to R'}$ with a ROPErep made of the $B_a(u)$ and $V_{ab}$ obtained in Theorem~\ref{theo:stability} is equal to
		\begin{equation}\label{eq:horizvertmodel:eigenweight}
			G_{R \to R'}(x,y,w,z) = \Lambda^{PQ-P'Q'} G_R(x,y,w,z)
		\end{equation}
		where $PQ-P'Q'$ is the area of $R\setminus R'$.
	\end{itemize}
\end{prop}
The "morphism of ROPE" property means in practice that relations such as
\begin{align*}
	&\Phi^{E,(N_2,M_1),(0,0)}_{\infty_E,\infty_S}\left(
		\begin{tikzpicture}[guillpart,yscale=1.8,xscale=2.4]
			\draw[guillsep] (0,0)--(0,1)--(4,1) (1,0)--(1,1) (2,0)--(2,1) (3,0)--(3,1);
			\draw[guillsep,dotted] (0,0)--(4,0)--(4,1);
			\node at (0.5,0.5) {$B_S(x_1)$};
			\node at (1.5,0.5) {$\ldots$};
			\node at (2.5,0.5) {$B_S(x_k)$};
			\node at (3.5,0.5) {$V_{SE}$};
		\end{tikzpicture}_{\ca{B}^{(R')}}
	\right)
\\
	&= \begin{tikzpicture}[guillpart,xscale=6.2,yscale=1.8]
		\draw[guillsep] (0,0)--(0,1)--(4,1) (1,0)--(1,1) (2,0)--(2,1) (3,0)--(3,1);
		\draw[guillsep,dotted] (0,0)--(4,0)--(4,1);
		\node at (0.5,0.5) {$\Phi_{1,\infty_S}^{N_2,0}(B_S(x_1))$};
		\node at (1.5,0.5) {$\ldots$};
		\node at (2.5,0.5) {$\Phi_{1,\infty_S}^{N_2,0}(B_S(x_1))$};
		\node at (3.5,0.5) {$\Phi^{(N_2,M_1),(0,0)}_{\infty_E,\infty_S}(V_{SE})$};
	\end{tikzpicture}_{\ca{B}}
\\
	&= \Lambda^{N_2(M_1+k)}
	\begin{tikzpicture}[guillpart,yscale=1.8,xscale=2.5]
		\draw[guillsep] (0,0)--(0,1)--(4,1) (1,0)--(1,1) (2,0)--(2,1) (3,0)--(3,1);
		\draw[guillsep,dotted] (0,0)--(4,0)--(4,1);
		\node at (0.5,0.5) {$A_S(x_1)$};
		\node at (1.5,0.5) {$\ldots$};
		\node at (2.5,0.5) {$A_S(x_k)$};
		\node at (3.5,0.5) {$U_{SE}$};
	\end{tikzpicture}_{\ca{B}}
\end{align*}
hold on any corner with an arbitrary number of actions of side elements $B_a(u)$.
\begin{proof}
	The first two points are a simple exercise in linear algebra using the fact that $a_l$, $a_r$, $b_l$ and $b_r$ are eigenvectors of $\MarkovWeight{A}$ and $\MarkovWeight{B}$.
	
	The last point is obtained either by the "morphism of ROPE" property or by direct computation. The exponent $PQ-P'Q'$ is obtained by summing the exponents of $\Lambda$ for each element in the ROPErep of $g_{R'}$: there are $P'$ (resp. $P'$, $Q'$ and $Q'$) elements $B_S(x)$ (resp. $B_N(y)$, $B_W(w)$ and $B_E(z)$), each with an exponent of $\Lambda$ equal to $M_1$ (resp. $M_2$, $N_1$ and $N_2$) when mapped to $A_S(x)$ (resp. $A_N(y)$, $A_W(w)$ and $A_E(z)$). The corner elements produce power of $\Lambda$ equal to $N_1 M_1$, $N_1 M_2$, $N_2 M_1$ and $N_2 M_2$. One observes that the exponent of $\Lambda$ corresponds to the number of weight faces used to produce the $B_a(u)$ out of the $A_a(u)$ in Theorem~\ref{theo:stability}. One then has:
	\[
		P' (M_1+M_2) + Q'(N_1+N_2) + N_1 M_1+ N_1 M_2+ N_2 M_1 N_2 M_2 = PQ-P'Q'
	\]
	by using $P=N_1+P'+N_2$ and $Q=M_1+Q'+M_2$.
\end{proof}

\subsubsection{Remarks}  One observes on \eqref{eq:horizvertmodel:eigenweight} that the weights~\eqref{eq:horizvertmodel:boundaryweight} behave as "eigen"-weights under the restriction to smaller rectangles, i.e. they acquire a \emph{global} factor given by $\Lambda^S$ where $S$ is the surface that is removed when shrinking rectangles, in the same way as, in dimension one, Perron-Frobenius eigenvectors acquire a factor $\Lambda^L$ where $L$ is the length removed when shrinking a segment to a smaller one. This is exactly what is expected from a generalization of "eigen"-boundary conditions.

It is interesting to notice that the factors $\Lambda^S$ where $S$ is a surface are already present at the \emph{local} level of the operators of the ROPErep of $G_{R}$, $G_{R'}$ and $G_{R\to R'}$. The previous computations on such a simple model already exhibit the key feature of the present section and definitions~\ref{def:eigenalgebrauptomorphims} and \ref{def:cornereigensemigroups} below: obtaining back the operators $A_a(u)$ and their factor $\Lambda^S$ from the $B_a(u)$ requires a composition by new morphisms $\Phi_\bullet^\bullet$.

\subsection{The second example: independent oblique lines.}

The previous model of horizontal and vertical independent processes on lines was very simple since the lines are parallel to the directions of the cuts. We now present a second model made of independent lines with other directions so that the algebraic is more involved.

We consider again finite sets $S_1$ and $S_2$ and two rectangular matrices $\MarkovWeight{C}$ and $\MarkovWeight{D}$ such that:
\begin{itemize}
\item the size of $\MarkovWeight{C}$ (resp. $\MarkovWeight{D}$) is $|S_2| \times |S_1|$ (resp. $|S_1|\times |S_2|$);
\item both matrices are indexed by elements of $S_1$ and $S_2$;
\item both matrices have positive coefficients and the products $\MarkovWeight{C}\MarkovWeight{D}$ and $\MarkovWeight{D}\MarkovWeight{C}$ are irreducible.
\end{itemize} 
We now build the face weight $\MarkovWeight{F}$ out of these matrices using 
\begin{equation}\label{eq:trivialmodelDIAG}
\MarkovWeight{F}(x,y,w,z) = \MarkovWeight{C}_{x,w}\MarkovWeight{D}_{z,y}
\end{equation}
for any $(x,y,w,z)\in S_1^2\times S_2^2$ and consider, on any rectangle $R$, the two-dimensional Markov process $(X_e)_{e\in\Edges{R}}$ with all face weights equal to $\MarkovWeight{F}$.

The structure of $\MarkovWeight{W}$ builds couplings only along oblique lines (see figure~\ref{fig:trivialmodels}). More precisely, we define an \emph{oblique line} $\ca{L}$ as an infinite sequence of edges $(e_n)_{n\in\setZ}$ such that:
\begin{itemize}
\item for any even (resp. odd) $n$, $e_n$ is an horizontal (resp. vertical) edge;
\item for any even $n$, there exists a face $f$ such that $e_n$ is the South boundary of $f$ and $e_{n+1}$ is the West boundary of $f$;
\item for any odd $n$, there exists a face $f$ such that $e_n$ is the East boundary of $f$ and $e_{n+1}$ is the North boundary of $f$.
\end{itemize} 

\begin{figure}
\begin{center}
Basic tiles:

\begin{tikzpicture}
\begin{scope}[scale=0.4,xshift=-4cm]
	\draw (0,0) rectangle (4,4);
	\draw[ultra thick,red] (2,0)--(2,4);
	\draw[ultra thick,orange] (0,2)--(1.9,2) (2.1,2)--(4,2);
	\node at (2,0) [anchor = north] {$x$};
	\node at (2,4) [anchor = south] {$y$};
	\node at (0,2) [anchor = east] {$w$};
	\node at (4,2) [anchor = west] {$z$};
	\node at (1,2) [anchor = south,orange] {$\MarkovWeight{B}_{w,z}$};
	\node at (2,1) [anchor = west,red] {$\MarkovWeight{A}_{x,y}$};
\end{scope}
\begin{scope}[scale=0.4,xshift=4cm]
	\draw (0,0) rectangle (4,4);
	\draw[ultra thick,blue] (2,0)--(0,2);
	\draw[ultra thick,green] (4,2)--(2,4);
	\node at (2,0) [anchor = north] {$x$};
	\node at (2,4) [anchor = south] {$y$};
	\node at (0,2) [anchor = east] {$w$};
	\node at (4,2) [anchor = west] {$z$};
	\node at (1,1) [anchor = west,blue] {$\MarkovWeight{C}_{x,w}$};
	\node at (3,3) [anchor = north,green] {$\MarkovWeight{D}_{z,y}$};
\end{scope}
\end{tikzpicture}

Tiles glued to form larger rectangles:

\begin{tikzpicture}
\begin{scope}[xshift=-4cm,scale=0.5]
\foreach \x in {0,1,...,8}
	\foreach \y in {0,1,...,6} {
		\begin{scope}[xshift=\x cm,yshift=\y cm]		
			\draw (0,0) rectangle (1,1);
			\draw[ultra thick,red] (0.5,0)--(0.5,1);
			\draw[ultra thick,orange] (0,0.5)--(0.4,0.5) (0.6,0.5)--(1,0.5);
		\end{scope}
		}
\end{scope}
\begin{scope}[xshift=4cm,scale=0.5]
\foreach \x in {0,1,...,8}
	\foreach \y in {0,1,...,6} {
		\begin{scope}[xshift=\x cm,yshift=\y cm]		
			\draw (0,0) rectangle (1,1);
			\draw[ultra thick,blue] (0.5,0)--(0,0.5);
			\draw[ultra thick,green] (1,0.5)--(0.5,1);
		\end{scope}
		}
\end{scope}
\end{tikzpicture}
\end{center}

\caption{\label{fig:trivialmodels}Pictures illustrating the two trivial models of Section~\ref{sec:trivialfactorizedcase}. On the top pictures, the structure of the face weights is illustrated: a solid internal line indicates a non-zero coupling between two edges. On the bottom, couplings are illustrated on a larger rectangle: processes on distinct lines are independent. On the left, horizontal lines are monochrome, which corresponds to powers of the matrices $\MarkovWeight{A}$ and $\MarkovWeight{B}$; on the right, each line has two alternating colours corresponding to a product of alternating $\MarkovWeight{C}$ and $\MarkovWeight{D}$ matrices.}
\end{figure}

Then, for any oblique line $\ca{L}$ intersecting a rectangle $R$, we define $X^\ca{L}$ as the restriction of the process $X$ on the intersection of $\ca{L}$ with the edges inside $R$. From a probabilistic point of view, the collection of processes $(X^\ca{L}_e)_{e\in\ca{L}}$ indexed by the oblique lines that intersect $R$ forms a sequence of \emph{independent} processes such that, for any $\ca{L}$, $X^{\ca{L}}$ is a one-dimensional Markov process with values alternatively in $S_1$ and $S_2$ and with alternating weights $\MarkovWeight{C}$ and $\MarkovWeight{D}$. 

Compared to the previous situation, we have again independent processes living on $P+Q$ lines, where the rectangle $R$ have shape $(P,Q)$: however, instead of having $Q$ horizontal and $P$ vertical lines with intersections on each face, there are $P+Q$ parallel oblique lines.  The partition function for $Q\leq P$ is again easy to compute as a product of one-dimensional partition functions:
\begin{equation}\label{eq:trivialoblique:Z}
Z_R(\MarkovWeight{F}; x,y,w,z) = 
\left(\prod_{k=1}^Q  (\MarkovWeight{C}(\MarkovWeight{D}\MarkovWeight{C})^{k-1})_{x_k,w_k}\right)
\left(\prod_{k=Q+1}^P  ((\MarkovWeight{C}\MarkovWeight{D})^{Q})_{x_k,y_{k-Q}}\right)
\left(\prod_{k=1}^Q  (\MarkovWeight{D}(\MarkovWeight{C}\MarkovWeight{D})^{Q-k})_{z_k,y_{P-Q+k}}\right)
\end{equation}
The complementary case $P<Q$ has a similar structure.
Asymptotics of $Z_R$ for large rectangles as well as invariant boundary conditions are obtained by introducing the Perron-Frobenius eigenvectors of $\MarkovWeight{C}\MarkovWeight{D}$ and $\MarkovWeight{D}\MarkovWeight{C}$, which share the same largest eigenvalue $\Lambda$:
\begin{align*}
\MarkovWeight{CD}v^{(1)}_R &= \Lambda v^{(1)}_R
&
v^{(1)}_L \MarkovWeight{CD} &=  v^{(1)}_L \Lambda
\\
\MarkovWeight{DC}v^{(2)}_R &= \Lambda v^{(2)}_R
&
v^{(2)}_L \MarkovWeight{DC} &=  v^{(2)}_L \Lambda
\\
\scal{v^{(1)}_L}{v^{(1)}_R} &= 1
&
\scal{v^{(2)}_L}{v^{(2)}_R} &= 1
\end{align*}
and one has necessarily the relations:
\begin{align*}
\MarkovWeight{D}v_R^{(1)} &= s_1 v_R^{(2)} 
&
\MarkovWeight{C}v_R^{(2)} &= s_2 v_R^{(1)} 
\\
v_L^{(1)} \MarkovWeight{C} &= s_2 v_L^{(2)} 
&
v_L^{(2)} \MarkovWeight{D} &= s_1 v_L^{(1)} 
\end{align*}
with $s_1s_2 =\Lambda$.
For any fixed $0<\alpha<1$ and sequence $Q_P=\lfloor \alpha P \rfloor$ and any sequence of boundary conditions, an easy computation (Cesar\`o mean) proves 
\[
\frac{1}{P Q_P} \log Z_{[P,Q_P]}(\MarkovWeight{F}; x^{(P)},y^{(P)},w^{(P)},z^{(P)}) \xrightarrow{P\to\infty} \log \Lambda 
\]
This free energy density $f=\log \Lambda$ can be exactly achieved for any rectangle size $(P,Q)$ by the factorized boundary conditions:
\begin{equation}
G_R(x,y,w,z) = \prod_{i=1}^{P} v_L^{(1)}(x_i)\prod_{i=1}^{P} v_R^{(1)}(y_i)\prod_{j=1}^Q v_R^{(2)}(w_j)\prod_{j=1}^Q v_L^{(2)}(z_j)
\end{equation}
for which one has the exact partition function
\begin{equation}
Z^{\boundaryweights}_R(\MarkovWeight{F};G_R) = \Lambda^{PQ}.
\end{equation}

It is interesting to recast these trivial calculations in the language of ROPE representation and of Theorem~\ref{theo:stability} in order to illustrate the underlying necessity of morphisms of $\Guill$-algebras on the boundary operators.

The initial ROPE is again given $\ca{B}_{p,q}=\setR$ for any colour $(p,q)$ and the ROPEs $\ca{B}'_{\PatternShapes(\patterntype{fp}^*)}$ and $\ca{B}^{(R')}_{\PatternShapes(\patterntype{fp}^*)}$ required by Theorem~\ref{theo:stability} and Corollary~\ref{coro:stab:smallerROPE} are the same as for the previous model. The ROPErep of $g_R$ is now given by the corner elements $U_{ab}=1$ and the side operators $A_S(x) = v_L^{(1)}(x)$, $A_N(y)=v_R^{(1)}(y)$, $A_W(w)=v_R^{(2)}(w)$ and $A_E(z)=v_L^{(2)}(z)$. 

We now present detailed computations only on the South case and let the reader adapt them to the three other sides. After the use of Corollary~\ref{coro:stab:smallerROPE} , the operators $B_S(x)$ may be seen as elements of $\End(V(S_2))^{\otimes M_1}$ and we introduce the canonical indices such that
$B_S(x)_{i_1,i_2,\ldots,i_{M_1}; j_1,\ldots,j_{M_1}}$ where $(i_k,j_k)$ corresponds to the $k$-th factor of the tensor product. To make equations shorter, we write $\mathbf{i} = (i_1,\ldots,i_{M_1})$ and $\mathbf{j}=(j_1,\ldots,j_{M_1})$.  Here, a direct evaluation of $B_S(x)$ from \eqref{eq:def:boundarySouthobjectrestriction} gives:
\begin{equation}
B_S(x)_{\mathbf{i},\mathbf{j}} = \MarkovWeight{D}_{j_1,x} \left( \prod_{k=1}^{M_1-1} (\MarkovWeight{DC})_{j_{k+1},i_k} \right) (v_L^{(1)}\MarkovWeight{C})_{i_{M_1}}
\end{equation}
The one-dimensional structure of these operators is less obvious than in the previous case of crossing horizontal and vertical lines. Indeed, one has:
\begin{equation*}
\left(B_S(x)B_S(y)\right)_{\mathbf{i},\mathbf{j}} = \MarkovWeight{D}_{j_1,y} (\MarkovWeight{DCD})_{j_2,x} \left( \prod_{k=1}^{M_1-2} (\MarkovWeight{CD})_{j_{k+2},i_k} \right) ( v_L^{(1)}\MarkovWeight{CDC})_{i_{M_1-1}}(v_L^{(1)}\MarkovWeight{C})_{i_{M_1}}
\end{equation*}
and similar formulae for larger products. A direct computation gives:
\begin{align*}
	B_S(x) \left( \bigotimes_{k=1}^{M_1} v_R^{(2)} \right) &= \Lambda^{M_1} \scal{v_L^{(1)}}{e_x} \bigotimes_{k=1}^{M_1} v_R^{(2)}
	\\
	\left(\bigotimes_{k=1}^{M_1} v_L^{(2)} \right) B_S(x)  &= s_2\Lambda^{M_1-1} (\MarkovWeight{D}e_x) \otimes \bigotimes_{k=2}^{M_1} v_L^{(2)}  
\end{align*}
which shows that all the $B_S(x)$ have a triangular structure
\[
\begin{pmatrix}
    \Lambda^{M_1}v_L^{(1)}(x) & \diamond_2(x)
	\\
	0 & \diamond_1(x)
\end{pmatrix}
\] on a direct decomposition $V(S_2)^{\otimes M_1} = \setK {v_L^{(2)}}^{\otimes M_1}\oplus W $ with non-trivial blocks $\diamond_i(x)$. The same linear maps as \eqref{eq:horizvertmodel:morphisms} with projectors on the eigenspace of $\MarkovWeight{DC}$ map $B_S(x)$ to $\Lambda^{M_1}A_S(x)$ and are morphisms of algebra because of the triangular block structures. On the three other sides, the morphisms of algebra are defined in the same way.

On the South-East corner, we have trivially \[
V_{SE} = \Lambda^{N_2M_1}(v_L^{(1)})^{\otimes N_2}\otimes (v_L^{(2)})^{\otimes M_1}
\]
from the action of the matrices $\MarkovWeight{C}$ and $\MarkovWeight{D}$ on the vectors that corresponds to $A_S(x)$ and $A_E(z)$. The right part corresponds to the eigenvector of the $B_S(x)$ and thus the action of $B_S(x)$ is one-dimensional. The same is true for the action of $B_E(z)$ on the left part of the tensor product. The one-dimensional space can then trivially mapped to $\setK$ through a canonical morphism $\Phi^{(N_2,M_1),(0,0)}_{\infty_E,\infty_S}$ as in the previous model of horizontal and vertical lines.

The South-West corner is a bit more involved. The expression of $V_{SW}$ is given by
\[
V_{SW} = \Lambda^{N_1(M_1-N_1)} (v_L^{(1)})^{\otimes N_1-M_1}\otimes
\sum_{
\substack{
x_1,\ldots,x_{M_1} \in S_1 \\
w_1,\ldots,w_{M_1} \in S_2 }}
\prod_{i=1}^{M_1} (\MarkovWeight{D}(\MarkovWeight{C}\MarkovWeight{D})^{i-1})_{w_i,x_i}
 \left(\otimes_{i=1}^{M_1} e^*_{x_{M_1+1-i}}\right) \otimes 
  \left(\otimes_{i=1}^{M_1} e^*_{w_i}\right)
\]
Each factor $\MarkovWeight{D}(\MarkovWeight{C}\MarkovWeight{D})^{i-1}$ corresponds to an oblique line that connects both sides and goes through $2i-1$ faces. Contrary to the South-East corner, it is not a tensor product of the Perron-Frobenius eigenvectors $B_S(x)$ and $B_W(w)$. 

However, it has a block structure compatible with the previously described block structure of $B_S(x)$ and $B_W(w)$. We now define the following morphism
\begin{align*}
	\Phi_{\infty_W,\infty_S}^{(N_1,M_1),(0,0)} : \ca{B}^{(R')}_{\infty_W,\infty_S} & \to \ca{B}_{\infty_W,\infty_S}
	\\
	\bigotimes_{i=1}^{N_1} u_i \otimes \bigotimes_{j=1}^{M_1} v_i 
	& \mapsto 
	\prod_{i=1}^{N_1} \scal{u_i}{v_R^{(1)} } \prod_{j=1}^{M_1} \scal{v_L^{(2)}}{v_i}
\end{align*}
It is then an exercise to check the following property of morphism of ROPE:
\begin{equation}\begin{split}
	&\Phi_{\infty_W,\infty_S}^{SW,(N_1,M_1),(0,0)}\left(
	\begin{tikzpicture}[guillpart,yscale=1.3,xscale=3.]
		%\fill[guillfill] (0,0) rectangle (4,4);
		\draw[guillsep, dotted] (0,4)--(0,0)--(4,0);
		\draw[guillsep] (4,0)--(4,1)--(1,1)--(1,4)--(0,4);
		\draw[guillsep] (1,0)--(1,1) (2,0)--(2,1) (3,0)--(3,1) (0,1)--(1,1) (0,2)--(1,2) (0,3)--(1,3);
		\node at (0.5,0.5) {$V_{SW}$};
		\node at (1.5,0.5) {$B_S(x_1)$};
		\node at (2.5,0.5) {$\ldots$};
		\node at (3.5,0.5) {$B_S(x_k)$};
		\node at (0.5,1.5) {$B_W(w_1)$};
		\node at (0.5,2.5) {$\ldots$};
		\node at (0.5,3.5) {$B_W(w_l)$};
	\end{tikzpicture}
	\right)
	\\
	&= \begin{tikzpicture}[guillpart,yscale=1.6,xscale=6.5]
		%\fill[guillfill] (0,0) rectangle (4,4);
		\draw[guillsep, dotted] (0,4)--(0,0)--(4,0);
		\draw[guillsep] (4,0)--(4,1)--(1,1)--(1,4)--(0,4);
		\draw[guillsep] (1,0)--(1,1) (2,0)--(2,1) (3,0)--(3,1) (0,1)--(1,1) (0,2)--(1,2) (0,3)--(1,3);
		\node at (0.5,0.5) {$\Phi^{(N_1,M_1),(0,0)}_{\infty_W,\infty_S}(V_{SW})$};
		\node at (1.5,0.5) {$\Phi_{1,\infty_S}^{N_1,0}(A_S(x_1))$};
		\node at (2.5,0.5) {$\ldots$};
		\node at (3.5,0.5) {$\Phi_{1,\infty_S}^{N_1,0}(A_S(x_k))$};
		\node at (0.5,1.5) {$\Phi_{\infty_W,1}^{M_1,0}(A_W(w_1))$};
		\node at (0.5,2.5) {$\ldots$};
		\node at (0.5,3.5) {$\Phi_{\infty_W,1}^{M_1,0}(A_W(w_l))$};
	\end{tikzpicture}
\end{split}
\end{equation}
and moreover, one has $\Phi_{\infty_W,\infty_S}^{(N_1,M_1),(0,0)}(V_{SW})=\Lambda^{N_1M_1} U_{SW}$ as for the previous model.

As for the previous model, the local operators $A_a(x)$ satisfy an "eigen"-property up to morphisms: the operator $B_a(x)$ obtained by gluing $S$ faces can be mapped back to $A_a(x)$ up to a factor $\Lambda^{S}$ by a morphism of ROPE.

\section{Generalized eigen-generators on South half-strips}

\subsection{A guiding principle and first definitions.}

\subsubsection{Associativities and inductive sequences of algebras.}

We consider a fixed 2D-semigroup $(\MarkovWeight{F}_{p,q})_{(p,q)\in\PatternShapes(\patterntype{r})}$ and we define the following linear maps, for any $p,r\in\setL^*$ and $q\in\setL^*\cup\{\infty_S\}$:
\begin{equation}\label{eq:def:actionmorphism:psi}
	\begin{split}
		\psi_p^{\MarkovWeight{F},S,q,r} : \ca{A}_{p,q} & \to \ca{A}_{p,q+r}  \\
		A & \mapsto \begin{tikzpicture}[guillpart,xscale=2]
			\fill[guillfill] (0,0) rectangle (1,2);
			\draw[guillsep] (0,0)--(0,2) -- (1,2) -- (1,0)
			(0,1)--(1,1);
			\node at (0.5,0.5) {$A$};
			\node at (0.5,1.5) {$\MarkovWeight{F}_{p,r}$};
		\end{tikzpicture}
	\end{split}
\end{equation}
with the same convention $\infty_S+r=\infty_S$ as before.

\begin{lemm}\label{lemm:morphismsfrom2Dsemigroup}
	For any $q\in\setL^*\cup\{\infty_S\}$ and $r\in\setL^*$, the maps $\psi_\bullet^{\MarkovWeight{F},S,q,r}$ define a morphism of $\Guill_1$-algebras from $\ca{A}_{\bullet,q}$ to $\ca{A}_{\bullet,q+r}$. Moreover, they satisfy the composition rule:
	\begin{equation}\label{eq:2Dsemigroup:directedmorph}
		\psi_p^{\MarkovWeight{F},S,q+r,r'} \circ \psi_p^{\MarkovWeight{F},S,q,r} = \psi_p^{\MarkovWeight{F},S,q,r+r'}
	\end{equation}
\end{lemm}
\begin{proof}
	A morphism of $\Guill_1$-algebra $\psi_\bullet^{\MarkovWeight{F},S,q,r}$ corresponds to a set of maps $\psi_{p}^{\MarkovWeight{F},S,q,r}: \ca{A}_{p,q}\to\ca{A}_{p,q+r}$ such that:
	\[
	\psi_{p_1+p_2}^{\MarkovWeight{F},S,q,r}\left(
	\begin{tikzpicture}[guillpart,scale=1.5]
		\fill[guillfill] (0,0) rectangle (2,1);
		\draw[guillsep] (0,0)--(0,1)--(2,1)--(2,0) (1,0)--(1,1);
		\draw[guillsep,dashed] (0,0)--(2,0);
		\node at (0.5,0.5) {$s_1$};
		\node at (1.5,0.5) {$s_2$};
	\end{tikzpicture}
	\right)
	= 
	\begin{tikzpicture}[guillpart,yscale=1.5,xscale=4.]
		\fill[guillfill] (0,0) rectangle (2,1);
		\draw[guillsep] (0,0)--(0,1)--(2,1)--(2,0) (1,0)--(1,1);
		\draw[guillsep,dashed] (0,0)--(2,0);
		\node at (0.5,0.5) {$\psi_{p_1}^{\MarkovWeight{F},S,q,r}(s_1)$};
		\node at (1.5,0.5) {$\psi_{p_2}^{\MarkovWeight{F},S,q,r}(s_2)$};
	\end{tikzpicture}
	\]
	for any $p_1,p_2\in\setL^*$ and $s_1\in\ca{A}_{p_1,q}$ and $s_2\in\ca{A}_{p_2,q}$ (the dashed line indicates a line or not depending on whether $q$ is finite or not). In the present case, it is fulfilled from the square associativity~\eqref{eq:guill2:interchangeassoc} of the guillotine products and \eqref{eq:horizsemigroup} for the 2D-semi-group $\MarkovWeight{F}_{\bullet}$.
	
	The composition identity~\eqref{eq:2Dsemigroup:directedmorph} is fulfilled using the vertical associativity~\eqref{eq:guill2:vertassoc} of the guillotine products and \eqref{eq:vertsemigroup} for the 2D-semi-group $\MarkovWeight{F}_\bullet$.
\end{proof}

The interest of this proof relies on the symmetry-breaking between the three associativities \eqref{eq:guill2:listassoc} when one considers an action on one side: the horizontal associativity defines the algebraic $\Guill_1$-structure of the space $\ca{A}_{\bullet,q}$ on which $\MarkovWeight{F}_\bullet$ acts, the square associativity define morphisms of $\Guill_1$-structures, the vertical associativity provides an action that is similar to the one-dimensional situation of matrices $A^{q'-q}$ acting on a vector $v$ through the composition rule~\eqref{eq:2Dsemigroup:directedmorph}.

We have similar constructions on the three other sides North, West and East using the action of $\MarkovWeight{F}_{\bullet}$ on the three other boundary spaces . One can also combine opposite sides and show that, for $(q,q')\neq(\infty_S,\infty_N)$, the morphism
\[
\begin{split}
	\psi_{p}^{\MarkovWeight{F},S,N,q,r,q'} : \ca{A}_{p,q}\otimes \ca{A}_{p,q'} \mapsto \ca{A}_{p,q+r+q'}
	\\	
	A\otimes B \mapsto \begin{tikzpicture}[guillpart,xscale=3]
		\fill[guillfill] (0,0) rectangle (1,3);
		\draw[guillsep] (0,0)--(0,3) (1,0)-- (1,3)  (0,1)--(1,1) (0,2)--(1,2);
		\draw[guillsep,dashed] (0,0)--(1,0) (0,3)--(1,3);;
		\node at (0.5,0.5) {$A$};
		\node at (0.5,1.5) {$\MarkovWeight{F}_{p,r}$};
		\node at (0.5,2.5) {$B$};
	\end{tikzpicture}
\end{split}
\]
satisfy, for any partition $r=r_1+r_2$,
\[
\psi_{p}^{\MarkovWeight{F},S,N,q,r,q'} = \psi_{p}^{\MarkovWeight{F},S,q,r_1}\circ \psi_{p}^{\MarkovWeight{F},N,q',r_2}
= \psi_{p}^{\MarkovWeight{F},N,q',r_2} \circ \psi_{p}^{\MarkovWeight{F},S,q,r_1}
\]
and is thus again a morphism of $\Guill_1$-algebra. In the case $(q,q')=(\infty_S,\infty_N)$, a base point $x$ have to be introduced and we have a collection, indexed by the base point, of morphisms $\ca{A}_{p,\infty_S}\otimes\ca{A}_{p,\infty_N}\to\ca{A}_{p,\infty_{SN}}$.

\begin{prop}\label{prop:towerofguillalg}
	Any sub-$\Guill_1$-algebra $\ca{S}_\bullet$ of $\ca{A}_{\bullet,q}$ for $q\in\setL^*\cup\{\infty_S\}$ defines a 
	directed structure of $\Guill_1$-algebras $\ca{S}_\bullet^{S,\MarkovWeight{F},r}=\psi_\bullet^{\MarkovWeight{F},S,q,r}(\ca{S}_\bullet)$ with $r\in \setL^*$ whose morphisms are given, for $r\leq r'$, by $\psi_\bullet^{\MarkovWeight{F},S,q+r,r'-r}$.
\end{prop}

This proposition illustrates the effect of the change of dimension from one to two by promoting successive images of a vector space $\ca{S}$ to successive images of an algebra using the transverse direction and the square associativity~\eqref{eq:guill2:interchangeassoc}: in dimension one, $\ca{S}_\bullet$ is just a vector space $\ca{S}$ of a vector space $\ca{A}_{\infty_L}=\ca{V}^*$ and the action of a weight $\MarkovWeight{A}$ is given by $v\mapsto v\MarkovWeight{A}$. There is no product on $\ca{V}^*$ and hence no notion of morphism. The vector spaces $\ca{S}^{L,\MarkovWeight{A},r}$ are the images of $\ca{S}$ by $\MarkovWeight{A}^r$ (acting on the left). In dimension two, even in the simple case of the canonical structure $\ca{T}_{p,q}$ for a finite state space Markov process, the presence of a transverse direction induces that the spaces $\ca{S}^{S,\MarkovWeight{F},r}_p$ necessarily belong to different subspaces of the graded space $\ca{T}_{p,\infty_S}$. However, the product in the transverse direction allows for the presence of non-trivial $\Guill_1$-morphisms between the various algebras $\ca{S}^{S,\MarkovWeight{F},r}_\bullet$.

Now, most constructions below will be inspired from the following principle: instead of imposing a direct algebraic relation between the elements of the various $\ca{S}_\bullet^{S,\MarkovWeight{F},r}$ as in dimension one, we may impose only these algebraic relation only up to $\Guill_1$-morphisms that map the $\ca{S}_\bullet^{S,\MarkovWeight{F},r}$ to a common space or to each other. We will see below two first examples with an eigenvalue-like algebraic relation and a commutation-like algebraic relation.

\subsubsection{Scaling maps}

In order to define eigen-elements, a notion of multiplication by a scalar that is compatible with the set of colours is required.

\begin{defi}[geometric sequence]
	A $\setK$-valued sequence $(u_p)_{p\in\setL^*}$ is called geometric whenever, for all $p_1,p_2\in\setL^*$, $u_{p_1+p_2}=u_{p_1}u_{p_2}$.
\end{defi}

In the discrete space case $\setL=\setN$, it coincides with the traditional notion of geometric sequence and such a sequence is characterized by the value $u_1$, so that, $u_p=u_1^p$. In the continuous space case $\setL=\setR_+$, this is not the case since it would require non-integer powers, which are not allowed for any $\setK$. In the interesting case $\setK=\setR_+$ for Markov processes, such a geometric sequence corresponds to $u_p = \exp(p f)$ for some $f\in\setR$.

\begin{prop}
Let $\ca{A}_{\PatternShapes(\patterntype{hs}_S)}$ be a $\Guill_2^{(\patterntype{hs}_S)}$-algebra. Let $\lambda=(\lambda_p)_{p\in\setL^*}$ be a geometric sequence.
The scaling maps $\mu_\bullet^{(u)}$ defined by:
\begin{align*}
	\mu_p^{(\lambda)}: \ca{A}_{p,\infty_S} &\to \ca{A}_{p,\infty_S}	\\
	a & \mapsto u_p a
\end{align*}
are morphisms of $\Guill_1$-algebras and commute with the action of $\ca{A}_{\PatternShapes(\patterntype{r})}$.
\end{prop}
\begin{proof}
The scaling maps are morphisms of $\Guill_1$-algebra since we have for any $p_1,p_2\in\setL^*$ and any $a_1\in\ca{A}_{p_1,\infty_S}$, $a_2\in\ca{A}_{p_2,\infty_S}$ since we have by linearity:
\begin{align*}
	\mu^{(\lambda)}_{p_1+p_2}\left(
	\begin{tikzpicture}[guillpart,scale=1.3]
		\fill[guillfill] (0,0) rectangle (2,1);
		\draw[guillsep] (0,0)--(0,1)--(2,1)--(2,0) (1,0)--(1,1);
		\node at (0.5,0.5) {$a_1$};
		\node at (1.5,0.5) {$a_2$};
	\end{tikzpicture}
	\right)
	&= \lambda_{p_1+p_2}
	\begin{tikzpicture}[guillpart,scale=1.3]
		\fill[guillfill] (0,0) rectangle (2,1);
		\draw[guillsep] (0,0)--(0,1)--(2,1)--(2,0) (1,0)--(1,1);
		\node at (0.5,0.5) {$a_1$};
		\node at (1.5,0.5) {$a_2$};
	\end{tikzpicture}
	=\lambda_{p_1}\lambda_{p_2}
	\begin{tikzpicture}[guillpart,xscale=1.5,yscale=1.3]
		\fill[guillfill] (0,0) rectangle (2,1);
		\draw[guillsep] (0,0)--(0,1)--(2,1)--(2,0) (1,0)--(1,1);
		\node at (0.5,0.5) {$a_1$};
		\node at (1.5,0.5) {$a_2$};
	\end{tikzpicture}
	\\
	&=\begin{tikzpicture}[guillpart,yscale=1.3,xscale=2.5]
		\fill[guillfill] (0,0) rectangle (2,1);
		\draw[guillsep] (0,0)--(0,1)--(2,1)--(2,0) (1,0)--(1,1);
		\node at (0.5,0.5) {$\lambda^{p_1}a_1$};
		\node at (1.5,0.5) {$\lambda^{p_2}a_2$};
	\end{tikzpicture}
	=\begin{tikzpicture}[guillpart,yscale=1.3,xscale=3.]
		\fill[guillfill] (0,0) rectangle (2,1);
		\draw[guillsep] (0,0)--(0,1)--(2,1)--(2,0) (1,0)--(1,1);
		\node at (0.5,0.5) {$\mu_{p_1}^{(\lambda)}(a_1)$};
		\node at (1.5,0.5) {$\mu_{p_2}^{(\lambda)}(a_2)$};
	\end{tikzpicture}
\end{align*}
using the geometric property of the sequence $\lambda$. We have moreover for any $p\in\setL^*$, any $q\in\setL^*$, any $a\in\ca{A}_{p,\infty_S}$ and any $\MarkovWeight{T}\in\ca{A}_{p,q}$:
\begin{align*}
	\mu_{p}^{(\lambda)}\left(
		\begin{tikzpicture}[guillpart,yscale=1.2,xscale=2]
			\fill[guillfill] (0,0) rectangle (1,2);
			\draw[guillsep] (0,0)--(0,2)--(1,2)--(1,0) (0,1)--(1,1);
			\node at (0.5,0.5) {$a$};
			\node at (0.5,1.5) {$\MarkovWeight{T}$};
		\end{tikzpicture}
	\right)
	&= \lambda_{p} 		\begin{tikzpicture}[guillpart,yscale=1.2,xscale=2]
		\fill[guillfill] (0,0) rectangle (1,2);
		\draw[guillsep] (0,0)--(0,2)--(1,2)--(1,0) (0,1)--(1,1);
		\node at (0.5,0.5) {$a$};
		\node at (0.5,1.5) {$\MarkovWeight{T}$};
	\end{tikzpicture}
	= \begin{tikzpicture}[guillpart,yscale=1.2,xscale=2]
		\fill[guillfill] (0,0) rectangle (1,2);
		\draw[guillsep] (0,0)--(0,2)--(1,2)--(1,0) (0,1)--(1,1);
		\node at (0.5,0.5) {$\lambda_p a$};
		\node at (0.5,1.5) {$\MarkovWeight{T}$};
	\end{tikzpicture}
	= \begin{tikzpicture}[guillpart,yscale=1.2,xscale=2.5]
		\fill[guillfill] (0,0) rectangle (1,2);
		\draw[guillsep] (0,0)--(0,2)--(1,2)--(1,0) (0,1)--(1,1);
		\node at (0.5,0.5) {$\mu^{(\lambda)}_p(a)$};
		\node at (0.5,1.5) {$\MarkovWeight{T}$};
	\end{tikzpicture}
\end{align*}
\end{proof}
We moreover have the composition property:
$\mu_{p}^{(\lambda\lambda')} = \mu_p^{(\lambda)}\circ\mu_p^{(\lambda')}$
where $\lambda\lambda'$ is the geometric sequence $(\lambda_p\lambda'_p)_{p\in\setL^*}$.

\subsubsection{A simplistic definition of eigen-generators up to morphisms of a \emph{single} morphism of \texorpdfstring{$\Guill_1$}{Guill1}-algebras.}

Proposition~\ref{prop:towerofguillalg} provides a whole collection of morphisms between various $\Guill_1$-algebras. We first concentrate on a single one of them and leave temporarily the dimension two.

\begin{defi}[system of eigen-generators up to morphisms for South algebras]\label{def:eigenalgebra:single}
	Let $(\ca{A}_{p})_{p\in\setL^*}$ and $(\ca{B}_{p})_{p\in\setL^*}$ be two $\Guill_1$-algebras. Let $\Psi_\bullet: \ca{A}_\bullet \to \ca{B}_\bullet$ be a morphism of $\Guill_1$-algebras, i.e. a collection of morphisms $\Psi_p : \ca{A}_p\to\ca{B}_p$ compatible with the products in the operad $\Guill_1$.
	
	A collection of sub-spaces $(\ca{V}_{p})_{p\in\setL^*}$ of $\ca{A}_\bullet$ with $\ca{V}_p\subset\ca{A}_p$ is a \emph{a system of eigen-generators up to morphisms} of $\Psi_\bullet$ with eigenvalue geometric sequence $(\lambda_{p})_{p\in\setL^*}$ if and only if there exist linear maps $\Phi_p : \ca{B}_{p}\to \ca{A}_p$, $p\in\setL^*$ such that, for all $N\in\setN^*$, for any sequence $(p_i)_{1\leq i\leq N}$ in $\setL^*$ with $P=\sum_i p_i$, for any sequence $(v_i)_{1\leq i\leq N}$ with $v_i\in \Psi_{p_i}(\ca{V}_{p_i})$, it holds:
		\begin{align}\label{eq:eigenalg:single:products}
			\Phi_{P}\left(
				m_{P;p_1,\ldots,p_N}(\Psi_{p_1}(v_1),\ldots,\Psi_{p_N}(v_N))
			\right)
			&= 
				m_{P;p_1,\ldots,p_N}(\Phi_{p_1}(\Psi_{p_1}(v_1)),\ldots,\Phi_{p_N}(\Psi(v_N)))
				\\
			&= \lambda_{P} m_{P;p_1,\ldots,p_N}(v1,\ldots,v_N)	
			\label{eq:eigenalg:single:eigenval}
		\end{align}
\end{defi}

At the level of $\ca{B}_\bullet$, the maps $\Phi_\bullet$ are only linear and not necessarily $\Guill_1$-morphisms: only their restrictions to the sub-$\Guill_1$-algebra generated by the image under $\Psi_\bullet$ of the $\ca{V}_p$ are $\Guill_1$-morphisms. The first equation of the definition introduce strong constraints on the possible maps $\Phi_p$, which are new with respect to the one-dimensional case for which there are no product.

The second item is very similar to the one-dimensional definition of an eigen-space. In particular, if the map $\Psi_\bullet$ corresponds to the one-dimension case of a map $A\in\End(\ca{W})$ acting on a vector space $\ca{W}$ through $w\mapsto Aw$, an eigenspace $\ca{V}$ for an eigenvalue $\lambda$ is the equalizer of the following diagram:
\begin{equation}\label{eq:eigenspaceegalizer}
	\begin{tikzpicture}[baseline={(current bounding box.center)}]
		\matrix (m) [matrix of math nodes,row sep=3em,column sep=1em,minimum width=1em]
		{
			\ca{V}
			\subset 
			\ca{W} 
			& & \ca{W}
			\\
		};
		\path[-stealth]	([yshift=0.5ex] m-1-1.east) edge node [above] { $A$ }  ([yshift=0.5ex] m-1-3.west);
		\path[-stealth] ([yshift=-0.5ex] m-1-1.east) edge node [below] { $\lambda$ } ([yshift=-0.5ex] m-1-3.west);
	\end{tikzpicture}
\end{equation}
which is the same as \eqref{eq:eigenalg:single:eigenval} without the colours $p\in\setL^*$ and without the map $\Phi_\bullet$. 

We have considered here a system of eigen-generators $\ca{V}_\bullet$ but, if one considers the $\Guill_1$-algebra $\ca{S}_\bullet$ generated by $\ca{V}_\bullet$, the first item of the previous definition can be modified to require that the $\Phi_\bullet$ restricted to the spaces $\ca{S}_\bullet$ are morphisms of $\Guill_1$-algebras. We will see below that it is often more practical to consider a system of generators rather than the whole generated algebra. In the present case, the definitions are equivalent and we may switch implicitly from one to the other.

\subsection{Definition of South eigen-elements up to morphisms}

We now consider the particular case $q=\infty_S$ in Lemma~\ref{lemm:morphismsfrom2Dsemigroup}, for which the maps $\psi_\bullet^{\MarkovWeight{F},S,\infty_S,r}$ are $\Guill_1$-endomorphisms of $\ca{A}_{\bullet,\infty_S}$. Given a sub-$\Guill_1$-algebra $\ca{S}_\bullet$ of $\ca{A}_{\bullet,\infty_S}$, all the spaces $\ca{S}_{\bullet}^{S,r}$ are subspaces of $\ca{A}_{\bullet,\infty_S}$ and it is natural to wish to compare them (up to morphisms if needed) through the multiplication by scalar values, in order to generalize the notion of eigen-element.

The morphisms $\psi_\bullet^{S,\MarkovWeight{F},\infty_S,r} : \ca{A}_{\bullet,\infty_S}\to\ca{A}_{\bullet,\infty_S}$ generalize the one-dimensional maps $v\mapsto vA^r$ and preserve the $\Guill_1$-products in the transverse direction and, following Definition~\ref{def:eigenalgebra:single}, we may look for a system of eigen-$\Guill_1$-generators $(\ca{V}_p)_{p\in\setL^*}$ \emph{common} to all these morphisms $\psi_\bullet^{S,\MarkovWeight{F},\infty_S,r}$ with eigenvalue sequences $(\lambda_{rp})_{p\in\setL^*}$ and up to morphisms $\Phi_\bullet^{(r)}$. This would correspond to West-East-product-preserving maps
\[
\Phi^{(r)}_p\left(
\begin{tikzpicture}[guillpart,xscale=1.5,yscale=1.2]
	\fill[guillfill] (0,0) rectangle (1,2);
	\draw[guillsep] (0,0)--(0,2)--(1,2)--(1,0) (0,1)--(1,1) ;
	\node at (0.5,0.5) {$v_p$};
	\node at (0.5,1.5) {$\MarkovWeight{F}_{p,r}$};
\end{tikzpicture}
\right)
=\lambda_{pr} v_p
\]
which corresponds to a generalization of \eqref{eq:1D:eigenval} and satisfy the requirement~\eqref{eq:2D:wishedreplacements:first}. This has however three main drawbacks: 
\begin{itemize}
	\item it forgets the module structure of $\ca{A}_{\bullet,\infty_S}$ on which $\ca{A}_{\PatternShapes(\patterntype{r})}$ acts from the North;
	\item it is not obvious whether composition rules have to be fixed between the various morphisms $\Phi^{(r)}_\bullet$;
	\item the previous morphisms do not allow for a simplification rule of expressions like~\eqref{eq:2D:wishedreplacements:second}.
\end{itemize}
These three drawbacks are related and the $\Guill_1$-structure of $\ca{A}_{\bullet,\infty_S}$ has to be enhanced to a South module-like structure, i.e. contains consistent actions from the North of the spaces $\ca{A}_{\PatternShapes(\patterntype{r})}$. The main requirement is the possibility of a simplification rule of \eqref{eq:2D:wishedreplacements:second} using suitable morphisms. We therefore use the following definition.

\begin{defi}[eigen-$\Guill_2^{(\patterntype{hs}_S)}$-generators up to morphisms of a 2D-semi-group]\label{def:eigenalgebrauptomorphims}
	Let $\ca{A}_{\PatternShapes(\patterntype{hs}_S)}$ be a $\Guill_2^{(\patterntype{hs}_S)}$-algebra and $\MarkovWeight{F}_{\bullet,\bullet}$ a 2D-semi-group in $\ca{A}_{\PatternShapes(\patterntype{r})}$. A collection of subspaces $(\ca{V}_p)_{p\in\setL^*}$ of $(\ca{A}_{p,\infty_S})_{p\in\setL^*}$ is \emph{a system of eigen-$\Guill_1$-generators up to morphisms} of $\MarkovWeight{F}_{\bullet,\bullet}$ with eigenvalue geometric sequence $(\lambda_p)_{p\in\setL^*}$ if and only if there exists a collection of linear maps $\Phi_p^{S,q,r} : \ca{A}_{p,\infty_S} \to \ca{A}_{p,\infty_S}$, $p\in\setL^*$, $q\in\setL^*$, $r\in\setL$ such that, for any $N\in\setN^*$, for any sizes $q,r\in\setL$, any sequence of sizes $(p_i)_{1\leq i\leq N}$ in $\setL^*$, for any sequence of South elements $(v_i)_{1\leq i\leq N}$ with $v_i\in\ca{V}_{p_i,\infty_S}$, for any element $\MarkovWeight{T}\in\ca{A}_{P,r}$, it holds, using  $P=\sum_{i=1}^N p_i$:
		\begin{equation}\label{eq:defeigenSouth:withM}
			\Phi_{P}^{S,q,r}\left(
			\begin{tikzpicture}[guillpart,yscale=1.3,xscale=2.]
				\fill[guillfill] (0,0) rectangle (3,3);
				\draw[guillsep] (0,0)--(0,3)--(3,3)--(3,0) (1,0)--(1,2) (2,0)--(2,2) (0,1)--(3,1) (0,2)--(3,2);
				\node at (0.5,0.5) {$v_1$};
				\node at (1.5,0.5) {$\dots$};
				\node at (2.5,0.5) {$v_N$};
				\node at (0.5,1.5) {$\MarkovWeight{F}_{p_1,q}$};
				\node at (1.5,1.5) {$\dots$};
				\node at (2.5,1.5) {$\MarkovWeight{F}_{p_N,q}$};
				\node at (1.5,2.5) {$\MarkovWeight{T}$};
			\end{tikzpicture}
			\right)
			= \lambda_{Pq} 	\begin{tikzpicture}[guillpart,xscale=1.5]
				\fill[guillfill] (0,0) rectangle (3,2);
				\draw[guillsep] (0,0)--(0,2)--(3,2)--(3,0) (1,0)--(1,1) (2,0)--(2,1) (0,1)--(3,1);
				\node at (0.5,0.5) {$v_1$};
				\node at (1.5,0.5) {$\dots$};
				\node at (2.5,0.5) {$v_N$};
				\node at (1.5,1.5) {$\MarkovWeight{T}$};
			\end{tikzpicture}
		\end{equation}
\end{defi}

The morphisms $\Phi_p^{S,q,r}$ are defined as linear maps on the whole space $\ca{A}_{p,\infty_S}$ but only their restrictions on the subspaces generated by the $\ca{V}_\bullet$, the 2D-semigroup and the action of the elements $T$ are relevant and their precise values on the other subspaces can be taken arbitrary at this point. 

This definition contains as many equations as equations~\eqref{eq:2D:wishedreplacements:first} (first item) and \eqref{eq:2D:wishedreplacements:second} (second item) where suitable colours $r\in\setL$ and $q\in\setL$ are added for $T$ and $\MarkovWeight{F}_{\bullet}$ with new morphisms and additional colours added for the transverse $\Guill_1$-structure.

Equations \eqref{eq:defeigenSouth:withM} for varying sizes $r$ and $q$ are not independent and the existence of some of the morphisms can be deduced from other ones. In particular, we have the following lemma.

\begin{lemm}\label{lemm:eigenalg:sufficientmorphisms}
	Let $\ca{V}_{p\in\setL^*}$ be a collection of subspaces of $(\ca{A}_{p,\infty_S})$ such that there exists $q_0\in\setL^*$ such that \eqref{eq:defeigenSouth:withM} is satisfied for some morphisms $\Phi_\bullet^{S,q_0,r}$ for all $r\in\setL^*$ but only $q_0\in\setL^*$. Then there exists morphisms $\Phi_\bullet^{S,q,r}$ for all $r\in\setL$ and $q=nq_0$ with $n\in\setN^*$ that satisfy \eqref{eq:defeigenSouth:withM}.
\end{lemm}
\begin{proof}
Let $q=n q_0$ with $n\in\setN^*$ and $r\in\setL^*$ be fixed. If $n=1$, the result is trivial. We now consider $n\geq 2$. Let $\MarkovWeight{T}\in \ca{A}_{P,r}$. We first use \eqref{eq:defeigenSouth:withM} with $\MarkovWeight{T}$ replaced by $\MarkovWeight{T}'=m_{SN}(\MarkovWeight{F}_{p,(n-1)q_0},\MarkovWeight{T})\in\ca{A}_{P,r+(n-1)q_0}$ and then use vertical associativities and the 2D-semi-group property of $\MarkovWeight{F}_\bullet$ to obtain
\begin{align*}
	\Phi_{P}^{S,q_0,r+(n-1)q_0}&\left(
	\begin{tikzpicture}[guillpart,yscale=1.25,xscale=2.3]
		\fill[guillfill] (0,0) rectangle (3,3);
		\draw[guillsep] (0,0)--(0,3)--(3,3)--(3,0) (1,0)--(1,2) (2,0)--(2,2) (0,1)--(3,1) (0,2)--(3,2);
		\node at (0.5,0.5) {$v_1$};
		\node at (1.5,0.5) {$\dots$};
		\node at (2.5,0.5) {$v_N$};
		\node at (0.5,1.5) {$\MarkovWeight{F}_{p_1,nq_0}$};
		\node at (1.5,1.5) {$\dots$};
		\node at (2.5,1.5) {$\MarkovWeight{F}_{p_N,nq_0}$};
		\node at (1.5,2.5) {$\MarkovWeight{T}$};
	\end{tikzpicture}
	\right)
	=
	\Phi_{P}^{S,q_0,r+(n-1)q_0}\left(
	\begin{tikzpicture}[guillpart,yscale=1.25,xscale=2.1]
		\fill[guillfill] (0,0) rectangle (3,3);
		\draw[guillsep] (0,0)--(0,3)--(3,3)--(3,0) (1,0)--(1,2) (2,0)--(2,2) (0,1)--(3,1) (0,2)--(3,2);
		\node at (0.5,0.5) {$v_1$};
		\node at (1.5,0.5) {$\dots$};
		\node at (2.5,0.5) {$v_N$};
		\node at (0.5,1.5) {$\MarkovWeight{F}_{p_1,q_0}$};
		\node at (1.5,1.5) {$\dots$};
		\node at (2.5,1.5) {$\MarkovWeight{F}_{p_N,q_0}$};
		\node at (1.5,2.5) {$\MarkovWeight{T}'$};
	\end{tikzpicture}
	\right)
	\\
	&= \lambda_{Pq_0} 	\begin{tikzpicture}[guillpart,xscale=1.5]
		\fill[guillfill] (0,0) rectangle (3,2);
		\draw[guillsep] (0,0)--(0,2)--(3,2)--(3,0) (1,0)--(1,1) (2,0)--(2,1) (0,1)--(3,1);
		\node at (0.5,0.5) {$v_1$};
		\node at (1.5,0.5) {$\dots$};
		\node at (2.5,0.5) {$v_N$};
		\node at (1.5,1.5) {$\MarkovWeight{T}'$};
	\end{tikzpicture}
	= \lambda_{Pq_0} 	\begin{tikzpicture}[guillpart,yscale=1.25,xscale=3.6]
		\fill[guillfill] (0,0) rectangle (3,3);
		\draw[guillsep] (0,0)--(0,3)--(3,3)--(3,0) (1,0)--(1,2) (2,0)--(2,2) (0,1)--(3,1) (0,2)--(3,2);
		\node at (0.5,0.5) {$v_1$};
		\node at (1.5,0.5) {$\dots$};
		\node at (2.5,0.5) {$v_N$};
		\node at (0.5,1.5) {$\MarkovWeight{F}_{p_1,(n-1)q_0}$};
		\node at (1.5,1.5) {$\dots$};
		\node at (2.5,1.5) {$\MarkovWeight{F}_{p_N,(n-1)q_0}$};
		\node at (1.5,2.5) {$\MarkovWeight{T}$};
	\end{tikzpicture}
\end{align*}
Using the recursion, we now use the existence of the morphism $\Phi_P^{S,(n-1)q_0,r}$ and apply it on the previous expression to obtain:
\begin{align*}
	&\Phi_{P}^{S,(n-1)q_0,r}\circ\Phi_{P}^{S,q_0,r+(n-1)q_0}\left(
	\begin{tikzpicture}[guillpart,yscale=1.25,xscale=2.3]
		\fill[guillfill] (0,0) rectangle (3,3);
		\draw[guillsep] (0,0)--(0,3)--(3,3)--(3,0) (1,0)--(1,2) (2,0)--(2,2) (0,1)--(3,1) (0,2)--(3,2);
		\node at (0.5,0.5) {$v_1$};
		\node at (1.5,0.5) {$\dots$};
		\node at (2.5,0.5) {$v_N$};
		\node at (0.5,1.5) {$\MarkovWeight{F}_{p_1,nq_0}$};
		\node at (1.5,1.5) {$\dots$};
		\node at (2.5,1.5) {$\MarkovWeight{F}_{p_N,nq_0}$};
		\node at (1.5,2.5) {$\MarkovWeight{T}$};
	\end{tikzpicture}\right)
\\
	&=\lambda_{Pq_0} 	\Phi_P^{S,(n-1)q_0,r}\left(\begin{tikzpicture}[guillpart,yscale=1.25,xscale=3.6]
		\fill[guillfill] (0,0) rectangle (3,3);
		\draw[guillsep] (0,0)--(0,3)--(3,3)--(3,0) (1,0)--(1,2) (2,0)--(2,2) (0,1)--(3,1) (0,2)--(3,2);
		\node at (0.5,0.5) {$v_1$};
		\node at (1.5,0.5) {$\dots$};
		\node at (2.5,0.5) {$v_N$};
		\node at (0.5,1.5) {$\MarkovWeight{F}_{p_1,(n-1)q_0}$};
		\node at (1.5,1.5) {$\dots$};
		\node at (2.5,1.5) {$\MarkovWeight{F}_{p_N,(n-1)q_0}$};
		\node at (1.5,2.5) {$\MarkovWeight{T}$};
	\end{tikzpicture}\right) 
	= \lambda_{Pq_0}\lambda_{P(n-1)q_0} 
	\begin{tikzpicture}[guillpart,xscale=1.5]
		\fill[guillfill] (0,0) rectangle (3,2);
		\draw[guillsep] (0,0)--(0,2)--(3,2)--(3,0) (1,0)--(1,1) (2,0)--(2,1) (0,1)--(3,1);
		\node at (0.5,0.5) {$v_1$};
		\node at (1.5,0.5) {$\dots$};
		\node at (2.5,0.5) {$v_N$};
		\node at (1.5,1.5) {$\MarkovWeight{T}$};
	\end{tikzpicture}
\end{align*}
Using the geometric property of the sequence $(\lambda_p)$, we define $\Phi_P^{S,nq_0,r}=\Phi_{P}^{S,(n-1)q_0,r}\circ\Phi_P^{S,q_0,r+(n-1)q_0}$ and it provides the expected morphisms. For $r=0$, the same construction holds without the additional term $\MarkovWeight{T}$ on top of the vertical products.
\end{proof}
\begin{coro}
	In the case of discrete state space $\setL=\setN$, it is sufficient to prove the existence of the morphisms $\Phi_\bullet^{S,1,r}$ for any $r\in\setL$.
	
	In the case of continuous state space $\setL=\setN$, it is sufficient to prove the existence of the morphisms $\Phi_\bullet^{S,q,r}$ for any $r\in\setL$ and any $q\leq \epsilon$ for some $\epsilon \geq 0$.
\end{coro}

This corollary is the 2D generalization of the 1D result that states that diagonalizing a matrix $A$ provides a diagonalization of $A^k$ for any $k\in\setN^*$. In the case of continuous space and a semi-group $(Q_t)_{t\geq 0}$, it is enough to diagonalize the matrices $Q_t$ in a neighbourhood of $0^+$.

The reader may wonder whether the dependence of the morphisms on $r$ and the necessity to include the elements $\MarkovWeight{T}$ create practical difficulties. We encourage the reader to proceed directly to Section~\ref{sec:eigensouth:canonicalMarkov} to verify that, at least in simple "factorized" cases as in the present paper, the dependence of the morphisms on $r$ is trivial and should not be a matter of worry.

\subsection{A digression on the one-dimensional case and an important remark on eigenvalues and morphisms}\label{sec:onedimdiagowithmorph}

Our proposition of introduction of morphisms in the definition of eigen-elements for $\Guill_2$-algebras seems to be a new ingredient absent in the one-dimensional case but this is not precisely the case: we may use similar definitions in dimension one but we will show that they will degenerate back to the standard definition.

We consider the $\Guill_1$-algebra $\ca{A}_p = \End(V)$ (indeed a $\Ass$-algebra) for all $p\in\setN^*$, where $V$ is a finite-dimensional real Hermitian vector space and fix a semi-group $(F_p)$ given by $F_p= A^p$ for some $A\in\End(V)$. We also introduce $\ca{A}_{\infty_L}=V$ and $\ca{A}_{\infty_R}=V$ as boundary spaces and $\ca{A}_{\infty_{LR}}=\setR$, with a pairing $\ca{A}_{\infty_L}\otimes\ca{A}_{\infty_R} \to \ca{A}_{\infty_{LR}}$ given by the scalar product.

Given $v\in\ca{A}_{\infty_R}$ ---a vector in this case---, we may consider the collection of elements $A^q v$ and $T A^q v$, with $T\in\ca{A}_{r}=\End(V)$ for any $q,r\in\setN^*$. Mimicking Definition~\ref{def:eigenalgebrauptomorphims}, we may say that $v$ is an eigenvector up to morphism with eigenvalue $\lambda$ if and only if there exist morphisms $\Phi^{R,q,r} : \ca{A}_{\infty_R}\to\ca{A}_{\infty_R}$, $r\in\setN$ such that$\Phi^{R,q,r}(T A^q v) = \lambda^q T v$
Here, $\ca{A}_{\infty_R}$ is only a vector space without any transverse $\Guill_1$-structure so the morphisms are just linear maps. We now study the linear map on $V$ given by $\Phi^{R,q,r}(u)= U_{q,r} u$ for some $U_{q,r} \in\End(V)$ which satisfy $U_{q,r} TA^qv=\lambda^q Tv$ for any $T\in\ca{A}_r=\End(V)$. 

We now show that the only possible choice of $U_{q,r}$ is $U_{q,r}=\gamma^q \id$ for some $\gamma\in\setR$. This proof is a nice 1D exercise in preparation to the 2D proof of Theorem~\ref{theo:reductionofmorphisms:canostruct} below. 

We first use the fact that the maps $\Phi$ commute with the action of any $T$. For any rank-one projector $T x = \scal{w'}{x} w$, we have
\[
\scal{w'}{A^q v}U_{q,r} w = \lambda^q \scal{w'}{v} w
\]
For any $w \in V$, there exists thus $c_{q,r,u}$ such that $U_{q,r}=c_{q,r,u} u$ and thus $U_{q,r}= c_{q,r} \id$ for some $c_{q,r}\in\setR$. We then have, for any $w\in V$,
$
c_{q,r}\scal{w'}{A^q v} =  \lambda^q \scal{w'}{v}
$
and thus $c_{q,r} A^q v =  \lambda^q v$ for all $r\in\setN^*$. The vector is assumed to be non-zero and thus, for all $r\in\setN^*$, $c_{q,r}=\ti{c}_q$ for some $\ti{c}_q\in\setR$. In order to study the dependence on $q$, we now consider, for any $q'\in\setN^*$, the element $T=T'A^{q'}$ with $T'\in\ca{A}_r$ and thus $T\in\ca{A}_{r+q'}$. We now have, for all $q,q'\in\setN^*$,
\[
\ti{c}_{q} T' A^{q+q'} v = \Phi^{R,q,q'+r}( T A^q v) = \lambda^{q} T  v = \lambda^q T' A^{q'}v
\]
We may further compose this identity by $\Phi^{R,q',r}$ and obtain
\[
\ti{c}_{q}\ti{c}_{q'} T' A^{q+q'} v = \lambda^{q+q'} T' v
\]
On the other hand, we may also have applied $\Phi^{R,q+q',r}$ on $T'A^{q+q' v}$ and obtain
\[
\ti{c}_{q+q'} T' A^{q+q'} v = \lambda^{q+q'} T' v
\]
Since the vectors on the r.h.s. may be chosen to be non-zero, we now obtain that $\ti{c}_{q} = \gamma^q$ for some $\gamma\in\setC$. 

One may remark that we cannot go any further by requiring $\gamma=1$. The parameter $\gamma$ is left free and may absorb the eigenvalue $\lambda$ by redefining the morphisms $\ti{\Phi}^{R,q,r} = \lambda^{-q}\Phi^{R,q,r}$ and we see that morphisms and eigenvalues are intrinsically related.

\subsection{Morphisms on the canonical structure for Markov processes}\label{sec:eigensouth:canonicalMarkov}

\subsubsection{Notations}
We now consider a two-dimensional Markov process in discrete space with a canonical space $\ca{T}_{\PatternShapes(\patterntype{fp}^*)}$ and a boundary weight with a ROPErep over a ROPE $\ca{B}_{\PatternShapes(\patterntype{fp}^*)}$ as in the theorems of structural stability~\ref{theo:stability}, for which all the computations took place in the $\Guill_2^{(\patterntype{fp}^*)}$-algebra $\ca{E}_{p,q}=\ca{T}_{p,q}\otimes\ca{B}_{p,q}$. We discuss only the discrete space case to have easier notations for morphisms closer to the traditional linear algebra.

We recall that the structure of the South half-strip spaces \[
\ca{E}_{p,\infty_S}= V(S_1)^{\otimes p}\otimes \left(\oplus_{r\in\setL} \End(V(S_2))^{\otimes r}\right)\otimes \ca{B}_{p,\infty_S}
\]
The West-East-products act on these spaces in the following way: the elements of $V(S_1)$ are concatenated, i.e. tensorized and thus it generates a free algebra over the basis vectors of $V(S_1)$, the elements of $\End(V(S_2))^{\otimes r}$ are multiplied with the usual composition product, the elements of $\ca{B}_{p,\infty_S}$ are multiplied using the West-East product of the ROPE.

The 2D-semi-group $\MarkovWeight{F}_\bullet$ corresponds to the surface powers $\MarkovWeight{F}_{p,q}=\MarkovWeight{F}^{[p,q]}$ of a fixed element \[
\MarkovWeight{F}\in\ca{T}_{1,1}=\End(V(S_1))\otimes\End(V(S_2))
\]

The morphism $\psi_p^{S,\MarkovWeight{F},\infty_S,r}$ maps an element in the subspace $V(S_1)^p\otimes \End(V(S_2))^{q}\otimes \ca{B}_{p,\infty_S}$ to an element of the subspace $V(S_1)^p\otimes \End(V(S_2))^{q+r}\otimes \ca{B}_{p,\infty_S}$ and acts as the identity on the factor $\ca{B}_{p,\infty_S}$. We will now see that the morphisms $\Phi_p^{S,q,r}$ must leave invariant various factors of $\ca{E}_{p,\infty_S}$ and can be expressed out of a small number of simpler morphisms acting on $\End(V(S_2))\otimes \ca{B}_{p,\infty_S}$.

\subsubsection{Reduction of morphisms}

In the embedding of the ROPErep of a boundary weight in the structure $\ca{E}_{\bullet}$, we consider the element~\eqref{eq:fromROPEreptoEpq:South} in the subspace 
\[
V(S_1) \otimes \setK \otimes \ca{B}_{1,\infty_S} \simeq V(S_1) \otimes  \ca{B}_{1,\infty_S}
\]
where $\setK$ corresponds to the convention $\End(V(S_2))^{\otimes 0}=\setK$. We are thus interested in the one-dimensional system of generators given by
\begin{align*}
	\ca{V}_1 &=\setK \ha{A}_S & \ca{V}_p &= \{0\}, \quad p>1
\end{align*} 
where $\ha{A}_S= \sum_{x\in S_2} e^*_x \otimes A(x)$ with $A(x)\in\ca{B}_{1,\infty_S}$. The unknown elements are precisely the $|S_2|$ elements $A(x)$.

We introduce the following notation in $\ca{E}_{p,\infty_S}$ to simplify equations of the following theorem. Any element $v \in \ca{E}_{p,\infty_S}$ can be decomposed as a finite sum
\[
v = \sum_{x\in S_1^p} e^*_x \otimes v(x),
\qquad
v(x) \in \left(\oplus_{r\in\setL} \End(V(S_2))^{\otimes r}\right)\otimes \ca{B}_{p,\infty_S}
\]
The factors $v(x)$ are then represented graphically as
\begin{equation}\label{eq:factornotationSouth}
	v(x) =
	\begin{tikzpicture}[guillpart,yscale=1,xscale=1]
		\fill[guillfill] (0,0) rectangle (1,1.5);
		\draw[guillsep] (0,0)--(0,1.5)--(1,1.5)--(1,0);
		\node at (0.5,0.75) {$v$};
		\node at (0.5,1.5) [anchor = south] {$x$};
		\node at (0.5,1.5) [inner sep=2pt, diamond, fill]{};
	\end{tikzpicture}
\end{equation}
We now state the main reduction theorem.

\begin{theo}[reduction of morphisms on half-strips for canonical structures of Markov processes]\label{theo:reductionofmorphisms:canostruct}
In the $\Guill_2^{(\patterntype{hs}_S)}$-algebra $\ca{E}_{\PatternShapes(\patterntype{hs}_S)}$, let $(\ca{V}_{p})_{p\in\setL^*}$ be a collection of subspaces of the spaces $(\ca{E}_{p,\infty_S})_{p\in\setL^*}$ such that there exists $q_0\in\setL$ such that, for all $p\in\setL^*$,
\begin{equation}\label{eq:hyp:leveleigengenerator}
	\ca{V}_p \subset V(S_1)^{\otimes p}\otimes\End(V(S_2))^{\otimes q_0}\otimes \ca{B}_{1,\infty_S} =: \ca{E}_{p,\infty_S}^{[q_0]}
\end{equation}
The spaces $\ca{V}_\bullet$ is a system of eigen-generators of a 2D-semi-group $\MarkovWeight{F}^{[p,q]}$ up to morphisms $\Phi_\bullet^{S,q,r}$ with eigenvalue sequence $\lambda_{\bullet}$ if and only if there exist linear maps
\[
\ha{\phi}_p :\End(V(S_2))\otimes \End(V(S_2))^{\otimes q_0} \otimes \ca{B}_{p,\infty_S} \to \End(V(S_2))^{\otimes q_0} \otimes\ca{B}_{p,\infty_S}
\]
such that, for all $y\in S_1^p$, for all $v_1,\ldots,v_p\in\ca{V}_p$,
\begin{equation}\label{eq:reducedmorphism:eigeneq}
\ha{\phi}_p\left(
\begin{tikzpicture}[guillpart,yscale=1.5,xscale=1.75]
	\fill[guillfill] (0,0) rectangle (3,2);
	\draw[guillsep] (0,0)--(0,2)--(3,2)--(3,0) (0,1)--(3,1) (1,0)--(1,1) (2,0)--(2,1);
	\node at (0.5,0.5) {$v_1$};
	\node at (1.5,0.5) {$\dots$};
	\node at (2.5,0.5) {$v_p$};
	\node at (1.5,1.5) {$\MarkovWeight{F}^{[p,1]}$};
	\node at (1.5,2) [anchor=south] {$y$};
	\node at (1.5,2) [inner sep=2pt, diamond, fill]{};
\end{tikzpicture}
\right)
=\lambda_p 
\begin{tikzpicture}[guillpart,yscale=1.5,xscale=1.5]
	\fill[guillfill] (0,0) rectangle (3,1);
	\draw[guillsep] (0,0)--(0,1)--(3,1)--(3,0)  (1,0)--(1,1) (2,0)--(2,1);
	\node at (0.5,0.5) {$v_1$};
	\node at (1.5,0.5) {$\dots$};
	\node at (2.5,0.5) {$v_p$};
	\node at (1.5,1) [anchor=south] {$y$};
	\node at (1.5,1) [inner sep=2pt, diamond, fill]{};
\end{tikzpicture}
\end{equation}
\end{theo}
The proof below is constructive and builds the morphisms $\Phi_{\bullet}^{S,q,r}$ out of the maps $\ha{\phi}_p$.

\begin{proof}
We first start by characterizing acceptable morphisms in Definition~\ref{def:eigenalgebrauptomorphims} in the case of the canonical structure of Markov processes and a ROPErep. In all the proof below, an element $\ov{v}_p$ stands for an arbitrary element of $\ca{E}_{p,\infty_S}$ such that there exist generators $v_1,\ldots,v_p \in \ca{V}_1$ such that
\[
\ov{v}_p = \begin{tikzpicture}[guillpart,yscale=1.5,xscale=1.5]
	\fill[guillfill] (0,0) rectangle (3,1);
	\draw[guillsep] (0,0)--(0,1)--(3,1)--(3,0)  (1,0)--(1,1) (2,0)--(2,1);
	\node at (0.5,0.5) {$v_1$};
	\node at (1.5,0.5) {$\dots$};
	\node at (2.5,0.5) {$v_p$};
\end{tikzpicture}
\]

Hypothesis \eqref{eq:hyp:leveleigengenerator} implies that, for all $\MarkovWeight{T}\in \ca{T}_{p,r}$ and all $\ov{v}_p$,
\[
	\begin{tikzpicture}[guillpart,yscale=1.,xscale=1.]
		\fill[guillfill] (0,0) rectangle (2,3);
		\draw[guillsep] (0,0)--(0,3)--(2,3)--(2,0) (0,1)--(2,1) (0,2)--(2,2);
		\node at (1,0.5) {$\ov{v}_p$};
		\node at (1,1.5) {$\MarkovWeight{F}_{p,q}$};
		\node at (1,2.5) {$\MarkovWeight{T}$};
	\end{tikzpicture}
	\in V(S_1)^{\otimes p}\otimes\End(V(S_2))^{\otimes q_0+q+r}\otimes \ca{B}_{1,\infty_S} := \ca{E}_{p,\infty_S}^{[q_0+q+r]}
\]
The morphisms $\Phi_p^{S,q,r}$ on $\ca{E}_{p,\infty}$ must then satisfy 
\[
\Phi_p^{S,q,r}(\ca{E}_{p,\infty_S}^{[q_0+q+r]})\subset  \ca{E}_{p,\infty_S}^{[q_0]}
\]
and their actions on the other subspaces $\ca{E}_{p,\infty_S}^{[q']}$ may be chosen arbitrarily and are not relevant. A linear map $\Phi_p^{S,q,r} : \ca{E}_{p,\infty_S}^{[q_0+q+r]}\to  \ca{E}_{p,\infty_S}^{[q_0]}$ (we keep the same notation for the restriction of the morphism) can always be written as:
\begin{equation}\label{eq:morphismparam}
\Phi_p^{S,q,r}\left( e^*_{x} \otimes A_r \otimes B_{q}\otimes B'_{q_0}  \otimes b \right) = \sum_{y\in S_1^p} e^*_y \otimes \sum_{\alpha \in I}  Q_{\alpha,x,y}^{(1)} A_r Q_{\alpha,x,y}^{(2)} \otimes \phi_{\alpha,x,y}\left( B_{q}\otimes B'_{q_0}\otimes b \right)
\end{equation}
for all $x\in S_1^p$, $A_r\in\End(V(S_2))^r$, $B_{q}\in \End(V(S_2))^{\otimes q}$, $B'_{q_0}\in \End(V(S_2))^{\otimes q_0}$ and $b\in\ca{B}_{p,\infty_S}$, where $I$ is a large enough finite set, the $Q_{\alpha,x,y}^{(i)}$ are elements of $\End(V(S_2))^r$ and the $\phi_{\alpha,x,y}$ are linear maps 
\[
\phi_{\alpha,x,y} : \End(V(S_2))^{\otimes q+q_0}\otimes \ca{B}_{p,\infty_S} \to \End(V(S_2))^{\otimes q_0}\otimes \ca{B}_{p,\infty_S}
\]
Since $\MarkovWeight{T}$ can be chosen arbitrarily in $\ca{T}_{p,r}$, we may consider any choice
\[
\MarkovWeight{T} =  u_1 u_2^* \otimes A_r
\]
with $u_1,u_2\in V(S_2)^{\otimes p}$ and $A_r\in\End(V(S_2))^{\otimes r}$ and plug it into Definition~\ref{def:eigenalgebrauptomorphims}. We then obtain that the right part of \eqref{eq:morphismparam} does not depend on $x$ and $y$ and the action on $A_r$ must be the identity. We then obtain the following simplification
\[
\Phi_p^{S,q,r}\left( e^*_{x} \otimes A_r \otimes B_{q}\otimes B'_{q_0}  \otimes b \right) = e^*_x \otimes  A_r  \otimes \varphi_{p,q,r}\left( B_{q}\otimes B'_{q_0}\otimes b \right)
\]
for some linear map $\varphi_{p,q,r}:\End(V(S_2))^{\otimes q+q_0}\otimes \ca{B}_{p,\infty_S} \to \End(V(S_2))^{\otimes q_0}\otimes \ca{B}_{p,\infty_S}$.

We now write the constraints on the linear maps $\varphi_{p,q,r}$. To this purpose, we first write down in terms of coordinates the action of the 2D-semi-group:
\begin{align*}
\begin{tikzpicture}[guillpart,yscale=1.,xscale=1.5]
	\fill[guillfill] (0,0) rectangle (1,2);
	\draw[guillsep] (0,0)--(0,2)--(1,2)--(1,0) (0,1)--(1,1);
	\node at (0.5,0.5) {$\ov{v}_p$};
	\node at (0.5,1.5) {$\MarkovWeight{F}_{p,q}$};
\end{tikzpicture}
&=\sum_{x\in S_1^p} e^{*}_x\otimes \sum_{(w,z)\in S_2^q\times S_2^q}\sum_{u\in S_1^p} \MarkovWeight{F}_{p,q}(u,x,w,z) E_{w,z}\otimes v(u) 
=
\sum_{x\in S_1^p} e^*_x \otimes
\begin{tikzpicture}[guillpart,yscale=1.,xscale=1.75]
	\fill[guillfill] (0,0) rectangle (1,2);
	\draw[guillsep] (0,0)--(0,2)--(1,2)--(1,0) (0,1)--(1,1);
 	\node at (0.5,0.5) {$\ov{v}_p$};
	\node at (0.5,1.5) {$\MarkovWeight{F}_{p,q}$};
	\node at (0.5,2) [inner sep=2pt,diamond, fill] {};
	\node at (0.5,2) [anchor=south] {$x$};
\end{tikzpicture}
\end{align*}
where $E_{w,z}$ is the matrix with null coefficients excepted the coefficient $(w,z)$ which is equal to one. Given a decomposition \eqref{eq:factornotationSouth} of $v$, the definition of eigen-generators is equivalent to \eqref{eq:reducedmorphism:eigeneq} with linear map $\phi_{p,q,r}$ instead of $\ha{\phi}_p$.

One observes in particular that, for all $r\in\setL^*$, the condition remains the same. If one knows the morphisms $\Phi_p^{S,q,r}$ then one obtains a whole family of linear maps $\phi_{p,q,r}$. Conversely, if one knows one $\varphi_{p,1}$ (no dependence on $r$) satisfying this equation for $q=1$, we may now use Lemma~\ref{lemm:eigenalg:sufficientmorphisms} to reduce the construction of $\phi_{p,q,r}$ and hence the whole morphisms $\Phi_p^{S,q,r}$ out of this single $\varphi_{p,1}$. The previous theorem is then proved.
\end{proof}

\subsubsection{Computational considerations, eager vs lazy ROPES and general dimensional considerations}\label{sec:computational-considerations} This theorem has important computational consequences. First, the presence of a large number of linear maps $\Phi_{\bullet}^{S,q,r}$ in Definition~\ref{def:eigenalgebrauptomorphims} can be reduced to a much smaller of linear maps $\ha{\phi}_p$ acting only the "vertical" spaces made of tensor products of $\End(V(S_2))$ and $\ca{B}_{p,\infty_S}$. As soon as the spaces $\ca{B}_{p,\infty_S}$ are identified to some explicit spaces, the maps $\ha{\phi}_p$ and the vectors $v\in\ca{V}_p$ can be described and studied as explicit linear maps and vectors and one may start solving concrete equations. In particular the reduction to the case from a general $q$ to the case $q=1$ using Lemma~\ref{lemm:eigenalg:sufficientmorphisms} is close to the one-dimensional case: diagonalization of all $A^q$ (even in the sense of Section~\ref{sec:onedimdiagowithmorph}) can be obtained from the diagonalization of the single element $A$.

A major apparent worry about Definition~\ref{def:eigenalgebrauptomorphims} and its Markovian realization in the previous theorem~\ref{theo:reductionofmorphisms:canostruct} is that there is a whole collection of linear maps $\ha{\phi}_p$ \emph{without} any relation between them. The $\Guill_1$-morphism property on the restriction to the spaces $\ca{V}_\bullet$ is indeed a direct consequence of Definition~\ref{def:eigenalgebrauptomorphims}. In practice, there may be much too many degrees of freedom behind these maps $\ha{\phi}_p$ without any further assumptions on the the ROPE structure.

We illustrate this issue on the particular case of $\ca{B}_{p,\infty_S} =\Mat_{n_p,n_p}(\setR)$ with suitable $\Guill_1$-products (not necessarily the traditional matrix product) for some integer $n_p\geq 1$ and $\ca{V}_1 = \Vect(u)$ (one-dimensional space). The unknowns are thus the $|S_1|$ matrices 
\[
u(x) = \begin{tikzpicture}[guillpart,yscale=1.5,xscale=1.5]
	\fill[guillfill] (0,0) rectangle (1,1);
	\draw[guillsep] (0,0)--(0,1)--(1,1)--(1,0);
	\node at (0.5,0.5) {$u$};
	\node at (0.5,1) [inner sep=2pt,diamond, fill] {};
	\node at (0.5,1) [anchor=south] {$x$};
\end{tikzpicture}  \in \Mat_{n_1,n_1}(\setR),
\]
the eigenvalue $\lambda_1$ and the morphisms $\ha{\phi}_p$ for all $p\in\setN^*$. For every $p\in\setN^*$, \eqref{eq:reducedmorphism:eigeneq} corresponds to $|S_1|^p$ equations on $\Mat_{n_p,n_p}(\setR)$ where $\ha{\phi}_p$ is a linear map from $\Mat_{n_p |S_2|, n_p |S_2|}$ to $\Mat_{n_p,n_p}$.

If one considers for $\ca{B}_{p,\infty_S} =\Mat_{n_p,n_p}(\setR)$ a lazy ROPE (see Definition~\ref{def:eagerandlazy}) with $n_p=n_1^p$ then morphisms can always be found and the previous definition looks like empty in practice. If one considers an eager ROPE (see Definition~\ref{def:eagerandlazy}) and thus impose a constant morphisms $\ha{\phi}_p=\phi$, then the constraints \eqref{eq:reducedmorphism:eigeneq} on $\phi$ and the elements $u(x)$ are very restrictive. Since, at this level, the ROPE structure itself is also unknown, one may thus wish to work with a lazy ROPE and thus obtain many solutions. However, for practical probabilistic purposes, this approach is vain since corners have to be introduced to describe the Gibbs measure: Theorem~\ref{theo:1D:removingcoloursonboundaries} then enters into play and crushes any additional degrees of freedom introduced by lazy ROPES on the eager ROPE defined by $\Hom(\ca{B}_{\infty_{b},\infty_a},\ca{B}_{\infty_{b},\infty_a})$, $a\in\{S,N\}$, $b\in\{W,E\}$. This is the reason why the examples described in this paper and in subsequent papers such as \cite{BodiotSimon} are all described by eager ROPEs.

Why then bother with other possible ROPEs than eager ones? There are three main reasons, one geometric, one algebraic and one probabilistic and all three of them go beyond the scope of the present paper. 

On the geometrical aspect, an eager ROPE corresponds to a construction of spaces $\ca{E}_{p,\infty_S}^{[0]}$ very similar to the construction of $\ca{T}_{p,q}$: each of the two directions carries an associative algebra, products along this direction is the internal product of the algebra and products in the transverse direction are given by tensor products (concatenation). There is no any additional subtle operations (excepted for the morphisms) mixing both dimensions and this may not be surprising in the guillotine framework. We however believe that for general domains of $\setZ^2$ or other lattices than the square lattices (for example the triangular lattices with three translation axes in dimension two), non-trivial constraints must exist between the structures along each dimension of the lattice. In this case, it is not clear whether eager ROPEs are well-suited or if the dependence of boundary spaces of the ROPEs on the geometrical colours may be relaxed.

On the probabilistic side, the Gibbs measure on $\setZ^2$ is related to the construction of the whole ROPE. Partial objects on half-strips are related to the action or an arbitrary large number of face weights only in one of the four directions as illustrated in Theorem~\ref{theo:stability} and in Definition~\ref{def:eigenalgebrauptomorphims} and this must not be a surprise if they capture only information coming from infinity in the corresponding direction and not all the information of the Gibbs measure; an additional corner element can be interpreted as providing information coming from infinity in the transverse direction. In models with a non-unique Gibbs measure in which extremal Gibbs measures break the symmetry between the two-dimensions, we expect that the constraints on a half-strip element may not be sufficient to distinguish extremal Gibbs measures related to different asymptotic behaviour in the transverse direction. 

On the algebraic side, the relation between $\Guill_1$- and $\Ass$-algebras given by Theorem~\ref{theo:1D:removingcoloursonboundaries} can be interpreted at the light of representation theory. When representing a group or an algebra on a vector space, part of the properties are inherited from the group or algebra structure itself whereas the other properties are inherited from the representation space itself. Here, properties on half-strip eigen-elements inherited from Definition~\ref{def:eigenalgebrauptomorphims} only (for any ROPE, lazy or eager or any other choice)  are valid for any corner choice whereas the choice of an eigen-corner space as defined below introduces additional properties on the corresponding eager ROPE.

More mathematical formalization of these intuitions deserve to be addressed in future papers by considering various concrete models and a deeper algebraic understanding of the extended $\Guill$-operads.

\subsection{From a traditional elementary exercise to the Yang-Baxter equation using morphisms}\label{sec:commutuptomorph}

A traditional basic exercise in linear algebra is the proof that, if $AB=BA$, then the eigen-spaces of $A$ are stable under the action of $B$. This result is widely use all through the mathematical literature to diagonalize $A$. For example, if $B$ encodes some symmetry then the eigenvectors of $A$ must be consistent with this symmetry. Another example is the theory of integrable systems, for which one looks for a whole family of commuting operators in order to have an easier codiagonalization. It is interesting to see if such a simple result admits such a generalization in our operadic framework.

\subsubsection{Commutation morphisms}
The trick consists once again to define the notion of commutation up to morphisms. 

\begin{defi}\label{def:commutuptomorph}
	 In the discrete space case $\setL=\setN$, a South-commutation morphism of $\MarkovWeight{A}$ over $\MarkovWeight{B}$ in $\ca{A}_{1,1}$ is a collection of linear maps $C_{p}^{r,s} : \ca{A}_{p,r+2+s}\to\ca{A}_{p,r+2+s}$, $p\in\setL^*$, $s\in\setL$, $r\in\setL\cup\{\infty_S\}$, such that, for all $M_1\in\ca{A}_{p,r}$ and $M_2\in\ca{A}_{p,s}$, 
	\begin{equation}\label{eq:def:commutuptmorph:single}
		C_{p}^{r,s}\left( 
			\begin{tikzpicture}[guillpart,yscale=1.5,xscale=2]
				\fill[guillfill] (0,0) rectangle (1,4);
				\draw[guillsep] (0,0)--(0,4) (1,0)--(1,4) (0,1)--(1,1) (0,2)--(1,2) (0,3)--(1,3)(0,4)--(1,4);
				\draw[guillsep, dashed] (0,0)--(1,0);
				\node at (0.5,0.5) {$M_1$};
				\node at (0.5,1.5) {$\MarkovWeight{B}^{[p,1]}$};
				\node at (0.5,2.5) {$\MarkovWeight{A}^{[p,1]}$};
				\node at (0.5,3.5) {$M_2$};
			\end{tikzpicture}
		\right)
		= \begin{tikzpicture}[guillpart,yscale=1.5,xscale=2]
			\fill[guillfill] (0,0) rectangle (1,4);
			\draw[guillsep] (0,0)--(0,4) (1,0)--(1,4) (0,1)--(1,1) (0,2)--(1,2) (0,3)--(1,3)(0,4)--(1,4);
			\draw[guillsep, dashed] (0,0)--(1,0) ;
			\node at (0.5,0.5) {$M_1$};
			\node at (0.5,2.5) {$\MarkovWeight{B}^{[p,1]}$};
			\node at (0.5,1.5) {$\MarkovWeight{A}^{[p,1]}$};
			\node at (0.5,3.5) {$M_2$};
		\end{tikzpicture}
	\end{equation}
\end{defi}
\begin{rema}
	In this definition, contrary to the one-dimensional case, there is no symmetry between $\MarkovWeight{A}$ and $\MarkovWeight{B}$ since the maps $C_p^{r,q,q',s}$ may not be invertible. If there are invertible, then $\MarkovWeight{A}$ and $\MarkovWeight{B}$ are said to commute up to morphisms and symmetry is restored.
\end{rema}
\begin{lemm}
	If there exists a South-commutation morphisms of $\MarkovWeight{A}$ over $\MarkovWeight{B}$ (in $\ca{A}_{1,1}$), then there exist linear maps $C_p^{r,q,q',s}$ from $\ca{A}_{p,r+q+q'+s}$ to itself, with $p,q,q'\in\setL^*$, $s\in\setL$ and $r\in\setL\cup\{\infty_S\}$ such that:
	\begin{equation}\label{eq:commutuptomorph}
		C_{p}^{r,q,q',s}\left( 
		\begin{tikzpicture}[guillpart,yscale=1.5,xscale=2]
			\fill[guillfill] (0,0) rectangle (1,4);
			\draw[guillsep] (0,0)--(0,4) (1,0)--(1,4) (0,1)--(1,1) (0,2)--(1,2) (0,3)--(1,3) (0,4)--(1,4);
			\draw[guillsep, dashed] (0,0)--(1,0);
			\node at (0.5,0.5) {$M_1$};
			\node at (0.5,1.5) {$\MarkovWeight{B}^{[p,q]}$};
			\node at (0.5,2.5) {$\MarkovWeight{A}^{[p,q']}$};
			\node at (0.5,3.5) {$M_2$};
		\end{tikzpicture}
		\right)
		= \begin{tikzpicture}[guillpart,yscale=1.5,xscale=2]
			\fill[guillfill] (0,0) rectangle (1,4);
			\draw[guillsep] (0,0)--(0,4) (1,0)--(1,4) (0,1)--(1,1) (0,2)--(1,2) (0,3)--(1,3) (0,4)--(1,4);
			\draw[guillsep, dashed] (0,0)--(1,0);
			\node at (0.5,0.5) {$M_1$};
			\node at (0.5,2.5) {$\MarkovWeight{B}^{[p,q]}$};
			\node at (0.5,1.5) {$\MarkovWeight{A}^{[p,q']}$};
			\node at (0.5,3.5) {$M_2$};
		\end{tikzpicture}
	\end{equation}
\end{lemm}
\begin{proof}
	The proof is similar to the one of Lemma~\ref{lemm:eigenalg:sufficientmorphisms} obtained by recursion on $q$ and $q'$ and by replacing elements $M_i$ by $M_i$ multiplied by $\MarkovWeight{A}^{[p,q'-1]}$ or $\MarkovWeight{B}^{[p,q-1]}$ and successive composition of the morphisms $C_p^{r,s}$.
\end{proof}

\subsubsection{Commutation morphisms and eigen-generators.}
\begin{prop}
	Let $\MarkovWeight{A}_\bullet$ and $\MarkovWeight{B}_\bullet$ be two 2D-semi-groups in discrete space generated by some elements $\MarkovWeight{A}$ and $\MarkovWeight{B}$ in $\ca{A}_{1,1}$ such that there exists a South-commutation morphisms of $\MarkovWeight{A}$ over $\MarkovWeight{B}$. Let $(\ca{V}_p)_{p\in\setL^*}$ be a system of South eigen-generators up to morphisms with eigenvalue sequence $\Lambda_\bullet$ of the 2D-semi-group $\MarkovWeight{A}_\bullet$. Then, for any $q\in\setL^*$, the spaces $(\ca{V}^{\MarkovWeight{B},(q)}_p)_{p\in\setL^*}$ defined by 
	\[
		\ca{V}^{\MarkovWeight{B},(q)}_p = \left\{
			\begin{tikzpicture}[guillpart,scale=1.5]
				\fill[guillfill] (0,0) rectangle (1,2);
				\draw[guillsep] (0,0)--(0,2)--(1,2)--(1,0) (0,1)--(1,1);
				\node at (0.5,0.5) {$v$};
				\node at (0.5,1.5) {$\MarkovWeight{B}_{p,q}$};
			\end{tikzpicture}
			; v\in\ca{V}_{p}
		\right\}
	\]
	is also a system of South eigen-generators up to morphisms of $\MarkovWeight{A}_\bullet$ with the same eigenvalue eigenvalue sequence $\Lambda_\bullet$.
\end{prop}
\begin{proof}
We consider a sequence $(v_i)_{1\leq i\leq n}$ with $v_i\in\ca{V}_{p_i}$ and define $P=\sum_i p_i$. We now use \eqref{eq:commutuptomorph} where $M_1$ is the $\Guill_1$-product of the $v_i$ and $M_2=\MarkovWeight{T}$ is an arbitrary element of $\ca{A}_{P,s}$ and compose it by $\Phi_{P}^{(\MarkovWeight{A}),S,q',q+s}$ and obtain:
\begin{align*}
	\Phi_{P}^{(\MarkovWeight{A}),S,q',q+s}
	\circ C_p^{r,q,q',s}\left(
	\begin{tikzpicture}[guillpart,yscale=1.25,xscale=1.25]
		\fill[guillfill] (0,0) rectangle (3,4);
		\draw[guillsep] (0,0)--(0,4) (3,0)--(3,4) (0,1)--(3,1) (0,2)--(3,2) (0,3)--(3,3) (0,4)--(3,4) (1,0)--(1,1) (2,0)--(2,1);
		\node at (0.5,0.5) {$v_1$};
		\node at (1.5,0.5) {$\dots$};
		\node at (2.5,0.5) {$v_n$};
		\node at (1.5,1.5) {$\MarkovWeight{B}^{[p,q]}$};
		\node at (1.5,2.5) {$\MarkovWeight{A}^{[p,q']}$};
		\node at (1.5,3.5) {$\MarkovWeight{T}$};
	\end{tikzpicture}
	\right)
	= \Lambda^{pq'}\begin{tikzpicture}[guillpart,yscale=1.25,xscale=1.25]
		\fill[guillfill] (0,0) rectangle (3,3);
		\draw[guillsep] (0,0)--(0,3) (3,0)--(3,3) (0,1)--(3,1) (0,2)--(3,2)  (0,3)--(3,3) (1,0)--(1,1) (2,0)--(2,1);
		\node at (0.5,0.5) {$v_1$};
		\node at (1.5,0.5) {$\dots$};
		\node at (2.5,0.5) {$v_n$};
		\node at (1.5,1.5) {$\MarkovWeight{B}^{[p,q]}$};
		\node at (1.5,2.5) {$\MarkovWeight{T}$};
		\end{tikzpicture}
\end{align*}
Splitting $\MarkovWeight{B}^{[p,q]}$ in a horizontal product $\MarkovWeight{B}^{[p_i,q]}$ and applying each $\MarkovWeight{B}^{[p_i,q]}$ vertically on $v_i$ provide the theorem.
\end{proof}

\subsubsection{Relation with Yang-Baxter equation.}

An interesting observation is the description of South-commutation morphisms of $\MarkovWeight{A}$ over $\MarkovWeight{B}$ in the case of the canonical $\Guill_2$-structure associated to two-dimensional Markov processes. The same scenario as the previous one for morphisms of eigen-generators holds: all the morphisms $C_p^{r,s}$ from $\End(V(S_1))^{\otimes p}\otimes \End(V(S_2))^{\otimes r+2+s}$ to itself are built out a single morphism $\ti{C} : \End(V(S_2))^{\otimes 2} \to \End(V(S_2))^{\otimes 2}$ through
\[
	C_p^{r,s}\left(
	 S^{(1)}\otimes \left( \bigotimes_{k=1}^{r+2+s} T^{(2)}_k\right) \right)
	 =
	 S^{(1)}\otimes \left( \bigotimes_{k=1}^{r} T^{(2)}_k\right)\otimes \ti{C}\left( 
 	T^{(2)}_{r+1} \otimes T^{(2)}_{r+2}
 	 \right)\otimes \left( \bigotimes_{k=r+3}^{r+2+s} T^{(2)}_k\right)
\]
We are thus left to study the morphism of (associative) algebra $\ti{C}$ from $\End(V(S_2))^{\otimes 2}$ to itself. The typical case corresponds to the existence of an invertible operator $R\in \End(V(S_2))^{\otimes 2}$ such that
\[
\ti{C}(T\otimes T')  = R (T\otimes T') R^{-1}
\]
such that one has
\[
\ti{C}(T_1T_2\ldots T_p \otimes T'_1T'_2\ldots T'_p)= \ti{C}(T_1\otimes T'_1)\ldots \ti{C}(T_p\otimes T'_p)
\]
The commutation relation up to morphisms \eqref{eq:commutuptomorph} then becomes 
\begin{equation}\label{eq:protoYangBaxter}
	R_{23} B_{12} A_{13} R_{23}^{-1} = A_{12} B_{13} 
\end{equation}
where we have enumerated the spaces $V(S_1)\otimes V(S_2)\otimes V(S_2)$ from left to right from $1$ to $3$. This is precisely the algebraic prototype of the \emph{Yang-Baxter equation} extensively used in integrable statistical mechanics to codiagonalize families of weight faces. 

Our framework does not provide solutions to the Yang-Baxter equation but places it in a wider operadic context where the relation with diagonalization is made deeper since it appears as one of the simplest two-dimensional generalization up to morphisms of the relation $AB=BA$.

\subsection{Generators vs generated algebras: an important remark}\label{sec:generators:remark}

Starting from a system of South eigen-generators $(\ca{V}_{p})$, one may build the following sub-algebras of $\ca{A}_{\bullet,\infty_S}$.
\begin{itemize}
	\item a sub-$\Guill_1$-algebra $\ca{S}_\bullet$ defined as the smallest sub-$\Guill_1$-algebra of $\ca{A}_{\bullet,\infty_S}$ that contains all the spaces $\ca{V}_p$.
	\item for each level $q\in\setL^*$, a sub-$\Guill_1$-algebra $\ca{S}^{(q)}_\bullet$ defined as the smallest sub-$\Guill_1$-algebra of $\ca{A}_{\bullet,\infty_S}$ that contains all the images of the spaces $\ca{V}_p$ under $\psi^{S,\MarkovWeight{F},\infty_S,q}$.
	\item the sub-$\Guill_1$-algebra $\ti{\ca{S}}_{\bullet,\infty_S}$ generated by all the generator spaces $\ca{V}_p$ and the action of the 2D-semi-group.
\end{itemize}
The first item describes the simplest algebra to be considered when considering the spaces $(\ca{V}_p)$ alone. The third item describes the (South part) of the smallest $\Guill_2^{(\patterntype{hs}_S)}$-algebra $\ti{\ca{S}}_{\PatternShapes(\patterntype{hs}_S)}$ that contains the $\ca{V}_p$ and the elements of the semi-group $\MarkovWeight{F}_\bullet$. Due to the absorbing property of the colour $\infty_S$, i.e. $\infty_S+q=\infty_S$ for any finite $q$, we have again $\ti{\ca{S}}_{p,q}=\setK\MarkovWeight{F}_{p,q}$ for finite $p$ and $q$ but the spaces $\ti{\ca{S}}_{p,\infty_S}$ contain at least all the elements $v_p$ and all the elements 
\begin{equation}\label{eq:actionFonV:lengthq}
\begin{tikzpicture}[guillpart,scale=1.5]
	\fill[guillfill] (0,0) rectangle (1,2);
	\draw[guillsep] (0,0)--(0,2)--(1,2)--(1,0) (0,1)--(1,1);
	\node at (0.5,0.5) {$v_p$};
	\node at (0.5,1.5) {$\MarkovWeight{F}_{p,q}$};
\end{tikzpicture}
\end{equation}
for $v_p\in\ca{V}_p$ and $q\in\setL^*$. In this sense, the algebra $\ti{\ca{S}}_{\bullet,\infty_S}$ is thus very natural. This is however not the spaces used in the previous definition, where the second item, which provides a collection of spaces indexed by $q\in\setL^*$, is preferred. This is why we prefer to formulate the definitions in terms of spaces of eigen-generators instead of some "eigen-algebras" since it may be trickier to define and to use. Moreover, for the study of concrete Markov processes, we are interested in the computational task of obtaining a finite set of generators from a well-defined set of equations.

A major difference between the second item and the third item is the possibility, in the third item, to consider elements of the type:
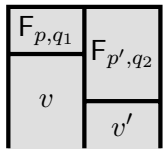
\[
\begin{tikzpicture}[guillpart,yscale=1.25,xscale=2.]
	\fill[guillfill] (0,0) rectangle (2,3);
	\draw[guillsep] (0,0)--(0,3)--(2,3)--(2,0) (1,0)--(1,3) (0,2)--(1,2) (1,1)--(2,1);
	\node at (0.5,1) {$v$};
	\node at (1.5,0.5) {$v'$};
	\node at (0.5,2.5) {$\MarkovWeight{F}_{p,q_1}$};
	\node at (1.5,2) {$\MarkovWeight{F}_{p',q_2}$};
\end{tikzpicture}
\]
with $v\in\ca{V}_p$ and $v'\in\ca{V}_{p'}$ and different $q_1$ and $q_2$. These hook shapes of elements of the 2D-semi-group as well as the gluing of two generators $v$ and $v'$ at different heights are proscribed in the second item and the previous Definition~\ref{def:eigenalgebrauptomorphims}. We however believe that acquiring a deeper understanding of these types of objects in terms of notions "up to morphisms" unlocks the door toward the study of domains with arbitrary shapes (i.e. not only rectangles with guillotine cuts).

\subsection{From one side to the four sides}

Our definitions on the South side can be directly adapted to any other North, West or East side by moving all the diagrams with a suitable action of dihedral group. The role of the vertical and horizontal associativities may be permuted as well the left or right action of the element $\MarkovWeight{F}_{p,q}$ in each of the two directions but the commutative diagrams and the role of morphisms remain the same. We leave the reader transpose precisely the previous definitions on its own.

Given two system of eigen-generators on two boundary side out of the four ones, we now present how to combine them to arrive at a full boundary structure. There are two situations: either the two boundary sides are adjacent and corner eigen-generators have to be introduced (see Section~\ref{sec:addingcornertoeigenSouth}), or the two boundary sides are opposite and the corresponding generators can be paired to strip elements with doubly-infinite shapes (see Section~\ref{sec:pairingeigentostrips}).

\section{Adding a corner between adjacent eigen-generators up to morphisms.}\label{sec:addingcornertoeigenSouth}

We consider in this section a fixed $\Guill_2^{(\patterntype{c}_{SW})}$-algebra $\ca{A}_{\PatternShapes(\patterntype{c}_{SW})}$. 

\subsection{Preliminary remarks on the corner geometry}

In the previous section on eigen-elements for South half-strips, we had a one-dimensional action from the North ---as in traditional linear algebra-- on which we had an additional transverse associative structure. Adding a corner is a non-trivial consequence of dimension two and we may expect several new features, without any equivalent feature in dimension one. Before any further computation, we emphasize on the following observations:
\begin{itemize}
	\item there is \emph{no} action of $\ca{A}_{\PatternShapes(\patterntype{r})}$ on the corner spaces $\ca{A}_{\infty_b,\infty_a}$.
	\item instead, there are left or right actions from two different $\Guill_1$-algebras (related to half-strips), one with a West-East product $m_{WE}$ and one with a South-North product $m_{SN}$
	\item each of these actions is thus related to a $\Guill_{1,L}$- or $\Guill_{1,R}$-algebras 
	\item both actions are not necessarily commuting (no bimodule a priori) but are compatible (square associativity~\eqref{eq:guill2:interchangeassoc}) through the actions of $\ca{A}_{\PatternShapes(\patterntype{r})}$ on the half-strip $\Guill_1$-algebras, using guillotine partitions of the following type:
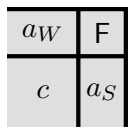
\[
	\begin{tikzpicture}[guillpart,yscale=1.25,xscale=1.25]
		\fill[guillfill] (0,0) rectangle (2.5,2.5);
		\draw[guillsep] (0,2.5)--(2.5,2.5)--(2.5,0) (1.5,0)--(1.5,2.5) (0,1.5)--(2.5,1.5);
		\node at (0.75,0.75) {$c$};
		\node at (2,0.75) {$a_S$};
		\node at (0.75,2) {$a_W$};
		\node at (2,2) {$\MarkovWeight{F}$};
	\end{tikzpicture}
\]
\end{itemize}

\subsection{Adding spaces on the left: generic vertical semi-groups on the West side}

In dimension one, $\Guill_1$-algebras $(\ca{A}_p)_{p\in\setL^*}$ can be enriched with a half-line space $\ca{A}_{\infty_L}$ to form a $\Guill_{1,L}$-algebra $(\ca{A}_p)_{p\in\setL^*\cup\{\infty_L\}}$, generalizing in the coloured case the notion of module over an algebra. 

The present section generalizes definition \ref{def:eigenalgebrauptomorphims} in order to include West spaces. For a fixed vertical colour, we have the following horizontal structures.
\begin{itemize}
	\item each (horizontal) $\Guill_1$-algebra $\ca{A}_{\bullet,q}$, for any fixed $q\in\setL^*$ is promoted to a $\Guill_{1,L}$-algebra by adding the West space $\ca{A}_{\infty_W,q}$.
	
	\item the South $\Guill_1$-algebra $\ca{A}_{\bullet,\infty_S}$ is promoted to a $\Guill_{1,L}$-algebra by adding the corner space $\ca{A}_{\infty_W,q}$.
\end{itemize}

We consider a 2D-semigroup $(\MarkovWeight{F}_{p,q})$ as in the previous section. For a fixed $p\in\setL^*$, each vertical semi-group $\MarkovWeight{F}_{p,\bullet}$ act on a South space $\ca{A}_{p,\infty_S}$. We must now add a new vertical semi-group $(M_q^W)$ on the West in order to act vertically on the corner space $\ca{A}_{\infty_W,\infty_S}$, i.e. a collection of elements $M_q^W\in\ca{A}_{\infty_W,q}$ such, for all $q_1,q_2\in\setL^*$:
\begin{equation}\label{eq:vertsemigroup:West}
	\begin{tikzpicture}[guillpart,scale=1.5]
		\fill[guillfill] (0,0) rectangle (1,2);
		\draw[guillsep] (0,0)--(1,0)--(1,2)--(0,2) (0,1)--(1,1);
		\node at (0.5,0.5) { $M^{W}_{q_1}$ };
		\node at (0.5,1.5) { $M^{W}_{q_2}$ };
	\end{tikzpicture}
	= M^{W}_{q_1+q_2}
\end{equation}
In the horizontal direction, there is a novelty due to the fact that $\infty_W$ is an absorbing colour, i.e. $\infty_W+p=\infty_W$ for finite $p$: the choice of $M_\bullet^W$ can not introduce any constraint on the 2D-semi-group $\MarkovWeight{F}_{\bullet}$ itself since $m_{WE}(M_q^W,\MarkovWeight{F}_{p,q})$ belongs to $\ca{A}_{\infty_W,q}$. We split the discussion in two steps: we first consider the case of a West semi-group $M_\bullet^W$ chosen arbitrarily, i.e. without any relation with $\MarkovWeight{F}_\bullet$ and, afterwards, we add the hypothesis that the elements $M_q^W$ form a system of West eigen-generators of $\MarkovWeight{F}_\bullet$.

Due to the absorbing colour $\infty_W$, the West elements $M^W_q$ and the 2D-semi-group $\MarkovWeight{F}_\bullet$ together generate subspaces
\begin{equation}
	\Vect_\setK\left(
	M^W_q, 	\begin{tikzpicture}[guillpart,yscale=1.5,xscale=1.6]
		\fill[guillfill] (0,0) rectangle (2,1);
		\draw[guillsep] (0,0)--(2,0)--(2,1)--(0,1) (1,0)--(1,1);
		\node at (0.5,0.5) {$M_q^{W}$};
		\node at (1.5,0.5) {$\MarkovWeight{F}_{p,q}$};
	\end{tikzpicture} ; p \in\setL^*
	\right)
\end{equation}
of $\ca{A}_{\infty_W,q}$ but, as stated in Section~\ref{sec:generators:remark}, these spaces may not be adequate. We prefer instead to consider the collection, indexed by $p\in\setL$, of vertical $\Guill_1$-algebras $\ca{M}^{W,(p)}_\bullet$ defined by:
\begin{align}\label{eq:eigencorner:westmodule:generators}
	\ca{M}^{W,(0)}_q &= \setK M_q^W 
	&
	\ca{M}^{W,(p)}_q &= \setK \begin{tikzpicture}[guillpart,yscale=1.5,xscale=1.6]
		\fill[guillfill] (0,0) rectangle (2,1);
		\draw[guillsep] (0,0)--(2,0)--(2,1)--(0,1) (1,0)--(1,1);
		\node at (0.5,0.5) {$M_q^{W}$};
		\node at (1.5,0.5) {$\MarkovWeight{F}_{p,q}$};
	\end{tikzpicture}
\end{align}
for $p\in\setL^*$. For each $p\in\setL$, these spaces $\ca{M}^{W,(0)}_q$ are one-dimensional and correspond to the elements of a vertical semi-group, hence mimicking and generalizing on the West the situation of a 2D-semi-group.

Using the notion of generators and the correct gradings, Definition~\ref{def:eigenalgebrauptomorphims} can then be directly generalized in the following way.

\begin{defi}\label{def:eigenalgebrauptomorph:leftextended}
	Let $\MarkovWeight{F}_{\bullet,\bullet}$ be a 2D-semi-group of $\ca{A}_{\PatternShapes(\patterntype{r})}$ and $M^W_\bullet$ a vertical semi-group of $\ca{A}_{\infty_W,\bullet}$. A collection of non-trivial sub-spaces $(\ca{V}_{p})_{p\in\setL^*\cup\{\infty_W\}}$ of $\ca{A}_{p,\infty_S})_{p\in\setL^*\cup\{\infty_W\}}$ is a left-extended system of eigen-$\Guill_1^L$-generators up to morphisms of $(\MarkovWeight{F}_{\bullet,\bullet},M^W_\bullet)$ with eigenvalue geometric sequences $(\lambda_r)_{r\in\setL^*}$ and $(\sigma^W_q)_{q\in\setL^*}$ if and only if there exists a collection of linear morphisms $\Phi_p^{S,q,r}: \ca{A}_{p,\infty_S}\to\ca{A}_{p,\infty_S}$, $p\in\setL\cup\{\infty_W\}$, $q\in\setL$, $r\in\setL$ such that:
		\begin{enumerate}[(i)]
			\item the spaces $(\ca{V}_{p})_{p\in\setL^*}$ are a system of eigen-$\Guill_1$-generators up to morphisms of $\MarkovWeight{F}_{\bullet,\bullet}$ with eigenvalue sequence $\lambda_\bullet$ up to the morphisms $\Phi_p^{S,q,r}$ with finite $p$.
			
			\item for any $N\in\setN$, for any $r\in\setL$, any sequence $(p_i)_{1\leq i\leq N}$ in $\setL^*$, for any $c\in\ca{V}_{\infty_W,\infty_S}$, for any sequence $(v_i)_{1\leq i\leq N}$ with $v_i\in\ca{V}_{p_i,\infty_S}$, for any $T_W\in \ca{A}_{\infty_W,r}$,
			\begin{equation}\label{eq:cornermorphism:def}
				\Phi_{\infty_W}^{S,q,r}\left(
					\begin{tikzpicture}[guillpart,yscale=1.25,xscale=2.]
						\fill[guillfill] (0,0) rectangle (4,3);
						\draw[guillsep] (4,0)--(4,3)--(0,3) (0,1)--(4,1) (0,2)--(4,2) (1,0)--(1,2) (2,0)--(2,2) (3,0)--(3,2);
						\node at (0.5,0.5) {$c$};
						\node at (1.5,0.5) {$v_1$};
						\node at (2.5,0.5) {$\dots$};
						\node at (3.5,0.5) {$v_N$};
						\node at (0.5,1.5) {$M^W_q$};
						\node at (1.5,1.5) {$\MarkovWeight{F}_{p_1,q}$};
						\node at (2.5,1.5) {$\dots$};
						\node at (3.5,1.5) {$\MarkovWeight{F}_{p_N,q}$};
						\node at (2,2.5) {$T_W$};
					\end{tikzpicture}
				\right)
				=\sigma^W_q\lambda_{Pq} 
				\begin{tikzpicture}[guillpart,yscale=1.25,xscale=1.25]
					\fill[guillfill] (0,0) rectangle (4,2);
					\draw[guillsep] (4,0)--(4,2)--(0,2) (0,1)--(4,1) (1,0)--(1,1) (2,0)--(2,1) (3,0)--(3,1);
					\node at (0.5,0.5) {$c$};
					\node at (1.5,0.5) {$v_1$};
					\node at (2.5,0.5) {$\dots$};
					\node at (3.5,0.5) {$v_N$};
					\node at (2,1.5) {$T_W$};
				\end{tikzpicture}
			\end{equation}
			with $P=\sum_{i=1}^N p_i$.
		\end{enumerate}
\end{defi}
In particular, for $N=0$, the sequences $(p_i)$ and $(v_i)$ are empty and equation~\eqref{eq:cornermorphism:def} becomes
\[
	\Phi_{\infty_W}^{S,q,r}\left(
	\begin{tikzpicture}[guillpart,yscale=1.25,xscale=1.5]
		\fill[guillfill] (0,0) rectangle (1,3);
		\draw[guillsep] (1,0)--(1,3)--(0,3) (0,1)--(1,1) (0,2)--(1,2) ;
		\node at (0.5,0.5) {$c$};
		\node at (0.5,1.5) {$M^W_q$};
		\node at (0.5,2.5) {$S$};
	\end{tikzpicture}
	\right)
	=\sigma^W_q
	\begin{tikzpicture}[guillpart,yscale=1.25,xscale=1.25]
		\fill[guillfill] (0,0) rectangle (1,2);
		\draw[guillsep] (1,0)--(1,2)--(0,2) (0,1)--(1,1);
		\node at (0.5,0.5) {$c$};
		\node at (0.5,1.5) {$S$};
	\end{tikzpicture}
\]
and it emphasizes that the eigenvalue sequence $\sigma^W_\bullet$ is associated to the vertical semi-group $M^W_\bullet$, in the same way as $(\lambda_{pq})_{q\in\setL^*}$ is associated to the vertical semi-group $\MarkovWeight{F}_{p,\bullet}$.

\begin{prop}\label{prop:shiftingWestgenerators}
	Let $\ca{F}_{\bullet,\bullet}$ be a 2D-semi-group of $\ca{A}_{\PatternShapes(\patterntype{r})}$ and $M^W_\bullet$ a vertical semi-group of $\ca{A}_{\infty_W,\bullet}$. Let $(\ca{V}_{p})_{p\in\setL^*\cup\{\infty_W\}}$ be a system of eigen-$\Guill_1^L$-generators up to morphisms of $(\MarkovWeight{F}_{\bullet,\bullet},M^W_\bullet)$ with eigenvalue geometric sequences $(\lambda_r)_{r\in\setL}$ and $(\sigma^W_q)_{q\in\setL}$. Let $p_0\in\setL^*$ be fixed. Then, the collection of spaces
	\begin{align}
		\ca{V}^{(p_0)}_{\infty_W} &= \left\{
			\begin{tikzpicture}[guillpart]
				\fill[guillfill] (0,0) rectangle (2,1);
				\draw[guillsep] (2,0)--(2,1)--(0,1) (1,0)--(1,1);
				\node at (0.5,0.5) {$c$};
				\node at (1.5,0.5) {$v$};
			\end{tikzpicture}
		; c\in\ca{V}_{\infty_W}, p\in\ca{V}_{p_0}
		\right\}
		\label{eq:shiftcornerspace}
		\\
		\ca{V}^{(p_0)}_{p} &= \ca{V}_p
	\end{align}
 	is a system of eigen-$\Guill_1^L$-generators up to morphisms of $\MarkovWeight{F}_{\bullet,\bullet}$ and the vertical semi-groups
 	\[ M^{W,(p_0)}_q=
	 	\begin{tikzpicture}[guillpart,yscale=1.5,xscale=1.5]
 		\fill[guillfill] (0,0) rectangle (2,1);
 		\draw[guillsep] (0,0)--(2,0)--(2,1)--(0,1) (1,0)--(1,1);
 		\node at (0.5,0.5) {$M_q^{W}$};
 		\node at (1.5,0.5) {$\MarkovWeight{F}_{p,q}$};
 		\end{tikzpicture}
 	\]
 	 with eigenvalue geometric sequences $(\lambda_r)_{r\in\setL^*}$ and $(\sigma^{W,(p_0)}_q)_{q\in\setL^*}$ with $\sigma^{W,(p_0)}_q=\sigma^W_q \lambda_{p_0 q}$ and the same morphisms.
\end{prop}
\begin{proof}
	This is a direct consequence of the previous definition by considering the equation \eqref{eq:cornermorphism:def} for $N\geq 1$ and split $N=1+N'$ with $N'\geq 0$. The first elements $v_1$ with $p_1=p_0$ are then used to form the new space $\ca{V}^{(p_0)}_{\infty_W}$.
\end{proof}
This result illustrates the link between the various vertical $\Guill_1$-algebras $\ca{M}^{W,(p_0)}_\bullet$ defined  in \eqref{eq:eigencorner:westmodule:generators}, obtained by the action from the East of the 2D-semi-group: the eigenvalue  sequence $\sigma^W_\bullet$ associated to the elements $M^W_\bullet$ are shifted by the eigenvalue sequence $\lambda_{p_0\bullet}$ but the morphisms remain the same for all of them.

\subsection{The case of a left vertical semi-group given as a system of West eigen-generators.}

We may now assume that the vertical semi-group $M^W_\bullet$ is not any semi-group but is related to the 2D-semi-group $\MarkovWeight{F}_{\bullet,\bullet}$. The simplest idea is to recognize in definition~\eqref{eq:eigencorner:westmodule:generators} the structure of a possible system of West eigen-generators.

We now assume that the spaces $\ca{M}_\bullet^{W,(0)}$ form a system of West eigen-generators with eigenvalue sequence $(\lambda'_r)_{r\in\setL}$ up to morphisms $\Phi^{W,p,r}_{q}$, which are linear maps from $\ca{A}_{\infty_W,q}$ to itself and satisfy
\begin{equation}
	\Phi_{q}^{W,p,r}\left(
	\begin{tikzpicture}[guillpart,yscale=1.5,xscale=1.5]
		\fill[guillfill] (0,0) rectangle (3,1);
		\draw[guillsep] (0,0)--(3,0)--(3,1)--(0,1) (1,0)--(1,1) (2,0)--(2,1);
		\node at (0.5,0.5) {$M^W_q$};
		\node at (1.5,0.5) {$\MarkovWeight{F}_{p,q}$};
		\node at (2.5,0.5) {$S$};
	\end{tikzpicture}
	\right)
	=
		\Phi_{q}^{W,p,r}\left(
	\begin{tikzpicture}[guillpart,yscale=1.5,xscale=1.5]
		\fill[guillfill] (0,0) rectangle (3,1);
		\draw[guillsep] (0,0)--(3,0)--(3,1)--(0,1) (2,0)--(2,1) ;
		\node at (1,0.5) {$M^{W,(p_0)}_q$};
		\node at (2.5,0.5) {$S$};
	\end{tikzpicture}
	\right)
	= \lambda'_{pq} 	
	\begin{tikzpicture}[guillpart,yscale=1.5,xscale=1.5]
		\fill[guillfill] (0,0) rectangle (2,1);
		\draw[guillsep] (0,0)--(2,0)--(2,1)--(0,1) (1,0)--(1,1) ;
		\node at (0.5,0.5) {$M^W_q$};
		\node at (1.5,0.5) {$S$};
	\end{tikzpicture}
\end{equation}
for any $p,q\in\setL^*$, $r\in\setL$ and $S\in\ca{A}_{r,q}$ (with the convention that $S$ is absent for $r=0$). This equation corresponds to \eqref{eq:defeigenSouth:withM} degenerated to the case where the spaces of the eigen-generators are one-dimensional and provided by a semi-group.

This provides isomorphisms (up to normalizations) between the various vertical semi-groups $M^{W,(p_0)}_\bullet$ and the natural next step is to impose a compatible relation between the spaces $\ca{V}^{(p_0)}_\bullet$ defined in Proposition~\ref{prop:shiftingWestgenerators} by introducing additional morphisms $\Phi_{\infty_S}^{W,p,r}$. This is the same procedure as for the left extension of system of eigenvectors in Definition~\ref{def:eigenalgebrauptomorph:leftextended} but now on the other side adjacent to the corner. One must be careful however about the fact that 
$\Phi_{\infty_S}^{W,p_0,0}$ should reduce any element of $\ca{V}^{(p_0)}_{\infty_W}$ as in \eqref{eq:shiftcornerspace} to $\ti{\sigma}^S_{p_0}(v) c$. The dependence of $\ti{\sigma}^S_{p_0}(v)$ on $v$ is a trickier question and should be formulated in terms of an operadic generalization of characters in the general case. However, if one considers system of South-eigen-generators given by a horizontal semi-group $M^S_p$, then this dependence is necessarily geometric and one has $\ti{\sigma}^S_{p_0}(M^S_{p_0})=\sigma^S_p$ where $\sigma^S_p$ is a geometric sequence.

\removable{
If one wishes to avoid the complication of the definition of generalized characters in the operadic framework ---we postpone it to a further work---, the easiest corner definitions are formulated for systems of South and East eigen-generators that are \emph{both} also horizontal and vertical semi-groups. This is only a small restriction in most cases of Markov processes where one expects a non-degeneracy of the largest eigenvalue related to the uniqueness of a Gibbs measure and to an eventual generalized Perron-Frobenius theorem. In case of degeneracies like in the low-temperature Ising model at zero magnetic fields, one expects the non-extremal Gibbs measures to be decomposed in the algebraic framework as direct sums of ROPEreps with the same sequences $\sigma^a_p$.
}

\subsection{Eigen-corner elements for adjacent systems of eigen-generators of a given semi-groups}

We now assume that we are given two semi-groups $M^W_\bullet$ and $M^S_\bullet$, one vertical and one horizontal such that the spaces
\begin{align*}
	\ca{M}^W_q &=\setK M^W_q 
	&
	\ca{M}^S_p &= \setK M^S_p
\end{align*}
form two systems of eigen-generators of a 2D-semi-group $\MarkovWeight{F}_{\bullet}$ with eigenvalue sequences $\lambda^W_\bullet$ and $\lambda^S_\bullet$ up to respective morphisms $\Phi_q^{W,p,r}$ and $\Phi_p^{S,q,r}$. We now add a corner space of generators $\ca{M}_{\infty_W,\infty_S}$ such that we have simultaneously:
\begin{itemize}
	\item $(\ca{M}_{\infty_W,\infty_S},\ca{M}^S_\bullet)$ is a system of left-extended $\Guill_1^L$-generators of $(\MarkovWeight{F}_\bullet,M^W_\bullet)$ 
	\item  $(\ca{M}_{\infty_W,\infty_S},\ca{M}^W_\bullet)$ is also a system of left-extended $\Guill_1^L$-generators of $(\MarkovWeight{F}_\bullet,M^S_\bullet)$.
\end{itemize}
In this case, there are various ways to reduce the following element, with $c\in\ca{M}_{\infty_W,\infty_S}$,
\[
\begin{tikzpicture}[guillpart,yscale=1.3,xscale=1.5]
	\fill[guillfill] (0,0) rectangle (3,3);
	\draw[guillsep] (3,0)--(3,3)--(0,3) (0,1)--(3,1) (1,0)--(1,3);
	\node at (0.5,0.5) {$c$};
	\node at (2.,0.5) {$M^S_P$};
	\node at (0.5,2) {$M^W_Q$};
	\node at (2,2) {$\MarkovWeight{F}_{P,Q}$};
\end{tikzpicture}
\]
to an element 
\[
\begin{tikzpicture}[guillpart,yscale=1.3,xscale=1.5]
	\fill[guillfill] (0,0) rectangle (2,2);
	\draw[guillsep] (2,0)--(2,2)--(0,2) (0,1)--(2,1) (1,0)--(1,2);
	\node at (0.5,0.5) {$c$};
	\node at (1.5,0.5) {$M^S_p$};
	\node at (0.5,1.5) {$M^W_q$};
	\node at (1.5,1.5) {$\MarkovWeight{F}_{p,q}$};
\end{tikzpicture}
\]
with $p\leq P$ and $q\leq Q$ by using a sequence of morphisms
\[
\Phi^{S,*,*}_*\circ \Phi^{W,*,*}_*\circ \Phi^{S,*,*}_*\circ\ldots\circ \Phi^{W,*,*}_*
\]
where the stars $*$ have to be replaced by any sequence of consistent indices. For each sequence of reduction morphisms, one obtain a scalar $\sigma^W_{Q-q}\sigma^S_{P-p}\lambda^W_{s_1}\Lambda^S_{s_2}$ 
where the values $(s_1,s_2)$ satisfy $s_1+s_2=PQ-pq$ but depend on the sequence of morphisms: $s_1$ is the area of the faces erased by the morphisms $\Phi^{S,*,*}_*$ and $s_2$ is the area of the faces erased by the morphisms $\Phi^{W,*,*}_*$. 

Such a dependence of the sequence, when generalized to the four corners of a full $\Guill_2^{(\patterntype{fp}_*)}$-algebra associated to a Markov process, leads (see below) to a badly-defined partition function and free energy since such quantities may depend on the way elements are simplified using the morphisms, hence violating the generalized associativities~\eqref{eq:guill2:listassoc} which provide evaluation of global quantities independent of the order of the gluing of the geometric shapes. Hence, it is natural to require $\lambda^W_\bullet=\lambda^S_\bullet$ when considering corner generators and we present the following final definition on corners (restricted to lateral semi-groups) and its immediate corollary about simplification rules on corners.

\begin{defi}[corner system of eigen-semi-groups]\label{def:cornereigensemigroups}
	A South-West corner system of eigen-semi-groups of a 2D-semi-group $\MarkovWeight{F}_\bullet$ up to morphisms with eigenvalue sequences $(\lambda_\bullet,\sigma^W_\bullet,\sigma^S_\bullet)$ is a triplet $(M^W_\bullet,M^S_\bullet,\ca{V}_{\infty_W,\infty_S})$ such that:
	\begin{enumerate}[(i)]
		\item $M^S_\bullet$ and $\ca{V}_{\infty_W,\infty_S}$ define a left-extended system of South-eigen-$\Guill_1^L$-generators of $(\MarkovWeight{F}_\bullet,M^W_\bullet)$ with eigenvalue sequences $(\lambda_\bullet,\sigma^W_\bullet)$ up to morphisms $\Phi^{S,q,r}_p$.
		\item $M^W_\bullet$ and $\ca{V}_{\infty_W,\infty_S}$ define a left-extended system of West-eigen-$\Guill_1^L$-generators of $(\MarkovWeight{F}_\bullet,M^S_\bullet)$ with eigenvalue sequences $(\lambda_\bullet,\sigma^S_\bullet)$ up to morphisms $\Phi^{W,p,r}_q$.
	\end{enumerate}
\end{defi}
\begin{prop}\label{prop:eigencorner:commutmorphisms}
	For all $p_1,q_1\in\setL$, $p_2,q_2\in\setL$, for all $c\in\ca{V}_{\infty_W,\infty_S}$ and all $T\in\ca{A}_{p_2,q_2}$, it holds:
	\begin{align*}
		\Phi_{\infty_S}^{W,p_1,p_2}\circ\Phi_{\infty_W}^{S,q_1,q_2}\left(
		\begin{tikzpicture}[guillpart,yscale=1.25,xscale=1.5]
			\fill[guillfill] (0,0) rectangle (4,4);
			\draw[guillsep] (4,0)--(4,4)--(0,4) (1,0)--(1,4) (2.5,0)--(2.5,4) (0,1)--(4,1) (0,2.5)--(4,2.5);
			\node at (0.5,0.5) {$c$};
			\node at (1.75,0.5) {$M^S_{p_1}$};
			\node at (3.25,0.5) {$M^S_{p_2}$};
			\node at (0.5,1.75) {$M^W_{q_1}$};
			\node at (0.5,3.25) {$M^W_{q_2}$};
			\node at (1.75,1.75) {$\MarkovWeight{F}_{p_1,q_1}$};
			\node at (1.75,3.25) {$\MarkovWeight{F}_{p_1,q_2}$};
			\node at (3.25,1.75) {$\MarkovWeight{F}_{p_2,q_1}$};
			\node at (3.25,3.25) {$T$};
		\end{tikzpicture}
	\right)
		&=	\Phi_{\infty_W}^{S,q_1,q_2} \circ \Phi_{\infty_S}^{W,p_1,p_2}\left(
		\begin{tikzpicture}[guillpart,yscale=1.25,xscale=1.5]
			\fill[guillfill] (0,0) rectangle (4,4);
			\draw[guillsep] (4,0)--(4,4)--(0,4) (1,0)--(1,4) (2.5,0)--(2.5,4) (0,1)--(4,1) (0,2.5)--(4,2.5);
			\node at (0.5,0.5) {$c$};
			\node at (1.75,0.5) {$M^S_{p_1}$};
			\node at (3.25,0.5) {$M^S_{p_2}$};
			\node at (0.5,1.75) {$M^W_{q_1}$};
			\node at (0.5,3.25) {$M^W_{q_2}$};
			\node at (1.75,1.75) {$\MarkovWeight{F}_{p_1,q_1}$};
			\node at (1.75,3.25) {$\MarkovWeight{F}_{p_1,q_2}$};
			\node at (3.25,1.75) {$\MarkovWeight{F}_{p_2,q_1}$};
			\node at (3.25,3.25) {$T$};
		\end{tikzpicture}
	\right)
	\\
		&= \sigma^S_{p_1}\sigma^W_{q_1}\lambda_{p_1q_1+p_1q_2+p_2q_1}
		\begin{tikzpicture}[guillpart,yscale=1.25,xscale=1.5]
			\fill[guillfill] (0,0) rectangle (2.5,2.5);
			\draw[guillsep] (2.5,0)--(2.5,2.5)--(0,2.5) (1,0)--(1,2.5) (0,1)--(2.5,1);
			\node at (0.5,0.5) {$c$};
			\node at (1.75,0.5) {$M^S_{p_2}$};
			\node at (0.5,1.75) {$M^W_{q_2}$};
			\node at (1.75,1.75) {$T$};
		\end{tikzpicture}
	\end{align*}
	The composition of the two morphisms \emph{restricted} to this type of elements will be noted shortly $\Phi^{SW,(p_1,q_1),(p_2,q_2)}_{\infty_W,\infty_S}$.
\end{prop}
\begin{proof}
	The proof is a simple exercise with Definition~\ref{def:eigenalgebrauptomorph:leftextended} applied on the West and on the South. In particular, the fact that the eigenvalue sequence $\lambda_\bullet$ is the same on both side plays a crucial role.
\end{proof}

This lemma provides precisely the simplification rules on the corners required in Section~\ref{par:requiredsimplificationrules} to described infinite volume Gibbs measures out of measures on finite domains with suitable ROPEreps of boundary weights. 
\begin{rema}
Definition~\ref{def:cornereigensemigroups} does not provide any simplification rule through morphisms for elements such as:
\[
\begin{tikzpicture}[guillpart,yscale=1.25,xscale=1.5]
	\fill[guillfill] (0,0) rectangle (4,4);
	\draw[guillsep] (4,0)--(4,4)--(0,4) (1,0)--(1,2.5) (2.5,0)--(2.5,4) (0,1)--(2.5,1) (0,2.5)--(4,2.5);
	\node at (0.5,0.5) {$c$};
	\node at (1.75,0.5) {$M^S_{p_1}$};
	\node at (0.5,1.75) {$M^W_{q_1}$};
	\node at (1.75,1.75) {$\MarkovWeight{F}_{p_1,q_1}$};
	\node at (1.25,3.25) {$T_W$};
	\node at (3.25,1.25) {$T_S$};
	\node at (3.25,3.25) {$\MarkovWeight{T}$};
\end{tikzpicture}
\]
for generic elements $T_W$, $T_S$ and $\MarkovWeight{T}$.
\end{rema}
One must then be careful when dealing with Definition~\ref{def:cornereigensemigroups} about the expressions that can be simplified and the ones that cannot without further assumptions: the zoology of simplifiable expressions and unsimplifiable ones becomes richer and richer as pattern shapes grow to strips, half-planes, etc.

\begin{rema}[absence of computational tools and the Yang-Baxter case]
Except for simple models (the ones of Section~\ref{sec:trivialfactorizedcase}, the Gaussian case in \cite{BodiotSimon} and the six-vertex model in \cite{SimonSixV}), we do not have yet any generic computational tools to describe corner eigen-spaces as in Definition~\ref{def:cornereigensemigroups}. The closest perspective is the extension of Section~\ref{sec:commutuptomorph} to the corner case: in particular, it would be very interesting for exactly-solvable systems to build the corner eigen-spaces from a Yang-Baxter-like equation, with its sprawling connections with different domains of algebra.
\end{rema}

\subsection{Simplification for morphisms in the case of discrete Markov processes}\label{sec:concreteeigencorner}

\subsubsection{A first simplification without any detailed ROPE structure.}
We consider here the precise case of an $\Guill_2$-algebra given by the canonical structure $\ca{T}_{\PatternShapes(\patterntype{fp}^*)}$ of a Markov process and a ROPE $\ca{B}_{\PatternShapes(\patterntype{fp}^*)}$, both assembled into the larger structure $\ca{E}_{\PatternShapes(\patterntype{fp}^*)}$ as introduced in the proof of Theorem~\ref{theo:stability}.

\begin{theo}[reduction of morphisms on corners for canonical structures of Markov processes]\label{theo:eigencorner:realizationMarkov}
	Let $\MarkovWeight{F}_{1,1} \in T_{1,1}(V(S_1),V(S_2))$ be the element associated to a face weight of a discrete-space Markov process with values in $S_1$ and $S_2$. Let $\ca{B}$ be a ROPE. Let $(M_\bullet^W,M_\bullet^S,\ca{V}_{\infty_W,\infty_S})$ be a triplet of elements such that $M_\bullet^W$ (resp $M_\bullet^S$) is a vertical semigroup in $\ca{E}_{\infty_W,\bullet}$ (resp. in $\ca{E}_{\bullet,\infty_S}$) and such that
	\begin{subequations}\label{eq:eigencorner:cano:levelassumptions}
	\begin{align}
		M^W_1 &\in \End(V(S_1))^{\otimes 0} \otimes {V(S_2)^*}^{\otimes 1} \otimes \ca{B}_{\infty_W,1} \simeq V(S_2)^*\otimes \ca{B}_{\infty_W,1}
	 \\
		M^S_1 &\in {V(S_1)^*}^{\otimes 1} \otimes \End(V(S_2))^{\otimes 0} \otimes \ca{B}_{1,\infty_S} \simeq V(S_1)^* \otimes \ca{B}_{1,\infty_S}
	\\
		\ca{V}_{\infty_W,\infty_S} &\subset V(S_1)^{\otimes 0} \otimes \End(V(S_2))^{\otimes 0}\otimes \ca{B}_{\infty_W,\infty_S}  \simeq \ca{B}_{\infty_W,\infty_S}
	\end{align}
	\end{subequations}
	The triplet $(M_\bullet^W,M_\bullet^S,\ca{V}_{\infty_W,\infty_S})$ is a corner system of eigen-semi-groups with eigenvalue sequences $(\lambda_\bullet,\sigma^W_\bullet,\sigma^S_\bullet)$ up to morphisms $\Phi_p^{S,q,r}$ and $\Phi_q^{W,p,r}$ (following Definition~\ref{def:cornereigensemigroups}) if and only if there exists a collection of linear maps
	\begin{align*}
	\varphi_{(y)}^{W,p,r}  :\ca{B}_{\infty_W,\infty_S} &\to \ca{B}_{\infty_W,\infty_S}
	\\
	\varphi_{(z)}^{S,q,r}  :\ca{B}_{\infty_W,\infty_S} &\to \ca{B}_{\infty_W,\infty_S}
	\end{align*}
	with $p,q\in\setN^*$, $r\in\setN$, $y\in S_1^p$, $z\in S_2^q$,
	such that, for all $z\in S_2^q$, for all $t\in\ca{B}_{r,\infty_S}$ and for all $c\in\ca{V}_{\infty_W,\infty_S}\simeq \ca{B}_{\infty_W,\infty_S}$, it holds,
	\begin{subequations}
		\label{eq:eigencorner:simplifiedmorphisms}
	\begin{equation}
		\sum_{(x,y,w)\in S_1^p\times S_1^p\times S_2^q}
		\begin{tikzpicture}[guillpart,yscale=1.75,xscale=1.75]
			\fill[guillfill] (0,0) rectangle (1,1);
			\draw[guillsep] (0,0) rectangle (1,1);
			\node at (0.5,0.5) {$\MarkovWeight{F}_{p,q}$};
			\node at (0.5,0.) [anchor=north] {$x$};
			\node at (0.5,1.) [anchor=south] {$y$};
			\node at (0.,0.5) [anchor=east] {$w$};
			\node at (1,0.5) [anchor=west] {$z$};
				\node at (0.5,1) [inner sep=2pt,diamond, fill] {};
				\node at (0.5,0) [inner sep=2pt,diamond, fill] {};
				\node at (0,0.5) [inner sep=2pt,diamond, fill] {};
				\node at (1,0.5) [inner sep=2pt,diamond, fill] {};
		\end{tikzpicture}
	\varphi_{(y)}^{W,p,r}\left(
		\begin{tikzpicture}[guillpart,yscale=1.5,xscale=1.5]
			\draw[guillsep] (4,0)--(4,1)--(1,1)--(1,3)--(-1,3) (1,0)--(1,3) (1,1)--(4,1) (3,0)--(3,1) (-1,1)--(1,1);
			\draw[guillsep,dashed] (4,0)--(-1,0)--(-1,3);
			\node at (0.,0.5) {$c$};
			\node at (2,0.5) {$M_p^S(x)$};
			\node at (3.5,0.5) {$t$};
			\node at (0.,2) {$M_q^W(w)$};
		\end{tikzpicture}
	\right)
		= 
		\sigma_p^S\lambda_{pq} 
		\begin{tikzpicture}[guillpart,yscale=1.5,xscale=1.5]
			\draw[guillsep] (2,0)--(2,1)--(1,1)--(1,3)--(-1,3) (1,0)--(1,3) (-1,1)--(2,1);
			\draw[guillsep,dashed] (2,0)--(-1,0)--(-1,3);
			\node at (0.,0.5) {$c$};
			\node at (1.5,0.5) {$t$};
			\node at (0.,2) {$M_q^W(z)$};
		\end{tikzpicture}
	\end{equation}
	and, for all $y\in S_1^p$, for all $t\in\ca{B}_{\infty_W,r}$ and for all $c\in\ca{V}_{\infty_W,\infty_S}\simeq \ca{B}_{\infty_W,\infty_S}$, it holds,
	\begin{equation}
		\sum_{(x,w,z)\in S_1^p\times S_2^q \times S_2^q}
		\begin{tikzpicture}[guillpart,yscale=1.75,xscale=1.75]
			\fill[guillfill] (0,0) rectangle (1,1);
			\draw[guillsep] (0,0) rectangle (1,1);
			\node at (0.5,0.5) {$\MarkovWeight{F}_{p,q}$};
			\node at (0.5,0.) [anchor=north] {$x$};
			\node at (0.5,1.) [anchor=south] {$y$};
			\node at (0.,0.5) [anchor=east] {$w$};
			\node at (1,0.5) [anchor=west] {$z$};
			\node at (0.5,1) [inner sep=2pt,diamond, fill] {};
			\node at (0.5,0) [inner sep=2pt,diamond, fill] {};
			\node at (0,0.5) [inner sep=2pt,diamond, fill] {};
			\node at (1,0.5) [inner sep=2pt,diamond, fill] {};
		\end{tikzpicture}
		\varphi_{(z)}^{W,p,r}\left(
		\begin{tikzpicture}[guillpart,yscale=1.5,xscale=1.5]
			\draw[guillsep] (3,0)--(3,1)--(1,1)--(1,4)--(-1,4) (1,0)--(1,1) (1,1)--(3,1) (3,0)--(3,1) (-1,1)--(1,1) (-1,3)--(1,3);
			\draw[guillsep,dashed] (3,0)--(-1,0)--(-1,4);
			\node at (0.,0.5) {$c$};
			\node at (2,0.5) {$M_p^S(x)$};
			\node at (0.,3.5) {$t$};
			\node at (0.,2) {$M_q^W(w)$};
		\end{tikzpicture}
		\right)
		= 
		\sigma_q^W\lambda_{pq} 
		\begin{tikzpicture}[guillpart,yscale=1.5,xscale=1.5]
			\draw[guillsep] (3,0)--(3,1)--(1,1)--(1,2)--(0,2) (1,0)--(1,2) (0,1)--(3,1);
			\draw[guillsep,dashed] (3,0)--(0,0)--(0,2);
			\node at (0.5,0.5) {$c$};
			\node at (0.5,1.5) {$t$};
			\node at (2.,0.5) {$M_p^S(y)$};
		\end{tikzpicture}
	\end{equation}
	\end{subequations}
	with the following conventions, which corresponds to the description of an element of $\ca{E}_{\bullet,\bullet}$ by tensorization with elements of $\ca{B}_{\bullet,\bullet}$:
	\begin{align*}
		M_p^{S} &=  \sum_{y\in S_1^p} e^*_y \otimes M_p^{S}(y) , \qquad M_p^{S}(y)\in\ca{B}_{p,\infty_S}
		\\
		M_q^{W} &=  \sum_{z\in S_2^q} e^*_z \otimes M_q^{W}(z) , \qquad M_q^{W}(z)\in\ca{B}_{\infty_W,q}
	\end{align*}
\end{theo}
Assumptions \eqref{eq:eigencorner:cano:levelassumptions} are satisfied in the the construction of Theorem~\ref{theo:stability}. The first two assumptions correspond to the case $q_0=0$ of Theorem~\ref{theo:eigencorner:realizationMarkov}, which we consider here to make notations easier (but a similar result holds for non-zero $q_0$). 
\begin{proof}
The proof follows closely the proof of Theorem~\ref{theo:reductionofmorphisms:canostruct} and consists in the separate study of the two items of Definition~\ref{def:cornereigensemigroups}: we then only work with one of them ---item~(i) here~--- and let the reader consider the second case in a similar way.

We consider \eqref{eq:cornermorphism:def} for arbitrary elements $T_W$. From the assumptions on $M_q^W$ and $\ca{V}_{\infty_W,\infty_S}$, it is sufficient to consider only the case of arbitrary elements
\[
T_W \in \End(V(S_1))^{\otimes P} \otimes {V(S_2)^*}^{\otimes r}\otimes \ca{B}_{\infty_W,r}.
\]
As in the previous proof, we may first choose $T_W$ equal to
\[
T_W = e_{x}e^*_y \otimes e^*_u \otimes t
\]
for arbitrary elements $x,y\in S_1^P$, $u\in S_2^r$ and $t\in \ca{B}_{\infty_W,r}$. The second step is the description of a generic linear map from the subspaces 
\[
\phi : {V(S_1)^*}^{\otimes P}\otimes {V(S_2)^*}^{\otimes q+r} \otimes \ca{B}_{\infty_W,r} \to {V(S_1)^*}^{\otimes P}\otimes {V(S_2)^*}^{\otimes r} \otimes \ca{B}_{\infty_W,r}
\]
through
\[
\phi( e^*_x\otimes e^*_w\otimes b ) = \sum_{(y,z)} e^*_y\otimes e^*_z\otimes \ti{\phi}_{x,y;w,z}(b)
\]
Combining arbitrary $T_W$ with this parametrization of the morphisms provides the expected result after heavy but easy computations, as in the previous proof.
\end{proof}

The previous theorem essentially relies on the careful study of the "matrix" part $\ca{T}_{\PatternShapes(\patterntype{fp}^*)}$ of the canonical structure $\ca{E}_{\PatternShapes(\patterntype{fp}^*)}$ through the careful choices of elements $T_W$. Any further simplification beyond \eqref{eq:eigencorner:simplifiedmorphisms} requires an insight on the ROPE $\ca{B}_{\PatternShapes(\patterntype{fp}^*)}$ in order to exploit the arbitrary choice of $t$ in $T_W$. From the restricted point of view of the present paper, the choice of the ROPE and of the arbitrary element $t$ in $T_W$ is not obvious at all since we care only about \emph{one} eigen-boundary condition of \emph{one}-semigroup in the rectangular guillotine geometry: we think that the generalization of our constructions may lead to a very large of choices of elements $t \in \ca{B}_{\infty_W,r}$ to achieve full generality. The easiest thing we may imagine is the introduction of arbitrary observables on segments between the faces, which may then replace the semi-group elements by modified values. We only provide here an example of further simplification in order to illustrate how \eqref{eq:eigencorner:simplifiedmorphisms} can be realized.

\subsubsection{Further simplification in the case of a elementary ROPE} As already discussed in Section~\ref{sec:eagerandlazyROPEs}, action of $\Guill_1$-algebras on boundary extensions "removes" the colours and corresponds to a representation of the boundary elements $\pi_W : \ca{B}_{\infty_W,q}\to\Hom(\ca{B}_{\infty_W,\infty_S})$ and $\pi_S:\ca{B}_{p,\infty_S} \to \Hom(\ca{B}_{\infty_W,\infty_S})$. In this case, the first equation~\eqref{eq:eigencorner:simplifiedmorphisms}, corresponds to:
\begin{equation}
	\begin{split}
\sum_{(x,y,w)\in S_1^p\times S_1^p\times S_2^q}&
\begin{tikzpicture}[guillpart,yscale=1.75,xscale=1.75]
	\fill[guillfill] (0,0) rectangle (1,1);
	\draw[guillsep] (0,0) rectangle (1,1);
	\node at (0.5,0.5) {$\MarkovWeight{F}_{p,q}$};
	\node at (0.5,0.) [anchor=north] {$x$};
	\node at (0.5,1.) [anchor=south] {$y$};
	\node at (0.,0.5) [anchor=east] {$w$};
	\node at (1,0.5) [anchor=west] {$z$};
	\node at (0.5,1) [inner sep=2pt,diamond, fill] {};
	\node at (0.5,0) [inner sep=2pt,diamond, fill] {};
	\node at (0,0.5) [inner sep=2pt,diamond, fill] {};
	\node at (1,0.5) [inner sep=2pt,diamond, fill] {};
\end{tikzpicture}
\varphi_{(y)}^{W,p,r}\left(
c \left( \pi_W(M_q^W(w))\otimes \pi_S(M_p^S(x))\pi_S(t) \right)
\right)
\\
&= 
\sigma_p^S\lambda_{pq} 
c \left( \pi_W(M_q^W(z)) \otimes \pi_S(t) \right)
\end{split}
\end{equation}
which holds for $c\in\ca{V}_{\infty_W,\infty_S}$ and is interpreted as an action of $\Hom(\ca{B}_{\infty_W,\infty_S},\ca{B}_{\infty_W,\infty_S})$ (the tensor product is present 
only to remind reader the bimodule structures of the corner space of a ROPE).

If one further assumes that the only elements which commutes with $\pi_S(t)$ are the scalar constant $\alpha \id$ (as a consequence some Schur's lemma for example) and that, $u\pi_S(t)=0$ for all $t$ implies $u=0$ , then the linear maps $\varphi_{(z)}^{W,p,r}$ necessarily acts from the West as some elements $K_{(z)}^{W,p,r}\in \Hom(\ca{B}_{\infty_W,\infty_S},\ca{B}_{\infty_W,\infty_S})$ and one obtains 
\begin{equation}\label{eq:concreteeigencorner}
	\begin{split}
	\sum_{(x,y,w)\in S_1^p\times S_1^p\times S_2^q}&
	\begin{tikzpicture}[guillpart,yscale=1.75,xscale=1.75]
		\fill[guillfill] (0,0) rectangle (1,1);
		\draw[guillsep] (0,0) rectangle (1,1);
		\node at (0.5,0.5) {$\MarkovWeight{F}_{p,q}$};
		\node at (0.5,0.) [anchor=north] {$x$};
		\node at (0.5,1.) [anchor=south] {$y$};
		\node at (0.,0.5) [anchor=east] {$w$};
		\node at (1,0.5) [anchor=west] {$z$};
		\node at (0.5,1) [inner sep=2pt,diamond, fill] {};
		\node at (0.5,0) [inner sep=2pt,diamond, fill] {};
		\node at (0,0.5) [inner sep=2pt,diamond, fill] {};
		\node at (1,0.5) [inner sep=2pt,diamond, fill] {};
	\end{tikzpicture}
	c \left( \pi_W(M_q^W(w))K_{(y)}^{W,p,r} \otimes \pi_S(M_p^S(x)) \right)
	\\
	&= 
	\sigma_p^S\lambda_{pq} 
	c \left( \pi_W(M_q^W(z))  \right)
	\end{split}
\end{equation}
This equation, together with its South equivalent, takes place as linear maps on $\ca{B}_{\infty_W,\infty_S}$ and it is the equation to solve in practice, once $\ca{B}_{\infty_W,\infty_S}$ is identified: the unknown variables are the eigenvalues $\sigma^S_\bullet$ (which can be absorbed as before in a redefinition of the elements $K_{(y)}^{W,p,r}$ or $M_p^{S}(w)$)), the elements $c$ and the linear maps $K_{(y)}^{W,p,r}$, which are the new ingredients related to the definition up to morphisms.

Such detailed computations for concrete models will be presented in the companion papers \cite{BodiotSimon,SimonSixV}.

\section{From one-sided infinite shapes to two-sided infinite shapes}\label{sec:pairingeigentostrips}

This section extends the previous definitions to the full plane geometry and, in particular, definitions must be adapted to pointed doubly-infinite guillotine partitions. Since it does not correspond to the basic computational building blocks of invariant ROPEreps, it can be skipped in a first reading. This is however necessary to glue together opposite eigen-elements in the generic setting. 

\subsection{Preliminary remarks on the geometry and the operadic structures}

We forget temporarily in this paragraph the previous corner case and considers a $\Guill_2^{(\patterntype{s}_{SN})}$-algebra $(\ca{A}_{p,q})_{(p,q)\in \PatternShapes(\patterntype{s}_{SN})}$, with $\PatternShapes(\patterntype{s}_{SN})=\setL^* \times (\setL^*\cup\{\infty_S,\infty_N,\infty_{SN}\})$. As seen before, the products 
\[
\begin{tikzpicture}[guillpart]
	\fill[guillfill] (0,0) rectangle (1,3);
	\draw[guillsep] (0,0)--(0,3) (1,0)--(1,3) (0,1)--(1,1) ;
	\node at (0.5,0.5) {$1$};
	\node at (0.5,1.8) {$2$};
	\node at (0,2.4) [circle, fill, inner sep=0.5mm] {};
	\node at (1,2.4) [circle, fill, inner sep=0.5mm] {};
	\draw [dotted] (0,2.4) -- (1,2.4);
	\draw [->] (1.3,1) -- node [midway, right] {$y$} (1.3,2.4);	
\end{tikzpicture}
: \ca{A}_{p,\infty_S} \otimes \ca{A}_{p,\infty_N} \to \ca{A}_{p,\infty_{SN}}
\]
contain a base point information in order to break translational invariance. Furthermore the spaces $(\ca{A}_{p,\infty_{SN}})_{p\in\setL^*}$ have a $\Guill_1$-algebra structure.

Following the same approach as in Section~\ref{sec:addingcornertoeigenSouth}, we first consider a North horizontal semi-group  $(M^N_p)_{p \in\setL^*}$ with $M^N_p\in \ca{A}_{p,\infty_N}$ for all $p\in\setL^*$ such that
\begin{equation}\label{eq:semigroup:north}
	\begin{tikzpicture}[guillpart,scale=1.5]
		\fill[guillfill] (0,0) rectangle (2,1);
		\draw[guillsep] (0,1)--(0,0)--(2,0)--(2,1) (1,0)--(1,1);
		\node at (0.5,0.5) {$M^N_{p_1}$};
		\node at (1.5,0.5) {$M^N_{p_2}$};
	\end{tikzpicture}
=M^N_{p_1+p_2}
\end{equation}
for all $p_1,p_2\in\setL^*$. 

We then define the following linear maps
\begin{align*}
\ti{\psi}_p^{N,M^N,q} :  \ca{A}_{p,q} & \to \ca{A}_{p,\infty_N}  
\\
\MarkovWeight{T} 
& \mapsto 
\begin{tikzpicture}[guillpart,yscale=1.,xscale=1.5]
	\fill[guillfill] (0,0) rectangle (1,2);
	\draw[guillsep] (0,2)--(0,0)--(1,0)--(1,2) (0,1)--(1,1);
	\node at (0.5,0.5) {$\MarkovWeight{T}$};
	\node at (0.5,1.5) {$M^N_p$};
\end{tikzpicture}
\\
\ti{\psi}_p^{N,M^N,\infty_S,(y)} : \ca{A}_{p,\infty_S} & \to \ca{A}_{p,\infty_{SN}}
\\
A &\mapsto
 \begin{tikzpicture}[guillpart,xscale=1.5]
	\fill[guillfill] (0,0) rectangle (1,3);
	\draw[guillsep] (0,0)--(0,3) (1,0)--(1,3) (0,1)--(1,1) ;
	\node at (0.5,0.5) {$A$};
	\node at (0.5,1.8) {$M^N_p$};
	\node at (0,2.4) [circle, fill, inner sep=0.5mm] {};
	\node at (1,2.4) [circle, fill, inner sep=0.5mm] {};
	\draw [dotted] (0,2.4) -- (1,2.4);
	\draw [->] (1.3,1) -- node [midway, right] {$y$} (1.3,2.4);	
\end{tikzpicture}
\end{align*}
which, using square associativity~\eqref{eq:guill2:interchangeassoc}, are morphisms of horizontal $\Guill_1$-algebra. There are additional suitable vertical compatibilities when one uses a 2D-semi-group $\MarkovWeight{F}_\bullet$:
\begin{equation}\label{eq:strip:algebramorph}
	\psi_p^{N,M^N,\infty_S,y} \circ \psi^{S,\MarkovWeight{F},q,0}(A) 
	= 
	\begin{tikzpicture}[guillpart,yscale=0.8,xscale=1.5]
		\fill[guillfill] (0,-1) rectangle (1,3);
		\draw[guillsep] (0,-1)--(0,3) (1,-1)--(1,3) (0,1)--(1,1) (0,0)--(1,0);
		\node at (0.5,-0.5) {$A$};
		\node at (0.5,0.5) {$\MarkovWeight{F}_{p,q}$};
		\node at (0.5,1.8) {$M^N_p$};
		\node at (0,2.4) [circle, fill, inner sep=0.5mm] {};
		\node at (1,2.4) [circle, fill, inner sep=0.5mm] {};
		\draw [dotted] (0,2.4) -- (1,2.4);
		\draw [->] (1.3,1) -- node [midway, right] {$y$} (1.3,2.4);	
	\end{tikzpicture}
	= 
	 \begin{tikzpicture}[guillpart,yscale=0.8,xscale=2.3]
		\fill[guillfill] (0,0) rectangle (1,3);
		\draw[guillsep] (0,0)--(0,3) (1,0)--(1,3) (0,1)--(1,1) ;
		\node at (0.5,0.5) {$A$};
		\node at (0.5,1.8) {$M^{N,(q)}_p$};
		\node at (0,2.4) [circle, fill, inner sep=0.5mm] {};
		\node at (1,2.4) [circle, fill, inner sep=0.5mm] {};
		\draw [dotted] (0,2.4) -- (1,2.4);
		\draw [->] (1.3,1) -- node [midway, right] {$y+q$} (1.3,2.4);	
	\end{tikzpicture}
	=
	\psi_p^{N,M^{N,(q)},\infty_S,y+q}(A)
\end{equation}
where the horizontal semi-group $M_\bullet^{N,(q)}$ is defined as
\[
M^{N,(q)}_p = \begin{tikzpicture}[guillpart,yscale=1.5,xscale=1.5]
	\fill[guillfill] (0,0) rectangle (1,2);
	\draw[guillsep] (0,2)--(0,0)--(1,0)--(1,2) (0,1)--(1,1);
	\node at (0.5,0.5) {$\MarkovWeight{F}_{p,q}$};
	\node at (0.5,1.5) {$M^N_p$};
\end{tikzpicture}
\]

\subsection{Action of North and South eigen-elements: simplification rules}

\subsubsection{Action of arbitrary North elements on South eigen-generators.}
We now consider a system $\ca{V}^S_\bullet$ of South-eigen-generators up to morphisms of a 2D-semi-group $\MarkovWeight{F}_\bullet$ and a North semi-group $M^N_\bullet$. 

We recall from Section~\ref{sec:generators:remark} that such a system $\ca{V}^S_\bullet$ generates a tower of isomorphic (up to scalar eigenvalues) sub-$\Guill_1$-algebras $\ca{S}^{(q)}_\bullet$ of $\ca{A}_{\bullet,\infty_S}$ with morphisms $\psi^{S,\MarkovWeight{F},\infty_S,q'-q}$ for $q\leq q'$ and reciprocal morphism $\Phi^{S,q'-q,0}$. Using the North semi-group $M^N_\bullet$ and the morphisms $\ti{\psi}^{N,M^N,\infty_S,(y)}$, we now have a collection of sub-$\Guill_1$-algebra $\ti{\ca{S}}_\bullet^{(q,y)}$ of $\ca{A}_{\bullet,\infty_{SN}}$ generated by elements \eqref{eq:strip:algebramorph} with $A\in\ca{S}_\bullet$. It may not be obvious a priori to choose which type of natural morphisms to introduce on these new sub-algebras of $\ca{A}_{\bullet,\infty_{SN}}$ since no particular object acts on them.

The solution lies again is the simplification rules required in Section~\ref{par:requiredsimplificationrules}. The idea is to have a simplification rule of any element of the type:
\begin{equation}\label{eq:simplifcationpattern:strip}
\begin{tikzpicture}[guillpart,yscale=1.1,xscale=1.5]
	\fill[guillfill] (0,0) rectangle (1,4);
	\draw[guillsep] (0,0)--(0,4) (1,0)--(1,4) (0,1)--(1,1) (0,2)--(1,2) ;
	\node at (0.5,0.5) {$A$};
	\node at (0.5,1.5) {$\MarkovWeight{F}_{p,q}$};
	\node at (0.5,2.8) {$\MarkovWeight{T}'$};	
	\node at (0,3.4) [circle, fill, inner sep=0.5mm] {};
	\node at (1,3.4) [circle, fill, inner sep=0.5mm] {};
	\draw [dotted] (0,3.4) -- (1,3.4);
	\draw [->] (1.3,2) -- node [midway, right] {$y$} (1.3,3.4);
\end{tikzpicture}
\longmapsto
\lambda^S_{pq}
\begin{tikzpicture}[guillpart,yscale=1.1,xscale=1.5]
	\fill[guillfill] (0,0) rectangle (1,3);
	\draw[guillsep] (0,0)--(0,3) (1,0)--(1,3) (0,1)--(1,1)  ;
	\node at (0.5,0.5) {$A$};
	\node at (0.5,1.8) {$\MarkovWeight{T}'$};
	\node at (0,2.4) [circle, fill, inner sep=0.5mm] {};
	\node at (1,2.4) [circle, fill, inner sep=0.5mm] {};
	\draw [dotted] (0,2.4) -- (1,2.4);
	\draw [->] (1.3,1) -- node [midway, right] {$y$} (1.3,2.4);
\end{tikzpicture}
\end{equation}
Due to the presence of the arbitrary element $\MarkovWeight{T}'\in\ca{A}_{p,\infty_N}$, the semi-group property of $M^N_\bullet$ does not play any particular role. We have instead an extension of Definition~\ref{def:eigenalgebrauptomorphims} where the constraint~\eqref{eq:defeigenSouth:withM} has to be extended to the case $r = \infty_N$. This case was obviously absent in this definition since we were considering only a $\Guill_2^{(\patterntype{hs}_S)}$-algebra. 

When considering a $\Guill_2^{(\patterntype{s}_{SN}^*)}$-algebra instead of a $\Guill_2^{(\patterntype{hs}_S)}$-algebra, equations~\eqref{eq:defeigenSouth:withM} (and subsequent ones), equations~\eqref{eq:def:commutuptmorph:single} and \eqref{eq:commutuptomorph} (and subsequent ones) have to be modified so that colours $r\in\setL$ of the North elements in the products are allowed to be equal to $\infty_N$. In this case, an additional base point have to be introduced. We then write $r\in\setL \cup\{\infty_N^*\}$ which corresponds to either $r\in\setL$ or $r=(\infty_N,y)$ where $y$ is a base point in the doubly-infinite direction.

\begin{defi}[completed system of eigen-generators]\label{def:extensioneigenalgtostrip}
	Using the same notations as definitions~\ref{def:eigenalgebrauptomorphims} with now a $\Guill_2^{(\patterntype{s}_{SN}^*)}$-algebra $\ca{A}_\bullet$ and \ref{def:eigenalgebrauptomorph:leftextended} with now a $\Guill_2^{(\patterntype{hp}_{W,SN}^*)}$, we define \emph{completed} systems of eigen-generators in both cases by assuming that equations~\eqref{eq:defeigenSouth:withM} and \eqref{eq:cornermorphism:def}, which are valid for $r\in\setL^*$, are \emph{also} valid for any $r=(\infty_N,y)$ where $y\in\setP$ is the base point used in the guillotine partition of vertical strip and half-planes, with additional \emph{external} morphisms
	$\Phi_{p}^{S,q,(\infty_N,y)} : \ca{A}_{p,\infty_{SN}} \to \ca{A}_{p,\infty_{SN}}$
	for $p\in\setL^*$ or $p\in\setL^*\cup\{\infty_W\}$.
\end{defi}

These additional morphisms play a different role than the previous one since they do not act of the same spaces: the maps $\Phi_p^{S,q,r)}$ with finite $r\in\setL^*$ are all linear maps on the same space $\ca{V}^S_p\subset\ca{A}_{p,\infty_S}$ and ignore the North spaces $\ca{A}_{p,\infty_N}$. The maps $\Phi_{p}^{S,q,(\infty_N,y)}$ take into account these opposite spaces through the pairing to the spaces $\ca{A}_{p,\infty_{SN}}$; however, Definition~\ref{def:extensioneigenalgtostrip} precisely state that simplification rules are not modified by this extra information.

\subsubsection{Combining North and South eigen-generators.}

We now assume that we are given two \emph{completed} systems $\ca{V}_\bullet^N$ and $\ca{V}_\bullet^S$, one on the North and one on the South, of eigen-generators of a 2D-semi-group $\MarkovWeight{F}_\bullet$ up to morphisms $\Phi_\bullet^{N,q,r}$ and $\Phi_\bullet^{S,q,r}$ with eigenvalue sequences $\lambda^N_\bullet$ and $\lambda^S_\bullet$ (following Definition~\ref{def:eigenalgebrauptomorphims} amended by Definition~\ref{def:extensioneigenalgtostrip}). This introduces simplification rules on elements:
\[
	\begin{tikzpicture}[guillpart,yscale=1.1,xscale=1.8]
		\fill[guillfill] (0,0) rectangle (1,5);
		\draw[guillsep] (0,0)--(0,5) (1,0)--(1,5) (0,1)--(1,1) (0,2)--(1,2) (0,3)--(1,3) (0,4)--(1,4);
		\node at (0.5,0.5) {$v$};
		\node at (0.5,1.5) {$\MarkovWeight{F}_{p,q_S}$};
		\node at (0.5,2.5) {$\MarkovWeight{T}$};
		\node at (0.5,3.5) {$\MarkovWeight{F}_{p,q_N}$};
		\node at (0.5,4.5) {$v'$};
	\end{tikzpicture}
\longmapsto \lambda^S_{p(q_S-q'_S)} \lambda^{N}_{p(q_N-q'_N)} 	
\begin{tikzpicture}[guillpart,yscale=1.1,xscale=1.8]
	\fill[guillfill] (0,0) rectangle (1,5);
	\draw[guillsep] (0,0)--(0,5) (1,0)--(1,5) (0,1)--(1,1) (0,2)--(1,2) (0,3)--(1,3) (0,4)--(1,4);
	\node at (0.5,0.5) {$v$};
	\node at (0.5,1.5) {$\MarkovWeight{F}_{p,q'_S}$};
	\node at (0.5,2.5) {$\MarkovWeight{T}$};
	\node at (0.5,3.5) {$\MarkovWeight{F}_{p,q'_N}$};
	\node at (0.5,4.5) {$v'$};
\end{tikzpicture}
\]
(base point skipped) with $q'_S\leq q_S$ and $q'_N\leq q_N$ both on the North and on the South using the two types of morphisms $\Phi^N$ and $\Phi^S$, which commute when correct choices of base points are made.

However, if $\MarkovWeight{T}=\MarkovWeight{F}_{p,r}$ for some $r\in\setL^*$, then there are various choices of sequence of simplification morphisms $\Phi^N$ and $\Phi^S$ depending on how the middle element is split between a North and South part using the vertical semi-group property of $\MarkovWeight{F}_\bullet$: in a similar manner as for the corner case, depending on the sequence of morphisms, all final expressions differ only by a multiplicative constant if $\lambda^S_\bullet\neq \lambda^N_\bullet$. We must also correctly take into account the presence of a base point.

\begin{prop}\label{prop:simplificationonstrips}
	Let $\ca{A}_{\PatternShapes(\patterntype{s}_{SN}^*)}$ be a $\Guill_2^{(\patterntype{s}_{SN}^*)}$-algebra and let $\MarkovWeight{F}$ be a $2D$-semigroup.	Let $\ca{V}^N_\bullet$ and $\ca{V}_\bullet^S$ be completed systems of North and South eigen-generators of $\MarkovWeight{F}_\bullet$ up to morphisms $\Phi_\bullet^{N,q,r}$ and $\Phi_\bullet^{N,q,(\infty_S,y)}$ and $\Phi_\bullet^{S,q,r}$ and $\Phi_\bullet^{S,q,(\infty_N,y)}$ with the \emph{same} eigenvalue sequence $\lambda_\bullet$ such that, for any $A_S\in\ca{V}_p^S$, $A_N\in\ca{V}_p^N$, $T_S\in\ca{A}_{p,\infty_S}$, $T_N\in\ca{A}_{p,\infty_N}$
	\begin{align}
		\label{eq:strip:NorthandSouthmorph}
		\Phi_\bullet^{S,q,(\infty_N,y)}\left(
		\begin{tikzpicture}[guillpart,yscale=1.1,xscale=1.5]
			\fill[guillfill] (0,0) rectangle (1,4);
			\draw[guillsep] (0,0)--(0,4) (1,0)--(1,4) (0,1)--(1,1) (0,2)--(1,2) ;
			\node at (0.5,0.5) {$A_S$};
			\node at (0.5,1.5) {$\MarkovWeight{F}_{p,q}$};
			\node at (0.5,2.8) {$\MarkovWeight{T}_N$};	
			\node at (0,3.4) [circle, fill, inner sep=0.5mm] {};
			\node at (1,3.4) [circle, fill, inner sep=0.5mm] {};
			\draw [dotted] (0,3.4) -- (1,3.4);
			\draw [->] (1.3,2) -- node [midway, right] {$y$} (1.3,3.4);
		\end{tikzpicture}\right)
		&= \lambda_{pq}
		\begin{tikzpicture}[guillpart,yscale=1.1,xscale=1.5]
			\fill[guillfill] (0,0) rectangle (1,3);
			\draw[guillsep] (0,0)--(0,3) (1,0)--(1,3) (0,1)--(1,1)  ;
			\node at (0.5,0.5) {$A$};
			\node at (0.5,1.8) {$\MarkovWeight{T}_N$};
			\node at (0,2.4) [circle, fill, inner sep=0.5mm] {};
			\node at (1,2.4) [circle, fill, inner sep=0.5mm] {};
			\draw [dotted] (0,2.4) -- (1,2.4);
			\draw [->] (1.3,1) -- node [midway, right] {$y$} (1.3,2.4);
		\end{tikzpicture}
	&
		\Phi_\bullet^{N,q',(\infty_N,z)}\left(
		\begin{tikzpicture}[guillpart,yscale=-1.1,xscale=1.5,]
			\fill[guillfill] (0,0) rectangle (1,4);
			\draw[guillsep] (0,0)--(0,4) (1,0)--(1,4) (0,1)--(1,1) (0,2)--(1,2) ;
			\node at (0.5,0.5) {$A_N$};
			\node at (0.5,1.5) {$\MarkovWeight{F}_{p,q'}$};
			\node at (0.5,2.8) {$\MarkovWeight{T}_S$};	
			\node at (0,3.4) [circle, fill, inner sep=0.5mm] {};
			\node at (1,3.4) [circle, fill, inner sep=0.5mm] {};
			\draw [dotted] (0,3.4) -- (1,3.4);
			\draw [->] (1.3,2) -- node [midway, right] {$z$} (1.3,3.4);
		\end{tikzpicture}\right)
		&= \lambda_{pq'}
		\begin{tikzpicture}[guillpart,yscale=-1.1,xscale=1.5,]
			\fill[guillfill] (0,0) rectangle (1,3);
			\draw[guillsep] (0,0)--(0,3) (1,0)--(1,3) (0,1)--(1,1)  ;
			\node at (0.5,0.5) {$A_N$};
			\node at (0.5,1.8) {$\MarkovWeight{T}_S$};
			\node at (0,2.4) [circle, fill, inner sep=0.5mm] {};
			\node at (1,2.4) [circle, fill, inner sep=0.5mm] {};
			\draw [dotted] (0,2.4) -- (1,2.4);
			\draw [->] (1.3,1) -- node [midway, right] {$z$} (1.3,2.4);
		\end{tikzpicture}	
	\end{align}
	The following relation holds for any $p,q,q',r\in\setL^*$, $y\in\setP$, $T\in\ca{A}_{p,r}$ and $A_a\in\ca{V}^a_p$, $a\in\{S,N\}$:
	\begin{equation}\begin{split}
		\Phi_{p}^{N,q',(\infty_S,y-r)}\circ\Phi_{p}^{S,q,(\infty_N,y)}\left(
		\begin{tikzpicture}[guillpart,yscale=1.15,xscale=1.8]
			\fill[guillfill] (0,0) rectangle (1,5);
			\draw[guillsep] (0,0)--(0,5) (1,0)--(1,5) (0,1)--(1,1) (0,2)--(1,2) (0,3)--(1,3) (0,4)--(1,4);
			\node at (0.5,0.5) {$A_S$};
			\node at (0.5,1.5) {$\MarkovWeight{F}_{p,q}$};
			\node at (0.5,2.5) {$\MarkovWeight{T}$};
			\node at (0.5,3.5) {$\MarkovWeight{F}_{p,q'}$};
			\node at (0.5,4.5) {$A_N$};
			\node at (0,3.4) [circle, fill, inner sep=0.5mm] {};
			\node at (1,3.4) [circle, fill, inner sep=0.5mm] {};
			\draw [dotted] (0,3.4) -- (1,3.4);
			\draw [->] (1.3,2) -- node [midway, right] {$y$} (1.3,3.4);
		\end{tikzpicture}\right)
		=&
		\Phi_{p}^{S,q,(\infty_N,y)}\circ\Phi_{p}^{N,q',(\infty_S,y-r)}\left(
		\begin{tikzpicture}[guillpart,yscale=1.15,xscale=1.8]
			\fill[guillfill] (0,0) rectangle (1,5);
			\draw[guillsep] (0,0)--(0,5) (1,0)--(1,5) (0,1)--(1,1) (0,2)--(1,2) (0,3)--(1,3) (0,4)--(1,4);
			\node at (0.5,0.5) {$A_S$};
			\node at (0.5,1.5) {$\MarkovWeight{F}_{p,q}$};
			\node at (0.5,2.5) {$\MarkovWeight{T}$};
			\node at (0.5,3.5) {$\MarkovWeight{F}_{p,q'}$};
			\node at (0.5,4.5) {$A_N$};
			\node at (0,3.4) [circle, fill, inner sep=0.5mm] {};
			\node at (1,3.4) [circle, fill, inner sep=0.5mm] {};
			\draw [dotted] (0,3.4) -- (1,3.4);
			\draw [->] (1.3,2) -- node [midway, right] {$y$} (1.3,3.4);
		\end{tikzpicture}
		\right)
		\\
		=& \lambda_{p(q+q')}
		\begin{tikzpicture}[guillpart,yscale=1.15,xscale=1.8]
			\fill[guillfill] (0,1) rectangle (1,4);
			\draw[guillsep] (0,1)--(0,4) (1,1)--(1,4)  (0,2)--(1,2) (0,3)--(1,3);
			\node at (0.5,1.5) {$A_S$};
			\node at (0.5,2.5) {$\MarkovWeight{T}$};
			\node at (0.5,3.5) {$A_N$};
			\node at (0,3.4) [circle, fill, inner sep=0.5mm] {};
			\node at (1,3.4) [circle, fill, inner sep=0.5mm] {};
			\draw [dotted] (0,3.4) -- (1,3.4);
			\draw [->] (1.3,2) -- node [midway, right] {$y$} (1.3,3.4);
		\end{tikzpicture}
	\end{split}
	\end{equation}
	The composition of the two commuting morphisms restricted to the North and South spaces of eigen-generators is then written $\Phi_{p}^{SN,q,q',r,(y)}$.
\end{prop}
\begin{proof}
	The only technical point is the presence of the base point. Using repeatedly Proposition~\ref{prop:stripalgebra} and the commutative diagram~\eqref{eq:commutdiag:pairing} shifts the base point used for the pairing so that the elements $\MarkovWeight{T}_N$ and $\MarkovWeight{T}_S$ can be chosen arbitrarily to apply formulae~\eqref{eq:strip:NorthandSouthmorph}. We thus have:
	\[
	\begin{tikzpicture}[guillpart,yscale=1.1,xscale=1.8]
		\fill[guillfill] (0,0) rectangle (1,5);
		\draw[guillsep] (0,0)--(0,5) (1,0)--(1,5) (0,1)--(1,1) (0,2)--(1,2) (0,3)--(1,3) (0,4)--(1,4);
		\node at (0.5,0.5) {$A_S$};
		\node at (0.5,1.5) {$\MarkovWeight{F}_{p,q}$};
		\node at (0.5,2.5) {$\MarkovWeight{T}$};
		\node at (0.5,3.5) {$\MarkovWeight{F}_{p,q'}$};
		\node at (0.5,4.5) {$A_N$};
		\node at (0,3.4) [circle, fill, inner sep=0.5mm] {};
		\node at (1,3.4) [circle, fill, inner sep=0.5mm] {};
		\draw [dotted] (0,3.4) -- (1,3.4);
		\draw [->] (1.3,2) -- node [midway, right] {$y$} (1.3,3.4);
	\end{tikzpicture}
	=
	\begin{tikzpicture}[guillpart,yscale=1.1,xscale=1.8]
		\fill[guillfill] (0,0) rectangle (1,5);
		\draw[guillsep] (0,0)--(0,5) (1,0)--(1,5) (0,1)--(1,1) (0,2)--(1,2) (0,3)--(1,3) (0,4)--(1,4);
		\node at (0.5,0.5) {$A_S$};
		\node at (0.5,1.5) {$\MarkovWeight{F}_{p,q}$};
		\node at (0.5,2.5) {$\MarkovWeight{T}$};
		\node at (0.5,3.5) {$\MarkovWeight{F}_{p,q'}$};
		\node at (0.5,4.5) {$A_N$};
		\node at (0,3.4) [circle, fill, inner sep=0.5mm] {};
		\node at (1,3.4) [circle, fill, inner sep=0.5mm] {};
		\draw [dotted] (0,3.4) -- (1,3.4);
		\draw [->] (1.3,3) -- node [midway, right] {$y-r$} (1.3,3.4);
	\end{tikzpicture}
	\]
	We may then use the l.h.s. (resp. r.h.s.) to apply first $\Phi_p^{S,q,(\infty_N,y)}$ (resp. $\Phi_p^{N,q',(\infty_S,y-r)}$) and eliminate the lower (upper) term $\MarkovWeight{F}_{p,q}$ (resp. $\MarkovWeight{F}_{p,q'}$) by correctly choosing $T_N$ (resp. $T_S$). Shifting again the base point with the same proposition allows one to apply the second morphism and hence the result the eigenvalue sequence is the same on the North and South side and can then be combined to produce $\lambda_{p(q+q')}$. 
\end{proof}

The same property holds between West and East systems of eigen-generators after use of dihedral group.

\subsection{Combining consecutive corners to half-planes}

The same mechanism as for Proposition~\ref{prop:simplificationonstrips} holds when extending the eigen-generators on the West or on the East as in Definition~\ref{def:cornereigensemigroups}: this is the content of the following Proposition~\ref{prop:NWandSWcorners:gluing:I}. There are however some non-trivial novelties when switching to gluing of corners.

First a new constraint is due to the fact that the North and South eigen-generators have to share the \emph{same} West (or East) semi-groups. In particular, this implies that the North and South eigenvalue sequences $\lambda_\bullet$ associated to the 2D-semi-group must coincide: it was only an assumption in Proposition~\ref{prop:simplificationonstrips} to obtain commutativity. Following the same steps as for Proposition~\ref{prop:simplificationonstrips}, we obtain the following proposition.

\begin{prop}\label{prop:NWandSWcorners:gluing:I}
Let $(M^W_\bullet,M^S_\bullet,\ca{V}_{\infty_W,\infty_S})$ and $(M^W_\bullet,M^N_\bullet,\ca{V}_{\infty_W,\infty_N})$ be a South-West and a North-West completed corner system of eigen-semi-groups of a 2D-semi-group $\MarkovWeight{F}_\bullet$ with the \emph{same} vertical semi-group $M^W_\bullet$ and eigenvalues sequences $(\lambda_\bullet^{SW},\sigma_\bullet^{W,(1)},\sigma_\bullet^{S})$ an $(\lambda_\bullet^{NW},\sigma_\bullet^{W,(2)},\sigma_\bullet^{N})$. 
If $\lambda_\bullet^{SW}=\lambda_{\bullet}^{NW}$ and $\sigma_\bullet^{W,(1)}=\sigma_\bullet^{W,(2)}$, then the morphisms $\Phi_{\infty_W}^{S,q,(\infty_N,y)}$ and $\Phi_{\infty_W}^{N,q',(\infty_N,y-r)}$ obtained by extending Definition~\ref{def:eigenalgebrauptomorph:leftextended} and equation~\eqref{eq:cornermorphism:def} to the vertical colours $\infty_N$ and $\infty_S$ with a base point satisfy, for all $q,q',r\in\setL^*$, $y\in\setP$, $v_{NW}\in\ca{V}_{\infty_W,\infty_N}$, $v_{SW}\in\ca{V}_{\infty_W,\infty_S}$ and $T_W\in\ca{A}_{\infty_W,r}$,
	\begin{equation}\begin{split}
		\Phi_{\infty_W}^{N,q',(\infty_S,y-r)}&\circ\Phi_{\infty_W}^{S,q,(\infty_N,y)}\left(
		\begin{tikzpicture}[guillpart,yscale=1.25,xscale=1.8]
			\fill[guillfill] (-1,0) rectangle (1,5);
			\draw[guillsep] (0,0)--(0,2) (0,3)--(0,5) (1,0)--(1,5) (-1,1)--(1,1) (-1,2)--(1,2) (-1,3)--(1,3) (-1,4)--(1,4);
			\node at (-0.5,0.5) {$v_{SW}$};
			\node at (-0.5,1.5) {$M^W_q$};
			\node at (-0.5,3.5) {$M^W_{q'}$};
			\node at (-0.5,4.5) {$v_{NW}$};
			\node at (0.5,0.5) {$M^S_p$};
			\node at (0.5,1.5) {$\MarkovWeight{F}_{p,q}$};
			\node at (0.,2.5) {$T_W$};
			\node at (0.5,3.5) {$\MarkovWeight{F}_{p,q'}$};
			\node at (0.5,4.5) {$M^W_p$};
			\node at (1,3.4) [circle, fill, inner sep=0.5mm] {};
			\draw [dotted] (-1,3.4) -- (1,3.4);
			\draw [->] (1.3,2) -- node [midway, right] {$y$} (1.3,3.4);
		\end{tikzpicture}\right)
		=
		\Phi_{\infty_W}^{S,q,(\infty_N,y)}\circ\Phi_{\infty_W}^{N,q',(\infty_S,y-r)}\left(
		\begin{tikzpicture}[guillpart,yscale=1.25,xscale=1.8]
		\fill[guillfill] (-1,0) rectangle (1,5);
		\draw[guillsep] (0,0)--(0,2) (0,3)--(0,5) (1,0)--(1,5) (-1,1)--(1,1) (-1,2)--(1,2) (-1,3)--(1,3) (-1,4)--(1,4);
		\node at (-0.5,0.5) {$v_{SW}$};
		\node at (-0.5,1.5) {$M^W_q$};
		\node at (-0.5,3.5) {$M^W_{q'}$};
		\node at (-0.5,4.5) {$v_{NW}$};
		\node at (0.5,0.5) {$M^S_p$};
		\node at (0.5,1.5) {$\MarkovWeight{F}_{p,q}$};
		\node at (0.,2.5) {$T_W$};
		\node at (0.5,3.5) {$\MarkovWeight{F}_{p,q'}$};
		\node at (0.5,4.5) {$M^W_p$};
		\node at (1,3.4) [circle, fill, inner sep=0.5mm] {};
		\draw [dotted] (-1,3.4) -- (1,3.4);
		\draw [->] (1.3,2) -- node [midway, right] {$y$} (1.3,3.4);
	\end{tikzpicture}
		\right)
		\\
		=& \lambda_{p(q+q')}\sigma^W_{q+q'}
		\begin{tikzpicture}[guillpart,yscale=1.25,xscale=1.8]
			\fill[guillfill] (-1,1) rectangle (1,4);
			\draw[guillsep] (0,1)--(0,2) (0,3)--(0,4) (1,1)--(1,4)  (-1,2)--(1,2) (-1,3)--(1,3);
			\node at (-0.5,1.5) {$v_{SW}$};
			\node at (-0.5,3.5) {$v_{NW}$};
			\node at (0.5,1.5) {$M^S_p$};
			\node at (0.,2.5) {$T_W$};
			\node at (0.5,3.5) {$M^N_p$};
			\node at (1,3.4) [circle, fill, inner sep=0.5mm] {};
			\draw [dotted] (-1,3.4) -- (1,3.4);
			\draw [->] (1.3,2) -- node [midway, right] {$y$} (1.3,3.4);
		\end{tikzpicture}
	\end{split}
\end{equation}
The composition of the two commuting morphisms is then written $\Phi_{\infty_W}^{SN,q,q',r,(y)}$.
\end{prop}
The proof is a simple exercise with the previous definitions.

Compared to the gluing of half-strips into strips, the second novelty with corner eigen-structures up to morphisms is the presence of morphisms in the transverse direction $\Phi_{\infty_S}^{W,p,p'}$ and $\Phi_{\infty_N}^{W,p,p'}$ such that:
\begin{subequations}
	\label{eq:neighbouringcornersmorphisms}
\begin{align}
\Phi_{\infty_N}^{W,p,p'}\left(
\begin{tikzpicture}[guillpart,yscale=-1.5,xscale=1.6]
	\fill[guillfill] (0,0) rectangle (3,2);
	\draw[guillsep] (3,0)--(3,2)--(0,2) (2,0)--(2,2) (1,0)--(1,2) (0,1)--(2,1);
	\node at (0.5,0.5) {$v_{NW}$};
	\node at (0.5,1.5) {$M_{q_2}^W$};
	\node at (1.5,0.5) {$M_{p}^N$};
	\node at (1.5,1.5) {$\MarkovWeight{F}_{p,q_2}$};
	\node at (2.5,1) {$\MarkovWeight{T_N}$};
\end{tikzpicture}
\right)
&=\lambda^{NW}_{pq}\sigma^S_p 
\begin{tikzpicture}[guillpart,yscale=-1.5,xscale=1.6]
	\fill[guillfill] (0,0) rectangle (2,2);
	\draw[guillsep] (2,0)--(2,2)--(0,2) (1,0)--(1,2) (0,1)--(1,1);
	\node at (0.5,0.5) {$v_{NW}$};
	\node at (0.5,1.5) {$M_{q_2}^W$};
	\node at (1.5,1) {$\MarkovWeight{T}_N$};
\end{tikzpicture}
\\
\Phi_{\infty_S}^{W,p,p'}\left(
\begin{tikzpicture}[guillpart,yscale=1.5,xscale=1.6]
	\fill[guillfill] (0,0) rectangle (3,2);
	\draw[guillsep] (3,0)--(3,2)--(0,2) (2,0)--(2,2) (1,0)--(1,2) (0,1)--(2,1);
	\node at (0.5,0.5) {$v_{SW}$};
	\node at (0.5,1.5) {$M_{q_1}^W$};
	\node at (1.5,0.5) {$M_{p}^S$};
	\node at (1.5,1.5) {$\MarkovWeight{F}_{p,q_1}$};
	\node at (2.5,1) {$\MarkovWeight{T_S}$};
\end{tikzpicture}
\right)
&=\lambda^{SW}_{pq}\sigma^S_p 
\begin{tikzpicture}[guillpart,yscale=1.5,xscale=1.6]
	\fill[guillfill] (0,0) rectangle (2,2);
	\draw[guillsep] (2,0)--(2,2)--(0,2) (1,0)--(1,2) (0,1)--(1,1);
	\node at (0.5,0.5) {$v_{SW}$};
	\node at (0.5,1.5) {$M_{q_1}^W$};
	\node at (1.5,1) {$\MarkovWeight{T}_S$};
\end{tikzpicture}
\end{align}
\end{subequations}
which must be compatible with gluing along a horizontal guillotine cut to form half-plane and strips.

\subsection{Combining half-planes with strips and opposite half-planes}.

In order to fulfil all wished simplification rules \eqref{eq:2D:wishedreplacements} starting from a full plane element, we must now introduce a way to simplify an expression such as
\[
\begin{tikzpicture}[guillpart,yscale=1.2,xscale=1.6]
	\fill[guillfill] (0,0) rectangle (3,3.5);
	\draw[guillsep] (3,0)--(3,3.5)  (2,0)--(2,3.5) (1,0)--(1,3.5) (0,1)--(2,1)  (0,2.5)--(2,2.5);
	\node at (0.5,0.5) {$v_{SW}$};
	\node at (0.5,3.) {$v_{NW}$};
	\node at (0.5,1.75) {$M_{q}^W$};
	\node at (1.5,0.5) {$M_{p}^S$};
	\node at (1.5,3.) {$M_{p}^N$};
	\node at (1.5,1.75) {$\MarkovWeight{F}_{p,q}$};
	\node at (2.5,1.75) {$T_{SN}$};
\end{tikzpicture}
\]
by replacing the middle part that contains the semi-group and the North and South eigen-semi-groups by $\lambda_{pq}\sigma^S_p\sigma^N_p$ up to some morphisms.

The situation is indeed very simple since it corresponds to a collection of one-dimensional situations, as described in Section~\ref{sec:onedimdiagowithmorph}: each vector \[
\begin{tikzpicture}[guillpart,yscale=1.2,xscale=1.6]
	\fill[guillfill] (0,0) rectangle (1,3.5);
	\draw[guillsep] (1,0)--(1,3.5) (0,1)--(1,1)  (0,2.5)--(1,2.5);
	\node at (0.5,0.5) {$v_{SW}$};
	\node at (0.5,3.) {$v_{NW}$};
	\node at (0.5,1.75) {$M_{q}^W$};
\end{tikzpicture}
\in \ca{A}_{\infty_W,\infty_{SN}}
\]
(with the specification of a base point)
is an eigenvector of the horizontal semigroup
\[
\begin{tikzpicture}[guillpart,yscale=1.2,xscale=1.6]
	\fill[guillfill] (0,0) rectangle (1,3.5);
	\draw[guillsep] (1,0)--(1,3.5) (0,1)--(1,1)  (0,2.5)--(1,2.5) (0,0)--(0,3.5);
	\node at (0.5,0.5) {$M^S_p$};
	\node at (0.5,3.) {$M^{N}_p$};
	\node at (0.5,1.75) {$\MarkovWeight{F}_{p,q}$};
\end{tikzpicture}
\in \ca{A}_{\infty_p,\infty_{SN}}
\]
(using the specification of the same base point)
with eigenvalue sequence $\mu^{SN}_p$ up to morphisms to $\Phi^{SN,p,r}$, where, as before, $r$ stands for the space $\ca{A}_{r,\infty_{SN}}$ to which $T_{SN}$ belongs. 

The only remaining task is to connect this collection of eigenvectors and their morphisms with the corner morphisms~\eqref{eq:neighbouringcornersmorphisms} as well as the compatibility conditions of Proposition~\ref{prop:NWandSWcorners:gluing:I}. 

\begin{defi}[compatible corners]\label{def:compatiblecorners:twobytwo}
	Let $\MarkovWeight{F}_{\bullet,\bullet}$ be a 2D-semigroup. Let the two triplets of spaces $(M^W_\bullet, M^S_\bullet, \ca{V}_{\infty_W,\infty_S})$ and $(M^W_\bullet,M^N_\bullet,\ca{V}_{\infty_W,\infty_N})$ be a South-West and a North-West completed corner system of eigen-semi-groups of a 2D-semi-group $\MarkovWeight{F}_\bullet$ with the \emph{same} vertical semi-group $M^W_\bullet$ and eigenvalues sequences $(\lambda_\bullet^{SW},\sigma_\bullet^{W,(1)},\sigma_\bullet^{S})$ an $(\lambda_\bullet^{NW},\sigma_\bullet^{W,(2)},\sigma_\bullet^{N})$. 
	
	These two corner systems are \emph{compatible} if and only if all the following conditions are fulfilled:
	\begin{enumerate}[(i)]
		\item the eigenvalue sequences $\lambda^{SW}_\bullet$ and $\lambda^{NW}_\bullet$ are equal;
		\item  the eigenvalue sequences $\sigma^{W,(1)}_\bullet$ and $\sigma^{W,(2)}_\bullet$ are equal;
		\item \label{item:stripeigenuptomorph} for all $q\in\setL^*$ and all $s\in\setP$, the spaces $\ca{V}_{\infty_W,\infty_{SN}}^{(q,s)}$ generated by
		\[
		\begin{tikzpicture}[guillpart,yscale=1.1,xscale=1.6]
			\fill[guillfill] (0,0) rectangle (1,3.5);
			\draw[guillsep] (1,0)--(1,3.5) (0,1)--(1,1)  (0,2.5)--(1,2.5);
			\node at (0.5,0.5) {$v_{SW}$};
			\node at (0.5,3.) {$v_{NW}$};
			\node at (0.5,1.75) {$M_{q}^W$};
						\node at (1,2.3) [circle, fill, inner sep=0.5mm] {};
			\draw [dotted] (0,2.3) -- (1,2.3);
			\draw [->] (1.3,1) -- node [midway, right] {$s$} (1.3,2.3);
		\end{tikzpicture}
		\]
		with $v_{SW}\in\ca{V}_{SW}$ and $v_{NW}\in\ca{V}_{NW}$ 
		are non-trivial eigen-spaces of the 1D-semi-group $M^{SN,(q,s)}_\bullet$ of $\ca{A}_{\bullet,\infty_{SN}}$
		defined by
		\[
		M^{SN,(q,s)}_p= \begin{tikzpicture}[guillpart,yscale=1.1,xscale=1.6]
			\fill[guillfill] (0,0) rectangle (1,3.5);
			\draw[guillsep] (1,0)--(1,3.5) (0,1)--(1,1)  (0,2.5)--(1,2.5) (0,0)--(0,3.5);
			\node at (0.5,0.5) {$M^S_p$};
			\node at (0.5,3.) {$M^N_p$};
			\node at (0.5,1.75) {$\MarkovWeight{F}_{p,q}$};
			\node at (1,2.3) [circle, fill, inner sep=0.5mm] {};
			\draw [dotted] (0,2.3) -- (1,2.3);
			\draw [->] (1.3,1) -- node [midway, right] {$s$} (1.3,2.3);
		\end{tikzpicture}
		\]
		with eigenvalue sequence $(\lambda_{pq}\sigma^{S}_p\sigma^{N}_p)_{p\in\setL^*}$ up to morphisms $\Phi^{W,p,r,(q,s)}_{\infty_W,\infty_{SN}}$ for $r\in\setL^*\cup\{\infty_E\}$.
		\item the morphisms $\Phi^{W,p,r,(q,s)}_{\infty_W,\infty_{SN}}$ satisfy, for any $r\in\setL\cup\{\infty_E^*\}$, for any $T_S\in\ca{A}_{r,\infty_S}$ and any $T_N\in\ca{A}_{r,\infty_S}$ and any partition $q=q_1+q_2$ in $\setL^*$,
		\begin{equation}\label{eq:eigenhalfplane:compatmorphisms}
		\begin{split}
		\Phi^{W,p,r,(q,s)}_{\infty_W,\infty_{SN}}&\left(
		\begin{tikzpicture}[guillpart,yscale=1.5,xscale=1.6]
			\fill[guillfill] (-1,0) rectangle (2,4);
			\draw[guillsep] (1,0)--(1,4) (-1,1)--(1,1)  (-1,3)--(1,3) (0,0)--(0,4) (2,0)--(2,4) (1,2)--(2,2) (-1,2)--(1,2);
			\node at (-0.5,0.5) {$v_{SW}$};
			\node at (-0.5,3.5) {$v_{NW}$};
			\node at (-0.5,1.5) {$M_{q_1}^W$};
			\node at (-0.5,2.5) {$M_{q_2}^W$};
			\node at (0.5,0.5) {$M^S_p$};
			\node at (0.5,3.5) {$M^N_p$};
			\node at (0.5,1.5) {$\MarkovWeight{F}_{p,q_1}$};
			\node at (0.5,2.5) {$\MarkovWeight{F}_{p,q_2}$};
			\node at (1.5,3) {$T_{N}$};
			\node at (1.5,1) {$T_{S}$};
			\node at (2,2.8) [circle, fill, inner sep=0.5mm] {};
			\draw [dotted] (-1,2.8) -- (2,2.8);
			\draw [->] (2.3,1) -- node [midway, right] {$s$} (2.3,2.8);
		\end{tikzpicture}
		\right)
		\\
		&= \begin{tikzpicture}[guillpart,yscale=1.5,xscale=1.5]
			\fill[guillfill] (0,0) rectangle (1,2);
			\draw[guillsep] (1,0)--(1,2) (0,1)--(1,1);
			\node at (0.5,0.5) {$1$};
			\node at (0.5,1.5) {$2$};
			\node at (1,1.8) [circle, fill, inner sep=0.5mm] {};
			\draw [dotted] (0,1.8) -- (1,1.8);
			\draw [->] (1.15,1) -- node [midway, right] {{\footnotesize $s-q_1$}} (1.15,1.8);
		\end{tikzpicture}\left(
	\Phi_{\infty_W,\infty_S}^{W,p,r}\left(
		\begin{tikzpicture}[guillpart,yscale=1.5,xscale=1.5]
			\fill[guillfill] (0,0) rectangle (3,2);
			\draw[guillsep] (3,0)--(3,2)--(0,2) (2,0)--(2,2) (1,0)--(1,2) (0,1)--(2,1);
			\node at (0.5,0.5) {$v_{SW}$};
			\node at (1.5,0.5) {$M_p^S$};
			\node at (0.5,1.5) {$M_{q_1}^W$};
			\node at (1.5,1.5) {$\MarkovWeight{F}_{p,q_1}$};
			\node at (2.5,1) {$T_S$};
		\end{tikzpicture}
		\right),
			\Phi_{\infty_W,\infty_N}^{W,p,r}\left(
			    \begin{tikzpicture}[guillpart,yscale=1.5,xscale=1.5]
			    	\fill[guillfill] (0,0) rectangle (3,2);
			    	\draw[guillsep] (0,0)--(3,0)--(3,2) (1,0)--(1,2) (2,0)--(2,2) (0,1)--(2,1);
			    	\node at (0.5,0.5) {$M_{q_2}^W$};
			    	\node at (1.5,0.5) {$\MarkovWeight{F}_{p,q_2}$};
			    	\node at (0.5,1.5) {$v_{NW}$};
			    	\node at (1.5,1.5) {$M_p^N$};
			    	\node at (2.5,1) {$T_N$};
			    \end{tikzpicture}
			\right)
		\right)
		\end{split}
		\end{equation}
	\end{enumerate}
\end{defi}

It is interesting to see that, if $\ca{A}_{p,\infty_{SN}}$ contains all the linear maps on $\ca{B}_{\infty_W,\infty_{SN}}$ then the morphisms $\Phi^{SN,p,r,(q,s)}$ are just multiplications by constants (which can be avoided or absorb the eigenvalue) following Section~\ref{sec:onedimdiagowithmorph}: in this simple case, \eqref{eq:eigenhalfplane:compatmorphisms} states that the corner morphisms let the pairing from two corners to a half-plane invariant up to some scalar value.

The four corner can then be glued together by the introduction of linear maps on the full plane which can be decomposed the previous corners.

\begin{prop}[full compatibility of four corners]\label{prop:compatiblecorners:full}
	Let $(M^W_\bullet,M^S_\bullet,\ca{V}_{\infty_W,\infty_S})$,  $(M^W_\bullet,M^N_\bullet,\ca{V}_{\infty_W,\infty_N})$, 
	$(M^E_\bullet,M^S_\bullet,\ca{V}_{\infty_E,\infty_S})$ and 
	$(M^E_\bullet,M^N_\bullet,\ca{V}_{\infty_E,\infty_N})$ be four completed corner eigen-spaces of a 2D-semigroup $\MarkovWeight{F}_{\bullet,\bullet}$ such that any two eigen corner systems are compatible, according to Definition~\ref{def:compatiblecorners:twobytwo}. Then, the following properties hold:
	\begin{enumerate}
		\item for any $v_{ab}\in\ca{V}_{\infty_b,\infty_a}$, any $r\in\setL$, any $T_{SN}\in\ca{A}_{r,\infty_{SN}}$, using \[
		a= \begin{tikzpicture}[guillpart,yscale=1.25,xscale=1.6]
			\fill[guillfill] (0,0) rectangle (5,4);
			\draw[guillsep] (1,0)--(1,4) (2,0)--(2,4) (3,0)--(3,4) (4,0)--(4,4);
			\draw[guillsep] (0,1)--(2,1) ( (0,3)--(2,3) (3,1)--(5,1)  (3,3)--(5,3);
			\node at (0.5,0.5) {$v_{SW}$};
			\node at (1.5,0.5) {$M_{p_1}^S$};
			\node at (3.5,0.5) {$M_{p_2}^S$};
			\node at (4.5,0.5) {$v_{SE}$};
			\node at (0.5,2) {$M^W_{q}$};
			\node at (1.5,2) {$\MarkovWeight{F}_{p_1,q}$};
			\node at (2.5,2) {$T_{SN}$};
			\node at (3.5,2) {$\MarkovWeight{F}_{p_2,q}$};
			\node at (4.5,2) {$M^E_{q}$};
			\node at (0.5,3.5) {$v_{NW}$};
			\node at (1.5,3.5) {$M_{p_1}^N$};
			\node at (3.5,3.5) {$M_{p_2}^N$};
			\node at (4.5,3.5) {$v_{NE}$};	
			\node at (2.75,2.5) [circle, fill, inner sep=0.5mm] {};
			\draw [dotted] (0,2.5) -- (5,2.5);
			\draw [dotted] (2.75,0) -- (2.75,4);
			\draw [->] (5.15,1.) -- node [midway, right] {{\footnotesize $y$}} (5.15,2.5);
			\draw [->] (2,4.15) -- node [midway, above] {{\footnotesize $x$}} (2.75,4.15);
		\end{tikzpicture} \in \ca{A}_{\infty_{WE},\infty_{SN}}
		\]
		we have:
		\begin{subequations}\label{eq:morphonfullplanefromhalfplanes}
		\begin{equation}
			\begin{split}
\Phi_{\infty_W,\infty_{SN}}^{W,p_1,(\infty_E,x),(q,y)}\circ \Phi_{\infty_E,\infty_{SN}}^{E,p_2,(\infty_W,x-r),(q,y)} 	\left(	a	\right)
&=
\Phi_{\infty_E,\infty_{SN}}^{E,p_2,(\infty_W,x-r),(q,y)}\circ \Phi_{\infty_W,\infty_{SN}}^{W,p_1,(\infty_W,x),(q,y)} 	\left(	a	\right)
\\
&=
\sigma^S_{p_1+p_2}\sigma^N_{p_1+p_2}\lambda_{(p_1+p_2)q}	\begin{tikzpicture}[guillpart,yscale=1.25,xscale=1.6]
	\fill[guillfill] (1,0) rectangle (4,4);
	\draw[guillsep]  (2,0)--(2,4) (3,0)--(3,4);
	\draw[guillsep] (1,1)--(2,1) (1,3)--(2,3) (3,1)--(4,1)  (3,3)--(4,3);
	\node at (1.5,0.5) {$v_{SW}$};
	\node at (3.5,0.5) {$v_{SE}$};
	\node at (1.5,2) {$M^W_{q}$};
	\node at (2.5,2) {$T_{SN}$};
	\node at (3.5,2) {$M^E_{q}$};
	\node at (1.5,3.5) {$v_{NW}$};
	\node at (3.5,3.5) {$v_{NE}$};	
	\node at (2.75,2.5) [circle, fill, inner sep=0.5mm] {};
	\draw [dotted] (1,2.5) -- (4,2.5);
	\draw [dotted] (2.75,0) -- (2.75,4);
	\draw [->] (4.15,1.) -- node [midway, right] {{\footnotesize $y$}} (4.15,2.5);
	\draw [->] (2,4.15) -- node [midway, above] {{\footnotesize $x$}} (2.75,4.15);
\end{tikzpicture} 
\end{split}
		\end{equation}
	\item for any $v_{ab}\in\ca{V}_{\infty_b,\infty_a}$, any $r\in\setL$, any $T_{WE}\in\ca{A}_{\infty_{WE},r}$, using \[
	b= \begin{tikzpicture}[guillpart,yscale=1.5,xscale=1.6,rotate=90]
		\fill[guillfill] (0,0) rectangle (5,4);
		\draw[guillsep] (1,0)--(1,4) (2,0)--(2,4) (3,0)--(3,4) (4,0)--(4,4);
		\draw[guillsep] (0,1)--(2,1) ( (0,3)--(2,3) (3,1)--(5,1)  (3,3)--(5,3);
		\node at (0.5,0.5) {$v_{SE}$};
		\node at (1.5,0.5) {$M_{p_1}^S$};
		\node at (3.5,0.5) {$M_{p_2}^S$};
		\node at (4.5,0.5) {$v_{NE}$};
		\node at (0.5,2) {$M^S_{p}$};
		\node at (1.5,2) {$\MarkovWeight{F}_{p,q_1}$};
		\node at (2.5,2) {$T_{WE}$};
		\node at (3.5,2) {$\MarkovWeight{F}_{p,q_2}$};
		\node at (4.5,2) {$M^N_{p}$};
		\node at (0.5,3.5) {$v_{SW}$};
		\node at (1.5,3.5) {$M_{q_1}^W$};
		\node at (3.5,3.5) {$M_{q_2}^W$};
		\node at (4.5,3.5) {$v_{NW}$};	
		\node at (2.75,2.5) [circle, fill, inner sep=0.5mm] {};
		\draw [dotted] (0,2.5) -- (5,2.5);
		\draw [dotted] (2.75,0) -- (2.75,4);
		\draw [->] (5.15,1.) -- node [midway, above] {{\footnotesize $x$}} (5.15,2.5);
		\draw [->] (2,4.15) -- node [midway, left] {{\footnotesize $y$}} (2.75,4.15);
	\end{tikzpicture} \in \ca{A}_{\infty_{WE},\infty_{SN}}
	\]
	we have
	\begin{equation}
		\begin{split}
			\Phi_{\infty_S,\infty_{WE}}^{S,q_1,(\infty_N,y),(p,x)}\circ \Phi_{\infty_N,\infty_{WE}}^{N,q_2,(\infty_S,y-r),(p,x)} 	\left(	b	\right)
			&=
			\Phi_{\infty_N,\infty_{WE}}^{N,q_2,(\infty_S,y-r),(p,x)}
			\circ
			\Phi_{\infty_S,\infty_{WE}}^{S,q_1,(\infty_N,y),(p,x)}
			\left(	b	\right)
			\\
			&=
			\sigma^W_{q_1+q_2}\sigma^E_{q_1+q_2}\lambda_{p(q_1+q_2)}	\begin{tikzpicture}[guillpart,yscale=1.5,xscale=1.6,rotate=90]
				\fill[guillfill] (1,0) rectangle (4,4);
				\draw[guillsep]  (2,0)--(2,4) (3,0)--(3,4);
				\draw[guillsep] (1,1)--(2,1) (1,3)--(2,3) (3,1)--(4,1)  (3,3)--(4,3);
				\node at (1.5,0.5) {$v_{SE}$};
				\node at (3.5,0.5) {$v_{NE}$};
				\node at (1.5,2) {$M^S_{p}$};
				\node at (2.5,2) {$T_{WE}$};
				\node at (3.5,2) {$M^N_{p}$};
				\node at (1.5,3.5) {$v_{SW}$};
				\node at (3.5,3.5) {$v_{NW}$};	
				\node at (2.75,2.5) [circle, fill, inner sep=0.5mm] {};
				\draw [dotted] (1,2.5) -- (4,2.5);
				\draw [dotted] (2.75,0) -- (2.75,4);
				\draw [->] (4.15,1.) -- node [midway, above] {{\footnotesize $x$}} (4.15,2.5);
				\draw [->] (2,4.15) -- node [midway, left] {{\footnotesize $y$}} (2.75,4.15);
			\end{tikzpicture} 
		\end{split}
	\end{equation}
	\end{subequations}
	\end{enumerate}
\end{prop}
Equations~\eqref{eq:morphonfullplanefromhalfplanes} show that the morphisms on opposite half-planes commute and their composition define linear maps of the full plane space $\ca{A}_{\infty_{WE},\infty_{SN}}$. For Markov processes, we have seen that both for the canonical structure and for ROPEs ---and thus for $\ca{E}_{\infty_{WE},\infty_{SN}}$---, this space is just the scalar field: the linear maps thus obtained on the full plane space are just scalar multiplications.
\begin{proof}
	We prove only in detail the identity \[
\Phi_{\infty_E,\infty_{SN}}^{E,p_2,(\infty_W,x-r),(q,y)}\circ 	\Phi_{\infty_W,\infty_{SN}}^{W,p_1,(\infty_E,x),(q,y)}	\left(	a	\right)
	=
	\sigma^S_{p}\sigma^N_{p}\lambda_{pq}	\begin{tikzpicture}[guillpart,yscale=1.25,xscale=1.6]
		\fill[guillfill] (1,0) rectangle (4,4);
		\draw[guillsep]  (2,0)--(2,4) (3,0)--(3,4);
		\draw[guillsep] (1,1)--(2,1) (1,3)--(2,3) (3,1)--(4,1)  (3,3)--(4,3);
		\node at (1.5,0.5) {$v_{SW}$};
		\node at (3.5,0.5) {$v_{SE}$};
		\node at (1.5,2) {$M^W_{q}$};
		\node at (2.5,2) {$T_{SN}$};
		\node at (3.5,2) {$M^E_{q}$};
		\node at (1.5,3.5) {$v_{NW}$};
		\node at (3.5,3.5) {$v_{NE}$};	
		\node at (2.75,2.5) [circle, fill, inner sep=0.5mm] {};
		\draw [dotted] (1,2.5) -- (4,2.5);
		\draw [dotted] (2.75,0) -- (2.75,4);
		\draw [->] (4.15,1.) -- node [midway, right] {{\footnotesize $y$}} (4.15,2.5);
		\draw [->] (2,4.15) -- node [midway, above] {{\footnotesize $x$}} (2.75,4.15);
	\end{tikzpicture} 
	\]
	with $p=p_1+p_2$, since the three other ones follow the same computations. We first consider apply $\Phi_{\infty_W,\infty_{SN}}^{W,p_1,(\infty_E,x),(q,y)}$ using item~(\ref{item:stripeigenuptomorph}) of Definition~\ref{def:compatiblecorners:twobytwo} applied to the two compatible corners on the West with $r=(\infty_E,x)$ and an element $T \in \ca{A}_{\infty_E,\infty_{SN}}$ given by:
	\[
	 T = \begin{tikzpicture}[guillpart,yscale=1.25,xscale=1.6]
	 	\fill[guillfill] (2,0) rectangle (5,4);
	 	\draw[guillsep]  (2,0)--(2,4) (3,0)--(3,4) (4,0)--(4,4);
	 	\draw[guillsep]   (3,1)--(5,1)  (3,3)--(5,3);
	 	\node at (3.5,0.5) {$M_{p_2}^S$};
	 	\node at (4.5,0.5) {$v_{SE}$};
	 	\node at (2.5,2) {$T_{SN}$};
	 	\node at (3.5,2) {$\MarkovWeight{F}_{p_2,q}$};
	 	\node at (4.5,2) {$M^E_{q}$};
	 	\node at (3.5,3.5) {$M_{p_2}^N$};
	 	\node at (4.5,3.5) {$v_{NE}$};	
	 	\node at (2.,2.5) [circle, fill, inner sep=0.5mm] {};
	 	\draw [dotted] (2,2.5) -- (5,2.5);
	 	\draw [->] (5.15,1.) -- node [midway, right] {{\footnotesize $y$}} (5.15,2.5);
	 \end{tikzpicture}
	\]
	We then obtain $\sigma^{S,(W)}_{p_1}\sigma^{N,(W)}_{p_1}\lambda^{(W)}_{p_1 q}$ times the same shape with the vertical strip made of $M_{p_1}^N$, $\MarkovWeight{F}_{p_1,q}$ and $M^N_{p_1}$, where $\sigma^{a,(b)}_\bullet$ stands for the eigenvalue sequences of $M_\bullet^{a}$ relatively with the West corner.
	
	Using the compatibility of the two East corners, we may then apply the morphisms $\Phi_{\infty_E,\infty_{SN}}^{E,p_2,(\infty_W,x-r),(q,y)}$ with an element $T'\in\ca{A}_{\infty_W,\infty_{SN}}$ now given by:
	\[
	T' = \begin{tikzpicture}[guillpart,yscale=1.25,xscale=1.6]
		\fill[guillfill] (1,0) rectangle (3,4);
		\draw[guillsep]  (2,0)--(2,4) (3,0)--(3,4);
		\draw[guillsep] (1,1)--(2,1) (1,3)--(2,3) ;
		\node at (1.5,0.5) {$v_{SW}$};
		\node at (1.5,2) {$M^W_{q}$};
		\node at (2.5,2) {$T_{SN}$};
		\node at (1.5,3.5) {$v_{NW}$};
		\node at (3.,2.5) [circle, fill, inner sep=0.5mm] {};
		\draw [dotted] (1,2.5) -- (3,2.5);
		\draw [->] (3.15,1.) -- node [midway, right] {{\footnotesize $y$}} (3.15,2.5);
	\end{tikzpicture} 
	\]
	and obtain $\sigma^{S,(E)}_{p_2}\sigma^{N,(E)}_{p_2}\lambda^{(E)}_{p_2 q}$ times the expected result. Using the compatibility of the North and South corners, we obtain the equality of the eigenvalue sequences on the West and on the East sides and the final results.
\end{proof}

Moreover, both items of the previous proposition can be composed one with the other to introduce morphisms in both directions to achieve full simplifications of expressions such as (base point skipped)
\[
x=\begin{tikzpicture}[guillpart,yscale=1.5,xscale=1.85]
	\fill[guillfill] (0,0) rectangle (5,5);
	\draw[guillsep] (1,0)--(1,5)  (2,0)--(2,5) (3,0)--(3,5) (4,0)--(4,5);
	\draw[guillsep] (0,1)--(5,1) (0,2)--(5,2) (0,3)--(5,3) (0,4)--(5,4);
	\node at (0.5,0.5) {$v_{SW}$};
		\node at (1.5,0.5) {$M^S_{p_1}$};
		\node at (2.5,0.5) {$M^S_{r}$};
		\node at (3.5,0.5) {$M^S_{p_2}$};
		\node at (4.5,0.5) {$v_{SE}$};
	\node at (0.5,1.5) {$M_{q_1}^W$};
		\node at (1.5,1.5) {$\MarkovWeight{F}_{p_1,q_1}$};
		\node at (2.5,1.5) {$\MarkovWeight{F}_{r,q_1}$};
		\node at (3.5,1.5) {$\MarkovWeight{F}_{p_2,q_1}$};
		\node at (4.5,1.5) {$M^E_{q_1}$};
	\node at (0.5,2.5) {$M^W_{s}$};
		\node at (1.5,2.5) {$\MarkovWeight{F}_{p_1,s}$};
		\node at (2.5,2.5) {$\MarkovWeight{T}$};
		\node at (3.5,2.5) {$\MarkovWeight{F}_{p_2,s}$};
		\node at (4.5,2.5) {$M^E_{s}$};
	\node at (0.5,3.5) {$M^W_{q_2}$};
		\node at (1.5,3.5) {$\MarkovWeight{F}_{p_1,q_2}$};
		\node at (2.5,3.5) {$\MarkovWeight{F}_{r,q_2}$};
		\node at (3.5,3.5) {$\MarkovWeight{F}_{p_2,q_2}$};
		\node at (4.5,3.5) {$M^E_{q_2}$};
	\node at (0.5,4.5) {$v_{NW}$};
		\node at (1.5,4.5) {$M^N_{p_1}$};
		\node at (2.5,4.5) {$M^N_{r}$};
		\node at (3.5,4.5) {$M^N_{p_2}$};
		\node at (4.5,4.5) {$v_{NE}$};
		\end{tikzpicture}
\]
as required in Section~\ref{par:requiredsimplificationrules}.

\begin{coro}\label{coro:fullplane:eigen} Let $(M^W_\bullet,M^S_\bullet,\ca{V}_{\infty_W,\infty_S})$, $(M^W_\bullet,M^N_\bullet,\ca{V}_{\infty_W,\infty_N})$, 
$(M^E_\bullet,M^S_\bullet,\ca{V}_{\infty_E,\infty_S})$ and $(M^E_\bullet,M^N_\bullet,\ca{V}_{\infty_E,\infty_N})$ be four compatible completed corner eigen-spaces of a 2D-semigroup $\MarkovWeight{F}_{\bullet,\bullet}$.
If $\ca{A}_{\infty_{WE},\infty_{SN}}=\setK$ then there exists a constant $C\in\setK^*$ such that
\begin{equation}\label{eq:fullplane:eigen}
\begin{tikzpicture}[guillpart,yscale=1.2,xscale=1.6]
	\fill[guillfill] (0,0) rectangle (4,3.5);
	\draw[guillsep] (1,0)--(1,3.5)   (3,0)--(3,3.5);
	\draw[guillsep] (0,1)--(4,1)  (0,2.5)--(4,2.5);
	\node at (0.5,0.5) {$v_{SW}$};
	\node at (3.5,0.5) {$v_{SE}$};
	\node at (0.5,3.) {$v_{NW}$};
	\node at (3.5,3.) {$v_{NE}$};
	\node at (2,0.5) {$M_p^{S}$};
	\node at (2,3.) {$M_p^{N}$};
	\node at (0.5,1.75) {$M_q^{W}$};
	\node at (3.5,1.75) {$M_q^{E}$};
	\node at (2,1.75) {$\MarkovWeight{F}_{p,q}$};
	\node at (2.5,2.1) [circle, fill, inner sep=0.5mm] {};
	\draw [dotted] (0,2.1) -- (4,2.1);
	\draw [dotted] (2.5,0) -- (2.5,3.5);
	\draw [->] (4.15,1.) -- node [midway, right] {{\footnotesize $y$}} (4.15,2.1);
	\draw [->] (1,3.65) -- node [midway, above] {{\footnotesize $x$}} (2.5,3.65);
\end{tikzpicture}
= C \lambda_{pq}\sigma^N_p\sigma^S_p\sigma^W_q\sigma^E_q
\end{equation}
\end{coro}
\begin{proof}
The morphisms on the West half-plane of the type
$\Phi_{\infty_W,\infty_{SN}}^{W,p',(\infty_E,x),(q,y)}$ as well as their counterpart on the three other sides are all linear maps on $\ca{A}_{\infty_{WE},\infty_{SN}}$. Under the hypothesis $\ca{A}_{\infty_{WE},\infty_{SN}}=\setK$, they are just multiplication maps $x\mapsto cx$ on $\setK$. 

Thus, up to a normalization constant, the l.h.s. of \eqref{eq:fullplane:eigen} is equal to its image by any of these morphisms. Using repeatedly Proposition~\ref{prop:compatiblecorners:full} with a suitable choice of the $T$ elements (or simply removing them), the l.h.s. of \eqref{eq:fullplane:eigen} is proportional to any similar elements with smaller $p$ and $q$. In the extremal case in which $p$ and $q$ are null, we obtain the eigenvalue in the r.h.s. times the element 
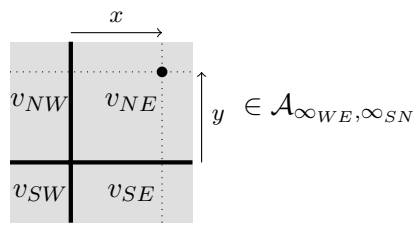
\[\begin{tikzpicture}[guillpart,yscale=1.2,xscale=1.6]
	\fill[guillfill] (0,0) rectangle (3,3);
	\draw[guillsep] (0,1)--(3,1) (1,0)--(1,3);
	\node at (0.5,0.5) {$v_{SW}$};
	\node at (2,0.5) {$v_{SE}$};
	\node at (0.5,2) {$v_{NW}$};
	\node at (2,2) {$v_{NE}$};
	\node at (2.5,2.5) [circle, fill, inner sep=0.5mm] {};
	\draw [dotted] (0,2.5) -- (3,2.5);
	\draw [dotted] (2.5,0) -- (2.5,3);
	\draw [->] (3.15,1.) -- node [midway, right] {{\footnotesize $y$}} (3.15,2.5);
	\draw [->] (1,3.15) -- node [midway, above] {{\footnotesize $x$}} (2.5,3.15);
\end{tikzpicture} \in \ca{A}_{\infty_{WE},\infty_{SN}}
\]
which do not depend on $p$ and $q$. Collecting the constants from the horizontal and the vertical simplification and the value of the gluing of the four corners alone, we obtain the expected result.
\end{proof}

Again, there is a nice interpretation of the r.h.s. of~\eqref{eq:fullplane:eigen}: the eigenvalue associated to the 2D-semigroup is a surface term with an exponent $pq$ equal to the area and each eigenvalue associated to a side 1D-semi-group is a lineic term with an exponent $p$ or $q$ equal to the length of the correspond side. The constant term corresponds to the four corners. Formula~\eqref{eq:fullplane:eigen} is the two-dimensional analogue of $v_\lambda^L A^pv_\lambda^R= c \lambda^p$ where $v_\lambda^{L}$ and $v_\lambda^R$ are the left and right eigen-vectors of $A$ for the eigenvalue $\lambda$.

\begin{rema}
	We do not spend too many lines about the description of the constant $C$ since it would produce, in the context of Markov processes, just an overall normalization of the partition functions. In particular, $C$ contains a normalization of the four corner elements (no dependence on $p$ and $q$) and may be changed by geometric sequences in $p$, $q$ and $pq$ by a redefinition of the morphisms and the eigenvalue sequences as in Section~\ref{sec:onedimdiagowithmorph}.
\end{rema}

\section{Full boundary eigen-structures: a summary and some questions.}

The previous sections are lengthy and introduce a lot of notations, especially in the description of morphisms. We spend some lines here to provide heuristic descriptions on the computations.

\subsection{How to identify an element simplifiable through morphisms and eigenvalues}\label{sec:summaryeigensituations}
Morphisms are introduced for every guillotine partition of a shape such that:
\begin{itemize}
	\item at least one direction, let say $a$, is infinite:
	
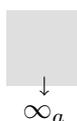
\[
	\begin{tikzpicture}[guillpart,yscale=0.8,xscale=0.8]
		\fill[guillfill] (0,0) rectangle (2,2);
		\draw[->] (1,0.25)--(1,-0.25);
		\node at (1.,-0.25) [below] {$\infty_a$};
	\end{tikzpicture}
	\]
	\item there are at least \emph{two} guillotine cuts in the transverse direction, cutting the domain in three regions:
	
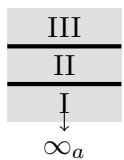
\[
	\begin{tikzpicture}[guillpart,yscale=1.,xscale=1.]
		\fill[guillfill] (0,0) rectangle (3,3);
		\draw[guillsep] (0,1)--(3,1) (0,2)--(3,2);
		\draw[->] (1.5,0.25)--(1.5,-0.25);
		\node at (1.5,-0.25) [below] {$\infty_a$};
		\node at (1.5,0.5) {I};
		\node at (1.5,1.5) {II};
		\node at (1.5,2.5) {III};
	\end{tikzpicture}
	\]
	\item in region I (infinite side), only the eigen 1D-semi-group $M^a_p$ with or without eigen corner space (depending on whether the transverse direction is finite or not)
	\item in region II, only the 2D-semi-group $\MarkovWeight{F}_{p,q}$ with or without eigen 1D-semigroup $M^b_q$ (depending on whether transverse direction is finite or not, see previous point)
	\item region III may be empty or contain any arbitrary element .
\end{itemize}
In the presence of two infinite directions with this structure, as in
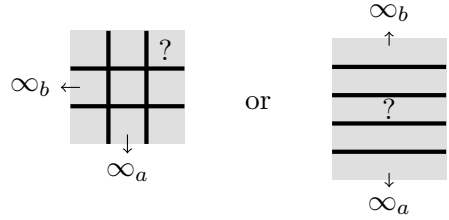
\[
\begin{tikzpicture}[guillpart,yscale=1.,xscale=1.]
	\fill[guillfill] (0,0) rectangle (3,3);
	\draw[guillsep] (0,1)--(3,1) (0,2)--(3,2);
	\draw[guillsep] (1,0)--(1,3) (2,0)--(2,3);
	\draw[->] (1.5,0.25)--(1.5,-0.25);
	\node at (1.5,-0.25) [below] {$\infty_a$};
	\draw[->] (0.25,1.5)--(-0.25,1.5);
	\node at (-0.25,1.5) [left] {$\infty_b$};
	\node at (2.5,2.5) {$?$};
\end{tikzpicture}
\qquad
\text{or}
\qquad
\begin{tikzpicture}[guillpart,yscale=0.75,xscale=1.]
	\fill[guillfill] (0,0) rectangle (3,5);
	\draw[guillsep] (0,1)--(3,1) (0,2)--(3,2) (0,3)--(3,3) (0,4)--(3,4);
	\draw[->] (1.5,0.25)--(1.5,-0.25);
	\node at (1.5,-0.25) [below] {$\infty_a$};
	\draw[->] (1.5,4.75)--(1.5,5.25);
	\node at (1.5,5.25) [above] {$\infty_b$};
	\node at (1.5,2.5) {$?$};
\end{tikzpicture}
\]
only the region $?$ may contain arbitrary elements. In both cases, simplifications through morphisms along the two directions can be performed in any order \emph{provided} all the adjacent eigen corner spaces present in the diagram are compatible.

\subsection{Morphisms commute with some operadic products}

When $T$ is itself a horizontal product of $N$ elements $T'_k$, \eqref{eq:defeigenSouth:withM} can be rewritten as
\[
		\Phi_{P}^{S,q,r}\left(
\begin{tikzpicture}[guillpart,yscale=1.1,xscale=2.]
	\fill[guillfill] (0,0) rectangle (3,3);
	\draw[guillsep] (0,0)--(0,3)--(3,3)--(3,0) (1,0)--(1,3) (2,0)--(2,3) (0,1)--(3,1) (0,2)--(3,2);
	\node at (0.5,0.5) {$v_1$};
	\node at (1.5,0.5) {$\dots$};
	\node at (2.5,0.5) {$v_N$};
	\node at (0.5,1.5) {$\MarkovWeight{F}_{p_1,q}$};
	\node at (1.5,1.5) {$\dots$};
	\node at (2.5,1.5) {$\MarkovWeight{F}_{p_N,q}$};
	\node at (0.5,2.5) {$T'_1$};
	\node at (1.5,2.5) {$\dots$};
	\node at (2.5,2.5) {$T'_N$};
\end{tikzpicture}
\right)
= \begin{tikzpicture}[guillpart,yscale=1.1,xscale=1.25]
	\fill[guillfill] (0,0) rectangle (3,1);
	\draw[guillsep] (0,0)--(0,1)--(3,1)--(3,0) (1,0)--(1,1) (2,0)--(2,1);
	\node at (0.5,0.5) {$1$};
	\node at (1.5,0.5) {$\ldots$};
	\node at (2.5,0.5) {$N$};
\end{tikzpicture}
\left(
\Phi_{p_1}^{S,q,r}\left(
\begin{tikzpicture}[guillpart,yscale=1.1,xscale=1.65]
	\fill[guillfill] (0,0) rectangle (1,3);
	\draw[guillsep] (0,0)--(0,3)--(1,3)--(1,0) (0,1)--(1,1) (0,2)--(1,2);
	\node at (0.5,0.5) {$v_1$};
	\node at (0.5,1.5) {$\MarkovWeight{F}_{p_1,q}$};
	\node at (0.5,2.5) {$T'_1$};
\end{tikzpicture}
\right),
\ldots,
\Phi_{p_N}^{S,q,r}\left(
\begin{tikzpicture}[guillpart,yscale=1.1,xscale=1.75]
	\fill[guillfill] (0,0) rectangle (1,3);
	\draw[guillsep] (0,0)--(0,3)--(1,3)--(1,0) (0,1)--(1,1) (0,2)--(1,2);
	\node at (0.5,0.5) {$v_N$};
	\node at (0.5,1.5) {$\MarkovWeight{F}_{p_N,q}$};
	\node at (0.5,2.5) {$T'_N$};
\end{tikzpicture}
\right)
\right)
\]
which tends to indicate that the morphisms have some commutation relations with the operadic products. The same type of relation holds on corners as it can be seen by rewriting \eqref{eq:cornermorphism:def} in the same way, as well as for half-planes and strips.

However, it is not clear to us which status these relations have, for three reasons. First, this phenomenon of commutation does not all hold for any morphism and any product: it requires a geometric configuration as in the previous Section~\ref{sec:summaryeigensituations}. It seems to be more related to the three elementary associativities \eqref{eq:guill2:listassoc} since at least two or three cuts are required than the $2$-ary products themselves. Second, there is a whole system of indices labelling the morphisms, which makes, at our present stage of understanding $\Guill_2$-operads and their boundary structure, hard to guess which precise geometric structure governs these morphisms. Third, we do not put any constraint on the action of the morphisms on other elements than the ones in the eigen-spaces (see section below) and thus it is hard to see at which precise level of generality the previous observed commutation relations shall be lifted.

We believe that the trails to follow are threefold: full computations with concrete models of probability theory and statistical mechanics, generalization of $\Guill_2$-operads and their boundary operads to non-rectangular shapes (in this case, the precise description of morphisms is less clear to us) and the creation of deeper relations with other tools of higher algebra.

\subsection{From one eigen-structures to several ones: a question about morphisms}

Most definitions in the present paper deal with 1D-semigroups on the four boundary sides. In the direction in which the 2d-semi-group acts, it is the equivalent of the 1D situation of a single eigenvector of an operator. In our probabilistic goal of infinite volume Gibbs measures, the equivalent of the Perron-Frobenius eigenvector is enough (see Section~\ref{sec:gibbsfromeigen}). From the perspective of linear algebra, this may be only just a starting point: in the one-dimensional traditional associative framework, beyond the single definition of an eigenvector, lie eigenspaces, Dunford decomposition, diagonalization of endomorphisms, the notion of spectrum and the whole yoga of commutative algebra.

Because of a lack of concrete situations we know how to treat until this point, we do not know yet how to put nicely together two generic eigen 1D-semi-groups on the same side with different eigenvalue sequences together. Indeed, in the present paper, each of them comes with its own set of morphisms whose action on the other 1D-semi-groups is not specified. We think that this is only a side-effect of our approach and we believe that the two set of morphisms must be compatible in some sense or obtained by restriction of some higher object which probably lives at the level of some generalized notion of spectrum of the 2D-semi-group. We have not managed to reach this level of generality due to a lack of imagination in the generalization of classical properties at a level "up to morphisms", beyond the particular case of Section~\ref{sec:commutuptomorph}.

	\section{Applications of boundary eigen-elements to infinite-volume Gibbs measures}\label{sec:gibbsfromeigen}

\subsection{From ROPErep full boundary eigen-elements to Gibbs measures through an extension theorem.}

We finally state the main result of the paper, which is the two-dimensional generalization of the construction of infinite-volume one-dimensional Gibbs measure from Perron-Frobenius eigenvectors and Kolmogorov's extension theorem.

\begin{theo}[invariant measures from boundary ROPE eigen-structure, from local to global]\label{theo:eigenROPErep:invmeas}
Let $S_1$ and $S_2$ be two finite sets and let 
$\MarkovWeight{W} : S_1\times S_1\times S_2 \times S_2 \to \setR_+$
be a fixed non-zero face weight, seen as an element of the space $T_{1,1}(V(S_1),V(S_2))=\ca{T}_{1,1}$ of the canonical $\Guill_2^{(\patterntype{fp}^*)}$-algebra $\ca{T}_{\PatternShapes(\patterntype{fp}^*)}$. Let $(g_{p,q})_{(p,q)\in\setL^*\times\setL^*}$ be a collection of non-identically null functions such that, for all $p,q\in\setL^*$, we have
\[
g_{p,q} : S_1^p\times S_1^p\times S_2^q\times S_2^q \to \setR_+
\]
(i.e. $g_{p,q}$ may be a boundary weight on a rectangle of size $(p,q)$).

If there exists a ROPE $\ca{B}_{\PatternShapes(\patterntype{fp}^*)}$ and ROPE elements $(A_S,A_N,A_W,A_E,U_{SW},U_{SE},U_{NW},U_{NE})$ of $(g_{p,q})_{(p,q)\in{\setL^*}^2}$ which satisfies the following conditions
\begin{enumerate}[(i)]
\item \label{item:hyp:roperepform}for all $(p,q)\in\setN^2$, $g_{p,q}$ admits the following ROPErep \begin{equation}\label{eq:boundaryweight:infiniteGibbs}
g_{p,q}(x,y,w,z)=
\begin{tikzpicture}[guillpart,yscale=1.55,xscale=2.8]
\draw[guillsep, dotted] (0,0) rectangle (6,5);
\draw[guillsep] (1,1) rectangle (5,4);
\draw[guillsep] 	(1,0)--(1,1)
				(2,0)--(2,1)
				(3,0)--(3,1)
				(4,0)--(4,1)
				(5,0)--(5,1);
\draw[guillsep] 	(1,5)--(1,4)
				(2,5)--(2,4)
				(3,5)--(3,4)
				(4,5)--(4,4)
				(5,5)--(5,4);
\draw[guillsep] 	(0,1)--(1,1)
				(0,2)--(1,2)
				(0,3)--(1,3)
				(0,4)--(1,4);
\draw[guillsep] 	(5,1)--(6,1)
				(5,2)--(6,2)
				(5,3)--(6,3)
				(5,4)--(6,4);
\node at (0.5,0.5) { $U_{SW}$ };
\node at (5.5,0.5) { $U_{SE}$ };
\node at (0.5,4.5) { $U_{NW}$ };
\node at (5.5,4.5) { $U_{NE}$ };
\node at (1.5,0.5) { $A_S(x_1)$ };
\node at (2.5,0.5) { $A_S(x_2)$ };
\node at (3.5,0.5) { $\ldots$ };
\node at (4.5,0.5) { $A_S(x_p)$ };
\node at (1.5,4.5) { $A_N(y_1)$ };
\node at (2.5,4.5) { $A_N(y_2)$ };
\node at (3.5,4.5) { $\ldots$ };
\node at (4.5,4.5) { $A_N(y_p)$ };
\node at (0.5,1.5) { $A_W(w_1)$ };
\node at (0.5,2.5) { $\vdots$ };
\node at (0.5,3.5) { $A_W(w_q)$ };
\node at (5.5,1.5) { $A_E(z_1)$ };
\node at (5.5,2.5) { $\vdots$ };
\node at (5.5,3.5) { $A_E(z_q)$ };
\end{tikzpicture}
\end{equation}

\item \label{item:hyp:boundaryeigenelmts} the four 1D-semi-groups $(\ha{A}^S_\bullet,\ha{A}^N_\bullet,\ha{A}^W_\bullet,\ha{A}^E_\bullet)$ generated by
\begin{subequations}
	\label{eq:ROPEreptosemigroups}
\begingroup
\allowdisplaybreaks	
\begin{align}
\ha{A}^S_1 &= \sum_{x\in S_1} e^*_x \otimes 1_\setK \otimes A_S(x) \in \ca{E}_{1,\infty_S}
\\
\ha{A}^N_1 &= \sum_{y\in S_1} e_y \otimes 1_\setK \otimes A_N(y) \in \ca{E}_{1,\infty_N}
\\
\ha{A}^W_1 &= \sum_{w\in S_2} 1_\setK\otimes e^*_w \otimes  A_W(w) \in \ca{E}_{\infty_W,1}
\\
\ha{A}^E_1 &= \sum_{z\in S_2} 1_\setK \otimes e_z  \otimes A_E(x) \in \ca{E}_{\infty_E,1}
\end{align}
\endgroup
and the corner elements
\begin{align}
\ha{U}_{a,b} &=  1_\setK  \otimes 1_\setK \otimes U_{a,b} \in \ca{E}_{\infty_b,\infty_a}
\end{align}
\end{subequations}
(see the canonical inclusions~\eqref{eq:fromROPEreptoEpq})
are such that the four triplets $(\ha{A}^S_\bullet,\ha{A}^W_\bullet,\setK U_{SW})$, 
$(\ha{A}^S_\bullet,\ha{A}^E_\bullet,\setK U_{SE})$, 
$(\ha{A}^N_\bullet,\ha{A}^W_\bullet,\setK U_{NW})$ and
$(\ha{A}^N_\bullet,\ha{A}^E_\bullet,\setK U_{NE})$ are four \emph{compatible completed corner eigen-spaces} of the 2D-semi-group generated by $\MarkovWeight{W}$ extended to $\ca{E}$ with eigenvalue sequences $\Lambda_\bullet$, $\sigma^S_\bullet$, $\sigma^N_\bullet$, $\sigma^W_\bullet$ and $\sigma^E_\bullet$ (see definitions~\ref{def:cornereigensemigroups} and \ref{def:compatiblecorners:twobytwo}),
\end{enumerate}
then there exists a two-dimensional Markov process $(X_e)_{e\in\Edges{\setZ^2}}$ on the full lattice $\setZ^2$ such that:
\begin{enumerate}[(a)]
\item $X$ is invariant in law by translation and all its restrictions to rectangles are Markov processes with face weights all equal to $\MarkovWeight{W}$;
\item for all $(p,q)\in\setN^*$, for any rectangle $R$ in $\setZ^2$ with shape $(p,q)$, the restriction $(X_e)_{e\in\Edges{R}}$ admits the function $g_{p,q}$ given by \eqref{eq:boundaryweight:infiniteGibbs} as boundary weight`, with a partition function given by:
\begin{equation}\label{eq:exactpartitionfuncrectangle}
Z_R^{\boundaryweights}(\MarkovWeight{W},g_{p,q}) = \kappa \sigma^S_p\sigma^N_p\sigma^W_q\sigma^E_q  \Lambda_{pq}
\end{equation}
(all sequences are geometric sequences) with a free energy density exactly equal to $\log \Lambda_1$.
\end{enumerate}
\end{theo}

\begin{proof}
Existence of the process on $\setZ^2$ is a direct consequence of Kolmogorov's extension theorem, with an index set $I=\Edges{\setZ^2}$ and a finite space $S=S_1\sqcup S_2$ with its complete $\sigma$-algebra. Any finite set $J\subset I$ can be embedded in a rectangle and it is enough to check the consistency condition on rectangles. We introduce the set of rectangles $[x_1,x_2]\times [y_1,y_2]$ in $\setZ^2$ with the partial order given by set inclusion. At the level of processes, there are, for any rectangle $R'\subset R$, natural restrictions $\pi_{R',R}( (x_e)_{e\in\Edges{R}} ) = (x_e)_{e\in\Edges{R'}}$ and the purpose is to construct laws $\mu_R$ of processes on every rectangle $R$ that satisfy the natural projections $\mu_{R'} = \mu_R \circ \pi_{R',R'}^{-1}$ of push-forward measures.

\paragraph*{Marginal law on a fixed rectangle of size \texorpdfstring{$(p,q)$}{(p,q)}.}

We consider now a fixed rectangle $R=[x_1,x_2]\times[y_1,y_2]$ with shape $(p,q)$ and consider the Markov process $(X^{(R)}_e)_{e\in \Edges{R}}$ with all face weights equal to $\MarkovWeight{W}$ and boundary weight $g_{p,q}$. This is made possible by the positivity condition on the ROPEreps of the functions $g_{p,q}$. Its probability law~\eqref{eq:MarkovLawDimTwo} contains a partition function given by \eqref{eq:proba:ZboundaryfromZdet}.

Using Proposition~\ref{prop:partitionfunc:homogeneous} together with the definition of the global $\Guill_2^{(\patterntype{fp}^*)}$-algebra $\ca{E}_{\PatternShapes(\patterntype{fp}^*)}$ in Lemma~\ref{lemm:globalGuillalgebraEpq} and the embeddings~\eqref{eq:ROPEreptosemigroups}, we obtain a partition function 
on $R$ given by \eqref{eq:fullplane:eigen} where the 2D- and 1D-semigroups are replaced by $\MarkovWeight{W}^{[p,q]}$ and the $\ha{A}^a_\bullet$. Corollary~\ref{coro:fullplane:eigen} then provides
\[
Z_R^{\boundaryweights}(\MarkovWeight{W},g_{p,q}) = \kappa \Lambda_{pq} \sigma_p^S \sigma^N_p \sigma^W_q\sigma^E_q
\]

\paragraph*{From a rectangle to a smaller one: marginal law and simplification through morphisms.}
Theorem~\ref{theo:stability} and its proof shows how a ROPErep of a boundary weight $g_{p,q}$ of a process $(X^{(R)}_e)_{e\in\Edges{R}}$ on a rectangle $R$ with shape $(p,q)$ provides a ROPErep \emph{on a different ROPE} of the boundary weight $g^{(p,q)}_{p',q'}$ of the law of the process $(X^{(R)}_e)_{e\in\Edges{R'}}$ restricted to a smaller rectangle $R'\subset R$ with shape $(p',q')$. We will now use the morphisms introduced with the boundary eigen-spaces to provide a second ROPErep of $g^{(p,q)}_{p',q'}$ on the \emph{same} ROPE as $g_{p,q}$ and then observe that $g^{(p,q)}_{p',q'}$ is proportional to $g_{p,q}$.

As in the proof of Theorem~\ref{theo:stability}, we characterize $g^{(p,q)}_{p',q'}$ by its action an arbitrary element of $Y\in\ca{E}_{p',q'}=\ca{T}_{p',q'}$ following \eqref{eq:bwrestrictedaction} (base point dropped for easier notation):
\[
\langle g^{(p,q)}_{p',q'},Y \rangle = 
\begin{tikzpicture}[guillpart,yscale=1.5,xscale=1.95]
	\fill[guillfill] (0,0) rectangle (5,5);
	\draw[guillsep] (1,0)--(1,5)  (2,0)--(2,5) (3,0)--(3,5) (4,0)--(4,5);
	\draw[guillsep] (0,1)--(5,1) (0,2)--(5,2) (0,3)--(5,3) (0,4)--(5,4);
	\node at (0.5,0.5) {$\ha{U}_{SW}$};
	\node at (1.5,0.5) {$\ha{A}^S_{p_1}$};
	\node at (2.5,0.5) {$\ha{A}^S_{p'}$};
	\node at (3.5,0.5) {$\ha{A}^S_{p_2}$};
	\node at (4.5,0.5) {$\ha{U}_{SE}$};
	\node at (0.5,1.5) {$\ha{A}_{q_1}^W$};
	\node at (1.5,1.5) {$\MarkovWeight{F}_{p_1,q_1}$};
	\node at (2.5,1.5) {$\MarkovWeight{F}_{p',q_1}$};
	\node at (3.5,1.5) {$\MarkovWeight{F}_{p_2,q_1}$};
	\node at (4.5,1.5) {$\ha{A}^E_{q_1}$};
	\node at (0.5,2.5) {$\ha{A}^W_{q'}$};
	\node at (1.5,2.5) {$\MarkovWeight{F}_{p_1,q'}$};
	\node at (2.5,2.5) {$Y\otimes 1$};
	\node at (3.5,2.5) {$\MarkovWeight{F}_{p_2,q'}$};
	\node at (4.5,2.5) {$\ha{A}^E_{q'}$};
	\node at (0.5,3.5) {$\ha{A}^W_{q_2}$};
	\node at (1.5,3.5) {$\MarkovWeight{F}_{p_1,q_2}$};
	\node at (2.5,3.5) {$\MarkovWeight{F}_{p',q_2}$};
	\node at (3.5,3.5) {$\MarkovWeight{F}_{p_2,q_2}$};
	\node at (4.5,3.5) {$\ha{A}^E_{q_2}$};
	\node at (0.5,4.5) {$\ha{U}_{NW}$};
	\node at (1.5,4.5) {$\ha{A}^N_{p_1}$};
	\node at (2.5,4.5) {$\ha{A}^N_{p'}$};
	\node at (3.5,4.5) {$\ha{A}^N_{p_2}$};
	\node at (4.5,4.5) {$\ha{U}_{NE}$};
\end{tikzpicture}
\]
In the present case, $\ca{E}_{\infty_{WE},\infty_{SN}}=\setR$ and thus, for any $a\in\setR$, the morphisms and eigenvalue sequences introduced in proposition~\eqref{prop:compatiblecorners:full} may be chosen so that
\[
\Phi_{\infty_W,\infty_{SN}}^{W,p_1,(\infty_E,x),(q,y)}\circ \Phi_{\infty_E,\infty_{SN}}^{E,p_2,(\infty_W,x-r),(q,y)} = \id
\]
(we have seen in Section~\ref{sec:onedimdiagowithmorph} that normalizations of the morphisms and eigenvalues are related to each other and we may use any gauge). Using the assumptions of compatible eigen corner spaces of the ROPErep of $g_{p,q}$, we may apply these horizontal morphisms to 
$\langle g^{(p,q)}_{p',q'} , Y \rangle $ and obtain from the previous equations:
\[
\langle g^{(p,q)}_{p',q'},Y \rangle = \Lambda^{(p_1+p_2)q}\sigma^N_{p_1+p_2}\sigma^S_{p_1+p_2}
\begin{tikzpicture}[guillpart,yscale=1.5,xscale=1.95]
	\fill[guillfill] (1,0) rectangle (4,5);
	\draw[guillsep]   (2,0)--(2,5) (3,0)--(3,5) ;
	\draw[guillsep] (1,1)--(4,1) (1,2)--(4,2) (1,3)--(4,3) (1,4)--(4,4);
	\node at (1.5,0.5) {$\ha{U}_{SW}$};
	\node at (2.5,0.5) {$\ha{A}^S_{p'}$};
	\node at (3.5,0.5) {$\ha{U}_{SE}$};
	\node at (1.5,1.5) {$\ha{A}_{q_1}^W$};
	\node at (2.5,1.5) {$\MarkovWeight{F}_{p',q_1}$};
	\node at (3.5,1.5) {$\ha{A}^E_{q_1}$};
	\node at (1.5,2.5) {$\ha{A}^W_{q'}$};
	\node at (2.5,2.5) {$Y\otimes 1$};
	\node at (3.5,2.5) {$\ha{A}^E_{q'}$};
	\node at (1.5,3.5) {$\ha{A}^W_{q_2}$};
	\node at (2.5,3.5) {$\MarkovWeight{F}_{p',q_2}$};
	\node at (3.5,3.5) {$\ha{A}^E_{q_2}$};
	\node at (1.5,4.5) {$\ha{U}_{NW}$};
	\node at (2.5,4.5) {$\ha{A}^N_{p'}$};
	\node at (3.5,4.5) {$\ha{U}_{NE}$};
\end{tikzpicture}
\]
Using again Proposition~\ref{prop:compatiblecorners:full} with vertical morphisms, we obtain
\[
\langle g^{(p,q)}_{p',q'},Y \rangle = \Lambda^{(p_1+p_2)q+(q_1+q_2)p'}\sigma^N_{p_1+p_2}\sigma^S_{p_1+p_2}\sigma^W_{q_1+q_2}\sigma^E_{q_1+q_2}
\begin{tikzpicture}[guillpart,yscale=1.5,xscale=1.95]
	\fill[guillfill] (1,1) rectangle (4,4);
	\draw[guillsep]   (2,1)--(2,4) (3,1)--(3,4) ;
	\draw[guillsep]  (1,2)--(4,2) (1,3)--(4,3) ;
	\node at (1.5,1.5) {$\ha{U}_{SW}$};
	\node at (2.5,1.5) {$\ha{A}^S_{p'}$};
	\node at (3.5,1.5) {$\ha{U}_{SE}$};
	\node at (1.5,2.5) {$\ha{A}^W_{q'}$};
	\node at (2.5,2.5) {$Y\otimes 1$};
	\node at (3.5,2.5) {$\ha{A}^E_{q'}$};
	\node at (1.5,3.5) {$\ha{U}_{NW}$};
	\node at (2.5,3.5) {$\ha{A}^N_{p'}$};
	\node at (3.5,3.5) {$\ha{U}_{NE}$};
\end{tikzpicture}
\]
which  identifies $g^{(p,q)}_{p',q'}$ to a multiple of $g_{p',q'}$ The probability law of the restricted process $(X_e)_{e\in\Edges{R'}}$ now contains a ratio
\[
\frac{g^{(p,q)}_{p',q'}}{Z^{\boundaryweights}_{R'}(\MarkovWeight{W},g_{p,q})} = \frac{\Lambda^{(p_1+p_2)q+(q_1+q_2)p'}\sigma^N_{p_1+p_2}\sigma^S_{p_1+p_2}\sigma^W_{q_1+q_2}\sigma^E_{q_1+q_2}}{\kappa \Lambda_{pq} \sigma_p^S \sigma^N_p \sigma^W_q\sigma^E_q}  g_{p,q}
\]
from our previous computation of the partition function. Using the geometric sequence property of the eigenvalue sequences, we now obtain directly
\[
\frac{\Lambda^{(p_1+p_2)q+(q_1+q_2)p'}\sigma^N_{p_1+p_2}\sigma^S_{p_1+p_2}\sigma^W_{q_1+q_2}\sigma^E_{q_1+q_2}}{\kappa \Lambda_{pq} \sigma_p^S \sigma^N_p \sigma^W_q\sigma^E_q}  g_{p,q}
= \frac{1}{\kappa \Lambda_{p'q'} \sigma_{p'}^S \sigma^N_{p'} \sigma^W_{q'}\sigma^E_{q'}} = \frac{1}{Z^{\boundaryweights}_{R'}(\MarkovWeight{W},g_{p,q})}
\]
Thus, the law of the restricted process $(X^{(R)}_{e})_{e\in\Edges{R'}}$ is equal to the law of the process $(X^{(R')}_e)_{e\in\Edges{R'}}$ directly defined out of $\MarkovWeight{W}$ and $g_{p,q}$. A direct application of Kolmogorov's extension on this consistent set of marginal laws on rectangles provides the existence of a global process on $\setZ^2$. 

Invariance in law by translation is a direct consequence of the homogeneity of the face weights and the fact that boundary weights depend only the shape of the rectangles and not their positions. Properties of the asymptotic free energy density is a direct consequence of the geometric nature of the eigenvalue sequences.

\end{proof}

\subsubsection{A simpler description of infinite-volume Gibbs measures.} One of the main difficulties in the construction of Gibbs measure is a simple description of the complete probability laws with the suitable restriction properties. As explained in \cite{GibbsVelenik}, both of the functional analysis approach and the DLR approach avoid the exact description of the boundary weights $g_R$, precisely because their determination is difficult. The present approach circumvents these difficulties and the only remaining hard point is the determination of the boundary eigen-structure.

\subsubsection{From local to global.} Theorem~\ref{theo:eigenROPErep:invmeas} induces a change of point of view on boundary weights. In dimension one, such a discussion is trivial: a boundary weight is a function on values above only two (maybe far away) points, which can be treated separately using a single "local" eigenvalue equation $Xu=\lambda u$ or $vX=\lambda v$. By the word of "local", we mean that it uses only the first power of $X$ and not $X^n$ where $n$ is the distance between the two boundary points. 

In dimension two however, the boundary of a rectangle is a one-dimensional curve: a boundary weight is a global object above the curve. When sticking to rectangles, reduction from a rectangle of sizes $(P,Q)$ to a smaller rectangle of size $(P-1,Q)$ (or $(P,Q-1)$) requires to use the matrix $\MarkovWeight{W}^{[1,Q]}$ (or $\MarkovWeight{W}^{[P,1]}$), which may be a large "global" object already hard to compute, and then solve an eigenvalue problem on a large matrix. Previous Theorem~\ref{theo:eigenROPErep:invmeas} solves this issue in two steps:
\begin{itemize}
\item first, the global object $g_{p,q}$ is built out of a finite number of local objects: the elements of the ROPE representation;
\item the elements of the ROPE representation solve  eigenvalue-like equations up to morphisms following the definitions of Section~\ref{sec:invariantboundaryelmts} on half-strips, corners and their compatibilities, which appear to be local equations (in the sense that they do not require the computation of any (large or not) surface powers $\MarkovWeight{W}^{[P,Q]}$.
\end{itemize}

%SIXVERTEX his approach will be performed in detail in a non-trivial case in Section~\ref{sec:sixvertex}.

\subsubsection{Some questions between probability and algebra} 

In dimension one, existence of Perron-Frobenius eigenvectors used to extend processes on the whole line can be obtained both with algebraic tools (cones) or probabilistic ones (excursion theory). It would be interesting to see how these echoing constructions can be extended in dimension two and larger.

From the point of view of probability theory, the reduction of dimension between the ROPEreps from Theorems~\ref{theo:stability} and \ref{theo:eigenROPErep:invmeas} through morphisms can be interpreted as some erasing of information along the transverse direction: we do not have yet nice examples to illustrate such constructions but it would be an important question in order to develop "higher" diagonalization tools.

\subsection{The simplex property for ROPEreps of invariant boundary weights.}

It is a well-known fact in the theory of Gibbs measures that, for a given model, the set of infinite-volume Gibbs measures is a simplex, i.e. any of them is a barycentre of extremal ones. Despite the fact that we do not know whether any infinite-volume Gibbs measure may arise from ROPEreps made of full boundary eigen-elements, it is interesting to check that this simplex property still holds in our case.

\begin{prop}[simplex property]
Let $\MarkovWeight{W}\in\ca{T}_{1,1}$ be a face weight. Let $\ca{B}_{\PatternShapes(\patterntype{fp}^*)}$ and $\ca{B}'_{\PatternShapes(\patterntype{fp}^*)}$ be two ROPE and two ROPEreps, respectively in $\ca{B}_{\PatternShapes(\patterntype{fp}^*)}$ and $\ca{B}'_{\PatternShapes(\patterntype{fp}^*)}$,
\begin{align*}
((A_S(x))_{x\in S_1},(A_N(y))_{y\in S_1},(A_W(w))_{w\in S_2},(A_E(z))_{z\in S_2},U_{SW},U_{NW},U_{SE},U_{NE})
\\
((A'_S(x))_{x\in S_1},(A'_N(y))_{y\in S_1},(A'_W(w))_{w\in S_2},(A'_E(z))_{z\in S_2},U'_{SW},U'_{NW},U'_{SE},U'_{NE})
\end{align*}
that satisfy the hypotheses \eqref{item:hyp:roperepform} and \eqref{item:hyp:boundaryeigenelmts} of Theorem~\ref{theo:eigenROPErep:invmeas} with the \emph{same} eigenvalue sequences.
There exists a ROPE $\ov{\ca{B}}_{\PatternShapes(\patterntype{fp}^*)}$ such that, for any 
$\alpha\in [0,1]$, there exists a ROPErep on $\ov{\ca{B}}_{\PatternShapes(\patterntype{fp}^*)}$
\[
((\ov{A}_S(x))_{x\in S_1},(\ov{A}_N(y))_{y\in S_1},(\ov{A}_W(w))_{w\in S_2},(\ov{A}_E(z))_{z\in S_2},\ov{U}_{SW},\ov{U}_{NW},\ov{U}_{SE},\ov{U}_{NE})
\]
which satisfies hypotheses \eqref{item:hyp:roperepform} and \eqref{item:hyp:boundaryeigenelmts} of Theorem~\ref{theo:eigenROPErep:invmeas} with the same eigenvalue sequences and such that the boundary weights built out of them satisfy on any rectangle
\begin{equation}\label{eq:simplexpropforboundaryweights}
\ov{g}_{p,q} =  \alpha g_{p,q} + (1-\alpha) g'_{p,q}
\end{equation}
The infinite-volume Gibbs measures built out of them satisfy also the same barycentric property.
\end{prop}
This is a straightforward consequence of the properties described in Section~\ref{sec:combinationofROPEreps}

\section{Fixed point up to morphisms}\label{sec:fixedpoints}

In the context of Markov processes, algebras over guillotine operads with a vector space structure and tensor products are directly relevant. However, as it is explained at the beginning of Section~\ref{sec:operadiclanguage}, algebras over guillotine operads can be formulated with only sets and Cartesian products (and we call it instead a Cartesian guillotine algebra"). In this simpler case, there is no action of the scalar field $\setK$ and no possibility of eigenvalue. However, it is possible to formulate definitions similar to the ones of eigen-elements on half-strips and corners. In this case, we replace the expression "\emph{eigen-element}/\emph{eigen-generators} up to morphisms" by "fixed point up to morphisms". This notion is particularly interesting in the case of Gaussian fields as in Section~\ref{sec:appli:gaussian} and \cite{BodiotSimon} and whenever the face weights are restricted to a family of weights parametrized by (maybe non-linear) objects with an internal structure of guillotine set algebra.

We sketch here how it works without rewriting all the definitions. The notion 2D-semi-group $F_{\bullet,\bullet}$ and the associated morphisms $\psi_p^{F,a,\bullet,\bullet}$ do not require any additional linear structure and all their operadic properties remain the same (Lemma~\ref{lemm:morphismsfrom2Dsemigroup}, Proposition~\ref{prop:towerofguillalg}). There is now no geometric sequence nor scaling maps. It is then possible to slightly shift Definition~\ref{def:eigenalgebra:single} by simply removing any mention of the eigenvalue geometric sequence and replacing tensor products by Cartesian products: all the operadic structure is kept unchanged. We can then proceed forward and shift in the same way all the following definitions. As an example, we provide an analogue of Definition~\ref{def:eigenalgebrauptomorphims} for fixed points and emphasize on the changes.

\begin{defi}[$\Guill_2^{(\patterntype{hs}_S)}$-fixed points up to morphisms of a 2D-semi-group]\label{def:fixedpointsouthuptomorphims}
	Let $\ca{A}_{\PatternShapes(\patterntype{hs}_S)}$ be a \emph{Cartesian} $\Guill_2^{(\patterntype{hs}_S)}$-algebra and $\MarkovWeight{F}_{\bullet,\bullet}$ a 2D-semi-group in $\ca{A}_{\PatternShapes(\patterntype{r})}$. A collection of \emph{subsets} $(\ca{V}_p)_{p\in\setL^*}$ of $(\ca{A}_{p,\infty_S})_{p\in\setL^*}$ is a system of \emph{$\Guill_1$-fixed points up to morphisms} of $\MarkovWeight{F}_{\bullet,\bullet}$ if and only if there exists a collection of maps $\Phi_p^{S,q,r} : \ca{A}_{p,\infty_S} \to \ca{A}_{p,\infty_S}$, $p,q,r\in\setL^*$ \emph{(no linearity imposed)} such that, for any $N\in\setN^*$, for any sizes $q$ and $r$, any sequence of sizes $(p_i)_{1\leq i\leq N}$ in $\setL^*$, any sequence of boundary elements $(v_i)_{1\leq i\leq N}$ with $v_i\in\ca{V}_{p_i,\infty_S}$, for any element $\MarkovWeight{T}\in\ca{A}_{P,r}$,
		\begin{equation}\label{eq:deffixedSouth:withM}
			\Phi_{\sum_i p_i}^{S,q,r}\left(
			\begin{tikzpicture}[guillpart,yscale=1,xscale=2.]
				\fill[guillfill] (0,0) rectangle (3,3);
				\draw[guillsep] (0,0)--(0,3)--(3,3)--(3,0) (1,0)--(1,2) (2,0)--(2,2) (0,1)--(3,1) (0,2)--(3,2);
				\node at (0.5,0.5) {$v_1$};
				\node at (1.5,0.5) {$\dots$};
				\node at (2.5,0.5) {$v_N$};
				\node at (0.5,1.5) {$\MarkovWeight{F}_{p_1,q}$};
				\node at (1.5,1.5) {$\dots$};
				\node at (2.5,1.5) {$\MarkovWeight{F}_{p_N,q}$};
				\node at (1.5,2.5) {$\MarkovWeight{T}$};
			\end{tikzpicture}
			\right)
			= \begin{tikzpicture}[guillpart,xscale=1.5]
				\fill[guillfill] (0,0) rectangle (3,2);
				\draw[guillsep] (0,0)--(0,2)--(3,2)--(3,0) (1,0)--(1,1) (2,0)--(2,1) (0,1)--(3,1);
				\node at (0.5,0.5) {$v_1$};
				\node at (1.5,0.5) {$\dots$};
				\node at (2.5,0.5) {$v_N$};
				\node at (1.5,1.5) {$\MarkovWeight{T}$};
			\end{tikzpicture}
		\end{equation}
\end{defi}
One indeed checks that the linear structure and the tensor products do not appear in this operadic definition. Lemma~\ref{lemm:eigenalg:sufficientmorphisms} purged from tensor products, eigenvalues and with scaling maps replaced by identities still hold.

Similar adaptations of definitions~\ref{def:eigenalgebrauptomorph:leftextended}, \ref{def:cornereigensemigroups}, \ref{def:extensioneigenalgtostrip} can be made by removing all eigenvalue sequences and scalar multiplication maps and propositions~\ref{prop:simplificationonstrips}, \ref{prop:NWandSWcorners:gluing:I} remains valid without eigenvalues. Compatibility between corners as in Definition~\ref{def:compatiblecorners:twobytwo} and Proposition~\ref{prop:compatiblecorners:full} are changed in the same way to describe compatible corner fixed points by keeping morphisms and removing eigenvalues.

We use this definitions of fixed points in order to simplify computations in Section~\ref{sec:appli:gaussian} but do not know whether it may have an interest in the theory of dynamical systems (see for example~\cite{SimonMixedRadix} for a simple model).

\chapter{Applications and concrete examples}\label{sec:applications}
\section{From abstract definitions down to practical computations}

The previous section has a strong operadic content, with the system of colours and morphisms. The purpose of the present section is to transform the previous results into a practical scheme of computations. Before entering models, we first list and count the unknown variables and their constraints. 

A boundary weight $g_R$ on a non-degenerate rectangle $R$ of size $(p,q)$ is a function $S_1^{2p}\times S_2^{2q}\to \setR_+$, thus has $|S_1|^{2p} |S_2|^{2q}$ parameters but is defined only up to multiplicative constant (which may be absorbed in the partition function $Z_R$). The constraints required for Kolmogorov's extension theorem in the translation invariant setting correspond to the erasure of a leftmost or rightmost column of the rectangle and the erasure of the topmost and bottommost lines of the rectangle. All of them takes the form of a linear relation with a rectangular matrix of size $|S_1|^{2p}|S_2|^{2q}\times |S_1|^{2p+2}|S_2|^{2q}$ such as
\[
g_{p,q}\left( 
\begin{tikzpicture}[scale=0.5,baseline={(current bounding box.center)}]
	\draw (0,0) rectangle (4,3);
	\node at (0.5,0) [anchor = north] {$x_1$};
	\node at (1.5,0) [anchor = north]{$x_2$};
	\node at (2.5,0) [anchor = north]{$\ldots$};
	\node at (3.5,0) [anchor = north]{$x_p$};
	\node at (0.5,3) [anchor = south] {$y_1$};
	\node at (1.5,3) [anchor = south]{$y_2$};
	\node at (2.5,3) [above]{$\ldots$};
	\node at (3.5,3) [above]{$y_p$};
	\node at (0,0.5) [anchor = east]{ $w_1$ };
	\node at (0,1.5) [anchor = east]{ $\vdots$ };
	\node at (0,2.5) [anchor = east]{ $w_q$ };
	\node at (4,0.5) [anchor = west]{ $z_1$ };
	\node at (4,1.5) [anchor = west]{ $\vdots$ };
	\node at (4,2.5) [anchor = west]{ $z_q$ };
\end{tikzpicture}
\right)  = \sum_{(x,y,z'_1,\ldots,z'_q)}
g_{p+1,q}\left( 
\begin{tikzpicture}[scale=0.5,baseline={(current bounding box.center)}]
	\draw (0,0) rectangle (5,3);
	\node at (0.5,0) [anchor = north] {$x_1$};
	\node at (1.5,0) [anchor = north]{$x_2$};
	\node at (2.5,0) [anchor = north]{$\ldots$};
	\node at (3.5,0) [anchor = north]{$x_p$};
	\node at (4.5,0) [anchor = north]{$x$};
	\node at (0.5,3) [anchor = south] {$y_1$};
	\node at (1.5,3) [anchor = south]{$y_2$};
	\node at (2.5,3) [above]{$\ldots$};
	\node at (3.5,3) [above]{$y_p$};
	\node at (4.5,3) [above]{$y$};
	\node at (0,0.5) [anchor = east]{ $w_1$ };
	\node at (0,1.5) [anchor = east]{ $\vdots$ };
	\node at (0,2.5) [anchor = east]{ $w_q$ };
	\node at (5,0.5) [anchor = west]{ $z'_1$ };
	\node at (5,1.5) [anchor = west]{ $\vdots$ };
	\node at (5,2.5) [anchor = west]{ $z'_q$ };
\end{tikzpicture}
\right)
\begin{tikzpicture}[guillpart]
	\fill[guillfill] (0,0) rectangle (1,3);
	\draw[guillsep] (0,0) rectangle (1,3);
	\draw[guillsep] (0,1) -- (1,1);
	\draw[guillsep] (0,2) -- (1,2);
	\node at (0.5,0.5) {$\MarkovWeight{W}$};
	\node at (0.5,1.5) {$\vdots$};
	\node at (0.5,2.5) {$\MarkovWeight{W}$};
	\node at (0.5,0) [anchor=north] {$x$};
	\node at (0.5,3) [anchor=south] {$y$};
	\node at (0,0.5) [anchor=east] {$z_1$};
	\node at (1,0.5) [anchor=west] {$z'_1$};
	\node at (0,1.5) [anchor=east] {$\vdots$};
	\node at (1,1.5) [anchor=west] {$\vdots$};
	\node at (0,2.5) [anchor=east] {$z_q$};
	\node at (1,2.5) [anchor=west] {$z'_q$};
\end{tikzpicture}
\]
There is a countable set of objects $(g_{p,q})$ with these equations. A first difficulty is to compute the rectangular matrix on the r.h.s. from the weight $\MarkovWeight{W}$ for arbitrary large $q$. A second difficulty, which is the major one, consists in guessing the solutions $(g_{p,q})$ and this is mainly the reason why Kolmogorov's extension theorem has not been used previously in dimension strictly larger than one.

We now assume that the family of boundary weights $g_{p,q}$ have a ROPErep such as \eqref{eq:boundaryweight:infiniteGibbs}, following notations of Theorem~\ref{theo:eigenROPErep:invmeas}. The previous rectangular consistency conditions are replaced in the formalism of this theorem by the eigen-element conditions. In the ROPErep, there is a \emph{finite} number of (generically) \emph{infinite-dimensional} objects, that we can list:
\begin{itemize}
	\item the $|S_1|$ elements $A_S(x)$, for $x\in S_1$
	\item the $|S_1|$ elements $A_N(y)$, for $y\in S_1$
	\item the $|S_2|$ elements $A_W(w)$, for $w\in S_2$
	\item the $|S_2|$ elements $A_E(z)$, for $z\in S_2$
	\item the four corner elements $U_{ab}$, for $a\in\{S,N\}$ and $b\in\{W,E\}$.
\end{itemize}
Once again, for a probabilistic use, they are defined only up to multiplicative constants that can be absorbed in the eigenvalue sequences $\sigma_\bullet^a$. On this finite set of objects, there is also a \emph{finite} set of constraints:
\begin{itemize}
	\item $|S_1|$ equations coupling linearly the operators $A_S(\bullet)$ of the type \eqref{eq:reducedmorphism:eigeneq} 
	\item $|S_1|$ equations coupling linearly the operators $A_N(\bullet)$ of the type \eqref{eq:reducedmorphism:eigeneq} 
	\item $|S_2|$ equations coupling linearly the operators $A_W(\bullet)$ of the type \eqref{eq:reducedmorphism:eigeneq} 
	\item $|S_2|$ equations coupling linearly the operators $A_W(\bullet)$ of the type \eqref{eq:reducedmorphism:eigeneq} 
	\item for each corner $ab$, a linear equation on $U_{ab}$ whose coefficients involve the operators $A_a(\bullet)$ and $A_b(\bullet)$ of type \eqref{eq:concreteeigencorner}.
\end{itemize}
All these equations are eigenvector-like equations with additional morphisms and described below for some models. In all these equations, there is only \emph{one} occurrence of the face weight $\MarkovWeight{W}$ and no computation of surface powers nor transfer matrices.

The major novelty is that the ROPE spaces $\ca{B}_\bullet$ and the morphisms $\Phi_\bullet^\bullet$ in these equations are also \emph{unknown}. As spaces and not elements of predefined spaces, they are not submitted to concrete equations but rather to the structural constraints through the morphism requirements. Moreover, as already explained, these spaces are generically infinite-dimensional and defined up to ROPE morphisms. The key difficulty in the present approach is that we do not have yet enough mathematical tools to build canonically or algorithmically these spaces; furthermore these difficulties are also related to the fact that, in dimension strictly larger than one, there is a wide variety of non-trivial Gibbs measures with possible phase transitions. One can reasonably expect that, as in the analytical approach of infinite-volume Gibbs measures, physical insights as well as mathematical inequalities put sufficiently strong constraints on these spaces so that their structure can be partly of fully guessed.

These apparent difficulties are not obstacles at all in many cases. Indeed, the boundary spaces can be partly guessed from the structure of the particular model under study. We list here some scenarios.

\subsection{Rigorous validation of a candidate solution (trivial case, six-vertex model, etc.)}

This corresponds to the simplest case in which a ROPErep is already known from other considerations, as in the case of the trivial models of Section~\ref{sec:trivialfactorizedcase}. In this case, the ROPE spaces are known and one only has to insert the candidate solution in the generalized eigen-vector equations and to check the formulae.

This also corresponds to less trivial situations where the ROPErep is only \emph{conjectured} from non-rigorous computations. This is the case of integrable systems with Bethe Ansatz solutions: some computations are made for finite system sizes using Bethe Ansatz approach and one expects, in the large size limit, to have condensation of Bethe roots on a curve in the complex plane. Proving this condensation is hard and lengthy and is still often an open problem, making the full rigorous study of phase diagram incomplete, as in the case of the elementary six-vertex model (see \cite{DuminilBethe,DuminilCondensation}). The present framework allows one to avoid this non-rigorous path: one first \emph{assumes} condensation of Bethe roots and directly describe the ROPE spaces and the ROPEreps from Bethe Ansatz computations (see \cite{AlcarazLazo,GolinelliMallick}) and then a direct check of the generalized eigenvector-like equations \emph{validates rigorously} the solution and hence the existence of the infinite-volume Gibbs measure. As seen below in Section~\ref{sec:appli:correlfn}, one further gains explicit representations of correlation functions directly in the thermodynamic limit without any need for approximations and controls of convergence. All this approach and the results are postponed to a separate paper \cite{SimonSixV}.

\subsection{Computations within given spaces: the Gaussian case}

For other models such as the Gaussian model described below and studied in detail in \cite{BodiotSimon}, the internal structure of the model (the Hilbert space isometry in the Gaussian case) prescribe canonical boundary spaces with canonical morphisms and a canonical form for the laws suitably parametrized (covariance matrix in the Gaussian case). Once the ROPE $\ca{B}_\bullet$ and the morphisms $\Phi_\bullet^\bullet$ are known, then all equations above, of type \eqref{eq:reducedmorphism:eigeneq} or \eqref{eq:concreteeigencorner}  become classical eigenvector equation to be solved. The simplification can go further as in the Gaussian case: if the space of parameters of the laws has itself a guillotine algebra structure (not necessarily linear) then the previous equations are lifted to a finite set of equations on parameters, which can be solved exactly or at least numerically.

\section[Line marginals and correlation functions]{Some uses of the full boundary eigen-structure: line marginals and correlation functions.} \label{sec:appli:correlfn}

We now assume that a boundary eigen-structure is fixed so that Theorem~\ref{theo:eigenROPErep:invmeas} holds and an infinite-volume Gibbs measure exists and we use the same notations. We are now interested to the marginal law of the process on a given horizontal or vertical segment of length $L$. We treat here the case of an horizontal segment made of $L$ consecutive edges:
\begin{equation}
	\mathbf{S} = \{ ((x+k,y),(x+k+1,y)); 0\leq k <L \}
\end{equation}
A algebraic interesting way is to consider the null-size commutative spaces $\ca{T}_{0,q}$ and $\ca{T}_{p,0}$ introduced for observables.

\begin{theo}[one-dimensional marginals from invariant ROPEs]\label{theo:onedimmarginaloutofinvROPErep}
	Under the hypotheses and notations of Theorem~\ref{theo:eigenROPErep:invmeas}, the marginal law  of the restriction of the process to a segment $(X_e)_{e\in\mathbf{S}}$ under the infinite-volume Gibbs measure defined by the ROPE eigen-elements is given by:
	\begin{subequations}
		\label{eq:onedim:marginal}
		\begin{align}
			\prob{(X_e)_{e\in\mathbf{S}}=(x_e)_{e\in\mathbf{S}} }
			&= \frac{1}{\kappa \sigma^S_L\sigma^N_L}
			\begin{tikzpicture}[guillpart,xscale=3.,yscale=1.75]
				\fill[guillfill] (0,0) rectangle (6,2);
				\draw[guillsep] (0,1)--(6,1) (1,0)--(1,2) (2,0)--(2,2) (3,0)--(3,2) (4,0)--(4,2) (5,0)--(5,2);
				\node at (1,1) [circle, fill, inner sep=0.5mm] {};
				\node at (0.5,0.5) { $\ha{U}_{SW}$ };
				\node at (0.5,1.5) { $\ha{U}_{NW}$ };
				\node at (1.5,0.5) { $\ha{A}_{S}$ };
				\node at (1.5,1.5) { $\ha{A}_{N}$ };
				\node at (2.5,0.5) { $\ha{A}_{S}$ };
				\node at (2.5,1.5) { $\ha{A}_{N}$ };
				\node at (3.5,0.5) { $\ldots$ };
				\node at (3.5,1.5) { $\ldots$ };
				\node at (4.5,0.5) { $\ha{A}_{S}$ };
				\node at (4.5,1.5) { $\ha{A}_{N}$ };
				\node at (5.5,0.5) { $\ha{U}_{SE}$ };
				\node at (5.5,1.5) { $\ha{U}_{NE}$ };
				\node at (1.5,1.) {$\indic{x_{e_1}}$};
				\node at (2.5,1.) {$\indic{x_{e_2}}$};
				\node at (3.5,1.) {$\ldots$};
				\node at (4.5,1.) {$\indic{x_{e_L}}$};
			\end{tikzpicture}_{\ca{E}}\label{eq:onedimmarginal:Epq}
			\\
			&=
			\frac{1}{\kappa \sigma^S_L\sigma^N_L}\begin{tikzpicture}[guillpart,yscale=1.5,xscale=2.8]
				\draw[guillsep, dotted] (0,0) rectangle (6,2);
				\draw[guillsep] (1,0)--(1,2)
				(2,0)--(2,2)
				(3,0)--(3,2)
				(4,0)--(4,2)
				(5,0)--(5,2)
				(0,1)--(6,1);
				\node at (0.5,0.5) { $U_{SW}$ };
				\node at (5.5,0.5) { $U_{SE}$ };
				\node at (0.5,1.5) { $U_{NW}$ };
				\node at (5.5,1.5) { $U_{NE}$ };
				\node at (1.5,0.5) { $A_S(x_{e_1})$ };
				\node at (2.5,0.5) { $A_S(x_{e_2})$ };
				\node at (3.5,0.5) { $\ldots$ };
				\node at (4.5,0.5) { $A_S(x_{e_L})$ };
				\node at (1.5,1.5) { $A_N(x_{e_1})$ };
				\node at (2.5,1.5) { $A_N(x_{e_2})$ };
				\node at (3.5,1.5) { $\ldots$ };
				\node at (4.5,1.5) { $A_N(x_{e_L})$ };
			\end{tikzpicture}_{\ca{B}}
			\label{eq:onedimmarginal:ROPErep}
			\\
			&= \frac{1}{\kappa \sigma^S_L\sigma^N_L}\langle u_{W}, A_{SN}(x_1)\ldots A_{SN}(x_L) u_{E} \rangle
			\label{eq:onedimmarginal:matrixansatz}
		\end{align}
	\end{subequations}
	where the edges are enumerated from left to right in the segment $\mathbf{S}$. The elements of the matrix product form~\eqref{eq:onedimmarginal:matrixansatz} take value in the space $\Hom(\ca{B}_{\infty_E,\infty_{SN}},\ca{B}_{\infty_E,\infty_{SN}})$ of the ROPE
	\begin{align*}
		u_W &=
		\pi_E\left(\begin{tikzpicture}[guillpart,xscale=2.,yscale=1.5]
			\draw[dotted] (0,0)rectangle(1,2);
			\draw[guillsep] (0,1)--(1,1)  (1,0)--(1,2);
			\node at (0.5,0.5) { $ U_{SW} $ };
			\node at (0.5,1.5) { $ U_{NW} $ };
		\end{tikzpicture}\right) \in \Hom(\ca{B}_{\infty_W,\infty_{SN}},\setR)
		&
		u_E &=
		\begin{tikzpicture}[guillpart,xscale=2.,yscale=1.5]
			\draw[dotted] (0,0)rectangle (1,2);
			\draw[guillsep] (0,1)--(1,1)  (0,0)--(0,2);
			\node at (0.5,0.5) { $ U_{SE} $ };
			\node at (0.5,1.5) { $ U_{NE} $ };
		\end{tikzpicture}
		\in \ca{B}_{\infty_E,\infty_{SN}}
		\\
		A_{SN}(x) &=
		\pi_E\left(\begin{tikzpicture}[guillpart,xscale=2.5,yscale=1.5]
			\draw[dotted] (0,0) rectangle(1,2);
			\draw[guillsep] (0,1)--(1,1)  (1,0)--(1,2) (0,0)--(0,2);
			\node at (0.5,0.5) { $ A_S(x) $ };
			\node at (0.5,1.5) { $ A_N(x) $ };
		\end{tikzpicture}\right)
		\in \Hom(\ca{B}_{\infty_E,\infty_{SN}},\ca{B}_{\infty_E,\infty_{SN}})
	\end{align*}
	where $\pi_E$ is the representation of $(\ca{B}_{p,\infty_{SN}})_{q\in\setL^*}$ on the corner space described in Theorem~\ref{theo:1D:removingcoloursonboundaries}.
\end{theo}

In practice, the matrices $A_{SN}(x)$ will often be a simple tensor product $A_S(x)\otimes A_N(x)$, or a linear transformation of it, and will be studied with standard tools of linear algebra, as for matrix product states in the physics literature.

This theorem has an important practical consequence: the same elementary and local objects --- the full boundary eigen-elements --- build as well interesting boundary conditions for Gibbs measures and also encode nicely correlation functions of observables, which are the quantities estimated in numerical simulations and of interest in applications.

\begin{proof}[Proof of theorem~\ref{theo:onedimmarginaloutofinvROPErep}]
	Taking the expectation over the boundary variables for any rectangle of size $(L,q_1+q_2)$ containing the segment $\mathbf{S}$ in its horizontal median leads to:
	\begin{equation}
		\Espi{g_R}{\prod_{k=1}^L \indic{x_{e_k}}(X_{e_k})}
		= \frac{1}{\kappa\Lambda_{L(q_1+q_2)}\sigma^{S}_L\sigma^N_L\sigma^W_{q_1+q_2}\sigma^E_{q_1+q_2}}
		\begin{tikzpicture}[guillpart,xscale=2,yscale=1.3]
			\fill[guillfill] (0,-3) rectangle (6,3);
			\draw[guillsep] 
			(0,-2) -- (6,-2)
			(0,0) -- (6,0)
			(0,2) -- (6,2)
			(1,-3)--(1,3)
			(5,-3)--(5,3);
			\node at (0.5,-2.5) {$\ha{U}_{SW}$};
			\node at (5.5,-2.5) {$\ha{U}_{SE}$};
			\node at (0.5,2.5) {$\ha{U}_{NW}$};
			\node at (5.5,2.5) {$\ha{U}_{NE}$};
			\node at (0.5,-1) {$\ha{A}^W_{q_1}$};
			\node at (0.5,1) {$\ha{A}^W_{q_2}$};
			\node at (5.5,-1) {$\ha{A}^E_{q_1}$};
			\node at (5.5,1) {$\ha{A}^E_{q_2}$};
			\node at (3,-2.5) {$\ha{A}^S_{L}$};
			\node at (3,2.5) {$\ha{A}^N_{L}$};
			\node at (3,-1.) {$\MarkovWeight{W}_{L,q_1}$};
			\node at (3,1.) {$\MarkovWeight{W}_{L,q_2}$};
			\node at (1.5,0) {$\indic{x_{e_1}}$};
			\node at (2.5,0) {$\indic{x_{e_2}}$};
			\node at (3.5,0) [below] {$\dots$};
			\node at (4.5,0) {$\indic{x_{e_L}}$};
			\draw (2,-0.1) -- (2,0.1)
			(3,-0.1) -- (3,0.1)
			(4,-0.1) -- (4,0.1);
		\end{tikzpicture}
	\end{equation}
	As in the previous proof of Corollary~\ref{coro:fullplane:eigen} and Theorem~\ref{theo:eigenROPErep:invmeas}, we may introduce first horizontal half-plane morphisms with suitable generic element $T$ chosen so that it incorporates both the 2D-semi-group and the diagonal matrices $D_{\indic{x_{e}}}$.
	
	Then, evaluating the South-North pairing with the diagonal matrices $D_{\indic{x_{e_k}}}$ in the $\ca{T}_{\PatternShapes(\patterntype{fp}^*)}$ factor of $\ca{E}_{\PatternShapes(\patterntype{fp}^*)}$ simplifies \eqref{eq:onedimmarginal:Epq} into \eqref{eq:onedimmarginal:ROPErep}. The final equation then corresponds to vertical guillotine cuts in\eqref{eq:onedimmarginal:ROPErep} to recover a one-dimensional $\Guill_1$-algebra with a matrix product state structure.
	
	The representation of the elements in $\Hom(\ca{B}_{\infty_E,\infty_{SN}},\ca{B}_{\infty_E,\infty_{SN}})$ follows directly from Theorem~\ref{theo:1D:removingcoloursonboundaries}.
\end{proof}

In order to understand these ROPErep of correlation functions, an interesting exercise consists in checking the consistency of formulae~\eqref{eq:onedim:marginal} in an inclusion-increasing sequence of segments. For any $p\in\setN^*$, for any $(x_1,\ldots,x_p)\in S_1^p$, the expectation value
\[
\Esp{\prod_{k=1}^p\indic{x_k}(X_{e_k})}
\]
is obtained from the same quantity with $p+p'$ variables by summing over the variables $x_{p+k}\in S_1$ with $1\leq k\leq p'$. This nesting corresponds on the algebraic side to the identity
\[
\langle u_W,  
A_{SN}(x_1)\ldots A_{SN}(x_p)
C_{SN}^{p'} u_E \rangle = \sigma^S_{p'}\sigma^{N}_{p'}\langle u_W,  
A_{SN}(x_1)\ldots A_{SN}(x_p)
u_E \rangle
\]
where \begin{equation}\label{eq:defCSN}
	C_{SN} = \sum_{u\in S_1} A_{SN}(u) = \pi_E\left(  
	\begin{tikzpicture}[guillpart,yscale=1.5,xscale=1.5]
		\fill[guillfill] (0,0) rectangle (1,2);
		\draw[guillsep] (0,0)--(0,2) (1,0)--(1,2) (0,1)--(1,1);
		\node at (0.5,0.5) {$\ha{A}^S_1$};
		\node at (0.5,1.5) {$\ha{A}^N_1$};
	\end{tikzpicture}
	\right)
\end{equation}

On the algebraic side, this formula ---~which states that $u_E$ behaves as an eigenvector of $C_{SN}$~--- is a direct consequence of the compatibility of consecutive corners and the fact that morphisms in dimension one (and thus for half-planes and strips) are reduced to scalar multiplication by a geometric sequence (see Section~\ref{sec:onedimdiagowithmorph}).

\begin{coro}\label{coro:correlationfunctionfromROPEreps}
	With the same notations as previously, for any segment $\mathbf{S}$ of length $L$, the horizontal one and two point correlations are given by:
	\begin{subequations}
		\begin{align}
			\Esp{ \indic{X_{e_1}=u} } &= 
			\frac{1}{\kappa\sigma^S_1\sigma^N_1}\langle u_W, A_{SN}(u) u_E \rangle
			\\
			\Esp{ \indic{X_{e_1}=u}\indic{X_{e_L}=v}} &= 
			\frac{1}{\kappa\sigma^S_L\sigma^N_L}\langle u_W, A_{SN}(u) C_{SN}^{L-2} A_{SN}(v) u_E \rangle
		\end{align}
	\end{subequations}
	with $C_{SN} = \sum_{u \in S_1} A_{SN}(u)$.
\end{coro}

Such formulae are well-known in statistical mechanics (see \cite{DEHP} or \cite{BlytheEvans} for illustrations) and are useful to derive the large distance asymptotics of correlations functions: since the distance $L$ appears only as a power of $C_{SN}$, it means that all the informations about the correlation decays are contained in the spectral properties of $C_{SN}$. 

\begin{rema}[correlation length from ROPErep of Gibbs boundary weights]
	Corollary~\ref{coro:correlationfunctionfromROPEreps} indicates that asymptotic correlations and thus correlation lengths can be read on the the spectrum of $C_{SN}$ given in~\eqref{eq:defCSN}. In particular, $\sigma^{S}_1\sigma^{N}_1$ is the dominant eigenvalue of $C_{SN}$ and the correlation length is the inverse of the logarithm of the second eigenvalue in presence of a spectral gap. If the spectrum is absolutely continuous near $\sigma^{S}_1\sigma^{N}_1$, the behaviour of the spectral density near 
	$\sigma^{S}_1\sigma^{N}_1$ provides the critical exponent of the decay of correlation functions.
\end{rema}

The study of eigen-algebra up to morphisms and eigen-corners up to morphisms shall contain a lot of physical information that are encapsulated in definitions~\ref{def:eigenalgebrauptomorphims} and \ref{def:cornereigensemigroups}
deserves a deeper study from the point of view of phase transitions and correlation functions.

The theory of phase transitions shows that a large palette of behaviours can be observed, from exponential decorrelation on a short length scale to first-order, second-order or higher-order phase transitions (or even Kosterlitz-Thouless transitions): if our approach is general enough to describe a large variety of infinite-volume Gibbs measures, then it implies that the ROPEreps of boundary weights may have a large variety of incarnations, related to the various physical properties of the phase transitions.

\section{The Gaussian model}\label{sec:appli:gaussian}	

This section presents an overview of \cite{BodiotSimon} where all detailed computations are presented and related to the existing literature. Gaussian variables are put on the edges of $\setZ^2$ and and we consider $S_1=\ca{H}_1$ and $S_2=\ca{H}_2$ two finite-dimensional Hilbert spaces. These two state spaces are not finite but are measurable spaces with a $\sigma$-finite Lebesgue reference measure and thus enters the framework of the present paper through the generalization presented below in Section~\ref{sec:higherdim}. We focus in this section on the operadic content and the interested reader can already understand the concepts and computations.

\subsection{From Gaussian isometries to face weight densities}

\subsubsection{Gaussian isometry} 
From the point of view of Gaussian processes, this corresponds to the specification of two vectors spaces $\ca{H}_1$ and $\ca{H}_2$ on horizontal and vertical edges and a Hermitian product $\Sigma_D$ on 
\[
\left(\bigoplus_{e\in \Edges{D}_h}\ca{H}_1 \right)\oplus\left(  \bigoplus_{e\in\Edges{D}_v} \ca{H}_2\right)
\]
for any finite domain $D\subset \setZ^2$. We note $\ca{H}_{D,\Sigma_D}$ the Hilbert completion of this space under $\Sigma$. A Gaussian field is thus an isometry $X$ from the Hilbert space $\ca{H}_{D,\Sigma_D}$ to a probability space $L^2(\Omega,\Prob,\nu)$ so that, for any $h\in \ca{H}_{D,\Sigma_D}$, $X(h)$ is a centred normal r.v. $\nu$ and, for any $h,h'\in H_D$, $\Esp{X(h)X(h')}=\Sigma_D(h,h')$. In particular, for any edge $e\in\Edges{D}_h$ (resp. $\Edges{D}_v$  )and $h\in \ca{H}_1$ (resp. $\ca{H}_2$), there is a canonical embedding $\iota_e(h)\in \ca{H}_{D,\Sigma_D}$ by inserting $h$ in the space $\ca{H}_i$ associated to $E$ and $0$ elsewhere. To any edge $e\in \Edges{D}$, we then have a collection of r.v. $X(\iota_e(h))$ indexed by $\ca{H}_1$ or $\ca{H}_2$. The Markov property can then rephrased as an identity between conditional covariance matrices, which can themselves be related to the blocks of the Hermitian inner product $\Sigma_D$. We may follow this machinery combined with the proof of Proposition~\ref{prop:MarkovOneDimSecond} to obtain elementary conditional covariance matrices and, then, face weights $\MarkovWeight{W}$ as exponential of suitable quadratic forms. This is however tedious and presents some traps. Moreover, constructing a product $\Sigma_D$ so that the face weights may be homogeneous is hard without further a priori on the law of the process. To circumvent this difficulty and still be able to apply the formalism of the present paper, we may also consider the definition of the process from its law or at least its density w.r.t. to reference measures, mimicking the definition of a finite-dimensional full-rank Gaussian r.v. from a density w.r.t. the Lebesgue measure.

\subsubsection{Face weight densities}
We introduce a bilinear form $Q$ on $\ca{H}_1\oplus \ca{H}_1\oplus \ca{H}_2\oplus \ca{H}_2$ with a block decomposition 
\[
Q = \begin{pmatrix}
	Q_{SS} & Q_{SN} & Q_{SW} & Q_{SE} \\
	Q_{NS} & Q_{NN} & Q_{NW} & Q_{NE} \\
	Q_{WS} & Q_{WN} & Q_{WW} & Q_{WE} \\
	Q_{ES} & Q_{EN} & Q_{EW} & Q_{EE}
\end{pmatrix}
\]
and the face weight density $\MarkovWeight{W}$ w.r.t.~the Lebesgue measure $H_1\oplus H_2\oplus H_2$ defined by:
\[\MarkovWeight{W}(x_S,x_N,x_W,x_E) = \exp\left(-\frac{1}{2}
\begin{pmatrix}
	x_S \\
	x_N \\
	x_W \\
	x_E 
\end{pmatrix}^*
Q
\begin{pmatrix}
	x_S \\
	x_N \\
	x_W \\
	x_E
\end{pmatrix}
\right)
\]
It is interesting to note that, in this case, integrating~w.r.t. a horizontal (resp. vertical) r.v. $x$ in the computation of partition functions or expectation values requires an integration w.r.t. the measure $\exp(-x^*(Q_{SS}+Q_{NN})x/2)dx$ (resp. with $Q_{WW}+Q_{EE}$). These integration steps appears as building blocks and are well-defined if $Q_{SS}+Q_{NN}$ and $Q_{WW}+Q_{EE}$ are positive definite operators. One may rewrite the model using a bilinear form \[
Q_c = Q-\Diag(Q_{SS},Q_{NN},Q_{WW},Q_{EE})
\]
and a weight density $\MarkovWeight{W}_c$ as above with $Q$ replaced with $Q_c$ but now these densities may be considered w.r.t. the normal laws $\ca{N}(0,(Q_{SS}+Q_{NN})^{-1})$ for horizontal edges and $\ca{N}(0,(Q_{WW}+Q_{EE})^{-1})$ for vertical edges, which can be chosen as reference measures $(\mu_e)$.

\begin{rema}If a restriction such as $x_S+x_E=x_N+x_0$ is imposed as it is the case for height processes such as a discrete version of the GFF (massive or not), then, even if the process is Gaussian, it is not possible any more to define reference measures on the edges. A solution beyond the scope of the present paper consists in putting reference measures not on individual edges but on shells of domains with suitable compatibility conditions. The operadic structure of the next section should however be modified and is not as easy to use.
\end{rema}

\subsubsection{Translation-invariant infinite-dimensional Gibbs measure}

Whenever $Q$ is positive semi-definite, it is easy to describe the unique translation-invariant Gibbs measures of the model (see \cite{GibbsGeorgii}) by Fourier transform through independent modes in $L^2(S^1\times S^1; \ca{H}_1\oplus\ca{H}_2)$: this construction is entirely non-local and the Markov property is hidden in the complex structure of some function $J$ on $\setC^*\times\setC^*$. The purpose of \cite{BodiotSimon} that we summarized here is to illustrate that all these results can be recovered from the present paper through explicit computations on ROPEreps.

\subsection{From the guillotine operad of weights to the guillotine Schur products}

Face weights as well as boundary weights and any marginal law are centred Gaussian weights on suitable Hilbert subspaces and are thus described by their covariance matrix or its inverse, the coupling matrix (such as $Q$ above). For any rectangle of size $(p,q)$, we associated the Hilbert space $\ca{H}_1^{2p}\oplus\ca{H}_2^{2q}$ and the set $\ca{Q}_{p,q}$ of bounded positive definite Hermitian operators on these space. When gluing rectangles, one integrates a product of Gaussian weights over the values on the guillotine cut and it produces again a Gaussian weight with a coupling matrix obtained by Schur complements of the two initial weights (see \cite{BodiotSimon}). 

\begin{lemm}[from \cite{BodiotSimon}]
	The products described in Theorem~\ref{theo:canonicalexampleGuill:continuous} restricted to the Gaussian weights parametrized by the sets $\ca{Q}_{\bullet}$ lift to a (non-linear) guillotine algebra structure on $\ca{Q}_{\bullet}$ given by Schur complements.
\end{lemm}

Since invariant boundary weights are expected to be Gaussian, it is natural to define a ROPE made of densities of Gaussian weights on suitable Hilbert spaces associated to the infinite half-lines that appear in the new pattern shapes (half-strips, corners, strips and half-planes). We now define the following four Hilbert spaces
\begin{align*}
\ca{H}_S &= L^2(\setZ_{-}^*; \ca{H}_2)
&
\ca{H}_N &= L^2(\setZ_{+}; \ca{H}_2)
\\
\ca{H}_W &= L^2(\setZ_{-}^*; \ca{H}_1)
&
\ca{H}_E &= L^2(\setZ_{+}; \ca{H}_2)
\end{align*}
which generalize the correspondence $\ca{H}_i^{k} \simeq L^2(\{1,\ldots,k\};\ca{H}_i)$ to infinite half-lines. For each non-finite pattern shape, we introduce the set $\ca{Q}_{\bullet}$ of Hermitian operators with suitable analytical properties (see \cite{BodiotSimon}) on the direct sum of these spaces; as an example $\ca{Q}_{\infty_E,q}$ is the space of Hermitian operators bounded and bounded away from zero on $\ca{H}_E\oplus \ca{H}_E \oplus \ca{H}_2^q$. Gluing a segment after a half-line produces the same half-line with a shift: this is encoded in shift isomorphisms on the four previous spaces:
\[
\begin{split}
\ca{H}_i \times L^2(\setZ_{+}; \ca{H}_i)  & \to L^2(\setZ_{+}; \ca{H}_i)
\\
(u,f) & \mapsto u\indic{1}(\cdot) + f(\cdot -1)\indic{>2}(\cdot)
\end{split}
\]
and hence a product $\ca{Q}_{1,q} \times\ca{Q}_{\infty_E,q}  \to \ca{Q}_{\infty_E,q}$ given by
\[
(Q,Q') \mapsto \begin{tikzpicture}[guillpart,yscale=1.2,xscale=1.5]
	\fill[guillfill] (0,0) rectangle (2,1);
	\draw[guillsep] (2,0)--(0,0)--(0,1)--(2,1) (1,0)--(1,1);
	\node at (0.5,0.5) {$Q$};
	\node at (1.5,0.5) {$Q'$};
\end{tikzpicture}_\ca{Q}
\]
with the same Gaussian integration procedure on the cut variable in $\ca{H}_2$, which produces a quadratic form on $(\ca{H}_1\oplus\ca{H}_E)^2\oplus \ca{H}_2$ and is identified then with $\ca{Q}_{\infty_E,1}$ with the previous shift morphism. Such products can be defined then on all shapes with the same procedure until reaching the following lemma.

\begin{lemm}[from \cite{BodiotSimon}]
	The spaces $\ca{Q}_{\PatternShapes(\patterntype{fp}^*)}$ with products on the boundaries involving the shift morphisms is a full plane $\Guill_2$-algebra.
\end{lemm}

As announced previously, the ROPE definition lifted at the level of parameters provides spaces $\ca{Q}_\bullet$ with predefined shift morphisms. In this case, the previous equations \eqref{eq:reducedmorphism:eigeneq} and \eqref{eq:concreteeigencorner} on $A_S(\bullet)$, $A_N(\bullet)$, $A_W(\bullet)$, $A_E(\bullet)$ and $U_{ab}$, are lifted to a set of $8$ fixed point equations on $8$ objects (one Hermitian matrix $\ov{Q}_a$ on each side and one Hermitian matrix $\ov{Q}_{ab}$ on each corner) as defined in Section~\ref{sec:fixedpoints} with \emph{already known} spaces and morphisms. We are then left with a traditional finite set of equations on a finite number of operator with the following result.

\begin{theo}[from \cite{BodiotSimon}]\label{theofromBS:Gaussian}
	The fixed point equations on $\ov{Q}_{a}$ and $\ov{Q}_{ab}$ are equivalent to \emph{well-defined recursions} on the matrix coefficients on these operators, with a unique solution (which can also be computed numerically and analytically). A full correspondence can be established between this solution for the operators $\ov{Q}_{a}$ and $\ov{Q}_{ab}$ and the infinite-volume Gibbs measure inherited from the Fourier transform as in \cite{GibbsGeorgii}.
\end{theo}

After the trivial models, this theorem is a second proof of concept that the present constructions can be useful in practice.

\section{The six-vertex model}\label{sec:appli:sixvertex}

\subsection{An overview of the model}

\paragraph*{"Exact" solutions}
The six-vertex model \cite{baxterbook,ReshetikhinSixV} is one of the paradigmatic model of two-dimensional statistical mechanics. On one hand side, it exhibits phase transitions and its scaling limit in the disordered phase is related to the Gaussian free field; on the other hand side, it is exactly solvable in the sense of integrable systems. However, this exact solvability does not mean that all the quantities are known nor that all the computations are easy or rigorous.

We present in this section the work presented in \cite{SimonSixV} by emphasizing on the contrast between the standard Bethe Ansatz approach and the operadic approach, both of them having their own interest.

\paragraph*{Definitions}

The six-vertex model corresponds to a spatial two-dimensional Markov process as defined in Section~\ref{sec:proba} for the state spaces $S_1=S_2=\{-1,1\}$. We define the face weight $\MarkovWeight{R}:\{-1,1\}^4\to\setR_+$ by
\begin{equation}
	\MarkovWeight{R}(x_S,x_N,x_W,x_E) =
	\begin{cases}
		\asixv & \text{if $(x_S,x_N,x_W,x_E)\in\{(1,1,1,1),(-1,-1,-1,-1)\}$} \\
		\bsixv & \text{if $(x_S,x_N,x_W,x_E)\in\{(-1,-1,1,1),(1,1,-1,-1)\}$} \\
		\csixv & \text{if $(x_S,x_N,x_W,x_E)\in\{(1,-1,-1,1),(-1,1,1,-1)\}$} \\
		0 & \text{else}
	\end{cases}
\end{equation}
for some positive numbers $\asixv$, $\bsixv$ and $\csixv$ as illustrated in figure~\ref{fig:sixvertexmodel}. In particular, all the non-zero weights correspond to configurations satisfying $x_W+x_S=x_N+x_E$. In the spin interpretation in which a $-1$ value on the South or on the North (resp. on the West or on the East) corresponds to a top-to-bottom vertical (resp. right-to-left horizontal) arrow and a $+1$ value to a bottom-to-top vertical (resp. left-to-right horizontal) arrow, a configuration on a subset of the lattice $\setZ^2$ with a non-zero weight have necessarily as many incoming and outgoing arrows around any elementary square face.

\begin{figure}
	\begin{center}
		\begin{tikzpicture}
			\begin{scope}
				\draw[thick] (0,0) rectangle (1,1);
				\draw[ultra thick,->] (0.5,0.7) -- (0.5,1.3);
				\draw[ultra thick,->] (0.5,-0.3) -- (0.5,0.3);
				\draw[ultra thick,->] (-0.3,0.5) -- (0.3,0.5);
				\draw[ultra thick,->] (0.7,0.5) -- (1.3,0.5);
				\node at (0.5,-0.3) [anchor = west,red] {$1$};
				\node at (0.5,1.3) [anchor = west,red] {$1$};
				\node at (-0.3,0.5) [anchor = south,red] {$1$};
				\node at (1.3,0.5) [anchor = south,red] {$1$};
				\node at (0.5,-0.6) {$\asixv$};
			\end{scope}
			\begin{scope}[xshift=2.2cm]
				\draw[thick] (0,0) rectangle (1,1);
				\draw[ultra thick,<-] (0.5,0.7) -- (0.5,1.3);
				\draw[ultra thick,<-] (0.5,-0.3) -- (0.5,0.3);
				\draw[ultra thick,<-] (-0.3,0.5) -- (0.3,0.5);
				\draw[ultra thick,<-] (0.7,0.5) -- (1.3,0.5);
				\node at (0.5,-0.3) [anchor = west,red] {$-1$};
				\node at (0.5,1.3) [anchor = west,red] {$-1$};
				\node at (-0.3,0.5) [anchor = north,red] {$-1$};
				\node at (1.3,0.5) [anchor = north,red] {$-1$};
				\node at (0.5,-0.6) {$\asixv$};
			\end{scope}
			\begin{scope}[xshift=4.4cm]
				\draw[thick] (0,0) rectangle (1,1);
				\draw[ultra thick,<-] (0.5,0.7) -- (0.5,1.3);
				\draw[ultra thick,<-] (0.5,-0.3) -- (0.5,0.3);
				\draw[ultra thick,->] (-0.3,0.5) -- (0.3,0.5);
				\draw[ultra thick,->] (0.7,0.5) -- (1.3,0.5);
				\node at (0.5,-0.3) [anchor = west,red] {$-1$};
				\node at (0.5,1.3) [anchor = west,red] {$-1$};
				\node at (-0.3,0.5) [anchor = south,red] {$1$};
				\node at (1.3,0.5) [anchor = south,red] {$1$};
				\node at (0.5,-0.6) {$\bsixv$};
			\end{scope}
			\begin{scope}[xshift=6.6cm]
				\draw[thick] (0,0) rectangle (1,1);
				\draw[ultra thick,->] (0.5,0.7) -- (0.5,1.3);
				\draw[ultra thick,->] (0.5,-0.3) -- (0.5,0.3);
				\draw[ultra thick,<-] (-0.3,0.5) -- (0.3,0.5);
				\draw[ultra thick,<-] (0.7,0.5) -- (1.3,0.5);
				\node at (0.5,-0.3) [anchor = west,red] {$1$};
				\node at (0.5,1.3) [anchor = west,red] {$1$};
				\node at (-0.3,0.5) [anchor = north,red] {$-1$};
				\node at (1.3,0.5) [anchor = north,red] {$-1$};
				\node at (0.5,-0.6) {$\bsixv$};
			\end{scope}
			\begin{scope}[xshift=8.8cm]
				\draw[thick] (0,0) rectangle (1,1);
				\draw[ultra thick,->] (0.5,0.7) -- (0.5,1.3);
				\draw[ultra thick,<-] (0.5,-0.3) -- (0.5,0.3);
				\draw[ultra thick,->] (-0.3,0.5) -- (0.3,0.5);
				\draw[ultra thick,<-] (0.7,0.5) -- (1.3,0.5);
				\node at (0.5,-0.3) [anchor = west,red] {$-1$};
				\node at (0.5,1.3) [anchor = west,red] {$1$};
				\node at (-0.3,0.5) [anchor = south,red] {$1$};
				\node at (1.3,0.5) [anchor = south,red] {$-1$};
				\node at (0.5,-0.6) {$\csixv$};
			\end{scope}
			\begin{scope}[xshift=11.cm]
				\draw[thick] (0,0) rectangle (1,1);
				\draw[ultra thick,<-] (0.5,0.7) -- (0.5,1.3);
				\draw[ultra thick,->] (0.5,-0.3) -- (0.5,0.3);
				\draw[ultra thick,<-] (-0.3,0.5) -- (0.3,0.5);
				\draw[ultra thick,->] (0.7,0.5) -- (1.3,0.5);
				\node at (0.5,-0.3) [anchor = west,red] {$1$};
				\node at (0.5,1.3) [anchor = west,red] {$-1$};
				\node at (-0.3,0.5) [anchor = north,red] {$-1$};
				\node at (1.3,0.5) [anchor = north,red] {$1$};
				\node at (0.5,-0.6) {$\csixv$};
			\end{scope}
		\end{tikzpicture}
	\end{center}
	\caption{The six-vertex model: configurations and weights. At the free fermion point, the values are given by $\asixv=\cos(u)$, $\bsixv=\sin(u)$ and $\csixv=1$. The ice model corresponds to $\asixv=\bsixv=\csixv=1$.}\label{fig:sixvertexmodel}
\end{figure}

The phase diagram is controlled by the reduced parameter
\[
\Delta = \frac{ \asixv^2+\bsixv^2-\csixv^2 }{ 2\asixv\bsixv }
\]
and the triplet $(\asixv,\bsixv,\csixv)$ is defined only up to a global multiplicative constant. The complete expected phase diagram is presented in figure~\ref{fig:sixv:phasediagram}: many parts of it are rigorous (see \cite{DuminilCondensation} for a review of the bibliography) but many computations are not since they still rely on precise asymptotics on the large size limit of Bethe Ansatz formulae for eigen-vectors of the periodic transfer matrix.
\begin{figure}
	\begin{center}
		\begin{tikzpicture}[scale=2]
			\draw[->] (0,0)--(2,0) ;
			\draw[->] (0,0)--(0,2) ;
			\node at (2,0) [right] {$a/c$};
			\node at (0,2) [left] {$b/c$};
			\draw[red, thick] (1,0)--(2,1);
			\draw[red, thick] (0,1)--(1,2);
			\draw[red, thick] (0,1)--(1,0);
			\node at (1.75,0.25) {$F$};
			\node at (0.25,1.75) {$F'$};
			\node at (1.5,1.5) {$D$};
			\node at (0.25,0.25) {AF};
			\node at (1,1) {$\bullet$};
			\node at (1,0) [below] {$1$};
			\node at (0,1) [left] {$1$};
			\draw[ultra thick, dashed, blue] (1,0) arc (0:90:1);
		\end{tikzpicture}
	\end{center}
	\caption{\label{fig:sixv:phasediagram}Expected phase diagram of the six-vertex model without external field. The domains $F$ and $F'$ correspond to the strongly ordered ferromagnetic phase $\Delta>1$. The phase AF corresponds to the strongly ordered antiferromagnetic phase with $\Delta<-1$. The phase $D$ corresponds to the disordered phase $-1<\Delta<1$ with long-range correlation. The point $\bullet$ with $a/c=b/c=1$ corresponds to the ice model studied by Lieb \cite{LiebIce}. The dashed line corresponds to the free fermion line $\Delta=0$.}
\end{figure}

\subsection{The traditional Bethe Ansatz approach: advantages and drawbacks}

\subsubsection{The Ansatz for finite size}
The traditional Bethe Ansatz consists in the exact diagonalization, \emph{for all horizontal size $P\in\setN^*$}, of the transfer matrix
\begin{equation}\label{eq:sixv:transfermatsizeP}
T_P = \begin{tikzpicture}[guillpart,yscale=1.5,xscale=1.5]
	\fill[guillfill] (0,0) rectangle (4,1);
	\draw[guillsep] (0,0)--(4,0)--(4,1)--(0,1)--(0,0) (1,0)--(1,1) (2,0)--(2,1) (3,0)--(3,1);
	\node at (0.5,0.5) {$\MarkovWeight{R}$};
	\node at (1.5,0.5) {$\MarkovWeight{R}$};
	\node at (2.5,0.5) {$\ldots$};
	\node at (3.5,0.5) {$\MarkovWeight{R}$};
	\node at (0,0.5) {$\bullet$};
	\node at (4,0.5) {$\bullet$};
\end{tikzpicture}
\end{equation}
acting on $(\setC^2)^{\otimes P}$ with periodic boundary conditions (see Section~\ref{sec:guill2:cylindershape}) by proposing a guess for the structure of the eigenvectors. By noting $1\leq x_1<x_2<\ldots<x_K<P$ the position of the $K$ values $+1$ on the bottom line (the number of up arrows on each line is conserved), Bethe's historical Ansatz \cite{Bethe} corresponds a vector $\psi\in (\setC^2)^{\otimes P}$ with 
\[
\psi = \sum_{1\leq x_1<x_2<\ldots<x_K<P} \psi(x_1,\ldots,x_K) e_{-1}^{\otimes x_1-1}\otimes e_{1}\otimes e_{-1}^{\otimes x_2-x_1-1}\otimesdots  e_{1} \otimes  e_{-1}^{\otimes P-x_K-1}
\]
in the canonical basis $(e_{-1},e_1)=( (0,1), (1,0) )$ with the following guessed structure
\begin{equation}
\label{eq:sixv:CBA}
\psi(x_1,\ldots,x_K) =  \sum_{\sigma\in \Sym{K}} A_\sigma \prod_{k=1}^K z_{\sigma(k)}^{x_k}
\end{equation}
for suitable coefficients $(A_\sigma)_{\sigma\in\Sym{K}}$ and complex numbers $(z_k)_{1\leq k\leq K}$. The condition $T_P\psi = \Lambda \psi$ determines the coefficients $A_\sigma$ up to an overall constant in terms of the coefficients $z_k$. These coefficients, called Bethe roots, satisfy the so-called non-linear algebraic Bethe equations
\begin{equation}
	\label{eq:sixv:betheeqs}
z_k^P = (-1)^{K-1} \prod_{k=1}^K \frac{1+z_kz_l-2\Delta z_k}{1+z_kz_l-2\Delta z_l} = \prod_{k=1}^K \Ssixv(z_k,z_l)
\end{equation}
and the eigenvalue is given by 
\begin{equation}\label{eq:sixv:Betheeigenval}
\Lambda_{P,K}((z_k)) = \asixv^P\prod_{k=1}^K \Lsixv(z_k) + \bsixv^P \prod_{k=1}^K \Msixv(z_k)
\end{equation}
for suitable homographic functions $\Lsixv$ and $\Msixv$  defined by
\begin{align*}
\Lsixv(z) &= \frac{\asixv\bsixv+(\csixv^2-\bsixv^2) z}{\asixv^2-\asixv\bsixv z}
&
\Msixv(z) &= \frac{(\asixv^2-\csixv^2)-\asixv\bsixv  z}{\asixv\bsixv-\bsixv^2 z}
\\
\Ssixv(z,w) &= \frac{\Msixv(w)\Lsixv(z)-1}{1-\Lsixv(w)\Msixv(z)} 
\end{align*}

The Bethe equations cannot be solved explicitly in the generic case and some of the hard tasks are
\begin{enumerate}
	\item find and describe the good solution $(z_k^*)_{1\leq k\leq K}$ that corresponds to the Perron-Frobenius eigenvector of $T_P$ with $K$ particles;
	\item show that all the roots are distinct and do not coincide with singularities of the factors of \eqref{eq:sixv:betheeqs}, of $\Lsixv$, $\Msixv$ and $\Ssixv$;
	\item show that the corresponding eigenvector $\psi$ is non-zero;
	\item extract the large $P$ behaviour of the Bethe roots associated to the Perron-Frobenius eigenvector (the so-called condensation, see below);
	\item extract, from the Bethe Ansatz expression of the eigenvector, formulae for correlation functions (in the large $P$ limit).
\end{enumerate}

Most of them are feasible and already performed by physicists in a non-rigorous way, while the rigorous proofs are often very technical steps (when they exist, see \cite{DuminilCondensation} for a review). The heuristic behaviour is that solutions of \eqref{eq:sixv:betheeqs} are enumerated by subsets of size $K$ of $P$-th roots of unity (or of $-1$) and, in the case of Perron-Frobenius case, by the ones with largest real parts with a symmetry under complex conjugation. In the limit $P\to\infty$ with $K/P\to \rho\in (0,1)$ (fixed density), it is expected that the Bethe roots get condensed on a one-dimensional curve parametrized by an interval $[-\pi\rho,\pi\rho]\subset  [-\pi,\pi]$ with an explicit equation derived from \eqref{eq:sixv:betheeqs}. More precisely, we expect the Bethe roots to take the asymptotic shape $z_k\simeq \zeta(2\pi I_k/P)$ with consecutive integers or half-integers $I_k$ ranging from $-\rho P /2$ to $\rho P/2$. These integers are related to the selection of correct determinations of the logarithms of the l.h.s.~and r.h.s.~of \eqref{eq:sixv:betheeqs} used to select the Perron-Frobenius solution. In the thermodynamic limit, these integers suitably normalized tend to a uniform density on $[-\pi\rho,\pi\rho]$. The function $\zeta : [-\pi \rho,\pi \rho] \to \setC$ is expected to satisfy the continuous equation:
\begin{equation}
	\zeta(\theta) e^{-i\theta} = \exp\left( \frac{1}{2\pi}\int_{-\pi\rho}^{\pi\rho} \log \Ssixv( \zeta(\theta), \zeta(\phi) ) d\phi \right) 
	\label{eq:condensedrootcurve}
\end{equation}
using the principal determination of the complex logarithm. This is equivalent to the traditional presentation in terms of the density $\rho(\zeta)$ of the roots along their one-dimensional curve.

Although all these heuristics are not always rigorous (or rigorous after long and technical proofs), this last equation is the only final information required to obtain the free energy density, given in the regime $\asixv>\bsixv$ by
\[
\frac{\log \Lambda_{P,K}((z_k))}{P} 
\to f(\rho) = \log \asixv + \frac{1}{2\pi}\int_{-\pi\rho}^{\pi\rho} \log \Lsixv(\zeta(\theta)) d\theta
\]
and a similar expression with $\bsixv$ and $\Msixv$ if $\bsixv>\asixv$ with a suitable chosen $\rho$ in the disordered phase $|\Delta|<1$. The target generalized eigenvalue we wish to obtain by the guillotine approach is then
\begin{equation}
	\label{eq:sixv:conjectureeigenval}
	\Lambda(\rho) = e^{f(\rho)} = \asixv \exp\left( \frac{1}{2\pi}\int_{-\pi\rho}^{\pi\rho} \log \Lsixv(\zeta(\theta)) d\theta
	\right)
\end{equation}

\subsubsection{Discussion of the Bethe Ansatz in relation with the guillotine operad}

We believe that some problems of rigour with the Bethe Ansatz are essentially related with the choice of boundary conditions and  be avoided by the operadic approach. We first present some of the basic assumptions at the heart of the Bethe Ansatz.

\paragraph*{Some limitations induced by the traditional approach}

The choice of the transfer matrix $T_P$ in \eqref{eq:sixv:transfermatsizeP} breaks the symmetry between the two dimensions of the plane: one of them (here the horizontal one) is used to a build a large-dimensional state space and the second one (here the vertical one) is used to multiply transfer matrices using standard linear algebra. It thus selects preferentially the products $m_{SN}$ and square associativity ---and all the associated tools like corner-eigen-elements and ROPEreps--- remains completely hidden in the computations. 

The situation is then worsened by the choice of periodic boundary conditions. Initially, the  purpose of periodic boundary conditions is to hide the question of choosing a suitable boundary condition and obtain more simpler computations, in particular in terms of Fourier transform (use of roots of unity). However, as discussed in Section~\ref{sec:guill2:cylindershape}, the cylindrical geometry avoids horizontal half-strips and \emph{thus} corners: the eigenvectors of $T_P$ are not expected to have particular relations with corner eigen-elements as defined in the present paper.

Beyond the question of eigen-elements, the choice of periodic conditions has an important practical consequence: it totally prevents the use of Kolmogorov's extension theorem to construct the infinite-volume Gibbs measures of the model. Indeed, the cylinders of size $P$ (and thus the matrices $T_P$) are not geometrically related ---no inclusion for example--- in the same way as rectangles can be nested. For example, no recursion can be implemented and no particular algebraic relation is expected to hold between the eigenvectors of the matrices $T_P$ when $P$ varies. Hence, for each quantity of interest $Q_P$, the only remaining path is the analytical one:
\begin{itemize}
	\item one guesses the limit $Q_{\infty}$ of the quantity of interest $(Q_P)$ (use of heuristics)
	\item one controls analytically, for large $P$, the error $|Q_P-Q_\infty|$ 
\end{itemize}
This corresponds to what is mostly done in the literature and, for most quantities such as the free energy density or the two-point correlations functions, an unavoidable bottleneck is the proof of condensation of the Bethe roots $(z_k)$ to the limiting curve \eqref{eq:condensedrootcurve}, as done for example in \cite{DuminilCondensation}. An inherent limit of this approach is that it is not possible to define and compute the limit $P\to\infty$ of the eigenvectors \eqref{eq:sixv:CBA} themselves. In practice, one must first compute the observable $Q_P$ for all $P$ and afterwards take the limit.

Finally, although the matrices $(T_P)_P$ are sufficiently well understood in the six-vertex model to use the Bethe Ansatz, their detailed computation remains hard since the state space is large and the number of matrices $\MarkovWeight{R}$ multiplied horizontally becomes larger and larger.

We claim that the operadic approach developed in \cite{SimonSixV} circumvents this difficulty in the following way.

\paragraph*{The operadic approach: the philosophy}

As already emphasized, the operadic approach allows one to construct the infinite-volume Gibbs measure using Kolmogorov's extension theorem. We now enter in the details.

The constraints put by this extension theorem are imposed on boundary weights and still require the computation of "annular" transfer matrices between nested rectangles. At this point, the transfer matrix approach with \eqref{eq:sixv:transfermatsizeP} still appears much simpler. When moving to ROPEreps and generalized half-strip and corner eigen-elements, the situation changes completely and we now explore the changes.

First, no power (horizontal or vertical) of $\MarkovWeight{R}$-matrices is required and the situation becomes closer to the 1D situation in which the invariant measure does not involve any power of the transition matrix. Only the action on $\MarkovWeight{R}$ on the boundary elements is used.

In a second time, the invariant boundary elements in the ROPErep of boundary weights are not specific to a given size and are already associated to the $P\to\infty$ limit. No analytical control of an error as $P\to\infty$ is required. For the six-vertex model, this corresponds to the \emph{absence} of finite size Bethe equations \eqref{eq:sixv:betheeqs} and the \emph{only} use of \eqref{eq:condensedrootcurve}. In particular, it circumvents a rigorous proof of the condensation of the roots.

The mechanism of this simplification is as follows: in opposition to the cylinders, the translation of the constraints of the extension theorem in the ROPErep language provides precise equations to be solved by the $P\to\infty$ limit. If the eigen-elements associated to the \emph{conjectured} limit such as \eqref{eq:condensedrootcurve} satisfy these constraints and are well-defined, then they become \emph{rigorous} building blocks of the boundary weights and of the infinite-volume limit. 

The situation can be compared to 1D recursions $u_{n+1}=f(u_n)$: the traditional Bethe Ansatz approach corresponds to the exact (but difficult) computation of $u_n =h(n)$, then the guess of the limit $l$ and finally the rigorous verification of the limit by the control of $|h(n)-l|$. The operadic approach corresponds to the direct computation of $l$ through the fixed point property $l=f(l)$.

Finally, the half-strip and corner structure does not break the invariance between the four directions in the plane as in the construction of $T_P$.

We do not pretend that the Bethe Ansatz is useless and we now show that it is the contrary. The form \eqref{eq:sixv:CBA} is adapted to the matrices $T_P$. We now see how variants of this form with the same physical insight can be adapted to the operadic language by precisely providing the spaces required in the (necessary) guessing of the ROPE structure.

\subsection{Operadic Bethe Ansatz}
\subsubsection{Recycling the basics of Bethe Ansatz}

\paragraph*{Fourier transform and directed walks}
The first ingredient is the Fourier transform behind the modes $z^x$ in \eqref{eq:sixv:CBA}. Each mode with $z^P=1$ corresponds to the eigenvector of a translation-invariant generator of a random walk on $\setZ/P\setZ$, with eigenvalue $\asixv^P\Lsixv(z)$ or $\bsixv^P\Msixv(z)$ (depending on the phase). On the plane, this picture persists: a up or right arrow in a sea of down and left arrows perform a directed walk from South-West to North-East. 

We show in \cite{SimonSixV} that the one-particle partition function of the oriented walk already has its own structure of $\Guill_2$-algebra, with associated generalized eigen-elements on half-strips and boundaries (without any guessing !). In particular, one recovers the normalized one-particle transfer matrix on $\setZ$, with its Fourier eigen-modes and eigenvalues $\Lsixv(e^{i\theta)}$ by gluing opposite half-strips elements. However, already at the one-particle level, we also recovers non-trivial corner eigen-elements.

\begin{figure}
	\begin{center}
		\begin{tikzpicture}[scale=0.5]
			\foreach \x in {0,1,2,3,11,12,13,14} {
				\draw (\x,0)--(\x,10);	
			}
			\foreach \x in {4,5,6,7,8,9,10} {
				\draw[ultra thick, violet] (\x,0)--(\x,10);	
			}
			\foreach \y in {0,1,2,8,9,10} {
				\draw (0,\y)--(14,\y);
			}
			\foreach \y in {3,4,5,6,7} {
				\draw[ultra thick, violet] (0,\y)--(14,\y);	
			}
			\foreach \x in {0,1,...,13} {
				\foreach \y in {0,1,...,10}
				\draw[->,red] (\x+0.5,\y+0.3)--(\x+0.5,\y-0.3);
			}
			\foreach \x in {0,1,...,14} {
				\foreach \y in {0,1,...,9}
				\draw[->,red] (\x+0.3,\y+0.5) -- (\x-0.3,\y+0.5);
			}
			
			\foreach \x/\y in { 1/0, 2/1, 6/2, 7/3, 7/4, 7/5, 8/6, 8/7, 9/8, 9/9, 12/10 } {
				\draw[->,blue,ultra thick] (\x+0.5,\y-0.3)--(\x+0.5,\y+0.3);
			}
			\foreach \x/\y in { 2/0, 3/1, 4/1, 5/1,6/1, 7/2, 8/5, 9/7, 10/9, 11/9,12/9 } {
				\draw[->,blue,ultra thick] (\x-0.3,\y+0.5)--(\x+0.3,\y+0.5);
			}
			\foreach \x/\y in { 3/1, 4/1,5/1, 7/3, 7/4, 8/6, 9/8, 10/9, 11/9} {
				\node[blue] at (\x+0.5,\y+0.5) {$\ti{\bsixv}$};	
			}
			\foreach \x/\y in { 1/0, 2/0, 2/1, 6/1, 6/2, 7/2, 7/5, 8/5, 8/7, 9/7, 9/9, 12/9} {
				\node[blue] at (\x+0.5,\y+0.5) {$\ti{\csixv}$};	
			}
			\fill[gray, opacity=0.3] 
			(0,0) rectangle (4,3)
			(10,7) rectangle (14,10)
			;
			\fill[pattern=horizontal lines,opacity=0.4] 
			(4,0) rectangle (7,3)
			(9,7) rectangle (10,10)
			;
			\fill[pattern=north east lines,opacity=0.4]
			(7,0) rectangle (8,3)
			(8,7) rectangle (9,10)
			;
		\end{tikzpicture}
	\end{center}
	\caption{\label{fig:sixv:defectline}Line of $+1$ values inserted in a reference configuration with $-1$ values everywhere. The normalized face weights along the line are defined by $\ti{\bsixv}=\bsixv/\asixv$ and $\ti{\csixv}=\csixv/\asixv$. The thick violet line indicate an example of guillotine cuts on the plane.}
\end{figure}

\paragraph*{The two-particle scattering}

The second ingredient of Bethe Ansatz is the factorization of an $n$-body dynamics into $2$-body interactions terms, as illustrated in the structure of the terms $A_\sigma$ and the $\Ssixv$ factors in \eqref{eq:sixv:betheeqs}. This $\Ssixv$ factor is already present when studying two interacting oriented walks in the plane.

We show in \cite{SimonSixV} that the two-particle partition functions together with the zero- and one-particle partition functions again have a structure of a $\Guill_2$-algebra. Again, the expression of $\Ssixv$ is recovered together when considering strips, together with new functions on half-strips and corners that are invisible to the cylindrical transfer matrix approach.

\paragraph*{The Fock space structure: a morphism}

The last ingredient of Bethe Ansatz is the Fock space structure in \eqref{eq:sixv:CBA} through the summation of permutations. Even away from the free fermion point $\Delta=0$, there is an underlying fermionic structure in the six-vertex model; one of the best way of considering it lies in \cite{AlcarazLazo}. It consists in the rewriting of \eqref{eq:sixv:CBA} as a Matrix Ansatz
\begin{equation}\label{eq:alcarazlazorep}
	\psi(\tau_1,\ldots,\tau_P) = \Tr_{\gr{F}_K}\left( Q_K A_K(\tau_1)\ldots A_K(\tau_P)\right) 
\end{equation}
where $\tau_i=1$ if there is an up arrow on the $i$-th edge of the line and $-1$ else. The operators $A_K(+1)$ and $A_K(-1)$ are precursors of the half-strip generalized eigen-elements; in particular, they depend on the number $K$ of particles but not on the size $P$ of the cylinder. These two operators have a fermionic structure and are combination of operators $\Gamma_a(Z)$ (related to the $z_k$), $\Gamma_a(\Sigma)$ (related to the factors $\Ssixv(z_k,z_l)$) and $b^+_k$ on a fermionic Fock space $\gr{F}_K=\Gamma_a(\ca{H}_K)$. The space $\ca{H}_K$ can be identified to $l^2(I_K)$ where $I_K$ is a well chosen interval of $\setU_P$ ($P$-th roots of unity), in relation with the Fourier structure. In the limit $K,P\to\infty$, $K/P\to \rho$, it is easy to \emph{guess the space} $\gr{F}_\rho=\Gamma_a( L^2(I_\rho) )$ where $I_\rho$ is the suitable interval in $S^1$, as well as the structure of the $A_\rho(\bullet)$ but the study of the convergence has the same lack of rigour as in standard Bethe Ansatz. 

The main advance in \cite{SimonSixV} consists in the rigorous definition a morphism of $\Guill_2$-algebra from the one-particle $\Guill_2$-algebra to a suitable Fock space $\Guill_2$-algebra at density $\rho$. This morphism uses the two-particle scattering factor $\Ssixv$, the condensation curve \eqref{eq:condensedrootcurve} and has the same structure as the one in \cite{AlcarazLazo}. However, no convergence of Bethe roots is required but rather we show directly that the boundary elements are generalized eigen-elements and the infinite-volume Gibbs measure follows from Theorem~\ref{theo:eigenROPErep:invmeas} with any approximation.

\paragraph*{Some remarks}

The method presented in \cite{SimonSixV} is yet-another variation around Bethe Ansatz (coordinate, algebraic, thermodynamic, etc.) that we choose to call "operadic Bethe Ansatz". The change of geometry w.r.t.~the cylindrical one introduces new interesting elements in the corners. 

Beyond circumventing the rigorous proof of condensation of Bethe roots, a second advantage is that the same building blocks ---the boundary generalized eigen-elements--- can be used to compute correlation functions in the infinite-volume limit as illustrated in Section~\ref{sec:appli:correlfn}, again without any convergence proof from finite cylinders.

\chapter{Generalizations to continuous space and larger dimension}\label{sec:generalizations}

All the previous sections deal with discrete two-dimensional space $\setZ^2$ and finite state spaces $S_1$ and $S_2$. This encompasses already many models of statistical mechanics and it illustrates best the various concepts by keeping 2D drawings and only finite sums and finite-dimensional spaces. However, many probabilistic models require measurable state spaces for the edge variables (see for example the Gaussian case above, Section~\ref{sec:appli:gaussian}). Moreover, the guillotine operad introduced in Section~\ref{sec:operad} is defined either on $\setZ^2$ or on the continuous space $\setR^2$ and can be generalized easily in larger dimensions $d\geq 2$ with additional guillotine cuts on the new additional dimensions. The purpose of the present section is to provide a general framework that encompasses all these additional features.

We start with the algebraic introduction of the $\Guill_d$-operad in arbitrary dimension $d$ both in the discrete and the continuous space. In a second time, we introduce the state spaces suitable for the definition of the Markov processes both in the discrete and the continuous space and realize the $\Guill_d$-operad on the associated canonical spaces.

\section{Higher-dimensional guillotine operads} \label{sec:higherdim}

\subsection{Discrete Markov processes in higher dimension: geometry}

In order to introduce the notations adapted to higher-dimensional processes, we first rewrite shortly the content of Section~\ref{sec:proba} for an arbitrary  dimension $d\geq 1$.

We consider the lattice $\setZ^d$ for any $d\geq 2$. Faces are replaced by the $d$-dimensional elementary cells given by the sets $\prod_{1\leq i\leq d} [k_i,k_i+1]$ with $(k_1,\ldots,k_d)\in\setZ^d$. The boundary of such a $d$-dimensional cell is made of $(d-1)$-dimensional cells given by \[
[k_1,k_1+1]\timesdots [k_{j-1},k_{j-1}+1]\times \{l_j\} \times [k_{j+1},k_{j+1}+1] \timesdots [k_d,k_d+1]
\]
where $j$ is any integer in $\{1,\ldots,d\}$ and $l_j\in\{k_j,k_j+1\}$. More generally, the $d'$-dimensional elementary cells, which are boundaries of $(d'+1)$-dimensional cells and have $(d'-1)$-dimensional cells as boundaries, are described by a subset $J\subset \{1,\ldots,d\}$ and a point of the lattice $\mathbf{k}=(k_1,\ldots,k_d)\in\setZ^d$ such that the cell is the Cartesian product of the $d$ sets $\{k_i\}$ if $i\in J$ and $[k_i,k_i+1]$ if i$\notin J$. 

An elementary $d'$-cell $c$ described by $(J,\mathbf{k})$ with $|J|=d-d'$ has $2d'$ boundary $d'-1$-cells obtained by considering, for any $j\notin J$ and $\delta\in\{0,1\}$, the cell $(J\cup\{j\},\mathbf{k}+ \delta e_j)$ where $(e_k)$ is the canonical basis of $\setR^d$ and which is noted $ b_{j,\delta}(c)$.  This generalizes the notion of faces, boundary edges and points seen in dimension two.

2D rectangles are now replaced by $d$-boxes: a $d$-box in $\setR^d$ is a Cartesian product $R=\prod_{1\leq i \leq d} [a_i,b_i]$ with $a_i<b_i$ integers. The \emph{shape} of the box is defined as the vector $(L_1,\ldots,L_d)$ with $L_i = b_i-a_i$. The boundary of this box $\partial R$ is the set of size $2d$ of all the $(d-1)$-boxes obtained by substituting $\{a_i\}$ or $\{b_i\}$ to one of the segment $[a_i,b_i]$. For a box $R$, the set of all the $L_1\ldots L_d$ elementary $d$-cells included in $R$ is written $\Cell(R)$.

Every $(d-1)$-cell that is part of the boundary of one of these $d$-cells either appears twice as a boundary of $d$-cells (in this case, it is said \emph{interior} to the box $R$) or belongs to the boundary of the box (in this case, it is said \emph{exterior}). The set of all the $(d-1)$-cells included in a box $R$ will be written $\Wall(R)$. The subset of all the $(d-1)$-cells included in a hyperplane $\setR^{j-1}\times \{k\} \times \setR^{d-j}$ for some $k\in\setZ$ will be written $\Wall_i(R)$. 

\begin{defi}\label{def:higherdim:markovlaw}
	Let $(S_1,\ldots,S_d)$ be $d$ finite sets. Let $R = \prod_{1\leq i \leq d} [a_i,b_i]$ be a $d$-box in $\setZ^d$. A $(S_1,\ldots,S_d)$-valued Markov process on $R$ is a process $(X_u)_{u\in\Wall(R)}$ with values in $\sqcup_{i} S_i$ such that:
	\begin{enumerate}[(i)]
		\item for any $1\leq i\leq d$, for any $u\in\Wall_i(R)$, $X_u\in S_i$ a.s.,
		\item there exist weights $(\MarkovWeight{A}_c)_{c\in\Cell(R)}$ with \[
		\MarkovWeight{A}_c  : S_1\times S_1\timesdots S_d\times S_d \to \setR_+,
		\]
		such that, conditionally on the boundary $\sigma$-algebra $\sigma(X_v; v \in \partial R)$, for any sequence $(x_u)_{u\in\Wall(R)}$ in $\cup_i S_i$,
		\begin{equation}\label{eq:anyd:markovlaw}
			\probc{ \prod_{u\in\Wall(R)} \indic{X_u=x_u} }{ X_v; v\in\partial R}
			= \frac{1}{Z(\MarkovWeight{A}_\bullet; (x_v)_{v\in\partial R} )} \prod_{c\in\Cell(R)} \MarkovWeight{A}_c(x_{\partial c})
		\end{equation}
		with the short notation \[
		x_{\partial c}=\left(x_{b_{1,0}(c)},x_{b_{1,1}(c)},\ldots,x_{b_{d,0}(c)},x_{b_{d,1}(c)} \right)
		\]
	\end{enumerate}
\end{defi}

In the case $d=1$ and $d=2$, we recover our definitions \ref{def:MarkovOneDimFirst} and \ref{def:MarkovTwoDimFirst} and we leave as an exercise to the reader to check that the law is equivalent to the Markov property as soon as $|\Cell(R)|\geq 3$. All the remarks made for $d\in \{1,2\}$ still hold, in particular the fact that every variable $x_u$ appears exactly twice and the fact that the partition function $Z$ is the sum of the products of the weights over all the possible values of the $(x_u)_{u\in\Cell(R)\setminus \partial R}$.

The law of such a Markov process $(X_u)_{u\in\Wall(R)}$ is then parametrized by the (unique up to diagonal gauges) weight $\MarkovWeight{A}_c$ and a boundary weight \[
g_R : S_1^{\partial R\cap \Wall_1(R)}  \times S_2^{\partial R\cap \Wall_1=2(R)} \timesdots S_d^{\partial R\cap \Wall_d(R)} \to  \setR_+ 
\]
For homogeneous processes (i.e. there exists $\MarkovWeight{A}$ such that, for all $c$, $\MarkovWeight{A}_c=\MarkovWeight{A}$), the same difficulties of construction of a translation invariant Gibbs measure of the full lattice $\setZ^d$ as in dimension two arise and maybe solved in the same way if one finds a coherent set of weights $\ha{g}_{L_1,\ldots,L_d}$ such that $g_R = \ha{g}_{L_1,\ldots,L_d}$ for any box $R$ with shape $(L_1,\ldots,L_d)$. We claim that boundary weights $\ha{g}_{L_1,\ldots,L_d}$ built out of local objects given by suitable extensions of a $\Guill_d$-operad and suitable notions of boundary eigen-objects up to morphisms as in Section~\ref{sec:invariantboundaryelmts} help to define infinite-volume Gibbs measure as in Theorem~\ref{theo:eigenROPErep:invmeas} and compute exactly free energy densities.

\subsection{Higher-dimensional guillotine partitions}

We mimic the cases of dimensions one and two to produce the following definitions. 
\begin{defi}
	Given a box $R=\prod_{1\leq i\leq d} [a_i,b_i] \subset \setP^d$ with shape $(b_1-a_1,\ldots,b_d-a_d)\in\setL^d$, a \emph{guillotine cut} is an element $(j,c)\in\{1,\ldots,d\}\times\setR$ with $a_j\leq c_j \leq b_j$ and it defines an \emph{elementary guillotine partition} $(R_1,R_2)$ with \begin{align*}
		R_1 &= \left(\prod_{1\leq i\leq j-1} [a_i,b_i] \right)\times [a_j,c] \times \left(\prod_{j+1\leq i\leq d} [a_i,b_i]  \right)\\
		\\
		R_2 &=\left(\prod_{1\leq i\leq j-1} [a_i,b_i] \right) \times [c,b_j] \times \left(\prod_{j+1\leq i\leq d} [a_i,b_i] \right)
	\end{align*}
\end{defi}
\begin{defi}
	A \emph{guillotine partition} of a box $R$ is a sequence of boxes $(R_1,\ldots,R_n)$ such that:
	\begin{enumerate}[(i)]
		\item if $n=1$, $R_1=R$;
		\item if $n=2$, $(R_1,R_2)$ is an elementary guillotine partition of $R$;
		\item if $n>2$, there exists $\sigma\in\Sym{n}$, an elementary guillotine partition $(R'_1,R'_2)$ of $R$ and an integer $1\leq k <n$, such that $(R_{\sigma(1)},\ldots,R_{\sigma(k)})$ is a guillotine partition of $R'_1$ and $(R_{\sigma(k+1)},\ldots,R_{\sigma(n)})$ is a guillotine partition of $R'_2$.
	\end{enumerate}
\end{defi}
\begin{defi}
	The coloured operad $\Guill_d$ is the operad with colour palette ($\BoxShapes$ stands for "Box Shape") $\BoxShapes_d = \setL^d$ where the set of $n$-ary operations $\Guill_d(c;c_1,\ldots,c_n)$ are the guillotine partitions \emph{up to translations} of a box with integer sizes $c=(L_1,\ldots,L_d)$ into boxes of integer sizes $c_1,\ldots,c_n$. The compositions are given by the nesting of partitions and the identities in $\Guill_d(c;c)$ are the trivial partitions of size $1$.
\end{defi}

All the proofs are direct generalizations of the two-dimensional case. The recursive definition of guillotine partitions again helps to identify a set of generators ---~the elementary guillotine partitions~---  together with their "associativity" conditions ---~successive cuts leading to the same partitions. The generators are the multiplications \begin{align*}
	m_j^{\mathbf{L},\mathbf{L}'} & \in \Guill_d(\mathbf{L}'';\mathbf{L},\mathbf{L}')\\
	\text{with }& \mathbf{L} = (L_1,\ldots,L_{j-1},L_j^{(1)},L_{j+1},\ldots,L_d) \\
	& \mathbf{L'} = (L_1,\ldots,L_{j-1},L_j^{(2)},L_{j+1},\ldots,L_d) \\
	& \mathbf{L''} = (L_1,\ldots,L_{j-1},L_j^{(1)}+L_j^{(2)},L_{j+1},\ldots,L_d) 
\end{align*}
As before, as long as the box dimensions can be guessed from context, these multiplications will be written simply $m_j$ where $j$ indicates in which dimension the gluing of two successive boxes takes place. 

The associativity conditions can be organized recursively on the dimension $d$: for any strict subset $J\subset \{1,\ldots,d\}$ and any sequence $(l_k)_{k\notin J}$, the multiplications $(m_j^{\mathbf{L},\mathbf{L}'})_{j\in J}$ with $L_k=L'_k=l_k$ for any $k\notin J$ generate a $\Guill_{|J|}$-operad: this corresponds to guillotine partitions with cuts $(j,c)$ for $j\in J$ only. Besides all the associativity conditions obtained for these subsets, there are additional overall associativity conditions corresponding to $d$ cuts with one in each directions executed in any order. More precisely, for any $\mathbf{L}$ and $\mathbf{L}'$ in $\BoxShapes_d$, the guillotine partition of the box $\prod_{1\leq i\leq d} [0,L_i+L'_i]$ into $2^d$ boxes $\prod_{1\leq i\leq d} [l_i,r_i]$ with $(l_i,r_i)=(0,L_i)$ or $(l_i,r_i)=[L_i,L'_i]$ is made of $d$ guillotine cuts $(j,L_j)$ with $1\leq i\leq d$ which can be performed in any order.

\begin{theo}[canonical linear algebra structure on hyper-cubic lattices] \label{theo:higherdim:canostruct}We consider the discrete setting $(\setP,\setL)=(\setZ,\setL)$ and a fixed dimension $d\geq 1$.
	Let $(V_i)_{1\leq i\leq d}$ be $d$ finite-dimensional vector spaces. For any $\mathbf{L}=(L_1,\ldots,L_d)\in\BoxShapes_d$, we define $\vol(\mathbf{L})=L_1\ldots L_d$ and $\vol_{i}=\prod_{j\neq i} L_j$ and the spaces
	\begin{equation}
		\label{eq:anyd:canonicalspace}T_{\mathbf{L}} = \otimes_{1\leq i\leq d} \End(V_i)^{\otimes \vol_i(\mathbf{L})}.
	\end{equation} By considering the operator product on each algebra and the concatenation products, the spaces $(T_{\mathbf{L}})_{\mathbf{L}\in\BoxShapes_d}$ can be endowed with a structure of algebra over the $\Guill_d$-operad.
\end{theo}

As it will be seen below, the partition functions on boxes of a Markov process on the faces of the hypercubes in $\setZ^d$ with values in finite spaces $S_i$ appearing in~\eqref{eq:anyd:markovlaw} can be seen as elements of the spaces $T_\mathbf{L}$ where $\mathbf{L}$ is the shape of the boxes. For any guillotine partition $(B_1,\ldots,B_n)$ of a box $B$, it holds $Z_B= m_{(B_1,\ldots,B_n)}(Z_{B_1},\ldots,Z_{B_n})$ where $Z_{B}$ is the partition function on the box $B$ built from the weights $(A_c)_{c\in\Cell(S)}$ of the Markov process.

\subsection{Boundary extensions of the guillotine operad in dimension $d$}

We may extend the guillotine operad $\Guill_d$ in any of the $2d$ canonical directions by considering guillotine partitions of infinite boxes $\prod_{1\leq i\leq d} [a_i,b_i]$ with $a_i\in \{-\infty\}\cup\setP$, $b_i\in \setP\cup\{+\infty\}$ with $a_i<b_i$. If, for a given direction, both $a_i$ and $b_i$ are infinite, we may consider pointed guillotine partitions in order to lift the translational ambiguity in the definition of the operads. There are much more shapes than in figure~\ref{fig:admissiblepatterns} and thus we will describe here directly $\Guill_d^{(\patterntype{f})}$ where $\patterntype{f}$ stands for \emph{full} space and corresponds to all the shapes generated by guillotine partitions of $\setP^d$.

We introduce the extended set of colours ($\BoxShapes$ stands for "Box Shapes") \[
\BoxShapes_d^{\patterntype{f}} = \left( \setL^* \cup \{ \infty_L,\infty_R,\infty_{LR} \} \right)^d
\] which corresponds to shapes of boxes obtained by cutting the full space into elementary boxes through guillotine cuts. The value $\infty_L$ (resp. $\infty_R$, $\infty_{LR}$) in a direction corresponds to an interval $(-\infty,b]$ with finite $b$ (resp. $[a,+\infty)$ with finite $a$, $(-\infty,+\infty)$). The extended guillotine operad $\Guill_d(\patterntype{f})$ is defined as the equivalence class under translations of guillotine partitions of such boxes (with a base point in each direction with a colour $\infty_{LR}$).

An algebra $(\ca{A}_c)_{c\in\BoxShapes_d^{\patterntype{f}}}$ over $\Guill_d^{\patterntype{f}}$ is thus made of:
\begin{itemize}
	\item an algebra $(\ca{A}_c)_{c\in\BoxShapes_d}$ over $\Guill_d$
	\item \emph{boundary spaces} $\ca{A}_c$ when $c$ contains at least one infinity value among $\{ \infty_L,\infty_R,\infty_{LR} \}$
\end{itemize}
The key remark is that boundary spaces on each side are themselves sub-$\Guill^{\patterntype{f}}_{d'}$-algebra. 

\begin{prop}
	Given a fixed subset $J\subset\{1,\ldots,d\}$ with $m$ elements $j_1<j_2<\ldots<j_m$ and a sequence $(l_j)_{j\in J}$ in $\{\infty_L,\infty_R\}$, the canonical inclusion
	\begin{align*}
		\iota_J :\BoxShapes_{d-m}^{\patterntype{f}}  & \to \BoxShapes_{d}^{\patterntype{f}}	\\
		(L_1,\ldots,L_{d-m}) & \mapsto (L_1,\ldots,L_{j_1-1},l_{j_1}, L_{j_1}, \ldots, L_{j_1+j_2-2},l_{j_2},\ldots)
	\end{align*}
	induces a $\Guill_{d-m}^{\patterntype{f}}$-algebra structure on the spaces $(\ca{A}_{\iota_J(c)})_{c\in\BoxShapes_{d-m}^{\patterntype{f}}}$.
\end{prop}
From a geometric point of view, this corresponds to avoid guillotine cuts in the directions $j\in J$ on the boxes with an infinity in the directions $j\in J$.

Applications to Markov processes are direct generalizations from our construction of Section~\ref{sec:invariantboundaryelmts} that we briefly summarize:
\begin{itemize}
	\item canonical spaces $(T_\mathbf{L})_{\mathbf{L}\in\BoxShapes_d^{\patterntype{f}}}$ can be defined through the same constructions as in Section~\ref{sec:canonicalboundarystructure}
	\item higher dimensional generalizations of ROPEs and ROPEreps can be defined by putting a trivial $\Guill_d$-algebra structure on $\setC$ (all the products are the standard commutative product) for finite boxes, which acts multiplicatively on boundary spaces and we now call them \emph{$\Guill_d$-operadic product state} (ROPEreps are $\Guill_2$-operadic product states)
	\item the stability theorem~\ref{theo:stability} holds for inclusion of boxes with boundary weights built from $\Guill_d$-operadic product states; 
	\item for each of the $2d$ planar parts of the boundary $\partial R$ of a box $R$, the same definition of eigen-element up to morphisms as in Definition~\ref{def:eigenalgebrauptomorphims} holds excepted that the morphisms are now morphisms of higher $\Guill^{\patterntype{f}}_{d-1}$-algebras, instead of morphisms of $\Guill_1$-algebras;
	\item in dimension two, corner spaces do not have any internal structure and carry only actions of boundary algebras; in larger dimensions, corner spaces are replaced by $\Guill^{\patterntype{f}}_{d-2}$-algebras on the intersection of two hyperplane parts of $\partial R$ and have now their own internal structure if $d>2$ and actions of the two incident $\Guill^{\patterntype{f}}_{d-1}$-algebras, with the suitable notion of morphisms for eigen-elements;
	\item there is now a whole tower of lower dimensional objects that are $\Guill^{\patterntype{f}}_{d-k}$-algebras on the intersections of the higher dimensional objects, with their own definitions of eigen-elements up to morphisms.
\end{itemize}
One may still list, as an painful exercise (see below), in dimension three, all the shapes and the all requirements for eigen-structures on the boundary $\Guill_2$-, $\Guill_1$- and $\Guill_0$-algebras but, without any powerful indexation sets and a more complete operadic toolbox, it becomes nearly impossible without a computer or a deeper understanding of the operadic up-to-morphisms definitions, to write all the equations needed to compute all the eigen-elements up to morphisms. We however believe that such a formalism is needed to describe properly boundary weights in larger dimensions. All the notions of morphisms and definitions of boundary eigen-spaces up to morphisms can be generalized straightfully: the only difficulty is the enumeration of equations to solve.

The easiest part in raising the dimension is the geometric one, i.e. describing the shapes and the palette of colours with the associated notions of volumes, areas, etc. (i.e. the metric notions). In particular, the partition functions of an homogeneous Markov process with cell weight $\MarkovWeight{W}$ on a box $B$ with size $(L_1,\ldots,L_d)$ with a boundary weight $g^\Lambda$ built out of completed compatible eigen-elements up to morphisms with eigenvalue sequences $\Lambda$ and $\sigma_{J,b}$ is exactly given by:
\begin{equation}\label{eq:partitionfunc:eigenvalues:gene}
	Z_{B}^{\boundaryweights}(\MarkovWeight{W}; g^\Lambda)
	= \Lambda^{L_1 \ldots L_d} \prod_{k=2}^{d}\prod_{\substack{J\subset \{1,\ldots,d\}\\ |J|=d-k }}\prod_{b\in\{R,L\}} \sigma_{J,b}^{\prod_{j\in J} L_j}
\end{equation}
where the coefficients $\sigma_{J,b}$ are associated to each sub-$\Guill_{d-k}$-algebra living on the corresponding $(d-k)$-dimensional boundary with one infinite size on the direction with $j\in J$ and generalize the coefficients $\sigma_{a}$, $a\in\{S,N,W,E\}$ and $\kappa$ (for $k=d$) of Corollary~\ref{coro:fullplane:eigen} in dimension two.

Theorem~\ref{theo:1D:removingcoloursonboundaries} shows how the introduction of a boundary removes colours in the corresponding direction. In dimension two, the boundary is one-dimensional with zero-dimensional corners and we are thus left with traditional linear algebra on matrices and vectors. Already in dimension three, the two-dimensional boundaries produce eigen-$\Guill_2$-algebras up to morphisms, which are already new structures without so many tools.

We think that, from a fundamental perspective and for the development of abstract tools, it would be more interesting to reformulate most constructions of the present paper using higher algebra tools as in \cite{LurieHA,CostelloGwilliam} in which one should add a colour palette and the associated metric notions on the underlying geometric category. Second, the way we have formulated our results in dimension two is very "computational", with the aims of solving equations and obtaining numbers: without a practical computation in mind, it does not deserve to write all the eigen-element equations; for such practical purposes, it would be interesting to have a combinatorial approach to the set of equations to write.

\subsection{A example of structure in dimension three}

We consider here the case $S_1=S_2=S_3=S$ of finite sets. $S$-valued random variables are associated to the (two-dimensional) faces of the elementary cubes of $\setZ^3$. A face weight is thus a function $\MarkovWeight{W}:S^6 \to \setR_+$ (faces oriented for example along: left, right, front, back, down, up). The face weights have a $\Guill_3$-structure inherited from the gluings along faces. Around a box of size $(p,q,r)$, a boundary configuration consists in six two-dimensional arrays $(x^{(a)}_{i,j})$, $a \in \{L,R,F,B,D,U\}$ with respective sizes $qr$, $qr$, $pr$, $pr$, $pq$ and $pq$.

Eigen-boundary conditions are then described by a collection $6|S|+12+8$ objects organized as follows:
\begin{itemize}
	\item $6|S|$ operators $A^{(a)}(x)$, $x\in\ S$ associated to the elementary faces on the side $a$, which, for each $a$ have a $\Guill_2$-structure.
	\item 12 operators $U^{(ab)}$ associated to each edge on the boundary between two faces $a$ and $b$, which, for each $(a,b)$, have a $\Guill_1$-structure.
	\item 8 octant operators $V^{(abc)}$ associated to the corners of the cube at the intersection between the three faces $a$, $b$ and $c$.
\end{itemize}
For example, a  weight $\MarkovWeight{W}$, a down element $A^{(D)}(x)$, a back element $A^{(B)}(y)$, a left-back edge element $U^{(BD)}$ and a left-back-down corner element $V^{(LBD)}$ are associated to the following shapes :
\[
\begin{tikzpicture}[guillpart,yscale=1.7,xscale=1.7]
	\draw[guillsep] (0,0,0)--(1,0,0)--(1,1,0)--(0,1,0)--cycle;
	\draw[guillsep] (0,0,1)--(1,0,1)--(1,1,1)--(0,1,1)--cycle;
	\draw[guillsep] (0,0,0)--(0,0,1);
	\draw[guillsep] (1,0,0)--(1,0,1);
	\draw[guillsep] (0,1,0)--(0,1,1);
	\draw[guillsep] (1,1,0)--(1,1,1);
	\node at (0.5,0.5,0.5) {$\MarkovWeight{W}$};
	\fill[gray,opacity=0.2] 
	(0,0,0)--(1,0,0)--(1,1,0)--(0,1,0)--cycle;
	\fill[gray,opacity=0.2]
	(0,0,1)--(1,0,1)--(1,1,1)--(0,1,1)--cycle;
	\fill[gray,opacity=0.2]
	(0,0,0)--(1,0,0)--(1,0,1)--(0,0,1)--cycle;
	\fill[gray,opacity=0.2]
	(0,1,0)--(1,1,0)--(1,1,1)--(0,1,1)--cycle;
	\fill[gray,opacity=0.2]
	(0,0,1)--(0,1,1)--(0,1,0)--(0,0,0)--cycle;
	\fill[gray,opacity=0.2]
	(1,0,1)--(1,1,1)--(1,1,0)--(1,0,0)--cycle;
\end{tikzpicture}
,\quad
\begin{tikzpicture}[guillpart,yscale=1.7,xscale=1.7]
	\draw[guillsep] (0,1.5,0)--(1,1.5,0)--(1,1.5,1)--(0,1.5,1)--cycle;
	\draw[guillsep] (0,1.5,0)--(0,0,0);
	\draw[guillsep] (1,1.5,0)--(1,0,0);
	\draw[guillsep] 	(0,1.5,1)--(0,0,1);
	\draw[guillsep] 	(1,1.5,1)--(1,0,1);
	\fill[gray,opacity=0.2] (0,1.5,0)--(1,1.5,0)--(1,1.5,1)--(0,1.5,1)--cycle;
	\fill[gray,opacity=0.2] 
	(0,1.5,0)--(0,1.5,1)--(0,0,1)--(0,0,0)--cycle;
	\fill[gray,opacity=0.2] 
	(1,1.5,0)--(1,1.5,1)--(1,0,1)--(1,0,0)--cycle;
	\fill[gray,opacity=0.2] 
	(0,1.5,0)--(1,1.5,0)--(1,0,0)--(0,0,0)--cycle;
	\fill[gray,opacity=0.2] 
	(0,1.5,1)--(1,1.5,1)--(1,0,1)--(0,0,1)--cycle;
	\node at (0.5,0.75,0.5) {$A^{(D)}(x)$};
\end{tikzpicture},
\quad
\begin{tikzpicture}[guillpart,yscale=2,xscale=2]
	\draw[guillsep] (0,0,0)--(1,0,0)--(1,1,0)--(0,1,0)--cycle;
	\draw[guillsep,] (0,0,0)--(0,0,-1.5);
	\draw[guillsep,] (1,0,0)--(1,0,-1.5);
	\draw[guillsep,] 	(0,1,0)--(0,1,-1.5);
	\draw[guillsep,] (1,1,0)--(1,1,-1.5);
	\fill[gray,opacity=0.2] (0,0,0)--(1,0,0)--(1,1,0)--(0,1,0)--cycle;
	\fill[gray,opacity=0.2] (1,0,0)--(1,0,-1.5)--(1,1,-1.5)--(1,1,0)--cycle;
	\fill[gray,opacity=0.2] (0,0,0)--(0,0,-1.5)--(0,1,-1.5)--(0,1,0)--cycle;
	\fill[gray,opacity=0.2] (0,1,0)--(1,1,0)--(1,1,-1.5)--(0,1,-1.5)--cycle;
	\fill[gray,opacity=0.2] (0,0,0)--(1,0,0)--(1,0,-1.5)--(0,0,-1.5)--cycle;
	\node at (0.5,0.5,-0.75) {$A^{(B)}(y)$};
\end{tikzpicture},
\quad 
\begin{tikzpicture}[guillpart,yscale=2,xscale=2]
	\draw[guillsep] (0,0,-1.5)--(0,0,0)--(1,0,0)--(1,0,-1.5) 
	(0,0,0)--(0,-1.5,0) (1,0,0)--(1,-1.5,0);
	\fill[gray,opacity=0.2] (0,0,0)--(0,0,-1.5)--(0,-1.5,-1.5)--(0,-1.5,0)--cycle;
	\fill[gray,opacity=0.2] (1,0,0)--(1,0,-1.5)--(1,-1.5,-1.5)--(1,-1.5,0)--cycle;
	\fill[gray,opacity=0.2] 
	(0,0,0)--(1,0,0)--(1,0,-1.5)--(0,0,-1.5)--cycle;
	\fill[gray,opacity=0.2]
	(0,0,0)--(1,0,0)--(1,-1.5,0)--(0,-1.5,0)--cycle;
	\node at (0.5,-0.75,-0.75) {$U^{(BD)}$};
\end{tikzpicture},
\quad
\begin{tikzpicture}[guillpart,yscale=2,xscale=2]
	\fill[gray,opacity=0.2] 
	(0,0,0) -- (0,-1.5,0)--(-1.5,-1.5,0)--(-1.5,0,0)--cycle;
	\fill[gray,opacity=0.2] 
	(0,0,0) -- (0,-1.5,0)--(0,-1.5,-1.5)--(0,0,-1.5)--cycle;
	\fill[gray,opacity=0.2] 
	(0,0,0) -- (-1.5,0,0)--(-1.5,0,-1.5)--(0,0,-1.5)--cycle;
	\draw[guillsep] (0,0,0)--(-1.5,0,0) (0,0,0)--(0,-1.5,0) (0,0,0)--(0,0,-1.5);
	\node at (-0.75,-0.75,-0.75) {$V^{(LBD)}$};
\end{tikzpicture}
\]
A face weight $\MarkovWeight{W}$ in a $\Guill_3$-algebra acting on the down $\Guill_2$-algebra corresponds to vertical gluings of $n$ face weights
\[
\ti{A}^{(D)}_n=\begin{tikzpicture}[guillpart,yscale=1.5,xscale=1.5]
	\draw[guillsep] (0,2,0)--(1,2,0)--(1,2,-1)--(0,2,-1)--cycle;
	\draw[guillsep] (0,1,0)--(1,1,0)--(1,1,-1)--(0,1,-1)--cycle;
	\draw[guillsep] (0,0,0)--(1,0,0)--(1,0,-1)--(0,0,-1)--cycle;
	\draw[guillsep] (0,-1,0)--(1,-1,0)--(1,-1,-1)--(0,-1,-1)--cycle;
	\draw[guillsep] (0,2,0)--(0,-2.5,0) (1,2,0)--(1,-2.5,0) (1,2,-1)--(1,-2.5,-1) (0,2,-1)--(0,-2.5,-1);
	\fill[gray,opacity=0.2] 
	(0,2,0)--(1,2,0)--(1,2,-1)--(0,2,-1)--cycle;
	\fill[gray,opacity=0.2] 
	(0,1,0)--(1,1,0)--(1,1,-1)--(0,1,-1)--cycle;
	\fill[gray,opacity=0.2] 
	(0,0,0)--(1,0,0)--(1,0,-1)--(0,0,-1)--cycle;
	\fill[gray,opacity=0.2] 
	(0,-1,0)--(1,-1,0)--(1,-1,-1)--(0,-1,-1)--cycle;
	\fill[gray,opacity=0.2]
	(0,2,0)--(1,2,0)--(1,-2.5,0)--(0,-2.5,0)--cycle;
	\fill[gray,opacity=0.2]
	(0,2,-1)--(1,2,-1)--(1,-2.5,-1)--(0,-2.5,-1)--cycle;
	\fill[gray,opacity=0.2]
	(0,2,0)--(0,2,-1)--(0,-2.5,-1)--(0,-2.5,0)--cycle;		
	\fill[gray,opacity=0.2]
	(1,2,0)--(1,2,-1)--(1,-2.5,-1)--(1,-2.5,0)--cycle;
	\node at (0.5,-0.5,-0.5) {$\MarkovWeight{W}$};
	\node at (0.5,0.5,-0.5) {$\vdots$};
	\node at (0.5,1.5,-0.5) {$\MarkovWeight{W}$};
	\node at (0.5,-1.75,-0.5) {$A^{(D)}$};
\end{tikzpicture}
\]
A volume power of the weight $\MarkovWeight{W}^{[p,q,r]}$ acts on the surface power $(A^{(D)})^{[p,q]}$ and, for example, we have by associativities that is equal to $\MarkovWeight{W}^{[p,q,r-k]}$ acting on $(\ti{A}_k^{(D)})^{[p,q]}$. A definition similar to Definition~\ref{def:eigenalgebrauptomorphims} corresponds to the existence of a morphism of $\Guill_2$-algebra (the one on the down side) that maps $(\ti{A}_n^{(D)})^{[p,q]}$ to $\Lambda^{npq} (A^{(D)})^{[p,q]}$.

\begin{rema}
	This is a first identified algebraic difficulty in higher dimension: whereas "standard" algebras are put on the edges of a rectangle in dimension two ---~hence allowing for "standard" tools of reduction of endomorphisms and morphisms (such as block extraction in Section~\ref{sec:ROPErep:fundamentalexample})---, we now need to understand in a deeper way the structure of morphisms in $\Guill_2$-algebras (i.e. the two products are preserved). 
\end{rema}

Once the six boundary face $\Guill_2$-eigen-elements of the $\Guill_3$-weight $\MarkovWeight{W}$ are understood. Then the elements $U^{(ab)}$ have to be identified. Definition~\ref{def:cornereigensemigroups} is now replaced by the existence of morphisms of $\Guill_1$-algebra (and not just linear maps) such that, for any $p\geq 1$, $q,r\geq 0$,
\[
\begin{tikzpicture}[guillpart,yscale=2,xscale=2]
	\draw[guillsep]
	(0,0,0)--(1,0,0) 
	(0,0,-3.)--(0,0,0)
	(1,0,0)--(1,0,-3.) 
	(0,0,0)--(0,-3.,0) 
	(1,0,0)--(1,-3.,0)
	(0,-3.,-1)--(0,0,-1)--(1,0,-1)--(1,-3.,-1)
	(0,-1,-3.)--(0,-1,0)--(1,-1,0)--(1,-1,-3.)
	(0,-1,-1)--(1,-1,-1);
	
	\fill[gray,opacity=0.2] 	(0,0,0)--(0,0,-3.)--(0,-3.,-3.)--(0,-3.,0)--cycle;
	\fill[gray,opacity=0.2] (1,0,0)--(1,0,-3.)--(1,-3.,-3.)--(1,-3.,0)--cycle;
	\fill[gray,opacity=0.2] 
	(0,0,0)--(1,0,0)--(1,0,-3.)--(0,0,-3.)--cycle;
	\fill[gray,opacity=0.2]
	(0,0,0)--(1,0,0)--(1,-3.,0)--(0,-3.,0)--cycle;
	
	\node at (0.5,-2.5,-2.75) {$(U^{(BD)})^{[p]}$};
	\node at (0.5,-0.5,-0.5) {$\MarkovWeight{W}^{[p,q,r]}$};
	\node at (0.5,-0.5,-2.75) {$(A^{(B)})^{[p,q]}$};
	\node at (0.5,-2.5,-0.5) {$(A^{(D)})^{[p,r]}$};
\end{tikzpicture}
\overset{\Phi}{\longmapsto}
\Lambda^{pqr} \sigma_{B}^{pq} \sigma_D^{pr}
\begin{tikzpicture}[guillpart,yscale=2,xscale=2]
	\draw[guillsep] (0,0,-1.5)--(0,0,0)--(1,0,0)--(1,0,-1.5) 
	(0,0,0)--(0,-1.5,0) (1,0,0)--(1,-1.5,0);
	\fill[gray,opacity=0.2] (0,0,0)--(0,0,-1.5)--(0,-1.5,-1.5)--(0,-1.5,0)--cycle;
	\fill[gray,opacity=0.2] (1,0,0)--(1,0,-1.5)--(1,-1.5,-1.5)--(1,-1.5,0)--cycle;
	\fill[gray,opacity=0.2] 
	(0,0,0)--(1,0,0)--(1,0,-1.5)--(0,0,-1.5)--cycle;
	\fill[gray,opacity=0.2]
	(0,0,0)--(1,0,0)--(1,-1.5,0)--(0,-1.5,0)--cycle;
	\node at (0.5,-0.75,-0.75) {$(U^{(BD)})^{[p]}$};
\end{tikzpicture}
\]
When $q=0$ or $r=0$, then the elements $\MarkovWeight{W}$ and the corresponding $A^{(a)}$ is just absent. As before, for discrete Markov processes and the canonical space, it is enough to find solutions as in Section~\ref{sec:concreteeigencorner} for $(q,r)\in\{(0,1),(1,0),(1,1)\}$ and arbitrary $p$. The combined action of the three elements $\MarkovWeight{W}$, $A^{(D)}$ and $A^{(B)}$ is just a multiplicative morphism on the $\Guill_1$-element $U^{(BD)}$, as well as the backward morphism that provides the eigenvalue: we now recover a situation similar to the one of Definition~\ref{def:eigenalgebrauptomorphims} and the associated tools inherited from linear algebra. In dimension three, the edge elements $U^{(ab)}$ combine both aspects of the side elements in dimension two ($\Guill_1$-structure) and the corner elements in dimension two.

Finally, the corner element $V^{(abc)}$ have a structure similar to a trimodule with respect to the $U^{(ab)}$, $U^{(ac)}$ and $U^{(bc)}$. The corner spaces are $\Guill_0$-algebra and thus do not have any internal structure. We only need the existence of linear maps on the corner space such that
\[
\begin{tikzpicture}[guillpart,yscale=3.,xscale=3]%z={(10mm, 10mm)}
	%face dessus
	\draw[guillsep] 
	(0,0,-3)--(0,0,0)--(-3,0,0)
	(-1,0,-3)--(-1,0,-1)--(-3,0,-1)
	(0,0,-1)--(-1,0,-1)--(-1,0,0)
	;
	\fill[gray,opacity=0.2]	
	(0,0,-3)--(0,0,0)--(-3,0,0)--(-3,0,-3)--cycle;
	%face horizontale en dessous	
	\draw[guillsep] 
	(0,-1,-3)--(0,-1,0)--(-3,-1,0)
	(-1,-1,-3)--(-1,-1,-1)--(-3,-1,-1)
	(0,-1,-1)--(-1,-1,-1)--(-1,-1,0)
	;
	\fill[gray,opacity=0.2]	
	(0,-1,-3)--(0,-1,0)--(-3,-1,0)--(-3,-1,-3)--cycle;
	%face à droite
	\fill[gray,opacity=0.2]
	(0,0,-3)--(0,0,0)--(0,-3,0)--(0,-3,-3)--cycle;
	\draw[guillsep]
	(0,0,0)--(0,-3,0)
	(0,0,-1)--(0,-3,-1)
	;
	% face verticale parallèle juste à gauche de la précédente
	\fill[gray,opacity=0.2]
	(-1,0,-3)--(-1,0,0)--(-1,-3,0)--(-1,-3,-3)--cycle;
	\draw[guillsep]
	(-1,0,0)--(-1,-3,0)
	(-1,0,-1)--(-1,-3,-1)
	;
	% face parallèle à la feuille
	\fill[gray,opacity=0.2]
	(0,0,0)--(0,-3,0)--(-3,-3,0)--(-3,0,0)--cycle;
	\fill[gray,opacity=0.2]
	(0,0,-1)--(0,-3,-1)--(-3,-3,-1)--(-3,0,-1)--cycle;
	% opérateurs:
	\node at (-0.5,-0.5,-0.5) {$\MarkovWeight{W}^{[p,q,r]}$};
	\node at (-0.5,-2,-0.5) {$(A^{(D)})^{[p,r]}$};
	\node at (-2,-0.5,-0.5) {$(A^{(L)})^{[q,r]}$};
	\node at (-0.5,-0.5,-2) {$(A^{(B)})^{[p,q]}$};
	\node at (-2,-0.5,-2) {$(U^{(LB)})^q$};
	\node at (-0.5,-2,-2) {$(U^{(BD)})^p$};
	\node at (-2,-2,-0.5) {$(U^{(LD)})^r$};
	\node at (-2,-2,-2) {$V^{(LBD)}$};
	
	\draw[<->,dashed] (0.1,0,-3) -- node [midway,anchor=east] {$q$} (0.1,-1,-3);
	\draw[<->,dashed] (0,0.1,-3) -- node [midway,anchor=north] {$p$} (-1,0.1,-3);
	\draw[<->,dashed] (-3.1,0,0) -- node [midway,anchor=south] {$r$} (-3.1,0,-1);
\end{tikzpicture}
\overset{\Phi}{\longmapsto }
\Lambda^{pqr} \sigma_{B}^{pq} \sigma_D^{pr} \sigma_L^{qr}
u_{LB}^q u_{BD}^p u_{LD}^r
\begin{tikzpicture}[guillpart,yscale=2,xscale=2]
	\fill[gray,opacity=0.2] 
	(0,0,0) -- (0,-1.5,0)--(-1.5,-1.5,0)--(-1.5,0,0)--cycle;
	\fill[gray,opacity=0.2] 
	(0,0,0) -- (0,-1.5,0)--(0,-1.5,-1.5)--(0,0,-1.5)--cycle;
	\fill[gray,opacity=0.2] 
	(0,0,0) -- (-1.5,0,0)--(-1.5,0,-1.5)--(0,0,-1.5)--cycle;
	\draw[guillsep] (0,0,0)--(-1.5,0,0) (0,0,0)--(0,-1.5,0) (0,0,0)--(0,0,-1.5);
	\node at (-0.75,-0.75,-0.75) {$V^{(LBD)}$};
\end{tikzpicture}
\]
As it can be seen, now the drawings become hard to interpret and careful description of the  3D guillotine partitions have to be provided. However, the previous drawing is a simple generalization of the 2D corner described above with more objects surrounding the 3D corner.

\begin{rema}
	Despite the length of the enumeration of all objects and eigen-constraints, the general structure in dimension $d$ is simple to guess. In order to encompass all the intermediate $(d-k)$-dimensional objects, a better description in terms of boundaries (as for simplicial sets) is required.
\end{rema}

\begin{rema}
	We have seen in Section~\ref{sec:commutuptomorph} that the notion of commutation up to morphism provides a relation with the Yang-Baxter: this is a fruit of our understanding of morphisms in linear algebra, i.e. $\Guill_1$-algebra here. In order to have a similar "integrability" structure in dimension three, it requires both to have a deeper understanding of $\Guill_2$-morphisms as well as to understand what should replace commutativity at a higher operadic level.
\end{rema}

\section[General state spaces in arbitrary dimension]{General state spaces for continuous and discrete space in arbitrary dimension}

\subsection{Discrete space and measurable state space}
Definition~\ref{def:higherdim:markovlaw} shows that raising the dimension for Markov processes with discrete space and finite state space introduces no technical novelty excepted that random variables are now attached to the faces of the elementary hypercubes $\prod_{i=1}^d [k_i,k_i+1]$. 

When moving from finite state space to measurable state space on $\setZ^d$, all the sums involved in partition functions and marginal laws can be seen as integral over the counting measure and, in general, have to be replaced by integral over a reference measure. Each set $S_i$ is replaced by a measured set $(S_i,\ca{S}_i,\nu_i)$. All associativity conditions come from Fubini's theorem for the product measure over all the elementary faces around boxes and thus $\sigma$-finite measures $\nu_i$ are required. The state space on the boundary of a box $B$ is then given by $\prod_{b\in \Wall_d(B)} S_i$ with the product $\sigma$-algebra of the $\ca{S}_i$ and the product measure of the $\nu_i$. All the operadic products are then similar to the 1D case described in Proposition~\ref{prop:markov:1Dmeasurable}.

\subsection{Generalized canonical realization of the $\Guill_d$-operad in both discrete and continuous space}

The spaces $V(S_i)^{\otimes p}$ correspond to functions $S_i^p\to\setK$ and the operators $A\in \End(V(S_i))^{\otimes p}$ correspond to functions $S_i^p\times S_i^p\to\setK$. For $\sigma$-finite measured spaces $(S_i,\ca{S}_i,\nu_i)$, we have to consider integrals w.r.t. to product measures on $S_i^p$ and thus introduce additional analytical requirements.

Given $d$ face state spaces $(S^{(i)}_c,\ca{S}^i_c,\nu^i_c)_{c\in \BoxShapes_{d-1}}$, $1\leq i\leq d$ in dimension $d$ and a box $B$ with shape $(p_1,\ldots,p_d)\in \BoxShapes_d$, we introduce the following product measured space
\[
S_{\Shell(B)}  = \prod_{k=1}^{d} \left( S^{(k)}_{(p_i)_{i\neq k})} \times S^{(k)}_{(p_i)_{i\neq k}}\right)
\]
to which belongs the $2d$-uplet of r.v. associated to the faces of $B$. The partition function can then be seen as an element $Z_B:S_{\Shell(B)} \to\setR_+$.

We now introduce the following spaces of complex functions, for any box $B$ with shape $(p_1,\ldots,p_d)\in \BoxShapes_d$:
\begin{align}
	L^+_{p_1,\ldots,p_d}&= L^+\left(  S_{\Shell(B)}\right)
	\\
	L^\infty_{p_1,\ldots,p_d}&= L^\infty\left(  S_{\Shell(B)}\right)
\end{align}
(with respect to the product measures of the corresponding $\nu_i$) which correspond to the set of, respectively, $\ov{\setR}_+$-valued measurable functions and $\setR$-valued measurable functions with a bounded essential supremum (w.r.t. the product reference measure). Only the second space is a vector space; however the first one is easier to manipulate for partition functions and positive random variables. For each $1\leq k\leq d$, we introduce suitable elementary products.

\begin{defi}\label{def:generalizedguillprod}
	Let $1\leq k\leq d$ and let $p$ and $p'$ be two elements of $\setL^*$. For any three sequences,
	 \begin{align*}
	 	\gr{p} &= (p_1,\ldots,p_{k-1},p,p_{k+1},\ldots,p_d)
	 	\\
	 	\gr{p}' &= (p_1,\ldots,p_{k-1},p',p_{k+1},\ldots,p_d)
	 	\\
	 	\gr{p}'' &= (p_1,\ldots,p_{k-1},p+p',p_{k+1},\ldots,p_d)
	 \end{align*}
	 the product $m_{\gr{p},\gr{p'}} : L^\infty_{\gr{p}}\times L^\infty_{\gr{p}'} \to L^\infty_{\gr{p}''}$ (or on the spaces $L^+_\bullet$), associated to the guillotine partition
	 \[
	 \left(  
	 	\prod_{i=1}^d [0,p_i]
	 ,
	 	\prod_{i=1}^{k-1} [0,p_i] \times [p,p+p'] \times \prod_{i=k+1}^d [0,p_i]
	 \right)
	 \]
	 with a cut on the $k$-th coordinate,
	 is defined, for any $x=(x_1,x'_1,\ldots,x_d,x'_d)$ by
	 \begin{equation}\label{eq:generalizedelemprod}
	 	m_{\gr{p},\gr{p'}}(f,g)(x_1,x'_1,\ldots,x_d,x'_d) = \int_{ S^{(k)}_{(l_i)_{i\neq k}}} f(\gamma_L(x,u))g(\gamma_R(x,u)) du
	 \end{equation}
	 in which $\gamma_L(x)$ and $\gamma_R(x)$ are given by
	 \begin{align}
	 	\gamma_L(x,u) &= \left( ( \pi_{1,p}(x_1), \pi_{1,p}(x'_1) ), \ldots, (x_k, u), \ldots, (\pi_{1,p}(x_d), \pi_{1,p}(x'_d))\right) 
	 	\\
	 	\gamma_R(x,u) &= \left( ( \pi_{p+1,p+p'}(x_1), \pi_{p+1,p+p'}(x'_1) ), \ldots, (u,x'_k), \ldots, (\pi_{p+1,p+p'}(x_d), \pi_{p+1,p+p'}(x'_d))\right) 
	 \end{align}
	 with $\pi_{i,j}(x)= (x_i,\ldots,x_j)$ for $i\leq j$.
\end{defi}

\begin{theo}[Guillotine algebra for Markov processes in continuous space] \label{theo:canonicalexampleGuill:continuous}
	The spaces $(L^+_\gr{p})_{\gr{p}\in\BoxShapes_d}$ (resp. $(L^\infty_\gr{p})_{\gr{p}\in\BoxShapes_d}$ with an additional assumption of finite mass for all the reference measures) endowed with the elementary guillotine products defined in Definition~\ref{def:generalizedguillprod} is a $\Guill_d$-algebra with a Cartesian (resp. tensor) structure.
\end{theo}

\begin{proof}
	The proof is essentially the same for Theorem~\ref{theo:canonicalexampleGuill} excepted that summations are replaced by integrals w.r.t the reference measures and the definition of face state spaces is used to transform integrals using the transfer lemma (through the coproducts $\gamma_\bullet^{(i)}$) and Fubini's theorem. 
	
	The finite mass assumption of the reference measures in the $L^{\infty}$ case is required to obtain that the result of the products~\eqref{eq:generalizedelemprod} have again a finite essential supremum.
\end{proof}

The same definition for boundary spaces can be performed using direct sums and we do not reproduce the whole theory to gain some place. It is then a simple exercise to reformulate all the results of sections~\ref{sec:operad} for measurable sets instead of finite sets.

\begin{prop}
		All the previous constructions \eqref{eq:productboundaryalgebra}, \eqref{eq:actionboundaryalgebra}, \eqref{eq:canonical:actionsoncorners} remain valid in the context of finite-mass measurable state spaces $(S^{(i)}_\bullet,\ca{S}^{(i)}_\bullet,\nu^{(i)}_\bullet)$ when the tensor product spaces $\otimes_k V(S_{i_k})$ are replaced by the spaces $L^\infty(\prod_k S_{i_k})$ or $L^\infty(A)$ (which allows for the integrals to be well-defined) where $A$ is a Cartesian product of state spaces and when all the sums over state spaces are replaced in the products by integrals with respect to the reference measures.
		
		\emph{Mutatis mutandis}, lemmata generalizing the ones of Section~\ref{sec:canonicalboundarystructure} and Theorem~\ref{theo:canonicalboundarystructure} are also valid in this context. 
\end{prop}
	
\section{Generalized boundary operator product representations in arbitrary dimension and measurable state spaces}

ROPE representations can be generalized easily to face guillotine state spaces. As for partition functions in Definition~\ref{def:higherdim:markovlaw} which cannot be reduced to elementary weights but rather form a $d$-dimensional semi-group, the definition of a ROPE representation goes along the same line and uses again the restriction coproducts of the guillotine state spaces.

Definition~\ref{def:ROPE} of a ROPE does not change in the present generalized setting: all the finite-shape spaces are equal to the field $\setR$ (or $\setC$) with products given by scalar multiplication. In dimension $d\geq 2$, a ROPE is replaced by a \emph{boundary operator product environment} with the same definition.

When juggling with boundaries, we need a description of all the embeddings of $(d-k)$-dimensional boxes in the boundary of a $d$-dimensional box. We thus introduce, for any shape  $\gr{p}=(p_1,\ldots,p_{d-k})\in\BoxShapes_{d-k}^{(\patterntype{f})}$, any increasing sequence $\gr{j}=(j_i)_{1\leq i\leq k}$ of $k$ elements in $\{1,\ldots,d\}$ and any sequence $\beta: J \to \{\infty_R,\infty_L\}$, the extended shape in $\BoxShapes_{d}(\patterntype{f})$ given by
\[
\iota_{\gr{j},\beta}(\gr{p}) = (p_1,\ldots,p_{j_1-1}, \beta(j_1), p_{j_1},\ldots, p_{j_2-1},\beta(j_2),p_{j_2},\ldots,\beta(j_k),p_{j_k},\ldots,p_{d-k})
\]

\begin{defi}[boundary operator product representation for guillotine state spaces]\label{def:ROPErep:SSS}
	Let $d$ be an arbitrary dimension larger than $1$.
	Let $(S^{(k)}_{\BoxShapes_{d-1}})$, $1\leq k\leq d$ be $d$ face $\Guill_{d-1}$-state spaces. Let $\ca{B}_\bullet$ be a boundary operator product environment with topological vector spaces allowing for integration. Let $I$ be a set and $(\gr{p}_i)_{i\in I}$ a sequence of shapes in $\BoxShapes_{d}$ and $(B_i)_{i\in I}$ associated boxes in $\setP^d$. 
	
	Let $(f_i)_{i\in I}$ be a collection of measurable functions such that $f_i: S_{\Shell(B_i)} \to \setK$. A \emph{homogeneous boundary operator product  representation} (BOPR) of the collection $(f_i)_{i\in I}$ over $\ca{B}_\bullet$ consists in:
	\begin{enumerate}[(i)]
		\item \label{item:genewordmatrixproduct}for any direction $j\in\{1,\ldots,d\}$, any shape $\gr{p}=(p_1,\ldots,p_{d-1})\in \BoxShapes_{d-1}$ of a $(d-1)$-dimensional box $B$ and any value $b\in\{\infty_R,\infty_L\}$,
		a collection of measurable functions 
		\[
		A_{\iota_{j,b}(\gr{p})} : S^{(j)}_{\gr{p}} \to \ca{B}_{\iota_{j,b}(\gr{p})}
		\]
		such that, for any guillotine partition $(B_1,\ldots,B_r)$ of $B$ with shapes $(\gr{p}_1,\ldots,\gr{p}_r)$, 
		\begin{equation}\label{eq:BOPR:morphismprop:face}
		A_{\iota_{j,b}(\gr{p})}( x )	
		=
		m_{\rho}\left( A_{\iota_{j,b}(\gr{p}_1)}( \gamma^{(j)}_{\rho}(x)_1 ),
		\ldots,
		A_{\iota_{j,b}(\gr{p}_r)}( \gamma^{(j)}_{\rho}(x)_r ) \right) 
		\end{equation}
		almost everywhere w.r.t.~the reference measures, with guillotine products in the $\Guill_{d-1}$-subalgebra $(\ca{B}_{\iota_{j,b}(\gr{p})})_{\gr{p}\in\BoxShapes_{d-1}}$ and where $\gamma_\rho^{(j)}$ is the extraction from the $n$-uplet associated to a face of $B$ to the corresponding face of the sub-box in the partition $\rho$.
		
		\item \label{item:genecorner} for any $k\in\{2,\ldots,d\}$ and and any sequence  $\gr{j}=(j_i)_{1\leq i\leq k}$ of $k$ elements in $\{1,\ldots,d\}$ and any strictly increasing sequence $\beta: J \to \{\infty_R,\infty_L\}$, a collection  of elements $(U_{\iota_{\gr{j},\beta}(\gr{p})})_{\gr{p}\in \BoxShapes_{d-k}}$ such that, for any guillotine partition $\rho=(B_1,\ldots,B_r)$ with shapes $(\gr{p}_1,\ldots,\gr{p}_r)$ of a box $B$ with shape $\gr{p}$,
		\begin{equation}
			U_{\iota_{\gr{j},\beta}(\gr{p})} = m_\rho\left(  
				U_{\iota_{\gr{j},\beta}(\gr{p}_1)},
				\ldots,
				U_{\iota_{\gr{j},\beta}(\gr{p}_r)}
			\right)
		\end{equation}
	\end{enumerate}
	such that, for all $i\in I$ and almost all $x=(x_1,x'_1,\ldots,x_d,x'_d)\in S_{\Shell(B_i)}$, 
	\begin{equation}
		\label{eq:generalizedBOPrep}
		f_i(x) = m_{\rho_{B_i}}\left(1_{\setK},
			A_{\infty_R,p_2,\ldots,p_d}(x_1),
			A_{\infty_L,p_2,\ldots,p_d}(x'_1),
			\ldots,
			A_{p_1,\ldots,p_{d-1},\infty_R}(x_d),
			A_{p_1,\ldots,p_{d-1},\infty_L}(x'_d),
			U_{\bullet}
		\right)
	\end{equation}
	where $\rho_{B_i}$ is the guillotine partition of $\setP^d$ (arbitrarily pointed) generated by the $2d$-hyperplanes containing the $2d$ faces of the box $B_i$ and all the elements $U_\bullet$ are associated to the $(d-k)$-dimensional boxes that appear in the partition $\rho_{B_i}$.
\end{defi}

 The definition involves much more complex notations than the discrete 2D case of Definition~\ref{def:ROPErep:FD} of a ROPE representation of functions on boundaries of rectangles. We now spend some lines to state the full equivalence as well as some distinctions. We have the Cartesian products $S^{(1)}_k=S_1^k$ and $S^{(2)}_k=S_2^k$. Item~\ref{item:genewordmatrixproduct} corresponds to the following identification of the four sides:
 \begin{align*}
 	 \iota_{1,\infty_{L}}(q) &= (\infty_W,q)
 	 &
 	 \iota_{1,\infty_{R}}(q) &= (\infty_E,q)
 	 \\
 	 \iota_{2,\infty_{L}}(p) &= (p,\infty_S)
 	 &
 	 \iota_{2,\infty_{R}}(p) &= (p,\infty_N)
 \end{align*}
 In the discrete case, all the elements $A_{\iota_{j,B}(p)}$ are generated through the semi-group property~\eqref{eq:BOPR:morphismprop:face} by elementary elements $A_{\iota_{j,\infty_a}(1)}(x)$ with $p=1$. This semi-group property then becomes, for example, on the South with the simpler notation $A_S(x)= A_{\iota_{2,\infty_{L}}(1)}(x)$,
 \[
 A_{\iota_{2,\infty_{L}}(p)}(x_1,\ldots,x_p) = 
 \begin{tikzpicture}[guillpart,yscale=1.5,xscale=3.]
 	\fill[guillfill] (0,0) rectangle (4,1);
 	\draw[guillsep] (0,0)--(0,1)--(4,1)--(4,0) (1,0)--(1,1) (2,0)--(2,1) (3,0)--(3,1);
 	\node at (0.5,0.5) { $A_{S}(x_1)$ };
 	\node at (1.5,0.5) { $A_{S}(x_2)$ };
 	\node at (2.5,0.5) { $\ldots$ };
 	\node at (3.5,0.5) { $A_{S}(x_p)$ };
 \end{tikzpicture}
 \]
In dimension two, item~\ref{item:genecorner} contains only the case $k=2$ with only one sequence $\gr{j}=(1,2)$ and two sequences $\beta:\{1,2\}\to\{\infty_R,\infty_L\}$, and the set of degenerate box shapes $\BoxShapes_{0}$ contains only one element $\{*\}$. We are thus left with four elements
\begin{align*}
U_{\iota_{(1,2),(\infty_L,\infty_L)}(*)} &= U_{SW}
&
U_{\iota_{(1,2),(\infty_L,\infty_R)}(*)} &= U_{NW}
\\
U_{\iota_{(1,2),(\infty_R,\infty_L)}(*)} &= U_{SE}
&
U_{\iota_{(1,2),(\infty_R,\infty_R)}(*)} &= U_{NE}
\end{align*}
The operad $\Guill_0$ is degenerate and there are no products excepted identities, so that the set of morphism requirements~\eqref{eq:generalizedBOPrep} is empty: this corresponds to the corner elements of Definition~\ref{def:ROPErep:FD}. For any shape $(p_i,q_i)$ associated to the function $f_i$ and the box $B_i=[0,p_i]\times [0,q_i]$, the partition $\rho_{B_i}$ in \eqref{eq:generalizedBOPrep} is generated by the four guillotine cuts $\{0\}\times\setR$,  $\{p_i\}\times\setR$, $\setR\times \{0\}$ and $\setR\times\{q\}$ and contains the nine boxes $B_i$ and
\begin{align*}
	[0,p]\times]-\infty,0],\;
	[0,p]\times[q,+\infty,0[,\;
	]-\infty,0]\times[0,q],\;
	[p,+\infty[\times[0,q],\;
	\\
	]-\infty,0]\times]-\infty,0],\;
	]-\infty,0]\times[q,+\infty[,\;
	[p,+\infty[\times]-\infty,0],\;
	[p,+\infty[\times[q,+\infty[
\end{align*}
The element $1_\setK$ is associated to $B_i$, elements $A_{\iota_{2,\infty_{L}}(p)}(x_1,\ldots,x_p)$ to the four next half-strips and the elements $U_{\iota_{(1,2),(\infty_a,\infty_n)}(*)}$ to the last four corners.

The novelty for dimension $d>2$ is that $(d-k)$-dimensional sub-boxes with $2\leq k\leq d-1$ acquire themselves an internal guillotine structure. It is still a mystery for us how to deal with this whole hierarchy of $\Guill_{d-k}$-structures in an efficient way.

Nonetheless, we conclude the present section by the following remark without surcharging the presentation with the complete proofs. 

\begin{rema}[generalization of the stability theorem]
The important point of the present construction is that Theorem~\ref{theo:stability} remains valid \emph{mutatis mutandis} with the previous definition of BOPR of boundary weights.
\end{rema}

\subsection{Generalized eigen-boundary elements up to morphisms}

The generalization from ROPE representation to generalized boundary operator product representations ---i.e.~from Definition~\ref{def:ROPErep:FD} to Definition~\ref{def:ROPErep:SSS}--- has to be put in front of the expected formula~\eqref{eq:partitionfunc:eigenvalues:gene}. The eigenvalue $\Lambda$ is attached to the elements $A_\bullet$ in the first item in Definition~\ref{def:ROPErep:SSS} seen as eigen-elements up to morphisms of the box partition functions. The hierarchy of boundary eigen-values $\sigma_{J,b}$ in \eqref{eq:partitionfunc:eigenvalues:gene} corresponds to the hierarchy of the $(d-k)$-dimensional sub-boxes of a $d$-dimensional box and the associated $\Guill_{d-k}$-elements $U_{\bullet}$ in the second item of Definition~\ref{def:ROPErep:SSS}.

All the definitions of Section~\ref{sec:invariantboundaryelmts} can be generalized to the boundary operator product representations of Definition~\ref{def:ROPErep:SSS} in order to reach a theorem equivalent to Theorem~\ref{theo:eigenROPErep:invmeas} in arbitrary dimension $d$. However, much heavier notations are needed to describe the hierarchy of morphisms in the definition of the eigen-elements $U_{\bullet}$. 

Indeed, in dimension two, equations on the $A_a$ depend only the face weight and $A_a$, whereas equations on the corner elements $U_{ab}$ depend on the face weight and the neighbouring elements $A_a$ and $A_b$. In higher dimension, the element $U_{a}$ associated to a $(d-k)$-dimensional box $B'$ included the boundary of a $d$-dimensional box $B$ is expected to satisfy generalized eigen-element equations up to morphisms that involve the box partition function on $B$ as well as all the $U_{c}$ associated to a $(d-l)$-dimensional box with $l\leq k$ included in the boundary of $B$ and containing $B'$.

It is still possible to list all the requirements in dimension $3$ or $4$ to validate a generalization of Theorem~\ref{theo:eigenROPErep:invmeas}. There are however many objects with many constraints: listing all of them is possible since they all have the same structure; it is even possible to formulate it in detail. However, solving the hierarchy of equations require a deeper and more systematic understanding of the hierarchical structure of the boundary operator product representation through eigen-elements and we postpone the general definition to a future work involving concrete computational tools.

As a long exercise left to the reader, it is however possible to check that these constructions can be performed exactly in the Gaussian framework of Section~\ref{sec:appli:gaussian} in arbitrary dimension $d$, hence providing a proof of concept of the utility of the present section.

\chapter{Conclusion and perspectives}\label{sec:openquestions}

\section{A brief summary}

We have shown in this paper how the use of suitable coloured operads can give firm algebraic grounds to the probabilistic Markov property and leads to a bunch of new notions both in algebra (with boundary extensions of our operads) and in probability with a local-to-global construction of invariant boundary rates. In particular, the colour palette of the coloured operads plays a strong role in many properties. A major novelty of the present paper is the new definition of eigen-elements for operads "with boundaries", which generalizes the classical definition of eigenvectors using equations \emph{up to morphisms}.

\subsection{From geometry and probability to algebra}

The Markov property is an example of an incursion of geometry into probability theory. The basic ingredient of gluing of domains (see Section~\ref{sec:operad}) with its associativities is then pushed to the algebraic side by describing how the probabilistic parameters are combined under gluings. This is done by the operadic formalism. 

Keeping only rectangles and guillotine partitions allows one to reduce the structure of gluings to only two products with three associativities and a palette of colours. It combines in a minimal way (square/interchange associativity) the two one-dimensional products associated to each dimension.

The geometry of the operad structure induces a natural theory of boundary elements, generalizing the left and right modules in associative algebras. The two-dimensionality induces four side module-like structure on sides of rectangles ---related to standard associativity---, as well as new corner structures on the vertices of rectangles, for which the square/interchange associativity plays a major role.

The advantage of the present construction is that the geometric intuition behind the algebraic definitions provides a whole set of natural definitions. In particular, one sees how the geometric dimensions of each object dictates its algebraic structure: matrices on 1D edges, four-legs tensors on elementary square faces, $(2p+2q)$-legs tensors on rectangle with sizes $(p,q)$, etc.

In contrast with the one-dimensional situation in which boundaries are described by vectors ---without any internal  
structure---, the ROPEreps are made of operators with a product structure associated to the one-dimensional gluing structure of edges. The emergence of the internal structure precisely opens the way to definitions "up to morphisms" as in Section~\ref{sec:invariantboundaryelmts}. These definitions are also minimal in the sense that already trivial models requires the presence of such morphisms.

\subsection{Kolmogorov's extension theorem: an algebraic alternative to the analytical method}

On the probabilistic side, existence of infinite-volume Gibbs measures is obtained by using Kolmogorov's extension theorem (see Theorem~\ref{theo:eigenROPErep:invmeas}). Up to now, this theorem was known to be difficult to apply in dimension larger that one since it requires a precise description of boundary weights. Instead, the standard approach is based on analysis: local observables studied on finite domains with \emph{prescribed well-chosen} boundary conditions, large size limit and Riesz representation theorem. The main problem is the requirement for a physical or intuitive input on the boundary conditions or a proof of its irrelevance.

The present paper rewrites the set of constraints required by Kolmogorov's extension theorem into constraints on a finite set of $2(|S_1|+|S_2|)+4$ boundary objects: the members of the ROPEreps of boundary weights. This is precisely what makes these constraints tractable in a way similar to the one-dimensional situation of eigenvectors and makes the extension theorem useful again.

The idea of having equations for the ROPErep bricks of boundary weights has various advantages. It avoids any arbitrary choice in the boundary conditions and makes possible an initial study of a new model without previous physical insight on the phase diagram. Moreover, as seen in the case of the six-vertex model, it can validate a conjectural solution by showing that it satisfies the eigen-element definition "up to morphisms".

\subsection{Computational power with given morphisms}

The definitions "up to morphisms" may scare the reader since the spaces of the invariant ROPErep with their morphisms are part of the unknowns. However, as seen in the Gaussian case and the six-vertex model, initial guesses can be obtained from the other properties (Hilbert $L^2$ space for Gaussian models, fermionic Fock spaces for the six-vertex model). Given spaces with their morphisms obtained such a structural point of view, the eigen-element definition then provides concrete equations that can be studied (and even solved in some cases). It would be very interesting to establish a full classification of the possible boundary spaces and a rigorous study of the equations within given boundary spaces.

As seen for Gaussian models \cite{BodiotSimon} and the six-vertex model \cite{SimonSixV}, it is possible to recover easily and quickly the already known results with a maximal degree of mathematical rigour.

\section{Perspectives and open questions}

The present construction paves the way to many studies in different directions. We try here to present some of them and sort them by traditional domains of mathematics, although the present paper emphasizes on a unified approach.

\subsection{The algebraic path}
The idea of present paper started with the observation of similarities between factorization algebras \cite{GinotQFTfacto} and matrix product states \cite{MPSreview}. Our care of and quest for practical numerical computations with concrete models led us to introduce definitions as close as possible to standard linear algebra with matrices but quite far from higher algebra, constructive quantum field theory and renormalization theory. 

However, it is obvious that many overlaps exist and that other tools may be imported to study our theory of eigen-elements. In particular, our guillotine partitions are quite close conceptually to the quilts in \cite{boxoperads} and our discussion on corners spaces in relation with the interchange associativity~\eqref{eq:guill2:interchangeassoc} is maybe related with the construction of corner categories in double categories \cite{mweber2015,doublecatfactosystem} or to some subset of \cite{LurieHA}.

Our only exploration of the study eigen-elements up to morphisms beyond the definition lies in Section~\ref{sec:commutuptomorph} and already shows relation with the Yang-Baxter equation. It would be interesting to explore further such tools. In particular, it would be interesting to extend our study to corner elements which constitute the real novelty w.r.t.~the one-dimensional case. Moreover, there are many traditional exercises in linear algebra around the diagonalization and it would be interesting to lift them in the present framework in order to bring new computation tools. In particular, there is yet no such thing in our case as Cayley-Hamilton theorem (which considers only eigenvalues without any eigenvectors) or the so-called \emph{lemme des noyaux}. 

In the same spirit, eigenvalues in classical associative algebras are related to the yoga of $C^*$ or von Neumann algebras whose commutative versions are known to be identified to algebras of functions (with suitable analytical properties) on their spectra. In our case, through the discussions of $\Guill_1$ and $\Guill_2$ algebras and with the commutative case described in Section~\ref{sec:eckmanhilton}, we have a clear hierarchy between commutative, associative, $\Guill_2$ and $\Guill_n$-algebras which are identified to $0$-, $1$-, $2$- and $n$-dimensional geometric objects, with the subtleties that a segment can be considered as $1$-dimensional in a direction and $0$-dimensional in the other directions. If a $\Guill_n$-algebra satisfy commutativity along $d<n$ directions together with suitable analytical properties, can it be identified to functions, characters or anything else on a set that we may call its spectrum?

\subsection{The scaling limit path: from discrete to continuous}

Extending the present formalism to a measurable state spaces is quite easy up to detailed computations of eigen-elements: describing morphisms and boundary algebras may require more advanced functional analysis (von Neumann algebras here) than just matrix algebras.

We have provided definitions in continuous space in order to show that the algebraic machinery remains the same however we have not illustrated them by examples.
Such examples are expected to be hard to define because of the distributional nature of the observables for generic models and the required renormalization procedures. However, in the regime of not-to-bad H\"older regularities, the similarity between ROPEreps and hidden Markov chains as illustrated in Proposition~\ref{prop:ROPErep:hiddenmarkov} combined with an rough path approach as in \cite{LopusanschiSimon} may give interesting starting points. In a similar spirit, the convergence of the height function of the six-vertex model to the Gaussian Free Field \cite{dominoGFF} may give interesting hints for the construction of continuous space analogues of our new tools. 

A key observation is that, with traditional trivial or periodic boundary conditions, double limits have to be taken to obtain scaling limits: the size of the domains must go to infinity and the elementary lattice scale should go to zero with renormalization of the macroscopic observables; now, considering full boundary eigen-elements removes the first limit since it corresponds already to the infinite-volume Gibbs measure and we are left with the second limit (lattice mesh going to zero). During the construction of such scaling limits, renormalization is expected to be required and it may be then very interesting to see how our structures may interact with other approaches such as \cite{Costello,CostelloGwilliam} or stochastic quantization with regularity structures.

\subsection{The computational path: solving models of statistical mechanics.} Knowing whether our new operadic approach may be efficient or not and acquiring intuition about it require essentially to solve models with this approach. Two works in progress on the general six-vertex model \cite{SimonSixV} and Gaussian models \cite{BodiotSimon} tend to indicate that finding full boundary eigen-elements is a feasible task with interesting probabilistic interpretations in many cases.

We would like to mention the special case of integrable systems, already studied throughout a huge algebraic literature. Since most integrable models of lattice statistical mechanics are Markovian, binding the two approaches should be fruitful. In particular, the standard approach to Yang-Baxter model with the transfer matrix method requires the symmetry breaking between the two directions: one direction carries the periodic boundary conditions or factorized one and the transfer matrix acts on the states on the other directions. At the level of the R-matrix itself, both directions play the same role. There must then be an intermediate level where computations of boundary eigen-elements may be performed. This is indeed a task in progress \cite{SimonSixV}.

\subsection{The barrier of undecidability.} We refer the reader to \cite{undecidabilityphysicsreview} for a wider description of the situation of undecidability in physics. However, by considering the logarithm of the face weights as (inverse temperature times) Hamiltonian acting on nearest neighbours, the study of invariant boundary weights is close to the study of translation-invariant fundamental states on Hamiltonians on the lattice (and coincide in the zero temperature limit). However this last problem is known to undecidable (in the Turing sense). We may thus worry about the decidability of ROPEreps of invariant boundary weights. In particular, we have seen that the boundary spaces are themselves part of the unknowns and infinite-dimensional in generic cases: this is an argument in favour of undecidability and it would be worth looking in detail at this question. We have seen however that in the Gaussian case as well as in the six-vertex case, there is no obstacle in the mathematical description of invariant boundary weights and thus undecidability is also not a strong obstacle to concrete computations and rigorous results.

\subsection{A combinatorial and probabilistic path: Perron-Frobenius property and enumeration of objects} In dimension one, existence of left and right eigenvectors with positive coefficients and real positive eigenvalues for matrices with positive entries is ensured by Perron-Frobenius theorem. Such eigenvectors also have nice probabilistic representations such as occupation times under excursion measures. On the other hand, in dimension larger than $2$, infinite-volume Gibbs measures are known to exist generically for our models. In between, the present paper introduces half-strip and corner eigen-elements up to morphisms from an algebraic perspective only: it would be very interesting to complete the picture by probabilistic representations of these elements, in the same way as, in dimension one, invariant measures are represented by local times under excursion laws. This may also produce interesting probabilistic "positive" representations useful for asymptotics: we think for example to \cite{CorteelWilliams} or \cite{ASEPBrownianExc} in the case of the ASEP Matrix Ansatz \cite{DEHP}.

\subsection{The analytical and numerical path: approximations and analytical bounds.} DMRG techniques \cite{DMRGMPS} show that matrix product states have very good approximation properties of ground states of Hamiltonians, together with efficient numerical techniques of linear algebra. On the other hand, in dimension $1$, the spectral radius of a matrix $A$ with positive entries is related to its Perron-Frobenius eigenvalue $\lambda_1$ and can also be obtained as a supremum of $\norm{ Ax}$ over the unit sphere $\norm{x}=1$, opening the way to bounds on $\lambda_1$ by suitable choices of vectors $x$. We would like to know whether such spectral radius approaches can be extended to full boundary eigen-elements to obtain bounds on the eigenvalue without computing all the eigen-elements.

In the same spirit, there are various iterative algorithms for converging estimations of leading eigenvalues and eigenvectors of a matrix $A$. The simplest of them consists in evaluating $u_{n+1}=Au_n/\norm{A u_n}$ for large $n$ to obtain $\lambda_1$. Can similar algorithms be developed for the estimation of full boundary eigen-elements? This question is non-trivial due to the identification of elements only up to morphisms. This would be particularly relevant in dimension larger than $3$ for which nearly no exact computational tools exist, even at the critical points.

\subsection{Operadic variations around the geometry}

\subsubsection{Relaxing the guillotine cut restriction.} When comparing sections~\ref{sec:proba} and \ref{sec:operad}, the striking difference is the restriction to guillotine cuts when partitioning domains. On one hand side, this restriction leads to an interesting mix of one-dimensional products with suitable associativities, which helps to understand to the boundary structures on rectangles. On the second hand side, this restriction is obviously only a first step towards a general theory since guillotine cuts do not play any role on the probabilistic side. It is not difficult to build a coloured operad associated to general partitions of domains as in Section~\ref{sec:proba}: the price to pay is essentially the definition of a good colour palette. Two tasks are much more difficult: first, the definition of an extended operad to deal with the boundaries and with general cuts of the whole plane and, second, the definition of good structure on the colour palette to obtain areas, perimeters and curvatures in the "surface powers" of eigenvalues.

\subsubsection{Varying the lattices.} A closely related question to the previous one is to formulate the theory for other lattices. For example, switching to the triangular lattice is quite easy since it is enough to consider guillotine cuts in three directions with $2\pi/3$-angles between them instead of two directions for the square lattice. The zoology of elementary shapes of figure~\ref{fig:admissiblepatterns} is richer but tractable; in particular, various gluings can be interpreted as $\Guill_1$-algebras with corner "bimodules" and areas and lengths are easy to define. Switching to other tilings of the plane, such as isoradial graphs considered in integrable statistical mechanics \cite{isoradial} for example, looks as difficult as relaxing the guillotine cut restriction on the square lattice.

\backmatter

\tableofcontents

\bibliographystyle{plain}
\bibliography{biblio_markov2D}

\begin{thebibliography}{10}

\bibitem{AlcarazLazo}
Francisco~C. Alcaraz and Matheus~J. Lazo.
\newblock Exact solutions of exactly integrable quantum chains by a matrix
  product ansatz.
\newblock {\em Journal of Physics A: Mathematical and General}, 37(14), 2004.

\bibitem{BagherzadehBremner}
F.~Bagherzadeh and M.~Bremner.
\newblock Commutativity in double interchange semigroups.
\newblock {\em Appl. Categor. Struct.}, 26:1185--1210, 2018.

\bibitem{BaxterCTM_2007}
R~J Baxter.
\newblock Corner transfer matrices in statistical mechanics.
\newblock {\em Journal of Physics A: Mathematical and Theoretical},
  40(42):12577--12588, oct 2007.

\bibitem{BaxterCTM_1981}
R.J. Baxter.
\newblock Corner transfer matrices.
\newblock {\em Physica A: Statistical Mechanics and its Applications},
  106(1-2):18--27, mar 1981.

\bibitem{baxterbook}
R.J. Baxter.
\newblock {\em Exactly Solved Models in Statistical Mechanics}.
\newblock Dover books on physics. Dover Publications, 2007.

\bibitem{Bethe}
H.~A. Bethe.
\newblock Zur {T}heorie der {M}etalle. i. {E}igenwerte und {E}igenfunktionen
  der linearen {A}tomkette.
\newblock {\em Zeit. f\"ur Physik}, 71:205, 1931.

\bibitem{BlytheEvans}
Richard~A. Blythe and Martin~R. Evans.
\newblock Nonequilibrium steady states of matrix-product form: a solver's
  guide.
\newblock {\em Journal of Physics A: Mathematical and Theoretical},
  40(46):R333--R441, 2007.

\bibitem{BodiotSimon}
Emilien Bodiot and Damien Simon.
\newblock Operadic structure of boundary conditions for two-dimensional markov
  gaussian random fields on the lattice.
\newblock https://arxiv.org/abs/2312.07230, 2024.

\bibitem{isoradial}
C{\'e}dric Boutillier and B{\'e}atrice de~Tili{\`e}re.
\newblock Statistical mechanics on isoradial graphs.
\newblock In Jean-Dominique Deuschel, Barbara Gentz, Wolfgang K{\"o}nig, Max
  von Renesse, Michael Scheutzow, and Uwe Schmock, editors, {\em Probability in
  Complex Physical Systems}, pages 491--512, Berlin, Heidelberg, 2012. Springer
  Berlin Heidelberg.

\bibitem{CorteelWilliams}
Sylvie Corteel and Lauren~K. Williams.
\newblock Staircase tableaux, the asymmetric exclusion process, and
  {Askey-Wilson} polynomials.
\newblock {\em Proceedings of the National Academy of Sciences},
  107(15):6726--6730, 2010.

\bibitem{Costello}
Kevin Costello.
\newblock {\em Renormalization and Effective Field Theory}, volume 170 of {\em
  Mathematical Surveys and Monographs}.
\newblock Amer. Math. Society, 2011.

\bibitem{CostelloGwilliam}
Kevin Costello and Owen Gwilliam.
\newblock {\em Factorization algebras in Quantum Field Theory, volumes {I and
  II}}.
\newblock New mathematical monographs. Cambridge University Press, 2016--2018.

\bibitem{vanicatZF}
N~Crampe, K~Mallick, E~Ragoucy, and M~Vanicat.
\newblock Open two-species exclusion processes with integrable boundaries.
\newblock {\em Journal of Physics A: Mathematical and Theoretical},
  48(17):175002, apr 2015.

\bibitem{CrampeRagoucySimon}
N~Crampe, E~Ragoucy, and D~Simon.
\newblock Matrix coordinate {Bethe Ansatz}: applications to {XXZ} and {ASEP}
  models.
\newblock {\em Journal of Physics A: Mathematical and Theoretical},
  44(40):405003, sep 2011.

\bibitem{ASEPBrownianExc}
Bernard Derrida, Camille Enaud, and Joel~L. Lebowitz.
\newblock The asymmetric exclusion process and {Brownian} excursions.
\newblock {\em J. Stat. Phys.}, 115:365--382, 2004.

\bibitem{DEHP}
Bernard Derrida, Martin~R. Evans, Vincent Hakim, and Vincent Pasquier.
\newblock Exact solution of a {1D} asymmetric exclusion model using a matrix
  formulation.
\newblock {\em Journal of Physics A: Mathematical and General},
  26(7):1493--1517, 1993.

\bibitem{LargeDevDensityASEP}
Bernard Derrida, Joel~L. Lebowitz, and Eugene~R. Speer.
\newblock Large deviation of the density drofile in the steady state of the
  symmetric simple exclusion process.
\newblock {\em J. Stat. Phys.}, 107:599--634, 2002.

\bibitem{DuminilBethe}
H.~Duminil-Copin, M.~Gagnebin, M.~Harel, I.~Manolescu, and V.~Tassion.
\newblock {The {Bethe Ansatz} for the six-vertex and {XXZ} models: An
  exposition}.
\newblock {\em Probability Surveys}, 15:102 -- 130, 2018.

\bibitem{DuminilCondensation}
Hugo Duminil-Copin, Karol~Kajetan Kozlowski, Dmitry Krachun, Ioan Manolescu,
  and Tatiana Tikhonovskaia.
\newblock On the six-vertex model's free energy.
\newblock {\em Communications in Mathematical PHysics}, 395(3):1383--1430, sep
  2022.

\bibitem{eckmannhilton}
B.~Eckmann and P.J. Hilton.
\newblock Group-like structures in general categories i: multiplications and
  comultiplicatins.
\newblock {\em Math. Ann.}, 145:227--255, 1962.

\bibitem{ZFFaddeev}
L.~D. Faddeev.
\newblock Quantum completely integrable models in field theory.
\newblock {\em Sov. Sci. Rev. C}, 1:107, 1980.

\bibitem{fressemoduleoveroperad}
B.~Fresse.
\newblock {\em Modules over operads and functors}, volume 1967 of {\em Lecture
  Notes in Mathematics}.
\newblock Springer-Verlag, 2009.

\bibitem{GibbsVelenik}
Sacha Friedli and Yvan Velenik.
\newblock {\em Statistical Mechanics of Lattice Systems}.
\newblock Cambridge University Press, 2017.

\bibitem{GibbsGeorgii}
Hans-Otto Georgii.
\newblock {\em Gibbs Measures and Phase Transitions}.
\newblock Studies in Mathematics. de Gruyter, second edition, 2011.

\bibitem{GinotQFTfacto}
Gr{\'e}gory Ginot.
\newblock {\em Notes on Factorization Algebras, Factorization Homology and
  Applications}, pages 429--552.
\newblock Springer International Publishing, Cham, 2015.

\bibitem{GolinelliMallick}
O~Golinelli and K~Mallick.
\newblock Derivation of a matrix product representation for the asymmetric
  exclusion process from the algebraic bethe ansatz.
\newblock {\em Journal of Physics A: Mathematical and General},
  39(34):10647--10658, aug 2006.

\bibitem{MPSreview}
Jutho Haegeman and Frank Verstraete.
\newblock Diagonalizing transfer matrices and matrix product operators: A
  medley of exact and computational methods.
\newblock {\em Annual Review of Condensed Matter Physics}, 8(1):355--406, 2017.

\bibitem{JIMBO_1989}
Michio Jimbo.
\newblock Introduction to the {Yang-Baxter} equation.
\newblock {\em International Journal of Modern Physics A}, 04(15):3759--3777,
  sep 1989.

\bibitem{dominoGFF}
Richard Kenyon.
\newblock Conformal invariance of domino tiling.
\newblock {\em The Annals of Probability}, 28(2), apr 2000.

\bibitem{Kock}
J.~Kock.
\newblock Note on commutativity in double semigroups and two-fold monoidal
  categories.
\newblock {\em J. Homot. Relat. Struct.}, 2(2):217--228, 2007.

\bibitem{Tetrahedron}
Atsuo Kuniba and Vincent Pasquier.
\newblock Matrix product solutions to the reflection equation from three
  dimensional integrability.
\newblock {\em Journal of Physics A: Mathematical and Theoretical},
  51(25):255204, may 2018.

\bibitem{darrick2}
Darrick Lee.
\newblock The surface signature and rough surfaces, 2024.
\newblock \url{https://arxiv.org/abs/2406.16857}.

\bibitem{darrick1}
Darrick Lee and Harald Oberhauser.
\newblock Random surfaces and higher algebra, 2024.
\newblock \url{https://arxiv.org/abs/2311.08366}.

\bibitem{LiebIce}
Elliott~H. Lieb.
\newblock Residual entropy of square ice.
\newblock {\em Phys. Rev.}, 162:162--172, Oct 1967.

\bibitem{LopusanschiSimon}
Olga Lopusanschi and Damien Simon.
\newblock Area anomaly in the rough path {Brownian} scaling limit of hidden
  markov walks.
\newblock {\em Bernoulli}, 26(4), nov 2020.

\bibitem{LurieHA}
Jacob Lurie.
\newblock Higher algebra.
\newblock \url{https://www.math.ias.edu/~lurie/papers/HA.pdf}, 2017.

\bibitem{stasheffoperad}
Martin Markl, Steve Shnider, and James~D. Stasheff.
\newblock {\em Operads in algebra, topology and physics}, volume~96 of {\em
  Mathematical Surveys and Monographs}.
\newblock American Mathematical Society, 2002.

\bibitem{mayoperad}
J.P. May.
\newblock {\em The geometry of iterated loop spaces}, volume 271 of {\em
  Lecture Notes in Mathematics}.
\newblock Springer-Verlag, 1972.

\bibitem{Nishino_1996}
Tomotoshi Nishino and Kouichi Okunishi.
\newblock {Corner Transfer Matrix Renormalization Group Method}.
\newblock {\em Journal of the Physical Society of Japan}, 65(4):891--894, apr
  1996.

\bibitem{Norris}
J.~R. Norris.
\newblock {\em Markov Chains}.
\newblock Cambridge Series in Statistical and Probabilistic Mathematics.
  Cambridge University Press, 1997.

\bibitem{ReshetikhinSixV}
Nikolai Reshetikhin.
\newblock Lectures on the integrability of the 6-vertex model.
\newblock In {\em Exact Methods in Low-Dimensional Statistical Physics and
  Quantum Computing}, pages 197--266. Oxford University Press, 2010.

\bibitem{DMRGMPS}
Ulrich Schollw\"ock.
\newblock The density-matrix renormalization group in the age of matrix product
  states.
\newblock {\em Annals of Physics}, 326(1):96--192, 2011.

\bibitem{SimonMixedRadix}
Damien Simon.
\newblock Mixed radix numeration bases: H\"orner's rule, yang-baxter equation
  and furstenberg's conjecture.
\newblock \url{https://arxiv.org/abs/2405.19798}, 2024.

\bibitem{SimonSixV}
Damien Simon.
\newblock Operadic approach to the six-vertex model: solving the model directly
  in the thermodynamic limit.
\newblock {\em in preparation}, 2024.

\bibitem{boxoperads}
Hoang~Dinh Van, Lander Hermans, and Wendy Lowen.
\newblock Box operads and higher {Gerstenhaber} brackets.
\newblock \url{https://arxiv.org/abs/2305.20036}, 2023.

\bibitem{VoronovSwissCheese}
A.A. Voronov.
\newblock The {Swiss} cheese operad.
\newblock In {\em Homotopy Invariant Algebraic Structures: A Conference in
  Honor of J. Michael Boardman : AMS Special Session on Homotopy Theory,
  January 7-10, 1998, Baltimore, MD}, Contemporary mathematics. American
  Mathematical Society, 1999.

\bibitem{mweber2015}
Mark Weber.
\newblock Internal algebra classifiers as codescent objects of crossed internal
  categories.
\newblock {\em Theory and Applications of Categories}, 30(50):1713--1792, 2015.

\bibitem{YauColoredOperads}
D.~Yau.
\newblock {\em Colored Operads}.
\newblock Graduate Studies in Mathematics. American Mathematical Society, 2016.

\bibitem{ZFZamolodchikov}
A.~B. Zamolodchikov.
\newblock Factorized {S-matrices} in two dimensions as the exact solutions of
  certain relativistic quantum field theory models.
\newblock {\em Ann. Phys.}, 120:253, 1979.

\bibitem{undecidabilityphysicsreview}
Álvaro Perales-Eceiza, Toby Cubitt, Mile Gu, David Pérez-García, and
  Michael~M. Wolf.
\newblock Undecidability in physics: a review.
\newblock \url{https://arxiv.org/abs/2410.16532}, 2024.

\bibitem{doublecatfactosystem}
Miloslav Štěpán.
\newblock Factorization systems and double categories.
\newblock \url{https://arxiv.org/abs/2305.06714}, 2023.

\end{thebibliography}

\end{document}